\documentclass[11pt,twoside,reqno,openany]{amsart}

\pdfoutput=1

\usepackage{amsbsy,amscd,amsfonts,amsmath,amssymb,amsthm,color,
fancybox,fancyhdr,footmisc,graphics,graphicx,ifthen,mathrsfs,
multicol,multirow,nccrules,pdfpages,pifont,rotating,shuffle,stmaryrd,
textcomp,times,wasysym,wrapfig}

\usepackage[dvipsnames,svgnames,x11names]{xcolor}

\usepackage[all]{xy}
\usepackage[utf8]{inputenc}
\usepackage[T1]{fontenc}
\usepackage{hyperref}

\usepackage{mathtools}
\newtagform{EngelLie}[\scriptstyle]{$}{$}
\makeatletter\newcommand{\leqnomode}{\tagsleft@true}
\newcommand{\reqnomode}{\tagsleft@false}\makeatother

\newtheorem{Theorem}[equation]{Theorem}

\newtheorem{Proposition}[equation]{Proposition}

\newtheorem{Lemma}[equation]{Lemma}

\newtheorem{Observation}[equation]{Observation}

\theoremstyle{definition}

\newtheorem{Definition}[equation]{Definition}
\newtheorem{Definition-Notation}[equation]{D\'efinition-Notation}

\newtheorem{Remark}[equation]{Remark}
\newtheorem{Principle}[equation]{Principle}
\newtheorem{Question}[equation]{Question}

\newtheorem{Problem}[equation]{Problem}
\newtheorem{SubProblem}[equation]{Sub-Problem}

\newcommand{\B}{\mathbb{B}}
\newcommand{\C}{\mathbb{C}}

\newcommand{\G}{\mathbb{G}}

\newcommand{\N}{\mathbb{N}}

\newcommand{\R}{\mathbb{R}}

\newcommand{\Z}{\mathbb{Z}}

\newcommand{\NN}{\text{\sc n}}

\newcommand{\qaux}{{\text{\usefont{T1}{qcs}{m}{sl}q}}}

\newcommand{\Iaux}{{\text{\usefont{T1}{qcs}{m}{sl}I}}}

\newcommand{\Kaux}{{\text{\usefont{T1}{qcs}{m}{sl}K}}}

\newcommand{\Paux}{{\text{\usefont{T1}{qcs}{m}{sl}P}}}

\newcommand{\Saux}{{\text{\usefont{T1}{qcs}{m}{sl}S}}}

\newcommand{\Waux}{{\text{\usefont{T1}{qcs}{m}{sl}W}}}

\definecolor{blue}{cmyk}{1.,1.,0.,0.63}
\definecolor{red}{cmyk}{0.,1.,1.,0.63}
\definecolor{green}{cmyk}{1.,0.,1.,0.63}
\definecolor{black}{cmyk}{1.,1.,1.,1.}

\newcommand{\blue}{\textcolor{blue}}
\newcommand{\green}{\textcolor{green}}
\newcommand{\red}{\textcolor{red}}

\makeatletter
\renewcommand{\@fnsymbol}[1]
{\ensuremath{\ifcase#1\or $*$\or $**$\or $***$\or $****$\or $*****$
\else\@ctrerr\fi}}
\makeatother

\newcommand{\HEAD}[2]{%
\pagestyle{fancy}
\fancyhead[RO]{\tiny\sf\thepage}
\fancyhead[CO]{{\tiny\sf #1}}
\fancyhead[LE]{\tiny\sf\thepage}
\fancyhead[CE]{{\tiny\sf #2}}
\fancyfoot{}}

\newcommand{\CITATION}[1]{\smallskip\hfill
\begin{minipage}[t]{13.5cm}\baselineskip=0.37cm\parindent=0.91cm
\blue{\scriptsize{\sf \!\!\!\!\!\!\!#1}}\end{minipage}\medskip}

\numberwithin{equation}{section}

\newcommand{\Section}[1]{
\renewcommand{\thesection}{\bf\arabic{section}}
\section{#1}
\renewcommand{\thesection}{\arabic{section}}}

\newcommand{\SectionHead}[2]{
\Section{\bf #1}
\label{#2}
\HEAD{\ref{#2}.~{\sf 
#1}}{
Julien {\sc Heyd} and 
Jo\"el {\sc Merker}, 
D\'epartement de Math\'ematiques d'Orsay, 
Universit\'e Paris-Saclay, France}}

\newcommand{\style}[1]{{\sf #1}}

\newcommand{\Aff}{\style{Aff}}

\newcommand{\dep}{\style{dep}}

\renewcommand{\det}{\style{det}}

\renewcommand{\dim}{\style{dim}}

\newcommand{\eqLF}{\style{eqLF}}

\newcommand{\eqLG}{\style{eqLG}}

\newcommand{\eqR}{\style{eqR}}

\newcommand{\eqS}{\style{eqS}}

\renewcommand{\exp}{\style{exp}}

\newcommand{\GL}{\style{GL}}

\newcommand{\Hessian}{\style{Hessian}}

\newcommand{\ind}{\style{ind}}

\newcommand{\Lie}{\style{Lie}}

\renewcommand{\lim}{\style{lim}}

\newcommand{\normal}{\style{normal}}

\newcommand{\order}{\style{order}}

\newcommand{\Proj}{\style{Proj}}

\newcommand{\rank}{\style{rank}}

\newcommand{\Saff}{\style{Saff}}

\newcommand{\SL}{\style{SL}}

\newcommand{\Span}{\style{Span}}

\newcommand{\stab}{\style{stab}}

\newcommand{\Sym}{\style{Sym}}

\newcommand{\Hall}{\Hall}

\renewcommand{\see}{{\em see}\:}

\newcommand{\smallbullet}{{\scriptscriptstyle{\bullet}}}

\newcommand{\smallsum}[1]{
\underset{#1}{\raisebox{1pt}{$\sum$\,}}
}

\newcommand{\vf}{\vfill\end{document}}

\newcommand{\zero}[1]{\underline{#1}_{\red{\circ}}}

\makeatletter
\def\hlinethick#1{%
\noalign{\ifnum0=`}\fi\hrule \@height #1 \futurelet
\reserved@a\@xhline}
\makeatother

\newlength{\myskip}
\begin{document}

\setcounter{section}{0}

\bigskip\bigskip

%\bigskip\bigskip\bigskip

\begin{center}

{\large\bf On Affinely Homogeneous Submanifolds:}

\medskip

{\large\bf The Power Series Method of Equivalence}

\label{power-series-equivalence-method}

\bigskip\bigskip

Julien {\sc Heyd}\footnotemark[1] 
and 
Jo\"el~{\sc Merker}\footnotemark[1]

\end{center}\bigskip

\footnotetext[1]{\,\,
D\'epartement de Math\'ematiques d'Orsay,
CNRS, Universit\'e Paris-Saclay, 91405 Orsay Cedex,
France, 
{\bf julien.heyd@universite-paris-saclay.fr},
{\bf joel.merker@universite-paris-saclay.fr}}

\footnotetext[1]{\,
This research was supported
in part by the Polish National Science Centre (NCN) 
via the grant number 2018/29/B/ST1/02583,
and by the Norwegian Financial Mechanism
2014--2021 via the project registration number 2019/34/H/ST1/00636.}

\begin{center}
\begin{minipage}[t]{12.5cm}
\parindent 0.53cm
\footnotesize
\noindent
{\sc Abstract}.
We determine all affinely homogeneous models for:

\smallskip\noindent$\square$\,
Surfaces $S^2 \subset \R^4$;

\smallskip\noindent
including the {\em simply transitive} models.

We employ an improved {\sl power series method of equivalence},
which captures invariants at the origin,
creates branches, and infinitesimalizes calculations.

We find several inequivalent terminal branches 
yielding each to some nonempty
moduli space of homogeneous models, 
sometimes parametrized by a certain invariant algebraic variety.

Three main features may be emphasized:

\smallskip\noindent{\bf 1.} 
Iterated {\em single-pointed} jet bundles;

\smallskip\noindent{\bf 2.} 
Cartan-enhanced power series method of equivalence;

\smallskip\noindent{\bf 3.}
Constant {\em ping-pong} between normal forms (nf)
and vector fields (vf).

\end{minipage}
\end{center}

%%%%%%%%%%%%%%%%%%%%%%%%%%%%%%%%%%%%%%%%%%%%%%%%%%%%%%%%%%%%%%%%%%%%%%
\SectionHead{Introduction}
{introduction-Sigma}
%%%%%%%%%%%%%%%%%%%%%%%%%%%%%%%%%%%%%%%%%%%%%%%%%%%%%%%%%%%%%%%%%%%%%%

Inspired by Peter Olver's theory of differential 
invariants~{\cite{Olver-2011, Olver-2018, 
Olver-Valiquette-2018}},
{\em cf.} the seminal article~{\cite{Fels-Olver-1999}} 
with Mark Fels, 
\"Orn Arnaldsson 
and Francis Valiquette~{\cite{Valiquette-2013,
Arnaldsson-Valiquette-2020}}, 
and the second author as well
joint with Zhangchi Chen~{\cite{Chen-Merker-2019, 
Chen-Merker-2020},
have endeavoured, in an affine context,
to explore {\em branches}
of geometric structures for which some (relative or absolute)
differential invariants either vanish identically, or are nonzero.

Either recursive recurrence 
formulas~{\cite{Olver-2011, Valiquette-2013}},
or explicit expressions of invariants (when accessible),
lead to "discover" certain sub-structures of a given geometry,
all having their own algebras of differential invariants.
This approach requires to work in jet bundles, as a whole.

On the other hand, although less represented in the literature,
the power series approach consists in expanding
(in power series)
at a fixed point a given geometric structure.
Roughly, once a complete normal form is attained, 
power series coefficients happen to be {\em values at the origin}
of differential invariants~{\cite{Olver-2018}}.

However, power series are essentially of no use 
in order to realize
what's happening in a branch where a collection of
(relative or absolute) differential invariants 
{\em vanish} (identically),
because any function which vanishes at one point
may certainly have nonzero values at nearby points.

So, because there is a true difficulty here,
we temporarily abandon the (ambitious) project of describing 
algebras of differential invariants in (all) sub-branches
of a given geometric structure.
Nonetheless, we hope to come back to 
this project in a near future.

Instead, the present memoir is devoted to 
2 illustrations of the power series
method of equivalence, with main 
focus on the determination of
{\em homogeneous models}.

Inspired by articles~{\cite{Eastwood-Ezhov-1999, 
Eastwood-Ezhov-2001-1, Eastwood-Ezhov-2001-2}}
of Mike Eastwood and Vladimir Ezhov in affine geometry,
we enrich the power series approach
by importing Cartan's technique of
step-by-step
{\em group reduction}.
But instead of proceeding
in jet bundles by reducing $G$-structures, 
as the celebrated Cartan method of equivalence does,
we proceed jet order by jet order,
always at a single fixed point,
and, iteratively order by order,
above fixed points or above orbit transversals 
(which sometimes are positive-dimensional).

Unexpectedly, and in contrast with Cartan's method
of equivalence, we are led at each jet order
to certain {\em linear representations of matrix Lie groups},
which are responsible for the creation of branches.
How? Just by examining, in parallel with 
the creation of normal forms at a jet order, 
the condition that a given vector field is
an infinitesimal symmetry of the geometric structure.

As a matter of fact, in the search for
homogeneous models (only), we reach in some way
a complete information about the collection of
(relative) invariants which vanish (identically,
by homogeneity) in any given branch.

This memoir has restricted scope:

\smallskip\noindent$\square$\,
affine classifications in small dimensions;

\smallskip\noindent$\square$\,
no closed form is shown for the found homogeneous models.

\smallskip

Other authors
were able to set up nice closed forms,
{\em see} {\em e.g.}~{\cite{Doubrov-Komrakov-1998,
Doubrov-Komrakov-Rabinovich-1996,
Eastwood-Ezhov-1999,
Eastwood-Ezhov-2001-1,
Eastwood-Ezhov-2001-2}}.
Is there a general criterion for the existence of nice closed forms?
We were not able to answer this question.

In this memoir, we abandon the task of setting up closed forms,
because our goal is to understand\big/to see the (pointwise) invariants
which distribute homogeneous models into inequivalent branches having
empty intersection (up to some finite group action).
Closed forms are nicer, but they do not show what their invariants
are, and often, when a closed form of a 
family of homogeneous model depends on
one or several parameters, the family "crosses" several distinct
branches.

This memoir has 2 parts, each studying affine submanifolds.
Computation files are available at~{\cite{Heyd-Merker-2023}}.

Part~I, consisting of 
Sections~{\ref{differential-invariants-general-methods}}
$\to$
{\ref{classification-invariants}},
is devoted to a general presentation of the
power series method of equivalence, 
with links with Peter Olver's approach
for determining the structures
of algebras of differential invariants
under the action of a finite-dimensional
Lie group $G$ acting on local graphed submanifolds
$\big\{ u = F(x) \big\}$.
In the spirit of Lie, 
{\em cf.}~{\cite{Engel-Lie-Merker-2015, Engel-Lie-1893}},
what we consider to be a classification must exhibit
explicit Lie algebras of vector fields, parametrized by
absolute invariants (if any).

Part~II, 
consisting of 
Sections~{\ref{S2-R4}}
$\to$
{\ref{2g-models}}, 
is devoted to surfaces $S^2 \subset \R^4$.
The creation of (non-overlapping) order 2 branches is
explained in detail.
The final classification appears at the end of the part,
with numerous (families of) homogeneous
models, in the spirit of 
Wermann's impressive Ph.D
(not published in a journal) 
on special affine homogeneous
hypersurfaces $S^3 \subset \R^4$, 
which followed the algebraic approach
of Doubrov-Komrakov-Rabinovich~{\cite{
Doubrov-Komrakov-Rabinovich-1996}}.

\medskip\noindent{\bf Acknowledgments.}
An anonymous referee is thanked for his {\em criticism}.
Authors benefited of the invaluable power of
Ian Anderson's Maple package 
{\footnotesize\sf DifferentialGeometry}.
Also, 
the power series method of equivalence owes a lot to 
several `{\sl private lectures}' of Pawe{\l} Nurowski
on the Cartan equivalence method during his many stays in Orsay.

And, above all, authors were astonished, attracted, seduced 
by Peter Olver's amazingly crystal clear monographs and articles.

%%%%%%%%%%%%%%%%%%%%%%%%%%%%%%%%%%%%%%%%%%%%%%%%%%%%%%%%%%%%%%%%%%%%%%
\SectionHead{Differential Invariants and Homogeneous Models}
{differential-invariants-general-methods}
%%%%%%%%%%%%%%%%%%%%%%%%%%%%%%%%%%%%%%%%%%%%%%%%%%%%%%%%%%%%%%%%%%%%%%

Consider a Lie group $G$ acting on a given 
type of geometric structure.
Examples are: 
Euclidean, affine, conformal, projective,
(pseudo-)Riemannian, symplectic, quaternionic, 
Cauchy-Riemann (CR), para-CR, 
\dots,
structures.
Other examples are: 
ordinary differential equations;
partial differential equations;
integrability systems; 
Pfaffian systems, \dots.

In his complete works, 
\'Elie Cartan often started by re-expressing
the considered geometric structure 
as being a specific
exterior differential system. 

On the other hand, as explained in Peter Olver's 
monographs and articles,
after transfer to 
an appropriate associated space 
({\em e.g.} a jet bundle),
several (local)
geometric structures with a (local) Lie group $G$ acting on them
can be expressed as (local) {\em graphs} 
$\{ u = F(x) \}$
in the associated space
equipped with a $G$-action.

In this memoir, we adopt the graph point of view.
Although our considerations 
are valid for infinite-dimensional Lie groups,
like the groups of diffeomorphisms,
of biholomorphisms, of CR-equivalences, \dots, 
we shall restrict ourselves to the finite-dimensional setting.
We shall work over $\R$ or $\C$.

Consider therefore a Lie group $G$ of finite dimension 
$1 \leqslant r < \infty$. 
Let $n \in \N_{\geqslant 1}$ and $c \in \N_{\geqslant 1}$.
In $\R^{n+c}$ with coordinates
$x = (x_1, \dots, x_n)$ and
$u = (u_1, \dots, u_c)$, 
consider a $c$-codimensional graph:
\[
u_j
\,=\,
F_j\big(x_1,\dots,x_n\big)
\eqno
{\scriptstyle{(1\,\leqslant\,j\,\leqslant\,c)}}.
\]
Throughout, our point of view will be local,
and the $F_j$ will be assumed to be {\em analytic}.
We will {\em not} introduce notations
for open sets, subsets, sub-subsets, \dots,
{\em cf.}~{\cite[Chap.~1]{Engel-Lie-Merker-2015}}.

Let the group $G$ act on $\R^{n+c}$, by analytic
diffeomorphisms. In this memoir,
$G$ will consist of {\em affine} transformations.
Also, an element $g$ 
of the group $G$ will always be explicitly
given by 
{\sl group parameters} $(g_1, \dots, g_r) \in \R^r$. 

Two general problems are of interest, about which
we will be more specific later,
{\em see} Problems~{\ref{Pbm-homogeneous-models}}
and~{\ref{Pbm-algebras-differential-invariants}}  {\em infra}.

\begin{Problem}
\label{Pbm-debut-diff-invts}
{\sl Describe algebras of differential invariants.}
\end{Problem}

\begin{Problem}
\label{Pbm-debut-hom-models}
{\sl Determine homogeneous models.}
\end{Problem}

These two problems are tightly linked with each other,
because most of the times, 
homogeneous models of a given geometric action
are `exceptional' objects in a wide universe of {\em non}symmetric
objects. The `exceptional' symmetric objects have {\em constant}
differential invariants, while the 'general' {\em non}symmetric
objects often have 
infinitely many {\em functional} differential invariants,
which share {\em complicated} differential-algebraic relations.

The Lie-Fels-Olver {\sl recurrence relations}
between differential invariants constitute a natural
`{\sl bridge}' between these two general problems.
Indeed, the effectiveness of Peter Olver's equivariant moving frame
approach lies in the powerful {\sl recurrence relations}, 
which produce
complete and explicit {\em differential-algebraic}
structures for the underlying algebras of
differential invariants\,\,---\,\,this without requiring explicit
coordinate expressions for either
the moving frame or the invariants. Evidently, 
differential invariants of homogeneous structures are 
{\em constant}, and it is a fact that
the {\em algebraic} relations 
between them retain major part of the recurrence relations.

In this memoir, we will focus mainly on the second 
Problem~{\ref{Pbm-debut-hom-models}}.
Our objective is to develop a computationally effective
method towards the determination of homogeneous models,
which is inspired by Cartan's equivalence method
and by Lie's original way of proceeding~{\cite{Engel-Lie-Merker-2015,
Engel-Lie-1893}}, {\em cf.} also
Tresse~{\cite{Tresse-1893}}.
We will manipulate in an algebraic way (constant)
differential invariants, seen as power series coefficients.

Our point of view will be different from the one 
in~{\cite{Chen-Merker-2019, Chen-Merker-2020}},
where homogenous models were determined by stating first
the Lie-Fels-Olver recurrence relations between
differential invariants before assuming them to be constant.

None of the articles~{\cite{Eastwood-Ezhov-1999, 
Eastwood-Ezhov-2001-1, Eastwood-Ezhov-2001-2}}, 
concerned with affinely homogeneous models,
mentions the 
presence of differential invariants.

%%%%%%%%%%%%%%%%%%%%%%%%%%%%%%%%%%%%%%%%%%%%%%%%%%%%%%%%%%%%%%%%%%%%%%
\SectionHead{Fibers Over Group Transversals
Versus Full Jet Bundle}
{fibers-over-transversals-versus-full-jet-bundles}
%%%%%%%%%%%%%%%%%%%%%%%%%%%%%%%%%%%%%%%%%%%%%%%%%%%%%%%%%%%%%%%%%%%%%%

Abbreviate $z := (x, u)$.
Denote the target coordinates as $\overline{z} := 
\big( \overline{x}, \overline{u} \big)$.
An element $g \in G$ in some neighborhood of the identity
sends the graph: 
\[
M 
\,:=\,
\{
u
=
F(x)
\}
\] 
to a similar graph:
\[
\overline{M}
\,:=\,
\big\{
\overline{u}
\,=\,
\overline{F}\big(\overline{x},g\big)
\big\},
\]
with certain analytic functions $\overline{F}_j$ which
depend on the group parameters.
 
The expressions of these $\overline{F}_j (\overline{x}, g)$
are difficult to write down,
highly nonlinear, 
often cumbersome. 
They in fact require the full strength
of the implicit function theorem.

Such transformations of graphs 
appear regularly in the original
complete works of Lie~{\cite{Engel-Lie-Merker-2015}}.
The monographs~{\cite{Olver-1986, Olver-1995}}
of Peter Olver,
to which nonexpert readers are referred,
offer invaluable modernized access 
to Lie's theory. Our general presentation 
will be close 
to~{\cite{Chen-Merker-2019, Merker-2020}}.

Let us write:
\[
g\cdot
z
\,=\,
g
\cdot
(x,u)
\,=:\,
\big(\overline{x},\overline{u}\big)
\,=\,
\overline{z}.
\]

We whall assume that the group $G$ acts transitively
on $\R^{n+c}$, and even, that $G$ contains all translations.
{\small (Non-transitive group actions are 
sometimes considered in Peter Olver's articles.)}
`Morally', the fact that $G$ acts {\em transitively} implies 
that all points are somewhat `equivalent'. 

\begin{center}
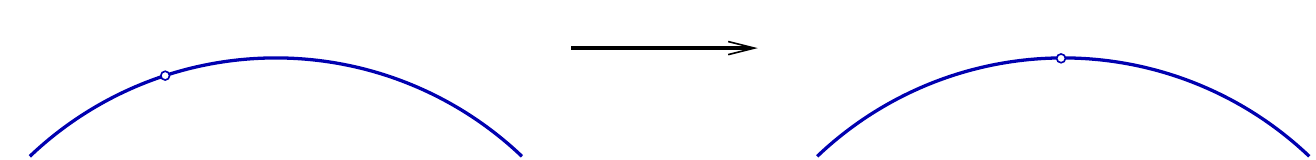
\end{center}

Therefore, any point $p_0 \in M$ can be `moved by $G$' to some
`central' point, $\overline{0} \in \R^{n+c}$,
the origin of the target 
coordinates $\overline{z}$.
Next, coordinates $z$
can be `re-centered' at $p_0$.

\begin{center}
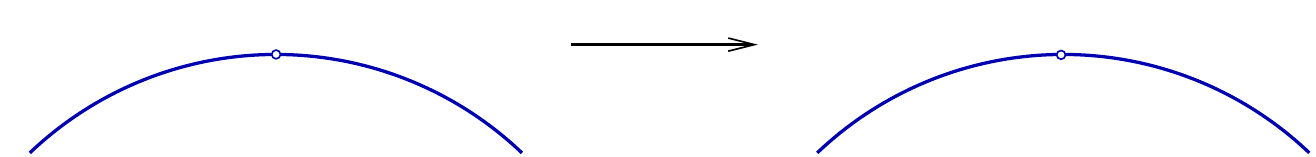
\end{center}

So both graphs $M$ and $\overline{M}$ pass through the origin.
And in fact, only the (isotropy) {\em subgroup} $G_\stab^0 \subset G$
of transformations $g \in G$ sending $0$ to $\overline{0}$
should be considered onward, as we will argue later.

To study invariants under $G$-actions
and to classify $G$-homogeneous geometries,
(roughly) two different (general) approaches exist:

\smallskip\noindent$\bullet$\,
Work within (full) jet bundles (Lie, Cartan, Olver, \dots);

\smallskip\noindent$\bullet$\,
Work with (truncated) power series centered at the origin
(Lagrange, Poincaré, Moser, \dots).

\smallskip

The second approach, 
less developed, 
has several defects.
One obvious defect is that differential invariants
of Lie type,
which require differentiation with respect to $x_1, \dots, x_n$,
cannot be computed by manipulating
power series only at $x_1 = \cdots = x_n = 0$!
Other defects will be discussed later.

The first steps of 
Lie's theory of differential invariants
consist in {\em prolongating}
the $G$-action to jet 
bundles~{\cite{Engel-Lie-Merker-2015, Olver-1986, Olver-1995}}. 
Sketching only key aspects, we will not present 
the complete 
details of~{\cite[Chaps. 3, 4, 5]{Chen-Merker-2019}}.

For a jet order $\kappa \in \N$, 
let $J_{n,c}^\kappa$ be the bundle of $\kappa$-jets
of $c$ functions of $n$ variables, at all base
points $\big(x, u(x) \big) \in M$. For instance,
$J_{n,c}^1$ has $n + c + n\,c$ 
independent coordinates corresponding
to the $x_i$, 
and to the $u_j$ together with all their first order derivatives
$u_{j,x_i}$.

\begin{center}
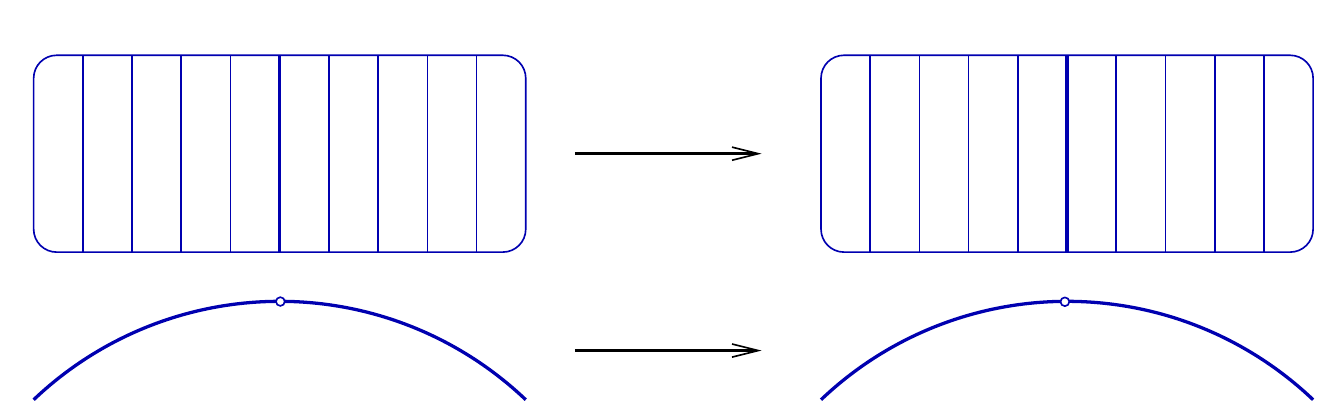
\end{center}

As is known, the $G$-action
uniquely {\em lifts} as a $G$-action
on {\em first jets} of graphs. 
This action is just the (differential) action on tangent spaces to
the two graphs at corresponding points. Illustrations
can be found in~{\cite{Chen-Merker-2019}}.
Here, we will introduce neither symbols nor downward arrows
for projections. 

Denote $z^1 = (x, u, u^1)$ and similarly
$\overline{z}^1 = \big( \overline{x}, \overline{u},
\overline{u}^1 \big)$. 
Although it is the same group $G$ that acts on $J_{n,c}^1$, 
denote its lifted action with the symbol $g^1$:
\[
g^1\cdot z^1
\,=:\,
\overline{z}^1
\eqno
{\scriptstyle{(g\,\in\,G)}}.
\] 

The $G$-action can be lifted to jet bundles $J_{n,c}^\kappa$
of any order $\kappa$. 
To represent the lifts of the $G$-action to the complete 
jet bundles of orders $1$ and $2$, 
here is a diagram.

\begin{center}
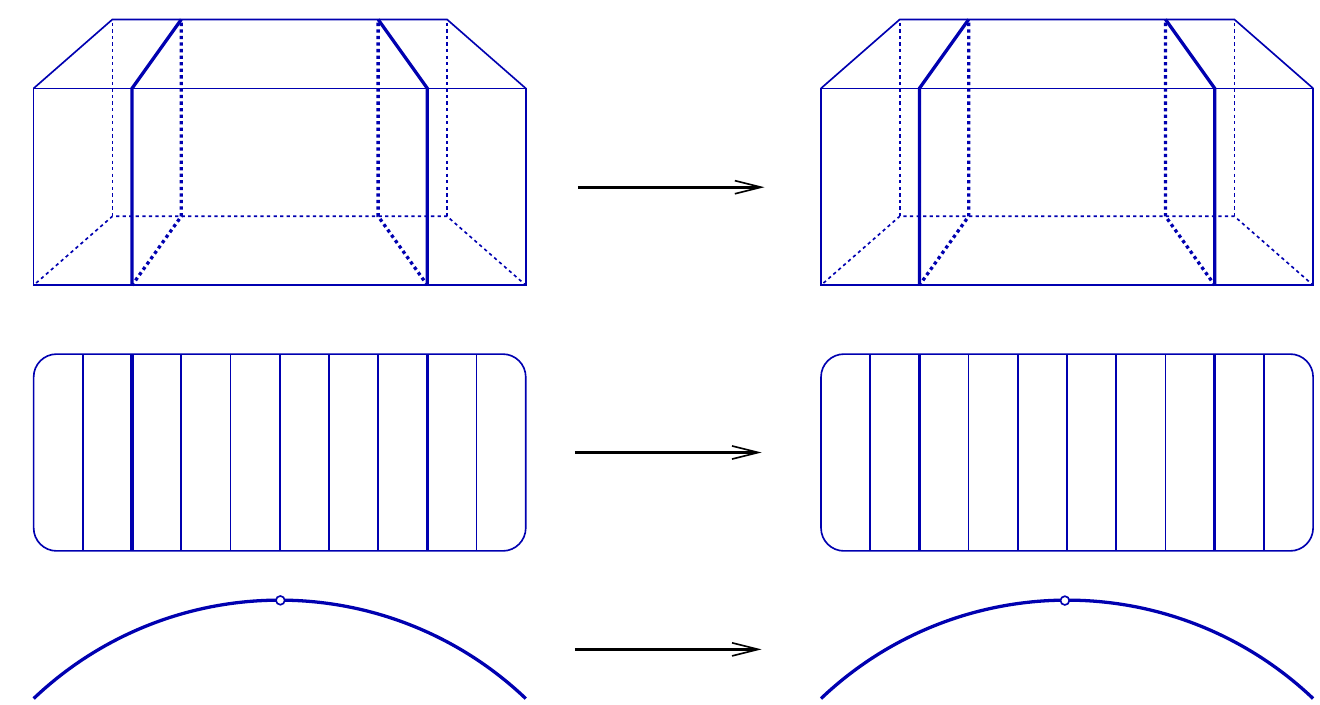
\end{center}

Roughly, differential invariants 
of order $\kappa$ are functions on $J_{n,c}^\kappa$
that are invariant under the lifted $G$-action.
{\em High difficulties are caused
by the unavoidable orbit dimension jumps of the $G$-actions on
$J_{n,c}^1$, on $J_{n,c}^2$, \dots, on $J_{n,c}^\kappa$, \dots}. 

General theoretical aspects are presented 
in the references mentioned above,
in Fels-Olver's seminal article~{\cite{Fels-Olver-1999}},
and in several other articles~{\cite{Olver-2007,
Olver-2011, 
Olver-2018,
Olver-Valiquette-2018}} 
of Peter Olver and his collaborators, 
\'Evelyne Hubert, 
Irina Kogan, 
Francis Valiquette, 
Masoud Sabzevari, \dots.

Our principal objective is to explain how some aspects 
of Lie's theory
of differential invariants can be `exported' (in part)
to the power series domain. 
Said differently, we aim at showing that 
Lagrange's foundational claim~{\cite{Lagrange-1797}}
that functions are plain power series,
although being false in general, 
is in fact very relevant and very adequate in this field.

In addition, an astonishing outcome occurs.
Namely that with power series (only), 
the computational 
exploration of algebras of differential
invariants happens to be more accessible
and more efficiently performed than ever before.
This we can assert
because we practiced intensive computations
on Lie-type jet bundles, 
and also because we practiced, 
joint with Pocchiola, Sabzevari, Foo, Nurowski,
the {\em parametric}
Cartan method of equivalence~{\cite{Merker-Pocchiola-2018, 
%%%Merker-2021,
Merker-Nurowski-2020-a, 
%%%Merker-Nurowski-2020-b,
Merker-Nurowski-2023}}.

So let us explain why and how power series fit.
Everything is rather simple and natural. 
An approach similar to the one developed in this memoir
appears in the article~{\cite{Merker-2020}},
devoted to a simplified construction of Moser 
normal forms and chains 
in CR geometry.

Unfortunately, in present times,
Lie's original works are not much read, not even
studied.
But the mathematical ideas and the geometric aspects 
that we are starting to present here are
essentially contained in 
Friedrich Engel and Sophus Lie's 
monographs~{\cite{Engel-Lie-Merker-2015, Engel-Lie-1893}}, 
and also in the memoir~{\cite{Tresse-1893}}
of Arthur Tresse.

\begin{center}
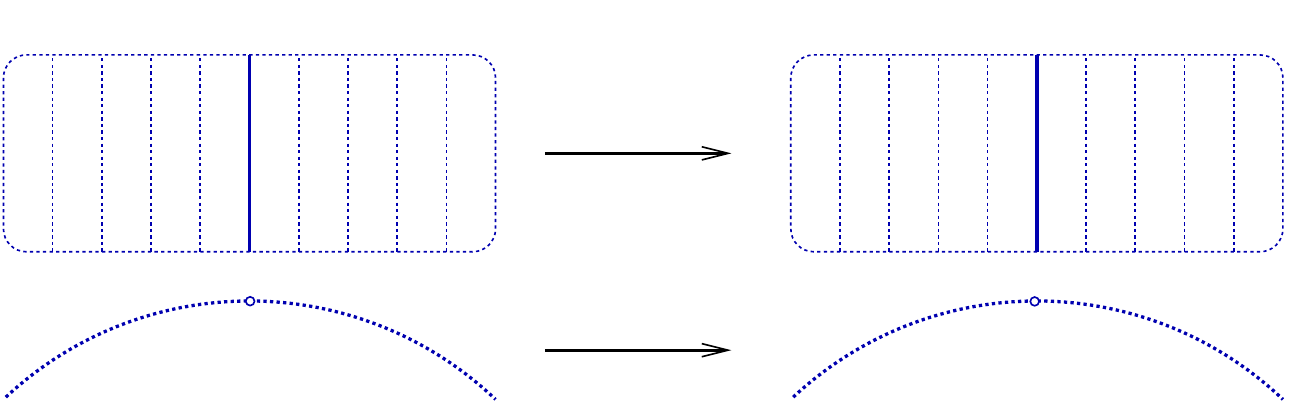
\end{center}

First of all, since the origin $0 \in M$ is sent to
the origin $\overline{0} \in \overline{M}$, 
the group action sends the first jet fiber 
$J_{n,c}^1 \big\vert_0$ 
over $0$ to the first jet fiber 
$\overline{J}_{n,c}^1 \big\vert_{\overline{0}}$
over $\overline{0}$. Of course,
we are considering only
group elements $g$ 
of the {\em subgroup} $G_\stab^0 \subset G$
fixing the origin, which we denote by $g\vert_0$.

\begin{center}
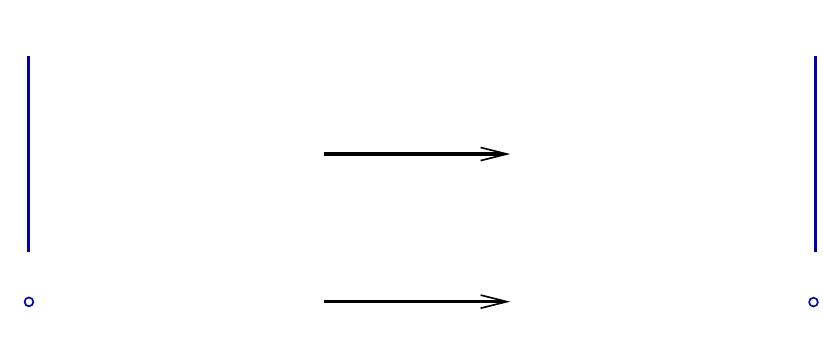
\end{center}

As a key decision here, we decide to {\em forget} 
other jet fibers (!).
Full bundles will not anymore be dealt with (!). 
When passing to higher jet orders, this decision of
restricting to selected fibers will be iterated.
Of course, there are prolongations $\big( g\vert_0 \big)^1$, 
$\big( g\vert_0 \big)^2$, \dots, to jet fibers
$J^1 \big\vert_0$, $J^2 \big\vert_0$, \dots,
and we will later show {\em formulas} for such prolongations,
which are simpler thatn the formulas in the full
jet bundle~{\cite{Olver-1986, Bluman-Kumei-1989, Olver-1995,
Merker-2008, Chen-Merker-2019}}

\begin{center}
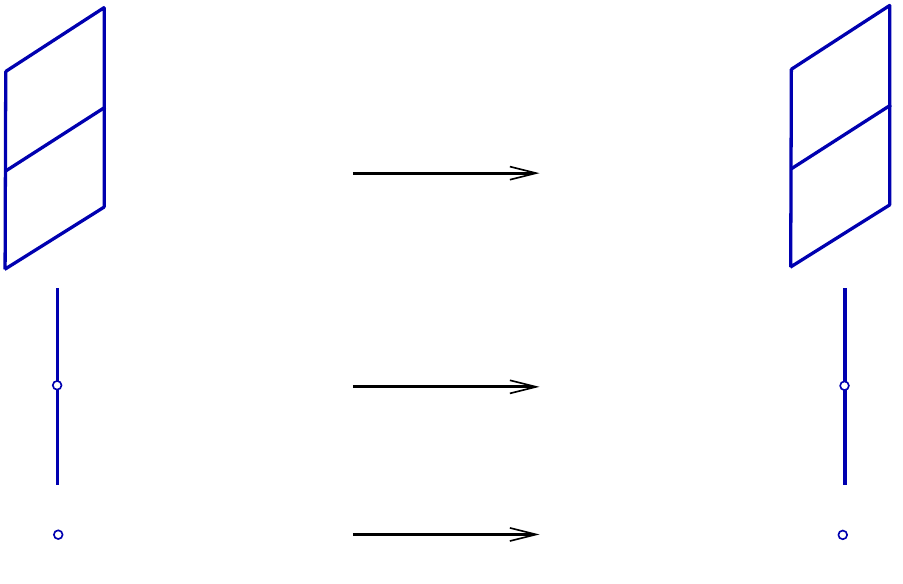
\end{center}

Several groups $G$, as {\em e.g.}
the affine or projective groups,
contain not only translations but also 
{\sl transvections}~{\cite[Sec.~14.3]{Chen-Merker-2019}},
namely maps of the form:
\[
v_j
\,=\,
u_j
+
\qaux_{j,1}\,x_1
+\cdots+
\qaux_{j,x}\,x_n
\eqno
{\scriptstyle{(1\,\leqslant\,j\,\leqslant\,c)}},
\]
with arbitrary $\qaux_{j,i} \in \R$.
Such maps enable to `strengthen' tangent spaces
of both $M$ at $0$ and $\overline{M}$ at $\overline{0}$
to be `horizontal', that is, to normalize to zero
all first order terms in the power series expansions:
\[
u
\,=\,
0
+
{\rm O}_{x_1,\dots,x_n}(2)
\ \ \ \ \ \ \ \ \ \ \ \ \ \ \ \ \ \ \ \
\text{and}
\ \ \ \ \ \ \ \ \ \ \ \ \ \ \ \ \ \ \ \
\overline{u}
\,=\,
0
+
{\rm O}_{\overline{x}_1,\dots,\overline{x}_n}(2),
\]
where of course:
\[
{\rm O}_{x_1,\dots,x_n}(2)
\,=\,
\sum_{i_1+\cdots+i_n\geqslant2}\,
x_1^{i_1}\cdots x_n^{i_n}\,
F_{i_1,\dots,i_n}.
\]

Precise formulas and
normalization equations can easily be written,
{\em cf.}~{\cite[Sec.~2]{Merker-2022}} for $G = \Aff(\R^{n+1})$.
Geometrically, this means that the $G$-action lifted
to the first jet bundle $J^1$ and {\em restricted} 
to its fiber $J^1 \big\vert_0$ over the origin $0$ only,
is {\em transitive}, and this means that the
{\em origin} $0^1 \in J^1 \big\vert_0$ is taken
as a {\em transversal} to the unique $G_\stab^0$-orbit
in $J^1 \big\vert_0$.   

\begin{center}
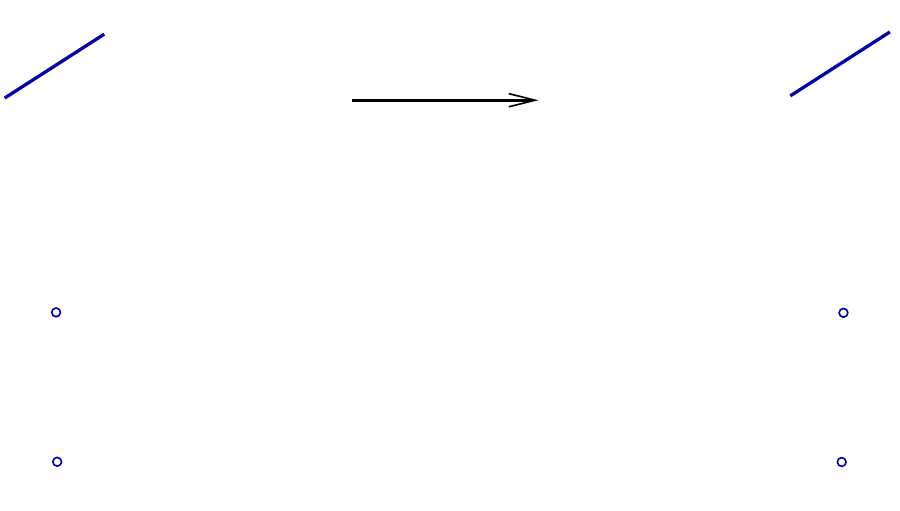
\end{center}

Therefore, not the whole second order jet fibers $J^2 \big\vert_0$ and
$\overline{J}^2 \big\vert_{\overline{0}}$ over the origins $0 \in M$
and $\overline{0} \in \overline{M}$ should be dealt with.
Instead, and precisely as it is drawn in the simplified diagram above,
one should consider only:

\smallskip\noindent$\bullet$\,
$J^2 \big\vert_{0^1} :=$ the part of $J^2$ 
over the origin $0^1$ of $J^1\big\vert_0$;

\smallskip\noindent$\bullet$\,
$\overline{J}^2 \big\vert_{\overline{0}^1} :=$
the part of $\overline{J}^2$ 
over the origin $\overline{0}^1$ of
$\overline{J}^1 \big\vert_{\overline{0}}$.

\smallskip\noindent
These two {\em smaller} subspaces 
are the respective two preimages of $0^1$ and of $\overline{0}^1$
under the (unwritten) projections
from the second floor to the first floor.

\smallskip

Furthermore, only the {\em subgroup} $G_\stab^1 \subset G_\stab^0
\subset G$ of transformations sending $0^1$ to $\overline{0}^1$ 
(hence sending $0$ to $\overline{0}$) 
should be dealt with. As in the figure above,
let us denote by $g^2 \big\vert_{0^1}$ 
the prolongation to $J^2 \big\vert_{0^1}$ of
group elements $g$ belonging to $G_\stab^1$.

Thus, exactly as in Cartan's method of equivalence,
there are here successive {\sl group reductions}.
None of the articles~{\cite{Eastwood-Ezhov-1999, 
Eastwood-Ezhov-2001-1, Eastwood-Ezhov-2001-2}} 
mentions the ubiquitous (and universal)
presence of group reduction.

\begin{center}
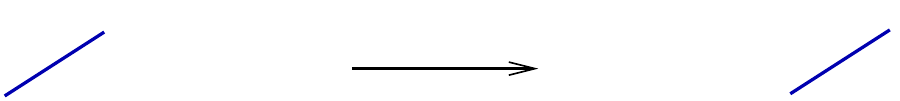
\end{center}

So again, there is an action on selected (reduced) fibers.
And again, the concerned fiber must be decomposed into
group orbits. 
Theorem~42 of Lie in~{\cite{Engel-Lie-Merker-2015}}\,\,---\,\,probably
the most complicated statement of the whole 
Volume~I of {\em Theorie der 
Transformationsgruppen}\,\,---\,\,explains 
in an algorithmic way how to decompose group actions into orbits,
applying an infinitesimal technique.

In the present memoir, 
this method of 
restricting to fibers, 
of decomposing
into reduced-group orbits,
and of finding appropriate transversals,
will be several times 
iterated in higher jet levels. 

As we intend to avoid
introducing sophisticated notation 
and to avoid showing a fully detailed
general theory, let us treat a simple example. 
The simpler, the better.

%%%%%%%%%%%%%%%%%%%%%%%%%%%%%%%%%%%%%%%%%%%%%%%%%%%%%%%%%%%%%%%%%%%%%%
\SectionHead{Example: Parabolic Surfaces $S^2 \subset \C^3$}
{parabolic-surfaces-S2-C3}
%%%%%%%%%%%%%%%%%%%%%%%%%%%%%%%%%%%%%%%%%%%%%%%%%%%%%%%%%%%%%%%%%%%%%%

With $G := \Aff (\C^3)$, 
in the left space, 
let $S^2 \subset \C^3 \ni (x ,y, u)$ be a graphed (analytic) surface:
\[
u
\,=\,
F(x,y)
\,=\,
0
+
0
+
F_{2,0}\,x^2
+
F_{1,1}\,x\,y
+
F_{0,2}\,y^2
+
{\rm O}_{x,y}(3),
\]
its 
constant term $0$ and its first order term $0$ being already 
normalized. Of course: 
\[
{\rm O}_{x,y}(3)
\,=\,
\sum_{i+j\geqslant3}\,
F_{i,j}\,x^i\,y^j.
\]
Clearly, $J^2 \big\vert_{0^1}$ is coordinatized by
$\big(F_{2,0}, F_{1,1}, F_{0,2} \big)$.

In the right space, let the target surface in $\C^3 \ni
(p, q, v)$ be similarly graphed as:
\[
v
\,=\,
G(p,q)
\,=\,
0
+
0
+
G_{2,0}\,p^2
+
G_{1,1}\,p\,q
+
G_{0,2}\,q^2
+
{\rm O}_{p,q}(3),
\]
with $\big( G_{2,0}, G_{1,1}, G_{0,2} \big)$ being coordinates
on $\overline{J}^2 \big\vert_{\overline{0}^1}$.

A general transformation of $\Aff(\C^3)$ writes:
\[
\aligned
p
&
\,:=\,
a_{1,1}\,x+a_{1,2}\,y+b_1\,u+\tau_1,
\\
q
&
\,:=\,
a_{2,1}\,x+a_{2,2}\,y+b_2\,u+\tau_2,
\\
v
&
\,:=\,
c_1\,x+c_2\,y+d\,u+\sigma,
\endaligned
\ \ \ \ \ \ \ \ \ \ \ \ \ \ \ \ \ \ \ \
\text{with}
\ \ \ \ \ \ \ \ \ \ \ \ \ \ \ \ \ \ \ \
0
\,\neq\,
\left\vert\!
\begin{array}{ccc}
a_{1,1} & a_{1,2} & b_1
\\
a_{2,1} & a_{2,2} & b_2
\\
c_1 & c_2 & d
\end{array}
\!\right\vert.
\]
But $0$ should be sent to $\overline{0}$, which holds if
and only if all translational parameters
$\tau_1 = \tau_2 = \sigma = 0$ vanish,
so that the transformation belongs to $\GL(\C^3)$:
\[
\left[
\begin{array}{c}
p \\ q \\ v
\end{array}
\right]
\,=\,
\left[\,
\begin{array}{ccc}
a_{1,1} & a_{1,2} & b_1
\\
a_{2,1} & a_{2,1} & b_1
\\
c_1 & c_2 & d
\end{array}
\,\right]
\left[
\begin{array}{c}
x \\ y \\ u
\end{array}
\right].
\]
Thus $G_\stab^0 = \GL(\C^3)$ here.

Furthermore, $0^1$ should be sent to 
$\overline{0}^1$, 
and the reader can verify that this 
corresponds to the group reduction
towards $G_\stab^1$:
\[
\left[\,
\begin{array}{ccc}
a_{1,1} & a_{1,2} & b_1
\\
a_{2,1} & a_{2,2} & b_2
\\
c_1 & c_2 & d
\end{array}
\,\right]^{\green{\bf 0}}
\ \ \ \ \
\leadsto
\ \ \ \ \
\left[\,
\begin{array}{ccc}
a_{1,1} & a_{1,2} & b_1
\\
a_{2,1} & a_{2,2} & b_2
\\
\red{\bf 0} & \red{\bf 0} & d
\end{array}
\,\right]^{\green{\bf 1}}.
\]
How? Simply by inspecting the {\sl fundamental equation:}
\leqnomode\usetagform{default}
\begin{align}
\label{intro-fund-eq-S2-C3}
0
&
\,\equiv\,
-\,c_1\,x
-
c_2\,y
-
d\,F(x,y)
\notag
\\
&
\ \ \ \ \ \
+
G\Big(
a_{1,1}\,x+a_{1,2}\,y+b_1\,F(x,y),\,\,\,
a_{2,1}\,x+a_{2,2}\,y+b_2\,F(x,y)
\Big),
\end{align}
which expresses that $\big\{ u = F(x,y) \big\}$ is mapped
to $\big\{ v = G(p,q) \big\}$.
This fundamental equation must hold {\em identically} in the ring
$\C\{x,y\}$ of convergent power series. Thus,
neglecting second and higher order terms:
\[
0
\,\equiv\,
-\,c_1\,x
-
c_2\,y
+
{\rm O}_{x,y}(2),
\]
we see that $0 = c_1 = c_2$, necessarily.
Visibly, in $G_\stab^1$, there remain 7 (isotropy) parameters.

And now, what is the action of $G_\stab^1$ 
on $J^2\big\vert_{0^1}$? How to prolong $G_\stab^1$
to second order jets?
Simply by looking at second order terms in the
fundamental equation! By hand or using a computer, we find:
\leqnomode\usetagform{default}
\begin{align}
\label{eq-x2-xy-y2}
0
&
\,\equiv\,
x^2\,
\Big[
a_{2,1}^2\,G_{0,2}
+
a_{1,1}\,a_{2,1}\,G_{1,1}
+
a_{1,1}^2\,G_{2,0}
-
d\,F_{2,0}
\Big]
\notag
\\
&
\ \ \ \ \
+
x\,y\,
\Big[
2\,a_{2,1}\,a_{2,2}\,G_{0,2}
+
a_{1,1}\,a_{2,2}\,G_{1,1}
+
a_{1,2}\,a_{2,1}\,G_{1,1}
+
2\,a_{1,1}\,a_{1,2}\,G_{2,0}
-
d\,F_{1,1}
\Big]
\notag
\\
&
\ \ \ \ \
+
y^2\,
\Big[
a_{2,2}^2\,G_{0,2}
+
a_{1,2}\,a_{2,2}\,G_{1,1}
+
a_{1,2}^2\,G_{2,0}
-
d\,F_{0,2}
\Big]
+
{\rm O}_{x,y}(3).
\end{align}
{\small (Another\,\,---\,\,less economic\,\,---\,\,way 
of doing would consist in applying
Lie's prolongation formulas~{\cite[Chap.~25]{Engel-Lie-Merker-2015}},
{\cite{Olver-1986, Bluman-Kumei-1989,
Olver-1995, Merker-2008, Chen-Merker-2019}}
of diffeomorphisms to the {\em full bundle} 
of second order jets,
before restricting these formulas to the considered fiber.)}

Since $G_\stab^1$ is a subgroup of $\GL(\C^3)$, its determinant
must be nonzero:
\[
0
\,\neq\,
\big(
a_{1,1}\,a_{2,2}
-
a_{2,1}\,a_{1,2}
\big)\,
d
\,=\,
\det\,
\left[\,
\begin{array}{ccc}
a_{1,1} & a_{1,2} & b_1
\\
a_{2,1} & a_{2,2} & b_2
\\
\red{\bf 0} & \red{\bf 0} & d
\end{array}
\,\right].
\]
Equating to zero the coefficients of $x^2$, of $x\,y$, of $y^2$,
and solving for $G_{2,0}$, $G_{1,1}$, $G_{0,2}$ gives
a {\em linear representation} on $\C^3$:
\[
\def\arraystretch{1.25}
\left[
\begin{array}{c}
G_{2,0}
\\
G_{1,1}
\\
G_{0,2}
\end{array}
\right]
\,=\,
\frac{1}{\left(a_{1,1}\,a_{2,2}-a_{2,1}\,a_{1,2}\right)^2}\,
\left[\,
\begin{array}{ccc}
a_{2,2}^2\,d & 
-\,a_{2,1}\,a_{2,2}\,d & 
a_{2,1}^2\,d
\\
-\,2\,a_{1,2}\,a_{2,2}\,d & 
a_{1,1}\,a_{2,2}\,d+a_{2,1}\,a_{1,2}\,d &
-\,2\,a_{1,1}\,a_{2,1}\,d
\\
a_{1,2}^2\,d &
-\,a_{1,1}\,a_{1,2}\,d &
a_{1,1}^2\,d
\end{array}
\,\right]
\left[
\begin{array}{c}
F_{2,0}
\\
F_{1,1}
\\
F_{0,2}
\end{array}
\right].
\]
This is the action of $G_\stab^1$ on $J^2 \big\vert_{0^1} = \C^3$,
and in fact, the action of 
the block-diagonal
subgroup $\GL(\C^2) \times \C^\ast \subset
\G_\stab^1$,
because $b_1$, $b_2$ are absent.

It is elementary to realize that this action is equivalent,
up to dilation, 
to the action of $\SL(\C^2)$ on binary quadrics, 
{\em cf.}~{\cite{Olver-1999}}, 
and to deduce that there are exactly 3 possible 
inequivalent normal forms at order 2:
\[
\aligned
\text{\green{\bf Branch 2a}}
\ \ \ \ \ \ \ \ \ \ \ \ \ \ \ \ \ \ \ \
u
&
\,=\,
0
\,\,\,\,
+
{\rm O}_{x,y}(3),
\\
\text{\green{\bf Branch 2b}}
\ \ \ \ \ \ \ \ \ \ \ \ \ \ \ \ \ \ \ \
u
&
\,=\,
x^2
\,\,
+
{\rm O}_{x,y}(3),
\\
\text{\green{\bf Branch 2c}}
\ \ \ \ \ \ \ \ \ \ \ \ \ \ \ \ \ \ \ \
u
&
\,=\,
x\,y
+
{\rm O}_{x,y}(3).
\endaligned
\]
Indeed, over the comblex numbers, both 
$x^2 + y^2$ and $x^2 - y^2$
are equivalent to $x\,y$. Geometrically, there are 3 group-orbits,
and there are 3\,\,---\,\,point-like, 
zero-dimensional\,\,---\,\,transversals

A quick way to recover this fact is to realize
by a direct computation that 
the Hessian at the origin
is a relative invariant:
\[
4\,G_{2,0}\,G_{0,2}
-
G_{1,1}^2
\,=\,
\frac{d^2}{a_{1,1}\,a_{2,2}-a_{2,1}\,a_{1,2}}\,
\Big[
4\,F_{2,0}\,F_{0,2}
-
F_{1,1}^2
\Big].
\]
Higher-dimensional Hessian matrices are also known to be
relatively invariant,
{\em see}~{\cite{Merker-2022}} for details from the power series
perspective.

\begin{Observation}
In all affine structures classified in this memoir,
at every jet order, there will 
always appear {\em explicit linear representations}
of subsequently reduced subgroups $G_\stab^{\kappa-1}$ 
on jet fibers $J^\kappa \big\vert_{T^{\kappa-1}}$ 
over certain group-transversals $T^{\kappa-1} \subset
J^{\kappa-1}$
from the jet level beneath.
\end{Observation}

None of the articles~{\cite{Eastwood-Ezhov-1999, 
Eastwood-Ezhov-2001-1, Eastwood-Ezhov-2001-2}} 
discovered the (universal!) occurrence of linear representations
on jet fibers at every jet order, 
like on $J^2 \big\vert_{0^1}$ here, 
identified with certain spaces of coefficients.
In~{\cite[App.~1]{Eastwood-Ezhov-2001-2}} about
hypersurfaces $H^3 \subset \C^4$ under $\Saff(\C^4)$,
a link with binary sextics is made, 
but without explicit linear representation,
and without explicit mention that the irreducible representations
of $\SL(\C^2)$ are (well)
known to be the spaces of binary quartics of any degree.

At the infinitesimal level, a general affine vector field:
\[
\aligned
L
&
\,=\,
\ \ 
\big(
T_1+A_{1,1}\,x+A_{1,2}\,y+B_1\,u
\big)\,\frac{\partial}{\partial x}
\\
&
\ \ \ \ \
+
\big(
T_2+A_{2,1}\,x+A_{2,2}\,y+B_2\,u
\big)\,\frac{\partial}{\partial y}
\\
&
\ \ \ \ \
+
\big(
U_0+C_1\,x+C_2\,y+D\,u
\big)\,\frac{\partial}{\partial u},
\endaligned
\]
is tangent to $\big\{u = F(x,y) \big\}$ 
if and only if:
\[
0
\,\equiv\,
L
\big(
-\,u+F(x,y)
\big)
\Big\vert_{u=F(x,y)},
\]
identically in $\C\{ x,y\}$.
With the normalization up to order 2
included:
\[
u
\,=\,
0
+
0
+
F_{2,0}\,x^2
+
F_{1,1}\,x\,y
+
F_{0,2}\,y^2
+
{\rm O}_{x,y}(3),
\]
these tangency equation reads:
\[
0
\,\equiv\,
-\,U_0
+
x\,
\big[
F_{1,1}\,T_2
+
2\,F_{2,0}\,T_1
-
C_1
\big]
+
y\,
\big[
2\,F_{0,2}\,T_2
+
F_{1,1}\,T_1
-
C_2
\big]
+
{\rm O}_{x,y}(2),
\]
whence necessarily:
\[
\aligned
U_0
&
\,:=\,
0,
\\
C_1
&
\,:=\,
F_{1,1}\,T_2
+
2\,F_{2,0}\,T_1,
\\
C_2
&
\,:=\,
2\,F_{0,2}\,T_2
+
F_{1,1}\,T_1.
\endaligned
\]

This corresponds to the group reduction to 
$G_\stab^1$ seen above, and this means that the general
infinitesimal generator of $G_\stab^1$ writes:
\[
\aligned
L_\stab^1
&
\,:=\,
\Big(
T_1
+
A_{1,1}\,x
+
A_{1,2}\,y
+
B_1\,u
\Big)\,\frac{\partial}{\partial x}
\\
&
\ \ \ \ \
+
\Big(
T_2
+
A_{2,1}\,x
+
A_{2,2}\,y
+
B_2\,u
\Big)\,\frac{\partial}{\partial y}
\\
&
\ \ \ \ \
+
\Big(
\big[
F_{1,1}\,T_2
+
2\,F_{2,0}\,T_1
\big]\,x
+
\big[
2\,F_{0,2}\,T_2
+
F_{1,1}\,T_1
\big]\,y
+
D\,u
\Big)\,\frac{\partial}{\partial u}.
\endaligned
\]

%%%%%%%%%%%%%%%%%%%%%%%%%%%%%%%%%%%%%%%%%%%%%%%%%%%%%%%%%%%%%%%%%%%%%%
\SectionHead{General Setting: Induction on Jet Order}
{general-setting-induction-jet-order}
%%%%%%%%%%%%%%%%%%%%%%%%%%%%%%%%%%%%%%%%%%%%%%%%%%%%%%%%%%%%%%%%%%%%%%

Now, come back to the general setting.
As before, we will avoid introducing 
sophisticated notation which would incorporate 
all the necessary indices.
Sometimes, we will employ the simple
notation ``${}_\smallbullet$'' to denote an index,
or even, several indices.

In the next paragraphs, we will not 
{\em prove}, in the rigorous sense, 
all of our assertions.
Rather, a bit in the spirit of
Cartan himself~{\cite{Cartan-1937}},
and as in the monograph~{\cite{Olver-1995}} 
of Peter Olver
on Cartan's method of equivalence,
we will {\em describe}
a method of equivalence,
emphasizing the diversity and the complexity
of what really occurs in various examples.
Especially, the present memoir will show a wealth
of examples.

Reasoning by induction up to p.~{\pageref{terminate-induction-kappa}}, 
let $\kappa \geqslant 1$ be the working jet order.
\label{start-induction-kappa}
We suppose 
that 
at the preceding orders $\leqslant \kappa-1$, there were 
a certain number
of transversals to the orbits of certain linear
representations of certain (reduced) subgroups
contained in $G$. 
What holds by induction at the preceding jet orders
will be clearer later.
Some of these
transversals may have consisted of a single point.
Some other transversals may have been positive-dimensional,
and in this case, they may have incorporated absolute invariants 
$\Iaux_\smallbullet$
which became parameters in subsequent computations.

The determination of transversals relies upon choices,
it is not canonical. At this preceding jet order $\kappa-1$,
as many branches were created as there were (reduced sugroup)
orbits, and also, as many as there were orbit transversals.
A similar situation will occur at the next jet order $\kappa$.

For instance, in the case of surfaces $S^2 \subset \C^3$, 
the creation of 3 branches at order 2 can be represented 
and labelled as:
\[
\def\arraystretch{1.25}
\begin{array}{rccccc}
\green{\bf 1}\,\,\,\,
\green{\downarrow}\,\,
& 
F_{2,0} & F_{1,1} & F_{0,2}
\\
\green{\bf 2a} & 
0 & 0 & 0
\\
\green{\bf 2b} & 
1 & 0 & 0
\\
\green{\bf 2c} & 
0 & 1 & 0
\end{array}
\]
These 3 transversals are zero-dimensional.
Several similar {\sl branch-creation arrays} 
will be set up 
later in this memoir, {\em see} 
{\em e.g.}~Sections~{\ref{classification-Hn-Rn-1}},
{\ref{branching-diagrams}}.

Pick one such transversal $T^{\kappa-1}$. 
Forget other transversals or
branches (which can be treated similarly).
For instance, pick \green{\bf 2b} above.

The stabilization
of the selected transversal $T^{\kappa-1}$
then requires to {\em reduce}
the group. Why? Because the normal form constructed
up to this point must be kept, be preserved, be conserved! 
And the normal form is the same, {\em on both sides!}

Denote by $G_\stab^{\kappa-1} \subset G$
the concerned
stabilization subgroup. This generalizes the 
subgroups $G_\stab^2 \subset G_\stab^1 \subset G$ introduced above.

In the $(x,u)$-space
and in the $(\overline{x}, \overline{u}$)-space as well, 
let the two normal forms be written as:
\[
u_j
\,=\,
{\sf N}_{j,\kappa-1}^\normal
\big(\Iaux_{\smallbullet},x\big)
+
\!\!\!\!\!\!\!
\sum_{i_1+\cdots+i_n\geqslant\kappa}
\!\!\!\!\!
F_{j,i_1,\dots,i_n}\,
x_1^{i_1}\cdots x_n^{i_n}
\ \ \ \ \ \ \
\text{and}
\ \ \ \ \ \ \
\overline{u}_j
\,=\,
{\sf N}_{j,\kappa-1}^\normal
\big(\Iaux_{\smallbullet},\overline{x}\big)
+
\!\!\!\!\!\!\!
\sum_{i_1+\cdots+i_n\geqslant\kappa}
\!\!\!\!\!
\overline{F}_{j,i_1,\dots,i_n}\,
\overline{x}_1^{i_1}\cdots\overline{x}_n^{i_n},
\]
with $1 \leqslant j \leqslant c$, 
where the following holds.

\medskip\noindent$\square$\,
The ${\sf N}_{j,\kappa-1}^\normal$ represent all $x$-monomials 
in the left space and all $\overline{x}$-monomials in the right
space, monomials which are
{\em normalized and finalized} up to order $\leqslant \kappa-1$.

\medskip\noindent$\square$\,
These normalized polynomials ${\sf N}_{j,\kappa-1}^\normal$ are
{\em exactly the same functions} on both 
sides\,\,---\,\,only 
the argument $x$ is changed to $\overline{x}$.

\medskip\noindent$\square$\,
The supplementary argument $\Iaux_\smallbullet$ 
(without indices, sometimes absent)
indicates that in some branches, there might remain a
certain number of {\em absolute invariants} found
in preceding orders,
namely function satisfying in this branch:
\[
\Iaux_\smallbullet
\Big(J^{\kappa-1}F\Big)
\,=\,
\Iaux_\smallbullet
\Big(\overline{J}^{\kappa-1}\overline{F}\Big),
\]
with on both sides {\em exactly the same functions}
$\Iaux_\smallbullet$ of the collection of
order $\leqslant \kappa-1$ 
power series coefficients\,\,---\,\,plainly 
denoted here with the notation $J^{\kappa-1}$.

\medskip

So now, how to determine $G_\stab^{\kappa-1}$? Just by requiring
that the normal form is preserved by
a transformation $g \in G$ up to order $\leqslant \kappa-1$.
In the example of $S^2 \subset \C^3$ under $\Aff(\C^3)$,
we saw the {\em fundamental 
equation}~{\eqref{intro-fund-eq-S2-C3}},
and we truncated it at order 1 to get
$G_\stab^1$ with $c_1 = c_2 = 0$.

In the general setting, the reduced group $G_\stab^{\kappa-1}
\subset G_\stab^{\kappa-2} \subset \cdots \subset G$ 
can be determined,
theoretically, as follows. At first, with $g \in G_\stab^{\kappa-2}$, 
let the group-dependent
diffeomorphism $(x,u) \longmapsto \big( \overline{x}, 
\overline{u} \big)$ be written as:
\[
\overline{x}
\,=\,
\overline{x}
\big(x,u,g\big),
\ \ \ \ \ \ \ \ \ \ \ \ \ \ \ \ \ \ \ \
\overline{u}
\,=\,
\overline{u}
\big(x,u,g\big).
\]
For $G = \Aff(\C^3)$, such formulas are explicit.
Such a diffeomorphism
maps $\big\{ u = F(x) \big\}$ to $\big\{\overline{u} =
\overline{F} (\overline{x}) \big\}$ if and only if:
\[
u
\,=\,
F(x)
\ \ \ \ \ \ \ \ \ \ \ \ \ \ \ \ \ \ \ \
\Longrightarrow
\ \ \ \ \ \ \ \ \ \ \ \ \ \ \ \ \ \ \ \
\overline{u}
\,=\,
\overline{F}
\big(\overline{x}\big),
\]
which yields the {\sl fundamental equations}, in the current branch:
\[
\aligned
0
&
\,\equiv\,
-\,
\overline{u}_j
\big(x,F(x),g\big)
+
\overline{F}_j
\Big(
\overline{x}
\big(x,F(x),g\big)
\Big)
\\
&
\,\equiv\,
\sum_{i_1+\cdots+i_n\geqslant 0}\,
{\sf E}_{j,i_1,\dots,i_n}^{\sf nf}
\Big(
\Iaux_\smallbullet,
F_\smallbullet,
\overline{F}_\smallbullet,
g
\Big)\,
x_1^{i_1}
\cdots
x_n^{i_n}.
\endaligned
\]

These $c$ equations for
$1 \leqslant j \leqslant c$
should be satisfied {\em identically} in
$\C\{x_1, \dots, x_n\}$. 
The upper index ${}^{\sf nf}$ in
${\sf E}_\smallbullet^{\sf nf}$ indicates that these equations
are involved in the production of {\em {\sf n}ormal {\sf f}orms}.
{\em Infra}, we will introduce 
other kinds of equations ${\sf E}_\smallbullet^{\sf vf}$ 
with the upper index ${}^{\sf vf}$,
indicating that they come from
tangential {\em {\sf v}ector {\sf f}ields}.

So all these $E_{j,i_1,\dots,i_n}^{\sf nf} = 0$ should
vanish.
Above, the lightened notation
$F_\smallbullet$ denotes a certain
finite collections of power series coefficient
$F_{j,i_1',\dots,i_n'}$, always with $i_1' + \cdots + i_n'
\leqslant i_1 + \cdots + i_n$,
and the same for $\overline{F}_\smallbullet$.
In practice, real formulas are challenging,
even for powerful symbolic computers.

By the induction hypothesis, since $g \in G_\stab^{\kappa-2}$,
all equations  $E_{j,i_1,\dots,i_n}^{\sf nf} = 0$ with
$1 \leqslant j \leqslant c$ and with
$i_1 + \cdots + i_n \leqslant \kappa-2$ are already fulfilled,
and it remains:
\[
0
\,\equiv\,
\sum_{i_1+\cdots+i_n=\kappa-1}\,
{\sf E}_{j,i_1,\dots,i_n}^{\sf nf}
\Big(
\Iaux_\smallbullet,
F_\smallbullet,
\overline{F}_\smallbullet,
g
\Big)\,
x_1^{i_1}
\cdots
x_n^{i_n}
+
{\rm O}_{x_1,\dots,x_n}(\kappa)
\eqno
{\scriptstyle{(1\,\leqslant\,j\,\leqslant\,c)}},
\]
whence:
\[
0
\,=\,
{\sf E}_{j,i_1,\dots,i_n}^{\sf nf}
\Big(
\Iaux_\smallbullet,
F_\smallbullet,
\overline{F}_\smallbullet,
g
\Big)
\eqno
{\scriptstyle{(\forall\,1\,\leqslant\,j\,\leqslant\,c,\,\,
\forall\,i_1+\cdots+i_n\,=\,\kappa-1)}}.
\]

Once $F_\smallbullet$
is chosen in a certain transversal $T^{\kappa-1}$
with (by invariancy)
the same choice for $\overline{F}_\smallbullet$,
these (algebraic) equations are used 
as supplementary constraints on $g \in G_\stab^{\kappa-2}$.
These equations therefore force $g$ to belong 
to a specific {\em reduced} subgroup 
$G_\stab^{\kappa-1} \subset G_\stab^{\kappa-2}$.

In this order $\kappa-1$ preceding the working order $\kappa$, 
because we reason by induction,
we have not yet explained how transversals 
$T^{\kappa-1}$ 
to $G_\stab^{\kappa-1}$-orbits 
were constructed\big/chosen. 
This aspect is more delicate.
{\em Infra}, 
at the next (working) order $\kappa$,
we will explain how to create
transversals $T^\kappa$.
At least for now, in our reasoning by induction,
we have explained
what we assume to be achieved 
at orders $\leqslant \kappa-1$

Once $G_\stab^{\kappa-1}$ is known, 
the next step is to {\em prolong} its action
to the space of $\kappa$-jets.
Remember that we do {\em not} work in full jet bundles,
which is a key trick to dominate the complexity of computations.
We work only above successive transversals.
This means that we work over the already normalized
power series coefficients, at orders $\leqslant \kappa-1$,
namely `over' ${\sf N}_{j,i_1,\dots,i_n}^\normal$,
symetrically on both left and right sides.

Also, this means that the relative fiber of the projection
from $\kappa$-jets to normalized jets of order $\leqslant \kappa-1$
is represented just by letting appear order $= \kappa$
power series coefficients:
\[
u_j
\,=\,
{\sf N}_{j,\kappa-1}^\normal
\big(\Iaux_{\smallbullet},x\big)
+
\!\!\!\!\!\!\!
\sum_{i_1+\cdots+i_n=\kappa}
\!\!\!\!\!
F_{j,i_1,\dots,i_n}\,
x_1^{i_1}\cdots x_n^{i_n}
+
{\rm O}_{x_1,\dots,x_n}(\kappa+1),
\]
and the same for:
\[
\overline{u}_j
\,=\,
{\sf N}_{j,\kappa-1}^\normal
\big(\Iaux_{\smallbullet},\overline{x}\big)
+
\!\!\!\!\!\!\!
\sum_{i_1+\cdots+i_n=\kappa}
\!\!\!\!\!
\overline{F}_{j,i_1,\dots,i_n}\,
\overline{x}_1^{i_1}\cdots\overline{x}_n^{i_n}
+
{\rm O}_{\overline{x}_1,\dots,\overline{x}_n}(\kappa+1).
\]
Of course, the appearing 
$F_{j,i_1,\dots,i_n}$ and $\overline{F}_{j,
i_1,\dots,i_n}$ are {\em a priori} different here
(while at orders $\leqslant \kappa-1$, they are equal
by construction).

The goal is to {\em normalize} these
$F_{j,i_1,\dots,i_n}$ and $\overline{F}_{j,
i_1,\dots,i_n}$, {\em i.e.} to
find appropriate orbit transversals. But for which group action?
It is at this precise step that things often happen to 
become delicate. 

Abbreviating:
\[
\aligned
J_\ast^\kappa F
&
\,:=\,
\big\{
F_{j,i_1,\dots,i_n}
\big\}_{i_1+\cdots+i_n=\kappa}^{1\leqslant j\leqslant c},
\\
J_\ast^\kappa\overline{F}
&
\,:=\,
\big\{
\overline{F}_{j,i_1,\dots,i_n}
\big\}_{i_1+\cdots+i_n=\kappa}^{1\leqslant j\leqslant c},
\endaligned
\]
the fundamental equation,
which is now identically satisfied up to all orders
$\leqslant \kappa-1$ when $g \in G_\stab^{\kappa-1}$, 
reads at order $\kappa$ as:
\[
0
\,=\,
{\sf E}_{j,i_1,\dots,i_n}^{\sf nf}
\Big(
\Iaux_\smallbullet,
J_\ast^\kappa F,
J_\ast^\kappa\overline{F},
g
\Big)
\eqno
{\scriptstyle{(\forall\,1\,\leqslant\,j\,\leqslant\,c,\,\,
\forall\,i_1+\cdots+i_n\,=\,\kappa)}}.
\]
Provided that $g \in G_\stab^{\kappa-1}$ lies in some
neighborhood of the identity, these algebraic equations,
of degree 1 with respect to $J_\ast^\kappa F$ and
to $J_\ast^\kappa \overline{F}$, may always be solved 
under the form:
\[
J_\ast^\kappa\overline{F}
\,=\,
\Lambda\,
\Big(
\Iaux_\smallbullet,J_\ast^\kappa F,g
\Big).
\]

But some key information may be missing.
This is a `defect' of the {\sf n}ormal {\sf f}orm
equations $0 = {\sf E}_\smallbullet^{\sf nf}$
which, by working only over the origin
$(x,u) = (0,0)$, are unable
{\em per se} to capture differentialo-geometric information.

%%%%%%%%%%%%%%%%%%%%%%%%%%%%%%%%%%%%%%%%%%%%%%%%%%%%%%%%%%%%%%%%%%%%%%
\SectionHead{Non-Reduced Linear Representation}
{non-reduced-linear-representation}
%%%%%%%%%%%%%%%%%%%%%%%%%%%%%%%%%%%%%%%%%%%%%%%%%%%%%%%%%%%%%%%%%%%%%%

At his point, before explaining further the general theory,
it is better to come back to our concrete example.
For surfaces $S^2 \subset \C^3$, in the branch $\green{\bf 2b}$:
\[
u
\,=\,
x^2
+
{\rm O}_{x,y}(3)
\ \ \ \ \ \ \ \ \ \ \ \ \ \ \ \ \ \ \ \
\text{and}
\ \ \ \ \ \ \ \ \ \ \ \ \ \ \ \ \ \ \ \
v
\,=\,
p^2
+
{\rm O}_{p,q}(3),
\]
the fundamental equation~{\eqref{intro-fund-eq-S2-C3}},
truncated at order $2 = \kappa-1$:
\[
0
\,=\,
x^2\,
\big[
a_{1,1}^2
-
d
\big]
+
xy\,
\big[
2\,a_{1,1}\,a_{1,2}
\big]
+
y^2\,
\big[
a_{1,2}^2
\big]
+
{\rm O}_{x,y}(3),
\]
resolved as:
\[
d
\,:=\,
a_{1,1}^2,
\ \ \ \ \ \ \ \ \ \ \ \ \ \ \ \ \ \ \ \
a_{1,2}
\,:=\,
0,
\]
imposes the group reduction towards $G_\stab^2$:
\[
\left[\,
\begin{array}{ccc}
a_{1,1} & a_{1,2} & b_1
\\
a_{2,1} & a_{2,1} & b_1
\\
\red{\bf 0} & \red{\bf 0} & d
\end{array}
\,\right]^{\green{\bf 1}}
\ \ \ \ \
\leadsto
\ \ \ \ \
\left[\,
\begin{array}{ccc}
a_{1,1} & \red{\bf 0} & b_1
\\
a_{2,1} & a_{2,1} & b_1
\\
\red{\bf 0} & \red{\bf 0} & \red{a_{1,1}^2}
\end{array}
\,\right]^{\green{\bf 2}},
\]
with nonzero determinant:
\[
0
\,\neq\,
a_{1,1}^3\,a_{2,2}.
\]
Hence, so is our $G_\stab^{\kappa-1}$.

From the infinitesimal side, 
after normalizing as above in this branch
\green{\bf 2b}:
\[
\aligned
U_0
&
\,:=\,
0,
\\
C_1
&
\,:=\,
2\,T_1,
\\
C_2
&
\,:=\,
0,
\endaligned
\]
whence:
\[
\aligned
L_\stab^1
&
\,:=\,
\Big(
T_1
+
A_{1,1}\,x
+
A_{1,2}\,y
+
B_1\,u
\Big)\,\frac{\partial}{\partial x}
\\
&
\ \ \ \ \
+
\Big(
T_2
+
A_{2,1}\,x
+
A_{2,2}\,y
+
B_2\,u
\Big)\,\frac{\partial}{\partial y}
\\
&
\ \ \ \ \
+
\Big(
2\,T_1\,x
+
D\,u
\Big)\,\frac{\partial}{\partial u},
\endaligned
\]
the tangency of $L_\stab^1$ to:
\[
u
\,=\,
x^2
+
F_{3,0}\,x^3
+
F_{2,1}\,x^2y
+
F_{1,2}\,xy^2
+
F_{0,3}\,y^3
+
{\rm O}_{x,y}(4),
\]
gives:
\leqnomode\usetagform{default}
\begin{align}
\label{F12-F03-zero}
0
&
\,\equiv\,
x^2\,
\Big[
3\,F_{3,0}\,T_1
+
F_{2,1}\,T_2
+
2\,A_{1,1}
-
\underline{D}
\Big]
\notag
\\
&
\ \ \ \ \
+
xy\,
\Big[
2\,F_{2,1}\,T_1
+
2\,F_{1,2}\,T_2
+
2\,\underline{A_{1,2}}
\Big]
\\
&
\ \ \ \ \
+
y^2\,
\Big[
F_{1,2}\,T_1
+
3\,F_{0,3}\,T_2
\Big]
+
{\rm O}_{x,y}(3),
\notag
\end{align}
so that the infinitesimal counterpart of group reduction 
is\,\,---\,\,observe the parallel between 
lower-case and upper-case letters\,\,---:
\[
\aligned
D
&
\,:=\,
3\,F_{3,0}\,T_1
+
F_{2,1}\,T_2
+
2\,A_{1,1},
\\
A_{1,2}
&
\,:=\,
-\,F_{2,1}\,T_1
-
F_{1,2}\,T_2.
\endaligned
\]
Replacing in $L_\stab^1$ gives $L_\stab^2$.
Temporarily, let us
disregard (or forget) the equation obtained 
from the coefficient of $y^2$ above.

Next, a map $g \in G_\stab^2$ between the two surfaces:
\[
\aligned
u
&
\,=\,
x^2
+
F_{3,0}\,x^3
+
F_{2,1}\,x^2y
+
F_{1,2}\,xy^2
+
F_{0,3}\,y^3
+
{\rm O}_{x,y}(4),
\\
v
&
\,=\,
p^2
+
G_{3,0}\,p^3
+
G_{2,1}\,p^2q
+
G_{1,2}\,pq^2
+
G_{0,3}\,q^3
+
{\rm O}_{p,q}(4),
\endaligned
\]
has fundamental equation at order 3:
\[
\aligned
0
&
\,\equiv\,
x^3\,
\Big[
-\,a_{1,1}^2\,F_{3,0}
+
a_{1,1}^3\,G_{3,0}
+
a_{2,1}^3\,G_{0,3}
+
a_{1,1}\,a_{2,1}^2\,G_{1,2}
+
a_{1,1}^2\,a_{2,1}\,G_{2,1}
+
2\,a_{1,1}\,\boxed{b_1}
\Big]
\\
&
\ \ \ \ \
+
x^2y\,
\Big[
-\,a_{1,1}^2\,F_{2,1}
+
a_{1,1}^2\,a_{2,2}\,G_{2,1}
+
3\,a_{2,1}^2\,a_{2,2}\,G_{0,3}
+
2\,a_{1,1}\,a_{2,1}\,a_{2,2}\,G_{1,2}
\Big]
\\
&
\ \ \ \ \
+
xy^2\,
\Big[
-\,a_{1,1}^2\,F_{1,2}
+
a_{1,1}\,a_{2,2}^2\,G_{1,2}
+
3\,a_{2,1}\,a_{2,2}^2\,G_{0,3}
\Big]
\\
&
\ \ \ \ \
+
y^3\,
\Big[
-\,a_{1,1}^2\,F_{0,3}
+
a_{2,2}^3\,G_{0,3}
\Big]
+
{\rm O}_{x,y}(4).
\endaligned
\]
The (boxed) free group parameter $b_1$ can be used to normalize
in the right space:
\[
G_{3,0}
\,:=\,
0,
\]
simply by setting:
\[
b_1
\,:=\,
-\,\tfrac{1}{2}\,
\frac{a_{2,1}^3}{a_{1,1}}\,
G_{0,3}
-
\tfrac{1}{2}\,
a_{2,1}^2\,G_{1,2}
-
\tfrac{1}{2}\,
a_{1,1}\,a_{2,1}\,
+
\tfrac{1}{2}\,a_{1,1}\,F_{3,0},
\]
the division by $a_{1,1}$ being
allowed since the determinant $a_{1,1}^3\, a_{2,2} \neq 0$.

Then, by restarting from the right
space, and by renaming $G_{\smallbullet, \smallbullet}$ 
as $F_{\smallbullet, \smallbullet}$ put on the left, 
we can assume that $F_{3,0} = 0$, 
and we can again normalize $G_{3,0} := 0$ using again $b_1$.
Lastly, to {\em stabilize} the symmetric normalization:
\[
G_{3,0}
\,=:\,
0
\,:=\,
F_{3,0},
\]
it is necessary and sufficient to set:
\[
b_1
\,:=\,
-\,\tfrac{1}{2}\,
\frac{a_{2,1}^3}{a_{1,1}}\,
G_{0,3}
-
\tfrac{1}{2}\,
a_{2,1}^2\,G_{1,2}
-
\tfrac{1}{2}\,
a_{1,1}\,a_{2,1}.
\]

Once this is done, 3 equations remain, the coefficients above of
$x^2y$, of $xy^2$, of $y^3$.
Solving for $\big( G_{2,1}, G_{1,2}, G_{0,3} \big)$ gives
a linear representation on $\C^3$:
\[
\def\arraystretch{1.25}
\left[
\begin{array}{c}
G_{2,1}
\\
G_{1,2}
\\
G_{0,3}
\end{array}
\right]
\,=\,
\frac{1}{a_{2,2}^3}\,
\left[\,
\begin{array}{ccc}
a_{2,2}^2 & -2a_{2,1}a_{2,2} & 3a_{2,1}^2
\\
0 & a_{1,1}a_{2,2} & -3a_{1,1}a_{2,1}
\\
0 & 0 & a_{1,1}^2
\end{array}
\,\right]
\left[
\begin{array}{c}
F_{2,1}
\\
F_{1,2}
\\
F_{0,3}
\end{array}
\right],
\]
with determinant:
\[
\frac{a_{1,1}^3}{a_{2,2}^6}
\,\neq\,
0.
\]

{\em However, this linear representation 
of $G_\stab^2$ 
is not the `right one'!}
But why?

Because in Branch~\green{\bf 2b}, there is 
a {\em degenerate}
geometric situation, in which the Hessian is of {\em non}maximal 
rank $1 < 2$. This assumption was made at the origin.
And the rank of the Hessian matrix:
\[
\Hessian_F(x,y)
\,:=\,
\left[
\begin{array}{cc}
F_{xx} & F_{xy}
\\
F_{yx} & F_{yy}
\end{array}
\right]
(x,y),
\]
is invariant under affine transformations at points
corresponding to each other through transformations.
This follows from the above considerations,
{\em see} also~{\cite[Sec.~2]{Merker-2022}}
where all details are written
for hypersurfaces $H^n \subset \R^{n+1}$ (same 
result over $\C$).

Of course, the rank of the Hessian matrix may well `jump'
from point to point. More generally, 
for any geometric structure, 
pointwise invariant hypotheses
may well `jump' from point to point,
as matrix ranks may do.
Singularity Theory studies such (`jumping') objects.

On the other hand, Differential Geometry in the spirit of Lie
and Cartan traditionally
decides to adopt {\em constant rank 
and constant geometric 
hypotheses} (only). In this memoir,
we definitely adopt this point of view.

%%%%%%%%%%%%%%%%%%%%%%%%%%%%%%%%%%%%%%%%%%%%%%%%%%%%%%%%%%%%%%%%%%%%%%
\SectionHead{Setting Up Dependent Jets}
{setting-up-dependent-jets}
%%%%%%%%%%%%%%%%%%%%%%%%%%%%%%%%%%%%%%%%%%%%%%%%%%%%%%%%%%%%%%%%%%%%%%

We are therefore led to assume that the
Hessian determinant is constantly of rank 1, at all point
$(x,y)$ in some open neighborhood of the origin $(0,0)$.
In the memoir~{\cite{Chen-Merker-2019}}
based on a function-theoretic partial
differential equations approach, this assumption was expressed
as:
\[
F_{xx}
\,\neq\,
0
\,\equiv\,
F_{xx}\,F_{yy}
-
F_{xy}^2,
\]
at all point $(x,y)$ near $(0,0)$. 
Then the so-called {\sl parabolic surfaces} were defined
by requiring that:
\[
F_{yy}(x,y)
\,\equiv\,
\frac{F_{xy}^2(x,y)}{F_{xx}(x,y)}.
\]

Plugging $F = x^2 + F_{3,0}\, x^3 + F_{2,1}\, x^2y + 
F_{1,2}\, xy^2 + F_{0,3}\, y^3 + \cdots$ into this {\sc pde}
then forces:
\leqnomode\usetagform{default}
\begin{align}
\label{Hessian-F12-F03-zero}
F_{1,2}
\,=\,
0,
\ \ \ \ \ \ \ \ \ \ \ \ \ \ \ \ \ \ \ \
F_{0,3}
\,=\,
0.
\end{align}
By invariancy, $G_{1,2} = G_{0,3} = 0$ also.
This shows that the above linear representation on $\C^3$
{\em reduces} to a (simpler) linear representation 
on $\C^1$:
\[
\Big[
G_{2,1}
\Big]
\,=\,
\Big[\tfrac{1}{a_{2,2}}\Big]\,
\Big[
F_{2,1}
\Big],
\]
which is the `right one'.
Visibly, $F_{2,1}$ is a (punctual) relative invariant.

In brief, the delicate step is to {\em reduce} the jet
fiber by taking account of the geometric assumptions
that run from previous jet orders in the current branch
under study.

%%%%%%%%%%%%%%%%%%%%%%%%%%%%%%%%%%%%%%%%%%%%%%%%%%%%%%%%%%%%%%%%%%%%%%
\SectionHead{Reduced Linear Representation and Branch Creation}
{reduced-linear-representation-branch-creation}
%%%%%%%%%%%%%%%%%%%%%%%%%%%%%%%%%%%%%%%%%%%%%%%%%%%%%%%%%%%%%%%%%%%%%%

Back to the general setting, with $g \in G_\stab^{\kappa-1}$, 
in the order $\kappa$ normal form equations:
\[
0
\,=\,
{\sf E}_{j,i_1,\dots,i_n}^{\sf nf}
\Big(
\Iaux_\smallbullet,
J_\ast^\kappa F,
J_\ast^\kappa\overline{F},
g
\Big)
\eqno
{\scriptstyle{(\forall\,1\,\leqslant\,j\,\leqslant\,c,\,\,
\forall\,i_1+\cdots+i_n\,=\,\kappa)}}.
\]
some jet coordinates in $J_\ast^\kappa F$ and, parallelly,
in $J_\ast^\kappa \overline{F}$, should disappear
due to the previous history within the branches created before.

\begin{center}
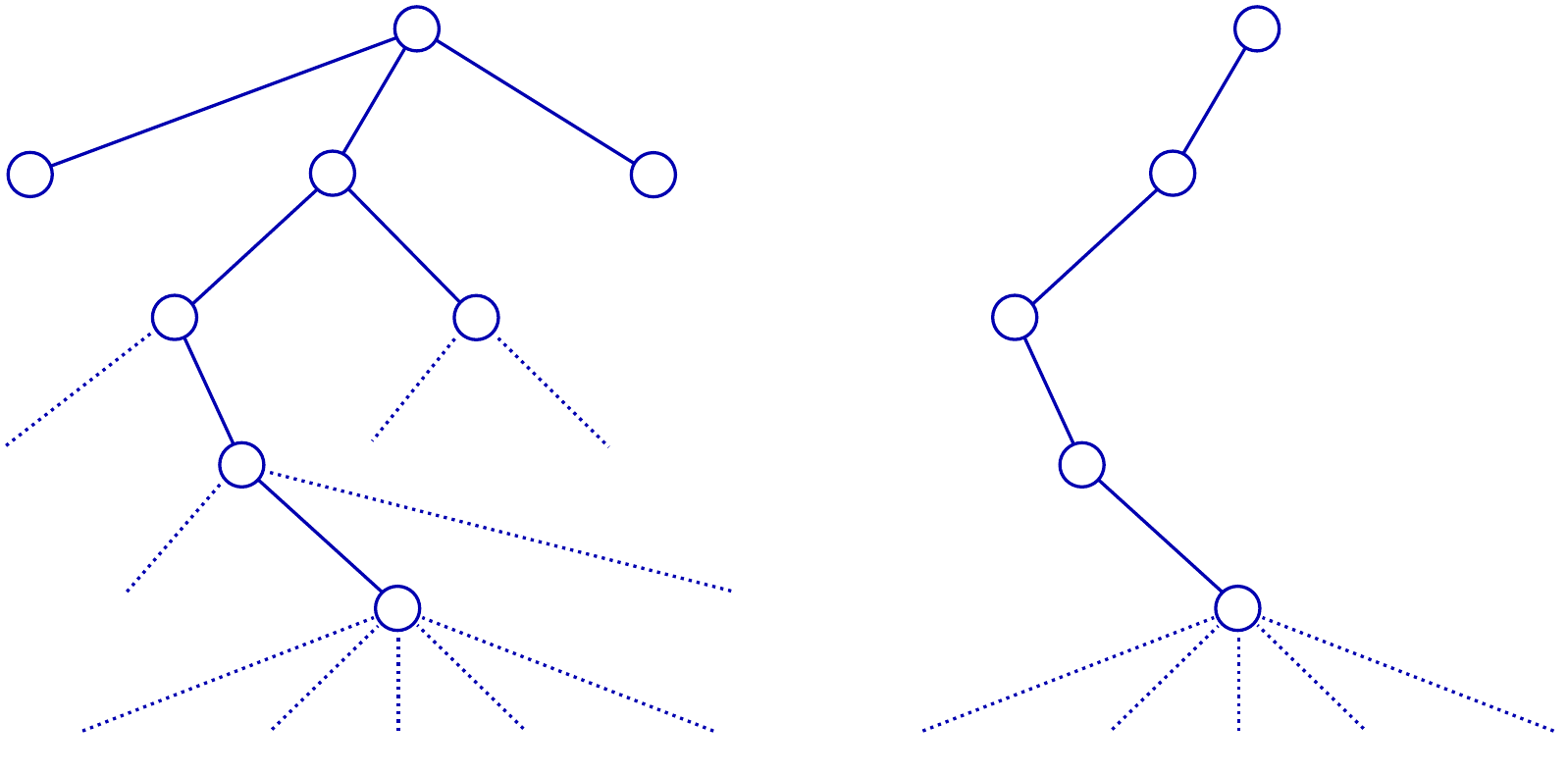
\end{center}

In the illustrating figure above, \green{\bf 5b} would be the
branch at order $\kappa-1 = 5$ at which considerations hold
(instead of $\kappa-1 = 3$ in the preceding
Sections~{\ref{non-reduced-linear-representation}} 
and~{\ref{setting-up-dependent-jets}}),
with nearby branches, and with the whole history
of preceding branches. Still,
the creation of order $6 = \kappa$ 
subsequent branches is not yet 
done.

\begin{center}
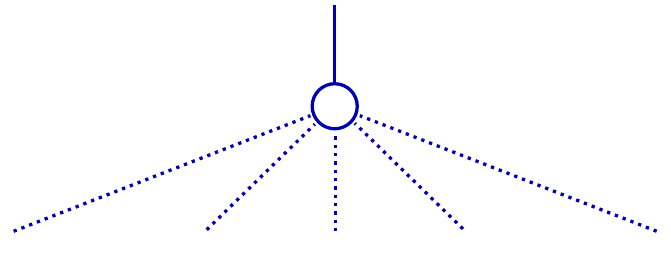
\end{center}

In the previous history, some
relative differential invariants, say
$\Kaux_1, \dots, \Kaux_t$, were encountered
which were assumed to be $\equiv 0$.
{\small 
(Some other relative differential invariants may have been
assumed to be nonzero and then normalized to $+1$ or to $-1$
with associated group reductions,
but such kinds of normalizations have no differential consequences.)}
These {\em invariant} differential relations:
\[
0
\,\equiv\,
\Kaux_1\big(J^{\kappa-1}F\big),\,\,\,
\dots\dots\dots,\,\,\,
0
\,\equiv\,
\Kaux_t\big(J^{\kappa-1}F\big),
\]
encountered at jet orders
$\leqslant \kappa-1$, do not depend on $J_\kappa^\ast F$.

But by differentiation with respect to $x_1, \dots, x_n$, 
these {\sc pde}s do (in general) provide resolutions of certain
{\sl dependent}
$J_{\ast,\dep}^\kappa F$ in terms of some other
{\em independent} $J_{\ast,\ind}^\kappa F$,
possibly with discussion of determinantal loci,
hence with creation of branches. 
For instance, from the parabolic surfaces differential
relation
$F_{yy} = \frac{F_{xy}^2}{F_{xx}}$ with $\kappa-1 = 2$,
it comes~{\cite{Chen-Merker-2019}}:
\[
\aligned
F_{xyy}
&
\,=\,
2\,
\frac{F_{xy}\,F_{xxy}}{F_{xx}}
-
\frac{F_{xy}^2\,F_{xxx}}{F_{xx}^2},
\\
F_{yyy}
&
\,=\,
3\,
\frac{F_{xy}^2\,F_{xxy}}{F_{xx}^2}
-
2\,
\frac{F_{xy}^3\,F_{xxx}}{F_{xx}^3}.
\endaligned
\]

Furthermore,
there can exist new invariant differential relations 
of orders $\leqslant \kappa-1$
obtained by cross differentiations,
{\em cf.}~{\cite{Chen-Merker-2020}}.
By invariancy, exactly the same (partial)
resolution formulas hold about 
$J_\ast^\kappa \overline{F}$.
The terminology {\sl dependent}\big/{\sl independent} jets
is employed in~{\cite{Chen-Merker-2019}}.

{\small There is still another aspect of complexity.
Relative differential invariants 
like $\Kaux_1, \dots, \Kaux_t$ above,
come each from a 1-dimensional
linear representation of the (trivial) group $\R^\ast$.
But sometimes, a linear representation which was
encountered at some jet order $\lambda \leqslant \kappa-1$
and which was responsible for the creation of branches
was {\em not} 1-dimensional, not diagonalizable,
not even triangularizable.
In this (delicate) situation, the differential consequences
of some degenerate branches might be difficult to set up.}

To express the differential consequences of
$0 \equiv \Kaux_1 \equiv \cdots \equiv \Kaux_t$, 
currently, two approaches exist.

The first approach, used {\em e.g.} 
in~{\cite{Chen-Merker-2019, Chen-Merker-2020}},
is to compute the explicit expressions
of these $\Kaux_\smallbullet \big( J_\ast^\kappa F \big)$.
It is practicable when dealing 
with small-sized invariants like the Hessian,
{\em cf.} also~{\cite{Merker-2020}}.
In the Cauchy-Riemann (CR) context 
and in the para-CR context,
the analog of the Hessian matrix is the Levi matrix, 
again a relative differential invariant (matrix)
of reasonable explicit size. 
Nonetheless, 
precise and fine
calculations of differential consequences 
happen to be surprisingly challenging, as may be seen 
in~{\cite{Merker-Pocchiola-2018, Foo-Merker-2019,
Foo-Merker-Ta-2019,
Merker-Nurowski-2020-a}},
%%%Merker-Nurowski-2020-b}},
{\em cf.} also the survey~{\cite{Merker-Nurowski-2023}}.

The second approach, 
{\em see} 
{\em e.g.}~{\cite{Arnaldsson-Valiquette-2020,
Olver-Sabzevari-Valiquette-2023}}, 
was recently developed
within the famous (partial)
moving (co)frames
theory developed by Peter Olver and his school.
It consists in apply partial, recursive,
invariant differentiations to the symbolic
(not explicit) invariant equations 
$0 \equiv \Kaux_\smallbullet$.
It is elegant, economic, efficient.

Unfortunately, by manipulating only 
power series centered at the origin,
it is impossible to reach differential
consequences of degeneracy hypotheses. 
Here lies the main defect of the {\em strict},
{\em i.e.} not simultaneously 
using the theory of differential invariants,
power series approach.

In summary, coming back to our power series,
let us admit that all the 
order $\kappa$ differential
consequences of the degeneracy assumptions 
encountered before in the current branch
are computable in some
`external' way and have been inserted in the
order $\kappa$ normal form equations:
\[
0
\,=\,
{\sf E}_{j,i_1,\dots,i_n}^{\sf nf}
\Big(
\Iaux_\smallbullet,
J_{\ast,\ind}^\kappa F,
J_{\ast,\ind}^\kappa\overline{F},
g
\Big)
\eqno
{\scriptstyle{(\forall\,1\,\leqslant\,j\,\leqslant\,c,\,\,
\forall\,i_1+\cdots+i_n\,=\,\kappa)}},
\]
with $g \in G_\stab^{\kappa-1}$.

In fact, since we will abandon 
the Differential Invariants Problem~{\ref{Pbm-debut-diff-invts}}
(but {\em see} Section~{\ref{creations-geometries}} {\em infra}),
and focus only on the Homogenous Models
Problem~{\ref{Pbm-debut-hom-models}},
we will develop a precise, elementary, and unambiguous
method for
determining
the explicit expressions of
the {\sl dependent} jets $J_{\ast,\dep}^\kappa F$,
together with some extra jet 
constraints required to construct homogeneous geometries,
{\em see} the explanations below.
This method will only use power series at the origin.

It seems that now, the appropriate linear representation
can be obtained by solving for 
$J_{\ast,\ind}^\kappa\overline{F}$.
But using some of the group parameters $g \in G_\stab^{\kappa-1}$,
some of the power series coefficients 
$J_{\ast,\ind}^\kappa\overline{F}$
may still be normalized, {\em e.g.} to $0$,
and then, associated group reductions must be set up.

Let us assume that such extra
normalizations have been made,
let us keep the same notation 
$J_{\ast,\ind}^\kappa F$ for the remaining independent jets, 
and let us keep the same notation $G_\stab^{\kappa-1}$
for the reduced group.

Once all these tasks are achieved, 
we can really solve:
\[
J_{\ast,\ind}^\kappa\overline{F}
\,=\,
\Lambda\,
\Big(
\Iaux_\smallbullet,J_{\ast,\ind}^\kappa F,g
\Big)
\eqno
{\scriptstyle{(g\,\in\,G_\stab^{\kappa-1})}}.
\]

\begin{Observation}
In all affine structures treated in this memoir,
in~{\cite{Foo-Merker-Nurowski-Ta-2021}}, 
and in other geometric structures as well,
at every jet order $\kappa$, 
these $\Lambda$-formulas {\em always}
were certain {\em explicit
linear matrix representations}
of a certain reduced Lie group $G_\stab^{\kappa-1} \subset
G_\stab^{\kappa-2} \subset \cdots \subset G$,
and even, always independent of the
absolute invariants $\Iaux_\smallbullet$ coming
from the preceding jet orders.
\end{Observation}

Consequently, to each node of the final branching tree
is attached a {\em linear representation} of a Lie group!

This is very analogous to
the existence of $G$-structures with
their successive reductions,
a central feature of Cartan's method of equivalence.
But there is an important difference:
$G$ structures have {\em functional} entries,
while our $\Lambda$-matrices always have {\em scalar}
entries, even when $G$ 
is a group of diffeomorphisms\,\,---\,\,is infinite-dimensional.

This is explained by our key decision not to work in
full jet bundles, but only above successively
selected points or transversals to group-orbits.

So quite unexpectedly for researchers like us who during several years
worked out the {\em very nonlinear} and {\sc pde}-theoretic
(parametric) Cartan equivalence method, {\em the theory of {\em
linear} representations of Lie groups became very useful, very
universal, and present at each step of the process, at every node of
every branching tree!}

To our knowledge, the observation
that linear representations of Lie groups
are universally present
has not been made
in the literature.

We can now terminate our induction reasoning.
The linear representation written above 
of $G_\stab^{\kappa-1}$
in the (finite-dimensional) 
vector space of the components of
$J_{\ast,\ind}^\kappa F$
then decomposes this vector space into a 
finite number of group-orbits.

Transversals $T_\smallbullet^\kappa$ to all these
group-orbits must then be appropriately chosen.

\begin{center}
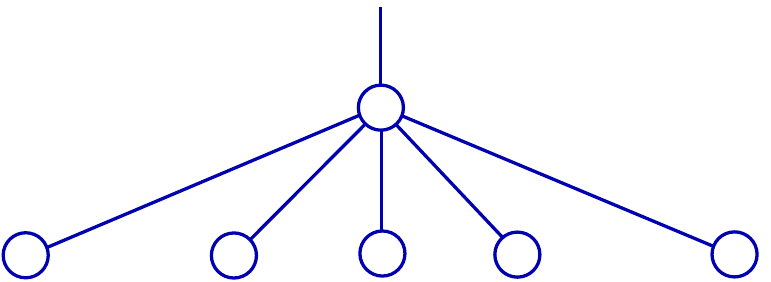
\end{center}

This is how we create the branches at the working order $\kappa$.

This terminates our description of the process,
by induction on $\kappa$, which we 
began on p.~{\pageref{start-induction-kappa}}.
\label{terminate-induction-kappa}

Of course, in the specific examples treated in this
memoir, details are presented,
especially, linear representations.
The reader is referred for instance to 
Sections~{\ref{linear-representations-branches}}.

%%%%%%%%%%%%%%%%%%%%%%%%%%%%%%%%%%%%%%%%%%%%%%%%%%%%%%%%%%%%%%%%%%%%%%
\SectionHead{Determination of Homogeneous Models}
{determination-homgeneous-models}
%%%%%%%%%%%%%%%%%%%%%%%%%%%%%%%%%%%%%%%%%%%%%%%%%%%%%%%%%%%%%%%%%%%%%%

However, some explanations are still missing.
Indeed, because in the present memoir,
we definitely abandon {\sc pde}s,
differential invariants,
functions defined in open sets,
and because we work only within fibers over selected transversals,
it remains to explain how we reach the delicate information about
the {\em dependent} jets $J_{\ast,\dep}^\kappa$. 
In Lagrange style, we work indeed only with (truncated) power series.

For instance, how to express that in Branch~{\green{\bf 2b}},
the Hessian is of constant rank $1$? 
Without using the {\em functional} Hessian?
By anticipation, the reader should
remember that we left untouched the equation:
\[
0
\,=\,
F_{1,2}\,T_1
+
3\,F_{0,3}\,T_2,
\] 
stemming from the
coefficient 
of $y^2$ in~{\eqref{F12-F03-zero}}.
It is from this equation that we will now deduce
$0 = F_{1,2} = F_{0,3}$ obtained
differently in~{\eqref{Hessian-F12-F03-zero}}.

Thus, we now focus our attention mainly on

\begin{Problem}
\label{Pbm-homogeneous-models}
{\bf [Determination of Homogeneous Models]}
{\sl Given a finite-dimensional local Lie group acting $G$ on graphed
submanifolds $M^n = \big\{ u = F(x) \big\}$ in $\R_{x,u}^{n+c}$,
find and classify all possible $M$ having a locally
transitive local automorphisms group $\Sym(M) \subset G$.}
\end{Problem} 
 
{\small (A precise, rigorous definition of a local Lie group,
in Lie's original spirit,
appears in~{\cite{Engel-Lie-Merker-2015, Olver-1986}}.)}
Here:
\[
\Sym(M)
\,=\,
\big\{
g\in G
\colon\,\,
g(M)
\subset
M
\big\},
\]
where we do not stipulate that open subsets 
$V \subset U \subset M$ should be chosen 
with $g(V) \subset U$ and that 
$g \in G$ should lie in some neighborhood of the identity. 

Local Lie groups, not often considered in the modern literature,
are easy to handle because they are well represented
(in a one-to-one manner)
by Lie algebras of vector fields. 

In fact, $\Sym\,M$ has Lie algebra:
\[
\Lie\,\Sym(M)
\,=\,
\mathfrak{sym}(M)
\,:=\,
\big\{
L\in\mathfrak{g}
\colon\,
L\big\vert_M\,\,
\text{is tangent to}\,\,
M
\big\},
\]
where $\mathfrak{g}$ denotes the Lie algebra of 
vector fields {\em inside $\R^{n+c}$}
obtained by diffentiating at the identity
the action of $G$ on $\R^{n+c}$.
For instance, when $G = \Aff(\R^{n+c})$:
\[
\mathfrak{g}
\,=\,
\Span\,
\Big(
\partial_{x_i},\,\,
\partial_{u_j},\,\,
x_{i_1}\,\partial_{x_{i_2}},\,\,
u_j\,\partial_{x_i},\,\,
x_i\,\partial_{u_j},\,\,
u_{j_1}\,\partial_{u_{j_2}}
\Big).
\]

Since all our considerations are {\em local}, 
after recentering the coordinates,
we can assume that everything takes place in some
neighborhood of the origin $0 \in M$.

\begin{Definition}
A $c$-codimensional submanifold $M^n \subset
\R^{n+c}$ is said to be 
{\sl (locally) affinely homogeneous} if:
\[
T_0M
\,=\,
\Span_\R\,
\big\{
L\big\vert_{0}
\colon\,
L
\in
\mathfrak{sym}(M)
\big\}.
\]
\end{Definition}

According to basic Lie theory, the $1$-parameter group
$p \longmapsto \exp(t\,L)(p)$ stabilizes $M$, and 
$\Sym(M)$ is then locally
transitive in a neighborhood of $0 \in M$.

As is known, the datum of
the Lie algebra $\mathfrak{sym} (M)$
enables (by exponentiation)
to reconstitute (a neighborhood of the
identity in) $\Sym(M)$.
But $\mathfrak{sym} (M)$ 
is much better handled than $\Sym (M)$,
thanks to its {\em linear} and
{\em infinitesimal} features.
Lie himself insisted~{\cite{Engel-Lie-Merker-2015}}
on the fact that 
{\em Lie algebras of vector fields}
are the right objects of study when classifying
continuous transformation group actions.
And all of Lie's classifications consist 
in {\em lists} of Lie algebras of {\em infinitesimal
transformations} (vector fields),
{\em see} {\em e.g.} on pages 
6, 17, 26, 57, 71, 106, 116, 139, 167, 203,
209, 214, 226, 246, 257,
271, 334, 370, 388, 384, 388, 391
of~{\cite{Engel-Lie-1893}}.

{\em We will adopt Lie's way of classifying geometries,
namely, by presenting explicit Lie algebras of vector fields.}

Now, in continuation with what precedes, set:
\[
\mathfrak{g}_\stab^{\kappa-1}
\,:=\,
\Lie\,
G_\stab^{\kappa-1}.
\]
Reasoning by induction on the jet order, 
assume that there are vector fields:
\[
e_1,\dots,e_n
\,\in\,
\mathfrak{g}_\stab^{\kappa-1},
\]
such that, at the origin $0 \in M$:
\[
\Span\,
\Big(
e_1\big\vert_0,\,
\dots,\,
e_n\big\vert_0
\Big)
\,=\,
T_0M.
\]
Certainly, $n \leqslant \dim\, \mathfrak{g}_\stab^{\kappa-1}
\leqslant \dim\, G$. 

Together with $e_1, \dots, e_n$, 
there are a certain number $\nu \geqslant 0$ of
{\sl isotropy} vector fields
$f_1, \dots, f_\nu \in \mathfrak{g}_\stab^{\kappa-1}$,
{\em i.e.} vector fields vanishing at 
the origin $(x,u) = (0,0)$, such that 
the general infinitesimal transformation
$L \in \mathfrak{g}_\stab^{\kappa-1}$ writes
for $1 \leqslant j \leqslant c$:
\[
L
\,=\,
T_1\,e_1
+\cdots+
T_n\,e_n
+
A_1\,f_1
+\cdots+
A_\nu\,f_\nu,
\]
with $n+\nu$ arbitrary parameters $T_m$ and $A_\mu$.

To guarantee local homogeneity
(transitivity), no linear relation can ever exist
between $T_1, \dots, T_n$.

The condition that $L$ be tangent to $M$
up to orders $\leqslant \kappa-1$, writes:
\[
\aligned
0
&
\,\equiv\,
L
\Big(
-\,u_j
+
F_j(x)
\Big)
\Big\vert_{u=F(x)}
\\
&
\,\equiv\,
\sum_{i_1+\cdots+i_n\leqslant\kappa-2}
x_1^{i_1}\cdots x_n^{i_n}\,
\underset{{\sf vanish}\atop{\sf by}\,{\sf induction}}{
\zero{\Big(
\cdots
\Big)}}
\\
&
\ \ \ \ \
+
\sum_{i_1+\cdots+i_n=\kappa-1}
x_1^{i_1}\cdots x_n^{i_n}\,
E_{j,i_1,\dots,i_n}^{\sf vf}
\Big(
\Iaux_\smallbullet,
J_\ast^\kappa F,
T_1,\dots,T_n,
A_1,\dots,A_\nu
\Big)
+
{\rm O}_{x_1,\dots,x_n}(\kappa),
\endaligned
\]
that is, after reorganization:
\[
\aligned
0
&
\,\equiv\,
\sum_{m=1}^n\,
T_m\,
\Big(
\Phi_{j,i_1,\dots,i_n,m}^{\sf vf}
\big(
\Iaux_\smallbullet,J_\ast^\kappa F
\big)
\Big)
+
\sum_{\mu=1}^\nu\,
A_\mu\,
\Big(
\Psi_{j,i_1,\dots,i_n,\mu}^{\sf vf}
\big(
\Iaux_\smallbullet,J_\ast^\kappa F
\big)
\Big)
\\
&
\ \ \ \ \ \ \ \ \ \ \ \ \ \ \ \ \ \ \ \ \ \ \ \ \ \ \ \ \ \ \ \ \ \ \
{\scriptstyle{(1\,\leqslant\,j\,\leqslant\,c,\,\,
i_1+\cdots+i_n\,=\,\kappa-1)}}.
\endaligned
\]
A few times below, we will abbreviate these equations as:
\[
0
\,=\,
{\sf E}_\smallbullet^{\sf vf}.
\]

Then exactly as in the example of the coefficient
of $y^2$ in~{\eqref{F12-F03-zero}},
whenever one of these equations,
say for some indices $\underline{j}, 
\underline{i}_1, \dots, \underline{i}_n$, 
does not incorporate any of 
the {\sl isotropy parameters}
$A_1, \dots, A_\nu$,
but incorporates only the {\sl transitivity parameters}
$T_1, \dots, T_n$, we receive $n$ equations:
\[
0
\,=\,
\Phi_{\underline{j},\underline{i}_1,\dots,\underline{i}_n,m}^{\sf vf}
\big(
\Iaux_\smallbullet,J_\ast^\kappa F
\big)
\eqno
{\scriptstyle{(1\,\leqslant\,m\,\leqslant\,n)}},
\]
which are of degree 1 with respect to $J_\ast^\kappa F$,
and which express constraints on
certain `dependent' jets $J_{\ast,\dep}^\kappa F$
to be resolved in terms
of certain other `independent' jets 
$J_{\ast,\ind}^\kappa F$.

Some of these `independent' jets may 
simultaneously become absolute invariants at order $\kappa$,
hence join the 
current collection $\Iaux_\smallbullet$ before
passing to order $\kappa+1$.

Sometimes even, some linear combinations between
these equations must be performed in some tricky way
in order to {\em eliminate} $A_1, \dots, A_\nu$,
so as to `discover' further {\sl transitivity equations}
which would reveal new constraints.
In many branches of our classification of affinely homogeneous
surfaces $S^2 \subset \R^4$, we were blocked 
for this reason.

{\small This method 
based on transitivity equations
has already been applied 
in~{\cite{Foo-Merker-Nurowski-Ta-2021}},
in a {\em degenerate} CR-geometric context,
for the infinite-dimensional group of biholomorphisms
of $\C^3$. However, no details of proof 
were given in~{\cite{Foo-Merker-Nurowski-Ta-2021}}.
A complete written proof would be about 50 pages long,
due to subtle computational aspects in degenerate branches.
Indeed, even for
a reduction to an {\em explicit}, 
{\sl parametric}, Cartan-type $\{e\}$-structure,
which is a preliminary step to determine homogeneous models, 
the calculations are 
long~{\cite{Merker-Pocchiola-2018, Foo-Merker-2019}}.
Similarly, in the degenerate para-CR context, by lack of space,
several computations used to determine homogeneous geometries
are not fully presented in~{\cite{Merker-Nurowski-2020-a, 
%%Merker-Nurowski-2020-b,
Merker-Nurowski-2023}}.}

This concludes our presentation of how we can
reach some differential information
in successive jet fibers over successive
group-transversals, 
even if we work only with truncated power series
at the origin, as Lagrange advocated to do.
In the affine context, a similar approach
though presented differently,
was developed in~{\cite{Eastwood-Ezhov-1999, 
Eastwood-Ezhov-2001-1, Eastwood-Ezhov-2001-2}}.

Contrary to what can be believed at first sight, 
this step of the process,
namely to extract transitivity equations,
is not (at all) straightfoward.
Experts know well how complicated algebras of differential
invariants can be, and here, when searching for homogeneous
models, the exploration
is confronted with such kinds of algebras.
Indeed, differential invariants being constant,
their invariant differentials in
the Lie-Fels-Olver recurrence formulas
vanish, but the other (complicated) algebraic 
terms in these formulas
are essentially
the same as the one we manipulate in this memoir.

Moreover, in many branches of our classifications,
the resolutions of these {\sl transitivity equations}
was delicate, required clever choices of variables
to be solved (proceeding by hand by solving
one variable at a time),
demanded to be able to go back to preceding
jet orders, 
{\em etc.}, exactly as in 
Cartan's equivalence method. 
In the literature, 
such delicate
aspects of computational explorations are almost never 
described or even mentioned.

Lastly, and importantly, at the end of the process,
we often obtain a collection of {\em algebraic} 
equations in the remaining absolute invariants
$\Iaux_\smallbullet$, 
some {\em key equations} whose zero-set defines
an {\em algebraic moduli space} of
a collection of homogeneous models, 
represented by a {\sl terminal leaf}
of the tree. 
The reader can have a look at
Model~{\green{\bf 2e3a4a}} in
Section~{\ref{2e-models}}.
We will
say more about terminal leaves in
Section~{\ref{termination-moduli-spaces-homogeneous-models}}.

We end up this lengthy 
description by briefly coming back
to Problem~{\ref{Pbm-debut-diff-invts}}.

%%%%%%%%%%%%%%%%%%%%%%%%%%%%%%%%%%%%%%%%%%%%%%%%%%%%%%%%%%%%%%%%%%%%%%
\SectionHead{Creations of Geometries}
{creations-geometries}
%%%%%%%%%%%%%%%%%%%%%%%%%%%%%%%%%%%%%%%%%%%%%%%%%%%%%%%%%%%%%%%%%%%%%%

By what precedes, the equations $0 = {\sf E}_\smallbullet^{\sf nf}$
together with 
the equations $0 = {\sf E}_\smallbullet^{\sf vf}$
are used to determine {\sl homogeneous geometries},
namely submanifolds $M \subset \R^{n+c}$ 
having (locally) transitive symmetry group 
$\Sym(M)$, jet order after jet order. 
These equations
are responsible for the creation of a certain 
{\sl branching tree}. In principle,
the {\sl terminal leaves} of this tree
correspond to (families of) homogeneous models.

To each node of the branching tree, there 
is associated a certain linear representation
of a certain subgroup $G' \subset G$ on a certain
vector space $V'$ coordinatized by
certain (independent)
jet coefficients of $F$.
As we already explained, 
from this node are born as many edges towards the
next jet order as there are 
transversals to $G'$-orbits in $V'$.

But instead of repeating in the next jet order
the use of the equations
$0 = {\sf E}_\smallbullet^{\sf nf}$ 
and $0 = {\sf E}_\smallbullet^{\sf vf}$
to continue to develop the branching
tree of homogeneous models,
we can {\em stop} the 
transitive analysis at this point.
We can take each created edge as the 
departure for a {\em new subgeometry}, 
without continuing the tree, even
without knowing what could happen next.

Indeed, in all the preceding jet orders,
there were certain (relative) differential invariants
which were assumed to be
zero at the origin. 
And to each one of these punctual invariants
there corresponded a (relative) differential invariant,
as is clearly explained in~{\cite{Olver-2018, Chen-Merker-2019,
Chen-Merker-2020, Olver-Sabzevari-Valiquette-2023}}.
Denote these differential invariants
as $\Kaux_1, \dots, \Kaux_t$.

So these $\Kaux_\smallbullet$ are assumed to vanish
at the origin $0 \in M$. Of course, 
they can take nonzero values nearby,
a situation that could be treated by Singularity Theory.
But as we decided to study only {\em constant-type}
geometries, adopting Lie's principle of 
thought~{\cite[Chap.~1]{Engel-Lie-Merker-2015}},
we are led to assume that:
\[
0
\,\equiv\,
\Kaux_1\big(x,J^\smallbullet F(x)\big)
\,\equiv\,\cdots\,\equiv\,
\Kaux_t\big(x,J^\smallbullet F(x)\big),
\]
for $x$ in some neighborhood of the origin.

Thus, we can stop the homogeneous geometries
process
$0 = {\sf E}_\smallbullet^{\sf nf} = {\sf E}_\smallbullet^{\sf vf}$
anywhere.

\begin{Principle}
{\bf [Creation of constant-type (degenerate) geometries]}
{\sl 
Given a group $G$ acting transitively 
on graphs $\big\{ u = F(x) \big\}$
in $\R^{n+c} \ni (x,u)$, with its prolonged actions
to jet bundles $J_{n,c}^1$, $J_{n,c}^2$, \dots, 
$J_{n,c}^\kappa$,
\dots, at each order $\kappa \geqslant 1$, 
at each node of the branching tree 
(even if incomplete)
which is constructed to determine homogeneous
models, create (introduce) new geometries,
of constant type, degenerate in a certain sense,
depending on the history of the node.}
\end{Principle}

Some of the nodes are such that all the power series coefficients
of $F$ 
are already uniquely determined,
especially the final nodes, {\em i.e.} the terminal leaves.

Some other nodes are such that there still remain
infinitely many power series $F$-coefficients
which are free, not normalized,
and then,
the $G'$-action must be prolonged to 
the jet (sub)bundle of this (sub)geometry,
in order to determine the corresponding algebras
of differential invariants.

In conclusion, 
{\em many new geometries having algebras of differential
invariants exist which should (can) be studied}.

Most of the times, the creation of constant-type geometries
is well known at jet order 2.
For instance, under the group
$\Aff(\C^{n+1})$ of
affine transformations
of $\C^{n+1}$\,\,---\,\,codimension $c = 1$\,\,---\,\,since 
the punctual rank of the Hessian matrix 
is invariant, inequivalent graphed normal forms are:
\[
u
\,=\,
x_1^2
+\cdots+
x_m^2
+
{\rm O}_{x_1,\dots,x_n}(3),
\]
with an invariant
integer $0 \leqslant m \leqslant n$, 
which produces $n+1$ different (inequivalent!)
geometric structures. 
Similarly, for hypersurfaces, 
the rank and the signature
of the Levi form are invariant under CR equivalences,
hence several order 2 geometries can be `created'.

Applying his theory of moving frames,
Peter Olver~{\cite{Olver-2007}} 
studied algebras of differential invariants
for elliptic and hyperbolic surfaces $S^2 \subset \R^3$
under Euclidean and Affine transformations,
{\em i.e.} with Hessians of maximal rank 2,
{\em see}~{\cite{Chen-Merker-2019, Arnaldsson-Valiquette-2020}}
for the Hessian rank 1 geometry.
To study only a single one of these 
constant Hessian rank affine geometries
from the point of view of differential invariants, 
for instance with $n = 5$ and $m = 3$
(a case probably never looked at), 
might already be a considerable task.

Constant type (degenerate) geometries
at jet order $\geqslant 3$ are not much
studied, but they are as legitimate
as the order 2 (degenerate) geometries.
The branching tree
in~{\cite{Merker-Nurowski-2020-a}}
shows certain degenerate para-CR 
geometries of jet order $> 3$,
{\em i.e.} beyond Levi form 
(which is of order 2)
and beyond $2$-nondegeneracy
(which of order 3).

\begin{Problem}
\label{Pbm-algebras-differential-invariants}
{\bf [Algebras of differential invariants
for degenerate geometries]}
{\sl Describe algebras of differential invariants
of constant-type degenerate geometries.
Find minimal sets of (differential) generators.
}
\end{Problem} 

We insist on the fact that we formulate this
general problem
for {\em all possible constant-type degenerate (sub)geometries}.
Even, 
the considered 
Lie group $G$ can be an infinite-dimensional Lie pseudo-group.

At the opposite are the {\em generic} geometries,
those for which it is allowed to assume that
some functions, some determinants, are nonzero, some
rank matrices are maximal, {\em etc.}
For some generic  geometries,
under some classical groups, 
Peter Olver~{\cite{Olver-2007}},
Hubert-Olver~{\cite{Hubert-Olver-2007}},
have established
remarkable theorems that a {\em single}
differential invariant is sufficient
to (differentially) generate the whole algebra
of differential invariants.

{\em But certainly, the genericity of a geometry is a {\em relative}
concept!} 
{\em Genericity also concerns subgeometries!}

Indeed, in any node at which a constant-type
(degenerate) geometry is created,
by assuming that all higher order encountered 
(relative) differential
invariants are nonvanishing (after restriction to open
subsets), by assuming in addition if it is convenient 
that some functions, some determinants, {\em etc.},
are nonzero, then a certain `{\sl generic}' (sub)geometry can 
be defined 
within the considered degenerate geometry.

For some degenerate 2D affine of 3D Cauchy-Riemann 
geometries
which are (relatively) generic,
Arnaldsson-Valiquette~{\cite{Arnaldsson-Valiquette-2020}},
Olver-Sabzevari-Valiquette~{\cite{Olver-Sabzevari-Valiquette-2023}},
have established the 
analogous theorems that a {\em single}
differential invariant is sufficient
to (differentially) generate the whole algebra
of differential invariants\,\,---\,\,an aspect
not touched in~{\cite{Chen-Merker-2019,
Chen-Merker-2020}}.

In forthcoming publications, we will endeavour to 
explore such algebras of differential invariants
on a number of further examples of {\em degenerate} geometries.
We will also explain how to obtain
differential consequences 
of assumptions like $0 \equiv \Kaux_1 \equiv \cdots \equiv 
\Kaux_t$, 
by working 
only with power series at the origin.

Also, we will explain how to obtain very quickly the
so-called Lie-Fels-Olver recurrence formulas
for differential invariants 
of some degenerate constant-type geometries,
{\em directly
from the {\sf v}ector {\sf f}ields 
equations
$0 = {\sf E}_\smallbullet^{\sf vf}$}.

%%%%%%%%%%%%%%%%%%%%%%%%%%%%%%%%%%%%%%%%%%%%%%%%%%%%%%%%%%%%%%%%%%%%%%
\SectionHead{General Method in Affine Context}
{general-method-affine-context}
%%%%%%%%%%%%%%%%%%%%%%%%%%%%%%%%%%%%%%%%%%%%%%%%%%%%%%%%%%%%%%%%%%%%%%

Let us now be more specific. In this memoir,
for certain small dimensions $n$ and small codimensions $c$,
we shall determine all affinely homogenous
submanifolds $M^n \subset \R^{n+c}$. 
Branches will be inequivalent, 
except in some cases, where 
equivalence migh hold up to the action of
some finite subgroup of the affine group,
as in most of Peter Olver's works.

In this section, we introdude general notations
and we express the equations 
$0 = {\sf E}_\smallbullet^{\sf nf} = 
{\sf E}_\smallbullet^{\sf vf}$ explicitly.

Let $n \geqslant 1$, let $x = (x_1, \dots, x_n) \in \R^n$,
let $c \geqslant 1$, and let 
$(u_1, \dots, u_c) \in \R^c$.
These are the {\sl source} coordinates in $\R^{n+c}$.
Similarly, let
$y = (y_1, \dots, y_n) \in \R^n$ 
and let $v = (v_1, \dots, v_c) \in \R^c$
be {\sl target} coordinates in $\R^{n+c}$\,\,---\,\,instead
of $\overline{x}$, $\overline{u}$.
Throughout, $\R$ can be replaced by $\C$. 

We consider equivalences of local analytic submanifolds
$M^n \subset \R_{x,u}^{n+c}$ and $N^n \subset \R_{y,v}^{n+c}$
graphed as:
\[
\aligned
u_j
&
\,=\,
F_j(x_1,\dots,x_n)
\ \ \ \ \ \ \ \ \ \ \ \ \ \ \ \ \ \ \ \
\text{and}
\ \ \ \ \ \ \ \ \ \ \ \ \ \ \ \ \ \ \ \
v_j
\,=\,
G_j(y_1,\dots,y_n),
\\
&
\ \ \ \ \ \ \ \ \ \
{\scriptstyle{(1\,\leqslant\,j\,\leqslant\,c)}}
\ \ \ \ \ \ \ \ \ \ \ \ \ \ \ \ \ \ \ \ \ \ \ \ \ \ \ \ \ \ \ \ \ \ \
\ \ \ \ \ \ \ \ \ \ \ \ \ \ \ \ \ \ \ \ \ \ \ \ \ \ \ \ \ \ \ \ 
{\scriptstyle{(1\,\leqslant\,j\,\leqslant\,c)}}
\endaligned
\]
under {\sl affine transformations} 
$\R^{n+c} \longrightarrow \R^{n+c}$:
\reqnomode\usetagform{EngelLie}
\begin{align}
y_i
&
\,=\,
\sum_{1\leqslant i'\leqslant n}\,
a_{i,i'}\,x_{i'}
+
\sum_{1\leqslant j'\leqslant c}\,
b_{i,j'}\,u_{j'}
+
\tau_i
\tag{(1\,\leqslant\,i\,\leqslant\,n),}
\\
v_j
&
\,=\,
\sum_{1\leqslant i'\leqslant n}\,
c_{j,i'}\,x_{i'}
+
\sum_{1\leqslant j'\leqslant c}\,
d_{j,j'}\,u_{j'}
+
\upsilon_j
\tag{(1\,\leqslant\,j\,\leqslant\,c),}
\end{align}
with translational constants $\tau_i$, $\upsilon_j$
(upsilon-j), 
where the $(n+c) \times (n+c)$ linear-part
$\big( \begin{smallmatrix} a & b \\ c & d \end{smallmatrix} \big)$
matrix belongs to 
$\GL(\R^{n+c})$, {\em i.e.} has nonzero determinant.
The collection of all these
transformations is the (continuous Lie transformation) group
$\Aff(\R^{n+c})$. 

Abbreviate:
\[
\aligned
y
&
\,=\,
a\,x
+
b\,u
+
\tau,
\\
v
&
\,=\,
c\,x
+
d\,u
+
\upsilon.
\endaligned
\]
Using the translation parameters $\tau_i$ and $\upsilon_j$,
we may assume that the origin is mapped to the origin,
whence $0 = \tau = \upsilon$, which reduces
the group $\Aff(\R^{n+c})$ to $\GL(\R^{n+c})$.

Then the fact that $\big\{ u = F(x) \big\}$ is mapped
to $\big\{ v = G(y) \big\}$ is equivalent to:
\[
0
\,\equiv\,
-\,c\,x
-
d\,F(x)
+
G
\big(
a\,x+b\,F(x)
\big),
\]
identically for all $x \in \R^n$ near the origin.

Expanding these $c$ equations in power series, we obtain:
\[
0
\,\equiv\,
\sum_{\kappa=0}^\infty\,
\sum_{i_1+\cdots+i_n=\kappa}\,
x_1^{i_1}\cdots x_n^{i_n}
\cdot
{\sf E}_{j,i_1,\dots,i_n}^{\sf nf}
\Big(
J^\kappa F,\,
J^\kappa G,\,
a,b,c,d
\Big).
\]
Thus, the {\sf n}ormal {\sf f}orm equations are:
\[
0
\,=\,
{\sf E}_{j,i_1,\dots,i_n}^{\sf nf}
\Big(
J^{i_1+\cdots+i_n}F,\,
J^{i_1+\cdots+i_n}F,\,
a,b,c,d
\Big)
\eqno
{\scriptstyle{(1\,\leqslant\,j\,\leqslant\,c,\,\,
i_1,\dots,i_n\,\geqslant\,0)}}.
\]
Though algebraic, they are highly nonlinear and often quite
difficult to manipulate.

One of Lie's deep ideas was to {\em infinitesimalize} 
calculations.

The Lie algebra $\mathfrak{aff} (\R^{n+c})$ of $\Aff(\R^{n+c})$
consists of the vector fields:
\[
\aligned
L
&
\,=\,
T_1\,
\tfrac{\partial}{\partial x_1}
+\cdots+
T_n\,
\tfrac{\partial}{\partial x_n}
+
U_1\,
\tfrac{\partial}{\partial u_1}
+\cdots+
U_c\,
\tfrac{\partial}{\partial u_c}
\,+
\\
&
\ \ \ \ \
+
\sum_{i=1}^n\,
\Big(
\smallsum{1\leqslant i'\leqslant n}\,
A_{i,i'}\,x_{i'}
+
\smallsum{1\leqslant j'\leqslant c}\,
B_{i,j'}\,u_{j'}
\Big)\,
\tfrac{\partial}{\partial x_i}
\\
&
\ \ \ \ \
+
\sum_{j=1}^c\,
\Big(
\smallsum{1\leqslant i'\leqslant n}\,
C_{j,i'}\,x_{i'}
+
\smallsum{1\leqslant j'\leqslant c}\
D_{j,j'}\,u_{j'}
\Big)\,
\tfrac{\partial}{\partial u_j},
\endaligned
\]
with free real constants $T_\smallbullet$, $U_\smallbullet$,
$A_{\smallbullet, \smallbullet}$,
$B_{\smallbullet, \smallbullet}$,
$C_{\smallbullet, \smallbullet}$,
$D_{\smallbullet, \smallbullet}$.

After an elementary affine transformation in the $(x,u)$-space,
we can assume that $M = \big\{ u = F(x) \big\}$ 
is graphed over $\R_x^n$ with {\em horizontal} tangent
plane at the origin $0 \in M$, 
so that $F(x) = {\rm O}_x(2)$.
Then the tangency of $L\big\vert_M$ to $M$ 
at the origin requires:
\[
U_1
\,=\,
\cdots
\,=\,
U_c
\,=\,
0.
\]

Therefore:
\[
L\Big\vert_0
\,=\,
T_1\tfrac{\partial}{\partial x_1}
\Big\vert_0
+\cdots+
T_n\tfrac{\partial}{\partial x_n}
\Big\vert_0,
\]
and in the search for locally 
affinely homogeneous models $M^n \subset \R^{n+c}$,
{\em no linear relation} can ever exist between
the {\sl transitivity parameters} $T_1, \dots, T_n$.

Then the tangency of $L\big\vert_M$ to $M$ is equivalent to:
\reqnomode\usetagform{EngelLie}
\begin{align}
0
&
\,\equiv\,
-\,
\smallsum{1\leqslant i'\leqslant n}\,
C_{j,i'}\,x_{i'}
-
\smallsum{1\leqslant j'\leqslant c}\,
D_{j,j'}\,F_{j'}(x)
\notag
\\
&
\ \ \ \ \
+
\sum_{i=1}^n\,
\Big(
T_i\,
+
\smallsum{1\leqslant i'\leqslant n}\,
A_{i,i'}\,x_{i'}\,
+
\smallsum{1\leqslant j'\leqslant c}\,
B_{i,j'}\,F_{j'}(x)
\Big)\,
F_{j,x_i}(x)
\tag{(1\,\leqslant\,j\,\leqslant\,c),}
\end{align}
identically for all $x \in \R^n$ near the origin.

Expanding these $c$ equations in power series, we obtain:
\[
0
\,\equiv\,
\sum_{\kappa=0}^\infty\,
\sum_{i_1+\cdots+i_n=\kappa}\,
x_1^{i_1}\cdots x_n^{i_n}
\cdot
{\sf E}_{j,i_1,\dots,i_n}^{\sf vf}
\Big(
J^{1+\kappa} F,\,
J^{1+\kappa} G,\,
T,\,
A,B,C,D
\Big).
\]
Thus, the {\sf v}ector {\sf f}ields equations are:
\[
0
\,=\,
{\sf E}_{j,i_1,\dots,i_n}^{\sf vf}
\Big(
J^{1+i_1+\cdots+i_n}F,\,
T,\,
A,B,C,D
\Big)
\eqno
{\scriptstyle{(1\,\leqslant\,j\,\leqslant\,c,\,\,
i_1,\dots,i_n\,\geqslant\,0)}}.
\]
These equations are all linear in $T$, $A$, $B$, $C$, $D$.

Special equations which depend only on $T_1, \dots, T_n$
are called {\sl transitivity equations},
for instance for some indices
$\underline{j}, \underline{i}_1, \dots, \underline{i}_n$:
\[
0
\,\equiv\,
\sum_{m=1}^n\,
T_m\,
\Big(
\Phi_{\underline{j},\underline{i}_1,\dots,\underline{i}_n,m}^{\sf vf}
\big(
J^{1+\underline{i}_1+\cdots+\underline{i}_n}F
\big)
\Big),
\]
because when searching for locally homogeneous
structures (our goal), we can deduce that the
$n$ coefficients of $T_1, \dots, T_n$ vanish:
\[
0
\,=\,
\Phi_{\underline{j},\underline{i}_1,\dots,\underline{i}_n,m}^{\sf vf}
\big(
J^{1+\underline{i}_1+\cdots+\underline{i}_n}F
\big)
\eqno
{\scriptstyle{(1\,\leqslant\,m\,\leqslant\,n)}}.
\]

In all affinely homogeneous models obtained in this memoir,
we apply the following algorithm, proceeding by induction
on the jet order $\kappa \geqslant 1$, playing ``ping-pong'' 
between $0 = {\sf E}_\smallbullet^{\sf nf}$
and $0 = {\sf E}_\smallbullet^{\sf vf}$.

\medskip\noindent{\footnotesize\sf Step~1.}
At order $\kappa-1$, analyze the {\sf v}ector {\sf f}ield
equations $0 = {\sf E}_\smallbullet^{\sf vf}$ 
to determine all {\sl transitivity equations},
so as to see which of the order $\kappa$ jet coefficients 
$J_\ast^\kappa F$ are dependent:
$J_{\ast,\dep}^\kappa F$.
Solve these $J_{\ast,\dep}^\kappa F$.
Only the remaining independent jets remain: 
$J_{\ast,\ind}^\kappa F$.

Of all the {\footnotesize\sf Steps~1, 2, 3, 4, 5},
this {\footnotesize\sf Step~1} is the most
time-consuming, from the computational point of view.
In some branches, 
{\em this {\footnotesize\sf\em Step~1} 
might happen to be very delicate,
not at all straightforward}.

\medskip\noindent{\footnotesize\sf Step~2.}
At the next jet order $\kappa$,
pass to the {\sf n}ormal {\sf f}orm
equations 
$0 = {\sf E}_\smallbullet^{\sf nf}$,
insert the transitivity constraints
obtained in {\footnotesize\sf Step~1}\,\,---\,\,which are
invariant!\,\,---~, simultaneously for the
source $\{ u = F(x) \}$ and for the 
source $\{ v = G(y) \}$.
Then some of the isotropy group parameters
which still remain free at this jet order $\kappa$,
namely some of the 
$a_{\smallbullet, \smallbullet}$,
$b_{\smallbullet, \smallbullet}$,
$c_{\smallbullet, \smallbullet}$,
$d_{\smallbullet, \smallbullet}$,
which have not been yet used at orders $\leqslant \kappa-1$,
can happen to be appropriate 
to normalize, at order $\kappa$,
some of the coefficients $J_\ast^\kappa G$,
and then symmetrically by invariancy,
to normalize {\em in exactly the same way}
the corresponding coefficients $J_\ast^\kappa F$. 

In some branches, 
{\em this Step~2 might be simply computationally blocked,
due to the nonlinear complexity of the equations
$0 = {\sf E}_\smallbullet^{\sf nf}$}.

\medskip\noindent{\footnotesize\sf Step~3.}
Provided that {\footnotesize\sf Step~1} and~{\footnotesize\sf Step~2}
have been achieved correctly, without any oblivion,
determine the {\sl stability} subgroup $G_\stab^{\kappa-1}
\subset G$ of affine transformations which stabilizes 
not only the normalizations done up to orders
$\leqslant \kappa-1$, but also the ones done
in~{\footnotesize\sf Step~2}. Then
obtain a certain {\em linear representation}
of $G_\stab^{\kappa-1}$ on the vector space
$\R^{\ell_\kappa}$ of the remaining unnormalized
coefficients $J_{\ast,\ind}^\kappa F$,
with target coefficients $J_{\ast,\ind}^\kappa G$:
\[
\Big(
J_{\ast,\ind}^\kappa G
\Big)
\,:=\,
\Big(
\Pi(a,b,c,d)
\Big)\,
\Big(
J_{\ast,\ind}^\kappa F
\Big),
\]
with a certain explicit $\ell_\kappa \times \ell_\kappa$
matrix $\Pi$ depending on the
(currently remaining) isotropy group parameters $(a,b,c,d) 
\in G_\stab^{\kappa-1}$. We have no {\em general proof}
that a linear representation necessarily appears,
but in all nodes of all trees treated in this memoir,
we always obtained linear representations,
possibly after corrections of mistakes in 
{\footnotesize\sf Steps~1} and~{\footnotesize\sf 2}

In some branches, 
this {\footnotesize\sf Step~3}
might lead to a linear representation 
which does not show its invariant subspaces,
hence a search for a change of basis
in $\C^{\ell_\kappa}$ adapted to the existing
invariant subspaces 
must be performed before going further.

\medskip\noindent{\footnotesize\sf Step~4.}
Decompose $\R^{\ell_\kappa}$ into orbits
under $G_\stab^{\kappa-1}$. Find `nice' transversals 
$T_\smallbullet^\kappa$ 
to all of these orbits. 

In some branches, 
this {\footnotesize\sf Step~4} cannot be passed
unless the linear representation has been
`simplified' in some way, {\em e.g.}
by changing the basis of $\R^{\ell_\kappa}$.

\medskip\noindent{\footnotesize\sf Step~5.}
Create as many branches as there are orbit transversals 
$T_\smallbullet^\kappa$.
Using the {\sf n}ormal {\sf f}orm equations
$0 = {\sf E}_\smallbullet^{\sf nf}$,
determine all stability subgroups
$G_\stab^\kappa \subset G_\stab^{\kappa-1}$ associated with all
transversals $T_\smallbullet^\kappa$. 
Restrict any further analysis to each one of these
$T_\smallbullet^\kappa$, that is,
work only with higher jets 
$J^{\kappa+1} \big\vert_{ T_\smallbullet^\kappa}$,
$J^{\kappa+2} \big\vert_{ T_\smallbullet^\kappa}$,
\dots,
restricted to these transversals. 

Then pass to {\footnotesize\sf Step~1} at order $\kappa$.

\medskip

The reader
interested in {\em termination} of the 
{\footnotesize\sf Steps~1-2-3-4-5} `algorithm'
can skip the next 
Section~{\ref{classification-affinely-homogeneous-S2-C3}}, 
and directly `jump' to 
Section~{\ref{termination-moduli-spaces-homogeneous-models}}.

%%%%%%%%%%%%%%%%%%%%%%%%%%%%%%%%%%%%%%%%%%%%%%%%%%%%%%%%%%%%%%%%%%%%%%
\SectionHead{A Classification of Affinely Homogeneous
$S^2 \subset \C^3$}
{classification-affinely-homogeneous-S2-C3}
%%%%%%%%%%%%%%%%%%%%%%%%%%%%%%%%%%%%%%%%%%%%%%%%%%%%%%%%%%%%%%%%%%%%%%

Since this memoir is devoted to higher (co)dimensional
affine homogeneity, it is advisable to explain,
for the well studied~{\cite{Guggenheimer-1963,
Jensen-1977,
Abdalla-Dillen-Vrancken-1997,
Doubrov-Komrakov-Rabinovich-1996,
Doubrov-Komrakov-1998,
Eastwood-Ezhov-1999,
Olver-2007,
Chen-Merker-2019,
Arnaldsson-Valiquette-2020,
Chen-Merker-2020}} 
case of surfaces $S^2 \subset \C^3$,
what would be the terminal leaves of a classification 
in Lie's spirit.

Remember that we abandon the search for closed forms
of homogeneous models, and that our goal is to obtain
Lie algebras of vector fields\,\,---\,\,sometimes
`parametrized' by a number of absolute invariants 
$\Iaux_\smallbullet$.

\begin{center}
\label{diag-arbre-S2-R3}
\includegraphics[scale=0.17]{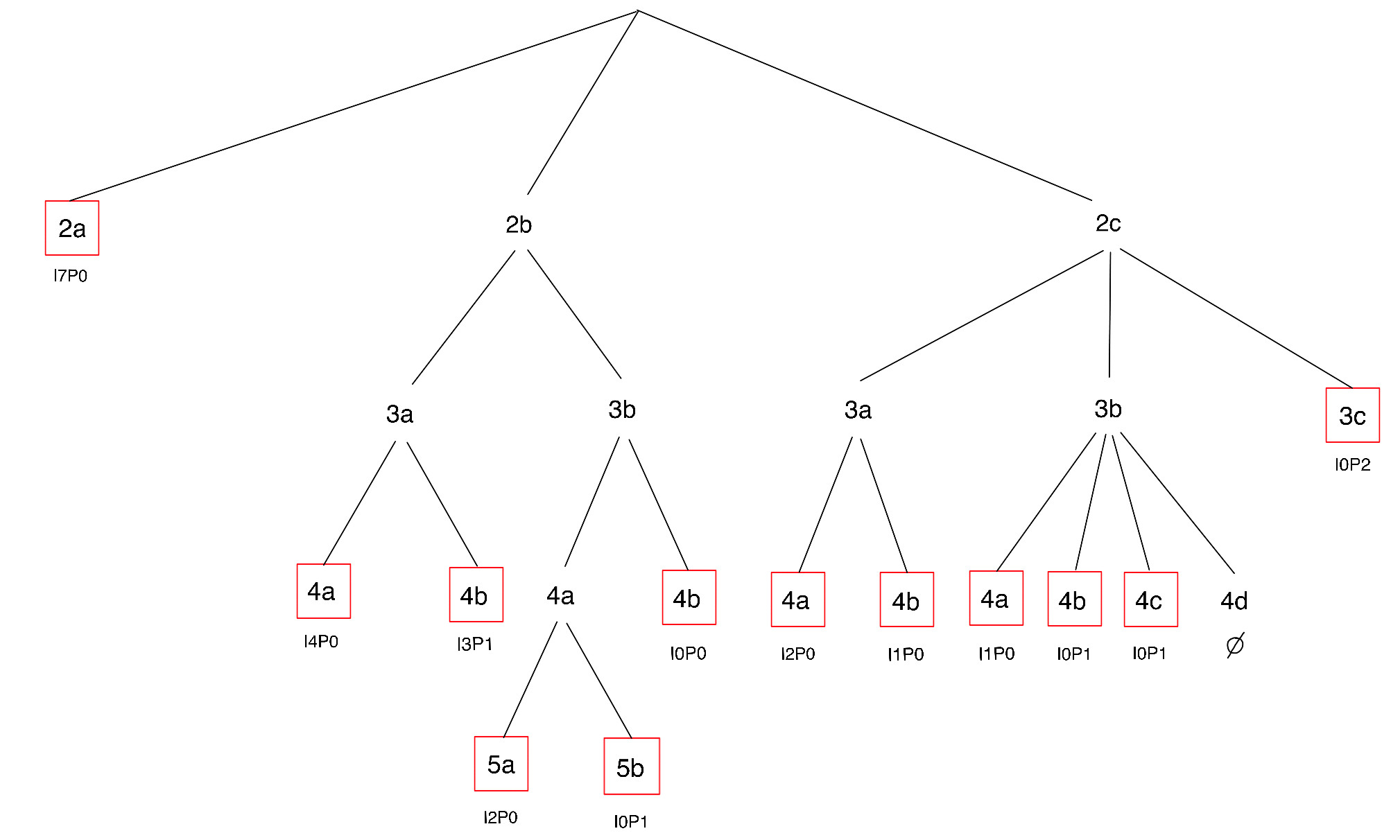}
\end{center}

The branching tree has 13 terminal leaves, with 1 of them,
{\green{\bf 2c3b4d}}, leading to
a contradiction (no homogeneous model).
In the notation I3P1, I3 means 3-dimensional Isotropy,
and 
P1 means
1-dimensional moduli space algebraic variety,
as the number of Parameters\,\,---\,\,often 
just $\C^1$. Further, P0 means
a single model, and rarely, a finite collection of
individual models, and certainly with zero parameter
(numerical constants only).

Remember the little array of creation 
of order 2 transversals:
\[
\def\arraystretch{1.25}
\begin{array}{rccc}
\green{\bf 1}\,\,\,\,
\green{\downarrow}\,\,
& 
F_{2,0} & F_{1,1} & F_{0,2}
\\
\green{\bf 2a} & 
0 & 0 & 0
\\
\green{\bf 2b} & 
1 & 0 & 0
\\
\green{\bf 2c} & 
0 & 1 & 0
\end{array}
\]
which can be translated as:
\[
\aligned
{}
&
{\green{\bf 2a}:}
&
\ \ \ \ \ \ \ \ \ \ \ \ \ \ \ \ \ \ \ \
u
&
\,=\,
0
+
0
+
{\rm O}_{x,y}(3),
\\
{}
&
{\green{\bf 2b}:}
&
\ \ \ \ \ \ \ \ \ \ \ \ \ \ \ \ \ \ \ \
u
&
\,=\,
0
+
x^2
+
{\rm O}_{x,y}(3),
\\
{}
&
{\green{\bf 2c}:}
&
\ \ \ \ \ \ \ \ \ \ \ \ \ \ \ \ \ \ \ \
u
&
\,=\,
0
+
x\,y
+
{\rm O}_{x,y}(3).
\endaligned
\]

As explained above, to each node of the branching diagram,
there corresponds a linear representation.
Before explaining computational details about 
how such linear representations appear in some
of the branches (not in all),
let us list the linear representations, together with the 
orbit-transversals they create.
At the final nodes consisting of boxed terminal leaves,
there are only trivial, identity linear representations.
Notation is as in Sections~{\ref{parabolic-surfaces-S2-C3}} 
and~{\ref{non-reduced-linear-representation}}.

As we already saw in Section~{\ref{setting-up-dependent-jets}},
at the node {\green{\bf 2b}}, with:
\[
u
\,=\,
x^2
+
F_{2,1}\,x^2\,y
+
{\rm O}_{x,y}(4)
\ \ \ \ \ \ \ \ \ \ \ \ \ \ \ \ \ \ \ \
\text{and}
\ \ \ \ \ \ \ \ \ \ \ \ \ \ \ \ \ \ \ \
v
\,=\,
q^2
+
G_{2,1}\,p^2\,q
+
{\rm O}_{p,q}(4),
\] 
the current reduced group consists of matrices of the form:
\[
\left[\,
\begin{array}{ccc}
a_{1,1} & \red{\bf 0} & b_1
\\
a_{2,1} & a_{2,1} & b_1
\\
\red{\bf 0} & \red{\bf 0} & \red{a_{1,1}^2}
\end{array}
\,\right],
\ \ \ \ \ \ \ \ \ \ \ \ \ \ \ \ \ \ \ \
0
\,\neq\,
a_{1,1}^3\,a_{2,2},
\]
and the linear representation is:
\[
\Big(
G_{2,1}
\Big)
\,=\,
\Big(
\tfrac{1}{a_{2,2}}
\Big)\,
\big(
F_{2,1}
\big).
\]

This means that $F_{2,1}$ is a {\em punctual} relative invariant,
{\see}~{\cite{Chen-Merker-2019}}
for the corresponding differential invariant, denoted 
there:
\[
\Saux
\,:=\,
\frac{F_{xx}\,F_{xxy}-F_{xy}\,F_{xxx}}{F_{xx}^2}.
\]

In this case, branch creation is immediate:
\[
\def\arraystretch{1.25}
\begin{array}{rc}
\green{\bf 2b}\,\,\,\,
\green{\downarrow}\,\,
& 
F_{2,1} 
\\
\green{\bf 3a} & 
0 
\\
\green{\bf 3b} & 
1 
\end{array}
\]
In Branch~{\green{\bf 2b3a}}, no group reduction occurs,
while in Branch~{\green{\bf 2b3b}}, after having normalized
$G_{2,1} := 1 =: F_{2,1}$, it necessarily comes the 
group reduction $a_{2,2} := 1$.

Next, at the node {\green{\bf 2b3a}}, with the order 3
power series coefficients immediately read off from the
above array, an analysis 
of the (here not delicate) equations
${\sf E}_\smallbullet^{\sf vf}$ shows that 
there are no other order 4 monomials
than $F_{4,0}\,x^4$:
\[
u
\,=\,
x^2
+
0
+
F_{4,0}\,x^4
+
{\rm O}_{x,y}(5)
\ \ \ \ \ \ \ \ \ \ \ \ \ \ \ \ \ \ \ \
\text{and}
\ \ \ \ \ \ \ \ \ \ \ \ \ \ \ \ \ \ \ \
v
\,=\,
q^2
+
0
+
G_{4,0}\,p^4
+
{\rm O}_{p,q}(5).
\] 
In~{\cite{Chen-Merker-2019}}, this is the branch
$\Saux \equiv 0$.

In this branch {\green{\bf 2b3a}}, the linear representation is:
\[
\Big(
G_{4,0}
\Big)
\,=\,
\Big(
\tfrac{1}{a_{1,1}^2}
\Big)\,
\big(
F_{4,0}
\big),
\]
hence branch creation is immediate since we work over $\C$:
\[
\def\arraystretch{1.25}
\begin{array}{rc}
\green{\bf 2b3a}\,\,\,\,
\green{\downarrow}\,\,
& 
F_{4,0} 
\\
\green{\bf 4a} & 
0 
\\
\green{\bf 4b} & 
1 
\end{array}
\]
In~{\cite{Chen-Merker-2019}}
the (Monge) differential invariant corresponding
to $F_{4,0}$ is denoted:
\[
\Paux
\,:=\,
\frac{1}{3}\,
\frac{-\,5\,F_{xxx}^2+3\,F_{xx}\,F_{xxxx}}{F_{xx}^2}.
\]

In Branch {\green{\bf 2b3a4a}}, 
the two differential invariants $0 \equiv \Saux \equiv \Paux$ 
vanish identically\,\,---\,\,{\em cf.} the equations we denoted
$0 \equiv \Kaux_1 \equiv \cdots \equiv \Kaux_t$ in the general
setting. Of course, 
this hypothesis imposes consequences
at higher order jet levels in the equations
$0 = {\sf E}_\smallbullet^{\sf vf}$.

In Branch {\green{\bf 2b3a4b}}, we may normalize
$F_{4,0} := 1 =: G_{4,0}$, and the equation:
\[
1
\,=\,
\tfrac{1}{a_{1,1}^2}
\cdot
1,
\]
imposes the group reduction $a_{1,1} := 1$,
even over $\R$, 
because we decided to work modulo the action of a finite
group. In fact, to neglect finite group actions
is `justified' by our goal of obtaining 
Lie algebras of vector fields parametrized by
absolute invariants (if present).

It happens that beyond~{\green{\bf 2b3a4a}} and
beyond~{\green{\bf 2b3a4b}}, there are no more group
reductions. This is a reason to declare that the 
{\footnotesize\sf Steps~1, 2, 3, 4, 5} algorithm stops,
{\em cf.} Section~{\ref{termination-moduli-spaces-homogeneous-models}} 
{\em infra}.
Nevertheless, 
the equations $0 = {\sf E}_\smallbullet^{\sf vf}$ 
must be still resolved up to some higher jets
in order to completely determine the 
concerned two Lie algebras. A list
of 12 homogeneous models with Lie algebras
appears at the end of the present
Section~{\ref{classification-affinely-homogeneous-S2-C3}}.

As already mentioned, 
in Branch {\green{\bf 2b3b}}, 
$G_{2,1} := 1 =: F_{2,1}$ as
well, and the equation:
\[
1
\,=\,
\tfrac{1}{a_{2,2}}
\cdot
1,
\]
imposes the group reduction $a_{2,2} := 1$. 
Then at order 4, after taking account of 
the equations $0 = {\sf E}_\smallbullet^{\sf vf}$, we have:
\[
u
\,=\,
x^2
+
x^2\,y
+
F_{3,1}\,x^3\,y
+
x^2\,y^2
+
{\rm O}_{x,y}(5)
\ \ \ \ \ \ \
\text{and}
\ \ \ \ \ \ \
v
\,=\,
p^2
+
p^2\,q
+
G_{3,1}\,p^3\,q
+
p^2\,q^2
+
{\rm O}_{p,q}(5),
\]
and the linear representation at order 4 is:
\[
\Big(
G_{3,1}
\Big)
\,=\,
\Big(
\tfrac{1}{a_{1,1}}
\Big)\,
\Big(
F_{3,1}
\Big),
\]
leading to the 2 order 4 branches:
\[
\def\arraystretch{1.25}
\begin{array}{rc}
\green{\bf 2b3b}\,\,\,\,
\green{\downarrow}\,\,
& 
F_{3,1} 
\\
\green{\bf 4a} & 
0 
\\
\green{\bf 4b} & 
1 
\end{array}
\]

To this
{\em punctual} relative invariant $F_{3,1}$ 
is associated the explicit:
\[
\Waux
\,:=\,
\frac{
F_{xx}^2\,F_{xxxy}
-
F_{xx}\,F_{xy}\,F_{xxxx}
+
2\,F_{xy}\,F_{xxx}^2
-
2\,F_{xx}\,F_{xxx}\,F_{xxy}}{
(F_{xx})^2\,\,
\big(
F_{xx}\,F_{xxy}
-
F_{xy}\,F_{xxx}
\big)^{2/3}},
\]
which was in~{\cite{Chen-Merker-2019}} an {\em absolute} 
differential invariant for the {\sl special} affine
group $\Saff (\C^3)$, and which is here is a {\em relative} 
(differential) invariant (with the same explicit expression
in terms of $J_{x,u}^4 F$) for the (full) affine group $\Aff (\C^3)$.

From now on, let us content ourselves with just {\em listing}
the linear representations 
together with the associated orbit-transversals
because in Theorem~{\ref{Thm-redo-S2-C3}}, the 
(truncated) normal forms
will be written anyway.
It remains to show what happens
at the four nodes\big/branches:
\[
\text{\green{\bf 2b3b4a}},
\ \ \ \ \ \ \ \ \ \ \ \ \ \ 
\text{\green{\bf 2c}},
\ \ \ \ \ \ \ \ \ \ \ \ \ \ 
\text{\green{\bf 2c3a}},
\ \ \ \ \ \ \ \ \ \ \ \ \ \ 
\text{\green{\bf 2c3b}}.
\]

In Branch~{\green{\bf 2b3b4a}}, the linear representation is:
\[
\Big(
G_{5,0}
\Big)
\,=\,
\Big(
\tfrac{1}{a_{1,1}^3}
\Big)\,
\Big(
F_{5,0}
\Big),
\]
leading to the orbit-transversals (sub-branches):
\[
\def\arraystretch{1.25}
\begin{array}{rc}
\green{\bf 2b3b4a}\,\,\,\,
\green{\downarrow}\,\,
& 
F_{5,0} 
\\
\green{\bf 5a} & 
0 
\\
\green{\bf 5b} & 
1 
\end{array}
\]

In Branch~{\green{\bf 2c}}, the linear representation is:
\[
\def\arraystretch{1.25}
\left[
\begin{array}{c}
G_{3,0}
\\
G_{0,3}
\end{array}
\right]
\,=\,
\left[
\begin{array}{cc}
\frac{a_{2,1}}{a_{1,1}^2} & 0
\\ 
0 & \frac{a_{1,1}}{a_{2,2}^2}
\end{array}
\right]\,
\left[
\begin{array}{c}
F_{3,0}
\\
F_{0,3}
\end{array}
\right],
\]
leading to:
\[
\def\arraystretch{1.25}
\begin{array}{rcc}
\green{\bf 2c}\,\,\,\,
\green{\downarrow}\,\,
& 
F_{3,0} & F_{0,3}
\\
\green{\bf 3a} & 
0 & 0
\\
\green{\bf 3b} & 
0 & 1
\\
\green{\bf 3c} & 
1 & 1  
\end{array}
\]

In Branch~{\green{\bf 2c3a}}, the linear representation is:
\[
\Big(
G_{2,2}
\Big)
\,=\,
\Big(
\tfrac{1}{a_{1,1}a_{2,2}}
\Big)\,
\Big(
F_{2,2}
\Big),
\]
leading to:
\[
\def\arraystretch{1.25}
\begin{array}{rc}
\green{\bf 2c3a}\,\,\,\,
\green{\downarrow}\,\,
& 
F_{2,2} 
\\
\green{\bf 4a} & 
0 
\\
\green{\bf 4b} & 
1 
\end{array}
\]

In Branch~{\green{\bf 2c3b}}, the linear representation is:
\[
\def\arraystretch{1.25}
\left[
\begin{array}{c}
G_{2,2}
\\
G_{1,3}
\\
G_{0,4}
\end{array}
\right]
\,=\,
\left[
\begin{array}{ccc}
\frac{1}{a_{2,2}^3} & 0 & 0
\\
0 & \frac{1}{a_{2,2}^2} & 0
\\
0 & 0 & \frac{1}{a_{2,2}}
\end{array}
\right]\,
\left[
\begin{array}{c}
F_{2,2}
\\
F_{1,3}
\\
F_{0,4}
\end{array}
\right],
\]
leading to:
\[
\def\arraystretch{1.25}
\begin{array}{rccc}
\green{\bf 2c3b}\,\,\,\,
\green{\downarrow}\,\,
& 
F_{2,2} & F_{1,3} & F_{0,4}
\\
\green{\bf 4a} & 
0 & 0 & 0
\\
\green{\bf 4b} & 
0 & 0 & 1
\\
\green{\bf 4c} & 
0 & 1  & F_{0,4}
\\
\green{\bf 4d} & 
1 & F_{1,3}  & F_{0,4}
\end{array}
\]
In the last two lines, $F_{0,4}$, and $F_{1,3}$, $F_{0,4}$ 
are absolute (punctual) invariants.

To present our classification principles on a known geometric
structure,
let us state a theorem which was 
essentially established by
Doubrov-Komrakov-Rabinovich~{\cite{Doubrov-Komrakov-Rabinovich-1996}}
and by Eastwood-Ezhov~{\cite{Eastwood-Ezhov-1999}}.
The theorem below is less precise
in the sense that we abandon
the search for closed forms, 
and we abandon to explore finite group equivalences.

Nevertheless, this theorem has the virtue of specifying
explicitly the way how homogeneous models
depend on the branching tree of punctual invariants,
in a purely algebraic manner, deeply related
with the Lie-Fels-Olver recurrence relations
between differential invariants.

Moreover, 
this somewhat restricted approach appears to be 
more suitable for generalizations 
to higher-dimensional geometric structures,
as {\em e.g.} the
classification of affinely homogeneous surfaces $S^2 \subset \R^4$
achieved in 
Sections~{\ref{S2-R4}}
$\to$
{\ref{2g-models}}, 
for which finding closed forms and\big/or
neutralizing finite equivalences
might be problematic,
especially 
in circumstances where much more branches exist,
and much more complicated simply transitive models exist.

\begin{Theorem}
\label{Thm-redo-S2-C3}
Up to the action of a finite subgroup of $\Aff(\C^3)$,
there are 12 (families of) 
affinely homogeneous
model surfaces $S^2 \subset \C^3$,
as represented by the diagram on 
p.~{\pageref{diag-arbre-S2-R3}}.

The list of (truncated) normal forms together with
their associated
transitive Lie algebras of vector fields is as follows.
\end{Theorem}

%%%%%%%%%%%%%%%%%%%%%%%%%%%%%%%%%%%%%%%%%%%%%%%%%%%%%%%%%%%%%%%%%%%%%%
%%%%%%%%%%%%%%%%%%%%%%%%%%%%%%%%%%%%%%%%%%%%%%%%%%%%%%%%%%%%%%%%%%%%%%
%%%%%%%%%%%%%%%%%%%%%%%%%%%%%%%%%%%%%%%%%%%%%%%%%%%%%%%%%%%%%%%%%%%%%%

\[
\text{\bf Model 2a} \ \ \ \ \
\Big\{
u = 0.
\]
%%%%%%%%%%%%%%%%%%%%%%%%%%%%%%%%%%%%%%%%%%%%%%%%%%%%%%%%%%%%%%%%%%%%%%
\[
\def\arraystretch{1.25}
\begin{array}{lllll}
e_1 := \partial_x, &
e_2 := \partial_y, &
e_3 := x\partial_x, &
e_4 := y\partial_x, &
e_5 := x\partial_y,
\\
e_6 := y\partial_y, &
e_7 := u\partial_x, &
e_8 := u\partial_y, &
e_9 := u\partial_u,
\end{array}
\]
%%%%%%%%%%%%%%%%%%%%%%%%%%%%%%%%%%%%%%%%%%%%%%%%%%%%%%%%%%%%%%%%%%%%%%
\[
\footnotesize
\def\arraystretch{1.25}
\begin{array}{c|ccccccccc}
{} & e_1 & e_2 & e_3 & e_4 & e_5 & e_6 & e_7 & e_8 & e_9 
\\
\hline
e_1 & 
0 & 0 & e_1 & 0 & e_2 & 0 & 0 &0  & 0
\\
e_2 &
0 & 0 & 0 & e_1 & 0 & e_2 & 0 & 0 & 0
\\
e_3 &
-e_1 & 0 &0 & -e_4 & e_5 & 0 & -e_7 & 0 & 0
\\
e_4 &
0 & -e_1 & e_4 & 0 & -e_3+e_6 & -e_4 & 0 & -e_7 & 0
\\
e_5 &
-e_2 & 0 & -e_5 & e_3-e_6 & 0 & e_5 & -e_8 & 0 & 0
\\
e_6 &
0 & -e_2 & 0 & e_4 & -e_5 & 0 & 0 & -e_8 & 0
\\
e_7 &
0 & 0 & e_7 & 0 & e_8 & 0 & 0& 0 & -e_7
\\
e_8 &
0 & 0 & 0 & e_7 & 0 & e_8 & 0 & 0 & -e_8
\\
e_9 &
0 & 0 & 0 & 0 & 0 & 0 & e_7 & e_8 & 0
\end{array}
\]

%%%%%%%%%%%%%%%%%%%%%%%%%%%%%%%%%%%%%%%%%%%%%%%%%%%%%%%%%%%%%%%%%%%%%%
%%%%%%%%%%%%%%%%%%%%%%%%%%%%%%%%%%%%%%%%%%%%%%%%%%%%%%%%%%%%%%%%%%%%%%
%%%%%%%%%%%%%%%%%%%%%%%%%%%%%%%%%%%%%%%%%%%%%%%%%%%%%%%%%%%%%%%%%%%%%%

\[
\text{\bf Model 2b3a4a}
\ \ \ \ \
\Big\{
\aligned
u
&
\,=\,
x^2,
\endaligned
\]
%%%%%%%%%%%%%%%%%%%%%%%%%%%%%%%%%%%%%%%%%%%%%%%%%%%%%%%%%%%%%%%%%%%%%%
\[
\def\arraystretch{1.25}
\begin{array}{llll}
e_1
\,:=\,
\partial_x+2x\partial_u
, &
e_2
\,:=\,
\partial_y
, &
& 
\\
e_3
\,:=\,
x\partial_x+2u\partial_u
, &
e_4
\,:=\,
x\partial_y
, &
e_5
\,:=\,
y\partial_y
, &
e_6
\,:=\,
u\partial_y,
\end{array}
\]

%%%%%%%%%%%%%%%%%%%%%%%%%%%%%%%%%%%%%%%%%%%%%%%%%%%%%%%%%%%%%%%%%%%%%%
\[
\footnotesize
\def\arraystretch{1.25}
\begin{array}{c|cccccc}
{} & e_1 & e_2 & e_3 & e_4 & e_5 & e_6 
\\
\hline
e_1 & 
0 & 0 & e_1 & e_2 & 0 & 2e_4
\\
e_2 &
0 & 0 & 0 & 0 & e_2 & 0
\\
e_3 &
-e_1 & 0 & 0 & e_4 & 0 & 2e_6
\\
e_4 &
-e_2 & 0 & -e_4 & 0 & e_4 & 0
\\
e_5 &
0 & -e_2 & 0 & -e_4 & 0 & -e_6
\\
e_6 &
-2e_4 & 0 & -2e_6 & 0 & e_6 & 0
\end{array}
\]

%%%%%%%%%%%%%%%%%%%%%%%%%%%%%%%%%%%%%%%%%%%%%%%%%%%%%%%%%%%%%%%%%%%%%%
%%%%%%%%%%%%%%%%%%%%%%%%%%%%%%%%%%%%%%%%%%%%%%%%%%%%%%%%%%%%%%%%%%%%%%
%%%%%%%%%%%%%%%%%%%%%%%%%%%%%%%%%%%%%%%%%%%%%%%%%%%%%%%%%%%%%%%%%%%%%%

\[
\text{\bf Model 2b3a4b}
\ \ \ \ \
\Big\{
\aligned
u
&
\,=\,
x^2+x^4+F_{5,0}x^5+\tfrac{5}{4}F_{5,0}^2x^6
+
\cdots,
\endaligned
\]
for any value for $F_{5,0}$,
%%%%%%%%%%%%%%%%%%%%%%%%%%%%%%%%%%%%%%%%%%%%%%%%%%%%%%%%%%%%%%%%%%%%%%
\[
\def\arraystretch{1.25}
\aligned
e_1 & := -\Big(\tfrac{5}{2}F_{5,0}x+2u-1\Big)\partial_x-\Big(5F_{5,0}u-2x\Big)\partial_u,\\
e_2 & := \partial_y,
\,\,\,\,\,\,\,\,\,\,
e_3 := x\partial_y,
\,\,\,\,\,\,\,\,\,\,
e_4 := y\partial_y,
\,\,\,\,\,\,\,\,\,\,
e_5 := u\partial_y,
\endaligned
\]

%%%%%%%%%%%%%%%%%%%%%%%%%%%%%%%%%%%%%%%%%%%%%%%%%%%%%%%%%%%%%%%%%%%%%%
\[
\footnotesize
\def\arraystretch{1.25}
\begin{array}{c|ccccc}
{} & e_1 & e_2 & e_3 & e_4 & e_5 
\\
\hline
e_1 & 
0 & 0 & e_2-\tfrac{5}{2}F_{5,0}e_3-2e_5 & 0 & 2e_3-5F_{5,0}e_5
\\
e_2 &
0 & 0 & 0 & e_2 & 0
\\
e_3 &
-e_2+\tfrac{5}{2}F_{5,0}e_3+2e_5 & 0 & 0 & e_3 & 0
\\
e_4 &
0 & -e_2 & -e_3 & 0 & -e_5
\\
e_5 &
5F_{5,0}e_5-2e_3 & 0 &0 & e_5 & 0
\end{array}
\]

%%%%%%%%%%%%%%%%%%%%%%%%%%%%%%%%%%%%%%%%%%%%%%%%%%%%%%%%%%%%%%%%%%%%%%
%%%%%%%%%%%%%%%%%%%%%%%%%%%%%%%%%%%%%%%%%%%%%%%%%%%%%%%%%%%%%%%%%%%%%%
%%%%%%%%%%%%%%%%%%%%%%%%%%%%%%%%%%%%%%%%%%%%%%%%%%%%%%%%%%%%%%%%%%%%%%

\[
\text{\bf Model 2b3b4a5a}
\ \ \ \ \
\Big\{
\aligned
u
&
\,=\,
x^2
+
x^2y
+
x^2y^2
+
x^2y^3
+
x^2y^4
+
x^2y^5
+
x^2y^6
+
\cdots,
\endaligned
\]
%%%%%%%%%%%%%%%%%%%%%%%%%%%%%%%%%%%%%%%%%%%%%%%%%%%%%%%%%%%%%%%%%%%%%%
\[
\def\arraystretch{1.25}
\aligned
e_1 & := (1-y)\partial_x+2x\partial_u,\\
e_2 & := (1-y)\partial_y+u\partial_u,\\
e_3 & := 2u\partial_u+x\partial_x,\\
e_4 & := -\tfrac{1}{2}u\partial_x+x\partial_y,
\endaligned
\]

%%%%%%%%%%%%%%%%%%%%%%%%%%%%%%%%%%%%%%%%%%%%%%%%%%%%%%%%%%%%%%%%%%%%%%
\[
\footnotesize
\def\arraystretch{1.25}
\begin{array}{c|cccc}
{} & e_1 & e_2 & e_3 & e_4
\\
\hline
e_1 & 
0 & e_1 & e_1 & e_2
\\
e_2 &
-e_1 & 0 & 0 & e_4
\\
e_3 &
-e_1 & 0 & 0 & e_4
\\
e_4 &
-e_2 & -e_4 & -e_4 & 0
\end{array}
\]

%%%%%%%%%%%%%%%%%%%%%%%%%%%%%%%%%%%%%%%%%%%%%%%%%%%%%%%%%%%%%%%%%%%%%%
%%%%%%%%%%%%%%%%%%%%%%%%%%%%%%%%%%%%%%%%%%%%%%%%%%%%%%%%%%%%%%%%%%%%%%
%%%%%%%%%%%%%%%%%%%%%%%%%%%%%%%%%%%%%%%%%%%%%%%%%%%%%%%%%%%%%%%%%%%%%%

\[
\text{\bf Model 2b3b4a5b}
\ \ \ \ \
\left\{
\aligned
u
&
\,=\,
x^2+x^2y+x^2y^2+x^2y^3+x^5+x^2y^4+4x^5y+x^2y^5+10x^5y^2
\\
&
\,\,\,\,\,\,
+
F_{7,0}x^7
+
x^2y^6+20x^5y^3+6F_{7,0}x^7y
+\tfrac{25}{8}x^8+\cdots,
\endaligned\right.
\]
for any value for $F_{7,0}$,
%%%%%%%%%%%%%%%%%%%%%%%%%%%%%%%%%%%%%%%%%%%%%%%%%%%%%%%%%%%%%%%%%%%%%%
\[
\def\arraystretch{1.25}
\aligned
e_1 & := (7F_{7,0}u-y+1)\partial_x-(14F_{7,0}x+5u)\partial_y+2x\partial_u,\\
e_2 & := -x\partial_x-(y-1)\partial_y-u\partial_u,
\endaligned
\]
%%%%%%%%%%%%%%%%%%%%%%%%%%%%%%%%%%%%%%%%%%%%%%%%%%%%%%%%%%%%%%%%%%%%%%
\[
\footnotesize
\def\arraystretch{1.25}
\begin{array}{c|cc}
{} & e_1 & e_2 
\\
\hline
e_1 & 
0 & 0
\\
e_2 &
0 & 0
\end{array}
\]

%%%%%%%%%%%%%%%%%%%%%%%%%%%%%%%%%%%%%%%%%%%%%%%%%%%%%%%%%%%%%%%%%%%%%%
%%%%%%%%%%%%%%%%%%%%%%%%%%%%%%%%%%%%%%%%%%%%%%%%%%%%%%%%%%%%%%%%%%%%%%
%%%%%%%%%%%%%%%%%%%%%%%%%%%%%%%%%%%%%%%%%%%%%%%%%%%%%%%%%%%%%%%%%%%%%%

\[
\text{\bf Model 2b3b4b}
\ \ \ \ \
\Big\{
\aligned
u
&
\,=\,
x^2+x^2y+x^3y+x^2y^2+x^5+x^2y^3+\cdots,
\endaligned
\]
%%%%%%%%%%%%%%%%%%%%%%%%%%%%%%%%%%%%%%%%%%%%%%%%%%%%%%%%%%%%%%%%%%%%%%
\[
\def\arraystretch{1.25}
\aligned
e_1 &:= -(-3x-1+y+5u)\partial_x-(15u-10x+3y)\partial_y+(6u+2x)\partial_u,\\
e_2 & := (-2x+3u\tfrac{1}{2})\partial_x+(4u-4x-y+1)\partial_y-3u\partial_u,
\endaligned
\]
%%%%%%%%%%%%%%%%%%%%%%%%%%%%%%%%%%%%%%%%%%%%%%%%%%%%%%%%%%%%%%%%%%%%%%
\[
\footnotesize
\def\arraystretch{1.25}
\begin{array}{c|cc}
{} & e_1 & e_2 
\\
\hline
e_1 & 
0 & -e_1-e_2
\\
e_2 &
e_1+e_2 & 0
\end{array}
\]

%%%%%%%%%%%%%%%%%%%%%%%%%%%%%%%%%%%%%%%%%%%%%%%%%%%%%%%%%%%%%%%%%%%%%%
%%%%%%%%%%%%%%%%%%%%%%%%%%%%%%%%%%%%%%%%%%%%%%%%%%%%%%%%%%%%%%%%%%%%%%
%%%%%%%%%%%%%%%%%%%%%%%%%%%%%%%%%%%%%%%%%%%%%%%%%%%%%%%%%%%%%%%%%%%%%%

\[
\text{\bf Model 2c3a4a}
\ \ \ \ \
\Big\{
\aligned
u
&
\,=\,
xy
.
\endaligned
\]
%%%%%%%%%%%%%%%%%%%%%%%%%%%%%%%%%%%%%%%%%%%%%%%%%%%%%%%%%%%%%%%%%%%%%%
\[
\def\arraystretch{1.25}
\aligned
e_1 & := \partial_x+y\partial_u,
\,\,\,\,\,\,\,\,\,\,\,\,
e_2  := \partial_y+x\partial_u,\\
e_3 & := x\partial_x+u\partial_u,
\,\,\,\,\,\,\,\,\,\,\,\,
e_4 := y\partial_y+u\partial_u,
\endaligned
\]

%%%%%%%%%%%%%%%%%%%%%%%%%%%%%%%%%%%%%%%%%%%%%%%%%%%%%%%%%%%%%%%%%%%%%%
\[
\footnotesize
\def\arraystretch{1.25}
\begin{array}{c|cccc}
{} & e_1 & e_2 & e_3 & e_4
\\
\hline
e_1 & 
0 & 0 & e_1 & 0
\\
e_2 &
0 & 0 & 0 & e_2
\\
e_3 &
-e_1 & 0 & 0 & 0
\\
e_4 &
0 & -e_2 &0 & 0
\end{array}
\]

%%%%%%%%%%%%%%%%%%%%%%%%%%%%%%%%%%%%%%%%%%%%%%%%%%%%%%%%%%%%%%%%%%%%%%
%%%%%%%%%%%%%%%%%%%%%%%%%%%%%%%%%%%%%%%%%%%%%%%%%%%%%%%%%%%%%%%%%%%%%%
%%%%%%%%%%%%%%%%%%%%%%%%%%%%%%%%%%%%%%%%%%%%%%%%%%%%%%%%%%%%%%%%%%%%%%

\[
\text{\bf Model 2c3a4b}
\ \ \ \ \
\Big\{
\aligned
u
&
\,=\,
xy
+
x^2y^2
+
2x^3y^3
+
\cdots,
\endaligned
\]
%%%%%%%%%%%%%%%%%%%%%%%%%%%%%%%%%%%%%%%%%%%%%%%%%%%%%%%%%%%%%%%%%%%%%%
\[
\def\arraystretch{1.25}
\begin{array}{ll}
e_1
&
\,:=\,
(1-2u)\partial_x+y\partial_u
,
\\
e_2
&
\,:=\,
(1-2u)\partial_y+x\partial_u
, 
\\
e_3
&
\,:=\,
x\partial_x-y\partial_y,
\end{array}
\]

%%%%%%%%%%%%%%%%%%%%%%%%%%%%%%%%%%%%%%%%%%%%%%%%%%%%%%%%%%%%%%%%%%%%%%
\[
\footnotesize
\def\arraystretch{1.25}
\begin{array}{c|ccc}
{} & e_1 & e_2 & e_3
\\
\hline
e_1 & 
0 & 2e_3 & e_1
\\
e_2 &
-2e_3 & 0 & -e_2
\\
e_3 &
-e_1 & e_2 & 0
\end{array}
\]

%%%%%%%%%%%%%%%%%%%%%%%%%%%%%%%%%%%%%%%%%%%%%%%%%%%%%%%%%%%%%%%%%%%%%%
%%%%%%%%%%%%%%%%%%%%%%%%%%%%%%%%%%%%%%%%%%%%%%%%%%%%%%%%%%%%%%%%%%%%%%
%%%%%%%%%%%%%%%%%%%%%%%%%%%%%%%%%%%%%%%%%%%%%%%%%%%%%%%%%%%%%%%%%%%%%%

\[
\text{\bf Model 2c3b4a}
\ \ \ \ \
\Big\{
\aligned
u
&
\,=\,
xy+y^3
+
\cdots,
\endaligned
\]
%%%%%%%%%%%%%%%%%%%%%%%%%%%%%%%%%%%%%%%%%%%%%%%%%%%%%%%%%%%%%%%%%%%%%%
\[
\def\arraystretch{1.25}
\begin{array}{ll}
e_1
&
\,:=\,
\partial_x+y\partial_u,
\\
e_2
&
\,:=\,
-3y\partial_x+\partial_y+x\partial_u,
\\
e_3
&
\,:=\,
2x\partial_x+y\partial_y+3u\partial_u,
\end{array}
\]

%%%%%%%%%%%%%%%%%%%%%%%%%%%%%%%%%%%%%%%%%%%%%%%%%%%%%%%%%%%%%%%%%%%%%%
\[
\footnotesize
\def\arraystretch{1.25}
\begin{array}{c|ccc}
{} & e_1 & e_2 & e_3
\\
\hline
e_1 &
0 & 0 & 2e_1 
\\
e_2 &
0 & 0 & e_2
\\
e_3 &
-2e_1 & -e_2 & 0
\end{array}
\]

%%%%%%%%%%%%%%%%%%%%%%%%%%%%%%%%%%%%%%%%%%%%%%%%%%%%%%%%%%%%%%%%%%%%%%
%%%%%%%%%%%%%%%%%%%%%%%%%%%%%%%%%%%%%%%%%%%%%%%%%%%%%%%%%%%%%%%%%%%%%%
%%%%%%%%%%%%%%%%%%%%%%%%%%%%%%%%%%%%%%%%%%%%%%%%%%%%%%%%%%%%%%%%%%%%%%

\[
\text{\bf Model 2c3b4b}
\ \ \ \ \
\Big\{
\aligned
u
&
\,=\,
xy
+
y^3
+
y^4+F_{5,0}y^5
+
\cdots,
\endaligned
\]
for any value for $F_{5,0}$,
%%%%%%%%%%%%%%%%%%%%%%%%%%%%%%%%%%%%%%%%%%%%%%%%%%%%%%%%%%%%%%%%%%%%%%
\[
\def\arraystretch{1.25}
\begin{array}{ll}
e_1
\,:=\,
\partial_x+y\partial_u,
\\
e_2
\,:=\,
-(10F_{0,5}x-12x+3y)\partial_x-(5F_{0,5}y-4y-1)\partial_y
-(15F_{0,5}u-16u-x)\partial_u,
\end{array}
\]
%%%%%%%%%%%%%%%%%%%%%%%%%%%%%%%%%%%%%%%%%%%%%%%%%%%%%%%%%%%%%%%%%%%%%%
\[
\footnotesize
\def\arraystretch{1.25}
\begin{array}{c|cc}
{} & e_1 & e_2 
\\
\hline
e_1 & 
0 & -(10F_{5,0}-12)e_1
\\
e_2 &
(10F_{5,0}-12)e_1 & 0
\end{array}
\]

%%%%%%%%%%%%%%%%%%%%%%%%%%%%%%%%%%%%%%%%%%%%%%%%%%%%%%%%%%%%%%%%%%%%%%
%%%%%%%%%%%%%%%%%%%%%%%%%%%%%%%%%%%%%%%%%%%%%%%%%%%%%%%%%%%%%%%%%%%%%%
%%%%%%%%%%%%%%%%%%%%%%%%%%%%%%%%%%%%%%%%%%%%%%%%%%%%%%%%%%%%%%%%%%%%%%

\[
\text{\bf Model 2c3b4c}
\ \ \ \ \
\Big\{
\aligned
u
&
\,=\,
xy+y^3+F_{0,4}y^4+xy^3+F_{0,4}xy^4+(\tfrac{6}{5}+\tfrac{6}{5}F_{0,4}^2)y^5
+
\cdots,
\endaligned
\]
for any value for $F_{0,4}$,
%%%%%%%%%%%%%%%%%%%%%%%%%%%%%%%%%%%%%%%%%%%%%%%%%%%%%%%%%%%%%%%%%%%%%%
\[
\def\arraystretch{1.25}
\aligned
e_1 & := (x+1)\partial_x+(u+y)\partial_u,\\
e_2 & := -(3u+3y)\partial_x-(2F_{0,4}y-1)\partial_y-(2F_{0,4}u-x)\partial_u,
\endaligned
\]
%%%%%%%%%%%%%%%%%%%%%%%%%%%%%%%%%%%%%%%%%%%%%%%%%%%%%%%%%%%%%%%%%%%%%%
\[
\footnotesize
\def\arraystretch{1.25}
\begin{array}{c|cc}
{} & e_1 & e_2 
\\
\hline
e_1 & 
0 & 0
\\
e_2 &
0 & 0
\end{array}
\]

%%%%%%%%%%%%%%%%%%%%%%%%%%%%%%%%%%%%%%%%%%%%%%%%%%%%%%%%%%%%%%%%%%%%%%
%%%%%%%%%%%%%%%%%%%%%%%%%%%%%%%%%%%%%%%%%%%%%%%%%%%%%%%%%%%%%%%%%%%%%%
%%%%%%%%%%%%%%%%%%%%%%%%%%%%%%%%%%%%%%%%%%%%%%%%%%%%%%%%%%%%%%%%%%%%%%

\[
\text{\bf Model 2c3c}
\ \ \ \ \
\left\{
\aligned
u
&
\,=\,
xy+x^3+y^3+F_{4,0}x^4+F_{3,1}x^3y
+
(\tfrac{2}{9}F_{1,3}F_{0,4}-\tfrac{1}{9}F_{3,1}F_{1,3}+\tfrac{8}{9}F_{4,0}F_{0,4})x^2y^2
+
\\
&
\,\,\,\,\,\,
+
F_{1,3}xy^3+F_{0,4}y^4
+
\cdots,
\endaligned\right.
\]
%%%%%%%%%%%%%%%%%%%%%%%%%%%%%%%%%%%%%%%%%%%%%%%%%%%%%%%%%%%%%%%%%%%%%%
\[
\def\arraystretch{1.25}
\aligned
e_1 &:= 
-
(-1+\tfrac{1}{3}xF_{1,3}+\tfrac{8}{3}xF_{4,0}-9u+\tfrac{4}{9}uF_{1,3}F_{0,4}-\tfrac{2}{9}uF_{3,1}F_{1,3}+\tfrac{16}{9}uF_{4,0}F_{0,4})\partial_x
+
\\
&
\,\,\,\,\,\,
-
(3x+\tfrac{4}{3}yF_{4,0}+\tfrac{2}{3}yF_{1,3}+3uF_{3,1})\partial_y-(F_{1,3}u+4F_{4,0}u-y)\partial_u,\\
e_2 &:= 
-
(\tfrac{4}{3}xF_{0,4}+\tfrac{2}{3}xF_{3,1}+3y+3uF_{1,3})\partial_x
-
(-1+\tfrac{1}{3}yF_{3,1}+\tfrac{8}{3}yF_{0,4}-9u+\tfrac{4}{9}uF_{1,3}F_{0,4}
+
\\
&
\,\,\,\,\,\,
-\tfrac{2}{9}uF_{3,1}F_{1,3}+\tfrac{16}{9}uF_{4,0}F_{0,4})\partial_y-(4F_{0,4}u+F_{3,1}u-x)\partial_u,
\endaligned
\]
with invariants $F_{0,4},F_{1,3},F_{3,1},F_{4,0}$ 
satisfying the 3 algebraic equations:
\[
\aligned
0 
& 
\,=\,
F_{0,4}F_{1,3}-F_{3,1}F_{4,0},
\\
0 
&
\,=\,
-32F_{0,4}^2F_{4,0}-16F_{0,4}F_{3,1}F_{4,0}
+F_{1,3}F_{3,1}^2+2F_{3,1}^2F_{4,0}+18F_{1,3}^2
+36F_{1,3}F_{4,0}+162F_{0,4}+81F_{3,1},
\\
0 
&
\,=\,
-32F_{0,4}F_{4,0}^2+F_{1,3}^2F_{3,1}+2F_{1,3}F_{3,1}F_{4,0}-16F_{3,1}F_{4,0}^2+36F_{0,4}F_{3,1}+18F_{3,1}^2+81F_{1,3}+162F_{4,0},
\endaligned
\]
%%%%%%%%%%%%%%%%%%%%%%%%%%%%%%%%%%%%%%%%%%%%%%%%%%%%%%%%%%%%%%%%%%%%%%
\[
\footnotesize
\def\arraystretch{1.25}
\begin{array}{c|cc}
{} & e_1 & e_2 
\\
\hline
e_1 &
0 & (-\tfrac{4}{3}F_{0,4}-\tfrac{2}{3}F_{3,1})e_1
+
(\tfrac{2}{3}F_{1,3}+\tfrac{4}{3}F_{4,0})e_2 
\\
e_2 &
-(-\tfrac{4}{3}F_{0,4}
-
\tfrac{2}{3}F_{3,1})e_1
-
(\tfrac{2}{3}F_{1,3}+\tfrac{4}{3}F_{4,0})e_2 & 0
\end{array}
\]

%%%%%%%%%%%%%%%%%%%%%%%%%%%%%%%%%%%%%%%%%%%%%%%%%%%%%%%%%%%%%%%%%%%%%%
\SectionHead{Termination: Moduli Spaces of Homogeneous Models}
{termination-moduli-spaces-homogeneous-models}
%%%%%%%%%%%%%%%%%%%%%%%%%%%%%%%%%%%%%%%%%%%%%%%%%%%%%%%%%%%%%%%%%%%%%%

Now, when, why, and how the 
{\footnotesize\sf Steps~1-2-3-4-5} `algorithm'
{\em terminates}?
What does its `{\em termination}' produce?
Before answering these questions, 
let us present some aspects
of the current state of the art.

First of all, as for Cartan's method of equivalence which 
is sometimes termed to be an `algorithm',
most of the times, as soon as the number $n$ of
independent variables $x_\smallbullet$ is $\geqslant 2$ or is
$\geqslant 3$, any `equivalence algorithm' 
can be `blocked' by computational complexity,
even with the help of powerful machines.

Indeed, the exploration of the branching tree 
of a given kind of homogeneous geometries requires
in some circumstances
to continue the computations until reaching
{\sl simply transitive} models,
{\em i.e.} those with: 
\[
\dim\, 
\mathfrak{sym}\,M
\,=\,
\dim\,M,
\]
and then in this case, the `algorithm' 
necessarily terminates.
The other homogeneous models $M$, 
those for which $\dim\, \mathfrak{sym} (M) > 
\dim\, M$, are
termed {\sl multiply transitive}.

As a matter of fact, it is (well) known in the literature that, for a
number of famous geometric structures, either simply transitive models
were never found yet, or were found by indirect methods, without
discovering the complete branching tree created by invariants together
with all the linear representations in the nodes.  Let us give 5
examples.

\medskip\noindent$\bullet$\,
For $(2, 3, 5)$ distributions $D^2$ in a five-manifold $M^5$, 
Cartan classified multiply transitive models,
{\em see}~{\cite{The-2022}}
(and the reference therein) for 
a recent synthesis based
on Cartan (parabolic) geometries,
and~{\em see} also~{\cite{Doubrov-Govorov-2013}}
for a complete classification,
including simply transitive models,
which is based on Lie algebraic techniques.
Beyond Cartan quartic types, it seems that
no complete picture exists for the branching
tree of order $\geqslant 5$ (punctual) invariants.

\medskip\noindent$\bullet$\,
For completely integrable second order {\sc pde} systems
in 2 dependent complex variables and 1 independent 
complex variable, 
the multiply transitive models
have been neatly classified 
in~{\cite{Doubrov-Medvedev-The-2019}},
but the complete branching tree of invariants is also missing,
and simply transitive models have not been determined yet.

\medskip\noindent$\bullet$\,
For CR-homogeneous Levi nondegenerate
hypersurfaces $M^5 \subset \C^3$,
the multiply transitive models
have been neatly classified 
in~{\cite{Doubrov-Medvedev-The-2020, Loboda-2020}},
the complete branching tree of invariants is also missing,
while the simply transitive models
have been determined by abstract Lie algebraic method,
{\em cf.}~{\cite{Loboda-2020, Doubrov-Merker-The-2020}}.

\medskip\noindent$\bullet$\,
For 4\textsuperscript{th} order {\sc ode}s
under point transformations, existing classifications
are not complete,
while classifications of 
homogeneous models 
3\textsuperscript{th} order {\sc ode}s
under fiber-preserving, point, contact, transformations
have been achieved by
Michal Godi{\'n}ski and 
Pawe{\l} Nurowski,
{\em cf.}~{\cite{Godlinski-2008, Godlinski-Nurowski-2009}}.

\medskip\noindent$\bullet$\,
For affinely homogeneous hypersurfaces $H^3 \subset \R^4$,
Eastwood-Ezhov in~{\cite{Eastwood-Ezhov-2001-1}}
do not reach simply transitive models.

\medskip

Now, let us come
back to the {\footnotesize\sf Steps~1-2-3-4-5} `algorithm'
in the general setting. Because: 
\[
n
\,\leqslant\,
\dim\,G 
\,<\,
\infty,
\]
it is clear
that the dimensions of the isotropy subgroups $G_\stab^\kappa 
\subset G$ at orders $\kappa = 0, 1, 2, 3, \dots$,
can decrease (strictly)
only a finite number of times,
in all branches. So in each one of the branches constructed
by induction, after a while, no more isotropy group reduction
can occur. This is when and why
the {\footnotesize\sf Steps~1-2-3-4-5} `algorithm' {\em terminates}.

And in fact, all boxed terminal leaves in the branching
trees shown in this memoir 
indicate {\em termination by
end-of-isotropy-reduction}.

But from the computational point of view, 
how termination does occur, concretely?
Namely, what really happens `at the end' 
of the `ping-pong' play
between the equations 
$0 = {\sf E}_\smallbullet^{\sf vf}$ and
$0 = {\sf E}_\smallbullet^{\sf nf}$?

First of all, 
after that {\footnotesize\sf Steps~1}
and~{\footnotesize\sf 2} have been passed, 
as soon as there is a nontrivial
linear representation in {\footnotesize\sf Step~3},
necessarily, 
there must be at least one further subbranch
which is accompanied with a nontrivial
group reduction\,\,---\,\,except for the mostly degenerate 
linear group-orbit: the origin in the vector space $\R^{\ell_\kappa}$,
{\em cf.} for instance 
Branch~{\green{\bf 2b3a4a}} in
Section~{\ref{classification-affinely-homogeneous-S2-C3}}.

Consequently, termination holds if and only if no (nontrivial)
linear representation occurs at {\footnotesize\sf Step~3},
at all higher jet orders $\kappa$, wich requires
to continue the `ping-pong' between 
$0 = {\sf E}_\smallbullet^{\sf vf}$ (firstly) as
{\footnotesize\sf Step~1} and
$0 = {\sf E}_\smallbullet^{\sf nf}$ 
(secondly) as {\footnotesize\sf Step~2}.

We did not attempt to prove or just state 
stabilization or pseudo-stabilization theorems
as in~{\cite{Olver-1995, Olver-2007b}}\,\,---\,\,which, 
we believe, can be done\,\,---, because there is here a simple 
alternative {\em and direct} way of realizing that the
process rigorously terminates, {\em see} 
Section~{\ref{termination-Jacobi-identities-zero-remainders}}
{\em infra}.

Concretely, the process stops if, after having
resolved the equations $0 = {\sf E}_\smallbullet^{\sf vf}$,
all equations $0 = {\sf E}_\smallbullet^{\sf nf}$ 
only show constant power series coefficients
or absolute invariants $G_\smallbullet = F_\smallbullet$, 
this, at every higher jet order.

Of course, it can happen that termination takes place
with isotropy dimension being stably constant and $> 0$,
whichever high is the jet order.

Thus, termination
holds when the equations 
$0 = {\sf E}_\smallbullet^{\sf nf}$ no more bring any
normalization of the $G_\smallbullet$ and $F_\smallbullet$
coefficients. However, the equations 
$0 = {\sf E}_\smallbullet^{\sf vf}$ still bring a lot of information!

\begin{Observation}
At all higher jet orders, 
when some absolute 
invariants $\Iaux_\smallbullet$ 
which come from preceding jet orders
are still present in computations,
the equations $0 = {\sf E}_\smallbullet^{\sf vf}$
do bring more and more algebraic equations 
in terms of $\Iaux_\smallbullet$
which coherently define a certain 
{\sl algebraic moduli space of homogeneous models}\,\,---\,\,unless 
some algebraic contradiction
occurs which indicates that no homogeneous model can exist
in the considered terminal leaf.
\end{Observation}

{\small (Contradictory terminal leaves are
indicated plainly
with the $\emptyset$ symbol, or even sometimes,
plainly erased.)}

An example of such an
{\sl algebraic moduli space of homogeneous models}
was already shown at the end of
Section~{\ref{classification-affinely-homogeneous-S2-C3}}, 
with the Branch~{\green{\bf 2c3c}}
for surfaces $S^2 \subset \C^3$. 
Certainly, the obtained algebraic equations are deeply related
with the Lie-Fels-Olver {\sl recurrence relations}
between differential invariants, {\em cf.}~{\cite{Chen-Merker-2020}}.

\begin{Observation}
To each boxed terminal leaf, there corresponds 
a family of homogeneous models parametrized by a certain
algebraic variety.
\end{Observation}

Quite often, a terminal leaf 
of a branching tree
is of the form IkP0, 
with $k$-dimensional isotropy Ik,
where P0
means that zero Parameter is present, 
so that the
concerned algebraic variety is just 1 point (or 2, 3, 4
points, never more in this memoir). 

As is known, 
every (complicated)
algebraic variety may always be decomposed
into a finite number of  
simpler disjoint smooth pieces, {\em e.g.},
by a process called {\sl stratifification}.

However, {\em conceptionally}, group reduction 
in the spirit of Lie and Cartan is 
of {\em different nature}, compared with
further explorations by stratifying
algebraic moduli spaces of homogeneous models.
In~{\cite{Eastwood-Ezhov-1999}}, 
both group reduction\,\,---\,\,in 
fact not explicitly 
mentioned there\,\,---\,\,and moduli space
stratification,
seem to be treated on equal footing,
{\em cf.} the flow diagrams on pages~67--69 there.
%%%Even, a certain stratification is presented with 
%%%positive-dimensional overlaps 
%%%in Theorem~2 of~{\cite{Eastwood-Ezhov-1999}}.

\CITATION{
Ma derni\`ere remarque g\'en\'erale concerne un aspect
de la Math\'ematique moderne en quelque sorte compl\'ementaire
de ses tendances unificatrices, \`a savoir
sa capacit\'e \`a dissocier ce qui \'etait 
ind\^ument confondu.
Jean {\sc Dieudonn\'e},~{\cite[p.~13]{Dieudonne-1964}}.}

As Dieudonn\'e writes, one must indeed `{\footnotesize\sf\em 
dissociate what
was unduly confused}'.

\begin{Observation}
In this memoir, we decided not to stratify the 
algebraic moduli spaces of homogeneous models 
that we obtained, after termination
of group reduction, at any terminal leaf.
\end{Observation}

Such a task could be endeavoured in a future publication.
Of course, stratifying an
algebraic moduli space of homogeneous
models would bring further sub-branches 
(of a different nature), 
devloping and branching after the boxed terminal
leaves.

In all of the 3 articles~{\cite{Eastwood-Ezhov-1999, 
Eastwood-Ezhov-2001-1, Eastwood-Ezhov-2001-2}},
only a single branch among all the branches
studied there does lead to a non-trivial algebraic
variety, namely what we call here Branch~{\green{\bf 2c3c}},
and what is called there `{\sl Nonvanishing Pick Invariant}'.

\CITATION{
In fact, if the system $E = F = G = H = 0$ is simply passed to the
‘solve’ routine of the computer algebra system MAPLE (Version V
Release 3), then the program returns the correct solutions as a set of
approximately 20 cases, in effect constructing its
own flow diagram)!~{\cite[pp.~32--33]{Eastwood-Ezhov-1999}}}

These 4 equations $E = F = G = H = 0$ are 
precisely equivalent to the
3 equations appearing in Branch~{\green{\bf 2c3c}} after
Theorem~{\ref{Thm-redo-S2-C3}}, and it is
indeed already difficult
to stratify their zero-locus.

All other branches treated in~{\cite{Eastwood-Ezhov-1999, 
Eastwood-Ezhov-2001-1, Eastwood-Ezhov-2001-2}}
directly lead either to individual
models with all coefficients $F_\smallbullet$ 
being numerical (hence uniquely prescribed),
or to the existence of a single real or complex
parameter (absolute invariant)
$\Iaux \in \R^1$ or $\Iaux \in \C^1$,
with no algebraic equation involved.

By contrast, in this memoir, several terminal leaves,
especially the {\em simply transitive} ones,
led us to certain quite complicated algebraic moduli spaces
of homogeneous models, far beyond
what was handled in~{\cite{Eastwood-Ezhov-1999, 
Eastwood-Ezhov-2001-1, Eastwood-Ezhov-2001-2}}.

Even for just one terminal leaf like {\em e.g.}:
\[
{\green{\bf 2f3a}},
\ \ \ \ \ \ \ 
\text{or}
\ \ \ \ \ \ \ 
{\green{\bf 2f3g}},
\ \ \ \ \ \ \
\text{or}
\ \ \ \ \ \ \
{\green{\bf 2g3a}},
\]
to set up a stratification could be a formidable task!
We believe that a similar algebraic complexity 
lies behind the simply transitive
affinely homogeneous hypersurfaces $H^3 \subset \R^4$,
not attained in~{\cite{Eastwood-Ezhov-2001-1}}.

{\small {\em Passim}, we would like to mention that
at the starting order 2, the 
other article~{\cite{Eastwood-Ezhov-2001-2}}
does not show a clear linear matrix group representation
(which appears implictly in Appendix~II).
Another point is that, even
at the starting order 2, 
group-transversals are not neatly stipulated, 
as the authors themselve recognize it:}

\CITATION{
Note that Table (2.5) is not a classification.
There are overlaps and repetitions under the action (2.2).
{\cite[p.~724]{Eastwood-Ezhov-2001-2}}}

In addition, no stratification in smooth neatly parametrized
pieces would be `canonical' in any 
sense\,\,---\,\,similarly as the choice of 
a group-transversal
is never `canonical'.

In conclusion, in this memoir, 
our classification approach decides to stop 
(to terminate)
once algebraic moduli spaces of homogeneous models
have been reached.

And now, we know the reason why, in the existing literature,
some classifications using the approach
with (differential) invariants are missing,
especially concerning the (difficult) 
simply transitive homogeneous models.

It is because the concerned algebraic varieties
which parametrize the sought 
(simply transitive) homogeneous 
models happen to be very complicated.

%%%%%%%%%%%%%%%%%%%%%%%%%%%%%%%%%%%%%%%%%%%%%%%%%%%%%%%%%%%%%%%%%%%%%%
\SectionHead{Termination: Jacobi Identities and Zero Remainders}
{termination-Jacobi-identities-zero-remainders}
%%%%%%%%%%%%%%%%%%%%%%%%%%%%%%%%%%%%%%%%%%%%%%%%%%%%%%%%%%%%%%%%%%%%%%

Lastly, we explain why it is not necessary to explore
the equations $0 = {\sf E}_\smallbullet^{\sf vf}$
at all jet orders $\kappa = 0, 1, 2, \dots$ up to infinity.

At a terminal leaf, no more group reduction takes place.
Consequently, from all the previous resolutions
of isotropy parameters 
$A$, $B$, $C$, $D$
present in the infinitesimal
affine transformation $L$, 
there are $n$ transitivity vector fields 
$e_1, \dots, e_n$ corresponding to the transitivity
parameters $T_1, \dots, T_n$, with:
\[
T_0M
\,=\,
\Span\,
\Big(
e_1\big\vert_0,\,
\dots,\,
e_n\big\vert_0
\Big),
\]
and there
are $\nu \geqslant 0$ isotropy vector fields
$f_1, \dots, f_\nu$, which vanish
at the origin $(x, u) = (0, 0)$:
\[
f_1\big\vert_0
\,=\,
\cdots
\,=\,
f_\nu\big\vert_0
\,=\,
0,
\]
such that the complete collection of vector fields:
\[
\big\{
e_1,\dots,e_n,\,
f_1,\dots,f_\nu
\big\},
\]
should constitute a {\em transitive
Lie algebra of vector fields},
provided the terminal collection of absolute invariants
$\Iaux_\smallbullet$ present in $e_1, \dots, e_n$ and
in $f_1, \dots, f_\nu$
satisfies
all the algebraic equations
equations $0 = {\sf E}_\smallbullet^{\sf vf} \big( 
\Iaux_\smallbullet \big)$,
at all orders $\kappa = 0, 1, 2, \dots$.

Let us admit the notational coincidence:
\[
e_{n+1}
\,=\,
f_1,\,\,\,
\dots,\,\,\,
e_{n+\nu}
\,=\,
f_\nu.
\]
Thus, each vector field $e_m$ with $1 \leqslant m \leqslant n+\nu$
is a linear combination of the $\partial_{x_i}$ and
of the $\partial_{u_j}$, with coefficients which
are of degree 1 (are affine) with respect to 
the $x_{i'}$ and to the $u_{j'}$, 
and which may also depend on a certain collection 
$\Iaux_\smallbullet$
of absolute invariants.

By knowing the explicit expressions of these $e_1, \dots, e_{n+\nu}$,
and by looking at a well chosen collection of 
$n + \nu$ monomials uniquely associated with
the isotropy generators $e_{n+1}, \dots, e_{n+\nu}$, 
we can easily determine (using computer programs)
the structure constants,
so that:
\[
\big[e_{m_1},\,e_{m_2}\big]
\,=\,
\sum_{s=1}^{n+\nu}\,
C_{m_1,m_2,s}\,
e_s.
\]

Then two conditions must be satisfied for $\big\{ e_m \big\}_{1 
\leqslant m \leqslant n+\nu}$ to really constitute a Lie algebra
of vector fields.

\medskip\noindent{\footnotesize\sf (Jac)}\,
Jacobi identities for the found structure constants.

\medskip\noindent{\footnotesize\sf (Null)}\,
Identically zero vector fields:
\[
0
\,\equiv\,
\big[e_{m_1},\,e_{m_2}\big]
-
\sum_{s=1}^{n+\nu}\,
C_{m_1,m_2,s}\,
e_s,
\]
while only a subset of these identities
was used to determine the structure constants.

\medskip

The interest of these two
(new) constraints is that now,
they are finite in number, because
the Lie brackets $\big[ e_{m_1}, e_{m_2} \big]$
are still of degree $\leqslant 1$ in
the $x_{i'}$ and the $u_{j'}$.

\begin{Observation}
In all terminal leaves of homogeneous models
attained in this memoir,
the algebraic equations $0 = {\sf E}_\smallbullet^{\sf vf}$
at a high enough order (but never $> 9$)
always defined the {\em same} algebraic variety
as the algebraic variety coming from 
{\footnotesize\sf (Jac)}
and
{\footnotesize\sf (Null)}.
\end{Observation}

Lastly, because $e_1, \dots, e_{n+\nu}$ are {\em analytic} vector
fields, the standard Frobenius theorem guarantees that there exists a
local analytic submanifold $M^n \subset \R^{n+c}$ (possibly
parametrized by some invariants $\Iaux_\smallbullet$ satisfying the
concerned algebraic equations), which is the sought locally affinely
homogeneous model.

%%%%%%%%%%%%%%%%%%%%%%%%%%%%%%%%%%%%%%%%%%%%%%%%%%%%%%%%%%%%%%%%%%%%%%
\SectionHead{Lie Algebras of Vector Fields Versus Closed Forms}
{lie-algebras-vector-fields-versus-closed-forms}
%%%%%%%%%%%%%%%%%%%%%%%%%%%%%%%%%%%%%%%%%%%%%%%%%%%%%%%%%%%%%%%%%%%%%%

In the literature, most of the times,
classifications of 
affinely or projectively homogenous
small-dimensional submanifolds in $\R^{n+c}$ 
(or in $\C^{n+c}$) attain {\em closed forms},
that is, equations $u = F(x)$ with
$F$ being expressed as a polynomial,
or\big/and in terms of elementary transcendental
functions: exponentials,
logarithms, trigonometric functions,
{\em see}~{\cite{Abdalla-Dillen-Vrancken-1997,
Doubrov-Govorov-2013,
Doubrov-Komrakov-1998,
Doubrov-Komrakov-Rabinovich-1996,
Doubrov-Medvedev-The-2019,
Doubrov-Medvedev-The-2020,
Doubrov-Merker-The-2020,
Eastwood-Ezhov-1999,
Eastwood-Ezhov-2001-1,
Eastwood-Ezhov-2001-2,
Godlinski-Nurowski-2009,
Loboda-2020,
Merker-Nurowski-2020-a,
Merker-Nurowski-2023,
Mozhei-2000,
The-2022,
Wermann-2001}}.

However, often, Lie algebras 
$\mathfrak{sym} (M)$ 
of infinitesimal symmetries are not shown,
simultaneously with the (nice) functions $F$.
And it is then a non-immediate task to determine
$\mathfrak{sym} (M)$ from a given closed graphed form $\big\{ 
u = F(x) \big\}$, especially when
some continuous parameters $\alpha$, $\beta$, 
\dots, are present
in $F = F_{\alpha, \beta}$.
Indeed, for various values of the parameters
in the closed form,
most probably,
the (graphed) 
manifold:
\[
M_{\alpha, \beta, \dots}
\,=\,
\big\{ 
u
= 
F_{\alpha,\beta,\dots} 
\big\},
\]
{\em crosses the branches}
of any invariant branching tree,
so that the dimensions of 
$\mathfrak{sym} \big( M_{\alpha, \beta, \dots} \big)$
`jump' in some way as the parameters 
$\alpha, \beta, \dots$ do vary.
Mathematical life is complicated!

Not only the Lie algebras 
$\mathfrak{sym} \big( M_{\alpha, \beta, \dots} \big)$
are not always shown in the literature,
but also, the branching trees created by invariants
are almost never constructed until reaching all terminal leaves,
except in some computationally simple cases, 
{\em e.g.} in the curve case $n = 1$.

Certainly, such invariant branching
trees are often at least partly known,
concerning relatively small order
differential invariants, 
for instance: 
Cartan's quartic 
for $(2, 3, 5)$ distributions;
or Chern-Moser's order 4 tensor
for CR hypersurfaces $M^5 \subset \C^3$. 

For the constant Hessian rank 1 affinely homogeneous
hypersurfaces $H^2 \subset \R^3$ and
$H^3 \subset \R^4$ that we treat 
in Sections~{\ref{surfaces-S2-R3}}
and~{\ref{threefolds-H3-R4}},
closed forms are known, 
{\em cf.}~{\cite{Mozhei-2000, 
Merker-Nurowski-2023}},
while we did not search
for a closed form 
representation of the single (up to sign)
model $H^4 \subset \R^5$ which we found 
in Theorem~{\ref{Thm-H4-R5}}. 

And for higher-dimensional constant Hessian rank 1
hypersurfaces $H^n \subset \R^{n+1}$, 
a quite unexpected fact was
established in~{\cite{Merker-2022}}, namely 
that in any dimension $n \geqslant 5$, 
there are no nonproduct homogeneous models at all!
{\em Passim}, let us raise a

\begin{Question}
{\sl Is it true, also, that 
in all high enough dimensions $n 
\geqslant \NN_2 \gg 1$,
there are, similarly, no nonproduct
constant Hessian rank 2
hypersurfaces $H^n \subset \R^{n+1}$ (or $\C^{n+1}$)
which are locally affinely homogeneous\,{\bf\em ?}}
\end{Question}

A similar question may be formulated for any 
fixed constant Hessian rank
$1 \leqslant r \leqslant n - 1$.
Also, the question may be considered with the
{\sl projective} group $\Proj (\R^{n+c})$ instead
of $\Aff (\R^{n+c})$.

Back to closed forms (that we will not seek
in this memoir), 
experts to whom we asked whether there exist theoretical
explanations\,\,---\,\,that could be read off from a given
Lie algebra of vector fields\,\,---\,\,why, 
when, how, closed forms (may)
exist, answered us that they ignore
what could be such reason(s), and that they obtained closed forms
with the help of Maple 
{\sc pde} integration programs.

It is therefore legitimate to raise a

\begin{Problem}
\label{Pbm-criteria-closed-forms}
{\sl Find criteria, if not necessary and sufficient conditions,
on given Lie algebras of vector fields 
that are symmetries of homogeneous models,
in order that, after a change of coordinates belonging
to the initial group $G$, the graphing function
$F(x)$ is either polynomial, or is expressed
in terms of usual transcendental functions:
exponentials, logarithms, trigonometric functions.}
\end{Problem}

Of course, the theorem of Frobenius guarantees
the existence of an analytic graph $M^n \subset \R^{n+c}$
which is simply the {\em orbit of the origin}
under the action of the found transitive  
Lie algebra of vector fields.
Such Lie algebras may be truncated,
especially when dealing with an infinite-dimensional
Lie group $G$ acting on $\R^{n+c}$
(or $\C^{n+c}$),
and again, the same problem appears to be meaningful.

In sum, there are 5 reasons why we did not seek closed forms
(for the moment).

\medskip\noindent$\bullet$\,
No general theory seems to exist around
Problem~{\ref{Pbm-criteria-closed-forms}},
and probably, there might exist certain
special homogeneous Lie algebras
of vector fields which would {\em not} be
elementarily integrable. 

\medskip\noindent$\bullet$\,
Lie's original principle of 
classification~{\cite{Engel-Lie-Merker-2015, Engel-Lie-1893}}, 
with which we agree, was to determine and to present
{\em Lie algebras of vector fields}, only.

\medskip\noindent$\bullet$\,
Punctual invariants of homogeneous models
are strongly related to algebras
of differential invariants, a research field
that we learned from Peter Olver's monographs and articles,
and in this field, branching by invariants is 
a natural process.

\medskip\noindent$\bullet$\,
Successive group reductions leading to linear
representations in all nodes seem to be
universal, although they were not discovered in 
the existing normal forms articles we know.

\medskip\noindent$\bullet$\,
Branching trees of invariants lie at the heart
matter, hence must be exhibited, even when quite ramified.

\medskip

In a near future,
we hope to contribute to 
Problem~{\ref{Pbm-criteria-closed-forms}}.
Certainly, the purely algebraic tools used in the present
memoir will appear to be largely insufficient.

%%%%%%%%%%%%%%%%%%%%%%%%%%%%%%%%%%%%%%%%%%%%%%%%%%%%%%%%%%%%%%%%%%%%%%
\SectionHead{Classification by Invariants}
{classification-invariants}
%%%%%%%%%%%%%%%%%%%%%%%%%%%%%%%%%%%%%%%%%%%%%%%%%%%%%%%%%%%%%%%%%%%%%%

Now, how a given concrete individual
$M^\ast = \big\{ u = F^\ast (x) \big\}$ 
can `{\em pass through}' a
finalized invariant branching tree\,?
How can it be submitted to our 
{\footnotesize\sf Steps~1, 2, 3, 4, 5} algorithm\,?
What are its invariants? How to compute them?

Sometimes, such questions are raised even
for {\em families} $M_{\alpha, \beta, \dots} 
= \big\{ u = F_{\alpha, \beta, \dots} \big\}$
of `individual' submanifolds\,\,---\,\,in fact
depending on continuous parameters.
Especially, when the question is to  
{\em compare} two (or more) different classifications.
Fortunately, 
we are not aware of the existence of any
affine classification of surfaces $S^2 \subset \R^4$
as done in Sections~{\ref{S2-R4}}
$\to$
{\ref{2g-models}},
hence in this case, the comparison question is void.

But the two different classifications
of affinely homogeneous $S^2 \subset \R^3$
(or $\C^3$) obtained 
in~{\cite{Doubrov-Komrakov-Rabinovich-1996}
and in~{\cite{Eastwood-Ezhov-1999}} 
required to produce normal forms for
certain
{\em parametrized} families of surfaces, 
{\em e.g.}:
\[
u
\,=\,
(1+x)^\alpha\,
(1+y)^\beta
\eqno
{\scriptstyle{(\alpha,\,\beta\,\in\,\R)}},
\]
viewed
at the origin, and this was nontrivial, computationally,
{\em see}~{\cite[Sec.~6]{Eastwood-Ezhov-1999}}.
Similar computations, though more complicated,
were done
in~{\cite[Sec.~6]{Doubrov-Merker-The-2020}},
using power series in 5 (instead of 2) variables, in order
to compare affine in $\R^3$
and Cauchy-Riemann in $\C^5$
homogeneous models
classifications. 

In any case, 
to determine how a given parametrized family
$M_{\alpha, \beta, \dots} 
= \big\{ u = F_{\alpha, \beta, \dots} \big\}$ of submanifolds
would {\sl pass through} a certain
finalized branching tree of invariants
forces to admit that the given family shall be cut into pieces.
And in any case
also, the construction of a complete branching tree of invariants
concerns the full family 
of $M = \big\{ u = F(x) \big\}$,
with all power series
coefficients $F_\smallbullet$ as parameters!
So the construction of a finalized branching tree
of invariants is in any case more general and more complicated
than for subfamilies!

Consequently, we will only treat the question of how a given
{\em individual} $M^\ast = \big\{ u = F^\ast (x) \big\}$,
with given {\em numerical} coefficients $F_\smallbullet^\ast$,
does pass through a branching tree,
which we assume to be already given, finalized, complete,
so as to determine the (punctual) invariants of $M^\ast$.

Let us illustrate this (simple) process in the case of surfaces
$S^2 \subset \C^3$, since the associated tree 
shown in Section~{\ref{classification-affinely-homogeneous-S2-C3}} 
is not too large.
For surfaces $S^2 \subset \R^4$, 
and for other geometric structures as well, 
the principle
of reasoning is identical.

The method proceeds order by order.
It consists in normalizing $\big\{ u = F^\ast (x,y) \big\}$
by means of the fundamental equation~{\eqref{intro-fund-eq-S2-C3}},
all along the (general) {\sl branched route} expressed
by the diagram on p.~{\ref{diag-arbre-S2-R3}}.
Notation is as in 
Section~{\ref{classification-affinely-homogeneous-S2-C3}}.

Since order 1 terms may always be normalized to 0
by means of a transvection
without changing the
values of the other power series coefficients,
we may start from:
\[
u
\,=\,
F_{2,0}^\ast\,
x^2
+
F_{1,1}^\ast\,
x\,y
+
F_{0,2}^\ast\,
y^2
+
{\rm O}_{x,y}(2).
\]
Here, the $F_{i,j}^\ast \in \C$ are thought of as being
{\em numerical} quantities, not formal variables.

According to the diagram on p.~{\ref{diag-arbre-S2-R3}},
there are 3 possible order 2 normalizations:
\[
v
\,=\,
\Big\{\,
0
\ \ \ \ \ 
\text{or}
\ \ \ \ \ 
p^2
\ \ \ \ \ 
\text{or}
\ \ \ \ \ 
p\,q
\,\Big\}
+
{\rm O}_{p,q}(3).
\]
We must determine to which of these 3 orbit-transversals
does the given numerical value
$\big( F_{2,0}^\ast, F_{1,1}^\ast, F_{0,2}^\ast \big)$ belong.

But we already know that the orbits are
distinguished by means of the Hessian. Hence
it suffices to insert these numerical
values into the Hessian matrix.

However, for general linear representations 
in arbitrary dimensions
and for groups with any number of parameters, 
the complexity of orbit-transversals can be high,
unlimited.
In other words, {\footnotesize\sf Step~4}
of the `algorithm' stated in 
Section~{\ref{general-method-affine-context}}
represents in itself an ample mathematical problem,
very much studied in the literature,
{\em cf.} 
Peter Olver's monograph~{\cite{Olver-1999}}
and the references therein.

\begin{Problem}
{\bf [Classical]}
{\sl For a given linear representation of a Lie group,
to determine invariant (polynomial or rational) 
functions which separate the orbits.}
\end{Problem}

Even for the affine group $\Aff(\R^{n+c})$ in general
dimension $n$ and codimension $c$, 
and for the linear representations which would appear
within the nodes of the associated branching diagrams,
this (sub)problem might be difficult in itself.

To determine explicit invariants of group actions,
either differential ones after prolongation to jet spaces,
or polynomial ones by restricting to linear representations,
requires some elimination computations. 
In fact, as was explained in~{\cite{Olver-2018,
Chen-Merker-2019, Chen-Merker-2020}}, 
explicit expressions of differential
invariants can be obtained by normalizing
power series at the origin,
keeping memory of each normalizing
transformation, 
before {\em performing the complete composition}
of the normalizing transformations. But such computations
are usually very demanding, and in this memoir,
our intention is to avoid them as much as possible.

\begin{SubProblem}
\label{Pbm-invariants-orbits-nodes}
{\sl For every linear representation appearing
in {\footnotesize\sf Step~4}
of the `algorithm' stated in 
Section~{\ref{general-method-affine-context}},
to determine invariant (polynomial or rational) 
functions which separate the partitions by orbit-transversals.}
\end{SubProblem}

A similar problem exists while classifying abstract 
Lie algebras~{\cite{Snobl-Winternitz-2014}},
and the so-called {\sl Casimir invariants}
are a classical tool to separate the (generic) 
orbits (only) under the action of a Lie group
on itself by inner automorphisms.

If such a task would be achieved for each node of 
a branching tree, then, by a straightforward
insertion of the values of the coefficients 
$F_{\ast,\ind}^\kappa$ into the computed invariants
for a node of order $\kappa$,
one would immediately recognize to which 
orbit $F_{\ast,\ind}^\kappa$ belongs.

However, in this memoir, it is not really necessary to 
solve Problem~{\ref{Pbm-invariants-orbits-nodes}},
which might lead to computational difficulties.
Indeed, since we take for granted that the complete 
branching tree of invariants has been constructed before
by a classification theorem,
it suffices to `{\sl test}' whether a given value
$\big( F_{2,0}^\ast, F_{1,1}^\ast, F_{0,2}^\ast \big)$
belongs to 1 of the 3 existing branches. How?

By just writing down the 
{\em fundamental equation}~{\eqref{intro-fund-eq-S2-C3}},
and more precisely, 
by writing the 3 equations~{\eqref{eq-x2-xy-y2}},
between 
$\big( F_{2,0}^\ast, F_{1,1}^\ast, F_{0,2}^\ast \big)$
and $(0,0,0)$ or $(1,0,0)$ or $(0,1,0)$, 
before `{\em testing}' whether they can
be realized with appropriate choices of
the appearing group parameters:
$a_{1,1}$, $a_{2,1}$, $a_{2,2}$.
By construction, and according to the classification
theorem,
only one normal form among
$(0,0,0)$ or $(1,0,0)$ or $(0,1,0)$ can fit.

To fix ideas, let us assume that the outcome
is that 
$\big( F_{2,0}^\ast, F_{1,1}^\ast, F_{0,2}^\ast \big)$
is equivalent to $(0,1,0)$.
Then solve the 3 equations~{\eqref{eq-x2-xy-y2}}
for some specific numerical values 
$a_{1,1}^*$, $a_{2,1}^*$, $a_{2,2}^*$.
Only a single choice is enough (while the general
choice is conjugated to the stability group
$G_\stab^2$). 

Once $a_{1,1}^*$, $a_{2,1}^*$, $a_{2,2}^*$ are chosen,
execute the associated affine (in fact linear)
transformation. How? 

With computer help, 
by truncating, say at order 10, the initial
numerical power series $F^\ast (x)$,
by writing the fundamental equation~{\eqref{intro-fund-eq-S2-C3}},
and by solving for the target power series
coefficients $G_\smallbullet^{3,\ast}, \dots,
G_\smallbullet^{10,\ast}$ at orders 
$3, \dots, 10$\,\,---\,\,after $G_\smallbullet^{2,\ast} 
= (0, 1, 0)$. These tasks are purely numerical,
hence instantly done by a computer.

At higher jet orders, the linear representations will be
diagonal, hence to determine to which branch
a numerical jet $J_{\ast, \ind}^\kappa$ does belong
happens to be a trivial task, without
the needs to address 
Sub-Problem~{\ref{Pbm-invariants-orbits-nodes}}.

After renaming $G_\smallbullet^\ast =: F_\smallbullet^\ast$, get:
\[
u
\,=\,
x\,y
+
F_{3,0}^\ast\,x^3
+
F_{2,1}^\ast\,x^2y
+
F_{1,2}^\ast\,xy^2
+
F_{0,3}^\ast\,y^3
+
{\rm O}_{x,y}(4).
\]
According to the branching diagram, 
in (node) Branch~{\green{\bf 2c}}
the (numerical) linear representation is:
\[
\def\arraystretch{1.25}
\left[
\begin{array}{c}
G_{3,0}^\ast
\\
G_{0,3}^\ast
\end{array}
\right]
\,=\,
\left[
\begin{array}{cc}
\frac{a_{2,1}}{a_{1,1}^2} & 0
\\ 
0 & \frac{a_{1,1}}{a_{2,2}^2}
\end{array}
\right]\,
\left[
\begin{array}{c}
F_{3,0}^\ast
\\
F_{0,3}^\ast
\end{array}
\right].
\]
Implicitly, this means that the
other two power series coefficients
$F_{2,1}^\ast$ and $F_{1,2}^\ast$ can 
{\em always} be normalized
to 0, by means of some free group parameters,
here $b_1$ and $b_2$.

Hence do normalize these $G_{2,1}^\ast := 0$ and
$G_{1,2}^\ast := 0$ in the fundamental
equation~{\eqref{intro-fund-eq-S2-C3}}
with some specific numerical values $b_1^\ast$ and 
$b_2^\ast$.
Again, with the same truncation order 10
and from the fundamental
equation~{\eqref{intro-fund-eq-S2-C3}},
solve (easily) on a computer for the target power series
coefficients $G_{3,0}^\ast$, $G_{0,3}^\ast$, and
$G_\smallbullet^{4,\ast}, \dots,
G_\smallbullet^{10,\ast}$.
Again, rename $G_\smallbullet^\ast =: F_\smallbullet^\ast$.

Since the linear representation is diagonal,
and since the normal forms at order 3 are:
\[
\def\arraystretch{1.25}
\begin{array}{rcc}
\green{\bf 2c}\,\,\,\,
\green{\downarrow}\,\,
& 
F_{3,0}^\ast & F_{0,3}^\ast
\\
\green{\bf 3a} & 
0 & 0
\\
\green{\bf 3b} & 
0 & 1
\\
\green{\bf 3c} & 
1 & 1  
\end{array}
\]
it is immediate to see to which of the 3 next branches does
$\big( F_{3,0}^\ast, F_{0,3}^\ast \big)$ belong.

To fix ideas, suppose that $F_{3,0}^\ast = 0$
while $F_{0,3}^\ast \neq 0$.
Again, look at the fundamental
equation~{\eqref{intro-fund-eq-S2-C3}} at order 3,
which is nothing but:
\[
\big(\,
1
\,=:\,
\big)\,\,
G_{0,3}^\ast
\,=\,
\frac{a_{1,1}}{a_{2,2}^2}\,
F_{0,3}^\ast,
\]
use one group parameter, for instance $a_{1,1}$, to 
normalize $G_{0,3}^\ast := 1$,
and again, with the same truncation order 10
and from the fundamental
equation~{\eqref{intro-fund-eq-S2-C3}},
solve for
$G_\smallbullet^{4,\ast}, \dots,
G_\smallbullet^{10,\ast}$.

After renaming $G_\smallbullet^\ast =: F_\smallbullet^\ast$, get:
\[
u
\,=\,
x\,y
+
y^3
+
F_{4,0}^\ast\,x^4
+
F_{3,1}^\ast\,x^3y
+
F_{2,2}^\ast\,x^2y^2
+
F_{1,3}^\ast\,xy^3
+
F_{0,4}^\ast\,y^4
+
{\rm O}_{x,y}(5).
\]

According to the branching diagram, 
in (node) Branch~{\green{\bf 2c3b}}
the (numerical) linear representation is:
\[
\def\arraystretch{1.25}
\left[
\begin{array}{c}
G_{2,2}^\ast
\\
G_{1,3}^\ast
\\
G_{0,4}^\ast
\end{array}
\right]
\,=\,
\left[
\begin{array}{ccc}
\frac{1}{a_{2,2}^3} & 0 & 0
\\
0 & \frac{1}{a_{2,2}^2} & 0
\\
0 & 0 & \frac{1}{a_{2,2}}
\end{array}
\right]\,
\left[
\begin{array}{c}
F_{2,2}^\ast
\\
F_{1,3}^\ast
\\
F_{0,4}^\ast
\end{array}
\right].
\]
Implicitly, this means that the two unwritten coefficients
$F_{4,0}^\ast$, $F_{3,1}^\ast$ have been normalized 
to 0\,\,---\,\,or are dependent jets! 

This is the case here, because the (punctual) invariant
$F_{3,0}^\ast = 0$ is assumed in
Branch~{\green{\bf 2c3b}} to vanish,
and because by homogeneity,
the corresponding differential invariant 
vanishes identically, hence by diffentiating
this differential invariant,
it follows that $F_{4,0} = 0$ and $F_{3,1} = 0$.

As we explained in Section~{\ref{general-method-affine-context}},
we reach the information about dependent jets
in an alternative way,
namely 
by using the {\sl transitivity equations}.

But here, because the surface has numerical coefficients,
nothing must be computed, and two cases may occur.

\medskip\noindent$\square$\,
Either one of the two numerical coefficients
$F_{4,0}^\ast \neq 0$ or $F_{3,1}^\ast \neq 0$ is nonzero,
which means that the given numerical 
surface $\big\{ u = F^\ast (x) \big\}$
is {\em not} affinely homogeneous\,\,---\,\,end of the story.
Equivalently, in the transitivity equations,
a nontrivial linear relation exists between
the transitivity parameters $T_1$ and $T_2$,
which is a contradiction in the search for
homogeneous models.

\medskip\noindent$\square$\,
Or the two numerical coefficients
$F_{4,0}^\ast = 0 = F_{3,1}^\ast$ are zero, which
leaves a chance for having a homogeneous model.

\medskip

Suppose that the second case holds.
Since the linear representation is diagonal, and
since the normal forms at order 3 are:
\[
\def\arraystretch{1.25}
\begin{array}{rccc}
\green{\bf 2c3b}\,\,\,\,
\green{\downarrow}\,\,
& 
F_{2,2}^\ast & F_{1,3}^\ast & F_{0,4}^\ast
\\
\green{\bf 4a} & 
0 & 0 & 0
\\
\green{\bf 4b} & 
0 & 0 & 1
\\
\green{\bf 4c} & 
0 & 1  & F_{0,4}^\ast
\\
\green{\bf 4d} & 
1 & F_{1,3}^\ast & F_{0,4}^\ast,
\end{array}
\]
it is immediate to see to which of the 4 next branches
does $\big( F_{2,2}^\ast, F_{1,3}^\ast, F_{0,4}^\ast \big)$
belong.

To fix ideas, suppose that $F_{2,2}^\ast = 0$ while
$F_{1,3}^\ast \neq 0$, which 
means that $\big\{ u = F^\ast (x) \big\}$
lands in Branch~{\green{\bf 2c3b4c}}.
Again, look at the fundamental equation~{\eqref{intro-fund-eq-S2-C3}},
at order 4, realize that the coefficient
of $xy^3$ in it gives the equation:
\[
\big(\,
1
\,=:\,
\big)\,\,\,
G_{1,3}^\ast
\,=\,
\frac{1}{a_{2,2}^2}\,
F_{1,3}^\ast,
\]
use the group parameter $a_{2,2} \in \C$ to normalize
$G_{1,3}^\ast := 1$, 
and again, with the same truncation order 10
and from the fundamental
equation~{\eqref{intro-fund-eq-S2-C3}},
solve for
$G_\smallbullet^{5,\ast}, \dots,
G_\smallbullet^{10,\ast}$.

After renaming $G_\smallbullet^\ast =: F_\smallbullet^\ast$, get:
\[
\aligned
u
&
\,=\,
x\,y
+
y^3
+
xy^3
+
F_{0,4}^\ast\,y^4
\\
&
\ \ \ \ \ 
+
F_{5,0}^\ast\,x^5
+
F_{4,1}^\ast\,x^4y
+
F_{3,2}^\ast\,x^3y^2
+
F_{2,3}^\ast\,x^2y^5
+
F_{1,4}^\ast\,xy^4
+
F_{0,5}^\ast\,y^5
+
{\rm O}_{x,y}(6).
\endaligned
\]
According to the branching diagram on
p.~{\ref{diag-arbre-S2-R3}},
the equivalence algorithm terminates,
in the sense that no more group reduction can happen.
In fact here, the existing homogeneous models
are simply transitive with 1-dimensional moduli
space algebraic variety, 
as is indicated by I0P1.

Here, two cases may occur.

\medskip\noindent$\square$\,
Either the obtained (normalized) numerical surface
$\big\{ u = F^\ast (x) \big\}$
admits 2 linearly independent tangential vector fields
as tested up to order 10 by the
{\sf v}ector {\sf f}ields equations 
$0 = {\sf E}_\smallbullet^{\sf vf}$.

\medskip\noindent$\square$\,
Or there is $0 = {\sf E}_\smallbullet^{\sf vf}$
at least one nontrivial linear relation
between the transitivity parameters $T_1$ and $T_2$,
whichs is a contradiction.

\medskip

Of course, if one starts from the beginning with 
a numerical surface $\big\{ u = F^\ast (x) \big\}$
which is known to be affinely homogeneous,
then whatever order higher than 10 is chosen to truncate,
non contradiction shall ever appear.

%%%%%%%%%%%%%%%%%%%%%%%%%%%%%%%%%%%%%%%%%%%%%%%%%%%%%%%%%%%%%%%%%%%%%%
\SectionHead{Classification of Affinely Homogeneous Surfaces 
$S^2 \subset \R^4$}
{S2-R4}
%%%%%%%%%%%%%%%%%%%%%%%%%%%%%%%%%%%%%%%%%%%%%%%%%%%%%%%%%%%%%%%%%%%%%%

The goal of this Part~II,
consisting of 
Sections~{\ref{S2-R4}}
$\to$
{\ref{2g-models}}, 
is to determine all affinely homogeneous
local hypersurfaces $S^2 \subset \R^4$, a problem
which seems to be still open.

In $\R^4 \ni (x,y, u,v)$, 
local analytic surfaces $S^2$ can be graphed, 
after an affine transformation, as:
\[
\aligned
u
&
\,=\,
F\big(x,y\big)
\,=\,
F_{2,0}\,x^2
+
F_{1,1}\,x\,y
+
F_{0,2}\,y^2
+
{\rm O}_{x,y}(3),
\\
v
&
\,=\,
G(x,y)
\,=\,
G_{2,0}\,x^2
+
G_{1,1}\,x\,y
+
G_{0,2}\,y^2
+
{\rm O}_{x,y}(3),
\endaligned
\]
with $F$, $G$ real-analytic at the origin. We show in 
Proposition~{\ref{Prp-parallel-quadrics-F2-G2}} that the property that
the two quadratic forms $F_2$ and $G_2$ are {\sl parallel}
(colinear) is affinely invariant.
Then in Section~{\ref{parallel-quadrics-order-2}}
and~{\ref{non-parallel-quadrics-order-2}}, 
we show that 7 inequivalent normalizations exist at order 2:
\[
\begin{array}{ccc}
{} & F_2 & G_2
\\
\green{\bf 2a} & 0 & 0
\\
\green{\bf 2b} & x^2 & 0
\\
\green{\bf 2c} & x\,y & 0
\\
\green{\bf 2d} & x^2+y^2 & 0
\\
\green{\bf 2e} & x\,y & x^2
\\
\green{\bf 2f} & x\,y & x^2+y^2
\\
\green{\bf 2g} & x\,y & x^2-y^2
\end{array}
\]

The power series method of equivalence 
can be applied to other geometric structures,
including equivalences under infinite-dimensional group actions,
{\em cf.}~{\cite{
Foo-Merker-Nurowski-Ta-2021}}.

%%%%%%%%%%%%%%%%%%%%%%%%%%%%%%%%%%%%%%%%%%%%%%%%%%%%%%%%%%%%%%%%%%%%%%
\SectionHead{Surfaces $S^2 \subset \R^4$ Under $\GL(4,\R)$}
{surfaces-S2-R4}
%%%%%%%%%%%%%%%%%%%%%%%%%%%%%%%%%%%%%%%%%%%%%%%%%%%%%%%%%%%%%%%%%%%%%%

After translation, an affine transformation of $\R^4$
fixes the origin. Consider therefore a linear map
$(x,y,u,v) \longmapsto (p,q,r,s)$:
\[
\aligned
p
&
\,:=\,
a_{1,1}\,x+a_{1,2}\,y+b_{1,1}\,u+b_{1,2}\,v,
\\
q
&
\,:=\,
a_{2,1}\,x+a_{2,2}\,y+b_{2,1}\,u+b_{2,2}\,v,
\\
r
&
\,:=\,
c_{1,1}\,x+c_{1,2}\,y+d_{1,1}\,u+d_{1,2}\,v,
\\
s
&
\,:=\,
c_{2,1}\,x+c_{2,2}\,y+d_{2,1}\,u+d_{2,2}\,v,
\endaligned
\ \ \ \ \ \ \ \ \ \ \ \ \ \ \ \ \ \ \ \
\text{with}
\ \ \ \ \ \ \ \ \ \ \ \ \ \ \ \ \ \ \ \
0
\,\neq\,
\left\vert\!
\begin{array}{cccc}
a_{1,1} & a_{1,2} & b_{1,1} & b_{1,2}
\\
a_{2,1} & b_{2,2} & b_{2,1} & b_{2,2}
\\
c_{1,1} & c_{1,2} & d_{1,1} & d_{1,2}
\\
c_{2,1} & c_{2,2} & d_{2,1} & d_{2,2}
\end{array}
\!\right\vert.
\]

Also, consider two local analytic surfaces passing through the origin,
graphed as:
\[
\aligned
u
&
\,=\,
F(x,y)
\ \ \ \ \ \ \ \ 
{\scriptstyle{(F(0,0)\,=\,0)}}
\\
v
&
\,=\,
G(x,y)
\ \ \ \ \ \ \ \ 
{\scriptstyle{(G(0,0)\,=\,0)}}
\endaligned
\ \ \ \ \ \ \ \ \ \ \ \
\text{and}
\ \ \ \ \ \ \ \ \ \ \ \
\aligned
r
\,=\,
R(p,q)
\ \ \ \ \ \ \ \ 
{\scriptstyle{(0\,=\,R(0,0))}},
\\
s
\,=\,
S(p,q)
\ \ \ \ \ \ \ \ 
{\scriptstyle{(0\,=\,S(0,0))}},
\endaligned
\]
with convergent series:
\[
\aligned
F
&
\,=\,
\sum_{i+j\geqslant 1}\,
F_{i,j}\,
x^i\,y^j,
\\
G
&
\,=\,
\sum_{i+j\geqslant 1}\,
G_{i,j}\,
x^i\,y^j,
\endaligned
\ \ \ \ \ \ \ \ \ \ \ \ \ \ \ \ \ \ \ \
\text{and}
\ \ \ \ \ \ \ \ \ \ \ \ \ \ \ \ \ \ \ \
\aligned
R
\,=\,
\sum_{k+l\geqslant 1}\,
R_{k,l}\,
p^k\,q^l,
\\
S
\,=\,
\sum_{k+l\geqslant 1}\,
S_{k,l}\,
p^k\,q^l.
\endaligned
\]

The concerned general linear transformation of $\GL(4,\R)$ 
maps the left surface $\{u = F, \, v = G\}$ 
to the right surface $\{r = R,\, s = S\}$ 
if and only if two {\sl fundamental equations:}
\leqnomode\usetagform{default}
\begin{align}
\label{eqR-eqS-S2-R4}
\aligned
0
&
\,\equiv\,
\eqR(x,y),
\\
0
&
\,\equiv\,
\eqS(x,y),
\endaligned
\end{align}
hold identically in $\R\{x,y\}$, where:
\[
\aligned
\eqR
&
\,:=\,
-\,c_{1,1}\,x-c_{1,2}\,y-d_{1,1}\,F(x,y)-d_{1,2}\,G(x,y)
\\
&
\ \ \ \ \ \ \
+
R
\Big(
a_{1,1}x+a_{1,2}y+b_{1,1}F(x,y)+b_{1,2}G(x,y),\,\,
a_{2,1}x+a_{2,2}y+b_{2,1}F(x,y)+b_{2,2}G(x,y)
\Big),
\\
\eqS
&
\,:=\,
-\,c_{2,1}\,x-c_{2,2}\,y-d_{2,1}\,F(x,y)-d_{2,2}\,G(x,y)
\\
&
\ \ \ \ \ \ \
+
S
\Big(
a_{1,1}x+a_{1,2}y+b_{1,1}F(x,y)+b_{1,2}G(x,y),\,\,
a_{2,1}x+a_{2,2}y+b_{2,1}F(x,y)+b_{2,2}G(x,y)
\Big),
\endaligned
\]
so that:
\[
\aligned
0
&
\,=\,
\eqR
\,=\,
\sum_{i,j\in\N}\,
\mathcal{R}_{i,j}
\Big(
a_{\smallbullet,\smallbullet},\,
b_{\smallbullet,\smallbullet},\,
c_{\smallbullet,\smallbullet},\,
d_{\smallbullet,\smallbullet},\,\,
F_{\smallbullet,\smallbullet},\,
G_{\smallbullet,\smallbullet},\,
R_{\smallbullet,\smallbullet}
\Big)\,
x^i\,y^j,
\\
0
&
\,=\,
\eqS
\,=\,
\sum_{i,j\in\N}\,
\mathcal{S}_{i,j}
\Big(
a_{\smallbullet,\smallbullet},\,
b_{\smallbullet,\smallbullet},\,
c_{\smallbullet,\smallbullet},\,
d_{\smallbullet,\smallbullet},\,\,
F_{\smallbullet,\smallbullet},\,
G_{\smallbullet,\smallbullet},\,
S_{\smallbullet,\smallbullet}
\Big)\,
x^i\,y^j.
\endaligned
\]

For all $i,j \in \N$, the coefficient of $x^i y^j$ in
$\eqR$ and in $\eqS$ can, in a standard way, be denoted as:
\[
\aligned
\big[x^i\,y^j\big]
\eqR
&
\,:=\,
\mathcal{R}_{i,j}
\,=\,
0,
\\
\big[x^i\,y^j\big]
\eqS
&
\,:=\,
\mathcal{S}_{i,j}
\,=\,
0,
\endaligned
\]
and we will indicate the corresponding indices 
$\green{\bf i, j}$ over the equal sign as:
\[
\aligned
&
0
\underset{\green{\bf R}}{\overset{\green{\bf i,j}}{\,\,=\,\,}}
\mathcal{R}_{i,j},
\\
&
0
\underset{\green{\bf S}}{\overset{\green{\bf i,j}}{\,\,=\,\,}}
\mathcal{S}_{i,j}.
\endaligned
\]
We will proceed by increasing:
\[
\order
\,:=\,
i+j.
\]

Two simple affine transformations of $\R_{x,y,u,v}^4$
and of $\R_{p,q,r,s}^4$ erase $\ast\, x^1 + \ast\, y^1$
and $\ast\, p^1 + \ast\, q^1$:
\[
\aligned
u
\,=\,
\red{\bf 0}
+
{\rm O}_{x,y}(2),
\\
v
\,=\,
\red{\bf 0}
+
{\rm O}_{x,y}(2),
\endaligned
\ \ \ \ \ \ \ \ \ \ \ \ \ \ \ \ \ \ \ \
\text{and}
\ \ \ \ \ \ \ \ \ \ \ \ \ \ \ \ \ \ \ \
\aligned
r
\,=\,
\red{\bf 0}
+
{\rm O}_{p,q}(2),
\\
s
\,=\,
\red{\bf 0}
+
{\rm O}_{p,q}(2),
\endaligned
\]
so that order 1 terms are {\sl normalized} to be $\red{\bf 0}$.

\begin{Lemma}
\label{Lm-stab-ordre-1}
Stabilization of order 1 terms holds if and only if 
$0 = c_{1,1} = c_{1,2} = c_{2,1} = c_{2,2}$:
\[
G_{\text{stab}}^{0}=
\left[
\begin{array}{cccc}
a_{1,1} & a_{1,2} & b_{1,1} & b_{1,2}
\\
a_{2,1} & a_{2,2} & b_{2,1} & b_{2,2}
\\
c_{1,1} & c_{1,2} & d_{1,1} & d_{1,2}
\\
c_{2,1} & c_{2,2} & d_{2,1} & d_{2,2}
\end{array}
\right]^{\green{\bf 0}}
\,\,\,\leadsto\,\,\,
G_{\text{\rm stab}}^{1}
=
\left[
\begin{array}{cccc}
a_{1,1} & a_{1,2} & b_{1,1} & b_{1,2}
\\
a_{2,1} & a_{2,2} & b_{2,1} & b_{2,2}
\\
\red{\bf 0} & \red{\bf 0} & d_{1,1} & d_{1,2}
\\
\red{\bf 0} & \red{\bf 0} & d_{2,1} & d_{2,2}
\end{array}
\right]^{\green{\bf 1}}.
\]
\end{Lemma}

\proof
Just read~{\eqref{eqR-eqS-S2-R4}} at order 1:
\begin{align}
0
&
\underset{\green{\bf R}}{\overset{\green{\bf 1,0}}{\,\,=\,\,}}
-\,c_{1,1},
\ \ \ \ \ \ \ \ \ \ \ \ \ \ \ \ \ \ \ \
0
\underset{\green{\bf R}}{\overset{\green{\bf 0,1}}{\,\,=\,\,}}
-\,c_{1,2},
\notag
\\
0
&
\underset{\green{\bf S}}{\overset{\green{\bf 1,0}}{\,\,=\,\,}}
-\,c_{2,1},
\ \ \ \ \ \ \ \ \ \ \ \ \ \ \ \ \ \ \ \
0
\underset{\green{\bf S}}{\overset{\green{\bf 0,1}}{\,\,=\,\,}}
-\,c_{2,2}.
\qedhere
\end{align}
\endproof

%%%%%%%%%%%%%%%%%%%%%%%%%%%%%%%%%%%%%%%%%%%%%%%%%%%%%%%%%%%%%%%%%%%%%%
\SectionHead{Parallel and NonParallel Quadratic Forms}
{parallel-non-parallel}
%%%%%%%%%%%%%%%%%%%%%%%%%%%%%%%%%%%%%%%%%%%%%%%%%%%%%%%%%%%%%%%%%%%%%%

At order 2:
\[
\aligned
u
&
\,=\,
F_{2,0}\,x^2
+
F_{1,1}\,xy
+
F_{0,2}\,y^2
+
{\rm O}_{x,y}(3)
\,=:\,
F_2
+\cdots,
\\
v
&
\,=\,
G_{2,0}\,x^2
+
G_{1,1}\,xy
+
G_{0,2}\,y^2
+
{\rm O}_{x,y}(3)
\,=:\,
G_2
+\cdots,
\endaligned
\]
two quadratic forms $F_2$, $G_2$ are present.
Any $2 \times 2$ linear transformation in the $(u,v)$-space,
is allowed, as it corresponds to the block 
$\big[\begin{smallmatrix} 
d_{1,1} & d_{1,2} \\ d_{2,1} & d_{2,2}
\end{smallmatrix}\big]$ above.

\begin{Proposition}
\label{Prp-parallel-quadrics-F2-G2}
The property that $F_2$ and $G_2$ are linearly dependent
is invariant under affine transformations.
\end{Proposition}

Of course, this property expresses as:
\[
0
\,=\,
\left\vert\!
\begin{array}{cc}
F_{2,0} & F_{1,1}
\\
G_{2,0} & G_{1,1}
\end{array}
\!\right\vert
\,=\,
\left\vert\!
\begin{array}{cc}
F_{2,0} & F_{0,2}
\\
G_{2,0} & G_{0,2}
\end{array}
\!\right\vert
\,=\,
\left\vert\!
\begin{array}{cc}
F_{1,1} & F_{0,2}
\\
G_{1,1} & G_{0,2}
\end{array}
\!\right\vert.
\]
More precisely, it is invariant under the subgroup of
$\GL(4,\R) \subset \Aff(4,\R)$ written 
in Lemma~{\ref{Lm-stab-ordre-1}}
above.

\proof
Using from now on a computer,
the equations~{\eqref{eqR-eqS-S2-R4}} at order 2 read as:
\[
\aligned
0
&
\underset{\green{\bf R}}{\overset{\green{\bf 2,0}}{\,\,=\,\,}}
a_{2,1}^2\,R_{0,2}
+
a_{1,1}a_{2,1}\,R_{1,1}
+
a_{1,1}^2\,R_{2,0}
-
d_{1,1}\,F_{2,0}
-
d_{1,2}\,G_{2,0},
\\
0
&
\underset{\green{\bf R}}{\overset{\green{\bf 1,1}}{\,\,=\,\,}}
2\,a_{2,1}a_{2,2}\,R_{0,2}
+
a_{1,1}a_{2,2}\,R_{1,1}
+
a_{1,2}a_{2,1}\,R_{1,1}
+
2\,a_{1,1}a_{1,2}\,R_{2,0}
-
d_{1,1}\,F_{1,1}
-
d_{1,2}\,G_{1,1},
\\
0
&
\underset{\green{\bf R}}{\overset{\green{\bf 0,2}}{\,\,=\,\,}}
a_{2,2}^2\,R_{0,2}
+
a_{1,2}a_{2,2}\,R_{1,1}
+
a_{1,2}^2\,R_{2,0}
-
d_{1,1}\,F_{0,2}
-
d_{1,2}\,G_{0,2},
\\
0
&
\underset{\green{\bf S}}{\overset{\green{\bf 2,0}}{\,\,=\,\,}}
a_{2,1}^2\,S_{0,2}
+
a_{1,1}\,a_{2,1}\,S_{1,1}
+
a_{1,1}^2\,S_{2,0}
-
d_{2,1}\,F_{2,0}
-
d_{2,2}\,G_{2,0},
\\
0
&
\underset{\green{\bf S}}{\overset{\green{\bf 1,1}}{\,\,=\,\,}}
2\,a_{2,1}a_{2,2}\,S_{0,2}
+
a_{1,1}a_{2,2}\,S_{1,1}
+
a_{1,2}a_{2,1}\,S_{1,1}
+
2\,a_{1,1}a_{1,2}\,S_{2,0}
-
d_{2,1}\,F_{1,1}
-
d_{2,2}\,G_{1,1},
\\
0
&
\underset{\green{\bf S}}{\overset{\green{\bf 0,2}}{\,\,=\,\,}}
a_{2,2}^2\,S_{0,2}
+
a_{1,2}a_{2,2}\,S_{1,1}
+
a_{1,2}^2S_{2,0}
-
d_{2,1}\,F_{0,2}
-
d_{2,2}\,G_{0,2}.
\endaligned
\]
From these 6 linear equations, solve the 6 target order 2 coefficients:
\[
\footnotesize
\aligned
R_{2,0}
&
\,=\,
\frac{1}{
(a_{1,1}a_{2,2}-a_{1,2}a_{2,1})^2}\,
\Big\{
a_{2,1}^2d_{1,1}F_{0,2}
-
a_{2,1}a_{2,2}d_{1,1}\,F_{1,1}
+
a_{2,2}^2d_{1,1}\,F_{2,0}
\\
&
\ \ \ \ \ \ \ \ \ \ \ \ \ \ \ \ \ \ \ \ \ \ \ \ \ \ \ \ \ \ \ \ \ \ \
\ \ \ \ \ \
+
a_{2,1}^2d_{1,2}\,G_{0,2}
-
a_{2,1}a_{2,2}d_{1,2}\,G_{1,1}
+
a_{2,2}^2d_{1,2}\,G_{2,0}
\Big\},
\\
R_{1,1}
&
\,=\,
\frac{1}{
(a_{1,1}a_{2,2}-a_{1,2}a_{2,1})^2}\,
\Big\{
-\,2\,a_{1,1}a_{2,1}d_{1,1}\,F_{0,2}
+
a_{1,1}a_{2,2}d_{1,1}\,F_{1,1}
+
a_{1,2}a_{2,1}d_{1,1}\,F_{1,1}
-
2\,a_{1,2}a_{2,2}d_{1,1}\,F_{2,0}
\\
&
\ \ \ \ \ \ \ \ \ \ \ \ \ \ \ \ \ \ \ \ \ \ \ \ \ \ \ \ \ \ \ \ \ \ \
\ \ \ \ \ \
-
2\,a_{1,1}a_{2,1}d_{1,2}\,G_{0,2}
+
a_{1,1}a_{2,2}d_{1,2}\,G_{1,1}
+
a_{1,2}a_{2,1}d_{1,2}\,G_{1,1}
-
2\,a_{1,2}a_{2,2}d_{1,2}\,G_{2,0}
\Big\},
\\
R_{0,2}
&
\,=\,
\frac{1}{
(a_{1,1}a_{2,2}-a_{1,2}a_{2,1})^2}\,
\Big\{
a_{1,1}^2d_{1,1}F_{0,2}
-
a_{1,1}a_{1,2}d_{1,1}\,F_{1,1}
+
a_{1,2}^2d_{1,1}\,F_{2,0}
\\
&
\ \ \ \ \ \ \ \ \ \ \ \ \ \ \ \ \ \ \ \ \ \ \ \ \ \ \ \ \ \ \ \ \ \ \
\ \ \ \ \ \
+
a_{1,1}^2d_{1,2}\,G_{0,2}
-
a_{1,1}a_{1,2}d_{1,2}\,G_{1,1}
+
a_{1,2}^2d_{1,2}\,G_{2,0}
\Big\},
\endaligned
\]
%%%%%%%%%%%%%%%%%%%%%%%%%%%%%%%%%%%%%%%%%%%%%%%%%%%%%%%%%%%%%%%%%%%%%%
\[
\footnotesize
\aligned
S_{2,0}
&
\,=\,
\frac{1}{
(a_{1,1}a_{2,2}-a_{1,2}a_{2,1})^2}\,
\Big\{
a_{2,1}^2d_{2,1}F_{0,2}
-
a_{2,1}a_{2,2}d_{2,1}\,F_{1,1}
+
a_{2,2}^2d_{2,1}\,F_{2,0}
\\
&
\ \ \ \ \ \ \ \ \ \ \ \ \ \ \ \ \ \ \ \ \ \ \ \ \ \ \ \ \ \ \ \ \ \ \
\ \ \ \ \ \
+
a_{2,1}^2d_{2,2}\,G_{0,2}
-
a_{2,1}a_{2,2}d_{2,2}\,G_{1,1}
+
a_{2,2}^2d_{2,2}\,G_{2,0}
\Big\},
\\
S_{1,1}
&
\,=\,
\frac{1}{
(a_{1,1}a_{2,2}-a_{1,2}a_{2,1})^2}\,
\Big\{
-\,2\,a_{1,1}a_{2,1}d_{2,1}\,F_{0,2}
+
a_{1,1}a_{2,2}d_{2,1}\,F_{1,1}
+
a_{1,2}a_{2,1}d_{2,1}\,F_{1,1}
-
2\,a_{1,2}a_{2,2}d_{2,1}\,F_{2,0}
\\
&
\ \ \ \ \ \ \ \ \ \ \ \ \ \ \ \ \ \ \ \ \ \ \ \ \ \ \ \ \ \ \ \ \ \ \
\ \ \ \ \ \
-
2\,a_{1,1}a_{2,1}d_{2,2}\,G_{0,2}
+
a_{1,1}a_{2,2}d_{2,2}\,G_{1,1}
+
a_{1,2}a_{2,1}d_{2,2}\,G_{1,1}
-
2\,a_{1,2}a_{2,2}d_{2,2}\,G_{2,0}
\Big\},
\\
S_{0,2}
&
\,=\,
\frac{1}{
(a_{1,1}a_{2,2}-a_{1,2}a_{2,1})^2}\,
\Big\{
a_{1,1}^2d_{2,1}F_{0,2}
-
a_{1,1}a_{1,2}d_{2,1}\,F_{1,1}
+
a_{1,2}^2d_{2,1}\,F_{2,0}
\\
&
\ \ \ \ \ \ \ \ \ \ \ \ \ \ \ \ \ \ \ \ \ \ \ \ \ \ \ \ \ \ \ \ \ \ \
\ \ \ \ \ \
+
a_{1,1}^2d_{2,2}\,G_{0,2}
-
a_{1,1}a_{1,2}d_{2,2}\,G_{1,1}
+
a_{1,2}^2d_{2,2}\,G_{2,0}
\Big\},
\endaligned
\]
which is allowed, since the determinant of the $4 \times 4$ 
matrix 
in Lemma~{\ref{Lm-stab-ordre-1}}
does not vanish:
\[
\big(
a_{1,1}a_{2,2}
-
a_{1,2}a_{2,1}
\big)\,
\big(
d_{1,1}d_{2,2}
-
d_{1,2}d_{2,1}
\big)
\,\neq\,
0.
\]

\begin{Observation}
\label{Obs-flat}
$0 = F_2 = G_2$ if and only if $R_2 = S_2 = 0$.\qed
\end{Observation}

\proof
The implication $\Longrightarrow$ is clear from the 6 formulas above,
while the reverse implication $\Longleftarrow$ stems from
analogous formulas expressing $F_{i,j}$, $G_{i,j}$ 
with $i + j = 2$, in terms of
$R_{k,l}$, $S_{k,l}$
with $k + l = 2$,
just by means of the {\em inverse} affine
(linear) transformation.
\endproof

Then by direct computations:
\[
\aligned
\left\vert\!
\begin{array}{cc}
R_{2,0} & R_{1,1}
\\
S_{2,0} & S_{1,1}
\end{array}
\!\right\vert
&
\,=\,
\frac{
\left\vert\begin{smallmatrix} 
d_{1,1} & d_{1,2} \\ d_{2,1} & d_{2,2} 
\end{smallmatrix}\right\vert
}{
\left\vert\begin{smallmatrix} 
a_{1,1} & a_{1,2} \\ a_{2,1} & a_{2,2} 
\end{smallmatrix}\right\vert^3
}
\left\{
a_{2,1}^2\,
\left\vert\!
\begin{array}{cc}
F_{1,1} & F_{0,2}
\\
G_{1,1} & G_{0,2}
\end{array}
\!\right\vert
-
2\,a_{2,1}a_{2,2}\,
\left\vert\!
\begin{array}{cc}
F_{2,0} & F_{0,2}
\\
G_{2,0} & G_{0,2}
\end{array}
\!\right\vert
-
a_{2,2}^2\,
\left\vert\!
\begin{array}{cc}
F_{2,0} & F_{1,1}
\\
G_{2,0} & G_{1,1}
\end{array}
\!\right\vert
\right\},
\\%%%%%
\left\vert\!
\begin{array}{cc}
R_{2,0} & R_{0,2}
\\
S_{2,0} & S_{0,2}
\end{array}
\!\right\vert
&
\,=\,
\frac{
\left\vert\begin{smallmatrix} 
d_{1,1} & d_{1,2} \\ d_{2,1} & d_{2,2} 
\end{smallmatrix}\right\vert
}{
\left\vert\begin{smallmatrix} 
a_{1,1} & a_{1,2} \\ a_{2,1} & a_{2,2} 
\end{smallmatrix}\right\vert^3
}
\left\{
-\,a_{1,1}a_{2,1}\,
\left\vert\!
\begin{array}{cc}
F_{1,1} & F_{0,2}
\\
G_{1,1} & G_{0,2}
\end{array}
\!\right\vert
+
\big(a_{1,1}a_{2,2}+a_{1,2}a_{2,1}\big)\,
\left\vert\!
\begin{array}{cc}
F_{2,0} & F_{0,2}
\\
G_{2,0} & G_{0,2}
\end{array}
\!\right\vert
\right.
\\
&
\ \ \ \ \ \ \ \ \ \ \ \ \ \ \ \ \ \ \ \ \ \ \ \ \ \ \ \ \ \ \ \ \ \ \
\ \ \ \ \ \ \ \ \ \ \ \ \ \ \ \ \ \ \ \ \ \ \ \ \ \ \ \ \ \ \ \ \ \ \
\ \ \ \ \ \ \ \ \ \ \ \ \ \ \ \ 
\left.
-\,
a_{1,2}a_{2,2}\,
\left\vert\!
\begin{array}{cc}
F_{2,0} & F_{1,1}
\\
G_{2,0} & G_{1,1}
\end{array}
\!\right\vert
\right\},
\\%%%%%
\left\vert\!
\begin{array}{cc}
R_{1,1} & R_{0,2}
\\
S_{1,1} & S_{0,2}
\end{array}
\!\right\vert
&
\,=\,
\frac{
\left\vert\begin{smallmatrix} 
d_{1,1} & d_{1,2} \\ d_{2,1} & d_{2,2} 
\end{smallmatrix}\right\vert
}{
\left\vert\begin{smallmatrix} 
a_{1,1} & a_{1,2} \\ a_{2,1} & a_{2,2} 
\end{smallmatrix}\right\vert^3
}
\left\{
a_{1,1}^2\,
\left\vert\!
\begin{array}{cc}
F_{1,1} & F_{0,2}
\\
G_{1,1} & G_{0,2}
\end{array}
\!\right\vert
-
2\,a_{1,1}a_{1,2}\,
\left\vert\!
\begin{array}{cc}
F_{2,0} & F_{0,2}
\\
G_{2,0} & G_{0,2}
\end{array}
\!\right\vert
-
a_{1,2}^2\,
\left\vert\!
\begin{array}{cc}
F_{2,0} & F_{1,1}
\\
G_{2,0} & G_{1,1}
\end{array}
\!\right\vert
\right\},
\endaligned
\]
Hence it is visible that:
\[
\left(\,\,\,
0
=
\left\vert\begin{smallmatrix}
F_{1,1} & F_{0,2}
\\
G_{1,1} & G_{0,2}
\end{smallmatrix}\right\vert
=
\left\vert\begin{smallmatrix}
F_{2,0} & F_{0,2}
\\
G_{2,0} & G_{0,2}
\end{smallmatrix}\right\vert
=
\left\vert\begin{smallmatrix}
F_{2,0} & F_{1,1}
\\
G_{2,0} & G_{1,1}
\end{smallmatrix}\right\vert
\,\,\,\right)
\ \ \ \ \ 
\Longrightarrow
\ \ \ \ \ 
\left(\,\,\,
\left\vert\begin{smallmatrix}
R_{1,1} & R_{0,2}
\\
S_{1,1} & S_{0,2}
\end{smallmatrix}\right\vert
=
\left\vert\begin{smallmatrix}
R_{2,0} & R_{0,2}
\\
S_{2,0} & S_{0,2}
\end{smallmatrix}\right\vert
=
\left\vert\begin{smallmatrix}
R_{2,0} & R_{1,1}
\\
S_{2,0} & S_{1,1}
\end{smallmatrix}\right\vert
=
0
\,\,\,\right).
\]
The reverse implication can be proved in exactly the same way,
using the inverse affine (linear) transformation.
\endproof

Once a normalization has been made on the right
space $\R^4 \ni (x,y,u,v)$,
always, it can also be made {\em exactly the same} on the left
space $(p,q,r,s) \in \R^4$.

\begin{Principle}
\label{Principle-FG-RS}
{\sl At each order, every performed normalization will always 
be instantly achieved on both hypersurfaces
$\{u = F,\, v=G\}$ and $\{r = R,\, s=S\}$.}\qed
\end{Principle}

Assume now that $F_2$ and $G_2$ are linearly dependent,
whence $R_2$ and $S_2$ are also linearly dependent. 
After simple linear 
transformations in the $(u,v)$-space and in the $(r,s)$-space,
we make $G_2 = 0$ and $S_2 = 0$. 

\begin{Proposition}
\label{Prp-parallel-rank-Hessian}
When $G_2 = 0 = S_2$, the rank of the Hessian matrix of $F$, $R$
at the origin is invariant:
\[
\rank\,
\left[
\begin{smallmatrix}
F_{xx}(0) & F_{xy}(0)
\\
F_{yx}(0) & F_{yy}(0)
\end{smallmatrix}
\right]
\,=\,
\rank\,
\left[
\begin{smallmatrix}
R_{xx}(0) & R_{xy}(0)
\\
R_{yx}(0) & R_{yy}(0)
\end{smallmatrix}
\right].
\]
\end{Proposition}

\proof
Assume therefore that $G_{2,0} = G_{1,1} = G_{0,2} = 0 = 
S_{2,0} = S_{1,1} = S_{0,2}$. 
Since Observation~{\ref{Obs-flat}} already showed that 
$\big[ \begin{smallmatrix} 
2\,F_{2,0} & F_{11} \\
F_{1,1} & 2\,F_{0,2} 
\end{smallmatrix} \big] 
= 
\left[ \begin{smallmatrix} 
0 & 0 \\ 0 & 0 
\end{smallmatrix} \right]$
if and only if
$\big[ \begin{smallmatrix} 
0 & 0 \\ 0 & 0 
\end{smallmatrix} \big]
=
\big[ \begin{smallmatrix} 
2\,R_{2,0} & R_{11} \\
R_{1,1} & 2\,R_{0,2} 
\end{smallmatrix} \big]$, we can assume that
$F_{2,0}$, $F_{1,1}$, $F_{0,2}$ are not all $0$.

The preceding resolutions then become:
\[
\footnotesize
\aligned
R_{2,0}
&
\,=\,
\frac{1}{
(a_{1,1}a_{2,2}-a_{1,2}a_{2,1})^2}\,
\Big\{
a_{2,1}^2d_{1,1}F_{0,2}
-
a_{2,1}a_{2,2}d_{1,1}\,F_{1,1}
+
a_{2,2}^2d_{1,1}\,F_{2,0}
\Big\},
\\
R_{1,1}
&
\,=\,
\frac{1}{
(a_{1,1}a_{2,2}-a_{1,2}a_{2,1})^2}\,
\Big\{
-\,2\,a_{1,1}a_{2,1}d_{1,1}\,F_{0,2}
+
a_{1,1}a_{2,2}d_{1,1}\,F_{1,1}
+
a_{1,2}a_{2,1}d_{1,1}\,F_{1,1}
-
2\,a_{1,2}a_{2,2}d_{1,1}\,F_{2,0}
\Big\},
\\
R_{0,2}
&
\,=\,
\frac{1}{
(a_{1,1}a_{2,2}-a_{1,2}a_{2,1})^2}\,
\Big\{
a_{1,1}^2d_{1,1}F_{0,2}
-
a_{1,1}a_{1,2}d_{1,1}\,F_{1,1}
+
a_{1,2}^2d_{1,1}\,F_{2,0}
\Big\},
\endaligned
\]
and a direct computation provides the nice formula:
\[
\left\vert\!
\begin{array}{cc}
R_{xx}(0) & R_{xy}(0)
\\
R_{yx}(0) & R_{yy}(0)
\end{array}
\!\right\vert
\,=\,
4\,R_{2,0}R_{0,2}
-
R_{1,1}^2
\,=\,
\frac{d_{1,1}^2}{
(a_{1,1}a_{2,2}-a_{1,2}a_{2,1})^2}\,
\Big(
4\,F_{2,0}F_{0,2}
-
F_{1,1}^2
\Big).
\]
But since on the other hand:
\[
\aligned
0
&
\underset{\green{\bf S}}{\overset{\green{\bf 2,0}}{\,\,=\,\,}}
-\,d_{2,1}\,F_{2,0},
\\
0
&
\underset{\green{\bf S}}{\overset{\green{\bf 1,1}}{\,\,=\,\,}}
-\,d_{2,1}\,F_{1,1},
\\
0
&
\underset{\green{\bf S}}{\overset{\green{\bf 0,2}}{\,\,=\,\,}}
-\,d_{2,1}\,F_{0,2},
\endaligned
\]
it is clear that $d_{2,1} = 0$,
whence the nonzero determinant 
of Lemma~{\ref{Lm-stab-ordre-1}}
becomes:
\[
\big(
a_{1,1}a_{2,2}
-
a_{1,2}a_{2,1}
\big)\,
d_{1,1}d_{2,2},
\]
guaranteeing that $d_{1,1} \neq 0$,
and in conclusion, the formula:
\[
\left\vert\!
\begin{array}{cc}
R_{xx}(0) & R_{xy}(0)
\\
R_{yx}(0) & R_{yy}(0)
\end{array}
\!\right\vert
\,\propto\,
\left\vert\!
\begin{array}{cc}
F_{xx}(0) & F_{xy}(0)
\\
F_{yx}(0) & F_{yy}(0)
\end{array}
\!\right\vert.
\]
shows the invariancy of the property that the Hessian has rank 2.

Beyond, one may verify that:
\[
\left[\!
\begin{array}{cc}
2\,R_{2,0} & R_{1,1}
\\
R_{1,1} & 2\,R_{0,2}
\end{array}
\!\right]
\,=\,
\frac{d_{1,1}}{(a_{1,1}a_{2,2}-a_{1,2}a_{2,1})^2}\,
\left[\!
\begin{array}{cc}
a_{2,2} & -a_{2,1}
\\
-a_{1,2} & a_{1,1}
\end{array}
\!\right]
\cdot
\left[\!
\begin{array}{cc}
2\,F_{2,0} & F_{1,1}
\\
F_{1,1} & 2\,F_{0,2}
\end{array}
\!\right]
\cdot
\left[\!
\begin{array}{cc}
a_{2,2} & -a_{1,2}
\\
-a_{2,1} & a_{1,1}
\end{array}
\!\right]
\]
and since the two $a_{i,j}$-matrices are invertible, it is clear that:
\[
\rank\,
\left[\!
\begin{array}{cc}
2\,R_{2,0} & R_{1,1}
\\
R_{1,1} & 2\,R_{0,2}
\end{array}
\!\right]
\,=\,
\rank\,
\left[\!
\begin{array}{cc}
2\,F_{2,0} & F_{1,1}
\\
F_{1,1} & 2\,F_{0,2}
\end{array}
\!\right].
\qedhere
\]
\endproof

%%%%%%%%%%%%%%%%%%%%%%%%%%%%%%%%%%%%%%%%%%%%%%%%%%%%%%%%%%%%%%%%%%%%%%
\SectionHead{Parallel Quadrics at Order 2}
{parallel-quadrics-order-2}
%%%%%%%%%%%%%%%%%%%%%%%%%%%%%%%%%%%%%%%%%%%%%%%%%%%%%%%%%%%%%%%%%%%%%%

When the two quadratic forms $F_2$, $G_2$ in
the left space are {\sl parallel},
{\em i.e.} colinear, so that $R_2$, $S_2$ in the right space
are also parallel, we can assume that $G_2 = 0$ and that $S_2 = 0$.
Then by Proposition~{\ref{Prp-parallel-rank-Hessian}},
the Hessians at the origin of $F$ and of $R$ have the same rank.

It is well known that, under the group 
$\GL(2,\R) \ni \big( \begin{smallmatrix} \alpha & \beta \\
\gamma & \delta \end{smallmatrix} \big)$, with $x \longmapsto
\alpha\,x + \beta\,y$ and $y \longmapsto \gamma\, x + \delta\,y$,
any quadratic form $F_{2,0}\, x^2 + F_{1,1}\, xy + F_{0,2}\, y^2$
can be normalized in 4 inequivalent ways as:
\[
\def\arraystretch{1.25}
\begin{array}{llll}
{\sf Zero} & {\sf Parabolic} & {\sf Hyperbolic} & {\sf Elliptic}
\\
0 & x^2 & x\,y & x^2+y^2
\end{array}
\]
the quadratic form $x\,y$, of signature $(1, 1)$,
being equivalent to $x^2 - y^2$, 

Visibly, the subgroup of $\GL(4,\R)$ written in 
Lemma~{\ref{Lm-stab-ordre-1}} which stabilizes 
the normalization of order 1 monomials, 
has an upper-left right $2 \times 2$ block
constituted of an arbitrary element of $\GL(2,\R)$, so that
we can use matrices of the subgroup form:
\[
\left(
\begin{smallmatrix}
a_{1,1} & a_{1,2} & 0 & 0
\\
a_{2,1} & a_{2,2} & 0 & 0
\\
0 & 0 & 1 & 0
\\
0 & 0 & 0 & 1
\end{smallmatrix}
\right)
\]
to normalize $F_2$ in the left space, and then do exactly
the same in the right space, {\em cf.} 
Principle~{\ref{Principle-FG-RS}}.

\begin{Proposition}
In the parallel case, 4 nonequivalent branches 
exist at order 2:
\[
\text{\bf Branch 2a}
\ \ \ \ \ \ \ \ \ \ \ \ \ \ \ \ \ \ \ \
\aligned
u
&
\,=\,
0
+
{\rm O}_{x,y}(3),
\\
v
&
\,=\,
0
+
{\rm O}_{x,y}(3);
\endaligned
\]
\[
\text{\bf Branch 2b}
\ \ \ \ \ \ \ \ \ \ \ \ \ \ \ \ \ \ \ \
\aligned
u
&
\,=\,
x^2
+
{\rm O}_{x,y}(3),
\\
v
&
\,=\,
0
+
{\rm O}_{x,y}(3);
\endaligned
\]
\[
\text{\bf Branch 2c}
\ \ \ \ \ \ \ \ \ \ \ \ \ \ \ \ \ \ \ \
\aligned
u
&
\,=\,
x\,y
+
{\rm O}_{x,y}(3),
\\
v
&
\,=\,
0
+
{\rm O}_{x,y}(3);
\endaligned
\]
\reqnomode\usetagform{EngelLie}
\begin{align}
\text{\bf Branch 2d}
\ \ \ \ \ \ \ \ \ \ \ \ \ \ \ \ \ \ \ \
u
&
\,=\,
x^2+y^2
+
{\rm O}_{x,y}(3),
\notag
\\
v
&
\,=\,
0
+
{\rm O}_{x,y}(3).
\tag{\qed}
\end{align}
\end{Proposition}

%%%%%%%%%%%%%%%%%%%%%%%%%%%%%%%%%%%%%%%%%%%%%%%%%%%%%%%%%%%%%%%%%%%%%%
\SectionHead{Non-Parallel Quadrics at Order 2}
{non-parallel-quadrics-order-2}
%%%%%%%%%%%%%%%%%%%%%%%%%%%%%%%%%%%%%%%%%%%%%%%%%%%%%%%%%%%%%%%%%%%%%%

It remains to treat the case when no linear relation
exists between $F_2$ and $G_2$.

\begin{Lemma}
In the non-parallel case, after a linear transformation 
in the group $\GL_2 \times \GL_2$
of matrices:
\[
\left(
\begin{smallmatrix}
a_{1,1} & a_{1,2} & 0 & 0
\\
a_{2,1} & a_{2,2} & 0 & 0
\\
0 & 0 & d_{1,1} & d_{1,2}
\\
0 & 0 & d_{2,1} & d_{2,2}
\end{smallmatrix}
\right),
\]
one can always normalize $F_2 = x\,y$.
\end{Lemma}

\proof
It suffices to show that 
one can normalize $F_2 = x\,y$ or $G_2 = x\,y$, since 
the exchange $u \longleftrightarrow v$ is allowed.

There are two cases:
\[
\text{Case~1:}
\ \ \ \ \
\big(
\begin{smallmatrix}
F_{2,0}
\\
G_{2,0}
\end{smallmatrix}
\big)
\,\neq\,
\big(
\begin{smallmatrix}
0
\\
0
\end{smallmatrix}
\big);
\ \ \ \ \ \ \ \ \ \ \ \ \ \ \ \ \ \ \ \
\ \ \ \ \ \ \ \ \ \ \ \ \ \ \ \ \ \ \ \
\text{Case~2:}
\ \ \ \ \
\big(
\begin{smallmatrix}
F_{2,0}
\\
G_{2,0}
\end{smallmatrix}
\big)
\,=\,
\big(
\begin{smallmatrix}
0
\\
0
\end{smallmatrix}
\big)
\,\neq\,
\big(
\begin{smallmatrix}
F_{1,1}
\\
G_{1,1}
\end{smallmatrix}
\big);
\]
and no more,
because $F_2 = 0 + 0 + F_{0,2}\,y^2$ would be parallel to 
$G_2 = 0 + 0 + G_{0,2}\,y^2$. 

In Case~1, after $u \longleftrightarrow v$, we may assume $G_{2,0}
\neq 0$.  By $v \longmapsto \frac{1}{G_{2,0}} v$, and by $u
\longmapsto u - F_{2,0}\,v$, we make:
\[
\aligned
u
&
\,=\,
\ \ \ \ \ \ \ \ \ 
F_{1,1}\,x\,y
+
F_{0,2}\,y^2,
\\
v
&
\,=\,
x^2
+
G_{1,1}\,x\,y
+
G_{0,2}\,y^2.
\endaligned
\]
where, to lighten the writing, we drop ``$+{\rm O}_{x,y}(3)$''.

Subcase 1-1: When $F_{1,1} \neq 0$, by $u \longmapsto 
\frac{1}{F_{1,1}}\,u$, and by introducing a new $x$, we are done:
\[
u
\,=\,
x\,y
+
\tfrac{F_{0,2}}{F_{1,1}}\,y^2
\,=\,
\big(
x
+
\tfrac{F_{0,2}}{F_{1,1}}\,y
\big)\,y
\,=:\,
x\,y.
\]

Subcase 1-2: If $F_{1,1} = 0$, whence $F_{0,2} \neq 0$ due to
non-parallelness, by $u \longmapsto \frac{1}{F_{0,2}}\, u$,
and by $v \longmapsto v - G_{0,2}\,u$, we are done again:
\[
v
\,=\,
x^2
+
G_{1,1}\,x\,y
\,=\,
x\,
\big(
x
+
G_{1,1}\,y
\big)
\,=:\,
x\,y.
\]

In Case~2, after $u \longleftrightarrow v$, we may assume $F_{1,1}
\neq 0$, and by $u \longmapsto \frac{1}{F_{1,1}}\, u$, we are
concluded:
\[
u
\,=\,
x\,y
+
\tfrac{F_{0,2}}{F_{1,1}}\,
y^2
\,=\,
\big(
x
+
\tfrac{F_{0,2}}{F_{1,1}}\,
y
\big)\,
y
\,=:\,
x\,y.
\qedhere
\]
\endproof

Therefore $F_2 = x\,y$, and by $v \longmapsto v - G_{1,1}\,u$, we
have:
\[
\aligned
u
&
\,=\,
\ \ \ \ \ \ \ \ \ \ \ \ \ \ \ \ \ 
x\,y
\ \ \ \ \ \ \ \ \ \ \ \ \ \ \ \ \ \ \ \
+\cdots,
\\
v
&
\,=\,
G_{2,0}\,x^2
\ \ \ \ \ \ \ \ \ \ \ \ \ \ \ 
+
G_{0,2}\,y^2
+\cdots,
\endaligned
\]
the $+\cdots$ being some ${\rm O}_{x,y}(3)$ remainder, of course.
By non-parallelness, $0 \neq G_2$, and after 
$x \longleftrightarrow y$,
we may assume $G_{2,0} \neq 0$. Next, by $v \longmapsto
\frac{1}{G_{2,0}}\, v$, we receive:
\[
\aligned
u
&
\,=\,
\ \ \ \ \ \ \ \
x\,y
\ \ \ \ \ \ \ \ \ \ \ \ \ \ \ \ \ \
+\cdots,
\\
v
&
\,=\,
x^2
\ \ \ \ \ \ \ \ \ \ \
+
G_{0,2}\,y^2
+\cdots,
\endaligned
\]

\begin{Question}
{\sl But then, what about the last remaining coefficient $G_{0,2}$?}
\end{Question}

Up to now, we have only used elementary affine transformations, which
is a (too) naive process, but in principle, all the $3 + 3$ equations:
\[
\aligned
0
&
\underset{\green{\bf R}}{\overset{\green{\bf 2,0}}{\,\,=\,\,}}
\mathcal{R}_{2,0},
&
\ \ \ \ \ \ \ \ \ \ \ \ \ \ \ \ \ \ \ \
0
&
\underset{\green{\bf R}}{\overset{\green{\bf 1,1}}{\,\,=\,\,}}
\mathcal{R}_{1,1},
&
\ \ \ \ \ \ \ \ \ \ \ \ \ \ \ \ \ \ \ \
0
&
\underset{\green{\bf R}}{\overset{\green{\bf 0,2}}{\,\,=\,\,}}
\mathcal{R}_{0,2},
\\
0
&
\underset{\green{\bf S}}{\overset{\green{\bf 2,0}}{\,\,=\,\,}}
\mathcal{S}_{2,0},
&
\ \ \ \ \ \ \ \ \ \ \ \ \ \ \ \ \ \ \ \
0
&
\underset{\green{\bf S}}{\overset{\green{\bf 1,1}}{\,\,=\,\,}}
\mathcal{S}_{1,1},
&
\ \ \ \ \ \ \ \ \ \ \ \ \ \ \ \ \ \ \ \
0
&
\underset{\green{\bf S}}{\overset{\green{\bf 0,2}}{\,\,=\,\,}}
\mathcal{S}_{0,2},
\endaligned
\]
should be dealt with.

\begin{Proposition}
Through an affine (linear) transformation which stabilizes 
the normalizations achieved up to now:
\[
\aligned
u
&
\,=\,
\ \ \ \ \ \ \ \
x\,y
\ \ \ \ \ \ \ \ \ \ \ \ \ \ \ \ \ \
+\cdots,
\\
v
&
\,=\,
x^2
\ \ \ \ \ \ \ \ \ \ \
+
G_{0,2}\,y^2
+\cdots,
\endaligned
\ \ \ \ \ \ \ \ \ \
\xrightarrow[{\rule[0pt]{50pt}{0pt}}]{\text{\rm Transformation}}
\ \ \ \ \ \ \ \ \ \
\aligned
r
&
\,=\,
\ \ \ \ \ \ \ \
p\,q
\ \ \ \ \ \ \ \ \ \ \ \ \ \ \ \ \ \
+\cdots,
\\
s
&
\,=\,
p^2
\ \ \ \ \ \ \ \ \ \ \
+
S_{0,2}\,q^2
+\cdots,
\endaligned
\]
it holds with $a_{1,1} \neq 0 \neq a_{2,2}$ that the only remaining
coefficient is a relative invariant:
\[
S_{0,2}
\,=\,
\frac{a_{1,1}^2}{a_{2,2}^2}\,
G_{0,2}.
\]
\end{Proposition}

\proof
So we write the $3 + 3$ equations in question:
\[
\aligned
0
&
\underset{\green{\bf R}}{\overset{\green{\bf 2,0}}{\,\,=\,\,}}
-\,d_{1,2}
+
a_{1,1}\,a_{2,1},
\\
0
&
\underset{\green{\bf R}}{\overset{\green{\bf 1,1}}{\,\,=\,\,}}
-\,d_{1,1}
+
a_{1,1}\,a_{2,2}
+
a_{2,1}\,a_{1,2},
\\
0
&
\underset{\green{\bf R}}{\overset{\green{\bf 0,2}}{\,\,=\,\,}}
-\,d_{1,2}\,G_{0,2}
+
a_{1,2}\,a_{2,2},
\\
0
&
\underset{\green{\bf S}}{\overset{\green{\bf 2,0}}{\,\,=\,\,}}
-\,d_{2,2}
+
(a_{2,1})^2\,S_{0,2}
+
(a_{1,1})^2,
\\
0
&
\underset{\green{\bf S}}{\overset{\green{\bf 1,1}}{\,\,=\,\,}}
-\,d_{2,1}
+
2\,a_{2,1}\,a_{2,2}\,S_{0,2}
+
2\,a_{1,1}\,a_{1,2},
\\
0
&
\underset{\green{\bf S}}{\overset{\green{\bf 0,2}}{\,\,=\,\,}}
-\,d_{2,2}\,G_{0,2}
+
(a_{2,2})^2\,S_{0,2}
+
(a_{1,2})^2.
\endaligned
\]

From the third equation, we solve using also the first equation:
\[
a_{1,2}
\,:=\,
\frac{a_{1,1}\,a_{2,1}}{a_{2,2}}\,
G_{0,2}.
\]
From the fourth equation, we solve:
\[
d_{2,2}
\,:=\,
(a_{2,1})^2\,S_{0,2}
+
(a_{1,1})^2.
\]
We replace these two values in the sixth equation, 
and we factorize it:
\[
0
\,=\,
\frac{a_{2,2}^2-a_{2,1}^2\,G_{0,2}}{a_{2,2}^2}\,\,
\Big(
-\,a_{1,1}^2\,G_{0,2}
+
a_{2,2}^2\,S_{0,2}
\Big).
\]

But we remind that the determinant of the matrix of 
Lemma~{\ref{Lm-stab-ordre-1}} is nonzero.
So after the replacement of $a_{1,2} :=
\frac{a_{1,1}\, a_{2,1}}{ a_{2,2}}\,
G_{0,2}$:
\[
\aligned
&
0
\,\neq\,
a_{1,1}a_{2,2}
-
a_{1,2}a_{2,1}
\\
&
\,=\,
\frac{a_{1,1}}{a_{2,2}}\,
\big(
a_{2,2}^2
-
a_{2,1}^2\,G_{0,2}
\big),
\endaligned
\]
we observe that the first factor above is nonvanishing, 
hence we conclude that:
\[
0
\,=\,
-\,a_{1,1}^2\,G_{0,2}
+
a_{2,2}^2\,S_{0,2}.
\qedhere
\]
\endproof

Since the dilation coefficient $\frac{a_{1,1}^2}{a_{2,2}^2}$
in the formula 
$S_{0,2}
\,=\,
\frac{a_{1,1}^2}{ a_{2,2}^2}\,
G_{0,2}$ is positive, and since $G_{0,2}, S_{0,2} \in \R$ are real,
3 normalizing values exist: $0$; $+1$, $-1$.

\begin{Proposition}
In the non-parallel case, 3 nonequivalent branches 
exist at order 2:
\[
\text{\bf Branch 2e}
\ \ \ \ \ \ \ \ \ \ \ \ \ \ \ \ \ \ \ \
\aligned
u
&
\,=\,
x\,y
+
{\rm O}_{x,y}(3),
\\
v
&
\,=\,
x^2
+
{\rm O}_{x,y}(3);
\endaligned
\]
\[
\text{\bf Branch 2f}
\ \ \ \ \ \ \ \ \ \ \ \ \ \ \ \ \ \ \ \
\aligned
u
&
\,=\,
x\,y
+
{\rm O}_{x,y}(3),
\\
v
&
\,=\,
x^2+y^2
+
{\rm O}_{x,y}(3);
\endaligned
\]
\reqnomode\usetagform{EngelLie}
\begin{align}
\text{\bf Branch 2g}
\ \ \ \ \ \ \ \ \ \ \ \ \ \ \ \ \ \ \ \
u
&
\,=\,
x\,y
+
{\rm O}_{x,y}(3),
\notag
\\
v
&
\,=\,
x^2-y^2
+
{\rm O}_{x,y}(3).
\tag{\qed}
\end{align}
\end{Proposition}

%%%%%%%%%%%%%%%%%%%%%%%%%%%%%%%%%%%%%%%%%%%%%%%%%%%%%%%%%%%%%%%%%%%%%%
\SectionHead{Infinitesimal Affine Transformations}
{infinitesimal-affine-transformations}
%%%%%%%%%%%%%%%%%%%%%%%%%%%%%%%%%%%%%%%%%%%%%%%%%%%%%%%%%%%%%%%%%%%%%%

A general affine vector field writes:
\[
\aligned
L
&
\,=\,
\ \ 
\big(
T_1+A_{1,1}\,x+A_{1,2}\,y+B_{1,1}\,u+B_{1,2}\,v
\big)\,\frac{\partial}{\partial x}
\\
&
\ \ \ \ \
+
\big(
T_2+A_{2,1}\,x+A_{2,2}\,y+B_{2,1}\,u+B_{2,2}\,v
\big)\,\frac{\partial}{\partial y}
\\
&
\ \ \ \ \
+
\big(
U_0+C_{1,1}\,x+C_{1,2}\,y+D_{1,1}\,u+D_{1,2}\,v
\big)\,\frac{\partial}{\partial u}
\\
&
\ \ \ \ \
+
\big(
V_0+C_{2,1}\,x+C_{2,2}\,y+D_{2,1}\,u+D_{2,2}\,v
\big)\,\frac{\partial}{\partial v}.
\endaligned
\]
It is tangent to $\big\{u = F(x,y),\,\, v = G(x,y) \big\}$ 
if and only if:
\[
\aligned
0
&
\,\equiv\,
\eqLF(x,y)
\,=:\,
L
\big(
-\,u+F(x,y)
\big)
\Big\vert_{u=F(x,y)\atop v=G(x,y)},
\\
0
&
\,\equiv\,
\eqLG(x,y)
\,=:\,
L
\big(
-\,v+G(x,y)
\big)
\Big\vert_{u=F(x,y)\atop v=G(x,y)},
\endaligned
\]
identically as power series in $\R\{x,y\}$.
With increasing orders $\mu = 0, 1, 2, 3, \dots$, 
these $\eqLF$ and $\eqLG$ may be expanded:
\[
\aligned
0
&
\,=\,
\eqLF
&
\,=\,
\sum_{\mu=0}^\infty\,\,
\sum_{i+j=\mu}\,
\mathcal{F}_{i,j}
\Big(
A_{\smallbullet,\smallbullet},\,
B_{\smallbullet,\smallbullet},\,
C_{\smallbullet,\smallbullet},\,
D_{\smallbullet,\smallbullet},\,\,
F_{\smallbullet,\smallbullet}
\Big)\,
x^i\,y^j,
\\
0
&
\,=\,
\eqLG
&
\,=\,
\sum_{\mu=0}^\infty\,\,
\sum_{i+j=\mu}\,
\mathcal{G}_{i,j}
\Big(
A_{\smallbullet,\smallbullet},\,
B_{\smallbullet,\smallbullet},\,
C_{\smallbullet,\smallbullet},\,
D_{\smallbullet,\smallbullet},\,\,
G_{\smallbullet,\smallbullet}
\Big)\,
x^i\,y^j.
\endaligned
\]
so that, for all $i, j \in \N$:
\[
\aligned
&
0
\underset{\green{\bf F}}{\overset{\green{\bf i,j}}{\,\,=\,\,}}
\mathcal{F}_{i,j},
\\
&
0
\underset{\green{\bf G}}{\overset{\green{\bf i,j}}{\,\,=\,\,}}
\mathcal{G}_{i,j}.
\endaligned
\]

Such a vector field $L$ is tangent to: 
\[
\aligned
u
&
\,=\,
{\rm O}_{x,y}(2),
\\
v
&
\,=\,
{\rm O}_{x,y}(2),
\endaligned
\]
if and only if:
\[
\aligned
0
&
\underset{\green{\bf F}}{\overset{\green{\bf 0,0}}{\,\,=\,\,}}
-\,U_0,
\\
0
&
\underset{\green{\bf G}}{\overset{\green{\bf 0,0}}{\,\,=\,\,}}
-\,V_0.
\endaligned
\]

The key constraint of transitivity:
\[
\Span\,
\big(
\tfrac{\partial}{\partial x},\,
\tfrac{\partial}{\partial y}
\big)
\,=\,
T_{\sf origin}S
\,=\,
\Span\,
L
\big\vert_{\sf origin}
\,=\,
\Span\,
\Big(
T_1\,
\tfrac{\partial}{\partial x}
+
T_2\,
\tfrac{\partial}{\partial y}
\Big),
\]
forces to always keep $T_1, T_2$ absolutely free\,\,---\,\,never
solved.

%%%%%%%%%%%%%%%%%%%%%%%%%%%%%%%%%%%%%%%%%%%%%%%%%%%%%%%%%%%%%%%%%%%%%%
\SectionHead{Branching diagrams}
{branching-diagrams}
%%%%%%%%%%%%%%%%%%%%%%%%%%%%%%%%%%%%%%%%%%%%%%%%%%%%%%%%%%%%%%%%%%%%%%

In this section, before entering any detailed description,
for clarity and readability, 
we show
all the branching diagrams at higher orders 3, 4, 5, \dots, 
whose terminal nodes, framed in red color, 
point out all the possible existing (families of)
affinely homogeneous models of surfaces $S^2 \subset \R^4$,
in our finallized classification.

\begin{Theorem}
\label{Thm-S2-R4}
Up to the action of a finite subgroup of $\Aff(\R^4)$,
all the possible 
affinely homogeneous
model surfaces $S^2 \subset \R^4$
are represented in a concise manner by 3 lists.

\smallskip\noindent$\bullet$\,
The list of branching trees appearing
{\em infra} in the present 
Section~{\ref{branching-diagrams}}.

\smallskip\noindent$\bullet$\,
The list of occuring linear representations appearing
in Section~{\ref{linear-representations-branches}}.

\smallskip\noindent$\bullet$\,
The lis of (truncated) normal forms together with
their associated
transitive Lie algebras of vector fields,
appearing in 
Sections~{\ref{2b-models}} $\to$
~{\ref{2g-models}}.

\end{Theorem}

We `captured' all
models, even if some of them are parametrized
by very complicated algebraic varieties\,\,---\,\,a novel aspect,
definitely, of our contribution\,\,---\,\,especially
in the simply transitive case.

Branch \green{\bf 2a} brings 1 homogeneous model,
whose details are presented in Section~{\ref{branch-2a}}.

\begin{center}
\includegraphics[scale=0.2]{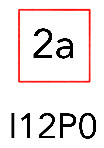}
\end{center}

\noindent
Here, the isotropy Lie subalgebra of (infinitesimal) transformations
fixing (vanishing at) the origin has dimension 12,
while the number of free parameters is 0, and more specifically,
there is a {\em single} model, framed in red.

\smallskip

Next, Branch \green{\bf 2b} brings 15 (families of) homogeneous model.

\begin{center}
\includegraphics[scale=0.19]{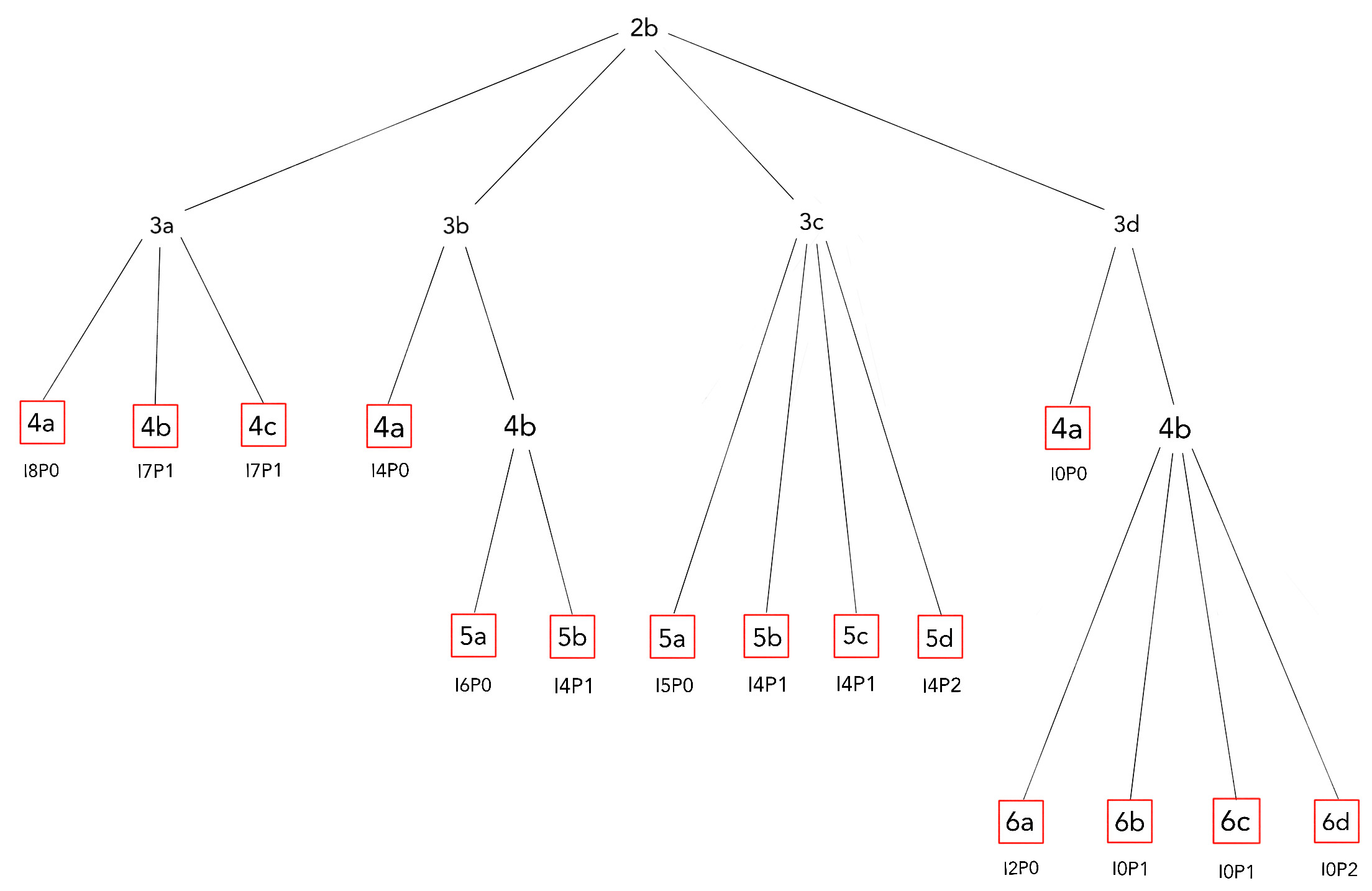}
\end{center}

\noindent
Here for instance, \green{\bf 2b3c5d} indicates a family of 
homogeneous models, parametrized by a
certain affine algebraic variety of
dimension 2, as is pointed out by "P2", all concerned models having
4-dimensional isotropy,
which is what "I4" point out.

\smallskip

Next, Branches~\green{\bf 2c} and~\green{\bf 2d} have less richness,
and one may realize that all the found homogeneous are
products with $\R^1$ of similar homogeneous models
of surfaces $S^2 \subset \R^3$, which is {\em not} true of
Branch~\green{\bf 2b}.

\begin{center}
\includegraphics[scale=0.19]{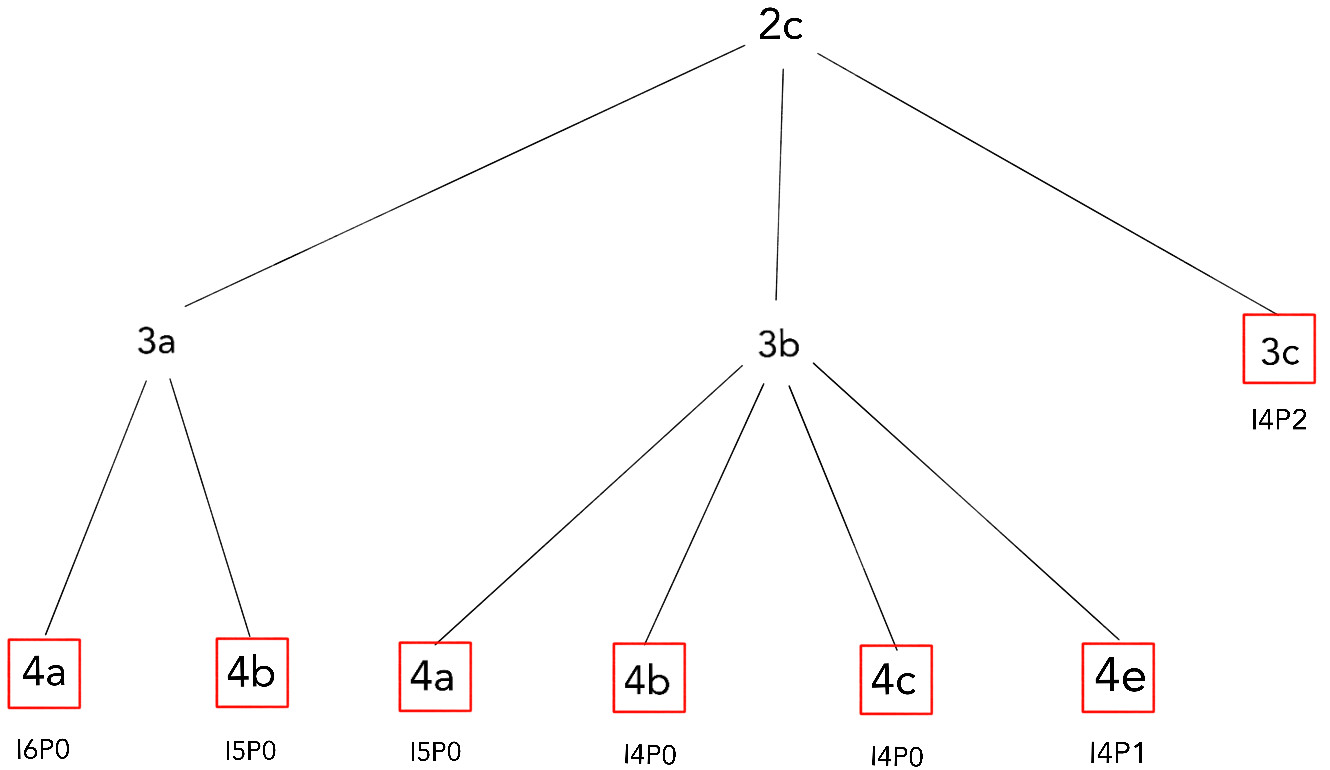}
\end{center}

\begin{center}
\includegraphics[scale=0.19]{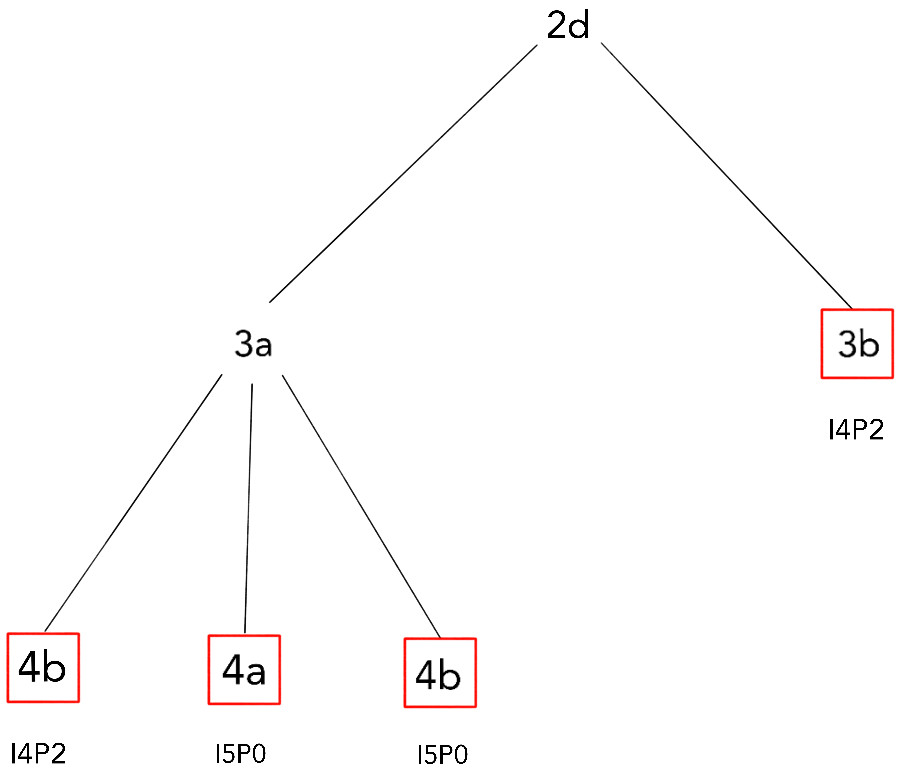}
\end{center}

Incredibly, Branch~\green{\bf 2e}
is surprisingly rich. None of the 27 (families of)
homogeneous models is a product with $\R^1$.

\hspace{-1.5cm}
\includegraphics[scale=0.19]{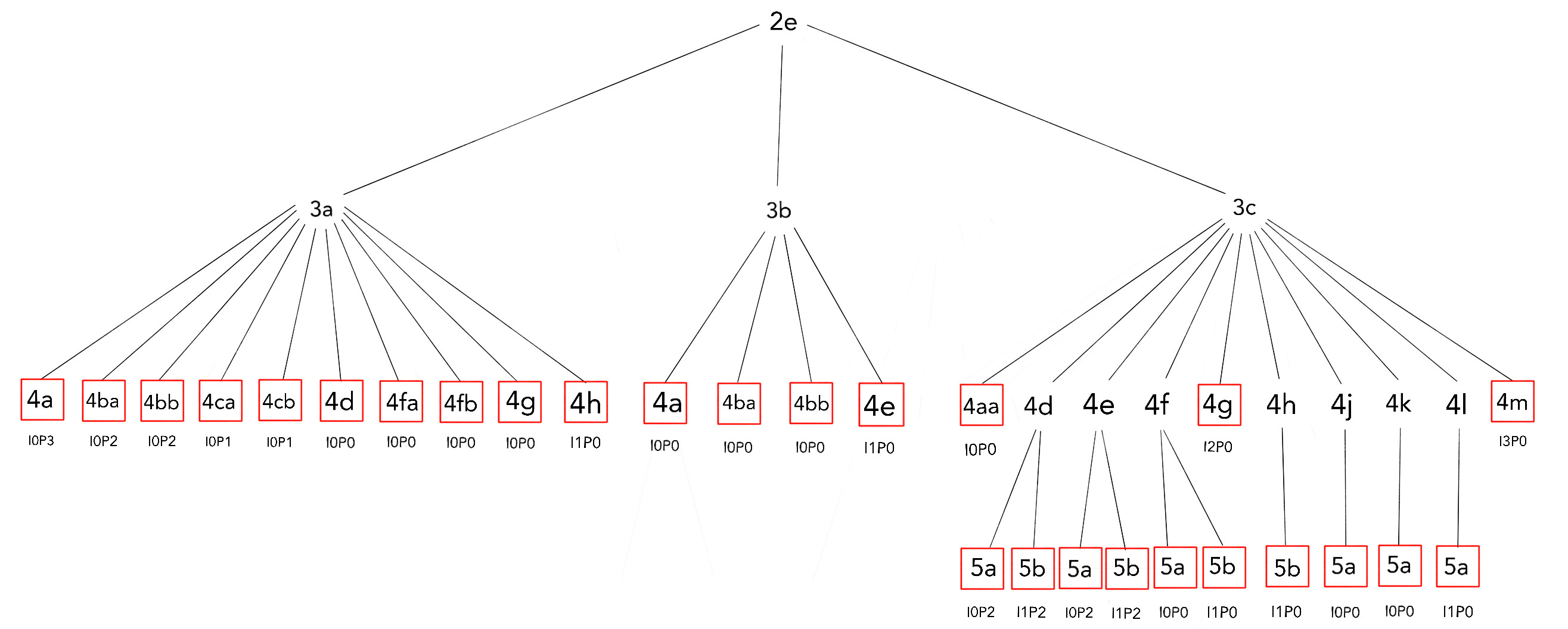}

Similarly, Branch~\green{\bf 2f} is amazing!

\hspace{-1.5cm}
\includegraphics[scale=0.19]{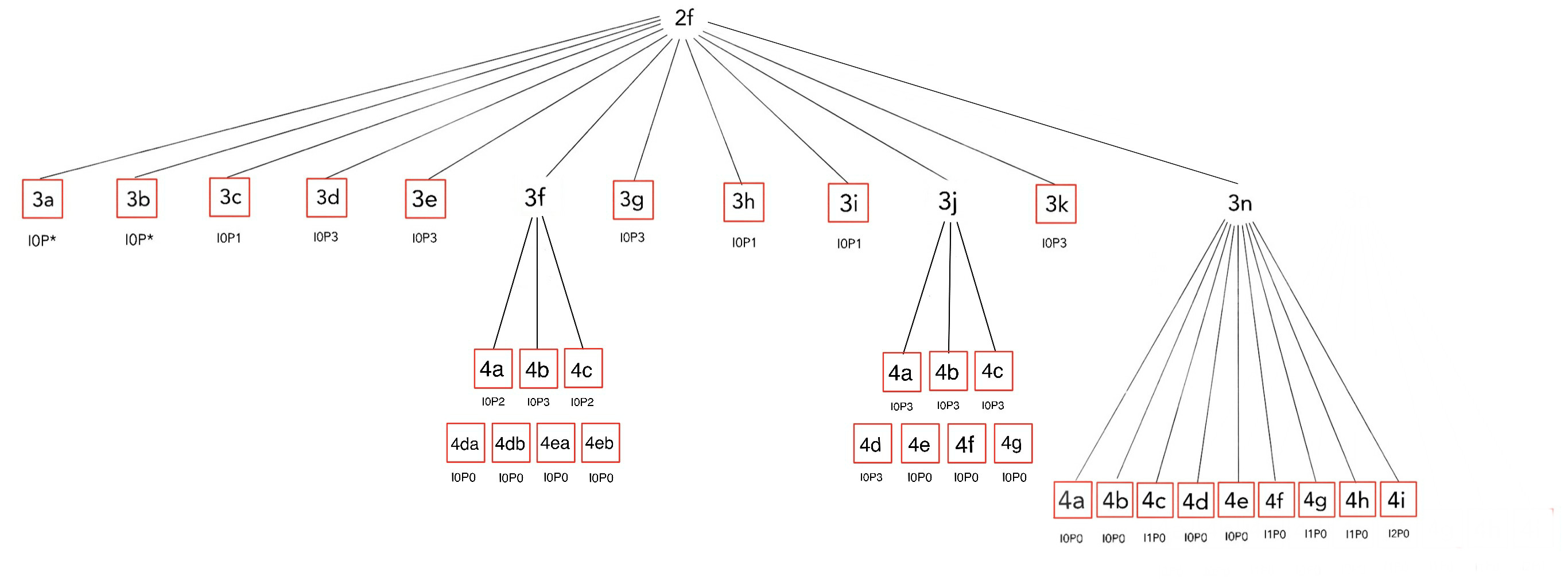}

And lastly, here is Branch~\green{\bf 2g}.

\begin{center}
\includegraphics[scale=0.19]{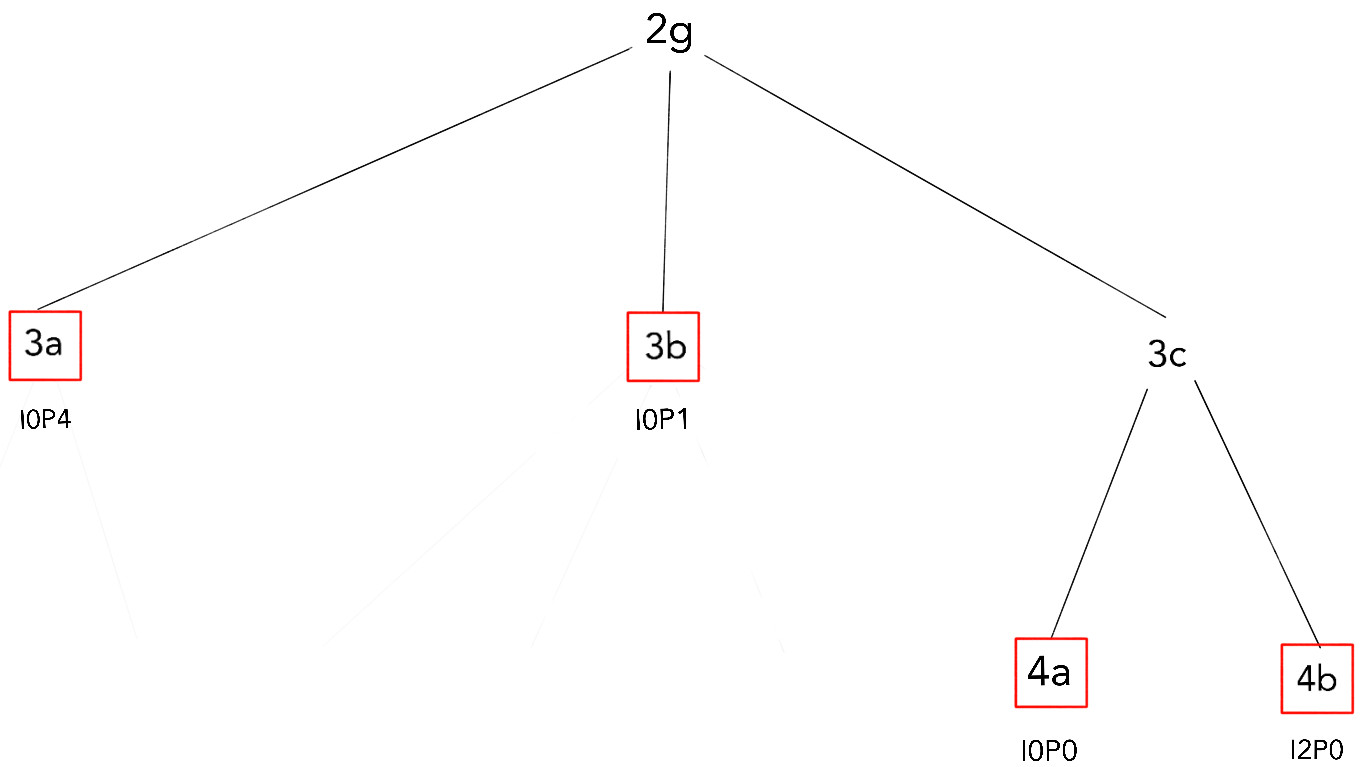}
\end{center}

%%%%%%%%%%%%%%%%%%%%%%%%%%%%%%%%%%%%%%%%%%%%%%%%%%%%%%%%%%%%%%%%%%%%%%
\SectionHead{Branch 2a}
{branch-2a}
%%%%%%%%%%%%%%%%%%%%%%%%%%%%%%%%%%%%%%%%%%%%%%%%%%%%%%%%%%%%%%%%%%%%%%

If an affine infinitesimal transformation $L$ is tangent to:
\[
\aligned
u
&
\,=\,
0
+
{\rm O}_{x,y}(3),
\\
v
&
\,=\,
0
+
{\rm O}_{x,y}(3),
\endaligned
\]
then at order 1:
\[
\aligned
0
&
\underset{\green{\bf F}}{\overset{\green{\bf 1,0}}{\,\,=\,\,}}
-\,C_{1,1},
\\
0
&
\underset{\green{\bf F}}{\overset{\green{\bf 0,1}}{\,\,=\,\,}}
-\,C_{1,2},
\\
0
&
\underset{\green{\bf G}}{\overset{\green{\bf 1,0}}{\,\,=\,\,}}
-\,C_{2,1},
\\
0
&
\underset{\green{\bf G}}{\overset{\green{\bf 0,1}}{\,\,=\,\,}}
-\,C_{2,2}.
\endaligned
\]

Next, after assigning all these $C_{i,j} := 0$,
the tangency of the affine infinitesimal transformation $L$ to:
\[
\aligned
u
&
\,=\,
0
+
F_{3,0}\,x^3
+
F_{2,1}\,x^2y
+
F_{1,2}\,xy^2
+
F_{0,3}\,y^3
+
{\rm O}_{x,y}(4),
\\
v
&
\,=\,
0
+
G_{3,0}\,x^3
+
G_{2,1}\,x^2y
+
G_{1,2}\,xy^2
+
G_{0,3}\,y^3
+
{\rm O}_{x,y}(4),
\endaligned
\]
gives at order 2:
\[
\aligned
0
&
\underset{\green{\bf F}}{\overset{\green{\bf 2,0}}{\,\,=\,\,}}
\big(3\,F_{3,0}\big)\,T_1
+
\big(F_{2,1}\big)\,T_2,
\\
0
&
\underset{\green{\bf F}}{\overset{\green{\bf 1,1}}{\,\,=\,\,}}
\big(2\,F_{2,1}\big)\,T_1
+
\big(2\,F_{1,2}\big)\,T_2,
\\
0
&
\underset{\green{\bf F}}{\overset{\green{\bf 0,2}}{\,\,=\,\,}}
\big(F_{1,2}\big)\,T_1
+
\big(3\,F_{0,3}\big)\,T_2,
\\
0
&
\underset{\green{\bf G}}{\overset{\green{\bf 2,0}}{\,\,=\,\,}}
\big(3\,G_{3,0}\big)\,T_1
+
\big(G_{2,1}\big)\,T_2,
\\
0
&
\underset{\green{\bf G}}{\overset{\green{\bf 1,1}}{\,\,=\,\,}}
\big(2\,G_{2,1}\big)\,T_1
+
\big(2\,G_{1,2}\big)\,T_2,
\\
0
&
\underset{\green{\bf G}}{\overset{\green{\bf 0,2}}{\,\,=\,\,}}
\big(G_{1,2}\big)\,T_1
+
\big(3\,G_{0,3}\big)\,T_2.
\endaligned
\]
Since the transitivity parameters $T_1$, $T_2$ must remain
absolutely free of any linear relations, it necessarily holds that:
\[
\aligned
0
&
\,=\,
F_{3,0}
\,=\,
F_{2,1}
\,=\,
F_{1,2}
\,=\,
F_{0,3},
\\
0
&
\,=\,
G_{3,0}
\,=\,
G_{2,1}
\,=\,
G_{1,2}
\,=\,
G_{0,3}.
\endaligned
\]

Next, the tangency of the affine infinitesimal transformation 
$L$ to:
\[
\aligned
u
&
\,=\,
0
+
F_{4,0}\,x^4
+
F_{3,1}\,x^3y
+
F_{2,2}\,x^2y^2
+
F_{1,3}\,xy^3
+
F_{0,4}\,y^4
+
{\rm O}_{x,y}(5),
\\
v
&
\,=\,
0
+
G_{4,0}\,x^4
+
G_{3,1}\,x^3y
+
G_{2,2}\,x^2y^2
+
G_{1,3}\,xy^3
+
G_{0,4}\,y^4
+
{\rm O}_{x,y}(5),
\endaligned
\]
gives at order 3:
\[
\aligned
0
&
\underset{\green{\bf F}}{\overset{\green{\bf 3,0}}{\,\,=\,\,}}
\big(4\,F_{4,0}\big)\,T_1
+
\big(F_{3,1}\big)\,T_2,
\\
0
&
\underset{\green{\bf F}}{\overset{\green{\bf 2,1}}{\,\,=\,\,}}
\big(3\,F_{3,1}\big)\,T_1
+
\big(2\,F_{2,2}\big)\,T_2,
\\
0
&
\underset{\green{\bf F}}{\overset{\green{\bf 1,2}}{\,\,=\,\,}}
\big(2\,F_{2,2}\big)\,T_1
+
\big(3\,F_{1,3}\big)\,T_2,
\\
0
&
\underset{\green{\bf F}}{\overset{\green{\bf 0,3}}{\,\,=\,\,}}
\big(F_{1,3}\big)\,T_1
+
\big(4\,F_{0,4}\big)\,T_2,
\\
0
&
\underset{\green{\bf G}}{\overset{\green{\bf 3,0}}{\,\,=\,\,}}
\big(4\,G_{4,0}\big)\,T_1
+
\big(G_{3,1}\big)\,T_2,
\\
0
&
\underset{\green{\bf G}}{\overset{\green{\bf 2,1}}{\,\,=\,\,}}
\big(3\,G_{3,1}\big)\,T_1
+
\big(2\,G_{2,2}\big)\,T_2,
\\
0
&
\underset{\green{\bf G}}{\overset{\green{\bf 1,2}}{\,\,=\,\,}}
\big(2\,G_{2,2}\big)\,T_1
+
\big(3\,G_{1,3}\big)\,T_2,
\\
0
&
\underset{\green{\bf G}}{\overset{\green{\bf 0,3}}{\,\,=\,\,}}
\big(G_{1,3}\big)\,T_1
+
\big(4\,G_{0,4}\big)\,T_2,
\endaligned
\]
whence necessarily:
\[
\aligned
0
&
\,=\,
F_{4,0}
\,=\,
F_{3,1}
\,=\,
F_{2,2}
\,=\,
F_{1,3}
\,=\,
F_{0,4},
\\
0
&
\,=\,
G_{4,0}
\,=\,
G_{3,1}
\,=\,
G_{2,2}
\,=\,
G_{1,3}
\,=\,
G_{0,4}.
\endaligned
\]

By induction, all $F_{i,j} = 0 = G_{i,j}$ for any $i + j \geqslant 3$,
and there are $2+12$ free infinitesimal coefficients:
\[
T_1,\ \ \
T_2,\ \ \
A_{1,1},\ \ \
A_{1,2},\ \ \
A_{2,1},\ \ \
A_{2,2},\ \ \
B_{1,1},\ \ \
B_{1,2},\ \ \
B_{2,1},\ \ \
B_{2,2},\ \ \
D_{1,1},\ \ \
D_{1,2},\ \ \
D_{2,1},\ \ \
D_{2,2}.
\] 

\begin{Proposition}
In Branch 2a, there is a single affinely homogeneous model:
\[
\aligned
u
&
\,=\,
0,
\\
v
&
\,=\,
0,
\endaligned
\]
with affine Lie algebra generated by:
\[
\aligned
e_1
&
\,:=\,
\partial_x,
&
\ \ \ \ \ \ \
e_2
&
\,:=\,
\partial_y,
\\
e_3
&
\,:=\,
x\,\partial_x,
&
\ \ \ \ \ \ \
e_4
&
\,:=\,
y\,\partial_x,
&
\ \ \ \ \ \ \
e_5
&
\,:=\,
x\,\partial_y,
&
\ \ \ \ \ \ \
e_6
&
\,:=\,
y\,\partial_y,
\\
e_7
&
\,:=\,
u\,\partial_x,
&
\ \ \ \ \ \ \
e_8
&
\,:=\,
v\,\partial_x,
&
\ \ \ \ \ \ \
e_9
&
\,:=\,
u\,\partial_y,
&
\ \ \ \ \ \ \
e_{10}
&
\,:=\,
v\,\partial_y,
\\
e_{11}
&
\,:=\,
u\,\partial_u,
&
\ \ \ \ \ \ \
e_{12}
&
\,:=\,
v\,\partial_u,
&
\ \ \ \ \ \ \
e_{13}
&
\,:=\,
u\,\partial_v,
&
\ \ \ \ \ \ \
e_{14}
&
\,:=\,
v\,\partial_v,
\endaligned
\]
having Lie-brackets table:
\[
\!\!\!\!\!\!\!\!\!\!\!\!\!\!\!
\footnotesize
\def\arraystretch{1.25}
\begin{array}{c|cccccccccccccc}
{} & e_1 & e_2 & e_3 & e_4 & e_5 & e_6 & e_7 & e_8 & e_9 & e_{10} &
e_{11} & e_{12} & e_{12} & e_{14}
\\
\hline
e_1 &
0 & 0 & e_1 & 0 & e_2 & 0 & 0 & 0 & 0 & 0 & 0 & 0 & 0 & 0
\\
e_2 &
0 & 0 & 0 & e_1 & 0 & e_2 & 0 & 0 & 0 & 0 & 0 & 0 & 0 & 0
\\
e_3 &
-e_1 & 0 & 0 & -e_4 & e_5 & 0 & -e_7 & -e_8 & 0 & 0 & 0 & 0 & 0 & 0
\\
e_4 &
0 & -e_1 & e_4 & 0 & -e_3+e_6 & -e_4 & 0 & 0 & -e_7 & -e_8 & 0 & 0 & 
0 & 0
\\
e_5 & 
-e_2 & 0 & -e_5 & e_3-e_6 & 0 & e_5 & -e_9 & -e_{10} & 0 & 0 & 0 & 
0 & 0 & 0
\\
e_6 &
0 & -e_2 & 0 & e_4 & -e_5 & 0 & 0 & 0 & -e_9 & -e_{10} & 0 & 0 & 0 & 0
\\
e_7 &
0 & 0 & e_7 & 0 & e_9 & 0 & 0 & 0 & 0 & 0 & -e_7 & -e_8 & 0 & 0
\\
e_8 &
0 & 0 & e_8 & 0 & e_{10} & 0 & 0 & 0 & 0 & 0 & 0 & 0 & -e_7 & -e_8
\\
e_9 &
0 & 0 & 0 & e_7 & 0 & e_9 & 0 & 0 & 0 & 0 & -e_9 & -e_{10} & 0 & 0
\\
e_{10} &
0 & 0 & 0 & e_8 & 0 & e_{10} & 0 & 0 & 0 & 0 & 0 & 0 & -e_9 & -e_{10}
\\
e_{11} &
0 & 0 & 0 & 0 & 0 & 0 & e_7 & 0 & e_9 & 0 & 0 & -e_{12} & e_{13} & 0
\\
e_{12} &
0 & 0 & 0 & 0 & 0 & 0 & e_8 & 0 & e_{10} & 0 & e_{12} & 0 & 
-e_{11}+e_{14} & -e_{12}
\\
e_{13} &
0 & 0 & 0 & 0 & 0 & 0 & 0 & e_7 & 0 & e_9 & -e_{13} & e_{11}-e_{14} &
0 & e_{13}
\\
e_{14} &
0 & 0 & 0 & 0 & 0 & 0 & 0 & e_8 & 0 & e_{10} & 0 & e_{12} &
-e_{13} & 0
\end{array}
\]
\end{Proposition}

%%%%%%%%%%%%%%%%%%%%%%%%%%%%%%%%%%%%%%%%%%%%%%%%%%%%%%%%%%%%%%%%%%%%%%
\SectionHead{Linear Representations Branch by Branch}
{linear-representations-branches}
%%%%%%%%%%%%%%%%%%%%%%%%%%%%%%%%%%%%%%%%%%%%%%%%%%%%%%%%%%%%%%%%%%%%%%

Here, after the (easy) Branch \green{\bf 2a},
we present all the other different branches, 
with their related linear representation nodes.

\medskip\noindent
{\bf Branch 2b.}
At order 3, 
the associated linear representation is:
\[
\left(\!
\begin{array}{c}
R_{2,1}
\\
S_{3,0}
\end{array}
\!\right)
\,=\,
\left(\!
\def\arraystretch{1.25}
\begin{array}{cc}
\tfrac{1}{a_{2,2}} & 0
\\
0 & \tfrac{d_{2,2}}{a_{1,1}^3}
\end{array}
\!\right)\,
\left(\!
\begin{array}{c}
F_{2,1}
\\
G_{3,0}
\end{array}
\!\right),
\]
Skipping details,
this leads us to the creation of $4$ branches,
mutually inequivalent and having empty intersection: 
\[
\def\arraystretch{1.25}
\begin{array}{rccc}
\green{\bf 2b}\,\,\,\,
\green{\downarrow}\,\,
& F_{2,1} & G_{3,0}
\\
\green{\bf 3a} & 
0 & 0
\\
\green{\bf 3b} & 
1 & 0
\\
\green{\bf 3c} & 
0 & 1 
\\
\green{\bf 3d} &
1 & 1
\end{array}
\]

\medskip

In branch \green{\bf 2b3a}:
\[
R_{4,0}
=
\frac{F_{4,0}}{a_{1,1}^2},
\]
leading to:
\[
\def\arraystretch{1.25}
\begin{array}{rcc}
\green{\bf 2b3a}\,\,\,\,
\green{\downarrow}\,\,
& F_{4,0}
\\
\green{\bf 4a} & 
0
\\
\green{\bf 4b} & 
1
\\
\green{\bf 4c} & 
-1 
\end{array}
\]

\medskip

In branch \green{\bf 2b3b}:
\[
R_{3,1}
=
\frac{F_{3,1}}{a_{1,1}},
\]
leading to:
\[
\def\arraystretch{1.25}
\begin{array}{rcc}
\green{\bf 2b3b}\,\,\,\,
\green{\downarrow}\,\,
& F_{3,1}
\\
\green{\bf 4a} & 
0
\\
\green{\bf 4b} & 
1
\end{array}
\]

\medskip

In branch \green{\bf 2b3b4a}:
\[
R_{5,0}
=
\frac{F_{5,0}}{a_{1,1}^3},
\]
leading to:
\[
\def\arraystretch{1.25}
\begin{array}{rcc}
\green{\bf 2b3b4a}\,\,\,\,
\green{\downarrow}\,\,
& F_{5,0}
\\
\green{\bf 5a} & 
0
\\
\green{\bf 5b} & 
1
\end{array}
\]

\medskip

In branch \green{\bf 2b3c}:
\[
\left(\!
\begin{array}{c}
R_{5,0}
\\
S_{5,0}
\end{array}
\!\right)
\,=\,
\left(\!
\def\arraystretch{1.25}
\begin{array}{cc}
\tfrac{1}{a_{1,1}^3} & 0
\\
0 & \tfrac{1}{a_{1,1}^2}
\end{array}
\!\right)\,
\left(\!
\begin{array}{c}
F_{5,0}
\\
G_{5,0}
\end{array}
\!\right),
\]
leading to:
\[
\def\arraystretch{1.25}
\begin{array}{rccc}
\green{\bf 2b3c}\,\,\,\,
\green{\downarrow}\,\,
& F_{5,0} & G_{5,0}
\\
\green{\bf 5a} & 
0 & 0
\\
\green{\bf 5b} & 
0 & 1
\\
\green{\bf 5c} & 
0 & -1
\\
\green{\bf 5d} & 
1 & G_{5,0}
\end{array}
\]

\medskip

In branch \green{\bf 2b3d}:
\[
R_{3,1}
=
\frac{F_{3,1}}{a_{1,1}},
\]
leading to:
\[
\def\arraystretch{1.25}
\begin{array}{rcc}
\green{\bf 2b3d}\,\,\,\,
\green{\downarrow}\,\,
& F_{3,1}
\\
\green{\bf 4a} & 
1
\\
\green{\bf 4b} & 
0
\end{array}
\]

\medskip

In branch \green{\bf 2b3d4b}:
\[
\left(\!
\begin{array}{c}
R_{6,0}
\\
S_{6,0}
\end{array}
\!\right)
\,=\,
\left(\!
\def\arraystretch{1.25}
\begin{array}{cc}
\tfrac{1}{a_{1,1}^4} & -\tfrac{a_{2,1}}{3a_{1,1}^4}
\\
0 & \tfrac{1}{a_{1,1}^3}
\end{array}
\!\right)\,
\left(\!
\begin{array}{c}
F_{6,0}
\\
G_{6,0}
\end{array}
\!\right),
\]
leading to:
\[
\def\arraystretch{1.25}
\begin{array}{rccc}
\green{\bf 2b3d4b}\,\,\,\,
\green{\downarrow}\,\,
& F_{6,0} & G_{6,0}
\\
\green{\bf 6a} & 
0 & 0
\\
\green{\bf 6b} & 
1 & 0
\\
\green{\bf 6c} & 
-1 & 0
\\
\green{\bf 6d} & 
0 & 1
\end{array}
\]

\medskip\noindent
{\bf Branch 2c.}
At order 3,
the associated linear representation is:
\[
\left(\!
\begin{array}{c}
R_{0,3}
\\
R_{3,0}
\end{array}
\!\right)
\,=\,
\left(\!
\def\arraystretch{1.25}
\begin{array}{cc}
\tfrac{a_{1,1}}{a_{2,2}^2} & 0
\\
0 & \tfrac{a_{2,2}}{a_{1,1}^2}
\end{array}
\!\right)\,
\left(\!
\begin{array}{c}
F_{0,3}
\\
F_{3,0}
\end{array}
\!\right),
\]
leading to:
\[
\def\arraystretch{1.25}
\begin{array}{rccc}
\green{\bf 2c}\,\,\,\,
\green{\downarrow}\,\,
& F_{0,3} & F_{3,0}
\\
\green{\bf 3a} & 
0 & 0
\\
\green{\bf 3b} & 
1 & 0
\\
\green{\bf 3c} & 
1 & 1
\end{array}
\]

\medskip

In branch \green{\bf 2c3a}:
\[
R_{2,2}
=
\frac{F_{2,2}}{a_{1,1}a_{2,2}},
\]
leading to:
\[
\def\arraystretch{1.25}
\begin{array}{rcc}
\green{\bf 2c3a}\,\,\,\,
\green{\downarrow}\,\,
& F_{2,2}
\\
\green{\bf 4a} & 
0
\\
\green{\bf 4b} & 
1
\end{array}
\]

\medskip

In branch \green{\bf 2c3b}:
\[
\left(\!
\begin{array}{c}
R_{0,4}
\\
R_{1,3}
\\
R_{2,2}
\end{array}
\!\right)
\,=\,
\left(\!
\def\arraystretch{1.25}
\begin{array}{ccc}
\tfrac{1}{a_{2,2}} & 0 & 0
\\
0 & \tfrac{1}{a_{2,2}^2} & 0
\\
0 & 0 & \tfrac{1}{a_{2,2}^3}
\end{array}
\!\right)\,
\left(\!
\begin{array}{c}
F_{0,4}
\\
F_{1,3}
\\
F_{2,2}
\end{array}
\!\right),
\]
leading to:
\[
\def\arraystretch{1.25}
\begin{array}{rcccc}
\green{\bf 2c3b}\,\,\,\,
\green{\downarrow}\,\,
& F_{0,4} & F_{1,3} & F_{2,2}
\\
\green{\bf 4a} & 
0 & 0 & 0
\\
\green{\bf 4b} & 
0 & 1 & 0
\\
\green{\bf 4c} & 
0 & -1 & 0
\\
\green{\bf 4d} & 
0 & F_{1,3} & 1
\\
\green{\bf 4e} & 
1 & F_{1,3} & F_{2,2}
\end{array}
\]

\medskip\noindent
{\bf Branch 2d.}
At order 3,
the associated linear representation is:
\[
\left(\!
\begin{array}{c}
R_{0,3}
\\
R_{3,0}
\end{array}
\!\right)
\,=\,
\left(\!
\def\arraystretch{1.25}
\begin{array}{cc}
\tfrac{a_{1,1}(a_{1,1}^2-3a_{1,2}^2)}{(a_{1,1}^2+a_{1,2}^2)^2} & \tfrac{a_{1,2}(3a_{1,1}^2-a_{1,2}^2)}{(a_{1,1}^2+a_{1,2}^2)^2}
\\
-\tfrac{a_{1,2}(3a_{1,1}^2-a_{1,2}^2)}{(a_{1,1}^2+a_{1,2}^2)^2} & \tfrac{a_{1,1}(a_{1,1}^2-3a_{1,2}^2)}{(a_{1,1}^2+a_{1,2}^2)^2}
\end{array}
\!\right)\,
\left(\!
\begin{array}{c}
F_{0,3}
\\
F_{3,0}
\end{array}
\!\right),
\]
leading to:
\[
\def\arraystretch{1.25}
\begin{array}{rccc}
\green{\bf 2d}\,\,\,\,
\green{\downarrow}\,\,
& F_{0,3} & F_{3,0}
\\
\green{\bf 3a} & 
0 & 0
\\
\green{\bf 3b} & 
0 & 1
\end{array}
\]

\medskip

In branch \green{\bf 2d3a}:
\[
R_{2,2}
=
\frac{F_{2,2}}{a_{1,1}^2+a_{1,2}^2},
\]
leading to:
\[
\def\arraystretch{1.25}
\begin{array}{rcc}
\green{\bf 2d3a}\,\,\,\,
\green{\downarrow}\,\,
& F_{2,2}
\\
\green{\bf 4a} & 
0
\\
\green{\bf 4b} & 
1
\\
\green{\bf 4c} & 
-1
\end{array}
\]

\medskip\noindent
{\bf Branch 2e.}
At order 3, 
the associated linear representation is:
\[
\left(\!
\begin{array}{c}
R_{1,2}
\\
R_{0,3}
\end{array}
\!\right)
\,=\,
\left(\!
\def\arraystretch{1.25}
\begin{array}{cc}
\tfrac{1}{a_{2,2}} & -\tfrac{3a_{2,1}}{a_{2,2}}
\\
0 & \tfrac{a_{1,1}}{a_{2,2}^2}
\end{array}
\!\right)\,
\left(\!
\begin{array}{c}
F_{1,2}
\\
F_{0,3}
\end{array}
\!\right),
\]
leading to:
\[
\def\arraystretch{1.25}
\begin{array}{rccc}
\green{\bf 2e}\,\,\,\,
\green{\downarrow}\,\,
& F_{1,2} & F_{0,3}
\\
\green{\bf 3a} & 
0 & 1
\\
\green{\bf 3b} & 
1 & 0
\\
\green{\bf 3c} & 
0 & 0
\end{array}
\]

\medskip

In branch \green{\bf 2e3a}:
\[
\left(\!
\begin{array}{c}
R_{4,0}
\\
S_{4,0}
\\
R_{3,1}
\\
S_{3,1}
\\
R_{2,2}
\\
R_{1,3}
\\
R_{0,4}
\end{array}
\!\right)
\,=\,
\left(\!
\def\arraystretch{1.25}
\begin{array}{ccccccc}
\tfrac{1}{a_{2,2}^5} & 0 & 0 & 0 & 0 & 0 & 0
\\
0 & \tfrac{1}{a_{2,2}^4} & 0 & 0 & 0 & 0 & 0
\\
0 & 0 & \tfrac{1}{a_{2,2}^4} & 0 & 0 & 0 & 0
\\
0 & 0 & 0 & \tfrac{1}{a_{2,2}^3} & 0 & 0 & 0
\\
0 & 0 & 0 & 0 & \tfrac{1}{a_{2,2}^3} & 0 & 0
\\
0 & 0 & 0 & 0 & 0 & \tfrac{1}{a_{2,2}^2} & 0
\\
0 & 0 & 0 & 0 & 0 & 0 & \tfrac{1}{a_{2,2}}
\end{array}
\!\right)\,
\left(\!
\begin{array}{c}
F_{4,0}
\\
G_{4,0}
\\
F_{3,1}
\\
G_{3,1}
\\
F_{2,2}
\\
F_{1,3}
\\
F_{0,4}
\end{array}
\!\right),
\]
leading to:
\[
\def\arraystretch{1.25}
\begin{array}{rccccccc}
\green{\bf 2e3a}\,\,\,\,
\green{\downarrow}\,\,
& F_{4,0} & G_{4,0} & F_{3,1} & G_{3,1} & F_{2,2} & F_{1,3} & F_{0,4}
\\
\green{\bf 4a} & 
1 & G_{4,0} & F_{3,1} & G_{3,1} & F_{2,2} & F_{1,3} & F_{0,4}
\\
\green{\bf 4ba} & 
0 & 1 & F_{3,1} & G_{3,1} & F_{2,2} & F_{1,3} & F_{0,4}
\\
\green{\bf 4bb} & 
0 & -1 & F_{3,1} & G_{3,1} & F_{2,2} & F_{1,3} & F_{0,4}
\\
\green{\bf 4ca} & 
0 & 0 & 1 & G_{3,1} & F_{2,2} & F_{1,3} & F_{0,4}
\\
\green{\bf 4cb} & 
0 & 0 & -1 & G_{3,1} & F_{2,2} & F_{1,3} & F_{0,4}
\\
\green{\bf 4d} & 
0 & 0 & 0 & 1 & F_{2,2} & F_{1,3} & F_{0,4}
\\
\green{\bf 4e} & 
0 & 0 & 0 & 0 & 1 & F_{1,3} & F_{0,4}
\\
\green{\bf 4fa} & 
0 & 0 & 0 & 0 & 0 & 1 & F_{0,4}
\\
\green{\bf 4fb} & 
0 & 0 & 0 & 0 & 0 & -1 & F_{0,4}
\\
\green{\bf 4g} & 
0 & 0 & 0 & 0 & 0 & 0 & 1
\\
\green{\bf 4h} & 
0 & 0 & 0 & 0 & 0 & 0 & 0
\\
\end{array}
\]

\medskip

In branch \green{\bf 2e3b}:
\[
\left(\!
\begin{array}{c}
R_{2,2}
\\
R_{3,1}
\\
S_{4,0}
\\
R_{4,0}
\end{array}
\!\right)
\,=\,
\left(\!
\def\arraystretch{1.25}
\begin{array}{cccc}
\tfrac{1}{a_{1,1}} & 0 & 0 & 0
\\
0 & \tfrac{1}{a_{1,1}^2} & 0 & 0
\\
0 & 0 & \tfrac{1}{a_{1,1}^2} & 0
\\
0 & 0 & 0 & \tfrac{1}{a_{1,1}^3}
\end{array}
\!\right)\,
\left(\!
\begin{array}{c}
F_{2,2}
\\
F_{1,3}
\\
G_{4,0}
\\
F_{4,0}
\end{array}
\!\right),
\]
leading to:
\[
\def\arraystretch{1.25}
\begin{array}{rcccc}
\green{\bf 2e3b}\,\,\,\,
\green{\downarrow}\,\,
& F_{2,2} & F_{3,1} & G_{4,0} & F_{4,0}
\\
\green{\bf 4a} & 
1 & F_{3,1} & G_{4,0} & F_{4,0}
\\
\green{\bf 4ba} & 
0 & 1 & G_{4,0} & F_{4,0}
\\
\green{\bf 4bb} & 
0 & -1 & G_{4,0} & F_{4,0}
\\
\green{\bf 4ca} & 
0 & 0 & 1 & F_{4,0}
\\
\green{\bf 4cb} & 
0 & 0 & -1 & F_{4,0}
\\
\green{\bf 4d} & 
0 & 0 & 0 & 1
\\
\green{\bf 4e} & 
0 & 0 & 0 & 0
\end{array}
\]

\medskip

In branch \green{\bf 2e3c}:
\[
\left(\!
\begin{array}{c}
S_{3,1}
\\
R_{3,1}
\\
S_{4,0}
\\
R_{4,0}
\end{array}
\!\right)
\,=\,
\left(\!
\def\arraystretch{1.25}
\begin{array}{cccc}
\tfrac{1}{a_{1,1}a_{2,2}} & 0 & 0 & 0
\\
-\tfrac{a_{2,1}}{2a_{1,1}^2a_{2,2}} & \tfrac{1}{a_{1,1}^2} & 0 & 0
\\
-\tfrac{a_{2,1}}{a_{1,1}^2a_{2,2}} & 0 & \tfrac{1}{a_{1,1}^2} & 0
\\
-\tfrac{a_{2,1}^2}{4a_{1,1}^3a_{2,2}} & -\tfrac{a_{2,1}}{a_{1,1}^3} & \tfrac{a_{2,1}}{a_{1,1}^3} & \tfrac{1}{a_{1,1}^3}
\end{array}
\!\right)\,
\left(\!
\begin{array}{c}
G_{3,1}
\\
F_{1,3}
\\
G_{4,0}
\\
F_{4,0}
\end{array}
\!\right),
\]
leading to:
\[
\def\arraystretch{1.25}
\begin{array}{rcccc}
\green{\bf 2e3c}\,\,\,\,
\green{\downarrow}\,\,
& G_{3,1} & F_{3,1} & G_{4,0} & F_{4,0}
\\
\green{\bf 4aa} & 
1 & 0 & 1 & F_{4,0}
\\
\green{\bf 4ab} & 
1 & 0 & -1 & F_{4,0}
\\
\green{\bf 4ba} & 
1 & 0 & 0 & 1
\\
\green{\bf 4bb} & 
1 & 0 & 0 & -1
\\
\green{\bf 4c} & 
1 & 0 & 0 & 0
\\
\green{\bf 4d} & 
0 & 1 & G_{4,0}^{\neq 1} & 0
\\
\green{\bf 4e} & 
0 & -1 & G_{4,0}^{\neq -1} & 0
\\
\green{\bf 4f} & 
0 & 1 & 1 & 1
\\
\green{\bf 4g} & 
0 & 1 & 1 & 0
\\
\green{\bf 4h} & 
0 & -1 & -1 & 1
\\
\green{\bf 4i} & 
0 & -1 & -1 & 0
\\
\green{\bf 4j} & 
0 & 0 & 1 & 0
\\
\green{\bf 4k} & 
0 & 0 & -1 & 0
\\
\green{\bf 4l} & 
0 & 0 & 0 & 1
\\
\green{\bf 4m} & 
0 & 0 & 0 & 0
\end{array}
\]

\medskip

In branch \green{\bf 2e3c4d}:
\[
R_{5,0}
=
a_{2,2}F_{5,0},
\]
leading to:
\[
\def\arraystretch{1.25}
\begin{array}{rcc}
\green{\bf 2e3c4d}\,\,\,\,
\green{\downarrow}\,\,
& F_{5,0}
\\
\green{\bf 5a} & 
1
\\
\green{\bf 5b} & 
0
\end{array}
\]

\medskip

In branch \green{\bf 2e3c4e}:
\[
R_{5,0}
=
a_{2,2}F_{5,0},
\]
leading to:
\[
\def\arraystretch{1.25}
\begin{array}{rcc}
\green{\bf 2e3c4e}\,\,\,\,
\green{\downarrow}\,\,
& F_{5,0}
\\
\green{\bf 5a} & 
1
\\
\green{\bf 5b} & 
0
\end{array}
\]

\medskip

In branch \green{\bf 2e3c4f}:
\[
R_{5,0}
-
F_{5,0}
+
\frac{1}{4}a_{2,1}G_{5,0}
=
0,
\]
leading to:
\[
\def\arraystretch{1.25}
\begin{array}{rcc}
\green{\bf 2e3c4f}\,\,\,\,
\green{\downarrow}\,\,
& G_{5,0}
\\
\green{\bf 5a} & 
G_{5,0}^{\neq 0}
\\
\green{\bf 5b} & 
0
\end{array}
\]

\medskip

In branch \green{\bf 2e3c4h}:
\[
R_{5,0}
-
F_{5,0}
+
\frac{1}{5}a_{2,1}F_{4,1}
=
0,
\]
leading to:
\[
\def\arraystretch{1.25}
\begin{array}{rcc}
\green{\bf 2e3c4h}\,\,\,\,
\green{\downarrow}\,\,
& F_{4,1}
\\
\green{\bf 5a} & 
1
\\
\green{\bf 5b} & 
0
\end{array}
\]

\medskip

In branch \green{\bf 2e3c4j}:
\[
R_{5,0}
=
a_{2,2}F_{5,0},
\]
leading to:
\[
\def\arraystretch{1.25}
\begin{array}{rcc}
\green{\bf 2e3c4j}\,\,\,\,
\green{\downarrow}\,\,
& F_{5,0}
\\
\green{\bf 5a} & 
1
\\
\green{\bf 5b} & 
0
\end{array}
\]

\medskip

In branch \green{\bf 2e3c4k}:
\[
R_{5,0}
=
a_{2,2}F_{5,0},
\]
leading to:
\[
\def\arraystretch{1.25}
\begin{array}{rcc}
\green{\bf 2e3c4k}\,\,\,\,
\green{\downarrow}\,\,
& F_{5,0}
\\
\green{\bf 5a} & 
1
\\
\green{\bf 5b} & 
0
\end{array}
\]

\medskip

In branch \green{\bf 2e3c4l}:
\[
R_{5,0}
=
a_{2,2}F_{5,0},
\]
leading to:
\[
\def\arraystretch{1.25}
\begin{array}{rcc}
\green{\bf 2e3c4l}\,\,\,\,
\green{\downarrow}\,\,
& F_{5,0}
\\
\green{\bf 5a} & 
1
\\
\green{\bf 5b} & 
0
\end{array}
\]

\medskip\noindent
{\bf Branch 2f.}
At order 3,
the associated linear representation is:
\[
\left(\!
\begin{array}{c}
R_{3,0}
\\
R_{0,3}
\\
S_{2,1}
\\
S_{1,2}
\end{array}
\!\right)
\,=\,
\tfrac{1}{(a_{1,1}^2-a_{2,1}^2)^3}
\,
\bf A
\,
\left(\!
\begin{array}{c}
F_{3,0}
\\
F_{0,3}
\\
G_{2,1}
\\
G_{1,2}
\end{array}
\!\right),
\]
where $\bf A$ is the following matrix:
\[
\!\!\!\!\!\!\!\!\!\!\!\!\!\!\!\!\!\!
\footnotesize
\aligned
\left(\!
\def\arraystretch{1.25}
\scriptsize
\begin{array}{cccc}
a_{1,1}(a_{1,1}^4+a_{1,1}^2a_{2,1}^2+2a_{2,1}^4) & -a_{2,1}(2a_{1,1}^4+a_{1,1}^2a_{2,1}^2) & -\tfrac{1}{2}a_{1,1}a_{2,1}^2(3a_{1,1}^2+a_{2,1}^2) & -\tfrac{1}{2}a_{1,1}^2a_{2,1}(a_{1,1}^2+3a_{2,1}^2)
\\
-a_{2,1}(2a_{1,1}^4+a_{1,1}a_{2,1}^2+a_{2,1}^4) & a_{1,1}(a_{1,1}^4+a_{1,1}^2a_{2,1}^2+2a_{2,1}^4) & \tfrac{1}{2}a_{1,1}^2a_{2,1}(a_{1,1}^2+3a_{2,1}^2) & -\tfrac{1}{2}a_{1,1}a_{2,1}^2(3a_{1,1}+a_{2,1}^2)
\\
-6a_{1,1}a_{2,1}^2(3a_{1,1}^2+a_{2,1}^2) & 6a_{1,1}^2a_{2,1}(a_{1,1}^2+3a_{2,1}^2) & a_{1,1}(a_{1,1}^4+7a_{1,1}^2a_{2,1}^2+4a_{2,1}^4) & -a_{2,1}(4a_{1,1}^4+7a_{1,1}^2a_{2,1}^2+a_{2,1}^4)
\\
6a_{1,1}^2a_{2,1}(a_{1,1}^2+3a_{2,1}^2) & -6a_{1,1}a_{2,1}^2(3a_{1,1}^2+a_{2,1}^2) & -a_{2,1}(4a_{1,1}^4+7a_{1,1}^2a_{2,1}^2+a_{2,1}^4) & a_{1,1}(a_{1,1}^4+7a_{1,1}^2a_{2,1}^2+4a_{2,1}^4)
\end{array}
\!\right),
\endaligned
\]
leading to:
\[
\def\arraystretch{1.25}
\begin{array}{rcccccc}
\green{\bf 2f}\,\,\,\,
\green{\downarrow}\,\,
& F_{3,0} & F_{0,3} & G_{2,1} & G_{1,2}
\\
\green{\bf 3a} & 
4-F_{0,3} & F_{0,3} & G_{2,1} & -8-G_{2,1}
\\
\green{\bf 3b} & 
-F_{0,3} & F_{0,3} & G_{2,1} & 16-G_{2,1} 
\\
\green{\bf 3c} & 
2-F_{0,3} & F_{0,3} & 4+6F_{0,3} & -6F_{0,3}
\\
\green{\bf 3d} & 
2 & 0 & -4 & 8
\\
\green{\bf 3e} & 
0 & 2 & 8 & -4
\\
\green{\bf 3f} & 
1 & 1 & 2 & 2
\\
\green{\bf 3g} & 
2-F_{0,3} & F_{0,3} & -4+6F_{0,3} & -8-6F_{0,3} 
\\
\green{\bf 3h} & 
2 & 0 & -12 & 0
\\
\green{\bf 3i} & 
0 & 2 & 0 & -12
\\
\green{\bf 3j} & 
1 & 1 & -6 & -6
\\
\green{\bf 3k} & 
2 & -2 & -4 & 4
\\
\green{\bf 3l} & 
1 & -1 & -6 & 6
\\
\green{\bf 3m} & 
1 & -1 & 2 & -2
\\
\green{\bf 3n} & 
0 & 0 & 0 & 0
\end{array}
\]
Through $x \longleftrightarrow -\,x$, 
branch~\green{\bf 2f3l} is equivalent to branch~\green{\bf 2f3j}
and
branch~\green{\bf 2f3m} is equivalent to branch~\green{\bf 2f3f}.

The complexity of the above linear representation matrix
${\bf A}$ makes it difficult to express group reductions
in the two branches~\green{\bf 2f3f}
and~\green{\bf 2f3j}. So to create subbranches,
{\em cf.} the diagram of branch~\green{\bf 2f},
we proceed by using only
the infinitesimal fundamental 
equations $\eqLF (x,y)$ and $\eqLG (x,y)$.
Details can be written.
 
\medskip

In branch \green{\bf 2f3n}:
\[
\left(\!
\begin{array}{c}
R_{2,2}
\\
S_{0,4}
\end{array}
\!\right)
\,=\,
\left(\!
\def\arraystretch{1.25}
\begin{array}{cc}
\tfrac{a_{1,1}^2+a_{2,1}^2}{(a_{1,1}^2-a_{2,1}^2)^2} & -\tfrac{6a_{1,1}a_{2,1}}{(a_{1,1}^2-a_{2,1}^2)^2}
\\
-\tfrac{2}{3}\tfrac{a_{1,1}a_{2,1}}{(a_{1,1}^2-a_{2,1}^2)^2} & \tfrac{a_{1,1}^2+a_{2,1}^2}{(a_{1,1}^2-a_{2,1}^2)^2}
\end{array}
\!\right)\,
\left(\!
\begin{array}{c}
F_{2,2}
\\
G_{0,4}
\end{array}
\!\right),
\]
leading to:
\[
\def\arraystretch{1.25}
\begin{array}{rccc}
\green{\bf 2f3n}\,\,\,\,
\green{\downarrow}\,\,
& F_{2,2} & G_{0,4}
\\
\green{\bf 4a} & 
0 & 2
\\
\green{\bf 4b} & 
6 & 0
\\
\green{\bf 4c} & 
3 & 1
\\
\green{\bf 4d} & 
-6 & 0
\\
\green{\bf 4e} & 
0 & -2
\\
\green{\bf 4f} & 
-3 & -1
\\
\green{\bf 4g} & 
-3 & 1
\\
\green{\bf 4h} & 
3 & -1
\\
\green{\bf 4i} & 
0 & 0
\end{array}
\]

\medskip\noindent
{\bf Branch 2g.}
At order 3,
the associated linear representation is:
\[
\left(\!
\begin{array}{c}
R_{3,0}
\\
R_{0,3}
\\
S_{2,1}
\\
S_{1,2}
\end{array}
\!\right)
\,=\,
\tfrac{1}{(a_{1,1}^2+a_{2,1}^2)^3}
\,
\bf A
\,
\left(\!
\begin{array}{c}
F_{3,0}
\\
F_{0,3}
\\
G_{2,1}
\\
G_{1,2}
\end{array}
\!\right),
\]
where $\bf A$ is the following matrix:
\[
\!\!\!\!\!\!\!\!\!\!\!\!\!\!\!\!\!\!\!\!
\footnotesize
\aligned
\left(\!
\def\arraystretch{1.25}
\scriptsize
\begin{array}{cccc}
a_{1,1}(a_{1,1}^4-a_{1,1}^2a_{2,1}^2+2a_{2,1}^4) & -a_{2,1}(2a_{1,1}^4-a_{1,1}^2a_{2,1}^2) & -\tfrac{1}{2}a_{1,1}a_{2,1}^2(3a_{1,1}^2-a_{2,1}^2) & -\tfrac{1}{2}a_{1,1}^2a_{2,1}(a_{1,1}^2-3a_{2,1}^2)
\\
-a_{2,1}(2a_{1,1}^4-a_{1,1}a_{2,1}^2+a_{2,1}^4) & a_{1,1}(a_{1,1}^4-a_{1,1}^2a_{2,1}^2+2a_{2,1}^4) & -\tfrac{1}{2}a_{1,1}^2a_{2,1}(a_{1,1}^2-3a_{2,1}^2) & \tfrac{1}{2}a_{1,1}a_{2,1}^2(3a_{1,1}-a_{2,1}^2)
\\
-6a_{1,1}a_{2,1}^2(3a_{1,1}^2-a_{2,1}^2) & 6a_{1,1}^2a_{2,1}(a_{1,1}^2-3a_{2,1}^2) & a_{1,1}(a_{1,1}^4-7a_{1,1}^2a_{2,1}^2+4a_{2,1}^4) & -a_{2,1}(4a_{1,1}^4-7a_{1,1}^2a_{2,1}^2+a_{2,1}^4)
\\
6a_{1,1}^2a_{2,1}(a_{1,1}^2-3a_{2,1}^2) & 6a_{1,1}a_{2,1}^2(3a_{1,1}^2-a_{2,1}^2) & a_{2,1}(4a_{1,1}^4-7a_{1,1}^2a_{2,1}^2+a_{2,1}^4) & a_{1,1}(a_{1,1}^4-7a_{1,1}^2a_{2,1}^2+4a_{2,1}^4)
\end{array}
\!\right),
\endaligned
\]
leading to:
\[
\def\arraystretch{1.25}
\begin{array}{rcccccc}
\green{\bf 2g}\,\,\,\,
\green{\downarrow}\,\,
& F_{3,0} & F_{0,3} & G_{2,1} & G_{1,2}
\\
\green{\bf 3a} & 
F_{3,0} & F_{0,3} & 6F_{3,0}-8 & -6F_{0,3}
\\
\green{\bf 3b} & 
1 & 0 & 6 & 0
\\
\green{\bf 3c} & 
0 & 0 & 0 & 0
\end{array}
\]

\medskip

In branch \green{\bf 2g3c}:
\[
\left(\!
\begin{array}{c}
R_{2,2}
\\
S_{0,4}
\end{array}
\!\right)
\,=\,
\left(\!
\def\arraystretch{1.25}
\begin{array}{cc}
\tfrac{a_{1,1}^2-a_{2,1}^2}{(a_{1,1}^2+a_{2,1}^2)^2} & \tfrac{6a_{1,1}a_{2,1}}{(a_{1,1}^2+a_{2,1}^2)^2}
\\
-\tfrac{2}{3}\tfrac{a_{1,1}a_{2,1}}{(a_{1,1}^2+a_{2,1}^2)^2} & \tfrac{a_{1,1}^2-a_{2,1}^2}{(a_{1,1}^2+a_{2,1}^2)^2}
\end{array}
\!\right)\,
\left(\!
\begin{array}{c}
F_{2,2}
\\
G_{0,4}
\end{array}
\!\right),
\]
leading to:
\[
\def\arraystretch{1.25}
\begin{array}{rccc}
\green{\bf 2g3c}\,\,\,\,
\green{\downarrow}\,\,
& F_{2,2} & G_{0,4}
\\
\green{\bf 4a} & 
1 & 0
\\
\green{\bf 4b} & 
0 & 0
\end{array}
\]

%%%%%%%%%%%%%%%%%%%%%%%%%%%%%%%%%%%%%%%%%%%%%%%%%%%%%%%%%%%%%%%%%%%%%%
\SectionHead{Model \green{\bf 2e3b4a}}
{model-2e3b4a}
%%%%%%%%%%%%%%%%%%%%%%%%%%%%%%%%%%%%%%%%%%%%%%%%%%%%%%%%%%%%%%%%%%%%%%

\begin{center}
\includegraphics[scale=0.19]{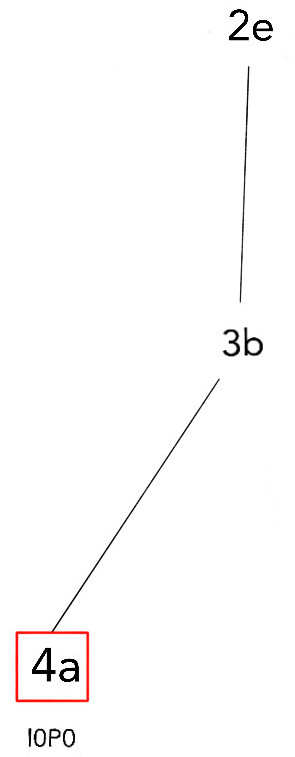}
\end{center}

In this section, we assume to be in the non-parallel case \green{\bf 2e}:
\begin{align}
\label{2e}
u
&
=
xy
+
{\rm O}_{x,y}(3),
\\
\notag
v
&
=
x^2
+
{\rm O}_{x,y}(3).
\end{align}
If an affine infinitesimal transformation $L$ is tangent to~{\eqref{2e}}, then at order $1$:
\[
\aligned
0
&
\underset{\green{\bf F}}{\overset{\green{\bf 1,0}}{\,\,=\,\,}}
-\,C_{1,1}+T_2,
\\
0
&
\underset{\green{\bf F}}{\overset{\green{\bf 0,1}}{\,\,=\,\,}}
-\,C_{1,2}+T_1,
\\
0
&
\underset{\green{\bf G}}{\overset{\green{\bf 1,0}}{\,\,=\,\,}}
-\,C_{2,1}+2T_1,
\\
0
&
\underset{\green{\bf G}}{\overset{\green{\bf 0,1}}{\,\,=\,\,}}
-\,C_{2,2},
\endaligned
\]
whence necessarily:
\[
C_{1,1}=T_2, \,\,\,\,\,\,\,\,\,\,\,\, 
C_{1,2}=T_1,  \,\,\,\,\,\,\,\,\,\,\,\,
C_{2,1}=2T_1, \,\,\,\,\,\,\,\,\,\,\,\, 
C_{2,2}=0.´
\]
Next, at order $2$, the two fundamental equations give:
\[
\aligned
0
&
\underset{\green{\bf R}}{\overset{\green{\bf 2,0}}{\,\,=\,\,}}
a_{1,1}a_{2,1}-d_{1,2},
\\
0
&
\underset{\green{\bf R}}{\overset{\green{\bf 1,1}}{\,\,=\,\,}}
a_{1,1}a_{2,2}+a_{1,2}a_{2,1}-d_{1,1},
\\
0
&
\underset{\green{\bf S}}{\overset{\green{\bf 2,0}}{\,\,=\,\,}}
a_{1,1}^2-d_{2,2},
\\
0
&
\underset{\green{\bf S}}{\overset{\green{\bf 1,1}}{\,\,=\,\,}}
2a_{1,1}a_{1,2}-d_{2,1},
\\
0
&
\underset{\green{\bf S}}{\overset{\green{\bf 0,2}}{\,\,=\,\,}}
a_{1,2}^2.
\endaligned
\]
Therefore, we obtain:
\begin{Lemma}\label{Lemma:G_{stab}^2}
The subgroup $G_{\text{\rm stab}}^{2}\subset G_{\text{\rm stab}}^{1}$ which sends
\[
\aligned
u
\,=\,
xy
+
{\rm O}_{x,y}(3),
\\
v
\,=\,
x^2
+
{\rm O}_{x,y}(3),
\endaligned
\ \ \ \ \ \ \ \ \ \ \ \ \ \ \ \ \ \ \ \
\text{to}
\ \ \ \ \ \ \ \ \ \ \ \ \ \ \ \ \ \ \ \
\aligned
r
\,=\,
pq
+
{\rm O}_{p,q}(3),
\\
s
\,=\,
p^2
+
{\rm O}_{p,q}(3),
\endaligned
\]
consists of matrices:
\[
\begin{pmatrix}
a_{1,1} & 0 & b_{1,1} & b_{1,2} \\
a_{2,1} & a_{2,2} & b_{2,1} & b_{2,2} \\
0 & 0 & a_{1,1}a_{2,2} & a_{1,1}a_{2,1} \\
0 & 0 & 0 & a_{1,1}^2
\end{pmatrix}
\in
G_{\text{\rm stab}}^{2}.
\]
\end{Lemma}
Next, the tangency of the affine infinitesimal transformation $L$ to:
\[
\aligned
u
&
=
xy
+
F_{3,0}x^3
+
F_{2,1}x^2y
+
F_{1,2}xy^2
+
F_{0,3}y^3
+
{\rm O}_{x,y}(4),
\\
v
&
=
x^2
+
G_{3,0}x^3
+
G_{2,1}x^2y
+
{\rm O}_{x,y}(4),
\endaligned
\]
gives, at order $2$:
\[
\aligned
0
&
\underset{\green{\bf F}}{\overset{\green{\bf 2,0}}{\,\,=\,\,}}
F_{2,1}T_2+3F_{3,0}T_1+A_{2,1}-D_{1,2},
\\
0
&
\underset{\green{\bf F}}{\overset{\green{\bf 1,1}}{\,\,=\,\,}}
2F_{1,2}T_2+2F_{2,1}T_1+A_{1,1}+A_{2,2}-D_{1,1},
\\
0
&
\underset{\green{\bf F}}{\overset{\green{\bf 0,2}}{\,\,=\,\,}}
3F_{0,3}T_2+F_{1,2}T_1+A_{1,2},
\\
0
&
\underset{\green{\bf G}}{\overset{\green{\bf 2,0}}{\,\,=\,\,}}
G_{2,1}T_2+3G_{3,0}T_1+2A_{1,1}-D_{2,2},
\\
0
&
\underset{\green{\bf G}}{\overset{\green{\bf 1,1}}{\,\,=\,\,}}
2G_{1,2}T_2+2G_{2,1}T_1+2A_{1,2}-D_{2,1},
\\
0
&
\underset{\green{\bf G}}{\overset{\green{\bf 0,2}}{\,\,=\,\,}}
3G_{0,3}T_2+G_{1,2}T_1,
\endaligned
\]
whence necessarily:
\[
\aligned
A_{1,2}
&
=
-3F_{0,3}T_2-F_{1,2}T_1,
\\
D_{1,1}
&
=
2F_{1,2}T_2
+
2F_{2,1}T_1
+
A_{1,1}
+
A_{2,2},
\\
D_{1,2}
&
=
F_{2,1}T_2
+
3F_{3,0}T_1
+
A_{2,1},
\\
D_{2,1}
&
=
-6F_{0,3}T_2-2F_{1,2}T_1+2G_{1,2}T_2+2G_{2,1}T_1,
\\
D_{2,2}
&
=
G_{2,1}T_2+3G_{3,0}T_1+2A_{1,1},
\\
G_{0,3}
&
=
0,
\\
G_{1,2}
&
=
0.
\endaligned
\]

Next, at order $3$, the fundamental equations give:
\[
\aligned
0
&
\underset{\green{\bf R}}{\overset{\green{\bf 3,0}}{\,\,=\,\,}}
R_{0,3}a_{2,1}^3+R_{1,2}a_{1,1}a_{2,1}^2+R_{2,1}a_{1,1}^2a_{2,1}+R_{3,0}a_{1,1}^3-F_{3,0}a_{1,1}a_{2,2}-G_{3,0}a_{1,1}a_{2,1}+a_{1,1}b_{2,2}
+
\\
&
\,\,\,\,\,\,\,\,\,\,
+a_{2,1}b_{1,2},
\\
0
&
\underset{\green{\bf R}}{\overset{\green{\bf 2,1}}{\,\,=\,\,}}
3R_{0,3}a_{2,1}^2a_{2,2}+2R_{1,2}a_{1,1}a_{2,1}a_{2,2}+R_{2,1}a_{1,1}^2a_{2,2}-F_{2,1}a_{1,1}a_{2,2}-G_{2,1}a_{1,1}a_{2,1}+a_{1,1}b_{2,1}+a_{2,1}b_{1,1}
+
\\
&
\,\,\,\,\,\,\,\,\,\,
+a_{2,2}b_{1,2},
\\
0
&
\underset{\green{\bf R}}{\overset{\green{\bf 1,2}}{\,\,=\,\,}}
3R_{0,3}a_{2,1}a_{2,2}^2+R_{1,2}a_{1,1}a_{2,2}^2-F_{1,2}a_{1,1}a_{2,2}+a_{2,2}b_{1,1},
\\
0
&
\underset{\green{\bf R}}{\overset{\green{\bf 0,3}}{\,\,=\,\,}}
R_{0,3}a_{2,2}^3-F_{0,3}a_{1,1}a_{2,2},
\\
0
&
\underset{\green{\bf S}}{\overset{\green{\bf 3,0}}{\,\,=\,\,}}
S_{2,1}a_{1,1}^2a_{2,1}+S_{3,0}a_{1,1}^3-G_{3,0}a_{1,1}^2+2a_{1,1}b_{1,2},
\\
0
&
\underset{\green{\bf S}}{\overset{\green{\bf 2,1}}{\,\,=\,\,}}
S_{2,1}a_{1,1}^2a_{2,2}-G_{2,1}a_{1,1}^2+2a_{1,1}b_{1,1}.
\endaligned
\]
By freeness of $b_{1,1}$ in $\underset{\green{\bf S}}{\overset{\green{\bf 2,1}}{\,=\,}}$, $G_{\text{\rm stab}}^{2}$ can act on the surface $\{u=F(x,y),\, v=G(x,y)\}$ by the element:
\[
g
=
\begin{pmatrix}
1 & 0 & \tfrac{a_{1,1}}{2}G_{2,1} & 0 \\
0 & 1 & 0 & 0 \\
0 & 0 & 1 & 0 \\
0 & 0 & 0 & 1
\end{pmatrix},
\]
to normalize $S_{2,1}=0$ and by equivalence:
\[
G_{2,1}
=
0
=
S_{2,1}.
\]
After that, $\underset{\green{\bf S}}{\overset{\green{\bf 2,1}}{\,=\,}}$ gives:
\[
b_{1,1}=0.
\]

Then, by looking at $\underset{\green{\bf S}}{\overset{\green{\bf 3,0}}{\,=\,}}$, and by an appropriate choice of $b_{1,2}$, we can normalize $S_{3,0}=0$, and by equivalence:
\[
G_{3,0}
=
0
=
S_{3,0}.
\]
After that, $\underset{\green{\bf S}}{\overset{\green{\bf 3,0}}{\,=\,}}$ gives:
\[
b_{1,2}
=
0.
\]

By the same process, in $\underset{\green{\bf R}}{\overset{\green{\bf 3,0}}{\,=\,}}$ , we normalize:
\[
F_{3,0}
=
0
=
R_{3,0},
\]
and we stabilize:
\[
b_{2,2}
=
-\tfrac{a_{2,1}}{a_{1,1}}^3R_{0,3}
-a_{2,1}^2R_{1,2}-a_{2,1}a_{1,1}R_{2,1}.
\]

Next, in $\underset{\green{\bf R}}{\overset{\green{\bf 2,1}}{\,=\,}}$, we normalize and stabilize:
\[
\aligned
F_{2,1}
&
=
0
=
R_{2,1},
\\
b_{2,1}
=
-3\tfrac{a_{2,1}}{a_{1,1}}^2 & a_{2,2}R_{0,3}
-2a_{2,1}a_{2,2}R_{1,2}.
\endaligned
\]
The two remaining equations $\underset{\green{\bf R}}{\overset{\green{\bf 1,2}}{\,=\,}}$ and $\underset{\green{\bf R}}{\overset{\green{\bf 0,3}}{\,=\,}}$ create an associated linear representation:
\[
\begin{pmatrix}
R_{1,2} \\
R_{0,3}
\end{pmatrix}
=
\begin{pmatrix}
\tfrac{1}{a_{2,2}} & -\tfrac{3a_{2,1}}{a_{2,2}} \\
0 & \tfrac{a_{1,1}}{a_{2,2}^2}
\end{pmatrix}
\begin{pmatrix}
F_{1,2} \\
F_{0,3}
\end{pmatrix}.
\]
\begin{Lemma}\label{Lemma:G_{stab}^3}
The subgroup $G_{\text{\rm stab}}^{3}\subset G_{\text{\rm stab}}^{2}$ which sends
\[
\aligned
u
&
\,=\,
xy
+
F_{1,2}xy^2
+
F_{0,3}y^3
+
{\rm O}_{x,y}(4),
\\
v
&
\,=\,
x^2
+
{\rm O}_{x,y}(4),
\endaligned
\ \ \ \ \ \ \ \ \ \ \ \ \ \ \ \ \ \ \ \
\text{to}
\ \ \ \ \ \ \ \ \ \ \ \ \ \ \ \ \ \ \ \
\aligned
r
&
\,=\,
pq
+
R_{1,2}pq^2
+
R_{0,3}q^3
+
{\rm O}_{p,q}(4),
\\
s
&
\,=\,
p^2
+
{\rm O}_{p,q}(4),
\endaligned
\]
consists of matrices:
\[
\begin{pmatrix}
a_{1,1} & 0 & 0 & 0 \\
a_{2,1} & a_{2,2} & -\tfrac{a_{2,1}a_{2,2}}{a_{1,1}}(3R_{0,3}a_{2,1}+2R_{1,2}a_{1,1}) & -\tfrac{a_{2,1}^2}{a_{1,1}}(R_{0,3}a_{2,1}+R_{1,2}a_{1,1}) \\
0 & 0 & a_{1,1}a_{2,2} & a_{1,1}a_{2,1} \\
0 & 0 & 0 & a_{1,1}^2
\end{pmatrix}
\in G_{\text{\rm stab}}^{3}.
\]
\end{Lemma}
Next, the tangency of the affine infinitesimal transformation $L$ to:
\[
\aligned
u
&
=
xy
+
xy^2
+
F_{4,0}x^4
+
F_{3,1}x^3y
+
F_{2,2}x^2y^2
+
xy^3
+
{\rm O}_{x,y}(5),
\\
v
&
=
x^2
+
G_{4,0}x^4
+
G_{3,1}x^3y
-
x^2y^2
+
{\rm O}_{x,y}(5),
\endaligned
\]
at order $3$, gives:
\[
\aligned
B_{1,1}
&
=
-G_{2,2}T_2-\tfrac{3}{2}G_{3,1}T_1,
\\
B_{1,2}
&
=
-\tfrac{1}{2}G_{3,1}T_2-2G_{4,0}T_1,
\\
B_{2,1}
&
=
-2F_{1,2}A_{2,1}+\tfrac{1}{2}G_{3,1}T_2-2F_{2,2}-3F_{3,1}T_1+2G_{4,0}T_1,
\\
B_{2,2}
&
=
-F_{3,1}T_2-4F_{4,0}T_1.
\endaligned
\]

Then, we open the branch $\green{\bf 2e3b}$:
\[
\aligned
F_{0,3}
&
=
0
=
R_{0,3},
\\
F_{1,2}
&
=
1
=
R_{1,2}.
\endaligned
\]
At order $3$, the fundamental equations give:
\[
a_{2,2}
=
1,
\]
and the tangency:
\[
\aligned
0
&
\underset{\green{\bf F}}{\overset{\green{\bf 1,2}}{\,\,=\,\,}}
A_{2,2}+(-G_{2,2}+3F_{1,3}-2)T_2+(-\tfrac{3}{2}G_{3,1}+2F_{2,2})T_1,
\\
0
&
\underset{\green{\bf F}}{\overset{\green{\bf 0,3}}{\,\,=\,\,}}
(F_{3,1}-1)T_1+4F_{0,4}T_2,
\\
0
&
\underset{\green{\bf G}}{\overset{\green{\bf 1,2}}{\,\,=\,\,}}
(2G_{2,2}+2)T_1+3G_{1,3}T_2,
\\
0
&
\underset{\green{\bf G}}{\overset{\green{\bf 0,3}}{\,\,=\,\,}}
4G_{0,4}T_2+G_{1,3}T_1,
\endaligned
\]
therefore:
\[
\aligned
A_{2,2}
&
=
-\,3F_{1,3}T_2+G_{2,2}T_2+2T_2+\tfrac{3}{2}G_{3,1}T_1-2F_{2,2}T_1,
\\
F_{1,3}
&
=
1
=
R_{1,3},
\,\,\,\,\,\,\,\,\,\,\,\,
G_{2,2}
=
-1
=
S_{2,2},
\\
F_{0,4}
&
=
0
=
R_{0,4},
\,\,\,\,\,\,\,\,\,\,\,\,
G_{0,4}
=
0
=
S_{0,4}.
\endaligned
\]

Next, at order $4$, the fundamental equations give:
\[
\aligned
0
&
\underset{\green{\bf R}}{\overset{\green{\bf 0,4}}{\,\,=\,\,}}
R_{2,2}a_{1,1}^2a_{2,1}^2+R_{3,1}a_{1,1}^3a_{2,1}+R_{4,0}a_{1,1}^4-a_{1,1}a_{2,1}^3-G_{4,0}a_{1,1}a_{2,1}-F_{4,0}a_{1,1},
\\
0
&
\underset{\green{\bf R}}{\overset{\green{\bf 3,1}}{\,\,=\,\,}}
2R_{2,2}a_{1,1}^2a_{2,1}+R_{3,1}a_{1,1}^3-G_{3,1}a_{1,1}a_{2,1}-3a_{1,1}a_{2,1}^2-F_{3,1}a_{1,1},
\\
0
&
\underset{\green{\bf R}}{\overset{\green{\bf 2,2}}{\,\,=\,\,}}
R_{2,2}a_{1,1}^2-F_{2,2}a_{1,1}-2a_{1,1}a_{2,1},
\\
0
&
\underset{\green{\bf S}}{\overset{\green{\bf 4,0}}{\,\,=\,\,}}
S_{3,1}a_{1,1}^3a_{2,1}+S_{4,0}a_{1,1}^4-a_{1,1}^2a_{2,1}^2-G_{4,0}a_{1,1}^2,
\\
0
&
\underset{\green{\bf S}}{\overset{\green{\bf 3,1}}{\,\,=\,\,}}
S_{3,1}a_{1,1}^3-G_{3,1}a_{1,1}^2-2a_{1,1}^2a_{2,1}.
\endaligned
\]
By freeness of $a_{2,1}$ in $\underset{\green{\bf S}}{\overset{\green{\bf 3,1}}{\,\,=\,\,}}$, we can normalize $S_{3,1}=0$ and by equivalence:
\[
G_{3,1}
=
0
=
S_{3,1}.
\]
After that, $\underset{\green{\bf S}}{\overset{\green{\bf 3,1}}{\,\,=\,\,}}$ gives:
\[
a_{2,1}
=
0.
\]
The remaining equations provide at order $4$ the linear representation:
\[
\begin{pmatrix}
R_{2,2} \\
R_{3,1} \\
S_{4,0} \\
R_{4,0}
\end{pmatrix}
=
\begin{pmatrix}
\tfrac{1}{a_{1,1}} & 0 & 0 & 0 \\
0 & \tfrac{1}{a_{1,1}^2} & 0 & 0 \\
0 & 0 & \tfrac{1}{a_{1,1}^2} & 0 \\
0 & 0 & 0 & \tfrac{1}{a_{1,1}^3}
\end{pmatrix}
\begin{pmatrix}
F_{2,2} \\
F_{3,1} \\
G_{4,0} \\
F_{4,0}
\end{pmatrix}.
\]
\begin{Lemma}
The subgroup $G_{\text{\rm stab}}^{4}\subset G_{\text{\rm stab}}^{3}$ which sends:
\[
\aligned
u
&
=
xy
+
xy^2
+
F_{4,0}x^4
+
F_{3,1}x^3y
+
F_{2,2}x^2y^2
+
xy^3
+
{\rm O}_{x,y}(5),
\\
v
&
=
x^2
+
G_{4,0}x^4
-
x^2y^2
+
{\rm O}_{x,y}(5),
\endaligned
\]
to:
\[
\aligned
r
&
=
pq
+
pq^2
+
R_{4,0}p^4
+
R_{3,1}p^3q
+
R_{2,2}p^2q^2
+
pq^3
+
{\rm O}_{p,q}(5),
\\
s
&
=
p^2
+
S_{4,0}p^4
-
p^2q^2
+
{\rm O}_{p,q}(5),
\endaligned
\]
consists of matrices:
\[
\begin{pmatrix}
a_{1,1} & 0 & 0 & 0 \\
0 & 1 & 0 & 0 \\
0 & 0 & a_{1,1}& 0\\
0 & 0 & 0 & a_{1,1}^2
\end{pmatrix}
\in G_{\text{\rm stab}}^{4}.
\]
\end{Lemma}
Next, the tangency to:
\[
\aligned
u
&
=
xy
+
xy^2
+
F_{4,0}x^4
+
F_{3,1}x^3y
+
F_{2,2}x^2y^2
+
xy^3
+
F_{5,0}x^5
+
F_{4,1}x^4y
+
F_{3,2}x^3y^2
+
\\
&
\,\,\,\,\,
+
F_{2,3}x^2y^3
+
F_{1,4}xy^4
+
F_{0,5}y^5
+
{\rm O}_{x,y}(6),
\\
v
&
=
x^2
+
G_{4,0}x^4
-
x^2y^2
+
G_{5,0}x^5
+
G_{4,1}x^4y
+
G_{3,2}x^3y^2
+
G_{2,3}x^2y^3
+
G_{1,4}xy^4
+
G_{0,5}y^5
+
{\rm O}_{x,y}(6),
\endaligned
\]
at order $4$, gives:
\[
A_{2,1}
=
F_{3,1}T_1+G_{3,2}T_2-2G_{4,0}T_1+2G_{4,1}T_1.
\]

Finally, we open the branch $\green{\bf 2e3b4a}$:
\[
F_{1,2}
=
1
=
R_{1,2}.
\]
At order $4$, we observe in $\underset{\green{\bf R}}{\overset{\green{\bf 2,2}}{\,\,=\,\,}}$:
\[
a_{1,1}
=
1,
\]
and in the associated tangency equation:
\[
A_{1,1}
=
-3F_{2,3}T_2+14F_{3,1}T_1-3F_{3,2}T_1+2G_{3,2}T_2-8G_{4,0}T_1+4G_{4,1}T_1+2T_1+10T_2.
\]

Moreover, the remaining equations $\underset{\green{\bf S}}{\overset{\green{\bf *,*}}{\,\,=\,\,}}$ give normalizations for $F_{4,0}$ and for:
\[
F_{i,j}
\eqno
{\scriptstyle{(i+j=5),}}
\]
\[
G_{k,l}
\eqno
{\scriptstyle{(k+l=5).}}
\]

At this step, all the coefficients:
\[
F_{i,j}
\eqno
{\scriptstyle{(1\leqslant i+j \leqslant 5),}}
\]
\[
G_{k,l}
\eqno
{\scriptstyle{(1\leqslant k+l\leqslant 5),}}
\]
are normalized, except $F_{3,1}$ and $G_{4,0}$.

Next, at order $5$, the fundamental equations no longer provide possible normalization, because:
\begin{Lemma}
The subgroup $G_{\text{\rm stab}}^{5}\subset G_{\text{\rm stab}}^{4}$ which sends:
\[
\footnotesize
\aligned
u
&
=
xy+xy^2+\big(-\tfrac{1}{12}F_{3,1}-\tfrac{4}{3}F_{3,1}G_{4,0}-\tfrac{1}{2}G_{4,0}+\tfrac{1}{12}F_{3,1}^2\big)x^4+F_{3,1}x^3y+x^2y^2+xy^3
+
\\
&
\,\,\,\,\,
+
\Big(-\tfrac{1}{3}F_{3,1}^2-\tfrac{42}{5}F_{3,1}(-\tfrac{1}{12}F_{3,1}-\tfrac{4}{3}F_{3,1}G_{4,0}-\tfrac{1}{2}G_{4,0}+\tfrac{1}{12}F_{3,1}^2)+\tfrac{1}{3}F_{3,1}G_{4,0}
+
\\
&
\,\,\,\,\,
+\tfrac{9}{5}(5F_{3,1}+2)(-\tfrac{1}{12}F_{3,1}-\tfrac{4}{3}F_{3,1}G_{4,0}-\tfrac{1}{2}G_{4,0}+\tfrac{1}{12}F_{3,1}^2)-\tfrac{48}{5}(-\tfrac{1}{12}F_{3,1}-\tfrac{4}{3}F_{3,1}G_{4,0}
+
\\
&
\,\,\,\,\,
-\tfrac{1}{2}G_{4,0}+\tfrac{1}{12}F_{3,1}^2)G_{4,0}+2G_{4,0}^2+\tfrac{2}{15}F_{3,1}+\tfrac{4}{5}G_{4,0}\Big)x^5
+
\Big(-\tfrac{27}{4}F_{3,1}^2+\tfrac{3}{2}F_{3,1}(5F_{3,1}+2)
+
\\
&
\,\,\,\,\,
-12F_{3,1}G_{4,0}-\tfrac{7}{4}F_{3,1}-\tfrac{13}{2}G_{4,0}\Big)x^4y+(5F_{3,1}+2)x^3y^2+3x^2y^3+xy^4
+
{\rm O}_{x,y}(6),
\\
v
&
=
x^2+G_{4,0}x^4-x^2y^2+\Big(-\tfrac{76}{15}F_{3,1}G_{4,0}+\tfrac{6}{5}(5F_{3,1}+2)G_{4,0}-\tfrac{32}{5}G_{4,0}^2+\tfrac{1}{30}F_{3,1}-\tfrac{3}{5}G_{4,0}
+
\\
&
\,\,\,\,\,
-\tfrac{1}{30}F_{3,1}^2\Big)x^5+6G_{4,0}x^4y-2x^3y^2-2x^2y^3
+
{\rm O}_{x,y}(6),
\endaligned
\]
to:
\[
\footnotesize
\aligned
r
&
=
pq+pq^2+\big(-\tfrac{1}{12}F_{3,1}-\tfrac{4}{3}F_{3,1}G_{4,0}-\tfrac{1}{2}G_{4,0}+\tfrac{1}{12}F_{3,1}^2\big)p^4+F_{3,1}p^3q+p^2q^2+pq^3
+
\\
&
\,\,\,\,\,
+
\Big(-\tfrac{1}{3}F_{3,1}^2-\tfrac{42}{5}F_{3,1}(-\tfrac{1}{12}F_{3,1}-\tfrac{4}{3}F_{3,1}G_{4,0}-\tfrac{1}{2}G_{4,0}+\tfrac{1}{12}F_{3,1}^2)+\tfrac{1}{3}F_{3,1}G_{4,0}
+
\\
&
\,\,\,\,\,
+\tfrac{9}{5}(5F_{3,1}+2)(-\tfrac{1}{12}F_{3,1}-\tfrac{4}{3}F_{3,1}G_{4,0}-\tfrac{1}{2}G_{4,0}+\tfrac{1}{12}F_{3,1}^2)-\tfrac{48}{5}(-\tfrac{1}{12}F_{3,1}-\tfrac{4}{3}F_{3,1}G_{4,0}
+
\\
&
\,\,\,\,\,
-\tfrac{1}{2}G_{4,0}+\tfrac{1}{12}F_{3,1}^2)G_{4,0}+2G_{4,0}^2+\tfrac{2}{15}F_{3,1}+\tfrac{4}{5}G_{4,0}\Big)p^5
+
\\
&
\,\,\,\,\,
+
\Big(-\tfrac{27}{4}F_{3,1}^2+\tfrac{3}{2}F_{3,1}(5F_{3,1}+2)
+
\\
&
\,\,\,\,\,
-12F_{3,1}G_{4,0}-\tfrac{7}{4}F_{3,1}-\tfrac{13}{2}G_{4,0}\Big)p^4q+(5F_{3,1}+2)p^3q^2+3p^2q^3+pq^4
+
\\
&
\,\,\,\,\,
+
{\rm O}_{p,q}(6),
\\
s
&
=
p^2+G_{4,0}p^4-p^2q^2+\Big(-\tfrac{76}{15}F_{3,1}G_{4,0}+\tfrac{6}{5}(5F_{3,1}+2)G_{4,0}-\tfrac{32}{5}G_{4,0}^2+\tfrac{1}{30}F_{3,1}-\tfrac{3}{5}G_{4,0}
+
\\
&
\,\,\,\,\,
-\tfrac{1}{30}F_{3,1}^2\Big)p^5+6G_{4,0}p^4q-2p^3q^2-2p^2q^3
+
\\
&
\,\,\,\,\,
+
{\rm O}_{p,q}(6),
\endaligned
\]
is reduced to the identity group $\{e\}$.
\end{Lemma}
Next, the tangency equations to:
\[
\aligned
u
&
=
xy+xy^2+\big(-\tfrac{1}{12}F_{3,1}-\tfrac{4}{3}F_{3,1}G_{4,0}-\tfrac{1}{2}G_{4,0}+\tfrac{1}{12}F_{3,1}^2\big)x^4+F_{3,1}x^3y+x^2y^2+xy^3
+
\\
&
\,\,\,\,\,
+
\Big(-\tfrac{1}{3}F_{3,1}^2-\tfrac{42}{5}F_{3,1}(-\tfrac{1}{12}F_{3,1}-\tfrac{4}{3}F_{3,1}G_{4,0}-\tfrac{1}{2}G_{4,0}+\tfrac{1}{12}F_{3,1}^2)+\tfrac{1}{3}F_{3,1}G_{4,0}
+
\\
&
\,\,\,\,\,
+\tfrac{9}{5}(5F_{3,1}+2)(-\tfrac{1}{12}F_{3,1}-\tfrac{4}{3}F_{3,1}G_{4,0}-\tfrac{1}{2}G_{4,0}+\tfrac{1}{12}F_{3,1}^2)-\tfrac{48}{5}(-\tfrac{1}{12}F_{3,1}-\tfrac{4}{3}F_{3,1}G_{4,0}
+
\\
&
\,\,\,\,\,
-\tfrac{1}{2}G_{4,0}+\tfrac{1}{12}F_{3,1}^2)G_{4,0}+2G_{4,0}^2+\tfrac{2}{15}F_{3,1}+\tfrac{4}{5}G_{4,0}\Big)x^5
+
\Big(-\tfrac{27}{4}F_{3,1}^2+\tfrac{3}{2}F_{3,1}(5F_{3,1}+2)
+
\\
&
\,\,\,\,\,
-12F_{3,1}G_{4,0}-\tfrac{7}{4}F_{3,1}-\tfrac{13}{2}G_{4,0}\Big)x^4y+(5F_{3,1}+2)x^3y^2+3x^2y^3+xy^4
+
F_{6,0}x^6
+
\\
&
\,\,\,\,\,
+
F_{5,1}x^5y
+
F_{4,2}x^4y^2
+
F_{3,3}x^3y^3
+
F_{2,4}x^2y^4
+
F_{1,5}xy^5
+
F_{0,6}y^6
+
{\rm O}_{x,y}(7),
\\
v
&
=
x^2+G_{4,0}x^4-x^2y^2+\Big(-\tfrac{76}{15}F_{3,1}G_{4,0}+\tfrac{6}{5}(5F_{3,1}+2)G_{4,0}-\tfrac{32}{5}G_{4,0}^2+\tfrac{1}{30}F_{3,1}-\tfrac{3}{5}G_{4,0}
+
\\
&
\,\,\,\,\,
-\tfrac{1}{30}F_{3,1}^2\Big)x^5+6G_{4,0}x^4y-2x^3y^2-2x^2y^3
+
+
G_{6,0}x^6
+
G_{5,1}x^5y
+
G_{4,2}x^4y^2
+
G_{3,3}x^3y^3
+
\\
&
\,\,\,\,\,
+
G_{2,4}x^2y^4
+
G_{1,5}xy^5
+
G_{0,6}y^6
+
{\rm O}_{x,y}(7),
\endaligned
\]
at order $5$ give normalizations for $F_{3,1},G_{4,0}$, and for all the $F_{i,j},G_{k,l}$ of order $6$.

Next, at orders $n\geqslant 7$, the tangency equations for $i+j=n$ and $k+l=n$ are of the form:
\[
\aligned
0
&
\underset{\green{\bf F}}{\overset{\green{\bf i,j}}{\,\,=\,\,}}
\big(
\alpha_{i,j}F_{i+1,j}+\beta_{i,j}
\big)\,
T_1
+
\big(
\gamma_{i,j}F_{i,j+1}+\delta_{i,j}
\big)\,
T_2,
\\
0
&
\underset{\green{\bf G}}{\overset{\green{\bf k,l}}{\,\,=\,\,}}
\big(
\epsilon_{k,l}G_{k+1,l}+\zeta_{k,l}
\big)\,
T_1
+
\big(
\eta_{k,l}G_{k,l+1}+\theta_{k,l}
\big)\,
T_2,
\endaligned
\]
with
$
\alpha_{i,j}, \,
\gamma_{i,j}, \,
\epsilon_{k,l}, \,
\eta_{k,l}, \in\R^*$
and
$
\beta_{i,j}, \,
\delta_{i,j}, \,
\zeta_{k,l}, \,
\theta_{k,l}\in \R$.
Therefore, all coefficients are determined at each order $n$.

Finally, we obtain the two generators of the Lie algebra:
\[
\aligned
e_1
&
:=
L\Big\vert_{T_1=1,T_2=0},
\\
e_2
&
:=
L\Big\vert_{T_1=0,T_2=1},
\endaligned
\]
which satisfy:
\[
[e_1,e_2]
=
-2e_1.
\]

\begin{Theorem}
In the branch $\green{\bf 2e3b4a}$, every affinely homogeneous surface is, in a unique way, equivalent to:
\[\left\{
\aligned
u
&
\,=\,
xy
+
xy^2
+
xy^3
+
x^2y^2
-
x^3y
+
\tfrac{3}{8}x^4
+
xy^4
+
3x^2y^3
-
3x^3y^2
+
\\
&
\ \ \ \ \
+
\tfrac{7}{8}x^4y
+
xy^5
+
6x^2y^4
-
4x^3y^3
-
\tfrac{21}{8}x^4y^2
+
\tfrac{9}{4}x^5y
-
\tfrac{3}{8}x^6
+
xy^6
+
\\
&
\ \ \ \ \
+
10x^2y^5
-
\tfrac{145}{8}x^4y^3
+
10x^5y^2
-
\tfrac{3}{2}x^6y
+
\cdots,
\\
v
&
\,=\,
x^2
-
x^2y^2
+
\tfrac{1}{4}x^4
-
2x^2y^3
-
2x^3y^2
+
\tfrac{3}{2}x^4y
-
\tfrac{1}{4}x^5
-
3x^2y^4
+
\\
&
\ \ \ \ \
-
8x^3y^3
+
\tfrac{7}{2}x^4y^2
-
\tfrac{1}{8}x^6
-
4x^2y^5
-
20x^3y^4
+
\tfrac{29}{4}x^5y^2
+
\\
&
\ \ \ \ \
-
3x^6y
+
\tfrac{3}{8}x^7
+\cdots
,
\endaligned \right.
\]
with $2$-dimensional Lie algebra:
\[
\aligned
e_1
&
\,:=\,
(
x
+
1
-
y
-
\tfrac{1}{2}v
)\partial_x
+
(
\tfrac{3}{2}x
-
2y
+
\tfrac{1}{2}u
-
\tfrac{3}{2}v
)\partial_y
+
(
y
-
u
+
\tfrac{3}{2}v
)\partial_u
+
\\
&
\ \ \ \ \
-
(
2u
-
2v
-
2x)\partial_v 
,
\\
e_2
&
\,:=\,
(
-
3x
+
u
)\partial_x
+
(
1
-
2x
-
2y
+
2u
+
v
)\partial_y
+
(
x
-
3u
-
2v
)\partial_u-6v\partial_v
,
\endaligned
\]
and Lie structure:
\[
[e_1,e_2]
=
-2e_1.
\]
\end{Theorem}

%%%%%%%%%%%%%%%%%%%%%%%%%%%%%%%%%%%%%%%%%%%%%%%%%%%%%%%%%%%%%%%%%%%%%%
\SectionHead{Model \green{\bf 2e3c4e5b}}
{model-2e3c4e5b}
%%%%%%%%%%%%%%%%%%%%%%%%%%%%%%%%%%%%%%%%%%%%%%%%%%%%%%%%%%%%%%%%%%%%%%

\begin{center}
\includegraphics[scale=0.19]{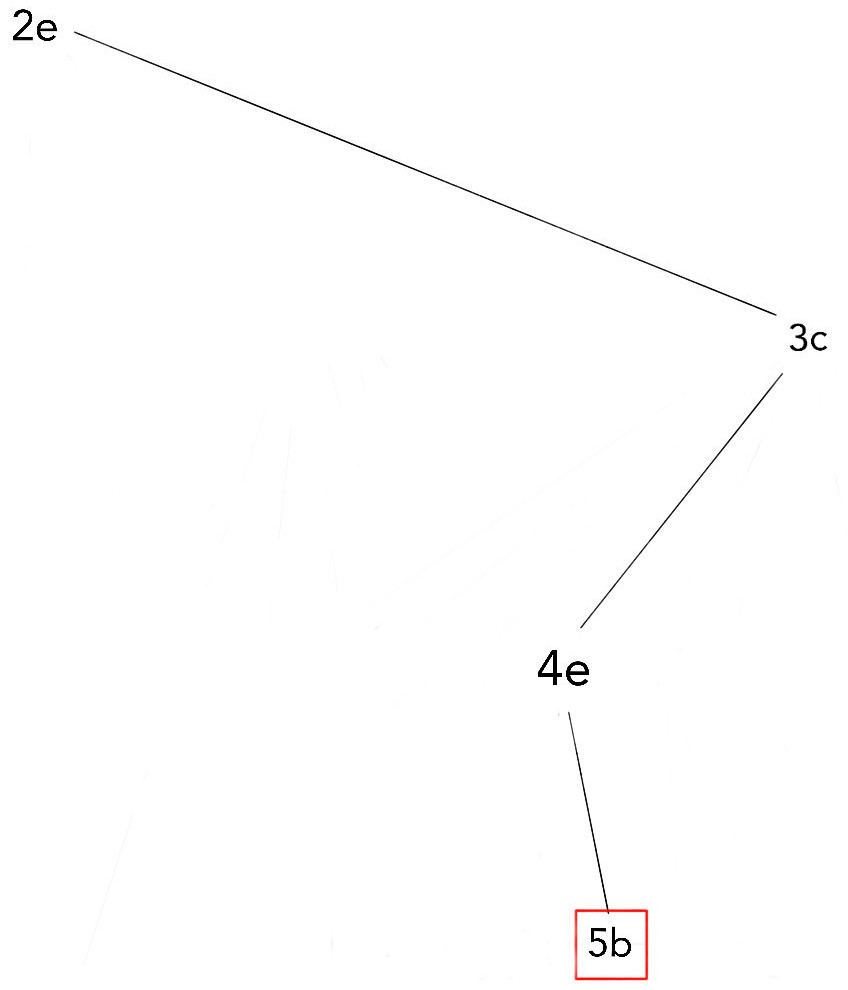}
\end{center}

Since this case is \green{\bf 2e3c4e5b}, we have the same result of \green{\bf 2e} before the opening the sub-branch \green{\bf 2e3b}: Lemmas~{\ref{Lemma:G_{stab}^2}},~{\ref{Lemma:G_{stab}^3}} and the normalizations for $L$:
\[
\aligned
C_{1,1} &=T_2, \,\,\,\,\,\,\,\,\,\,\,\,\,
C_{1,2}=T_1, \,\,\,\,\,\,\,\,\,\,\,\,\,\,
C_{2,1}=2T_1, \,\,\,\,\,\,\,\,\,\,\,\,\,\,
C_{2,2}=0,
\\
A_{1,2}
&
=
-3F_{0,3}T_2-F_{1,2}T_1,
\\
D_{1,1}
&
=
2F_{1,2}T_2
+
2F_{2,1}T_1
+
A_{1,1}
+
A_{2,2},
\\
D_{1,2}
&
=
F_{2,1}T_2
+
3F_{3,0}T_1
+
A_{2,1},
\\
D_{2,1}
&
=
-6F_{0,3}T_2-2F_{1,2}T_1+2G_{1,2}T_2+2G_{2,1}T_1,
\\
D_{2,2}
&
=
G_{2,1}T_2+3G_{3,0}T_1+2A_{1,1},
\\
G_{0,3}
&
=
0,
\\
G_{1,2}
&
=
0,
\\
B_{1,1}
&
=
-G_{2,2}T_2-\tfrac{3}{2}G_{3,1}T_1,
\\
B_{1,2}
&
=
-\tfrac{1}{2}G_{3,1}T_2-2G_{4,0}T_1,
\\
B_{2,1}
&
=
-2F_{1,2}A_{2,1}+\tfrac{1}{2}G_{3,1}T_2-2F_{2,2}-3F_{3,1}T_1+2G_{4,0}T_1,
\\
B_{2,2}
&
=
-F_{3,1}T_2-4F_{4,0}T_1.
\endaligned
\]

Then, we open the branch $\green{\bf 2e3c}$:
\[
\aligned
F_{0,3}
&
=
0
=
R_{0,3},
\\
F_{1,2}
&
=
0
=
R_{1,2}.
\endaligned
\]
At order 3, the tangency equations give:
\[
F_{1,3} = 0, \,\,\,\,\,\,
F_{2,2} = \tfrac{3}{4}G_{3,1}, \,\,\,\,\,\,
G_{1,3} = 0, \,\,\,\,\,\,
G_{2,2} = 0, \,\,\,\,\,\,
F_{0,4} = 0, \,\,\,\,\,\,
G_{0,4} = 0.
\]

Next, at order $4$, the fundamental equations give:
\[
\aligned
0
&
\underset{\green{\bf S}}{\overset{\green{\bf 3,1}}{\,\,=\,\,}}
S_{3,1}a_{1,1}^3a_{2,2}-G_{3,1}a_{1,1}^2,
\\
0
&
\underset{\green{\bf R}}{\overset{\green{\bf 3,1}}{\,\,=\,\,}}
-a_{1,1}a_{2,1}G_{3,1}-a_{1,1}a_{2,2}F_{3,1}+R_{3,1}a_{1,1}^3a_{2,2}+\tfrac{3}{2}S_{3,1}a_{1,1}^2a_{2,2}a_{2,1},
\\
0
&
\underset{\green{\bf S}}{\overset{\green{\bf 4,0}}{\,\,=\,\,}}
S_{3,1}a_{1,1}^3a_{2,1}+S_{4,0}a_{1,1}^4-G_{4,0}a_{1,1}^2,
\\
0
&
\underset{\green{\bf R}}{\overset{\green{\bf 4,0}}{\,\,=\,\,}}
-a_{1,1}a_{2,1}G_{4,0}-a_{1,1}a_{2,2}F_{4,0}+R_{4,0}a_{1,1}^4+R_{3,1}a_{1,1}^3a_{2,1}+\tfrac{3}{4}S_{3,1}a_{1,1}^2a_{2,1}^2,
\endaligned
\]
providing at order $4$ the linear representation:
\[
\left(\!
\begin{array}{c}
S_{3,1}
\\
R_{3,1}
\\
S_{4,0}
\\
R_{4,0}
\end{array}
\!\right)
\,=\,
\left(\!
\def\arraystretch{1.25}
\begin{array}{cccc}
\tfrac{1}{a_{1,1}a_{2,2}} & 0 & 0 & 0
\\
-\tfrac{a_{2,1}}{2a_{1,1}^2a_{2,2}} & \tfrac{1}{a_{1,1}^2} & 0 & 0
\\
-\tfrac{a_{2,1}}{a_{1,1}^2a_{2,2}} & 0 & \tfrac{1}{a_{1,1}^2} & 0
\\
-\tfrac{a_{2,1}^2}{4a_{1,1}^3a_{2,2}} & -\tfrac{a_{2,1}}{a_{1,1}^3} & \tfrac{a_{2,1}}{a_{1,1}^3} & \tfrac{1}{a_{1,1}^3}
\end{array}
\!\right)\,
\left(\!
\begin{array}{c}
G_{3,1}
\\
F_{1,3}
\\
G_{4,0}
\\
F_{4,0}
\end{array}
\!\right).
\]
\begin{Lemma}
The subgroup $G_{\text{\rm stab}}^{4}\subset G_{\text{\rm stab}}^{3}$ which sends:
\[
\aligned
u
&
=
xy+F_{4,0}x^4+F_{3,1}x^3y+\tfrac{3}{4}G_{3,1}x^2y^2
+
{\rm O}_{x,y}(5),
\\
v
&
=
x^2+G_{3,1}x^3y+G_{4,0}x^4
+
{\rm O}_{x,y}(5),
\endaligned
\]
to:
\[
\aligned
r
&
=
pq+R_{4,0}p^4+R_{3,1}p^3q+\tfrac{3}{4}S_{3,1}p^2q^2
+
{\rm O}_{p,q}(5),
\\
s
&
=
p^2+S_{3,1}p^3q+S_{4,0}p^4
+
{\rm O}_{p,q}(5),
\endaligned
\]
consists of matrices:
\[
\begin{pmatrix}
a_{1,1} & 0 & 0 & 0 \\
a_{2,1} & a_{2,2} & 0 & 0 \\
0 & 0 & a_{1,1}a_{2,2}& a_{1,1}a_{2,1}\\
0 & 0 & 0 & a_{1,1}^2
\end{pmatrix}
\in G_{\text{\rm stab}}^{4}.
\]
\end{Lemma}

Then, we open the branch $\green{\bf 2e3c4e}$:
\[
\aligned
G_{3,1}
&
=
0
=
S_{3,1},
\\ 
F_{3,1}
&
=
-1
=
R_{3,1},
\\
F_{4,0}
&
=
0
=
R_{4,0},
\\
G_{4,0}
&
\neq -1
\neq
S_{4,0}.
\endaligned
\]
At order $4$, the fundamental equations provide:
\[
a_{1,1}
=
1,
\,\,\,\,
a_{2,1}
=
0,
\]
and the associated tangency equations, at order $4$, to:
\[
\aligned
u
&
=
xy-x^3y
+
F_{5,0}x^5
+
F_{4,1}x^4y
+
F_{3,2}x^3y^2
+
F_{2,3}x^2y^3
+
F_{1,4}xy^4
+
F_{0,5}y^5
+
{\rm O}_{x,y}(6),
\\
v
&
=
x^2+G_{4,0}x^4
+
+
G_{5,0}x^5
+
G_{4,1}x^4y
+
G_{3,2}x^3y^2
+
G_{2,3}x^2y^3
+
G_{1,4}xy^4
+
G_{0,5}y^5
+
{\rm O}_{x,y}(6),
\endaligned
\]
give:
\[
\aligned
A_{2,1}
&
=
5\frac{F_{5,0}}{G_{4,0}+1}T_1+\frac{F_{4,1}}{G_{4,0}+1}T_2,
\\
A_{1,1}
&
=
2F_{4,1}T_1+F_{3,2}T_2,
\\
F_{1,4}
&
=
F_{2,3}
=
F_{3,2}
=
G_{1,4}
=
G_{2,3}
=
G_{3,2}
=
G_{4,1}
=
F_{0,5}
=
G_{0,5} 
= 0,
\\
G_{5,0}
&
=
-\tfrac{4}{5}F_{4,1}G_{4,0}.
\endaligned
\]

Next, at order $5$, the fundamental equations give the linear representation:
\[
R_{5,0}
=
a_{2,2}F_{5,0}.
\]
\begin{Lemma}
The subgroup $G_{\text{\rm stab}}^{5}\subset G_{\text{\rm stab}}^{4}$ which sends:
\[
\aligned
u
&
=
xy
-
x^3y
+
F_{5,0}x^5
+
F_{4,1}x^4y
+
{\rm O}_{x,y}(6),
\\
v
&
=
x^2
+
G_{4,0}x^4
-\tfrac{4}{5}F_{4,1}G_{4,0}x^5
+
{\rm O}_{x,y}(6),
\endaligned
\]
to:
\[
\aligned
r
&
=
pq
-
p^3q
+
R_{5,0}p^5
+
R_{4,1}p^4q
+
{\rm O}_{p,q}(6),
\\
s
&
=
p^2
+
S_{4,0}p^4
-\tfrac{4}{5}R_{4,1}S_{4,0}p^5
+
{\rm O}_{p,q}(6),
\endaligned
\]
consists of matrices:
\[
\begin{pmatrix}
1 & 0 & 0 & 0 \\
0 & a_{2,2} & 0 & 0 \\
0 & 0 & a_{2,2}& 0\\
0 & 0 & 0 & 1
\end{pmatrix}
\in G_{\text{\rm stab}}^{5}.
\]
\end{Lemma}

Then, we open the branch $\green{\bf 2e3c4e5b}$:
\[
F_{5,0}
=
0
=
R_{5,0}.
\]
The tangency equations at order $5$, to:
\[
\aligned
u
&
=
xy
-
x^3y
+
F_{4,1}x^4y
+
F_{6,0}x^6
+
F_{5,1}x^5y
+
F_{4,2}x^4y^2
+
F_{3,3}x^3y^3
+
F_{2,4}x^2y^4
+
F_{1,5}xy^5
+
\\
&
\,\,\,\,\,\,
+
F_{0,6}y^6
+
{\rm O}_{x,y}(7),
\\
v
&
=
x^2
+
G_{4,0}x^4
-
\tfrac{4}{5}F_{4,1}G_{4,0}x^5
+
G_{6,0}x^6
+
G_{5,1}x^5y
+
G_{4,2}x^4y^2
+
G_{3,3}x^3y^3
+
G_{2,4}x^2y^4
+
\\
&
\,\,\,\,\,\,
+
G_{1,5}xy^5
+
G_{0,6}y^6
+
{\rm O}_{x,y}(7),
\endaligned
\]
give:
\[
\aligned
F_{1,5} 
&
=
F_{2,4}
=
F_{3,3}
=
F_{4,2}
=
F_{6,0}
= 0,
\\
F_{5,1}
&
= -\tfrac{6}{5}F_{4,1}^2+\tfrac{2}{5}G_{4,0}^2-\tfrac{2}{5}G_{4,0}+\tfrac{6}{5},
\\
G_{1,5}
&
=
G_{2,4}
=
G_{3,3}
=
G_{4,2}
=
G_{5,1}
=
F_{0,6}
=
G_{0,6} 
=
0,
\\
G_{6,0} & = 2G_{4,0}^2+\tfrac{4}{5}F_{4,1}^2G_{4,0},
\endaligned
\]
and the last equation $\underset{\green{\bf F}}{\overset{\green{\bf 5,0}}{\,\,=\,\,}}$ provides a single algebraic relation between $G_{4,0},F_{4,1}$:
\[
0
=
\mathbb{B}
:=
2F_{4,1}^2G_{4,0}-2G_{4,0}^3+F_{4,1}^2-5G_{4,0}^2-4G_{4,0}-1,
\]
which is the {\em core algebraic equation}.

\begin{Remark}
In some other branches, we find several {\em core algebraic equations}. In such branches, we define $I_{\rm ca}$ the ideal generated by the associated polynomials of the core algebraic equations, and we denote by $\mathbb{B}_i$ its generators.
\end{Remark}

Next, at orders $n\geqslant 6$, the fundamental equations no longer provide possible normalizations, because the tangency equations are all of the form:
\[
\aligned
0
&
\underset{\green{\bf F}}{\overset{\green{\bf n,0}}{\,\,=\,\,}}
\big(
\alpha_{n,0}F_{n+1,0}
\big)\,
T_1
+
\big(
\gamma_{n,0}F_{n,1}+\delta_{n,0}
\big)\,
T_2,
\\
0
&
\underset{\green{\bf G}}{\overset{\green{\bf n,0}}{\,\,=\,\,}}
\big(
\epsilon_{n,0}G_{n+1,0}+\zeta_{n,0}
\big)\,
T_1
+
\big(
\eta_{n,0}G_{n,1}
\big)\,
T_2,
\endaligned
\]
with 
$\alpha_{n,0}, \,
\gamma_{n,0}, \,
\epsilon_{n,0}, \,
\eta_{n,0} \in\R^*$
and
$\delta_{n,0}, \,
\zeta_{n,0}
\in \R[G_{4,0},F_{4,1}]$.

Moreover, for $i+j=n$ and $k+l=n$, except $(n,0)$, 
the equations are all of the form:
\[
\aligned
0
&
\underset{\green{\bf F}}{\overset{\green{\bf i,j}}{\,\,=\,\,}}
\big(
\alpha_{i,j}F_{i+1,j}
\big)\,
T_1
+
\big(
\gamma_{i,j}F_{i,j+1}
\big)\,
T_2,
\\
0
&
\underset{\green{\bf G}}{\overset{\green{\bf k,l}}{\,\,=\,\,}}
\big(
\epsilon_{k,l}G_{k+1,l}
\big)\,
T_1
+
\big(
\eta_{k,l}G_{k,l+1}
\big)\,
T_2,
\endaligned
\]
with 
$
\alpha_{i,j}, \,
\gamma_{i,j}, \,
\epsilon_{i,j}, \,
\eta_{i,j} \in\R^*$.

Therefore, all coefficients $F_{i,j}, G_{k,l}$ of order $n+1$ are determined at all orders $n\geqslant 6$.

\begin{Lemma}
For $n\geqslant 5$, the subgroup $G_{\text{\rm stab}}^{n+1}\subset G_{\text{\rm stab}}^{n}$ which sends:
\[
\aligned
u
&
=
xy
-
x^3y
+
F_{4,1}x^4y
+
\sum_{i=5}^n
F_{i,1}x^iy
+
{\rm O}_{x,y}(n+2),
\\
v
&
=
x^2
+
G_{4,0}x^4
-\tfrac{4}{5}F_{4,1}G_{4,0}x^5
+
\sum_{j=6}^{n+1}
G_{j,0}x^j
+
{\rm O}_{x,y}(n+2),
\endaligned
\]
to:
\[
\aligned
r
&
=
pq
-
p^3q
+
R_{4,1}p^4q
+
\sum_{i=5}^n
R_{i,1}p^iq
+
{\rm O}_{p,q}(n+2),
\\
s
&
=
p^2
+
S_{4,0}p^4
-
\tfrac{4}{5}R_{4,1}S_{4,0}p^5
+
\sum_{j=6}^{n+1}
S_{j,0}p^j
+
{\rm O}_{p,q}(n+2),
\endaligned
\]
with $F_{i,1}G_{j,0}\in \R[G_{4,0},F_{4,1}]$ and $R_{i,1},S_{j,0}\in \R[S_{4,0},R_{4,1}]$ consists of matrices:
\[
\begin{pmatrix}
1 & 0 & 0 & 0 \\
0 & a_{2,2} & 0 & 0 \\
0 & 0 & a_{2,2}& 0\\
0 & 0 & 0 & 1
\end{pmatrix}
\in G_{\text{\rm stab}}^{n+1}.
\]
\end{Lemma}

Finally, we obtain the three generators of the Lie algebra:
\[
\aligned
e_1
&
:=
L\Big\vert_{T_1=1,T_2=0,A_{2,2}=0},
\\
e_2
&
:=
L\Big\vert_{T_1=0,T_2=1,A_{2,2}=0},
\\
e_3
&
:=
L\Big\vert_{T_1=0,T_2=0,A_{2,2}=1},
\endaligned
\]
which satisfy:
\[
[e_1,e_2]
=
\frac{F_{4,1}}{G_{4,0}+1} \, e_2,
\,\,\,\,\,\,\,\,\,\,\,\,\,\,\,\,
[e_2,e_3]
=
e_2.
\]

\begin{Theorem}
In the branch $\green{\bf 2e3c4e5b}$, 
every affinely homogeneous surface is, in a unique way, equivalent to:
\[
\left\{
\aligned
u
&
=
xy
-
x^3y
+
F_{4,1}x^4y
+
\sum_{i=5}^n
F_{i,1}x^iy
+
\cdots,
\\
v
&
=
x^2
+
G_{4,0}x^4
-\tfrac{4}{5}F_{4,1}G_{4,0}x^5
+
\sum_{j=6}^{n+1}
G_{j,0}x^j
+
\cdots,
\endaligned\right.
\]
with $F_{i,1},G_{j,0}\in \R[G_{4,0},F_{4,1}]$, and with $3$-dimensional Lie algebra:
\[
\aligned
e_1
&
\,:=\,
2F_{4,1}u\partial_u
+
4F_{4,1}v\partial_v
+
2F_{4,1}x\partial_x
-
2G_{4,0}v\partial_x
+
2G_{4,0}u\partial_y
+
y\partial_u
+
2u\partial_v
+
3u\partial_y
+
\partial_x
, 
\\
e_2
&
\,:=\,
\tfrac{1}{G_{4,0}+1}(F_{4,1}x+G_{4,0}v+G_{4,0}+v+1)\partial_y+
\tfrac{1}{G_{4,0}+1}(F_{4,1}v+G_{4,0}x+x)\partial_u
, 
\\
e_3
&
\,:=\,
\tfrac{1}{G_{4,0}+1}(G_{4,0}y+y)\partial_y+\tfrac{1}{G_{4,0}+1}(G_{4,0}u+u)\partial_u,
\endaligned
\]
and Lie structure:
\[
\footnotesize
\def\arraystretch{1.25}
\begin{array}{c|ccc}
{} & e_1 & e_2 & e_3
\\
\hline
e_1 & 0 & \tfrac{F_{4,1}}{G_{4,0}+1}e_2
\\
e_2 & -\tfrac{F_{4,1}}{G_{4,0}+1}e_2 & 0 & e_2
\\
e_3 & -\tfrac{F_{4,1}}{G_{4,0}+1}e_2 & 0 & e_2
\end{array}
\]
with $G_{4,0},F_{4,1}$ satisfying:
\[
0 = \mathbb{B} := 2F_{4,1}^2G_{4,0}-2G_{4,0}^3+F_{4,1}^2-5G_{4,0}^2-4G_{4,0}-1.
\]
\end{Theorem}

As explained in the introduction,
when reaching the terminal leaf of a branch,
we stop the process 
either once we obtain
a single homogeneous model, 
or once we obtain 
a moduli space of homogeneous models
parametrized by a certain algebraic variety.

For instance here, the ambient space is $\R^2$ with
coordinates $\big(G_{4,0}, F_{4,1}\big)$,
and the moduli space is the algebraic variety (hypersurface):
\[
Z_\B
\,:=\,
\big\{ 
\B 
=
0 
\big\}.
\]
The whole process of equivalence 
tells us that:

\medskip\noindent$\bullet$\,
every homogeneous model of $Z_\B$
is {\em inequivalent} to any
other homogeneous model found in any other terminal leaf
of the branching diagram;

\medskip\noindent$\bullet$\,
to each point $p \in Z_\B$ there corresponds a homogeneous
model;

\medskip\noindent$\bullet$\,
but still, two homogeneous models $p \neq q \in \Z_\B$ 
can be finitely equivalent one to another.

\medskip

Such a residual finite equivalence may be 
understood\big/set up explicitly,
by analyzing the remaining fundamental 
equations~{\eqref{eqR-eqS-S2-R4}} 
at the terminal leaf of this branch
\green{\bf 2e3c4e5b}.

However, because the number of existing branches 
that we have fully analyzed for {\em continuous} group equivalences,
is already quite large,
we will not endeavour any complete {\em finite equivalence} analysis.
In fact, several authors in Lie-Cartan theory 
usually disregard equivalences up to finite groups.
And in addition, during the whole process of branch creation,
we always considered {\em only the identity component}
of the concerned matrix stability groups, 
so that we did not consider the complete
finite equivalence in the complete matrix stability groups.

Anyhow, just for this example
\green{\bf 2e3c4e5b}
of terminal leaf, let us explain how one could go
a bit further by
decomposing the moduli space algebraic variety $Z_\B$
into convenient strata. By chance here, $\B$ factorizes:
\[
0
\,=\,
\B
\,=\,
\big(2\,G_{4,0}+1\big)\,
\big(F_{4,1}+1+G_{4,0}\big)\,
\big(F_{4,1}-1-G_{4,0}\big),
\]
and hence, still with $1+G_{4,0} \neq 0$, 
there are three "sub-cases":

\medskip\noindent{\bf (1)}\,
$G_{4,0} := -\,\frac{1}{2}$,

\medskip\noindent{\bf (2)}\,
$G_{4,0} := -\, 1 - F_{4,1}$, with $F_{4,1} \neq 0$,

\medskip\noindent{\bf (3)}\,
$G_{4,0} := -\, 1 +  F_{4,1}$, with $F_{4,1} \neq 0$,

\medskip\noindent
with corresponding vector fields, for {\bf (1)}:
\[
\footnotesize
\aligned
e_1
&
\,=\,
\big(
1+2F_{4,1}x+v
\big)\,\partial_x
+
\big(
2u
\big)\,\partial_y
+
\big(
y+2F_{4,1}u
\big)\,\partial_u
+
\big(2x+4F_{4,1}v\big)\,\partial_v,
\\
e_2
&
\,=\,
\big(
1+2F_{4,1}x+v
\big)\,\partial_y
+
\big(
x+2F_{4,1}v
\big)\,\partial_u,
\\
e_3
&
\,=\,
y\,\partial_y
+
u\,\partial_u,
\endaligned
\]
for {\bf (2)}:
\[
\footnotesize
\aligned
e_1
&
\,=\,
\big(
1+2F_{4,1}x+2v+2F_{4,1}v
\big)\,\partial_x
+
\big(
u-2F_{4,1}v
\big)\,\partial_y
+
\big(
y+2F_{4,1}u
\big)\,\partial_u
+
\big(2x+4F_{4,1}v\big)\,\partial_v,
\\
e_2
&
\,=\,
\big(
1-x+v
\big)\,\partial_y
+
\big(
x-v
\big)\,\partial_u,
\\
e_3
&
\,=\,
y\,\partial_y
+
u\,\partial_u,
\endaligned
\]
for {\bf (3)}:
\[
\footnotesize
\aligned
e_1
&
\,=\,
\big(
1+2F_{4,1}x+2v-2F_{4,1}v
\big)\,\partial_x
+
\big(
u+2F_{4,1}v
\big)\,\partial_y
+
\big(
y+2F_{4,1}u
\big)\,\partial_u
+
\big(2x+4F_{4,1}v\big)\,\partial_v,
\\
e_2
&
\,=\,
\big(
1+x+v
\big)\,\partial_y
+
\big(
x-v
\big)\,\partial_u,
\\
e_3
&
\,=\,
y\,\partial_y
+
u\,\partial_u,
\endaligned
\]
and with Lie structures:
\[
\footnotesize
\aligned
{\bf (1)}\,\,\,
\def\arraystretch{1.25}
\begin{array}{c|ccc}
{} & e_1 & e_2 & e_3
\\
\hline
e_1 & 0 & 2F_{4,1}e_2 & 0
\\
e_2 & -2F_{4,1}e_2 & 0 & e_2
\\
e_3 & 0 & -e_2 & 0
\end{array}
\ \ \ \ \ \ \ \ \ \ \ \ \ \ \
{\bf (2)}\,\,\,
\def\arraystretch{1.25}
\begin{array}{c|ccc}
{} & e_1 & e_2 & e_3
\\
\hline
e_1 & 0 & -e_2 & 0
\\
e_2 & e_2 & 0 & e_2
\\
e_3 & 0 & -e_2 & 0
\end{array}
\ \ \ \ \ \ \ \ \ \ \ \ \ \ \
{\bf (3)}\,\,\,
\def\arraystretch{1.25}
\begin{array}{c|ccc}
{} & e_1 & e_2 & e_3
\\
\hline
e_1 & 0 & e_2 & 0
\\
e_2 & -e_2 & 0 & e_2
\\
e_3 & 0 & -e_2 & 0
\end{array}
\endaligned
\]

Similar analyses\big/decompositions of all the algebraic
varieties appearing in 
Section~{\ref{2b-models}} $\to$ {\ref{2g-models}} 
can be made,
in order to simplify\big/compactify the presentations
of the concerned Lie algebras, to the price of creating
"sub-cases". But such "subcases" are {\em not} subbranches,
because they do not come from 
orbit decompositions\big/transversals in some
linear representation at some jet order.

In conclusion, we would like to emphasize once again that, 
{\em reaching core algebraic varieties that are moduli spaces
of homogeneous models is what the equivalence method does}.

%%%%%%%%%%%%%%%%%%%%%%%%%%%%%%%%%%%%%%%%%%%%%%%%%%%%%%%%%%%%%%%%%%%%%%
%%%%%%%%%%%%%%%%%%%%%%%%%%%%%%%%%%%%%%%%%%%%%%%%%%%%%%%%%%%%%%%%%%%%%%
%%%%%%%%%%%%%%%%%%%%%%%%%%%%%%%%%%%%%%%%%%%%%%%%%%%%%%%%%%%%%%%%%%%%%%

%%%%%%%%%%%%%%%%%%%%%%%%%%%%%%%%%%%%%%%%%%%%%%%%%%%%%%%%%%%%%%%%%%%%%%
\SectionHead{2b Models}
{2b-models}
%%%%%%%%%%%%%%%%%%%%%%%%%%%%%%%%%%%%%%%%%%%%%%%%%%%%%%%%%%%%%%%%%%%%%%

\[
\text{\bf Model 2b3a4a}\ \ \ \ \
\left\{
\aligned
u
&
\,=\,
x^2,
\\
v
&
\,=\,
0,
\endaligned\right.
\]
%%%%%%%%%%%%%%%%%%%%%%%%%%%%%%%%%%%%%%%%%%%%%%%%%%%%%%%%%%%%%%%%%%%%%%
\[
\def\arraystretch{1.25}
\begin{array}{llll}
e_1
\,:=\,
\partial_x+2\,x\,\partial_u, &
e_2
\,:=\,
\partial_y, &
& 
\\
e_3
\,:=\,
x\,\partial_x+2\,u\,\partial_u, &
e_4
\,:=\,
x\,\partial_y, &
e_5
\,:=\,
y\,\partial_y, &
e_6
\,:=\,
v\,\partial_x,
\\
e_7
\,:=\,
u\,\partial_y, &
e_8
\,:=\,
v\,\partial_y, &
e_9
\,:=\,
v\,\partial_u, &
e_{10}
\,:=\,
v\,\partial_v,
\end{array}
\]
%%%%%%%%%%%%%%%%%%%%%%%%%%%%%%%%%%%%%%%%%%%%%%%%%%%%%%%%%%%%%%%%%%%%%%
\[
\footnotesize
\def\arraystretch{1.25}
\begin{array}{c|cccccccccc}
{} & e_1 & e_2 & e_3 & e_4 & e_5 & e_6 & e_7 & e_8 & e_9 & e_{10}
\\
\hline
e_1 & 
0 & 0 & e_1 & e_2 & 0 & -\,2\,e_9 & 2\,e_4 & 0 & 0 & 0
\\
e_2 &
0 & 0 & 0 & 0 & e_2 & 0 & 0 & 0 & 0 & 0
\\
e_3 &
-\,e_1 & 0 & 0 & e_4 & 0 & -\,e_6 & 2\,e_7 & 0 & -\,2\,e_9 & 0
\\
e_4 &
-\,e_2 & 0 & -\,e_4 & 0 & e_4 & -\,e_8 & 0 & 0 & 0 & 0
\\
e_5 &
0 & -\,e_2 & 0 & -\,e_4 & 0 & 0 & -\,e_7 & -\,e_8 & 0 & 0
\\
e_6 &
2\,e_9 & 0 & e_6 & e_8 & 0 & 0 & 0 & 0 & 0 & -\,e_6
\\
e_7 &
-\,2\,e_4 & 0 & -\,2\,e_7 & 0 & e_7 & 0 & 0 & 0 & -\,e_8 & 0
\\
e_8 &
0 & 0 & 0 & 0 & e_8 & 0 & 0 & 0 & 0 & -\,e_8
\\
e_9 &
0 & 0 & 2\,e_9 & 0 & 0 & 0 & e_8 & 0 & 0 & -\,e_9
\\
e_{10} &
0 & 0 & 0 & 0 & 0 & e_6 & 0 & e_8 & e_9 & 0
\end{array}
\]

%%%%%%%%%%%%%%%%%%%%%%%%%%%%%%%%%%%%%%%%%%%%%%%%%%%%%%%%%%%%%%%%%%%%%%
%%%%%%%%%%%%%%%%%%%%%%%%%%%%%%%%%%%%%%%%%%%%%%%%%%%%%%%%%%%%%%%%%%%%%%
%%%%%%%%%%%%%%%%%%%%%%%%%%%%%%%%%%%%%%%%%%%%%%%%%%%%%%%%%%%%%%%%%%%%%%

\[
\text{\bf Model 2b3a4b}
\ \ \ \ \
\left\{
\aligned
u
&
\,=\,
x^2+x^4+F_{5,0}\,x^5
+
\big(
2+\tfrac{5}{4}\,F_{5,0}^2
\big)\,x^6
+
\big(
\tfrac{34}{7}\,F_{5,0}
+
\tfrac{25}{14}\,F_{5,0}^3
\big)\,x^7
+\cdots,
\\
v
&
\,=\,
0,
\endaligned\right.
\]
%%%%%%%%%%%%%%%%%%%%%%%%%%%%%%%%%%%%%%%%%%%%%%%%%%%%%%%%%%%%%%%%%%%%%%
\[
\def\arraystretch{1.25}
\begin{array}{llll}
e_1
\,:=\,
\big(1-\tfrac{5}{2}\,F_{5,0}\,x+2\,u\big)\,\partial_x
+
\big(2\,x-5\,F_{5,0}\,u\big)\,\partial_u, &
e_2
\,:=\,
\partial_y, & 
\\
e_3
\,:=\,
x\,\partial_y, &
e_4
\,:=\,
y\,\partial_y, &
e_5
\,:=\,
v\,\partial_x, &
e_6
\,:=\,
u\,\partial_y,
\\
e_7
\,:=\,
v\,\partial_y, &
e_8
\,:=\,
v\,\partial_u, &
e_9
\,:=\,
v\,\partial_v,
\end{array}
\]
%%%%%%%%%%%%%%%%%%%%%%%%%%%%%%%%%%%%%%%%%%%%%%%%%%%%%%%%%%%%%%%%%%%%%%
\[
\tiny
\def\arraystretch{1.25}
\begin{array}{c|ccccccccc}
{} & e_1 & e_2 & e_3 & e_4 & e_5 & e_6 & e_7 & e_8 & e_9 
\\
\hline
e_1 & 
0 & 0 & e_2-\frac{5}{2}\,F_{5,0}\,e_3-2\,e_6 & 0 & 
\frac{5}{2}\,F_{5,0}\,e_5-2\,e_8 & 2\,e_3-5\,F_{5,0}\,e_6 &
0 & 2\,e_5+5\,F_{5,0}\,e_8 & 0
\\
e_2 &
0 & 0 & 0 & e_2 & 0 & 0 & 0 & 0 & 0 
\\
e_3 &
-\,e_2+\frac{5}{2}\,F_{5,0}\,e_3+2\,e_6 & 0 & 0 & e_3 & 
-\,e_7 & 0 & 0 & 0 & 0 
\\
e_4 &
0 & -\,e_2 & -\,e_3 & 0 & 0 & -\,e_6 & -\,e_7 & 0 & 0 
\\
e_5 &
-\,\frac{5}{2}\,F_{5,0}\,e_5+2\,e_8 & 0 & e_7 & 0 & 0 & 0 & 0 &
0 & -\,e_5
\\
e_6 &
-\,2\,e_3+5\,F_{5,0}\,e_6 & 0 & 0 & e_6 & 0 & 0 & 0 & -\,e_7 & 0
\\
e_7 &
0 & 0 & 0 & e_7 & 0 & 0 & 0 & 0 & -\,e_7
\\
e_8 &
-\,2\,e_5-5\,F_{5,0}\,e_8 & 0 & 0 & 0 & 0 & e_7 & 0 & 0 & -\,e_8
\\
e_9 &
0 & 0 & 0 & 0 & e_5 & 0 & e_7 & e_8 & 0
\end{array}
\]

%%%%%%%%%%%%%%%%%%%%%%%%%%%%%%%%%%%%%%%%%%%%%%%%%%%%%%%%%%%%%%%%%%%%%%
%%%%%%%%%%%%%%%%%%%%%%%%%%%%%%%%%%%%%%%%%%%%%%%%%%%%%%%%%%%%%%%%%%%%%%
%%%%%%%%%%%%%%%%%%%%%%%%%%%%%%%%%%%%%%%%%%%%%%%%%%%%%%%%%%%%%%%%%%%%%%

\[
\text{\bf Model 2b3a4c}
\ \ \ \ \
\left\{
\aligned
u
&
\,=\,
x^2-x^4+F_{5,0}\,x^5
+
\big(
2-\tfrac{5}{4}\,F_{5,0}^2
\big)\,x^6
+
\big(
-\,\tfrac{34}{7}\,F_{5,0}
+
\tfrac{25}{14}\,F_{5,0}^3
\big)\,x^7
+\cdots,
\\
v
&
\,=\,
0,
\endaligned\right.
\]
for any value of $F_{5,0}$,
%%%%%%%%%%%%%%%%%%%%%%%%%%%%%%%%%%%%%%%%%%%%%%%%%%%%%%%%%%%%%%%%%%%%%%
\[
\def\arraystretch{1.25}
\begin{array}{llll}
e_1
\,:=\,
\big(1+\tfrac{5}{2}\,F_{5,0}\,x+2\,u\big)\,\partial_x
+
\big(2\,x+5\,F_{5,0}\,u\big)\,\partial_u, &
e_2
\,:=\,
\partial_y, & 
\\
e_3
\,:=\,
x\,\partial_y, &
e_4
\,:=\,
y\,\partial_y, &
e_5
\,:=\,
v\,\partial_x, &
e_6
\,:=\,
u\,\partial_y,
\\
e_7
\,:=\,
v\,\partial_y, &
e_8
\,:=\,
v\,\partial_u, &
e_9
\,:=\,
v\,\partial_v,
\end{array}
\]
%%%%%%%%%%%%%%%%%%%%%%%%%%%%%%%%%%%%%%%%%%%%%%%%%%%%%%%%%%%%%%%%%%%%%%
\[
\tiny
\def\arraystretch{1.25}
\begin{array}{c|ccccccccc}
{} & e_1 & e_2 & e_3 & e_4 & e_5 & e_6 & e_7 & e_8 & e_9 
\\
\hline
e_1 & 
0 & 0 & e_2+\frac{5}{2}\,F_{5,0}\,e_3+2\,e_6 & 0 & 
-\,\frac{5}{2}\,F_{5,0}\,e_5-2\,e_8 & 2\,e_3+5\,F_{5,0}\,e_6 &
0 & -\,2\,e_5-5\,F_{5,0}\,e_8 & 0
\\
e_2 &
0 & 0 & 0 & e_2 & 0 & 0 & 0 & 0 & 0 
\\
e_3 &
-\,e_2-\frac{5}{2}\,F_{5,0}\,e_3-2\,e_6 & 0 & 0 & e_3 & 
-\,e_7 & 0 & 0 & 0 & 0 
\\
e_4 &
0 & -\,e_2 & -\,e_3 & 0 & 0 & -\,e_6 & -\,e_7 & 0 & 0 
\\
e_5 &
\frac{5}{2}\,F_{5,0}\,e_5+2\,e_8 & 0 & e_7 & 0 & 0 & 0 & 0 &
0 & -\,e_5
\\
e_6 &
-\,2\,e_3-5\,F_{5,0}\,e_6 & 0 & 0 & e_6 & 0 & 0 & 0 & -\,e_7 & 0
\\
e_7 &
0 & 0 & 0 & e_7 & 0 & 0 & 0 & 0 & -\,e_7
\\
e_8 &
2\,e_5+5\,F_{5,0}\,e_8 & 0 & 0 & 0 & 0 & e_7 & 0 & 0 & -\,e_8
\\
e_9 &
0 & 0 & 0 & 0 & e_5 & 0 & e_7 & e_8 & 0
\end{array}
\]

%%%%%%%%%%%%%%%%%%%%%%%%%%%%%%%%%%%%%%%%%%%%%%%%%%%%%%%%%%%%%%%%%%%%%%
%%%%%%%%%%%%%%%%%%%%%%%%%%%%%%%%%%%%%%%%%%%%%%%%%%%%%%%%%%%%%%%%%%%%%%
%%%%%%%%%%%%%%%%%%%%%%%%%%%%%%%%%%%%%%%%%%%%%%%%%%%%%%%%%%%%%%%%%%%%%%

\[
\text{\bf Model 2b3b4a5a}
\ \ \ \ \
\left\{
\aligned
u
&
\,=\,
x^2+x^2y+x^2y^2+x^2y^3+x^2y^4+x^2y^5+x^2y^6+\cdots,
\\
v
&
\,=\,
0,
\endaligned\right.
\]
%%%%%%%%%%%%%%%%%%%%%%%%%%%%%%%%%%%%%%%%%%%%%%%%%%%%%%%%%%%%%%%%%%%%%%
\[
\def\arraystretch{1.25}
\begin{array}{llll}
e_1
\,:=\,
(1-y)\,\partial_x+2\,x\,\partial_u, &
e_2
\,:=\,
(1-y)\,\partial_y+u\,\partial_u, &
& 
\\
e_3
\,:=\,
x\,\partial_x+2\,u\,\partial_u, &
e_4
\,:=\,
-\,\frac{1}{2}\,u\,\partial_x+x\,\partial_y, &
e_5
\,:=\,
v\,\partial_x, &
e_6
\,:=\,
v\,\partial_y,
\\
e_7
\,:=\,
v\,\partial_u, &
e_8
\,:=\,
v\,\partial_v,
\end{array}
\]
%%%%%%%%%%%%%%%%%%%%%%%%%%%%%%%%%%%%%%%%%%%%%%%%%%%%%%%%%%%%%%%%%%%%%%
\[
\footnotesize
\def\arraystretch{1.25}
\begin{array}{c|cccccccc}
{} & e_1 & e_2 & e_3 & e_4 & e_5 & e_6 & e_7 & e_8 
\\
\hline
e_1 & 
0 & e_1 & e_1 & e_2 & -\,2\,e_7 & e_5 & 0 & 0
\\
e_2 &
-\,e_1 & 0 & 0 & e_4 & 0 & e_6 & -\,e_7 & 0
\\
e_3 &
-\,e_1 & 0 & 0 & e_4 & -\,e_5 & 0 & -\,2\,e_7 & 0
\\
e_4 &
-\,e_2 & -\,e_4 & -\,e_4 & 0 & -\,e_6 & 0 & \frac{1}{2}\,e_5 & 0
\\
e_5 &
2\,e_7 & 0 & e_5 & e_6 & 0 & 0 & 0 & -\,e_5
\\
e_6 &
-\,e_5 & -\,e_6 & 0 & 0 & 0 & 0 & 0 & -\,e_6
\\
e_7 &
0 & e_7 & 2\,e_7 & -\,\frac{1}{2}\,e_5 & 0 & 0 & 0 & -\,e_7
\\
e_8 &
0 & 0 & 0 & 0 & e_5 & e_6 & e_7 & 0
\end{array}
\]
%%%%%%%%%%%%%%%%%%%%%%%%%%%%%%%%%%%%%%%%%%%%%%%%%%%%%%%%%%%%%%%%%%%%%%
%%%%%%%%%%%%%%%%%%%%%%%%%%%%%%%%%%%%%%%%%%%%%%%%%%%%%%%%%%%%%%%%%%%%%%
%%%%%%%%%%%%%%%%%%%%%%%%%%%%%%%%%%%%%%%%%%%%%%%%%%%%%%%%%%%%%%%%%%%%%%

\[
\text{\bf Model 2b3b4a5b}
\ \ \ \ \
\left\{
\aligned
u
&
\,=\,
x^2+x^2y+x^2y^2+x^2y^3+x^5+x^2y^4+4\,x^5y
\\
&
\ \ \ \ \
+x^2y^5+10\,x^5y^2+F_{7,0}\,x^7+\cdots,
\\
v
&
\,=\,
0,
\endaligned\right.
\]
for any value of $F_{7,0}$,
%%%%%%%%%%%%%%%%%%%%%%%%%%%%%%%%%%%%%%%%%%%%%%%%%%%%%%%%%%%%%%%%%%%%%%
\[
\def\arraystretch{1.25}
\begin{array}{llll}
e_1
\,:=\,
\big(1-y+7\,F_{7,0}\,u\big)\,\partial_x
-\big(14\,F_{7,0}\,x+5\,u\big)\,\partial_y+2\,x\,\partial_u, &
e_2
\,:=\,
-\,x\,\partial_x+(1-y)\,\partial_y-u\,\partial_u, &
& 
\\
e_3
\,:=\,
v\,\partial_x, &
e_4
\,:=\,
v\,\partial_y, &
\\
e_5
\,:=\,
v\,\partial_u, &
e_6
\,:=\,
v\,\partial_v,
\end{array}
\]
%%%%%%%%%%%%%%%%%%%%%%%%%%%%%%%%%%%%%%%%%%%%%%%%%%%%%%%%%%%%%%%%%%%%%%
\[
\footnotesize
\def\arraystretch{1.25}
\begin{array}{c|cccccc}
{} & e_1 & e_2 & e_3 & e_4 & e_5 & e_6 
\\
\hline
e_1 & 
0 & 0 & 14\,F_{7,0}\,e_4-2\,e_5 & e_3 & -\,7\,F_{7,0}\,e_3+5\,e_4 &
0
\\
e_2 &
0 & 0 & e_3 & e_4 & e_5 & 0
\\
e_3 &
-\,14\,F_{7,0}\,e_4+2\,e_5 & -\,e_3 & 0 & 0 & 0 & -\,e_3
\\
e_4 &
-\,e_3 & -\,e_4 & 0 & 0 & 0 & -\,e_4
\\
e_5 &
7\,F_{7,0}\,e_3-5\,e_4 & -\,e_5 & 0 & 0 & 0 & -\,e_5
\\
e_6 &
0 & 0 & e_3 & e_4 & e_5 & 0
\end{array}
\]

%%%%%%%%%%%%%%%%%%%%%%%%%%%%%%%%%%%%%%%%%%%%%%%%%%%%%%%%%%%%%%%%%%%%%%
%%%%%%%%%%%%%%%%%%%%%%%%%%%%%%%%%%%%%%%%%%%%%%%%%%%%%%%%%%%%%%%%%%%%%%
%%%%%%%%%%%%%%%%%%%%%%%%%%%%%%%%%%%%%%%%%%%%%%%%%%%%%%%%%%%%%%%%%%%%%%

\[
\text{\bf Model 2b3b4b}
\ \ \ \ \
\left\{
\aligned
u
&
\,=\,
x^2+x^2y+x^2y^2+x^3y+x^2y^3+3\,x^3y^2+x^5+x^2y^4+6\,x^3y^3
\\
&
\ \ \ \ \
+\tfrac{9}{4}\,x^4y^2+3\,x^5y+x^6
+\cdots,
\\
v
&
\,=\,
0,
\endaligned\right.
\]
%%%%%%%%%%%%%%%%%%%%%%%%%%%%%%%%%%%%%%%%%%%%%%%%%%%%%%%%%%%%%%%%%%%%%%
\[
\def\arraystretch{1.25}
\begin{array}{llll}
e_1
\,:=\,
\big(1+3\,x-y-5\,u\big)\,\partial_x
+\big(10\,x-3\,y-15\,u\big)\,\partial_y
+\big(2\,x+6\,u\big)\,\partial_u, &
\\
e_2
\,:=\,
\big(-\,2\,x+\tfrac{3}{2}\,u\big)\,\partial_x
+
\big(1-4\,x-y+4\,u\big)\,\partial_y
-\,3\,u\,\partial_u, &
& 
\\
e_3
\,:=\,
v\,\partial_x, &
e_4
\,:=\,
v\,\partial_y, &
\\
e_5
\,:=\,
v\,\partial_u, &
e_6
\,:=\,
v\,\partial_v,
\end{array}
\]
%%%%%%%%%%%%%%%%%%%%%%%%%%%%%%%%%%%%%%%%%%%%%%%%%%%%%%%%%%%%%%%%%%%%%%
\[
\footnotesize
\def\arraystretch{1.25}
\begin{array}{c|cccccc}
{} & e_1 & e_2 & e_3 & e_4 & e_5 & e_6 
\\
\hline
e_1 & 
0 & -\,e_1-e_2 & -\,3\,e_3-10\,e_4-2\,e_5 & e_3+3\,e_4 & 
5\,e_3+15\,e_4-6\,e_5 & 0
\\
e_2 &
e_1+e_2 & 0 & 2\,e_3+4\,e_4 & e_4 &
-\frac{3}{2}\,e_3-4\,e_4+3\,e_5 & 0
\\
e_3 &
3\,e_3+10\,e_4+2\,e_5 & -2\,e_3-4\,e_4 & 0 & 0 & 0 & -\,e_3
\\
e_4 &
-\,e_3-3\,e_4 & -\,e_4 & 0 & 0 & 0 & -\,e_4
\\
e_5 &
-\,5\,e_3-15\,e_4+6\,e_5 & \frac{3}{2}\,e_3+4\,e_4-3\,e_5 & 0
& 0 & 0 & -\,e_5
\\
e_6 &
0 & 0 & e_3 & e_4 & e_5 & 0
\end{array}
\]

%%%%%%%%%%%%%%%%%%%%%%%%%%%%%%%%%%%%%%%%%%%%%%%%%%%%%%%%%%%%%%%%%%%%%%
%%%%%%%%%%%%%%%%%%%%%%%%%%%%%%%%%%%%%%%%%%%%%%%%%%%%%%%%%%%%%%%%%%%%%%
%%%%%%%%%%%%%%%%%%%%%%%%%%%%%%%%%%%%%%%%%%%%%%%%%%%%%%%%%%%%%%%%%%%%%%

\[
\text{\bf Model 2b3c5a}
\ \ \ \ \
\left\{
\aligned
u
&
\,=\,
x^2,
\\
v
&
\,=\,
x^3,
\endaligned\right.
\]
%%%%%%%%%%%%%%%%%%%%%%%%%%%%%%%%%%%%%%%%%%%%%%%%%%%%%%%%%%%%%%%%%%%%%%
\[
\def\arraystretch{1.25}
\begin{array}{llll}
e_1
\,:=\,
\partial_x+2\,x\,\partial_u+3\,u\,\partial_v, &
e_2
\,:=\,
\partial_y, &
& 
\\
e_3
\,:=\,
x\,\partial_x+2\,u\,\partial_u+3\,v\,\partial_v, &
e_4
\,:=\,
x\,\partial_y, &
e_5
\,:=\,
y\,\partial_y, &
e_6
\,:=\,
u\,\partial_y,
\\
e_7
\,:=\,
v\,\partial_y,
\end{array}
\]
%%%%%%%%%%%%%%%%%%%%%%%%%%%%%%%%%%%%%%%%%%%%%%%%%%%%%%%%%%%%%%%%%%%%%%
\[
\footnotesize
\def\arraystretch{1.25}
\begin{array}{c|ccccccc}
{} & e_1 & e_2 & e_3 & e_4 & e_5 & e_6 & e_7
\\
\hline
e_1 & 
0 & 0 & e_1 & e_2 & 0 & 2\,e_4 & 3\,e_6
\\
e_2 &
0 & 0 & 0 & 0 & e_2 & 0 & 0
\\
e_3 &
-\,e_1 & 0 & 0 & e_4 & 0 & 2\,e_6 & 3\,e_7
\\
e_4 &
-\,e_2 & 0 & -\,e_4 & 0 & e_4 & -\,e_6 & -\,e_7
\\
e_5 &
0 & -\,e_2 & 0 & -\,e_4 & 0 & -\,e_6 & -\,e_7
\\
e_6 &
-\,2\,e_4 & 0 & -\,2\,e_6 & 0 & e_6 & 0 & 0
\\
e_7 &
-\,3\,e_6 & 0 & -\,3\,e_7 & 0 & e_7 & 0 & 0
\end{array}
\]

%%%%%%%%%%%%%%%%%%%%%%%%%%%%%%%%%%%%%%%%%%%%%%%%%%%%%%%%%%%%%%%%%%%%%%
%%%%%%%%%%%%%%%%%%%%%%%%%%%%%%%%%%%%%%%%%%%%%%%%%%%%%%%%%%%%%%%%%%%%%%
%%%%%%%%%%%%%%%%%%%%%%%%%%%%%%%%%%%%%%%%%%%%%%%%%%%%%%%%%%%%%%%%%%%%%%

\[
\text{\bf Model 2b3c5b}
\ \ \ \ \
\left\{
\aligned
u
&
\,=\,
x^2-\tfrac{5}{9}\,x^6-\tfrac{10}{7}\,G_{6,0}\,x^7
+\big(-\,\tfrac{250}{189}-\tfrac{45}{14}\,G_{6,0}^2\big)\,x^8
+\cdots,
\\
v
&
\,=\,
x^3+x^5+G_{6,0}\,x^6
+\big(\tfrac{20}{21}+\tfrac{9}{7}\,G_{6,0}^2\big)\,x^7
+\big(\tfrac{15}{7}\,G_{6,0}+\tfrac{27}{14}\,G_{6,0}^3\big)\,x^8
+\cdots,
\endaligned\right.
\]
for any value of $G_{6,0}$,
%%%%%%%%%%%%%%%%%%%%%%%%%%%%%%%%%%%%%%%%%%%%%%%%%%%%%%%%%%%%%%%%%%%%%%
\[
\def\arraystretch{1.25}
\begin{array}{llll}
e_1
\,:=\,
\big(1-3\,G_{6,0}\,x-\tfrac{5}{3}\,u\big)\,\partial_x
+\big(2\,x-6\,G_{6,0}\,u-\tfrac{10}{3}\,v\big)\,\partial_u
+
\big(3\,u-9\,G_{6,0}\,v\big)\,\partial_v, &
e_2
\,:=\,
\partial_y, &
& 
\\
e_3
\,:=\,
x\,\partial_y, &
e_4
\,:=\,
y\,\partial_y, &
\\
e_5
\,:=\,
u\,\partial_y, &
e_6
\,:=\,
v\,\partial_y,
\end{array}
\]
%%%%%%%%%%%%%%%%%%%%%%%%%%%%%%%%%%%%%%%%%%%%%%%%%%%%%%%%%%%%%%%%%%%%%%
\[
\footnotesize
\def\arraystretch{1.25}
\begin{array}{c|cccccc}
{} & e_1 & e_2 & e_3 & e_4 & e_5 & e_6 
\\
\hline
e_1 & 
0 & 0 & e_2-3\,G_{6,0}\,e_3-\frac{5}{3}\,e_5 & 0 &
2\,e_3-6\,G_{6,0}\,e_5-\frac{10}{3}\,e_6 &
3\,e_5-9\,G_{6,0}\,e_6
\\
e_2 &
0 & 0 & 0 & e_2 & 0 & 0
\\
e_3 &
-\,e_2+3\,G_{6,0}\,e_3+\frac{5}{3}\,e_5 & 0 & 0 &
e_3 & 0 & 0
\\
e_4 &
0 & -\,e_2 & -\,e_3 & 0 & -\,e_5 & -\,e_6
\\
e_5 &
-\,2\,e_3+6\,G_{6,0}\,e_5+\frac{10}{3}\,e_6 & 0 & 0 &
e_5 & 0 & 0
\\
e_6 &
-\,3\,e_5+9\,G_{6,0}\,e_6 & 0 & 0 & e_6 & 0 & 0 
\end{array}
\]

%%%%%%%%%%%%%%%%%%%%%%%%%%%%%%%%%%%%%%%%%%%%%%%%%%%%%%%%%%%%%%%%%%%%%%
%%%%%%%%%%%%%%%%%%%%%%%%%%%%%%%%%%%%%%%%%%%%%%%%%%%%%%%%%%%%%%%%%%%%%%
%%%%%%%%%%%%%%%%%%%%%%%%%%%%%%%%%%%%%%%%%%%%%%%%%%%%%%%%%%%%%%%%%%%%%%

\[
\text{\bf Model 2b3c5c}
\ \ \ \ \
\left\{
\aligned
u
&
\,=\,
x^2-\tfrac{5}{9}\,x^6+\tfrac{10}{7}\,G_{6,0}\,x^7
+\big(\tfrac{250}{189}-\tfrac{45}{14}\,G_{6,0}^2\big)\,x^8
+\cdots,
\\
v
&
\,=\,
x^3-x^5+G_{6,0}\,x^6
+\big(\tfrac{20}{21}-\tfrac{9}{7}\,G_{6,0}^2\big)\,x^7
+\big(-\,\tfrac{15}{7}\,G_{6,0}+\tfrac{27}{14}\,G_{6,0}^3\big)\,x^8
+\cdots,
\endaligned\right.
\]
for any value of $G_{6,0}$,
%%%%%%%%%%%%%%%%%%%%%%%%%%%%%%%%%%%%%%%%%%%%%%%%%%%%%%%%%%%%%%%%%%%%%%
\[
\def\arraystretch{1.25}
\begin{array}{llll}
e_1
\,:=\,
\big(1+3\,G_{6,0}\,x+\tfrac{5}{3}\,u\big)\,\partial_x
+\big(2\,x+6\,G_{6,0}\,u+\tfrac{10}{3}\,v\big)\,\partial_u
+
\big(3\,u+9\,G_{6,0}\,v\big)\,\partial_v, &
e_2
\,:=\,
\partial_y, &
& 
\\
e_3
\,:=\,
x\,\partial_y, &
e_4
\,:=\,
y\,\partial_y, &
\\
e_5
\,:=\,
u\,\partial_y, &
e_6
\,:=\,
v\,\partial_y,
\end{array}
\]
%%%%%%%%%%%%%%%%%%%%%%%%%%%%%%%%%%%%%%%%%%%%%%%%%%%%%%%%%%%%%%%%%%%%%%
\[
\footnotesize
\def\arraystretch{1.25}
\begin{array}{c|cccccc}
{} & e_1 & e_2 & e_3 & e_4 & e_5 & e_6 
\\
\hline
e_1 & 
0 & 0 & e_2+3\,G_{6,0}\,e_3+\frac{5}{3}\,e_5 & 0 &
2\,e_3+6\,G_{6,0}\,e_5+\frac{10}{3}\,e_6 &
3\,e_5+9\,G_{6,0}\,e_6
\\
e_2 &
0 & 0 & 0 & e_2 & 0 & 0
\\
e_3 &
-\,e_2-3\,G_{6,0}\,e_3-\frac{5}{3}\,e_5 & 0 & 0 &
e_3 & 0 & 0
\\
e_4 &
0 & -\,e_2 & -\,e_3 & 0 & -\,e_5 & -\,e_6
\\
e_5 &
-\,2\,e_3-6\,G_{6,0}\,e_5-\frac{10}{3}\,e_6 & 0 & 0 &
e_5 & 0 & 0
\\
e_6 &
-\,3\,e_5-9\,G_{6,0}\,e_6 & 0 & 0 & e_6 & 0 & 0 
\end{array}
\]

%%%%%%%%%%%%%%%%%%%%%%%%%%%%%%%%%%%%%%%%%%%%%%%%%%%%%%%%%%%%%%%%%%%%%%
%%%%%%%%%%%%%%%%%%%%%%%%%%%%%%%%%%%%%%%%%%%%%%%%%%%%%%%%%%%%%%%%%%%%%%
%%%%%%%%%%%%%%%%%%%%%%%%%%%%%%%%%%%%%%%%%%%%%%%%%%%%%%%%%%%%%%%%%%%%%%

\[
\text{\bf Model 2b3c5d}
\ \ \ \ \
\left\{
\aligned
u
&
\,=\,
x^2+x^5+F_{6,0}\,x^6
+\big(\tfrac{65}{42}\,G_{5,0}+\tfrac{8}{7}\,F_{6,0}^2
+\tfrac{20}{63}\,F_{6,0}\,G_{5,0}^2
-\tfrac{100}{567}\,G_{5,0}^4\big)\,x^7
+\cdots,
\\
v
&
\,=\,
x^3+G_{5,0}\,x^5
+\big(\tfrac{7}{4}+\tfrac{2}{3}\,G_{5,0}\,F_{6,0}
+\tfrac{10}{27}\,G_{5,0}^3\big)\,x^6
\\
&
\ \ \ \ \
+\big(\tfrac{85}{42}\,G_{5,0}^2
+\tfrac{27}{14}\,F_{6,0}+\tfrac{4}{7}\,G_{5,0}\,F_{6,0}^2
+\tfrac{40}{63}\,G_{5,0}^3F_{6,0}+\tfrac{100}{567}\,G_{5,0}^5\big)\,x^7
+\cdots,
\endaligned\right.
\]
for any value of $G_{5,0},F_{6,0}$,
%%%%%%%%%%%%%%%%%%%%%%%%%%%%%%%%%%%%%%%%%%%%%%%%%%%%%%%%%%%%%%%%%%%%%%
\[
\aligned
e_1
&
\,:=\,
\big(
1-2\,F_{6,0}\,x-\tfrac{10}{9}\,G_{5,0}^2\,x-\tfrac{5}{3}\,G_{5,0}\,u
-\tfrac{5}{2}\,v
\big)\,\partial_x
\\
&
\ \ \ \ \ 
+
\big(
2\,x-4\,F_{6,0}\,u-\tfrac{20}{9}\,G_{5,0}^2\,u
-\tfrac{10}{3}\,G_{5,0}\,v
\big)\,\partial_u
\\
&
\ \ \ \ \ 
+
\big(
3\,u-6\,F_{6,0}\,v-\tfrac{10}{3}\,G_{5,0}^2\,v
\big)\,\partial_v,
\\
\ \ \ \ \ 
e_2
&
\,:=\,
\partial_y,
\\
e_3
&
\,:=\,
x\,\partial_y,
\ \ \ \ \ 
e_4
\,:=\,
y\,\partial_y,
\ \ \ \ \ 
e_5
\,:=\,
u\,\partial_y,
\ \ \ \ \ 
e_6
\,:=\,
v\,\partial_y,
\endaligned
\]
%%%%%%%%%%%%%%%%%%%%%%%%%%%%%%%%%%%%%%%%%%%%%%%%%%%%%%%%%%%%%%%%%%%%%%
\[
\aligned
{}
[e_1,e_3]
&
\,=\,
e_2
-
\big(\tfrac{10}{9}\,G_{5,0}^2+2\,F_{6,0}\big)\,e_3
-
\tfrac{5}{3}\,G_{5,0}\,e_5
-
\tfrac{5}{2}\,e_6,
\\
[e_1,e_5]
&
\,=\,
2\,e_3
-
\big(\tfrac{20}{9}\,G_{5,0}^2+4\,F_{6,0}\big)\,e_5
-
\tfrac{10}{3}\,G_{5,0}\,e_6,
\\
[e_1,e_6]
&
\,=\,
3\,e_5
-
\big(\tfrac{10}{3}\,G_{5,0}^2+6\,F_{6,0}\big)\,e_6,
\\
[e_4,e_5]
&
\,=\,
-\,e_5,
\\
[e_4,e_6]
&
\,=\,
-\,e_6.
\endaligned
\]
%%%%%%%%%%%%%%%%%%%%%%%%%%%%%%%%%%%%%%%%%%%%%%%%%%%%%%%%%%%%%%%%%%%%%%
%%%%%%%%%%%%%%%%%%%%%%%%%%%%%%%%%%%%%%%%%%%%%%%%%%%%%%%%%%%%%%%%%%%%%%
%%%%%%%%%%%%%%%%%%%%%%%%%%%%%%%%%%%%%%%%%%%%%%%%%%%%%%%%%%%%%%%%%%%%%%

\[
\text{\bf Model 2b3d4a}
\ \ \ \ \
\left\{
\aligned
u
&
\,=\,
x^2+x^2y+x^2y^2+x^3y+x^2y^3+3x^3y^2+x^2y^4+
\\
&
\ \ \ \ \ 
+
6x^3y^3+\tfrac{9}{4}x^4y^2-\tfrac{27}{32}x^5y-\tfrac{135}{256}x^6
+
\cdots
,
\\
v
&
\,=\,
x^3+2x^3y+3x^3y^2+\tfrac{9}{4}x^4y+4x^3y^3
+9x^4y^2+\tfrac{27}{64}x^6
+
\cdots,
\endaligned\right.
\]
%%%%%%%%%%%%%%%%%%%%%%%%%%%%%%%%%%%%%%%%%%%%%%%%%%%%%%%%%%%%%%%%%%%%%%
\[
\def\arraystretch{1.25}
\begin{array}{ll}
e_1
&
\,:=\,
(
3x+1-y-\tfrac{27}{4}u+\tfrac{81}{16}v
)\partial_x
+
(
\tfrac{81}{8}x-3y-\tfrac{81}{4}u+\tfrac{2997}{128}v
)\partial_y
+
\\
&
\ \ \ \ \ 
+
(
2x+6u-\tfrac{27}{8}v
)\partial_u
+
(
3u+9v
)\partial_v
,
\\
e_2
&
\,:=\,
(
-4x+\tfrac{9}{2}u-\tfrac{9}{4}v
)\partial_x
+
(
-9x+2-2y+\tfrac{27}{2}u-\tfrac{189}{16}v
)\partial_y
+
(
-6u+2v
)\partial_u
-8v\partial_v,
\end{array}
\]
%%%%%%%%%%%%%%%%%%%%%%%%%%%%%%%%%%%%%%%%%%%%%%%%%%%%%%%%%%%%%%%%%%%%%%
\[
\footnotesize
\def\arraystretch{1.25}
\begin{array}{c|cc}
{} & e_1 & e_2 
\\
\hline
e_1 & 
0 & -2e_1-\tfrac{3}{2}e_2
\\
e_2 &
2e_1+\tfrac{3}{2}e_2 & 0
\end{array}
\]

%%%%%%%%%%%%%%%%%%%%%%%%%%%%%%%%%%%%%%%%%%%%%%%%%%%%%%%%%%%%%%%%%%%%%%
%%%%%%%%%%%%%%%%%%%%%%%%%%%%%%%%%%%%%%%%%%%%%%%%%%%%%%%%%%%%%%%%%%%%%%
%%%%%%%%%%%%%%%%%%%%%%%%%%%%%%%%%%%%%%%%%%%%%%%%%%%%%%%%%%%%%%%%%%%%%%

\[
\text{\bf Model 2b3d4b6a}
\ \ \ \ \
\left\{
\aligned
u
&
\,=\,
x^2+x^2y+x^2y^2+x^2y^3+x^2y^4+x^2y^5+x^2y^6
+\cdots,
\\
v
&
\,=\,
x^3+2\,x^3y+3\,x^3y^2+4\,x^3y^3+5\,x^3y^4+6\,x^3y^5
+\cdots,
\endaligned\right.
\]
%%%%%%%%%%%%%%%%%%%%%%%%%%%%%%%%%%%%%%%%%%%%%%%%%%%%%%%%%%%%%%%%%%%%%%
\[
\def\arraystretch{1.25}
\begin{array}{llll}
e_1
\,:=\,
\big(1-y\big)\,\partial_x
+2\,x\,\partial_u+3\,u\,\partial_v, &
e_2
\,:=\,
\big(1-y\big)\,\partial_y+u\,\partial_u+2\,v\,\partial_v, &
& 
\\
e_3
\,:=\,
x\,\partial_x+2\,u\,\partial_u+3\,v\,\partial_v, &
e_4
\,:=\,
-\,\tfrac{2}{3}\,u\,\partial_x
+x\,\partial_y-\tfrac{1}{3}\,v\,\partial_u,
\end{array}
\]
%%%%%%%%%%%%%%%%%%%%%%%%%%%%%%%%%%%%%%%%%%%%%%%%%%%%%%%%%%%%%%%%%%%%%%
\[
\footnotesize
\def\arraystretch{1.25}
\begin{array}{c|cccc}
{} & e_1 & e_2 & e_3 & e_4
\\
\hline
e_1 & 
0 & e_1 & e_1 & e_2-\frac{1}{3}\,e_3
\\
e_2 &
-\,e_1 & 0 & 0 & e_4
\\
e_3 &
-\,e_1 & 0 & 0 & e_4
\\
e_4 &
-\,e_2+\frac{1}{3}\,e_3 & -\,e_4 & -\,e_4 & 0
\end{array}
\]

%%%%%%%%%%%%%%%%%%%%%%%%%%%%%%%%%%%%%%%%%%%%%%%%%%%%%%%%%%%%%%%%%%%%%%
%%%%%%%%%%%%%%%%%%%%%%%%%%%%%%%%%%%%%%%%%%%%%%%%%%%%%%%%%%%%%%%%%%%%%%
%%%%%%%%%%%%%%%%%%%%%%%%%%%%%%%%%%%%%%%%%%%%%%%%%%%%%%%%%%%%%%%%%%%%%%

\[
\text{\bf Model 2b3d4b6b}
\ \ \ \ \
\left\{
\aligned
u
&
\,=\,
x^2+x^2y+x^2y^2+x^2y^3+x^2y^4+x^6+x^2y^5+5\,x^6y
\\
&
\ \ \ \ \ 
+x^2y^6+15\,x^6y^2+F_{8,0}\,x^8+x^2y^7+35\,x^6y^3+7\,F_{8,0}
+\cdots,
\\
v
&
\,=\,
x^3+2\,x^3y+3\,x^3y^2+4\,x^3y^3+5\,x^3y^4+\tfrac{15}{7}\,x^7
\\
&
\ \ \ \ \ 
+6\,x^3y^5+\tfrac{90}{7}\,x^7y+7\,x^3y^6+45\,x^7y^2
+\tfrac{7}{3}\,F_{8,0}\,x^9+\cdots,
\endaligned\right.
\]
for any value of $F_{8,0}$,
%%%%%%%%%%%%%%%%%%%%%%%%%%%%%%%%%%%%%%%%%%%%%%%%%%%%%%%%%%%%%%%%%%%%%%
\[
\aligned
e_1
&
\,:=\,
\big(1-y+14\,F_{8,0}\,u\big)\,\partial_x
-\big(21\,F_{8,0}\,x+6\,v\big)\,\partial_y
+\big(2\,x+7\,F_{8,0}\,v\big)\,\partial_u
+3\,u\,\partial_v,
\\
e_2
&
\,:=\,
-\,x\,\partial_x
+\big(1-y\big)\,\partial_y
-u\,\partial_u-v\,\partial_v,
\endaligned
\]
%%%%%%%%%%%%%%%%%%%%%%%%%%%%%%%%%%%%%%%%%%%%%%%%%%%%%%%%%%%%%%%%%%%%%%
\[
\footnotesize
\def\arraystretch{1.25}
\begin{array}{c|cc}
{} & e_1 & e_2 
\\
\hline
e_1 & 
0 & 0
\\
e_2 &
0 & 0
\end{array}
\]

%%%%%%%%%%%%%%%%%%%%%%%%%%%%%%%%%%%%%%%%%%%%%%%%%%%%%%%%%%%%%%%%%%%%%%
%%%%%%%%%%%%%%%%%%%%%%%%%%%%%%%%%%%%%%%%%%%%%%%%%%%%%%%%%%%%%%%%%%%%%%
%%%%%%%%%%%%%%%%%%%%%%%%%%%%%%%%%%%%%%%%%%%%%%%%%%%%%%%%%%%%%%%%%%%%%%

\[
\text{\bf Model 2b3d4b6c}
\ \ \ \ \
\left\{
\aligned
u
&
\,=\,
x^2+x^2y+x^2y^2+x^2y^3+x^2y^4-x^6+x^2y^5-5\,x^6y
\\
&
\ \ \ \ \ 
+x^2y^6-15\,x^6y^2+F_{8,0}\,x^8+x^2y^7-35\,x^6y^3
+7\,F_{8,0}\,x^8y
+\cdots,
\\
v
&
\,=\,
x^3+2\,x^3y+3\,x^3y^2+4\,x^3y^3+5\,x^3y^4-\tfrac{15}{7}\,x^7
\\
&
\ \ \ \ \ 
+6\,x^3y^5-\tfrac{90}{7}\,x^7y+7\,x^3y^6-45\,x^7y^2
+\tfrac{7}{3}\,F_{8,0}\,x^9+\cdots,
\endaligned\right.
\]
for any value of $F_{8,0}$,
%%%%%%%%%%%%%%%%%%%%%%%%%%%%%%%%%%%%%%%%%%%%%%%%%%%%%%%%%%%%%%%%%%%%%%
\[
\aligned
e_1
&
\,:=\,
\big(1-y-14\,F_{8,0}\,u\big)\,\partial_x
+\big(21\,F_{8,0}\,x+6\,v\big)\,\partial_y
+\big(2\,x-7\,F_{8,0}\,v\big)\,\partial_u
+3\,u\,\partial_v,
\\
e_2
&
\,:=\,
-\,x\,\partial_x
+\big(1-y\big)\,\partial_y
-u\,\partial_u-v\,\partial_v,
\endaligned
\]
%%%%%%%%%%%%%%%%%%%%%%%%%%%%%%%%%%%%%%%%%%%%%%%%%%%%%%%%%%%%%%%%%%%%%%
\[
\footnotesize
\def\arraystretch{1.25}
\begin{array}{c|cc}
{} & e_1 & e_2 
\\
\hline
e_1 & 
0 & 0
\\
e_2 &
0 & 0
\end{array}
\]

%%%%%%%%%%%%%%%%%%%%%%%%%%%%%%%%%%%%%%%%%%%%%%%%%%%%%%%%%%%%%%%%%%%%%%
%%%%%%%%%%%%%%%%%%%%%%%%%%%%%%%%%%%%%%%%%%%%%%%%%%%%%%%%%%%%%%%%%%%%%%
%%%%%%%%%%%%%%%%%%%%%%%%%%%%%%%%%%%%%%%%%%%%%%%%%%%%%%%%%%%%%%%%%%%%%%

\[
\text{\bf Model 2b3d4b6d}
\ \ \ \ \
\left\{
\aligned
u
&
\,=\,
x^2+x^2y+x^2y^2+x^2y^3+x^2y^4+x^2y^5+F_{7,0}x^7+x^2y^6
+
\\
&
\ \ \ \ \ 
+
6F_{7,0}x^7y+(-\tfrac{3}{2}+\tfrac{7}{3}F_{7,0}G_{7,0})x^8
+
\cdots
,
\\
v
&
\,=\,
x^3+2x^3y+3x^3y^2+4x^3y^3+x^6+5x^3y^4+5x^6y+G_{7,0}x^7
+
\\
&
\ \ \ \ \ 
+
6x^3y^5+15x^6y^2+6G_{7,0}x^7y+(3F_{7,0}+\tfrac{7}{6}G_{7,0}^2)x^8
+
\cdots,
\endaligned\right.
\]
for any values of $F_{7,0}$, $G_{7,0}$,
%%%%%%%%%%%%%%%%%%%%%%%%%%%%%%%%%%%%%%%%%%%%%%%%%%%%%%%%%%%%%%%%%%%%%%
\[
\def\arraystretch{1.25}
\begin{array}{ll}
e_1
&
\,:=\,
\big(
1-\tfrac{7}{3}G_{7,0}x-y+14F_{7,0}u+6v
\big)\partial_x
-
\big(
21F_{7,0}x+12u
\big)\partial_y
+
\big(
2x-\tfrac{14}{3}G_{7,0}u+7F_{7,0}v
\big)\partial_u
+
\\
&
\ \ \ \ \ 
-
\big(
7G_{7,0}v-3u
\big)\partial_v
,
\\
e_2
&
\,:=\,
-x\partial_x-(y-1)\partial_y-u\partial_u-v\partial_v,
\end{array}
\]
%%%%%%%%%%%%%%%%%%%%%%%%%%%%%%%%%%%%%%%%%%%%%%%%%%%%%%%%%%%%%%%%%%%%%%
\[
\footnotesize
\def\arraystretch{1.25}
\begin{array}{c|cc}
{} & e_1 & e_2 
\\
\hline
e_1 & 
0 & 0
\\
e_2 &
0 & 0
\end{array}
\]

%%%%%%%%%%%%%%%%%%%%%%%%%%%%%%%%%%%%%%%%%%%%%%%%%%%%%%%%%%%%%%%%%%%%%%
%%%%%%%%%%%%%%%%%%%%%%%%%%%%%%%%%%%%%%%%%%%%%%%%%%%%%%%%%%%%%%%%%%%%%%
%%%%%%%%%%%%%%%%%%%%%%%%%%%%%%%%%%%%%%%%%%%%%%%%%%%%%%%%%%%%%%%%%%%%%%

%%%%%%%%%%%%%%%%%%%%%%%%%%%%%%%%%%%%%%%%%%%%%%%%%%%%%%%%%%%%%%%%%%%%%%
\SectionHead{2c Models}
{2c-models}
%%%%%%%%%%%%%%%%%%%%%%%%%%%%%%%%%%%%%%%%%%%%%%%%%%%%%%%%%%%%%%%%%%%%%%

%%%%%%%%%%%%%%%%%%%%%%%%%%%%%%%%%%%%%%%%%%%%%%%%%%%%%%%%%%%%%%%%%%%%%%
%%%%%%%%%%%%%%%%%%%%%%%%%%%%%%%%%%%%%%%%%%%%%%%%%%%%%%%%%%%%%%%%%%%%%%
%%%%%%%%%%%%%%%%%%%%%%%%%%%%%%%%%%%%%%%%%%%%%%%%%%%%%%%%%%%%%%%%%%%%%%

\[
\text{\bf Model 2c3a4a}
\ \ \ \ \
\left\{
\aligned
u
&
\,=\,
xy
,
\\
v
&
\,=\,
0,
\endaligned\right.
\]
%%%%%%%%%%%%%%%%%%%%%%%%%%%%%%%%%%%%%%%%%%%%%%%%%%%%%%%%%%%%%%%%%%%%%%
\[
\def\arraystretch{1.25}
\begin{array}{llll}
e_1
\,:=\,
\partial_x+y\partial_u
, &
e_2
\,:=\,
\partial_y+x\partial_u
, &
e_3
\,:=\,
x\partial_x+u\partial_u
, &
e_4
\,:=\,
y\partial_y+u\partial_u
, 
\\
e_5
\,:=\,
v\partial_x
, &
e_6
\,:=\,
v\partial_y
, &
e_7
\,:=\,
v\partial_u
, &
e_8
\,:=\,
v\partial_v,
\end{array}
\]
%%%%%%%%%%%%%%%%%%%%%%%%%%%%%%%%%%%%%%%%%%%%%%%%%%%%%%%%%%%%%%%%%%%%%%
\[
\footnotesize
\def\arraystretch{1.25}
\begin{array}{c|cccccccc}
{} & e_1 & e_2 & e_3 & e_4 & e_5 & e_6 & e_7 & e_8 
\\
\hline
e_1 & 0 & 0 & e_1 & 0 & 0 & -e_7 & 0 & 0
\\
e_2 & 0 & 0 & 0 & e_2 & -e_7 & 0 & 0 & 0
\\
e_3 & -e_1 & 0 & 0 & 0 & -e_5 & 0 & -e_7 & 0
\\
e_4 & 0 & -e_2 & 0 & 0 & 0 & -e_6 & -e_7 & 0
\\
e_5 & 0 & e_7 & e_5 & 0 & 0 & 0 & 0 & -e_5
\\
e_6 & e_7 & 0 & 0 & e_6 & 0 & 0 & 0 & -e_6
\\
e_7 & 0 & 0 & e_7 & e_7 & 0 & 0 & 0 & -e_7
\\
e_8 & 0 & 0 & 0 & 0 & e_5 & e_ 6 & e_7 & 0
\end{array}
\]

%%%%%%%%%%%%%%%%%%%%%%%%%%%%%%%%%%%%%%%%%%%%%%%%%%%%%%%%%%%%%%%%%%%%%%
%%%%%%%%%%%%%%%%%%%%%%%%%%%%%%%%%%%%%%%%%%%%%%%%%%%%%%%%%%%%%%%%%%%%%%
%%%%%%%%%%%%%%%%%%%%%%%%%%%%%%%%%%%%%%%%%%%%%%%%%%%%%%%%%%%%%%%%%%%%%%

\[
\text{\bf Model 2c3a4b}
\ \ \ \ \
\left\{
\aligned
u
&
\,=\,
xy+x^2y^2+2x^3y^3+5x^4y^4+14x^5y^5
+
\cdots
,
\\
v
&
\,=\,
0,
\endaligned\right.
\]
%%%%%%%%%%%%%%%%%%%%%%%%%%%%%%%%%%%%%%%%%%%%%%%%%%%%%%%%%%%%%%%%%%%%%%
\[
\def\arraystretch{1.25}
\begin{array}{llll}
e_1
\,:=\,
-
(
2u-1
)\partial_x
+
y\partial_u
, &
e_2
\,:=\,
-
(
2u-1
)\partial_y
+
x\partial_u
, &
e_3
\,:=\,
x\partial_x-y\partial_y
,
\\
e_4
\,:=\,
v\partial_x
, &
e_5
\,:=\,
v\partial_y
, &
e_6
\,:=\,
v\partial_u
,
\\
e_7
\,:=\,
v\partial_v,
\end{array}
\]
%%%%%%%%%%%%%%%%%%%%%%%%%%%%%%%%%%%%%%%%%%%%%%%%%%%%%%%%%%%%%%%%%%%%%%
\[
\footnotesize
\def\arraystretch{1.25}
\begin{array}{c|ccccccc}
{} & e_1 & e_2 & e_3 & e_4 & e_5 & e_6 & e_7
\\
\hline
e_1 & 0 & 2e_3 & e_1 & 0 & -e_6 & 2e_4 & 0
\\
e_2 & -2e_3 & 0 & -e_2 & -e_6 & 0 & 2e_5 & 0
\\
e_3 & -e_1 & e_2 & 0 & -e_4 & e_5 & 0 & 0
\\
e_4 & 0 & e_6 & e_4 & 0 & 0 & 0 & -e_4
\\
e_5 & e_6 & 0 & -e_5 & 0 & 0 & 0 & -e_5
\\
e_6 & -2e_4 & -2e_5 & 0 & 0 & 0 & 0 & -e_6
\\
e_7 & 0 & 0 & 0 & e_4 & e_5 & e_6 & 0
\end{array}
\]

%%%%%%%%%%%%%%%%%%%%%%%%%%%%%%%%%%%%%%%%%%%%%%%%%%%%%%%%%%%%%%%%%%%%%%
%%%%%%%%%%%%%%%%%%%%%%%%%%%%%%%%%%%%%%%%%%%%%%%%%%%%%%%%%%%%%%%%%%%%%%
%%%%%%%%%%%%%%%%%%%%%%%%%%%%%%%%%%%%%%%%%%%%%%%%%%%%%%%%%%%%%%%%%%%%%%

\[
\text{\bf Model 2c3b4a}
\ \ \ \ \
\left\{
\aligned
u
&
\,=\,
xy+y^3
,
\\
v
&
\,=\,
0,
\endaligned\right.
\]
%%%%%%%%%%%%%%%%%%%%%%%%%%%%%%%%%%%%%%%%%%%%%%%%%%%%%%%%%%%%%%%%%%%%%%
\[
\def\arraystretch{1.25}
\begin{array}{llll}
e_1
\,:=\,
\partial_x+y\partial_u
, &
e_2
\,:=\,
-3y\partial_x+\partial_y+x\partial_u
, &
e_3
\,:=\,
2x\partial_x+y\partial_y+3u\partial_u
,
\\
e_4
\,:=\,
v\partial_x
, &
e_5
\,:=\,
v\partial_y
, &
e_6
\,:=\,
v\partial_u
, &
e_7
\,:=\,
v\partial_v,
\end{array}
\]
%%%%%%%%%%%%%%%%%%%%%%%%%%%%%%%%%%%%%%%%%%%%%%%%%%%%%%%%%%%%%%%%%%%%%%
\[
\footnotesize
\def\arraystretch{1.25}
\begin{array}{c|ccccccc}
{} & e_1 & e_2 & e_3 & e_4 & e_5 & e_6 & e_7
\\
\hline
e_1 & 0 & 0 & 2e_1 & 0 & -e_6 & 0 & 0
\\
e_2 & 0 & 0 & e_2 & -e_6 & 3e_4 & 0 & 0
\\
e_3 & -2e_1 & -e_2 & 0 & -2e_4 & -e_5 & -3e_6 & 0
\\
e_4 & 0 & e_6 & 2e_4 & 0 & 0 & 0 & -e_4
\\
e_5 & e_6 & -3e_4 & e_5 & 0 & 0 & 0 & -e_5
\\
e_6 & 0 & 0 & 3e_6 & 0 & 0 & 0 & -e_6
\\
e_7 & 0 & 0 & 0 & e_4 & e_5 & e_6 & 0
\end{array}
\]
%%%%%%%%%%%%%%%%%%%%%%%%%%%%%%%%%%%%%%%%%%%%%%%%%%%%%%%%%%%%%%%%%%%%%%
%%%%%%%%%%%%%%%%%%%%%%%%%%%%%%%%%%%%%%%%%%%%%%%%%%%%%%%%%%%%%%%%%%%%%%
%%%%%%%%%%%%%%%%%%%%%%%%%%%%%%%%%%%%%%%%%%%%%%%%%%%%%%%%%%%%%%%%%%%%%%

\[
\text{\bf Model 2c3b4b}
\ \ \ \ \
\left\{
\aligned
u
&
\,=\,
xy+y^3+xy^3+\tfrac{6}{5}y^5+\tfrac{6}{5}xy^5+\tfrac{51}{35}y^7+\tfrac{51}{35}xy^7+\tfrac{62}{35}y^9
+
\cdots,
\\
v
&
\,=\,
0,
\endaligned\right.
\]
%%%%%%%%%%%%%%%%%%%%%%%%%%%%%%%%%%%%%%%%%%%%%%%%%%%%%%%%%%%%%%%%%%%%%%
\[
\def\arraystretch{1.25}
\begin{array}{lll}
e_1
\,:=\,
(x+1)\partial_x
+
(u+y)\partial_u
, &
e_2
\,:=\,
-(3u+3y)\partial_x+\partial_y+x\partial_u
, &
e_3
\,:=\,
v\partial_x,
\\
e_4
\,:=\,
v\partial_y
, &
e_5
\,:=\,
v\partial_u
, &
e_6
\,:=\,
v\partial_v,
\end{array}
\]
%%%%%%%%%%%%%%%%%%%%%%%%%%%%%%%%%%%%%%%%%%%%%%%%%%%%%%%%%%%%%%%%%%%%%%
\[
\footnotesize
\def\arraystretch{1.25}
\begin{array}{c|cccccc}
{} & e_1 & e_2 & e_3 & e_4 & e_5 & e_6 
\\
\hline
e_1 & 0 & 0 & -e_3 & -e_5 & -e_5 & 0
\\
e_2 & 0 & 0 & -e_5 & 3e_3 & 3e_3 & 0
\\
e_3 & e_3 & e_5 & 0 & 0 & 0 & -e_3
\\
e_4 & e_5 & -3e_3 & 0 & 0 & 0 & -e_4
\\
e_5 & e_5 & -3e_3 & 0 & 0 & 0 & -e_5
\\
e_6 & 0 & 0 & e_3 & e_4 & e_5 & 0
\end{array}
\]

%%%%%%%%%%%%%%%%%%%%%%%%%%%%%%%%%%%%%%%%%%%%%%%%%%%%%%%%%%%%%%%%%%%%%%
%%%%%%%%%%%%%%%%%%%%%%%%%%%%%%%%%%%%%%%%%%%%%%%%%%%%%%%%%%%%%%%%%%%%%%
%%%%%%%%%%%%%%%%%%%%%%%%%%%%%%%%%%%%%%%%%%%%%%%%%%%%%%%%%%%%%%%%%%%%%%

\[
\text{\bf Model 2c3b4c}
\ \ \ \ \
\left\{
\aligned
u
&
\,=\,
xy+y^3-xy^3-\tfrac{6}{5}y^5+\tfrac{6}{5}xy^5
+\tfrac{51}{35}y^7-\tfrac{51}{35}xy^7-\tfrac{62}{35}y^9
+
\cdots,
\\
v
&
\,=\,
0,
\endaligned\right.
\]
%%%%%%%%%%%%%%%%%%%%%%%%%%%%%%%%%%%%%%%%%%%%%%%%%%%%%%%%%%%%%%%%%%%%%%
\[
\def\arraystretch{1.25}
\begin{array}{lll}
e_1
\,:=\,
-(x-1)\partial_x-(u-y)\partial_u
, &
e_2
\,:=\,
(3u-3y)\partial_x+\partial_y+x\partial_u
, &
e_3
\,:=\,
v\partial_x
,
\\
e_4
\,:=\,
v\partial_y
, &
e_5
\,:=\,
v\partial_u
, &
e_6
\,:=\,
v\partial_v,
\end{array}
\]
%%%%%%%%%%%%%%%%%%%%%%%%%%%%%%%%%%%%%%%%%%%%%%%%%%%%%%%%%%%%%%%%%%%%%%
\[
\footnotesize
\def\arraystretch{1.25}
\begin{array}{c|cccccc}
{} & e_1 & e_2 & e_3 & e_4 & e_5 & e_6 
\\
\hline
e_1 & 0 & 0 & e_3 & -e_5 & e_5 & 0
\\
e_2 & 0 & 0 & -e_5 & 3e_3 & -3e_3 & 0
\\
e_3 & -e_3 & e_5 & 0 & 0 & 0 & -e_3
\\
e_4 & e_5 & -3e_3 & 0 & 0 & 0 & -e_4
\\
e_5 & -e_5 & 3e_3 & 0 & 0 & 0 & -e_5
\\
e_6 & 0 & 0 & e_3 & e_4 & e_5 & 0
\end{array}
\]

%%%%%%%%%%%%%%%%%%%%%%%%%%%%%%%%%%%%%%%%%%%%%%%%%%%%%%%%%%%%%%%%%%%%%%
%%%%%%%%%%%%%%%%%%%%%%%%%%%%%%%%%%%%%%%%%%%%%%%%%%%%%%%%%%%%%%%%%%%%%%
%%%%%%%%%%%%%%%%%%%%%%%%%%%%%%%%%%%%%%%%%%%%%%%%%%%%%%%%%%%%%%%%%%%%%%

\[
\text{\bf Model 2c3b4e}
\ \ \ \ \
\left\{
\aligned
u
&
\,=\,
xy+y^3+y^4+F_{1,3}xy^3+F_{0,5}y^5
+
(
\tfrac{5}{2}F_{0,5}F_{1,3}-3F_{1,3}^2-2F_{1,3}
)xy^4
+
\\
&
\ \ \ \ \ 
+
(
-3F_{1,3}^2+\tfrac{5}{2}F_{0,5}F_{1,3}^2-3F_{1,3}^3
)x^2y^3
+
\cdots,
\\
v
&
\,=\,
0,
\endaligned\right.
\]
%%%%%%%%%%%%%%%%%%%%%%%%%%%%%%%%%%%%%%%%%%%%%%%%%%%%%%%%%%%%%%%%%%%%%%

\noindent
Gr\"obner basis generators of 
moduli space core algebraic variety in $\R^2 \ni F_{1,3},F_{0,5}$:
\[
\aligned
\B_1 & := F_{1,3}^2(5F_{0,5}-6F_{1,3}-6),
\\
\B_2 & := F_{1,3}(5F_{0,5}-8+6F_{1,3})(5F_{0,5}-6F_{1,3}-6),
\endaligned
\]

%%%%%%%%%%%%%%%%%%%%%%%%%%%%%%%%%%%%%%%%%%%%%%%%%%%%%%%%%%%%%%%%%%%%%%
\[
\def\arraystretch{1.25}
\begin{array}{ll}
e_1
&
\,:=\,
-
(
5F_{0,5}F_{1,3}x-6F_{1,3}^2x-7F_{1,3}x-1
)\partial_x
-
(
-3F_{1,3}y+\tfrac{5}{2}F_{0,5}F_{1,3}y-3F_{1,3}^2y
)\partial_y
+
\\
&
\ \ \ \ \ 
-
(
-y-10F_{1,3}u+\tfrac{15}{2}F_{0,5}F_{1,3}u-9F_{1,3}^3u
)\partial_u
,
\\
e_2
&
\,:=\,
-
(
10F_{0,5}x+3F_{1,3}u-12F_{1,3}x-12x+3y
)\partial_x
-
(
5F_{0,5}y-6F_{1,3}y-4y-1
)\partial_y
\\
&
\ \ \ \ \ 
-
(
15F_{0,5}u-18F_{1,3}u-16u-x
)\partial_u,
\end{array}
\]
\[
\def\arraystretch{1.25}
\begin{array}{llll}
e_3
\,:=\,
v\partial_x
, &
e_4
\,:=\,
v\partial_y
, &
e_5
\,:=\,
v\partial_u
, &
e_6
\,:=\,
v\partial_v,
\end{array}
\]

%%%%%%%%%%%%%%%%%%%%%%%%%%%%%%%%%%%%%%%%%%%%%%%%%%%%%%%%%%%%%%%%%%%%%%
\[
\tiny
\def\arraystretch{1.25}
\begin{array}{c|cccccc}
{} & e_1 & e_2 & e_3 & e_4 & e_5 & e_6 
\\
\hline
e_1 & 0 &
\rotatebox[origin=c]{90}{
\begin{tabular}{p{3cm}} $(-10F_{0,5}+12F_{1,3}+12)e_1+(-3F_{1,3}+\tfrac{5}{2}F_{0,5}F_{1,3}-3F_{1,3}^2)e_2$
\end{tabular} 
}
& 
\rotatebox[origin=c]{90}{
\begin{tabular}{p{3cm}}
$F_{1,3}(5F_{0,5}-6F_{1,3}-7)e_3$
\end{tabular}
}
&
\rotatebox[origin=c]{90}{
\begin{tabular}{p{3cm}}
 $\tfrac{1}{2}F_{1,3}(-6+5F_{0,5}-6F_{1,3})e_4-e_5$
\end{tabular}
}
& 
\rotatebox[origin=c]{90}{
\begin{tabular}{p{3cm}}
$\tfrac{1}{2}F_{1,3}(-20+15F_{0,5}-18F_{1,3})e_5$
\end{tabular}
}
& 0
\\
e_2 & 
\rotatebox[origin=c]{90}{
\begin{tabular}{p{3cm}}
$-(-10F_{0,5}+12F_{1,3}+12)e_1-(-3F_{1,3}+\tfrac{5}{2}F_{0,5}F_{1,3}-3F_{1,3}^2)e_2$
\end{tabular}
}
&
0
&
\rotatebox[origin=c]{90}{
\begin{tabular}{p{3cm}}
$(-12+10F_{0,5}-12F_{1,3})e_3-e_5$
\end{tabular}
}
&
\rotatebox[origin=c]{90}{
\begin{tabular}{p{3cm}}
$3e_3+(5F_{0,5}-6F_{1,3}-4)e_4$
\end{tabular}
}
&
\rotatebox[origin=c]{90}{
\begin{tabular}{p{3cm}}
$3F_{1,3}e_3+(15F_{0,5}-18F_{1,3}-16)e_5$
\end{tabular}
}
& 0
\\
e_3 
&
\rotatebox[origin=c]{90}{
\begin{tabular}{p{3cm}}
$-F_{1,3}(5F_{0,5}-6F_{1,3}-7)e_3$
\end{tabular}
}
&
\rotatebox[origin=c]{90}{
\begin{tabular}{p{3cm}}
$-(-12+10F_{0,5}-12F_{1,3})e_3+e_5$
\end{tabular}
}
& 0 & 0 & 0 & -e_3
\\
e_4 
&
\rotatebox[origin=c]{90}{
\begin{tabular}{p{3cm}}
$-\tfrac{1}{2}F_{1,3}(-6+5F_{0,5}-6F_{1,3})e_4+e_5$
\end{tabular}
}
&
\rotatebox[origin=c]{90}{
\begin{tabular}{p{3cm}}
$-3e_3-(5F_{0,5}-6F_{1,3}-4)e_4$
\end{tabular}
}
& 0 & 0 & 0 & -e_4
\\
e_5 
&
\rotatebox[origin=c]{90}{
\begin{tabular}{p{3cm}}
$-\tfrac{1}{2}F_{1,3}(-20+15F_{0,5}-18F_{1,3})e_5$
\end{tabular}
}
&
\rotatebox[origin=c]{90}{
\begin{tabular}{p{3cm}}
$-3F_{1,3}e_3-(15F_{0,5}-18F_{1,3}-16)e_5$
\end{tabular}
}
& 0 & 0 & 0 & -e_5
\\
e_6 & 0 & 0 & e_3 & e_4 & e_5 & 0
\end{array}
\]

%%%%%%%%%%%%%%%%%%%%%%%%%%%%%%%%%%%%%%%%%%%%%%%%%%%%%%%%%%%%%%%%%%%%%%
%%%%%%%%%%%%%%%%%%%%%%%%%%%%%%%%%%%%%%%%%%%%%%%%%%%%%%%%%%%%%%%%%%%%%%
%%%%%%%%%%%%%%%%%%%%%%%%%%%%%%%%%%%%%%%%%%%%%%%%%%%%%%%%%%%%%%%%%%%%%%

\[
\text{\bf Model 2c3d}
\ \ \ \ \
\left\{
\aligned
u
&
\,=\,
xy+y^3+x^3+F_{4,0}y^4+F_{1,3}xy^3+
(
\tfrac{2}{9}F_{1,3}F_{0,4}-\tfrac{1}{9}F_{3,1}F_{1,3}+\tfrac{8}{9}F_{4,0}F_{0,4}
)x^2y^2
+
\\
&
\ \ \ \ \ 
+
F_{3,1}x^3y
+
F_{4,0}x^4
+
\cdots
,
\\
v
&
\,=\,
0,
\endaligned\right.
\]

%%%%%%%%%%%%%%%%%%%%%%%%%%%%%%%%%%%%%%%%%%%%%%%%%%%%%%%%%%%%%%%%%%%%%%

\noindent
Gr\"obner basis generators of 
moduli space core algebraic variety in 
$\R^4 \ni F_{1,3},F_{3,1},F_{0,4},F_{4,0}$:
\[
\aligned
\B_1 & := F_{0,4}F_{1,3}-F_{3,1}F_{4,0},
\\
\B_2 & := -2F_{0,4}F_{3,1}^2+F_{1,3}^3+2F_{1,3}^2F_{4,0}-F_{3,1}^3,
\\
\B_3 & := 2F_{0,4}^2F_{3,1}+F_{0,4}F_{3,1}^2-F_{1,3}^2F_{4,0}-2F_{1,3}F_{4,0}^2,
\\
\B_4 & := -32F_{0,4}F_{4,0}^2+F_{1,3}^2F_{3,1}+2F_{1,3}F_{3,1}F_{4,0}-16F_{3,1}F_{4,0}^2+36
F_{0,4}F_{3,1}+18F_{3,1}^2+81F_{1,3}+162F_{4,0},
\\
\B_5 & :=  -32F_{0,4}^2F_{4,0}-16F_{0,4}F_{3,1}F_{4,0}+F_{1,3}F_{3,1}^2+2F_{3,1}^2F_{4,0}+18
F_{1,3}^2+36F_{1,3}F_{4,0}+162F_{0,4}+81F_{3,1},
\endaligned
\]

%%%%%%%%%%%%%%%%%%%%%%%%%%%%%%%%%%%%%%%%%%%%%%%%%%%%%%%%%%%%%%%%%%%%%%
\[
\def\arraystretch{1.25}
\begin{array}{llll}
e_1
&
\,:=\,
-
(
-1
+
\tfrac{1}{3}F_{1,3}x
+
\tfrac{8}{3}F_{4,0}x
-
9u
+
\tfrac{4}{9}F_{1,3}F_{0,4}u
-
\tfrac{2}{9}F_{3,1}F_{1,3}
+
\tfrac{16}{9}F_{4,0}F_{0,4}u
)\partial_x
-
(
3x
+
\tfrac{4}{3}F_{4,0}y
\\
&
\ \ \ \ \ 
+
\tfrac{2}{3}F_{1,3}y
+
3F_{3,1}u
)\partial_y
-
(
F_{1,3}u
+
4F_{4,0}u
-
y
)\partial_u
,
\\
e_2
&
\,:=\,
-
(
\tfrac{4}{3}F_{0,4}x
+
\tfrac{2}{3}F_{3,1}x
+
3y
+
3F_{1,3}u
)\partial_x
-
(
-1
+
\tfrac{1}{3}F_{3,1}y
+
\tfrac{8}{3}F_{0,4}y
-
9u
+
\tfrac{4}{9}F_{1,3}F_{0,4}u
+
\\
&
\ \ \ \ \ 
-
\tfrac{2}{9}F_{3,1}F_{1,3}u
+
\tfrac{16}{9}F_{4,0}F_{0,4}u
)\partial_y
-
(
4F_{0,4}u
+
F_{3,1}u
-
x
)\partial_u
, 
\\
e_3
&
\,:=\,
v\partial_x
,
\ \ \ \ \
e_4
\,:=\,
v\partial_y
,
\ \ \ \ \
e_5
\,:=\,
v\partial_u
,
\ \ \ \ \
e_6
\,:=\,
v\partial_v,
\end{array}
\]

%%%%%%%%%%%%%%%%%%%%%%%%%%%%%%%%%%%%%%%%%%%%%%%%%%%%%%%%%%%%%%%%%%%%%%
\[
\footnotesize
\def\arraystretch{1.25}
\begin{array}{c|cccccc}
{} & e_1 & e_2 & e_3 & e_4 & e_5 & e_6 
\\
\hline
e_1 & 0 &
\rotatebox[origin=c]{90}{
\begin{tabular}{p{3cm}} $(-\tfrac{4}{3}F_{0,4}-\tfrac{2}{3}F_{3,1})e_1+(\tfrac{2}{3}F_{1,3}+\tfrac{4}{3}F_{4,0})e_2$
\end{tabular}
}
& 
\rotatebox[origin=c]{90}{
\begin{tabular}{p{3cm}}
$(\tfrac{1}{3}F_{1,3}+\tfrac{8}{3}F_{4,0})e_3+3e_4$
\end{tabular}
}
&
\rotatebox[origin=c]{90}{
\begin{tabular}{p{3cm}}
$(\tfrac{2}{3}F_{1,3}+\tfrac{4}{3}F_{4,0})e_4-e_5$
\end{tabular}
}
&
\rotatebox[origin=c]{90}{
\begin{tabular}{p{3cm}}
$(\tfrac{4}{9}F_{1,3}F_{0,4}+\tfrac{16}{9}F_{4,0}F_{0,4}-\tfrac{2}{9}F_{3,1}F_{1,3}-9)e_3+3F_{3,1}e_4+(F_{1,3}+4F_{4,0})e_5$
\end{tabular}
}
& 0
\\
e_2 & \rotatebox[origin=c]{90}{
\begin{tabular}{p{3cm}}
$(\tfrac{4}{3}F_{0,4}+\tfrac{2}{3}F_{3,1})e_1-(\tfrac{2}{3}F_{1,3}+\tfrac{4}{3}F_{4,0})e_2$
\end{tabular}
}
& 0 &
\rotatebox[origin=c]{90}{
\begin{tabular}{p{3cm}}
$(\tfrac{4}{3}F_{0,4}+\tfrac{2}{3}F_{3,1})e_3-e_5$
\end{tabular}
}
&
\rotatebox[origin=c]{90}{
\begin{tabular}{p{3cm}}
$3e_3+(\tfrac{1}{3}F_{3,1}+\tfrac{8}{3}F_{0,4})e_4$
\end{tabular}
}
& 
\rotatebox[origin=c]{90}{
\begin{tabular}{p{3cm}}
$3F_{1,3}e_3+(\tfrac{4}{9}F_{1,3}F_{0,4}+\tfrac{16}{9}F_{4,0}F_{0,4}-\tfrac{2}{9}F_{3,1}F_{1,3}-9)e_4+(4F_{0,4}+F_{3,1})e_5$
\end{tabular}
}
& 0
\\
e_3 
&
\rotatebox[origin=c]{90}{
\begin{tabular}{p{3cm}}
$-(\tfrac{1}{3}F_{1,3}+\tfrac{8}{3}F_{4,0})e_3-3e_4$
\end{tabular}
}
&
\rotatebox[origin=c]{90}{
\begin{tabular}{p{3cm}}
$-(\tfrac{4}{3}F_{0,4}+\tfrac{2}{3}F_{3,1})e_3+e_5$
\end{tabular}
}
& 0 & 0 & 0 & -e_3
\\
e_4 & 
\rotatebox[origin=c]{90}{
\begin{tabular}{p{3cm}}
$-(\tfrac{2}{3}F_{1,3}+\tfrac{4}{3}F_{4,0})e_4+e_5$
\end{tabular}
}
&
\rotatebox[origin=c]{90}{
\begin{tabular}{p{3cm}}
$-3e_3-(\tfrac{1}{3}F_{3,1}+\tfrac{8}{3}F_{0,4})e_4$
\end{tabular}
}
& 0 & 0 & 0 & - e_4
\\
e_5 
&
\rotatebox[origin=c]{90}{
\begin{tabular}{p{3cm}}
$-(\tfrac{4}{9}F_{1,3}F_{0,4}-\tfrac{2}{9}F_{3,1}F_{1,3}-9)e_3-3F_{3,1}e_4-(F_{1,3}+4F_{4,0})e_5$
\end{tabular}
}
& 
\rotatebox[origin=c]{90}{
\begin{tabular}{p{3cm}}
$-3F_{1,3}e_3-(\tfrac{4}{9}F_{1,3}F_{0,4}+\tfrac{16}{9}F_{4,0}F_{0,4}-\tfrac{2}{9}F_{3,1}F_{1,3}-9)e_4-(4F_{0,4}+F_{3,1})e5$
\end{tabular}
}
& 0 & 0 & 0 & -e_5
\\
e_6 & 0 & 0 & e_3 & e_4 & e_5 & 0
\end{array}
\]

%%%%%%%%%%%%%%%%%%%%%%%%%%%%%%%%%%%%%%%%%%%%%%%%%%%%%%%%%%%%%%%%%%%%%%
%%%%%%%%%%%%%%%%%%%%%%%%%%%%%%%%%%%%%%%%%%%%%%%%%%%%%%%%%%%%%%%%%%%%%%
%%%%%%%%%%%%%%%%%%%%%%%%%%%%%%%%%%%%%%%%%%%%%%%%%%%%%%%%%%%%%%%%%%%%%%

%%%%%%%%%%%%%%%%%%%%%%%%%%%%%%%%%%%%%%%%%%%%%%%%%%%%%%%%%%%%%%%%%%%%%%
\SectionHead{2d Models}
{2d-models}
%%%%%%%%%%%%%%%%%%%%%%%%%%%%%%%%%%%%%%%%%%%%%%%%%%%%%%%%%%%%%%%%%%%%%%

%%%%%%%%%%%%%%%%%%%%%%%%%%%%%%%%%%%%%%%%%%%%%%%%%%%%%%%%%%%%%%%%%%%%%%
%%%%%%%%%%%%%%%%%%%%%%%%%%%%%%%%%%%%%%%%%%%%%%%%%%%%%%%%%%%%%%%%%%%%%%
%%%%%%%%%%%%%%%%%%%%%%%%%%%%%%%%%%%%%%%%%%%%%%%%%%%%%%%%%%%%%%%%%%%%%%

\[
\text{\bf Model 2d3a4a}
\ \ \ \ \
\left\{
\aligned
u
&
\,=\,
x^2+y^2,
\\
v
&
\,=\,
0,
\endaligned\right.
\]
%%%%%%%%%%%%%%%%%%%%%%%%%%%%%%%%%%%%%%%%%%%%%%%%%%%%%%%%%%%%%%%%%%%%%%
\[
\def\arraystretch{1.25}
\begin{array}{llll}
e_1
\,:=\,
\partial_x
+
2x\partial_u
, &
e_2
\,:=\,
\partial_y
+
2y\partial_u
, &
e_3
\,:=\,
x\partial_x
+
y\partial_y
+
2u\partial_u
, &
e_4
\,:=\,
y\partial_x
-
x\partial_y
,
\\
e_5
\,:=\,
v\partial_x
, &
e_6
\,:=\,
v\partial_y
, &
e_7
\,:=\,
v\partial_u
, &
e_8
\,:=\,
v\partial_v,
\end{array}
\]

%%%%%%%%%%%%%%%%%%%%%%%%%%%%%%%%%%%%%%%%%%%%%%%%%%%%%%%%%%%%%%%%%%%%%%
\[
\footnotesize
\def\arraystretch{1.25}
\begin{array}{c|cccccccc}
{} & e_1 & e_2 & e_3 & e_4 & e_5 & e_6 & e_7 & e_8 
\\
\hline
e_1 & 0 & 0 & e_1 & -e_2 & -2e_7 & 0 & 0 & 0
\\
e_2 & 0 & 0 & e_2 & e_1 & 0 & -2e_7 & 0 & 0
\\
e_3 & -e_1 & -e_2 & 0 & 0 & -e_5 & -e_6 & -2e_7 & 0
\\
e_4 & e_2 & -e_1 & 0 & 0 & e_6 & -e_5 & 0 & 0
\\
e_5 & 2e_7 & 0 & e_5 & -e_6 & 0 & 0 & 0 & -e_5
\\
e_6 & 0 & 2e_7& e_6 & e_5 & 0 & 0 & 0 & -e_6
\\
e_7 & 0 & 0 & 2e_7 & 0 & 0 & 0 & 0 & -e_7
\\
e_8 & 0 & 0 & 0 & 0 & e_5 & e_6 & e_7 & 0
\end{array}
\]

%%%%%%%%%%%%%%%%%%%%%%%%%%%%%%%%%%%%%%%%%%%%%%%%%%%%%%%%%%%%%%%%%%%%%%
%%%%%%%%%%%%%%%%%%%%%%%%%%%%%%%%%%%%%%%%%%%%%%%%%%%%%%%%%%%%%%%%%%%%%%
%%%%%%%%%%%%%%%%%%%%%%%%%%%%%%%%%%%%%%%%%%%%%%%%%%%%%%%%%%%%%%%%%%%%%%

\[
\text{\bf Model 2d3a4b}
\ \ \ \ \
\left\{
\aligned
u
&
\,=\,
x^2+y^2
+
\tfrac{1}{2}y^4
+
x^2y^2
+
\tfrac{1}{2}x^4
+
\tfrac{1}{2}y^6
+
\tfrac{3}{2}x^2y^4
+
\tfrac{3}{2}x^4y^2
+
\tfrac{1}{6}x^6
+
\\
&
\ \ \ \ \ 
+
\tfrac{5}{8}y^8
+
\tfrac{5}{2}x^2y^6
+
\tfrac{15}{4}x^4y^4
+
\tfrac{5}{2}x^6y^2
+
\tfrac{5}{8}x^8
+
\cdots,
\\
v
&
\,=\,
0,
\endaligned\right.
\]
%%%%%%%%%%%%%%%%%%%%%%%%%%%%%%%%%%%%%%%%%%%%%%%%%%%%%%%%%%%%%%%%%%%%%%
\[
\def\arraystretch{1.25}
\begin{array}{llll}
e_1
\,:=\,
(1-u)\partial_x
+
2x\partial_u
, &
e_2
\,:=\,
(1-u)\partial_y
+
2y\partial_u
, &
e_3
\,:=\,
y\partial_x
-
x\partial_y
, 
\\
e_4
\,:=\,
v\partial_x
, &
e_5
\,:=\,
v\partial_y
, &
e_6
\,:=\,
v\partial_u
, &
e_7
\,:=\,
v\partial_v,
\end{array}
\]

%%%%%%%%%%%%%%%%%%%%%%%%%%%%%%%%%%%%%%%%%%%%%%%%%%%%%%%%%%%%%%%%%%%%%%
\[
\footnotesize
\def\arraystretch{1.25}
\begin{array}{c|ccccccc}
{} & e_1 & e_2 & e_3 & e_4 & e_5 & e_6 & e_7
\\
\hline
e_1 & 0 & 2e_3 & -e_2 & -2e_6 & 0 & e_4 & 0
\\
e_2 & -2e_3 & 0 & e_1 & 0 & -2e_6 & e_5 & 0
\\
e_3 & e_2 & -e_1 & 0 & e_5 & -e_4 & 0 & 0
\\
e_4 & 2e_6 & 0 & -e_5 & 0 & 0 & 0 & -e_4
\\
e_5 & 0 & 2e_6 & e_4 & 0 & 0 & 0 & -e_5
\\
e_6 & -e_4 & -e_5 & 0 & 0 & 0 & 0 & -e_6
\\
e_7 & 0 & 0 & 0 & e_4 & e_5 & e_6 & 0
\end{array}
\]

%%%%%%%%%%%%%%%%%%%%%%%%%%%%%%%%%%%%%%%%%%%%%%%%%%%%%%%%%%%%%%%%%%%%%%
%%%%%%%%%%%%%%%%%%%%%%%%%%%%%%%%%%%%%%%%%%%%%%%%%%%%%%%%%%%%%%%%%%%%%%
%%%%%%%%%%%%%%%%%%%%%%%%%%%%%%%%%%%%%%%%%%%%%%%%%%%%%%%%%%%%%%%%%%%%%%

\[
\text{\bf Model 2d3a4c}
\ \ \ \ \
\left\{
\aligned
u
&
\,=\,
x^2+y^2
-
\tfrac{1}{2}y^4
-
x^2y^2
-
\tfrac{1}{2}x^4
+
\tfrac{1}{2}y^6
+
\tfrac{3}{2}x^2y^4
+
\tfrac{3}{2}x^4y^2
+
\tfrac{1}{6}x^6
+
\\
&
\ \ \ \ \ 
-
\tfrac{5}{8}y^8
-
\tfrac{5}{2}x^2y^6
-
\tfrac{15}{4}x^4y^4
-
\tfrac{5}{2}x^6y^2
-
\tfrac{5}{8}x^8
+
\cdots,
\\
v
&
\,=\,
0,
\endaligned\right.
\]
%%%%%%%%%%%%%%%%%%%%%%%%%%%%%%%%%%%%%%%%%%%%%%%%%%%%%%%%%%%%%%%%%%%%%%
\[
\def\arraystretch{1.25}
\begin{array}{llll}
e_1
\,:=\,
(1+u)\partial_x
+
2x\partial_u
, &
e_2
\,:=\,
(1+u)\partial_y
+
2y\partial_u
, &
e_3
\,:=\,
y\partial_x
-
x\partial_y
, 
\\
e_4
\,:=\,
v\partial_x
, &
e_5
\,:=\,
v\partial_y
, &
e_6
\,:=\,
v\partial_u
, &
e_7
\,:=\,
v\partial_v,
\end{array}
\]

%%%%%%%%%%%%%%%%%%%%%%%%%%%%%%%%%%%%%%%%%%%%%%%%%%%%%%%%%%%%%%%%%%%%%%
\[
\footnotesize
\def\arraystretch{1.25}
\begin{array}{c|ccccccc}
{} & e_1 & e_2 & e_3 & e_4 & e_5 & e_6 & e_7
\\
\hline
e_1 & 0 & -2e_3 & -e_2 & -2e_6 & 0 & -e_4 & 0
\\
e_2 & 2e_3 & 0 & e_1 & 0 & -2e_6 & -e_5 & 0
\\
e_3 & e_2 & -e_1 & 0 & e_5 & -e_4 & 0 & 0
\\
e_4 & 2e_6 & 0 & -e_5 & 0 & 0 & 0 & -e_4
\\
e_5 & 0 & 2e_6 & e_4 & 0 & 0 & 0 & -e_5
\\
e_6 & e_4 & e_5 & 0 & 0 & 0 & 0 & -e_6
\\
e_7 & 0 & 0 & 0 & e_4 & e_5 & e_6 & 0
\end{array}
\]

%%%%%%%%%%%%%%%%%%%%%%%%%%%%%%%%%%%%%%%%%%%%%%%%%%%%%%%%%%%%%%%%%%%%%%
%%%%%%%%%%%%%%%%%%%%%%%%%%%%%%%%%%%%%%%%%%%%%%%%%%%%%%%%%%%%%%%%%%%%%%
%%%%%%%%%%%%%%%%%%%%%%%%%%%%%%%%%%%%%%%%%%%%%%%%%%%%%%%%%%%%%%%%%%%%%%

\[
\text{\bf Model 2d3b}
\ \ \ \ \
\left\{
\aligned
u
&
\,=\,
y^2
+
x^2
+
x^3
+
F_{0,4}y^4
+
F_{1,3}xy^3
+
F_{2,2}x^2y^2
+
F_{3,1}x^3y
+
F_{4,0}x^4
+
\cdots,
\\
v
&
\,=\,
0,
\endaligned\right.
\]
%%%%%%%%%%%%%%%%%%%%%%%%%%%%%%%%%%%%%%%%%%%%%%%%%%%%%%%%%%%%%%%%%%%%%%
\noindent
Gr\"obner basis generators of 
moduli space core algebraic variety in 
$\R^5 \ni F_{0,4},F_{4,0},F_{1,3},F_{3,1},F_{2,2}$:
\[
\aligned
\B_1 
&
:= 
32F_{0,4}F_{3,1}
+
4F_{1,3}F_{2,2}
+
136F_{1,3}F_{4,0}
-
20F_{2,2}
F_{3,1}
+
24F_{3,1}F_{4,0}
-
174F_{1,3}
-
27F_{3,1},
\\
\B_2
&
:= -4F_{0,4}F_{2,2}
+
8F_{0,4}F_{4,0}
+
F_{1,3}^2
+
5F_{1,3}F_{3,1}
-
2F_{2,2}^2
+
4F_{2,2}F_{4,0}
-
6F_{0,4}
-
3F_{2,2},
\\
\B_3
& 
:= 8F_{0,4}F_{1,3}
-
24F_{0,4}F_{3,1}
+
88F_{1,3}F_{4,0}
-
16F_{2,2}F_{3,1}
+
24F_{3,1}F_{4,0}
-
111F_{1,3}
-
27F_{3,1},
\\
\B_4
&
:= 8F_{0,4}^2
-
12F_{0,4}F_{2,2}
+
8F_{0,4}F_{4,0}
+
32F_{1,3}F_{3,1}
-
12F_{2,2}^2
+
28F_{2,2}F_{4,0}+
12F_{0,4}
-
21F_{2,2},
\endaligned
\]
%%%%%%%%%%%%%%%%%%%%%%%%%%%%%%%%%%%%%%%%%%%%%%%%%%%%%%%%%%%%%%%%%%%%%%
\[
\aligned
e_1
&
\,:=\,
(
2F_{2,2}x
-
4F_{4,0}x
+
3x
+
1
+
\tfrac{1}{3}yF_{1,3}
-
yF_{3,1}-uF_{2,2}
)\partial_x
-
(
\tfrac{1}{3}F_{1,3}x
-
xF_{3,1}
-
\tfrac{9}{2}y
-
2yF_{2,2}
+
\\
&
\ \ \ \ \ 
+
4yF_{4,0}
+
\tfrac{1}{2}uF_{1,3}
)\partial_y
+
(
4F_{2,2}u-8F_{4,0}u+9u+2x
)\partial_u
,
\\
e_2
&
\,:=\,
(
3F_{1,3}x
-
xF_{3,1}
+
\tfrac{4}{3}yF_{0,4}
-
\tfrac{2}{3}yF_{2,2}
-
\tfrac{3}{2}uF_{1,3}
)\partial_x
-
(
\tfrac{4}{3}xF_{0,4}
-
1
-
\tfrac{2}{3}F_{2,2}x
-
3yF_{1,3}
+
yF_{3,1}
+
\\
&
\ \ \ \ \ 
+
2uF_{0,4}
)\partial_y
+
(
6F_{1,3}u-2F_{3,1}u+2y
)\partial_u
, 
\\
e_3
&
\,:=\,
v\partial_x
,
\ \ \ \ \ 
e_4
\,:=\,
v\partial_y
,
\ \ \ \ \ 
e_5
\,:=\,
v\partial_u
, 
\ \ \ \ \ 
e_6
\,:=\,
v\partial_v,
\endaligned
\]

%%%%%%%%%%%%%%%%%%%%%%%%%%%%%%%%%%%%%%%%%%%%%%%%%%%%%%%%%%%%%%%%%%%%%%
\[
\footnotesize
\def\arraystretch{1.25}
\begin{array}{c|cccccc}
{} & e_1 & e_2 & e_3 & e_4 & e_5 & e_6 
\\
\hline
e_1 
&
0
&
\rotatebox[origin=c]{90}{
\begin{tabular}{p{4cm}}
$\tfrac{8}{3}F_{1,3}e_1
+
(
-
\tfrac{9}{2}
-
\tfrac{4}{3}F_{2,2}
+
4F_{4,0}
-
\tfrac{4}{3}F_{0,4}
)e_2$
\end{tabular}}
&
\rotatebox[origin=c]{90}{
\begin{tabular}{p{4cm}}
$
(
-
2F_{2,2}
+
4F_{4,0}
-
3
)e_3
+
(
\tfrac{1}{3}F_{1,3}
-
F_{3,1}
)e_4
-
2e_5$
\end{tabular}}
&
\rotatebox[origin=c]{90}{
	\begin{tabular}{p{4cm}}
$
(
-
\tfrac{1}{3}F_{1,3}
+
F_{3,1}
)e_3
+
(
-
2F_{2,2}
+
4F_{4,0}
-
\tfrac{9}{2}
)e_4$
\end{tabular}}
&
\rotatebox[origin=c]{90}{
	\begin{tabular}{p{4cm}}
$
F_{2,2}e_3
+
\tfrac{1}{2}F_{1,3}e_4
+
(
-
4F_{2,2}
+
8F_{4,0}
-
9
)e_5$
\end{tabular}}
&
0
\\
e_2 &
\rotatebox[origin=c]{90}{
	\begin{tabular}{p{4cm}}
$
-
\tfrac{8}{3}F_{1,3}e_1
-
(
-
\tfrac{9}{2}
-
\tfrac{4}{3}F_{2,2}
+
4F_{4,0}
-
\tfrac{4}{3}F_{0,4}
)e_2$
\end{tabular}}
&
0
&
\rotatebox[origin=c]{90}{
	\begin{tabular}{p{4cm}}
$
(
-
3F_{1,3}
+
F_{3,1}
)e_3
+
(
\tfrac{4}{3}F_{0,4}
-
\tfrac{2}{3}F_{2,2}
)e_4$
\end{tabular}}
&
\rotatebox[origin=c]{90}{
	\begin{tabular}{p{4cm}}
$
(
-
\tfrac{4}{3}F_{0,4}
+
\tfrac{2}{3}F_{2,2}
)e_3
+
(
-
3F_{1,3}
+
F_{3,1}
)e_4
-
2e_5$
\end{tabular}}
&
\rotatebox[origin=c]{90}{
	\begin{tabular}{p{4cm}}
$
\tfrac{3}{2}F_{1,3}e_3
+
2F_{0,4}e_4
+
(
-
6F_{1,3}
+
2F_{3,1}
)e_5$
\end{tabular}}
&
0
\\
e_3 &
\rotatebox[origin=c]{90}{
	\begin{tabular}{p{4cm}}
$
-
(
-
2F_{2,2}
+
4F_{4,0}
-
3
)e_3
-
(
\tfrac{1}{3}F_{1,3}
-
F_{3,1}
)e_4
+
2e_5$
\end{tabular}}
&
\rotatebox[origin=c]{90}{
	\begin{tabular}{p{4cm}}
$
-
(
-
3F_{1,3}
+
F_{3,1}
)e_3
-
(
\tfrac{4}{3}F_{0,4}
-
\tfrac{2}{3}F_{2,2}
)e_4$
\end{tabular}}
&
0
&
0
&
0
&
-e_3
\\
e_4 &
\rotatebox[origin=c]{90}{
	\begin{tabular}{p{4cm}}
$
-
(
-
\tfrac{1}{3}
F_{1,3}
+
F_{3,1}
)e_3
-
(
-
2F_{2,2}
+
4F_{4,0}
-
\tfrac{9}{2}
)e_4
$
\end{tabular}}
&
\rotatebox[origin=c]{90}{
	\begin{tabular}{p{4cm}}
$
-
(
-
\tfrac{4}{3}F_{0,4}
+
\tfrac{2}{3}
F_{2,2}
)e_3
-
(
-
3F_{1,3}
+
F_{3,1}
)e_4
+
2e_5
$
\end{tabular}}
& 
0
&
0
&
0
&
-e_4
\\
e_5 &
\rotatebox[origin=c]{90}{
	\begin{tabular}{p{4cm}}
$
-
F_{2,2}e_3
-
\tfrac{1}{2}F_{1,3}e_4
-
(
-
4F_{2,2}
+
8F_{4,0}
-
9
)e_5
$
\end{tabular}}
&
\rotatebox[origin=c]{90}{
	\begin{tabular}{p{4cm}}
$
-
\tfrac{3}{2}F_{1,3}e_3
-
2F_{0,4}e_4
-
(
-
6F_{1,3}
+
2F_{3,1}
)e_5
$
\end{tabular}}
&
0
&
0
&
0
&
-e_5
\\
e_6 & 0 & 0 & e_3 & e_4 & e_5 & 0
\end{array}
\]

%%%%%%%%%%%%%%%%%%%%%%%%%%%%%%%%%%%%%%%%%%%%%%%%%%%%%%%%%%%%%%%%%%%%%%
%%%%%%%%%%%%%%%%%%%%%%%%%%%%%%%%%%%%%%%%%%%%%%%%%%%%%%%%%%%%%%%%%%%%%%
%%%%%%%%%%%%%%%%%%%%%%%%%%%%%%%%%%%%%%%%%%%%%%%%%%%%%%%%%%%%%%%%%%%%%%

%%%%%%%%%%%%%%%%%%%%%%%%%%%%%%%%%%%%%%%%%%%%%%%%%%%%%%%%%%%%%%%%%%%%%%
\SectionHead{2e Models}
{2e-models}
%%%%%%%%%%%%%%%%%%%%%%%%%%%%%%%%%%%%%%%%%%%%%%%%%%%%%%%%%%%%%%%%%%%%%%

%%%%%%%%%%%%%%%%%%%%%%%%%%%%%%%%%%%%%%%%%%%%%%%%%%%%%%%%%%%%%%%%%%%%%%
%%%%%%%%%%%%%%%%%%%%%%%%%%%%%%%%%%%%%%%%%%%%%%%%%%%%%%%%%%%%%%%%%%%%%%
%%%%%%%%%%%%%%%%%%%%%%%%%%%%%%%%%%%%%%%%%%%%%%%%%%%%%%%%%%%%%%%%%%%%%%

\[
\text{\bf Model 2e3a4a}
\ \ \ \ \
\left\{
\aligned
u
&
\,=\,
xy
+
y^3
+
F_{0,4}y^4
+
F_{1,3}xy^3
+
F_{2,2}x^2y^2
+
F_{3,1}x^3y
+
x^4
+
\\
&
\ \ \ \ \ 
+
(
\tfrac{9}{10}F_{1,3}
-
\tfrac{9}{250}F_{0,4}F_{1,3}F_{3,1}
+
\tfrac{1}{25}F_{0,4}F_{1,3}G_{4,0}
+
\tfrac{32}{25}F_{0,4}^2
+
\\
&
\ \ \ \ \ 
-
\tfrac{1}{250}F_{0,4}G_{4,0}F_{3,1}G_{3,1}
+
\tfrac{2}{375}F_{0,4}G_{4,0}F_{3,1}F_{2,2}
+
\tfrac{1}{200}F_{0,4}G_{3,1}^2
+
\\
&
\ \ \ \ \ 
+
\tfrac{1}{75}F_{0,4}F_{2,2}^2
-
\tfrac{1}{60}F_{0,4}G_{3,1}F_{2,2}
+
\tfrac{1}{25}F_{0,4}F_{3,1}F_{5,0}
+
\tfrac{1}{250}F_{0,4}F_{3,1}^2G_{3,1}
+
\\
&
\ \ \ \ \ 
-
\tfrac{2}{375}F_{0,4}F_{3,1}^2F_{2,2})y^5
+
(
\tfrac{3}{4}G_{3,1}
+
F_{2,2}
+
\tfrac{1}{10}F_{0,4}F_{3,1}G_{3,1}
+
\\
&
\ \ \ \ \ 
-
\tfrac{2}{15}F_{0,4}F_{3,1}F_{2,2}
-
\tfrac{1}{10}F_{0,4}G_{4,0}G_{3,1}
+
\tfrac{2}{15}F_{0,4}G_{4,0}F_{2,2}
+
\tfrac{8}{5}F_{0,4}F_{1,3}
+
\\
&
\ \ \ \ \ 
+
F_{0,4}F_{5,0})xy^4
+
(
6F_{3,1}
-
4G_{4,0}
+
\tfrac{1}{10}G_{3,1}F_{1,3}F_{3,1}
-
\tfrac{2}{15}F_{2,2}F_{1,3}F_{3,1}
+
\\
&
\ \ \ \ \ 
-
\tfrac{1}{10}G_{3,1}F_{1,3}G_{4,0}
+
\tfrac{2}{15}F_{2,2}F_{1,3}G_{4,0}
+
\tfrac{3}{5}F_{1,3}^2
+
F_{1,3}F_{5,0}
-
F_{0,4}G_{3,1}
+
\\
&
\ \ \ \ \ 
+
\tfrac{4}{3}F_{2,2}F_{0,4}
)x^2y^3
+
(
4
+
\tfrac{1}{10}G_{3,1}F_{3,1}F_{2,2}
-
\tfrac{1}{10}G_{3,1}G_{4,0}F_{2,2}
-
\tfrac{1}{2}G_{3,1}F_{1,3}
+
\\
&
\ \ \ \ \ 
+
\tfrac{14}{15}F_{2,2}F_{1,3}
+
F_{2,2}F_{5,0}
-
\tfrac{2}{15}F_{3,1}F_{2,2}^2
+
\tfrac{2}{15}G_{4,0}F_{2,2}^2)x^3y^2
+
\\
&
\ \ \ \ \ 
+
(
-
\tfrac{1}{10}G_{4,0}F_{3,1}G_{3,1}
+
\tfrac{2}{15}G_{4,0}F_{3,1}F_{2,2}
+
\tfrac{1}{8}G_{3,1}^2
+
\tfrac{1}{3}F_{2,2}^2
-
\tfrac{5}{12}G_{3,1}F_{2,2}
+
\\
&
\ \ \ \ \ 
+
\tfrac{1}{10}F_{3,1}F_{1,3}
+
F_{3,1}F_{5,0}
+
\tfrac{1}{10}F_{3,1}^2G_{3,1}
-
\tfrac{2}{15}F_{3,1}^2F_{2,2})x^4y
+
F_{5,0}x^5
+
\cdots
,
\\
v
&
\,=\,
x^2
-
\tfrac{3}{2}y^4
+
G_{3,1}x^3y
+
G_{4,0}x^4
-
\tfrac{18}{5}F_{0,4}y^5
-
3F_{1,3}xy^4
+
\\
&
\ \ \ \ \ 
+
(
-
2F_{2,2}
+
3G_{3,1}
)x^2y^3
+
(
-
\tfrac{1}{10}G_{3,1}F_{3,1}F_{2,2}
+
\tfrac{1}{10}G_{3,1}G_{4,0}F_{2,2}
+
\\
&
\ \ \ \ \ 
+
\tfrac{1}{5}G_{3,1}F_{1,3}
+
\tfrac{3}{4}G_{3,1}F_{5,0}
+
\tfrac{3}{40}F_{3,1}G_{3,1}^2
-
\tfrac{3}{40}G_{4,0}G_{3,1}^2
)x^4y
+
\\
&
\ \ \ \ \ 
+
(
\tfrac{2}{25}G_{4,0}F_{3,1}G_{3,1}
-
\tfrac{8}{75}G_{4,0}F_{3,1}F_{2,2}
-
\tfrac{1}{10}G_{3,1}^2
+
\tfrac{2}{25}G_{4,0}F_{1,3}
+
\\
&
\ \ \ \ \ 
+
\tfrac{4}{5}G_{4,0}F_{5,0}
-
\tfrac{2}{25}G_{4,0}^2G_{3,1}
+
\tfrac{8}{75}G_{4,0}^2F_{2,2}
+
\tfrac{2}{15}G_{3,1}F_{2,2}
)x^5
+
\cdots
.
\endaligned\right.
\]
%%%%%%%%%%%%%%%%%%%%%%%%%%%%%%%%%%%%%%%%%%%%%%%%%%%%%%%%%%%%%%%%%%%%%%
\noindent
Gr\"obner basis generators of 
moduli space core algebraic variety in 
$\R^7 \ni F_{0,4},F_{2,2},F_{1,3},F_{3,1},G_{3,1},G_{4,0},F_{5,0}$:
\[
\aligned
\B_1 
&
:= 16F_{0,4}^2F_{3,1}G_{3,1}
-
48F_{0,4}F_{3,1}^2
+
96F_{0,4}F_{3,1}G_{4,0}
+
12F_{1,3}F_{2,2}F_{3,1}
-
24F_{1,3}F_{2,2}G_{4,0}
+
\\
&
\ \ \ \ \ 
+
30F_{2,2}^2G_{3,1}
-
15F_{2,2}G_{3,1}^2
-
24F_{0,4}G_{3,1},
\\
\B_2
&
:= 180F_{0,4}F_{1,3}^2G_{3,1}
+
576F_{0,4}F_{2,2}F_{3,1}
+
864F_{0,4}F_{2,2}G_{4,0}
+
528F_{0,4}F_{3,1}G_{3,1}
+
\\
&
\ \ \ \ \ 
-
1448F_{0,4}G_{3,1}G_{4,0}
+
60F_{1,3}F_{2,2}^2
-
60F_{1,3}F_{2,2}G_{3,1}
+
375F_{1,3}G_{3,1}^2
-
432F_{0,4}F_{1,3}
+
\\
&
\ \ \ \ \ 
-
4320F_{0,4}F_{5,0}
-
2880F_{3,1}^2
+
9240F_{3,1}G_{4,0}
-
6960G_{4,0}^2
+
5040F_{2,2}
-
6840G_{3,1},
\\
\B_3
&
:= 
-
288F_{0,4}F_{2,2}F_{3,1}
+
1728F_{0,4}F_{2,2}G_{4,0}
+
996F_{0,4}F_{3,1}G_{3,1}
-
2216F_{0,4}G_{3,1}G_{4,0}
+
\\
&
\ \ \ \ \ 
+
540F_{1,3}^2F_{3,1}
-
1080F_{1,3}^2G_{4,0}
-
120F_{1,3}F_{2,2}^2
-
285F_{1,3}F_{2,2}G_{3,1}
+
465F_{1,3}G_{3,1}^2
+
\\
&
\ \ \ \ \ 
-
3024F_{0,4}F_{1,3}
+
2160F_{0,4}F_{5,0}
-
2340F_{3,1}^2
+
8520F_{3,1}G_{4,0}
-
7680G_{4,0}^2
-
10080F_{2,2}
+
1530G_{3,1},
\\
\B_4
&
:= 320F_{0,4}F_{1,3}G_{3,1}^2
+
480F_{1,3}^3G_{3,1}
+
1488F_{1,3}F_{2,2}F_{3,1}
+
1632F_{1,3}F_{2,2}G_{4,0}
-
2676F_{1,3}F_{3,1}G_{3,1}
+
\\
&
\ \ \ \ \ 
+
696F_{1,3}G_{3,1}G_{4,0}
+
480F_{2,2}^3
+
1440F_{2,2}^2G_{3,1}
-
450F_{2,2}G_{3,1}^2
-
195G_{3,1}^3
-
1920F_{0,4}F_{2,2}
+
\\
&
\ \ \ \ \ 
+
1440F_{0,4}G_{3,1}
-
576F_{1,3}^2
-
12960F_{1,3}F_{5,0}
+
5040F_{3,1}
+
12960G_{4,0},
\\
\B_5
&
:= 240F_{0,4}^2F_{2,2}G_{3,1}
-
180F_{0,4}^2G_{3,1}^2
-
72F_{0,4}F_{2,2}F_{3,1}
+
432F_{0,4}F_{2,2}G_{4,0}
+
204F_{0,4}F_{3,1}G_{3,1}
+
\\
&
\ \ \ \ \ 
-
464F_{0,4}G_{3,1}G_{4,0}
+
60F_{1,3}F_{2,2}^2
+
525F_{1,3}F_{2,2}G_{3,1}
-
255F_{1,3}G_{3,1}^2
+
864F_{0,4}F_{1,3}
+
\\
&
\ \ \ \ \ 
-
2160F_{0,4}F_{5,0}
-
180F_{3,1}^2
+
240F_{3,1}G_{4,0}
+
240G_{4,0}^2
+
5040F_{2,2}
-
2790G_{3,1},
\\
\B_6 
&
:= 80F_{1,3}^2F_{2,2}G_{3,1}
-
360F_{1,3}F_{3,1}^2
+
1440F_{1,3}F_{3,1}G_{4,0}
-
1440F_{1,3}G_{4,0}^2
+
608F_{2,2}^2F_{3,1}
+
\\
&
\ \ \ \ \ 
-
208F_{2,2}^2G_{4,0}
-
4F_{2,2}F_{3,1}G_{3,1}
+
704F_{2,2}G_{3,1}G_{4,0}
-
159F_{3,1}G_{3,1}^2
-
291G_{3,1}^2G_{4,0}
+
\\
&
\ \ \ \ \ 
+
1440F_{0,4}F_{3,1}
-
2880F_{0,4}G_{4,0}
-
456F_{1,3}F_{2,2}
-
1344F_{1,3}G_{3,1}
-
2760F_{2,2}F_{5,0}
+
\\
&
\ \ \ \ \ 
+
1110F_{5,0}G_{3,1}
+
14400,
\\
\B_7 
&
:= 32F_{0,4}F_{1,3}F_{3,1}G_{3,1}
-
168F_{1,3}F_{3,1}^2
+
480F_{1,3}F_{3,1}G_{4,0}
-
288F_{1,3}G_{4,0}^2
+
96F_{2,2}^2F_{3,1}
+
\\
&
\ \ \ \ \ 
-
144F_{2,2}^2G_{4,0}
+
60F_{2,2}F_{3,1}G_{3,1}
+
192F_{2,2}G_{3,1}G_{4,0}
-
51F_{3,1}G_{3,1}^2
-
63G_{3,1}^2G_{4,0}
+
\\
&
\ \ \ \ \ 
+
288F_{0,4}F_{3,1}
-
576F_{0,4}G_{4,0}
+
24F_{1,3}F_{2,2}
-
288F_{1,3}G_{3,1}
-
360F_{2,2}F_{5,0}
+
270F_{5,0}G_{3,1}
+
2880,
\\
\B_8 
&
:= 2560F_{0,4}F_{1,3}F_{2,2}G_{4,0}
+
6840F_{1,3}F_{3,1}^2
-
15840F_{1,3}F_{3,1}G_{4,0}
+
4320F_{1,3}G_{4,0}^2
+
\\
&
\ \ \ \ \ 
-
288F_{2,2}^2F_{3,1}
+
25968F_{2,2}^2G_{4,0}
-
10836F_{2,2}F_{3,1}G_{3,1}
-
14784F_{2,2}G_{3,1}G_{4,0}
+
5229F_{3,1}G_{3,1}^2
+
\\
&
\ \ \ \ \ 
+
1161G_{3,1}^2G_{4,0}
-
27360F_{0,4}F_{3,1}
-
6720F_{0,4}G_{4,0}
+
216F_{1,3}F_{2,2}
+
6624F_{1,3}G_{3,1}
+
\\
&
\ \ \ \ \ 
-
3240F_{2,2}F_{5,0}
-
6210F_{5,0}G_{3,1}
-
43200,
\endaligned
\]
%%%%%%%%%%%%%%%%%%%%%%%%%%%%%%%%%%%%%%%%%%%%%%%%%%%%%%%%%%%%%%%%%%%%%%
\[
\aligned
e_1
&
\,:=\,
-
(
-
1
+
\tfrac{1}{5}F_{1,3}x
+
2xF_{5,0}
+
\tfrac{1}{5}xF_{3,1}G_{3,1}
-
\tfrac{4}{15}xF_{3,1}F_{2,2}
-
\tfrac{1}{5}xG_{4,0}G_{3,1}
+
\tfrac{4}{15}xG_{4,0}F_{2,2}
+
\\
&
\ \ \ \ \ 
+
\tfrac{3}{2}G_{3,1}u
+
2vG_{4,0}
)\partial_x
-
(
4v
-
\tfrac{1}{2}xG_{3,1}
+
\tfrac{2}{3}xF_{2,2}
+
\tfrac{3}{5}yF_{1,3}
+
yF_{5,0}
-
2uG_{4,0}
+
3uF_{3,1}
+
\\
&
\ \ \ \ \ 
+
\tfrac{1}{10}yF_{3,1}G_{3,1}
-
\tfrac{2}{15}yF_{3,1}F_{2,2}
-
\tfrac{1}{10}yG_{4,0}G_{3,1}
+
\tfrac{2}{15}yG_{4,0}F_{2,2}
)\partial_y
-
(
-
y
+
\tfrac{4}{5}uF_{1,3}
+
3uF_{5,0}
+
\\
&
\ \ \ \ \ 
+
\tfrac{3}{10}uF_{3,1}G_{3,1}
-
\tfrac{2}{5}uF_{3,1}F_{2,2}
-
\tfrac{3}{10}uG_{4,0}G_{3,1}
+
\tfrac{2}{5}uG_{4,0}F_{2,2}
-
\tfrac{1}{2}vG_{3,1}
+
\tfrac{2}{3}vF_{2,2}
)\partial_u
+
\\
&
\ \ \ \ \ 
-
(
-
2x
+
\tfrac{2}{5}vF_{1,3}
+
4vF_{5,0}
+
\tfrac{2}{5}vF_{3,1}G_{3,1}
-
\tfrac{8}{15}vF_{3,1}F_{2,2}
-
\tfrac{2}{5}vG_{4,0}G_{3,1}
+
\tfrac{8}{15}vG_{4,0}F_{2,2}
)\partial_v
, 
\\
e_2
&
\,:=\,
-
(
\tfrac{4}{5}xF_{0,4}
-
\tfrac{9}{25}F_{1,3}xF_{3,1}
+
\tfrac{2}{5}F_{1,3}xG_{4,0}
-
\tfrac{1}{25}xG_{4,0}F_{3,1}G_{3,1}
+
\tfrac{4}{75}xG_{4,0}F_{3,1}F_{2,2}
+
\tfrac{1}{20}xG_{3,1}^2
+
\\
&
\ \ \ \ \ 
+
\tfrac{2}{15}xF_{2,2}^2
-
\tfrac{1}{6}xG_{3,1}F_{2,2}
+
\tfrac{2}{5}xF_{3,1}F_{5,0}
+
\tfrac{1}{25}xF_{3,1}^2G_{3,1}
-
\tfrac{4}{75}xF_{3,1}^2F_{2,2}
+
3y
+
\tfrac{1}{2}vG_{3,1}
)\partial_x
+
\\
&
\ \ \ \ \ 
-
(
F_{1,3}x
-
1
+
\tfrac{12}{5}yF_{0,4}
+
\tfrac{1}{40}yG_{3,1}^2
+
\tfrac{1}{15}yF_{2,2}^2
-
\tfrac{1}{2}G_{3,1}u
+
2uF_{2,2}
+
vF_{3,1}
-
\tfrac{1}{50}yG_{4,0}F_{3,1}G_{3,1}
+
\\
&
\ \ \ \ \ 
+
\tfrac{2}{75}yG_{4,0}F_{3,1}F_{2,2}
-
\tfrac{9}{50}yF_{1,3}F_{3,1}
+
\tfrac{1}{5}yF_{1,3}G_{4,0}
-
\tfrac{1}{12}yG_{3,1}F_{2,2}
+
\tfrac{1}{5}yF_{3,1}F_{5,0}
+
\tfrac{1}{50}yF_{3,1}^2G_{3,1}
+
\\
&
\ \ \ \ \ 
-
\tfrac{2}{75}yF_{3,1}^2F_{2,2}
)\partial_y
-
(
-
x
+
\tfrac{16}{5}uF_{0,4}
-
\tfrac{27}{50}uF_{1,3}F_{3,1}
+
\tfrac{3}{5}uF_{1,3}G_{4,0}
-
\tfrac{3}{50}uG_{4,0}F_{3,1}G_{3,1}
+
\\
&
\ \ \ \ \ 
+
\tfrac{2}{25}uG_{4,0}F_{3,1}F_{2,2}
+
\tfrac{3}{40}uG_{3,1}^2
+
\tfrac{1}{5}uF_{2,2}^2
-
\tfrac{1}{4}G_{3,1}uF_{2,2}
+
\tfrac{3}{5}uF_{3,1}F_{5,0}
+
\tfrac{3}{50}uF_{3,1}^2G_{3,1}
+
\\
&
\ \ \ \ \ 
-
\tfrac{2}{25}uF_{3,1}^2F_{2,2}
+
vF_{1,3}
)\partial_u
-
(
6u
+
\tfrac{8}{5}vF_{0,4}
-
\tfrac{18}{25}vF_{1,3}F_{3,1}
+
\tfrac{4}{5}vF_{1,3}G_{4,0}
-
\tfrac{2}{25}vG_{4,0}F_{3,1}G_{3,1}
+
\\
&
\ \ \ \ \ 
+
\tfrac{8}{75}vG_{4,0}F_{3,1}F_{2,2}
+
\tfrac{1}{10}vG_{3,1}^2
+
\tfrac{4}{15}vF_{2,2}^2
-
\tfrac{1}{3}vG_{3,1}F_{2,2}
+
\tfrac{4}{5}vF_{3,1}F_{5,0}
+
\tfrac{2}{25}vF_{3,1}^2G_{3,1}
-
\tfrac{8}{75}vF_{3,1}^2F_{2,2}
)\partial_v,
\endaligned
\]

%%%%%%%%%%%%%%%%%%%%%%%%%%%%%%%%%%%%%%%%%%%%%%%%%%%%%%%%%%%%%%%%%%%%%%
\[
\footnotesize
\def\arraystretch{1.25}
\begin{array}{c|cc}
{} & e_1 & e_2 
\\
\hline
e_1 & 
0
&
\rotatebox[origin=c]{0}{
	\begin{tabular}{p{8cm}}
$(
-
\tfrac{4}{5}F_{0,4}
+
\tfrac{9}{25}F_{3,1}F_{1,3}
-
\tfrac{2}{5}G_{4,0}F_{1,3}
+
\tfrac{1}{25}G_{4,0}F_{3,1}G_{3,1}
-
\tfrac{4}{75}G_{4,0}F_{3,1}F_{2,2}
-
\tfrac{1}{20}G_{3,1}^2
-
\tfrac{2}{15}F_{2,2}^2
+
\tfrac{1}{6}G_{3,1}F_{2,2}
-
\tfrac{2}{5}F_{3,1}F_{5,0}
-
\tfrac{1}{25}F_{3,1}^2G_{3,1}
+
\tfrac{4}{75}F_{3,1}^2F_{2,2}
)e_1
+
(
\tfrac{1}{10}F_{3,1}G_{3,1}
-
\tfrac{2}{15}F_{3,1}F_{2,2}
-
\tfrac{1}{10}G_{4,0}G_{3,1}
+
\tfrac{2}{15}G_{4,0}F_{2,2}
-
\tfrac{2}{5}F_{1,3}
+
F_{5,0})e_2$
\end{tabular}}
\\
e_2 &
\rotatebox[origin=c]{0}{
	\begin{tabular}{p{8cm}}
$
-
(
-
\tfrac{4}{5}F_{0,4}
+
\tfrac{9}{25}F_{3,1}F_{1,3}
-
\tfrac{2}{5}G_{4,0}F_{1,3}
+
\tfrac{1}{25}G_{4,0}F_{3,1}G_{3,1}
-
\tfrac{4}{75}G_{4,0}F_{3,1}F_{2,2}
-
\tfrac{1}{20}G_{3,1}^2
-
\tfrac{2}{15}F_{2,2}^2
+
\tfrac{1}{6}G_{3,1}F_{2,2}
-
\tfrac{2}{5}F_{3,1}F_{5,0}
-
\tfrac{1}{25}F_{3,1}^2G_{3,1}
+
\tfrac{4}{75}F_{3,1}^2F_{2,2}
)e_1
-
(
\tfrac{1}{10}F_{3,1}G_{3,1}
-
\tfrac{2}{15}F_{3,1}F_{2,2}
-
\tfrac{1}{10}G_{4,0}G_{3,1}
+
\tfrac{2}{15}G_{4,0}F_{2,2}
-
\tfrac{2}{5}F_{1,3}
+
F_{5,0}
)e_2$
\end{tabular}}
&
0
\end{array}
\]

%%%%%%%%%%%%%%%%%%%%%%%%%%%%%%%%%%%%%%%%%%%%%%%%%%%%%%%%%%%%%%%%%%%%%%
%%%%%%%%%%%%%%%%%%%%%%%%%%%%%%%%%%%%%%%%%%%%%%%%%%%%%%%%%%%%%%%%%%%%%%
%%%%%%%%%%%%%%%%%%%%%%%%%%%%%%%%%%%%%%%%%%%%%%%%%%%%%%%%%%%%%%%%%%%%%%

\[
\text{\bf Model 2e3a4ba}
\ \ \ \ \
\left\{
\aligned
u
&
\,=\,
xy
+
y^3
+
F_{0,4}y^4
+
F_{1,3}xy^3
+
F_{2,2}x^2y^2
+
(
2
+
\tfrac{3}{256}G_{3,1}^4
-
\tfrac{7}{64}F_{1,3}G_{3,1}^2
+
\\
&
\ \ \ \ \ 
+
\tfrac{15}{128}G_{3,1}^2G_{5,0}
-
\tfrac{1}{64}G_{3,1}^3F_{2,2}
+
\tfrac{1}{3}G_{3,1}F_{0,4}
)x^3y
+
(
\tfrac{9}{10}F_{1,3}
-
\tfrac{7}{160}F_{0,4}F_{1,3}G_{3,1}
+
\\
&
\ \ \ \ \ 
+
\tfrac{6}{5}F_{0,4}^2
+
\tfrac{3}{64}G_{3,1}F_{0,4}G_{5,0}
+
\tfrac{3}{640}F_{0,4}G_{3,1}^3
-
\tfrac{1}{160}F_{0,4}G_{3,1}^2F_{2,2}
)y^5
+
(
\tfrac{3}{4}G_{3,1}
+
\\
&
\ \ \ \ \ 
+
\tfrac{3}{2}F_{0,4}F_{1,3}
+
\tfrac{5}{4}F_{0,4}G_{5,0}
+
\tfrac{1}{8}F_{0,4}G_{3,1}^2
-
\tfrac{1}{6}G_{3,1}F_{0,4}F_{2,2}
+
F_{2,2}
)xy^4
+
\\
&
\ \ \ \ \ 
+
(
8
+
\tfrac{9}{128}G_{3,1}^4
-
\tfrac{17}{32}F_{1,3}G_{3,1}^2
+
\tfrac{45}{64}G_{3,1}^2G_{5,0}
-
\tfrac{3}{32}G_{3,1}^3F_{2,2}
+
G_{3,1}F_{0,4}
+
\\
&
\ \ \ \ \ 
-
\tfrac{1}{6}F_{2,2}F_{1,3}G_{3,1}
+
\tfrac{1}{2}F_{1,3}^2
+
\tfrac{5}{4}F_{1,3}G_{5,0}
+
\tfrac{4}{3}F_{2,2}F_{0,4}
)x^2y^3
+
(
\tfrac{1}{8}G_{3,1}^2F_{2,2}
+
\\
&
\ \ \ \ \ 
-
\tfrac{1}{2}F_{1,3}G_{3,1}
+
\tfrac{5}{4}F_{2,2}G_{5,0}
-
\tfrac{1}{6}G_{3,1}F_{2,2}^2
+
\tfrac{5}{6}F_{2,2}F_{1,3}
)x^3y^2
+
(
\tfrac{1}{8}G_{3,1}^2
+
\\
&
\ \ \ \ \ 
-
\tfrac{1}{6}
(
2
+
\tfrac{3}{256}G_{3,1}^4
-
\tfrac{7}{64}F_{1,3}G_{3,1}^2
+
\tfrac{15}{128}G_{3,1}^2G_{5,0}
-
\tfrac{1}{64}G_{3,1}^3F_{2,2}
+
\\
&
\ \ \ \ \ 
+
\tfrac{1}{3}G_{3,1}F_{0,4}
)
)
G_{3,1}F_{2,2}
-
\tfrac{5}{12}G_{3,1}F_{2,2}
+
\tfrac{5}{4}
(
2
+
\tfrac{3}{256}G_{3,1}^4
-
\tfrac{7}{64}F_{1,3}G_{3,1}^2
+
\\
&
\ \ \ \ \ 
+
\tfrac{15}{128}G_{3,1}^2G_{5,0}
-
\tfrac{1}{64}G_{3,1}^3F_{2,2}
+
\tfrac{1}{3}G_{3,1}F_{0,4}
)
)
G_{5,0}
+
\tfrac{1}{8}
(
2
+
\tfrac{3}{256}G_{3,1}^4
+
\\
&
\ \ \ \ \ 
-
\tfrac{7}{64}F_{1,3}G_{3,1}^2
+
\tfrac{15}{128}G_{3,1}^2G_{5,0}
-
\tfrac{1}{64}G_{3,1}^3F_{2,2}
+
\tfrac{1}{3}G_{3,1}F_{0,4}
)
)
G_{3,1}^2
+
\tfrac{1}{3}F_{2,2}^2
)x^4y
+
\\
&
\ \ \ \ \ 
+
(
-
\tfrac{1}{10}
(
2
+
\tfrac{3}{256}G_{3,1}^4
-
\tfrac{7}{64}F_{1,3}G_{3,1}^2
+
\tfrac{15}{128}G_{3,1}^2G_{5,0}
-
\tfrac{1}{64}G_{3,1}^3F_{2,2}
+
\\
&
\ \ \ \ \ 
+
\tfrac{1}{3}G_{3,1}F_{0,4}
))
G_{3,1}
+
\tfrac{2}{15}
(
2
+
\tfrac{3}{256}G_{3,1}^4
-
\tfrac{7}{64}F_{1,3}G_{3,1}^2
+
\tfrac{15}{128}G_{3,1}^2G_{5,0}
+
\\
&
\ \ \ \ \ 
-
\tfrac{1}{64}G_{3,1}^3F_{2,2}
+
\tfrac{1}{3}G_{3,1}F_{0,4}
)
)
F_{2,2}
+
\tfrac{1}{10}G_{3,1}
-
\tfrac{2}{15}F_{2,2})x^5
+
\cdots
,
\\
v
&
\,=\,
x^2
-
\tfrac{3}{2}y^4
+
G_{3,1}x^3y
+
x^4
-
\tfrac{18}{5}F_{0,4}y^5
-
3F_{1,3}xy^4
+
(
-
2F_{2,2}
+
3G_{3,1}
)x^2y^3
+
\\
&
\ \ \ \ \ 
+
(
\tfrac{15}{16}G_{3,1}G_{5,0}
+
\tfrac{3}{32}G_{3,1}^3
-
\tfrac{1}{8}G_{3,1}^2F_{2,2}
+
\tfrac{1}{8}F_{1,3}G_{3,1})x^4y
+
G_{5,0}x^5
+
\cdots,
\endaligned\right.
\]
%%%%%%%%%%%%%%%%%%%%%%%%%%%%%%%%%%%%%%%%%%%%%%%%%%%%%%%%%%%%%%%%%%%%%%
\noindent
Gr\"obner basis generators of 
moduli space core algebraic variety in 
$\R^5 \ni F_{0,4},F_{1,3},F_{2,2},G_{3,1},G_{5,0}$:
\[
\aligned
\B_1 
&
:= 40F_{0,4}F_{1,3}F_{2,2}G_{3,1}
-
24F_{2,2}^3
-
98F_{2,2}^2G_{3,1}
+
119F_{2,2}G_{3,1}^2
+
63G_{3,1}^3
-
376F_{1,3}F_{2,2}
+
\\
&
\ \ \ \ \ 
+
1512F_{1,3}G_{3,1}
+
4020F_{2,2}G_{5,0}
+
840G_{3,1}G_{5,0},
\\
\B_2
&
:= 2560F_{0,4}^2F_{2,2}
-
1280F_{0,4}^2G_{3,1}
+
336F_{2,2}^3
-
1468F_{2,2}^2G_{3,1}
+
1514F_{2,2}G_{3,1}^2
-
297G_{3,1}^3
+
\\
&
\ \ \ \ \ 
+
8544F_{1,3}F_{2,2}
-
1668F_{1,3}G_{3,1}
+
12720F_{2,2}G_{5,0}
-
2310G_{3,1}G_{5,0},
\\
\B_3
&
:= 183840F_{2,2}^3G_{3,1}
-
1679176F_{1,3}F_{2,2}G_{3,1}
+
1714878F_{1,3}G_{3,1}^2
-
3867120F_{2,2}^2G_{5,0}
+
\\
&
\ \ \ \ \ 
+
5261580F_{2,2}G_{3,1}G_{5,0}
-
925695G_{3,1}^2G_{5,0}
-
2851328F_{0,4}F_{2,2}
+
1425664F_{0,4}G_{3,1}
+
\\
&
\ \ \ \ \ 
+
2373288F_{1,3}^2
-
6520200F_{1,3}G_{5,0}
+
1467450G_{5,0}^2,
\\
\B_4
&
:= 245120F_{2,2}^4
-
5966248F_{1,3}F_{2,2}G_{3,1}
+
2145774F_{1,3}G_{3,1}^2
+
13260240F_{2,2}^2G_{5,0}
+
\\
&
\ \ \ \ \ 
-
10330260F_{2,2}G_{3,1}G_{5,0}
+
3681465G_{3,1}^2G_{5,0}
+
17645056F_{0,4}F_{2,2}
-
8822528F_{0,4}G_{3,1}
+
\\
&
\ \ \ \ \ 
-
6649176F_{1,3}^2
+
161400F_{1,3}G_{5,0}
+
41153850G_{5,0}^2,
\\
\B_5
&
:= 
-
4802048F_{0,4}F_{2,2}G_{3,1}
+
7243008F_{0,4}G_{3,1}^2
+
10831752F_{1,3}^2F_{2,2}
-
3683292F_{1,3}^2G_{3,1}
+
\\
&
\ \ \ \ \ 
-
97567200F_{1,3}F_{2,2}G_{5,0}
+
15179940F_{1,3}G_{3,1}G_{5,0}
+
67614750F_{2,2}G_{5,0}^2
-
33807375G_{3,1}G_{5,0}^2
+
\\
&
\ \ \ \ \ 
-
58103808F_{0,4}F_{1,3}
+
145259520F_{0,4}G_{5,0}
-
164265600F_{2,2}
+
82132800G_{3,1},
\\
\B_6
&
:= 
-
9300648F_{0,4}F_{1,3}F_{2,2}
-
14461652F_{0,4}F_{1,3}G_{3,1}
-
37669380F_{0,4}F_{2,2}G_{5,0}
+
\\
&
\ \ \ \ \ 
+
18834690F_{0,4}G_{3,1}G_{5,0}
+
7055412F_{1,3}^3
-
11327820F_{1,3}^2G_{5,0}
-
15776775F_{1,3}G_{5,0}^2
+
\\
&
\ \ \ \ \ 
-
35143720F_{2,2}^2
-
38921092F_{2,2}G_{3,1}
+
26206992G_{3,1}^2
+
24473808F_{1,3}
-
61184520G_{5,0},
\endaligned
\]
%%%%%%%%%%%%%%%%%%%%%%%%%%%%%%%%%%%%%%%%%%%%%%%%%%%%%%%%%%%%%%%%%%%%%%
\[
\aligned
e_1
&
\,:=\,
(
1
-
\tfrac{5}{2}xG_{5,0}
-
\tfrac{1}{4}xG_{3,1}^2
+
\tfrac{1}{3}xG_{3,1}F_{2,2}
-
\tfrac{3}{2}G_{3,1}u
-
2v
)\partial_x
-
(
4u
-
\tfrac{1}{2}xG_{3,1}
+
\tfrac{2}{3}xF_{2,2}
+
\tfrac{1}{2}yF_{1,3}
+
\\
&
\ \ \ \ \ 
+
\tfrac{5}{4}yG_{5,0}
+
\tfrac{1}{8}yG_{3,1}^2
+
\tfrac{9}{256}uG_{3,1}^4
-
\tfrac{1}{6}yG_{3,1}F_{2,2}
-
\tfrac{21}{64}uF_{1,3}G_{3,1}^2
+
\tfrac{45}{128}uG_{3,1}^2G_{5,0}
-
\tfrac{3}{64}uG_{3,1}^3F_{2,2}
+
\\
&
\ \ \ \ \ 
+
G_{3,1}uF_{0,4}
)\partial_y
-
(
-
y
+
\tfrac{15}{4}uG_{5,0}
+
\tfrac{3}{8}uG_{3,1}^2
-
\tfrac{1}{2}G_{3,1}uF_{2,2}
+
\tfrac{1}{2}uF_{1,3}
-
\tfrac{1}{2}vG_{3,1}
+
\tfrac{2}{3}vF_{2,2}
)\partial_u
+
\\
&
\ \ \ \ \ 
+
(
-
5vG_{5,0}
+
2x
-
\tfrac{1}{2}vG_{3,1}^2
+
\tfrac{2}{3}vG_{3,1}F_{2,2}
)\partial_v
,
\\
e_2
&
\,:=\,
(
\tfrac{7}{16}xF_{1,3}G_{3,1}
-
\tfrac{15}{32}xG_{3,1}G_{5,0}
-
\tfrac{3}{64}xG_{3,1}^3
+
\tfrac{1}{16}xG_{3,1}^2F_{2,2}
-
3y
-
\tfrac{1}{2}vG_{3,1}
)\partial_x
-
(
2yF_{0,4}
+
xF_{1,3}
+
\\
&
\ \ \ \ \ 
+
2v
-
1
+
\tfrac{3}{128}yG_{3,1}^3
-
\tfrac{1}{2}G_{3,1}u
+
2uF_{2,2}
+
\tfrac{3}{256}vG_{3,1}^4
-
\tfrac{7}{32}yF_{1,3}G_{3,1}
+
\tfrac{15}{64}yG_{3,1}G_{5,0}
+
\\
&
\ \ \ \ \ 
-
\tfrac{1}{32}yG_{3,1}^2F_{2,2}
-
\tfrac{7}{64}vF_{1,3}G_{3,1}^2
+
\tfrac{15}{128}vG_{3,1}^2G_{5,0}
-
\tfrac{1}{64}vG_{3,1}^3F_{2,2}
+
\tfrac{1}{3}vG_{3,1}F_{0,4}
)\partial_y
+
\\
&
\ \ \ \ \ 
-
(
-
x
-
\tfrac{21}{32}uF_{1,3}G_{3,1}
+
\tfrac{45}{64}G_{3,1}uG_{5,0}
+
\tfrac{9}{128}uG_{3,1}^3
-
\tfrac{3}{32}uG_{3,1}^2F_{2,2}
+
2uF_{0,4}
+
vF_{1,3})\partial_u
+
\\
&
\ \ \ \ \ 
+
(
-
6u
+
\tfrac{7}{8}vF_{1,3}G_{3,1}
-
\tfrac{15}{16}vG_{3,1}G_{5,0}
-
\tfrac{3}{32}vG_{3,1}^3
+
\tfrac{1}{8}vG_{3,1}^2F_{2,2}
)\partial_v,
\endaligned
\]

%%%%%%%%%%%%%%%%%%%%%%%%%%%%%%%%%%%%%%%%%%%%%%%%%%%%%%%%%%%%%%%%%%%%%%
\[
\footnotesize
\def\arraystretch{1.25}
\begin{array}{c|cc}
{} & e_1 & e_2 
\\
\hline
e_1 & 
0
&
\rotatebox[origin=c]{0}{
\begin{tabular}{p{7cm}}
$
(
\tfrac{7}{16}F_{1,3}G_{3,1}
-
\tfrac{15}{32}G_{3,1}G_{5,0}
-
\tfrac{3}{64}G_{3,1}^3
+
\tfrac{1}{16}G_{3,1}^2F_{2,2})e_1
+
(
-
\tfrac{1}{6}G_{3,1}F_{2,2}
+
\tfrac{1}{8}G_{3,1}^2
-
\tfrac{1}{2}F_{1,3}
+
\tfrac{5}{4}G_{5,0}
)e_2$
\end{tabular}}
\\
e_2 &
\rotatebox[origin=c]{0}{
	\begin{tabular}{p{7cm}}
$
-
(
\tfrac{7}{16}F_{1,3}G_{3,1}
-
\tfrac{15}{32}G_{3,1}G_{5,0}
-
\tfrac{3}{64}G_{3,1}^3
+
\tfrac{1}{16}G_{3,1}^2F_{2,2}
)e_1
-
(
-
\tfrac{1}{6}G_{3,1}F_{2,2}
+
\tfrac{1}{8}G_{3,1}^2
-
\tfrac{1}{2}F_{1,3}
+
\tfrac{5}{4}G_{5,0}
)e_2$
\end{tabular}}
&
0
\end{array}
\]

%%%%%%%%%%%%%%%%%%%%%%%%%%%%%%%%%%%%%%%%%%%%%%%%%%%%%%%%%%%%%%%%%%%%%%
%%%%%%%%%%%%%%%%%%%%%%%%%%%%%%%%%%%%%%%%%%%%%%%%%%%%%%%%%%%%%%%%%%%%%%
%%%%%%%%%%%%%%%%%%%%%%%%%%%%%%%%%%%%%%%%%%%%%%%%%%%%%%%%%%%%%%%%%%%%%%

\[
\text{\bf Model 2e3a4bb}
\ \ \ \ \
\left\{
\aligned
u
&
\,=\,
xy
+
y^3
+
F_{0,4}y^4
+
F_{1,3}xy^3
+
F_{2,2}x^2y^2
+
(
-
2
+
\tfrac{3}{256}G_{3,1}^4
+
\tfrac{7}{64}G_{3,1}^2F_{1,3}
+
\\
&
\ \ \ \ \ 
+
\tfrac{15}{128}G_{3,1}^2G_{5,0}
-
\tfrac{1}{64}G_{3,1}^3F_{2,2}
+
\tfrac{1}{3}G_{3,1}F_{0,4})x^3y
+
(
\tfrac{9}{10}F_{1,3}
+
\tfrac{6}{5}F_{0,4}^2
+
\\
&
\ \ \ \ \ 
+
\tfrac{7}{160}G_{3,1}F_{0,4}F_{1,3}
+
\tfrac{3}{64}G_{3,1}F_{0,4}G_{5,0}
+
\tfrac{3}{640}F_{0,4}G_{3,1}^3
-
\tfrac{1}{160}F_{0,4}G_{3,1}^2F_{2,2})y^5
+
\\
&
\ \ \ \ \ 
+
(
\tfrac{3}{4}G_{3,1}
+
\tfrac{3}{2}F_{0,4}F_{1,3}
-
\tfrac{5}{4}F_{0,4}G_{5,0}
-
\tfrac{1}{8}F_{0,4}G_{3,1}^2
+
\tfrac{1}{6}G_{3,1}F_{0,4}F_{2,2}
+
F_{2,2}
)xy^4
+
\\
&
\ \ \ \ \ 
+
(
-
8
+
\tfrac{9}{128}G_{3,1}^4
+
\tfrac{17}{32}G_{3,1}^2F_{1,3}
+
\tfrac{45}{64}G_{3,1}^2G_{5,0}
-
\tfrac{3}{32}G_{3,1}^3F_{2,2}
+
G_{3,1}F_{0,4}
+
\\
&
\ \ \ \ \ 
+
\tfrac{1}{6}F_{2,2}G_{3,1}F_{1,3}
+
\tfrac{1}{2}F_{1,3}^2
-
\tfrac{5}{4}F_{1,3}G_{5,0}
+
\tfrac{4}{3}F_{2,2}F_{0,4}
)x^2y^3
+
(
-
\tfrac{1}{8}G_{3,1}^2F_{2,2}
+
\\
&
\ \ \ \ \ 
-
\tfrac{1}{2}G_{3,1}F_{1,3}
-
\tfrac{5}{4}F_{2,2}G_{5,0}
+
\tfrac{1}{6}G_{3,1}F_{2,2}^2
+
\tfrac{5}{6}F_{2,2}F_{1,3}
)x^3y^2
+
(
\tfrac{1}{8}G_{3,1}^2
+
\\
&
\ \ \ \ \ 
+
\tfrac{1}{6}
(
-
2
+
\tfrac{3}{256}G_{3,1}^4
+
\tfrac{7}{64}G_{3,1}^2F_{1,3}
+
\tfrac{15}{128}G_{3,1}^2G_{5,0}
-
\tfrac{1}{64}G_{3,1}^3F_{2,2}
+
\\
&
\ \ \ \ \ 
+
\tfrac{1}{3}G_{3,1}F_{0,4}
)
)
G_{3,1}F_{2,2}
-
\tfrac{5}{12}F_{2,2}G_{3,1}
-
(
\tfrac{5}{4}
(
-
2
+
\tfrac{3}{256}G_{3,1}^4
+
\tfrac{7}{64}G_{3,1}^2F_{1,3}
+
\\
&
\ \ \ \ \ 
+
\tfrac{15}{128}G_{3,1}^2G_{5,0}
-
\tfrac{1}{64}G_{3,1}^3F_{2,2}
+
\tfrac{1}{3}G_{3,1}F_{0,4}
)
)
G_{5,0}
-
(
\tfrac{1}{8}
(
-
2
+
\tfrac{3}{256}G_{3,1}^4
+
\\
&
\ \ \ \ \ 
+
\tfrac{7}{64}G_{3,1}^2F_{1,3}
+
\tfrac{15}{128}G_{3,1}^2G_{5,0}-
\tfrac{1}{64}G_{3,1}^3F_{2,2}
+
\tfrac{1}{3}G_{3,1}F_{0,4}
)
)
G_{3,1}^2
+
\tfrac{1}{3}F_{2,2}^2
)x^4y
+
\\
&
\ \ \ \ \ 
+
(
-
(
\tfrac{1}{10}
(
-
2
+
\tfrac{3}{256}G_{3,1}^4
+
\tfrac{7}{64}G_{3,1}^2F_{1,3}
+
\tfrac{15}{128}G_{3,1}^2G_{5,0}
-
\tfrac{1}{64}G_{3,1}^3F_{2,2}
+
\\
&
\ \ \ \ \ 
+
\tfrac{1}{3}G_{3,1}F_{0,4}
)
)
G_{3,1}
+
(
\tfrac{2}{15}
(
-
2
+
\tfrac{3}{256}G_{3,1}^4
+
\tfrac{7}{64}G_{3,1}^2F_{1,3}
+
\tfrac{15}{128}G_{3,1}^2G_{5,0}
+
\\
&
\ \ \ \ \ 
-
\tfrac{1}{64}G_{3,1}^3F_{2,2}
+
\tfrac{1}{3}G_{3,1}F_{0,4}
)
)
F_{2,2}
-
\tfrac{1}{10}G_{3,1}
+
\tfrac{2}{15}F_{2,2}
)x^5
+\cdots
,
\\
v
&
\,=\,
x^2
-
\tfrac{3}{2}y^4
+
G_{3,1}x^3y
-
x^4
-
\tfrac{18}{5}F_{0,4}y^5
-
3F_{1,3}xy^4
+
(
-
2F_{2,2}
+
3G_{3,1}
)x^2y^3
+
\\
&
\ \ \ \ \ 
+
(
-
\tfrac{15}{16}G_{3,1}G_{5,0}
-
\tfrac{3}{32}G_{3,1}^3
+
\tfrac{1}{8}G_{3,1}^2F_{2,2}
+
\tfrac{1}{8}G_{3,1}F_{1,3}
)x^4y
+
G_{5,0}x^5
+\cdots,
\endaligned\right.
\]
%%%%%%%%%%%%%%%%%%%%%%%%%%%%%%%%%%%%%%%%%%%%%%%%%%%%%%%%%%%%%%%%%%%%%%
\noindent
Gr\"obner basis generators of 
moduli space core algebraic variety in 
$\R^5 \ni F_{0,4},F_{1,3},F_{2,2},G_{3,1},G_{5,0}$:
\[
\aligned
\B_1 
&
:= 1701000F_{1,3}F_{2,2}G_{5,0}^2+1354976F_{2,2}^3+8802232F_{2,2}^2G_{3,1}-1660396F_{2,2}G_{3,1}^2-5460477G_{3,1}^3
+
\\
&
\ \ \ \ \ -21058704F_{1,3}F_{2,2}+69200028F_{1,3}G_{3,1}-192961080F_{2,2}G_{5,0}-21141810G_{3,1}G_{5,0},
\\
\B_2 
&
:= 1575F_{1,3}^2G_{3,1}G_{5,0}-56F_{2,2}^3+1118F_{2,2}^2G_{3,1}+1501F_{2,2}G_{3,1}^2-1398G_{3,1}^3-1056F_{1,3}F_{2,2}
+
\\
&
\ \ \ \ \ +6732F_{1,3}G_{3,1}-28620F_{2,2}G_{5,0}+3060G_{3,1}G_{5,0},
\\
\B_3 
&
:= 1915G_{3,1}^4+14976F_{1,3}F_{2,2}G_{3,1}-26928F_{1,3}G_{3,1}^2-13920F_{2,2}^2G_{5,0}+62880F_{2,2}G_{3,1}G_{5,0}
+
\\
&
\ \ \ \ \ +29880G_{3,1}^2G_{5,0}+33792F_{0,4}F_{2,2}-16896F_{0,4}G_{3,1}-23472F_{1,3}^2-54000F_{1,3}G_{5,0}+11700G_{5,0}^2,
\\
\B_4 
&
:= 45960F_{0,4}F_{1,3}G_{3,1}G_{5,0}+27944F_{1,3}F_{2,2}G_{3,1}-9618F_{1,3}G_{3,1}^2-75960F_{2,2}^2G_{5,0}
+
\\
&
\ \ \ \ \ -108150F_{2,2}G_{3,1}G_{5,0}+61425G_{3,1}^2G_{5,0}+32256F_{0,4}F_{2,2}-16128F_{0,4}G_{3,1}+30240F_{1,3}^2
+
\\
&
\ \ \ \ \ 
-64080F_{1,3}G_{5,0}-349200G_{5,0}^2,
\\
\B_5 
&
:= -4802048F_{0,4}F_{2,2}G_{3,1}+7243008F_{0,4}G_{3,1}^2+10831752F_{1,3}^2F_{2,2}-3683292F_{1,3}^2G_{3,1}
+
\\
&
\ \ \ \ \ 
+97567200F_{1,3}F_{2,2}G_{5,0}-15179940F_{1,3}G_{3,1}G_{5,0}+67614750F_{2,2}G_{5,0}^2-33807375G_{3,1}G_{5,0}^2
+
\\
&
\ \ \ \ \ 
+58103808F_{0,4}F_{1,3}+145259520F_{0,4}G_{5,0}+164265600F_{2,2}-82132800G_{3,1},
\\
\B_6 
&
:= -9300648F_{0,4}F_{1,3}F_{2,2}-14461652F_{0,4}F_{1,3}G_{3,1}+37669380F_{0,4}F_{2,2}G_{5,0}
+
\\
&
\ \ \ \ \ 
-18834690F_{0,4}G_{3,1}G_{5,0}+7055412F_{1,3}^3+11327820F_{1,3}^2G_{5,0}-15776775F_{1,3}G_{5,0}^2
+
\\
&
\ \ \ \ \ 
-35143720F_{2,2}^2-38921092F_{2,2}G_{3,1}+26206992G_{3,1}^2-24473808F_{1,3}-61184520G_{5,0},
\endaligned
\]
%%%%%%%%%%%%%%%%%%%%%%%%%%%%%%%%%%%%%%%%%%%%%%%%%%%%%%%%%%%%%%%%%%%%%%
\[
\aligned
e_1
&
\,:=\,
-
(
-
1
-
\tfrac{5}{2}xG_{5,0}
-
\tfrac{1}{4}xG_{3,1}^2
+
\tfrac{1}{3}xG_{3,1}F_{2,2}
+
\tfrac{3}{2}G_{3,1}u
-
2v
)\partial_x
-
(
-
4u
-
\tfrac{1}{2}xG_{3,1}
+
\tfrac{2}{3}xF_{2,2}
+
\\
&
\ \ \ \ \ 
+
\tfrac{1}{2}yF_{1,3}
-
\tfrac{5}{4}yG_{5,0}
-
\tfrac{1}{8}yG_{3,1}^2
+
\tfrac{9}{256}uG_{3,1}^4
+
\tfrac{1}{6}yG_{3,1}F_{2,2} 
+
\tfrac{21}{64}uF_{1,3}G_{3,1}^2
+
\tfrac{45}{128}uG_{3,1}^2G_{5,0}
+
\\
&
\ \ \ \ \ 
-
\tfrac{3}{64}uG_{3,1}^3F_{2,2}
+
uG_{3,1}F_{0,4}
)\partial_y
-
(
-
y
-
\tfrac{15}{4}uG_{5,0}
-
\tfrac{3}{8}uG_{3,1}^2
+
\tfrac{1}{2}uG_{3,1}F_{2,2}
+
\tfrac{1}{2}uF_{1,3}
-
\tfrac{1}{2}vG_{3,1}
+
\\
&
\ \ \ \ \ 
+
\tfrac{2}{3}vF_{2,2}
)\partial_u
-
(
-
2x
-
5vG_{5,0}
-
\tfrac{1}{2}vG_{3,1}^2
+
\tfrac{2}{3}vG_{3,1}F_{2,2}
)\partial_v
, 
\\
e_2
&
\,:=\,
-
(
\tfrac{7}{16}xG_{3,1}F_{1,3}
+
\tfrac{15}{32}xG_{3,1}G_{5,0}
+
\tfrac{3}{64}xG_{3,1}^3
-
\tfrac{1}{16}xG_{3,1}^2F_{2,2}
+
3y
+
\tfrac{1}{2}vG_{3,1}
)\partial_x
-
(
-
1
-
2v
+
\\
&
\ \ \ \ \ 
+
xF_{1,3}
+
2yF_{0,4}
+
\tfrac{3}{128}yG_{3,1}^3
-
\tfrac{1}{2}G_{3,1}u
+
2uF_{2,2}
+
\tfrac{3}{256}vG_{3,1}^4
+
\tfrac{7}{32}yG_{3,1}F_{1,3}
+
\tfrac{15}{64}yG_{3,1}G_{5,0}
+
\\
&
\ \ \ \ \ 
-
\tfrac{1}{32}yG_{3,1}^2F_{2,2}
+
\tfrac{7}{64}vG_{3,1}^2F_{1,3}
+
\tfrac{15}{128}vG_{3,1}^2G_{5,0}
-
\tfrac{1}{64}vG_{3,1}^3F_{2,2}
+
\tfrac{1}{3}vG_{3,1}F_{0,4}
)\partial_y
+
\\
&
\ \ \ \ \ 
-
(
-
x
+
\tfrac{21}{32}uG_{3,1}F_{1,3}
+
\tfrac{45}{64}uG_{3,1}G_{5,0}
+
\tfrac{9}{128}uG_{3,1}^3
-
\tfrac{3}{32}uG_{3,1}^2F_{2,2}
+
2uF_{0,4}
+
vF_{1,3}
)\partial_u
+
\\
&
\ \ \ \ \ 
-
(
6u
+
\tfrac{7}{8}vG_{3,1}F_{1,3}
+
\tfrac{15}{16}vG_{3,1}G_{5,0}
+
\tfrac{3}{32}vG_{3,1}^3
-
\tfrac{1}{8}vG_{3,1}^2F_{2,2}
)\partial_v,
\endaligned
\]
%%%%%%%%%%%%%%%%%%%%%%%%%%%%%%%%%%%%%%%%%%%%%%%%%%%%%%%%%%%%%%%%%%%%%%
\[
\footnotesize
\def\arraystretch{1.25}
\begin{array}{c|cc}
{} & e_1 & e_2 
\\
\hline
e_1 & 
0
&
\rotatebox[origin=c]{0}{
	\begin{tabular}{p{7cm}}
$
(
-
\tfrac{7}{16}G_{3,1}F_{1,3}
-
\tfrac{15}{32}G_{3,1}G_{5,0}
-
\tfrac{3}{64}G_{3,1}^3
+
\tfrac{1}{16}G_{3,1}^2F_{2,2})e_1
+
(
\tfrac{1}{6}F_{2,2}G_{3,1}
-
\tfrac{1}{8}G_{3,1}^2
-
\tfrac{5}{4}G_{5,0}
-
\tfrac{1}{2}F_{1,3})e_2
$
\end{tabular}}
\\
e_2 &
\rotatebox[origin=c]{0}{
	\begin{tabular}{p{7cm}}
$
-
(
-
\tfrac{7}{16}G_{3,1}F_{1,3}
-
\tfrac{15}{32}G_{3,1}G_{5,0}
-
\tfrac{3}{64}G_{3,1}^3
+
\tfrac{1}{16}G_{3,1}^2F_{2,2}
)e_1
-
(
\tfrac{1}{6}F_{2,2}G_{3,1}
-
\tfrac{1}{8}G_{3,1}^2
-
\tfrac{5}{4}G_{5,0}
-
\tfrac{1}{2}F_{1,3})e_2
$
\end{tabular}}
&
0
\end{array}
\]

%%%%%%%%%%%%%%%%%%%%%%%%%%%%%%%%%%%%%%%%%%%%%%%%%%%%%%%%%%%%%%%%%%%%%%
%%%%%%%%%%%%%%%%%%%%%%%%%%%%%%%%%%%%%%%%%%%%%%%%%%%%%%%%%%%%%%%%%%%%%%
%%%%%%%%%%%%%%%%%%%%%%%%%%%%%%%%%%%%%%%%%%%%%%%%%%%%%%%%%%%%%%%%%%%%%%

\[
\text{\bf Model 2e3a4ca}
\ \ \ \ \
\left\{
\aligned
u
&
\,=\,
xy
+
y^3
+
F_{0,4}y^4
+
F_{2,2}x^2y^2+x^3y
+\cdots,
\\
v
&
\,=\,
x^2
-
\tfrac{3}{2}y^4
+
2F_{2,2}x^3y
+
\cdots,
\endaligned\right.
\]
%%%%%%%%%%%%%%%%%%%%%%%%%%%%%%%%%%%%%%%%%%%%%%%%%%%%%%%%%%%%%%%%%%%%%%
\noindent
Gr\"obner basis generator of 
moduli space core algebraic variety in 
$\R^2 \ni F_{2,2},F_{0,4}$:
\[
\B_1 := 2F_{0,4}F_{2,2}-3,
\]
%%%%%%%%%%%%%%%%%%%%%%%%%%%%%%%%%%%%%%%%%%%%%%%%%%%%%%%%%%%%%%%%%%%%%%
\[
\aligned
e_1
&
\,:=\,
-
(
3F_{2,2}u
-
1
)\partial_x
+
(
\tfrac{1}{3}F_{2,2}x
-
3u
)\partial_y
+
(
y
+
\tfrac{1}{3}F_{2,2}v
)\partial_u
+
2x\partial_v
,
\\
e_2
&
\,:=\,
-
(
F_{2,2}v
+
3y
)\partial_x
-
(
2F_{0,4}y
+
F_{2,2}u
+
v
-
1
)\partial_y
-
(
2F_{0,4}u
-
x
)\partial_u-6u\partial_v,
\endaligned
\]

%%%%%%%%%%%%%%%%%%%%%%%%%%%%%%%%%%%%%%%%%%%%%%%%%%%%%%%%%%%%%%%%%%%%%%
\[
\footnotesize
\def\arraystretch{1.25}
\begin{array}{c|cc}
{} & e_1 & e_2 
\\
\hline
e_1 & 0 & 0
\\
e_2 & 0 & 0
\end{array}
\]

%%%%%%%%%%%%%%%%%%%%%%%%%%%%%%%%%%%%%%%%%%%%%%%%%%%%%%%%%%%%%%%%%%%%%%
%%%%%%%%%%%%%%%%%%%%%%%%%%%%%%%%%%%%%%%%%%%%%%%%%%%%%%%%%%%%%%%%%%%%%%
%%%%%%%%%%%%%%%%%%%%%%%%%%%%%%%%%%%%%%%%%%%%%%%%%%%%%%%%%%%%%%%%%%%%%%

\[
\text{\bf Model 2e3a4cb}
\ \ \ \ \
\left\{
\aligned
u
&
\,=\,
xy
+
y^3
+
F_{0,4}y^4
+
F_{2,2}x^2y^2
-
x^3y
+
\cdots,
\\
v
&
\,=\,
x^2
-
\tfrac{3}{2}y^4
+
2F_{2,2}x^3y
+
\cdots,
\endaligned\right.
\]
%%%%%%%%%%%%%%%%%%%%%%%%%%%%%%%%%%%%%%%%%%%%%%%%%%%%%%%%%%%%%%%%%%%%%%
\noindent
Gr\"obner basis generator of 
moduli space core algebraic variety in 
$\R^2 \ni F_{2,2},F_{0,4}$:
\[
\B_1= 2F_{0,4}F_{2,2}+3,
\]
%%%%%%%%%%%%%%%%%%%%%%%%%%%%%%%%%%%%%%%%%%%%%%%%%%%%%%%%%%%%%%%%%%%%%%
\[
\aligned
e_1
&
\,:=\,
-
(
3F_{2,2}u
-
1
)\partial_x
+
(
\tfrac{1}{3}F_{2,2}x
+
3u
)\partial_y
+
(
y
+
\tfrac{1}{3}F_{2,2}v
)\partial_u
+
2x\partial_v,
\\
e_2
&
\,:=\,
-
(
F_{2,2}v
+
3y
)\partial_x
-
(
2F_{0,4}y
+
F_{2,2}u
-
v
-
1
)\partial_y
-
(
2F_{0,4}u
-
x
)\partial_u
-6u\partial_v,
\endaligned
\]

%%%%%%%%%%%%%%%%%%%%%%%%%%%%%%%%%%%%%%%%%%%%%%%%%%%%%%%%%%%%%%%%%%%%%%
\[
\footnotesize
\def\arraystretch{1.25}
\begin{array}{c|cc}
{} & e_1 & e_2 
\\
\hline
e_1 & 0 & 0
\\
e_2 & 0 & 0
\end{array}
\]

%%%%%%%%%%%%%%%%%%%%%%%%%%%%%%%%%%%%%%%%%%%%%%%%%%%%%%%%%%%%%%%%%%%%%%
%%%%%%%%%%%%%%%%%%%%%%%%%%%%%%%%%%%%%%%%%%%%%%%%%%%%%%%%%%%%%%%%%%%%%%
%%%%%%%%%%%%%%%%%%%%%%%%%%%%%%%%%%%%%%%%%%%%%%%%%%%%%%%%%%%%%%%%%%%%%%

\[
\text{\bf Model 2e3a4d}
\ \ \ \ \
\left\{
\aligned
u
&
\,=\,
xy
+
y^3
+
\tfrac{1}{2}x^2y^2
+
\tfrac{5}{4}xy^4
+
\tfrac{3}{4}y^6
+
\tfrac{7}{12}x^3y^3
-
\tfrac{1}{1080}x^6
+
\tfrac{9}{4}x^2y^5
+
\\
&
\ \ \ \ \
-
\tfrac{1}{60}x^5y^2
+
\cdots,
\\
v
&
\,=\,
x^2
-
\tfrac{3}{2}y^4
+
x^3y
+
2x^2y^3
-
\tfrac{1}{30}x^5
+
\tfrac{5}{4}x^4y^2
-
\tfrac{15}{14}y^7
+
\tfrac{19}{4}x^3y^4
+
\\
&
\ \ \ \ \
-
\tfrac{1}{9}x^6y
+
\cdots,
\endaligned\right.
\]
%%%%%%%%%%%%%%%%%%%%%%%%%%%%%%%%%%%%%%%%%%%%%%%%%%%%%%%%%%%%%%%%%%%%%%
\[
\aligned
e_1
&
\,:=\,
-
(
-
1
+
\tfrac{3}{2}u
)\partial_x
+
\tfrac{1}{6}x\partial_y
+
(
y
+
\tfrac{1}{6}v
)\partial_u
+
2x\partial_v,
\\
e_2
&
\,:=\,
-
(
3y
+
\tfrac{1}{2}v
)\partial_x
-
(
-
1
+
\tfrac{1}{2}u
)\partial_y
+
x\partial_u
-
6u\partial_v,
\endaligned
\]
%%%%%%%%%%%%%%%%%%%%%%%%%%%%%%%%%%%%%%%%%%%%%%%%%%%%%%%%%%%%%%%%%%%%%%
\[
\footnotesize
\def\arraystretch{1.25}
\begin{array}{c|cc}
{} & e_1 & e_2 
\\
\hline
e_1 & 0 & 0
\\
e_2 & 0 & 0
\end{array}
\]

%%%%%%%%%%%%%%%%%%%%%%%%%%%%%%%%%%%%%%%%%%%%%%%%%%%%%%%%%%%%%%%%%%%%%%
%%%%%%%%%%%%%%%%%%%%%%%%%%%%%%%%%%%%%%%%%%%%%%%%%%%%%%%%%%%%%%%%%%%%%%
%%%%%%%%%%%%%%%%%%%%%%%%%%%%%%%%%%%%%%%%%%%%%%%%%%%%%%%%%%%%%%%%%%%%%%

\[
\text{\bf Model 2e3a4fa}
\ \ \ \ \
\left\{
\aligned
u
&
\,=\,
xy
+
y^3
+
F_{0,4}y^4
+
xy^3
+
(
\tfrac{9}{10}
+
\tfrac{6}{5}F_{0,4}^2
)y^5
+
2F_{0,4}xy^4
+
x^2y^3
+\cdots,
\\
v
&
\,=\,
x^2
-
\tfrac{3}{2}y^4
-
\tfrac{18}{5}F_{0,4}y^5
-
3xy^4
+
\cdots,
\endaligned\right.
\]
%%%%%%%%%%%%%%%%%%%%%%%%%%%%%%%%%%%%%%%%%%%%%%%%%%%%%%%%%%%%%%%%%%%%%%
for any value of $F_{0,4}$,
%%%%%%%%%%%%%%%%%%%%%%%%%%%%%%%%%%%%%%%%%%%%%%%%%%%%%%%%%%%%%%%%%%%%%%
\[
\aligned
e_1
&
\,:=\,
-
(
x
-
1
)\partial_x
-
y\partial_y
-
(
2u
-
y
)\partial_u
-
(
2v
-
2x
)\partial_v
,
\\
e_2
&
\,:=\,
-
3y\partial_x
-
(
2F_{0,4}y
+
x
-
1
)\partial_y
-
(2F_{0,4}u+v-x)\partial_u-6u\partial_v,
\endaligned
\]
%%%%%%%%%%%%%%%%%%%%%%%%%%%%%%%%%%%%%%%%%%%%%%%%%%%%%%%%%%%%%%%%%%%%%%
\[
\footnotesize
\def\arraystretch{1.25}
\begin{array}{c|cc}
{} & e_1 & e_2 
\\
\hline
e_1 & 0 & 0	
\\
e_2 & 0 & 0 
\end{array}
\]

%%%%%%%%%%%%%%%%%%%%%%%%%%%%%%%%%%%%%%%%%%%%%%%%%%%%%%%%%%%%%%%%%%%%%%
%%%%%%%%%%%%%%%%%%%%%%%%%%%%%%%%%%%%%%%%%%%%%%%%%%%%%%%%%%%%%%%%%%%%%%
%%%%%%%%%%%%%%%%%%%%%%%%%%%%%%%%%%%%%%%%%%%%%%%%%%%%%%%%%%%%%%%%%%%%%%

\[
\text{\bf Model 2e3a4fb}
\ \ \ \ \
\left\{
\aligned
u
&
\,=\,
xy
+
y^3
+
F_{0,4}y^4
-
xy^3
+
(
\tfrac{-9}{10}
+
\tfrac{6}{5}F_{0,4}^2
)y^5
-
2F_{0,4}xy^4
+
x^2y^3
+
\\
&
\ \ \ \ \ 
+
(
-
\tfrac{14}{5}F_{0,4}
+
\tfrac{8}{5}F_{0,4}^3
)y^6
+
(
\tfrac{27}{10}
-
\tfrac{18}{5}F_{0,4}^2
)xy^5
+
3F_{0,4}x^2y^4
-
x^3y^3
+
\\
&
\ \ \ \ \ 
+
(
\tfrac{3}{2}
-
\tfrac{46}{7}F_{0,4}^2
+
\tfrac{16}{7}F_{0,4}^4
)y^7
+
(
\tfrac{56}{5}F_{0,4}
-
\tfrac{32}{5}F_{0,4}^3
)xy^6
+
\\
&
\ \ \ \ \ 
+
(
\tfrac{-27}{5}
+
\tfrac{36}{5}F_{0,4}^2
)x^2y^5
-
4F_{0,4}x^3y^4
+
x^4y^3
+
\cdots
,
\\
v
&
\,=\,
x^2
-
\tfrac{3}{2}y^4
-
\tfrac{18}{5}F_{0,4}y^5
+
3xy^4
+
(
\tfrac{12}{5}
-
\tfrac{36}{5}F_{0,4}^2
)y^6
+
\tfrac{54}{5}F_{0,4}xy^5
+
\\
&
\ \ \ \ \ 
-
\tfrac{9}{2}x^2y^4
+
(
\tfrac{-96}{7}F_{0,4}^3
+
\tfrac{78}{7}F_{0,4}
)y^7
+
(
\tfrac{-48}{5}
+
\tfrac{144}{5}F_{0,4}^2
)xy^6
+
\\
&
\ \ \ \ \ 
-
\tfrac{108}{5}F_{0,4}x^2y^5
+
6x^3y^4
+
\cdots,
\endaligned\right.
\]
for any value of $F_{0,4}$,
%%%%%%%%%%%%%%%%%%%%%%%%%%%%%%%%%%%%%%%%%%%%%%%%%%%%%%%%%%%%%%%%%%%%%%
\[
\aligned
e_1
&
\,:=\,
(x+1)\partial_x+y\partial_y+(2u+y)\partial_u+(2v+2x)\partial_v
,\\
e_2
&
\,:=\,
-3y\partial_x-(2F_{0,4}y-x-1)\partial_y-(2F_{0,4}u-v-x)\partial_u-6u\partial_v,
\endaligned
\]
%%%%%%%%%%%%%%%%%%%%%%%%%%%%%%%%%%%%%%%%%%%%%%%%%%%%%%%%%%%%%%%%%%%%%%
\[
\footnotesize
\def\arraystretch{1.25}
\begin{array}{c|cc}
{} & e_1 & e_2 
\\
\hline
e_1 & 0 & 0
\\
e_2 & 0 & 0
\end{array}
\]

%%%%%%%%%%%%%%%%%%%%%%%%%%%%%%%%%%%%%%%%%%%%%%%%%%%%%%%%%%%%%%%%%%%%%%
%%%%%%%%%%%%%%%%%%%%%%%%%%%%%%%%%%%%%%%%%%%%%%%%%%%%%%%%%%%%%%%%%%%%%%
%%%%%%%%%%%%%%%%%%%%%%%%%%%%%%%%%%%%%%%%%%%%%%%%%%%%%%%%%%%%%%%%%%%%%%

\[
\text{\bf Model 2e3a4g}
\ \ \ \ \
\left\{
\aligned
u
&
\,=\,
xy
+
y^3
+
y^4
+
F_{0,5}y^5
+
(
\tfrac{5}{3}F_{0,5}^2
-
\tfrac{2}{3}F_{0,5}
)y^6
+
(
\tfrac{25}{7}F_{0,5}^3
-
\tfrac{10}{3}F_{0,5}^2
+
\\
&
\ \ \ \ \
+
\tfrac{16}{21}F_{0,5}
)y^7
+
\cdots
,
\\
v
&
\,=\,
x^2
-
\tfrac{3}{2}y^4
-
\tfrac{18}{5}y^5
+
(
-
4F_{0,5}
-
\tfrac{12}{5}y^6+
(
-
\tfrac{50}{7}F_{0,5}^2
-
\tfrac{20}{7}F_{0,5}
)y^7
+
\cdots,
\endaligned\right.
\]
for any value of $F_{0,5}$,
%%%%%%%%%%%%%%%%%%%%%%%%%%%%%%%%%%%%%%%%%%%%%%%%%%%%%%%%%%%%%%%%%%%%%%
\[
\aligned
e_1
&
\,:=\,
y\partial_u+2x\partial_v+\partial_x
,
\\
e_2
&
\,:=\,
-(10F_{0,5}x-12x+3y)\partial_x-(5F_{0,5}y-4y-1)\partial_y-(15F_{0,5}u-16u-x)\partial_u
+
\\
&
\ \ \ \ \ 
-(20F_{0,5}v+6u-24v)\partial_v
,
\endaligned
\]
%%%%%%%%%%%%%%%%%%%%%%%%%%%%%%%%%%%%%%%%%%%%%%%%%%%%%%%%%%%%%%%%%%%%%%
\[
\footnotesize
\def\arraystretch{1.25}
\begin{array}{c|cc}
{} & e_1 & e_2 
\\
\hline
e_1 & 0 & -(10F_{0,5}-12)e_1
\\
e_2 & (10F_{0,5}-12)e_1
\end{array}
\]

%%%%%%%%%%%%%%%%%%%%%%%%%%%%%%%%%%%%%%%%%%%%%%%%%%%%%%%%%%%%%%%%%%%%%%
%%%%%%%%%%%%%%%%%%%%%%%%%%%%%%%%%%%%%%%%%%%%%%%%%%%%%%%%%%%%%%%%%%%%%%
%%%%%%%%%%%%%%%%%%%%%%%%%%%%%%%%%%%%%%%%%%%%%%%%%%%%%%%%%%%%%%%%%%%%%%

\[
\text{\bf Model 2e3a4h}
\ \ \ \ \
\left\{
\aligned
u
&
\,=\,
y^3+xy
,
\\
v
&
\,=\,
x^2-\tfrac{3}{2}y^4
.
\endaligned\right.
\]
%%%%%%%%%%%%%%%%%%%%%%%%%%%%%%%%%%%%%%%%%%%%%%%%%%%%%%%%%%%%%%%%%%%%%%
\[
\def\arraystretch{1.25}
\begin{array}{lll}
e_1
\,:=\,
y\partial_u+2x\partial_v+\partial_x
, &
e_2
\,:=\,
x\partial_u-6u\partial_v-3y\partial_x+\partial_y
,
& 
e_3
\,:=\,
3u\partial_u+4v\partial_v+2x\partial_x+y\partial_y
.
\end{array}
\]

%%%%%%%%%%%%%%%%%%%%%%%%%%%%%%%%%%%%%%%%%%%%%%%%%%%%%%%%%%%%%%%%%%%%%%
\[
\footnotesize
\def\arraystretch{1.25}
\begin{array}{c|ccc}
{} & e_1 & e_2 & e_3
\\
\hline
e_1 & 0 & 0 & 2e_1
\\
e_2 & 0 & 0 & e_2
\\
e_3 & -2e_1 & -e_2 & 0
\end{array}
\]

%%%%%%%%%%%%%%%%%%%%%%%%%%%%%%%%%%%%%%%%%%%%%%%%%%%%%%%%%%%%%%%%%%%%%%
%%%%%%%%%%%%%%%%%%%%%%%%%%%%%%%%%%%%%%%%%%%%%%%%%%%%%%%%%%%%%%%%%%%%%%
%%%%%%%%%%%%%%%%%%%%%%%%%%%%%%%%%%%%%%%%%%%%%%%%%%%%%%%%%%%%%%%%%%%%%%

\[
\text{\bf Model 2e3b4a}
\ \ \ \ \
\left\{
\aligned
u
&
\,=\,
xy
+
xy^2
+
xy^3
+
x^2y^2
-
x^3y
+
\tfrac{3}{8}x^4
+
xy^4
+
3x^2y^3
-
3x^3y^2
+
\\
&
\ \ \ \ \
+
\tfrac{7}{8}x^4y
+
xy^5
+
6x^2y^4
-
4x^3y^3
-
\tfrac{21}{8}x^4y^2
+
\tfrac{9}{4}x^5y
-
\tfrac{3}{8}x^6
+
xy^6
+
\\
&
\ \ \ \ \
+
10x^2y^5
-
\tfrac{145}{8}x^4y^3
+
10x^5y^2
-
\tfrac{3}{2}x^6y
+
\cdots,
\\
v
&
\,=\,
x^2
-
x^2y^2
+
\tfrac{1}{4}x^4
-
2x^2y^3
-
2x^3y^2
+
\tfrac{3}{2}x^4y
-
\tfrac{1}{4}x^5
-
3x^2y^4
+
\\
&
\ \ \ \ \
-
8x^3y^3
+
\tfrac{7}{2}x^4y^2
-
\tfrac{1}{8}x^6
-
4x^2y^5
-
20x^3y^4
+
\tfrac{29}{4}x^5y^2
+
\\
&
\ \ \ \ \
-
3x^6y
+
\tfrac{3}{8}x^7
+\cdots
.
\endaligned\right.
\]
%%%%%%%%%%%%%%%%%%%%%%%%%%%%%%%%%%%%%%%%%%%%%%%%%%%%%%%%%%%%%%%%%%%%%%
\[
\aligned
e_1
&
\,:=\,
-
(
-
x
-
1
+
y
+
\tfrac{1}{2}v
)\partial_x
+
(
\tfrac{3}{2}x
-
2y
+
\tfrac{1}{2}u
-
\tfrac{3}{2}v
)\partial_y
-
(
-
y
+
u
-
\tfrac{3}{2}v
)\partial_u
+
\\
&
\ \ \ \ \
-
(
2u
-
2v
-
2x)\partial_v 
,
\\
e_2
&
\,:=\,
(
-
3x
+
u
)\partial_x
+
(
1
-
2x
-
2y
+
2u
+
v
)\partial_y
-
(
-
x
+
3u
+
2v
)\partial_u-6v\partial_v
.
\endaligned
\]
%%%%%%%%%%%%%%%%%%%%%%%%%%%%%%%%%%%%%%%%%%%%%%%%%%%%%%%%%%%%%%%%%%%%%%
\[
\footnotesize
\def\arraystretch{1.25}
\begin{array}{c|cc}
{} & e_1 & e_2 
\\
\hline
e_1 & 0 & -2e_1
\\
e_2 & 2e_1 & 0
\end{array}
\]

%%%%%%%%%%%%%%%%%%%%%%%%%%%%%%%%%%%%%%%%%%%%%%%%%%%%%%%%%%%%%%%%%%%%%%
%%%%%%%%%%%%%%%%%%%%%%%%%%%%%%%%%%%%%%%%%%%%%%%%%%%%%%%%%%%%%%%%%%%%%%
%%%%%%%%%%%%%%%%%%%%%%%%%%%%%%%%%%%%%%%%%%%%%%%%%%%%%%%%%%%%%%%%%%%%%%

\[
\text{\bf Model 2e3b4ba}
\ \ \ \ \
\left\{
\aligned
u
&
\,=\,
xy
+
xy^2
+
xy^3
+
x^3y
+
xy^4+
(
\tfrac{-8}{5}G_{4,0}
+
\tfrac{18}{5}x^3y^2
+
\tfrac{3}{25}
-
\tfrac{2}{5}G_{4,0}^2
+
\\
&
\ \ \ \ \
+
\tfrac{2}{5}G_{4,0}
(
\tfrac{4}{5}G_{4,0}
-
\tfrac{4}{5}+
\tfrac{7}{25}G_{4,0}
)x^5
+
\cdots,
\\
v
&
\,=\,
x^2
-
x^2y^2
+
G_{4,0}x^4
-
2x^2y^3
+
(
\tfrac{4}{5}G_{4,0}
-
\tfrac{4}{5}
)x^4y
+
\cdots,
\endaligned\right.
\]
%%%%%%%%%%%%%%%%%%%%%%%%%%%%%%%%%%%%%%%%%%%%%%%%%%%%%%%%%%%%%%%%%%%%%%
Gr\"obner basis generator of 
moduli space core algebraic variety in 
$\R^1 \ni G_{4,0}$:
\[
\B_1 := (4G_{4,0}+1)(G_{4,0}-1),
\]
%%%%%%%%%%%%%%%%%%%%%%%%%%%%%%%%%%%%%%%%%%%%%%%%%%%%%%%%%%%%%%%%%%%%%%
\[
\aligned
e_1
&
\,:=\,
-
(
2G_{4,0}v
+
y
-
1
)\partial_x
+
(
-
\tfrac{3}{5}x
-
\tfrac{2}{5}xG_{4,0}
-
\tfrac{9}{5}u
+
\tfrac{14}{5}uG_{4,0}
)\partial_y
-
(
-
y
+
\tfrac{3}{5}v
+
\tfrac{2}{5}G_{4,0}v
)\partial_u
+
\\
&
\ \ \ \ \
-
(
-
2x
+
2u
)\partial_v
,
\\
e_2
&
\,:=\,
(
\tfrac{8}{5}xG_{4,0}
-
\tfrac{8}{5}x
+
u
)\partial_x
-
(
-
1
+
2y
+
v
)\partial_y
+
(
x
+
\tfrac{8}{5}uG_{4,0}
-
\tfrac{8}{5}u
)\partial_u
+
(
\tfrac{16}{5}G_{4,0}v
-
\tfrac{16}{5}v
)\partial_v,
\endaligned
\]

%%%%%%%%%%%%%%%%%%%%%%%%%%%%%%%%%%%%%%%%%%%%%%%%%%%%%%%%%%%%%%%%%%%%%%
\[
\footnotesize
\def\arraystretch{1.25}
\begin{array}{c|cc}
{} & e_1 & e_2 
\\
\hline
e_1 & 0 & (-\tfrac{3}{5}+\tfrac{8}{5}G_{4,0})e_1
\\
e_2 & -(-\tfrac{3}{5}+\tfrac{8}{5}G_{4,0})e_1 & 0
\end{array}
\]

%%%%%%%%%%%%%%%%%%%%%%%%%%%%%%%%%%%%%%%%%%%%%%%%%%%%%%%%%%%%%%%%%%%%%%
%%%%%%%%%%%%%%%%%%%%%%%%%%%%%%%%%%%%%%%%%%%%%%%%%%%%%%%%%%%%%%%%%%%%%%
%%%%%%%%%%%%%%%%%%%%%%%%%%%%%%%%%%%%%%%%%%%%%%%%%%%%%%%%%%%%%%%%%%%%%%

\[
\text{\bf Model 2e3b4bb}
\ \ \ \ \
\left\{
\aligned
u
&
\,=\,
xy
+
xy^2
+
xy^3
-
x^3y
+
xy^4
+
(
-
\tfrac{8}{5}G_{4,0}
-
\tfrac{18}{5}
)x^3y^2
+
(
\tfrac{3}{25}
-
\tfrac{2}{5}G_{4,0}^2
+
\\
&
\ \ \ \ \
+
\tfrac{2}{5}G_{4,0}
(
\tfrac{4}{5}G_{4,0}
+
\tfrac{4}{5}
)
-
\tfrac{7}{25}G_{4,0}
)x^5
+
\cdots,
\\
v
&
\,=\,
x^2
-
x^2y^2
+
G_{4,0}x^4
-
2x^2y^3
+
(
\tfrac{4}{5}G_{4,0}
+
\tfrac{4}{5}
)x^4y
+
\cdots,
\endaligned\right.
\]
%%%%%%%%%%%%%%%%%%%%%%%%%%%%%%%%%%%%%%%%%%%%%%%%%%%%%%%%%%%%%%%%%%%%%%
Gr\"obner basis generator of 
moduli space core algebraic variety in 
$\R^1 \ni G_{4,0}$:
\[
\B_1 := 4G_{4,0}^2+3G_{4,0}-1,
\]
%%%%%%%%%%%%%%%%%%%%%%%%%%%%%%%%%%%%%%%%%%%%%%%%%%%%%%%%%%%%%%%%%%%%%%
\[
\aligned
e_1
&
\,:=\,
-
(
2G_{4,0}v
+
y
-
1
)\partial_x
+
(
\tfrac{3}{5}x
-
\tfrac{2}{5}xG_{4,0}
+
\tfrac{9}{5}u
+
\tfrac{14}{5}uG_{4,0}
)\partial_y
-
(
-
y
-
\tfrac{3}{5}v
+
\tfrac{2}{5}G_{4,0}v
)\partial_u
+
\\
&
\ \ \ \ \
-
(
-
2x
+
2u
)\partial_v
,
\\
e_2
&
\,:=\,
-
(
\tfrac{8}{5}xG_{4,0}
+
\tfrac{8}{5}x
-
u
)\partial_x
+
(
1
-
2y
+
v
)\partial_y
-
(
-
x
+
\tfrac{8}{5}uG_{4,0}
+
\tfrac{8}{5}u
)\partial_u
-
(
\tfrac{16}{5}G_{4,0}v
+
\tfrac{16}{5}v
)\partial_v,
\endaligned
\]
%%%%%%%%%%%%%%%%%%%%%%%%%%%%%%%%%%%%%%%%%%%%%%%%%%%%%%%%%%%%%%%%%%%%%%
\[
\footnotesize
\def\arraystretch{1.25}
\begin{array}{c|cc}
{} & e_1 & e_2 
\\
\hline
e_1 & 0 & (-\tfrac{3}{5}-\tfrac{8}{5}G_{4,0})e_1
\\
e_2 & -(-\tfrac{3}{5}-\tfrac{8}{5}G_{4,0})e_1 & 0
\end{array}
\]

%%%%%%%%%%%%%%%%%%%%%%%%%%%%%%%%%%%%%%%%%%%%%%%%%%%%%%%%%%%%%%%%%%%%%%
%%%%%%%%%%%%%%%%%%%%%%%%%%%%%%%%%%%%%%%%%%%%%%%%%%%%%%%%%%%%%%%%%%%%%%
%%%%%%%%%%%%%%%%%%%%%%%%%%%%%%%%%%%%%%%%%%%%%%%%%%%%%%%%%%%%%%%%%%%%%%

\[
\text{\bf Model 2e3b4e}
\ \ \ \ \
\left\{
\aligned
u
&
\,=\,
xy^6+xy^5+xy^4+xy^3+xy^2+xy
+\cdots
,
\\
v
&
\,=\,
-4x^2y^5-3x^2y^4-2x^2y^3-x^2y^2+x^2
+\cdots
.
\endaligned\right.
\]
%%%%%%%%%%%%%%%%%%%%%%%%%%%%%%%%%%%%%%%%%%%%%%%%%%%%%%%%%%%%%%%%%%%%%%
\[
\def\arraystretch{1.25}
\begin{array}{ll}
e_1
\,:=\,
-(y-1)\partial_x+y\partial_u-(2u-2x)\partial_v
, &
e_2
\,:=\,
u\partial_x-(2y-1)\partial_y+x\partial_u
, 
\\
e_3
\,:=\,
\partial_uu+2\partial_vv+\partial_xx
.
\end{array}
\]

%%%%%%%%%%%%%%%%%%%%%%%%%%%%%%%%%%%%%%%%%%%%%%%%%%%%%%%%%%%%%%%%%%%%%%
\[
\footnotesize
\def\arraystretch{1.25}
\begin{array}{c|ccc}
{} & e_1 & e_2 & e_3
\\
\hline
e_1 & 0 & e_1 & e_1
\\
e_2 & -e_1 & 0 & 0
\\
e_3 & -e_1 & 0 & 0
\end{array}
\]

%%%%%%%%%%%%%%%%%%%%%%%%%%%%%%%%%%%%%%%%%%%%%%%%%%%%%%%%%%%%%%%%%%%%%%
%%%%%%%%%%%%%%%%%%%%%%%%%%%%%%%%%%%%%%%%%%%%%%%%%%%%%%%%%%%%%%%%%%%%%%
%%%%%%%%%%%%%%%%%%%%%%%%%%%%%%%%%%%%%%%%%%%%%%%%%%%%%%%%%%%%%%%%%%%%%%

\[
\text{\bf Model 2e3c4aa}
\ \ \ \ \
\left\{
\aligned
u
&
\,=\,
xy
+
\tfrac{3}{4}x^2y^2
-
\tfrac{8}{9}x^4
+
\tfrac{1}{6}\sqrt{3}x^3y^2
+
\tfrac{2}{9}\sqrt{3}x^4y
+
\tfrac{9}{8}x^3y^3
+
\tfrac{1}{6}x^4y^2
-
\tfrac{8}{3}x^5y
+
\\
&
\ \ \ \ \
-
\tfrac{40}{27}x^6
+
\tfrac{5}{8}\sqrt{3}x^4y^3
+
\sqrt{3}x^5y^2
-
\tfrac{2}{9}\sqrt{3}x^6y
-
\tfrac{16}{27}\sqrt{3}x^7
+\cdots
,
\\
v
&
\,=\,
x^2
+
x^3y
+
x^4
+
\tfrac{1}{6}\sqrt{3}x^4y
+
\tfrac{2}{9}\sqrt{3}x^5
+
\tfrac{27}{16}x^4y^2
+
\tfrac{17}{6}x^5y
+
x^6
+
\tfrac{3}{4}\sqrt{3}x^5y^2
+
\\
&
\ \ \ \ \
+
\tfrac{5}{3}\sqrt{3}x^6y
+
\tfrac{8}{9}\sqrt{3}x^7 
+\cdots,
\endaligned\right.
\]
%%%%%%%%%%%%%%%%%%%%%%%%%%%%%%%%%%%%%%%%%%%%%%%%%%%%%%%%%%%%%%%%%%%%%%
\[
\aligned
e_1
&
\,:=\,
-
(
-
\tfrac{1}{3}x\sqrt{3}
-
1
+
\tfrac{3}{2}u
+
2v
)\partial_x
+
(
-
\tfrac{16}{9}x\sqrt{3}
-
y\sqrt{3}
+
2u
+
\tfrac{32}{9}v
)\partial_y
+
\\
&
\ \ \ \ \
-
(
-
y
+
\tfrac{2}{3}\sqrt{3}u
+
\tfrac{16}{9}v\sqrt{3}
)\partial_u
+
(
2x
+
\tfrac{2}{3}v\sqrt{3}
)\partial_v
,
\\
e_2
&
\,:=\,
-(
-
\tfrac{1}{4}x\sqrt{3}
+
\tfrac{1}{2}v
)\partial_x
-
(
\tfrac{2}{3}x\sqrt{3}
-
1
+
\tfrac{1}{4}y\sqrt{3}
+
u
)\partial_y
-
(
-
x
+
\tfrac{2}{3}v\sqrt{3}
)\partial_u
+
\tfrac{1}{2}v\sqrt{3}\partial_v
.
\endaligned
\]
%%%%%%%%%%%%%%%%%%%%%%%%%%%%%%%%%%%%%%%%%%%%%%%%%%%%%%%%%%%%%%%%%%%%%%
\[
\footnotesize
\def\arraystretch{1.25}
\begin{array}{c|cc}
{} & e_1 & e_2 
\\
\hline
e_1 & 0 & \tfrac{1}{4}\sqrt{3}e_1+\tfrac{1}{3}\sqrt{3}e_2
\\
e_2 & -\tfrac{1}{4}\sqrt{3}e_1-\tfrac{1}{3}\sqrt{3}e_2 & 0
\end{array}
\]

%%%%%%%%%%%%%%%%%%%%%%%%%%%%%%%%%%%%%%%%%%%%%%%%%%%%%%%%%%%%%%%%%%%%%%
%%%%%%%%%%%%%%%%%%%%%%%%%%%%%%%%%%%%%%%%%%%%%%%%%%%%%%%%%%%%%%%%%%%%%%
%%%%%%%%%%%%%%%%%%%%%%%%%%%%%%%%%%%%%%%%%%%%%%%%%%%%%%%%%%%%%%%%%%%%%%

\[
\text{\bf Model 2e3c4d5a}
\ \ \ \ \
\left\{
\aligned
u
&
\,=\,
xy
+
x^3y
+
F_{4,1}x^4y
+
x^5
+
(
\tfrac{6}{5}F_{4,1}^2
+
\tfrac{2}{5}G_{4,0}^2
+
\tfrac{2}{5}G_{4,0}
+
\tfrac{6}{5}
)x^5y
+
\\
&
\ \ \ \ \
+
F_{6,0}x^6
+
(
\tfrac{13}{5}F_{4,1}
+
\tfrac{4}{5}G_{4,0}^2F_{4,1}
+
\tfrac{6}{5}F_{4,1}G_{4,0}
+
\tfrac{8}{5}F_{4,1}^3)x^6y
+
\\
&
\ \ \ \ \
+
F_{7,0}x^7
+
\cdots
,
\\
v
&
\,=\,
x^2
+
G_{4,0}x^4
+
\tfrac{4}{5}F_{4,1}G_{4,0}x^5
+
(
2G_{4,0}^2
+
\tfrac{4}{5}F_{4,1}^2G_{4,0}
)x^6
+
\\
&
\ \ \ \ \
+
(
\tfrac{136}{35}G_{4,0}^2F_{4,1}
+
\tfrac{32}{35}F_{4,1}^3G_{4,0}
)x^7
+
\cdots,
\endaligned\right.
\]
%%%%%%%%%%%%%%%%%%%%%%%%%%%%%%%%%%%%%%%%%%%%%%%%%%%%%%%%%%%%%%%%%%%%%%
\[
\aligned
e_1
&
\,:=\,
-
(
2F_{4,1}G_{4,0}x
+
2G_{4,0}^2v
-
2F_{4,1}x
-
2G_{4,0}v
-
G_{4,0}
+
1
)\partial_x
-
(
12F_{4,1}G_{4,0}y
-
6F_{6,0}G_{4,0}y
+
\\
&
\ \ \ \ \
-
2G_{4,0}^2u
-
13F_{4,1}y
+
6F_{6,0}y
+
5G_{4,0}u
-
3u
-
5x
)\partial_y
-
(
14F_{4,1}G_{4,0}u
-
6F_{6,0}G_{4,0}u
+
\\
&
\ \ \ \ \
-
15F_{4,1}u
+
6F_{6,0}u
-
G_{4,0}y
-
5v
+
y
)\partial_u
-
(
4F_{4,1}G_{4,0}v
-
4F_{4,1}v
-
2G_{4,0}x
+
2x
)\partial_v
,
\\
e_2
&
\,:=\,
(
F_{4,1}x
+
G_{4,0}
-
1
+
\tfrac{2}{5}F_{4,1}^2G_{4,0}y
+
\tfrac{2}{5}G_{4,0}^3y
-
\tfrac{1}{5}F_{4,1}^2y
-
G_{4,0}^2y
+
\tfrac{4}{5}G_{4,0}y
-
\tfrac{1}{5}y
-
G_{4,0}v
+
\\
&
\ \ \ \ \
+
v
)\partial_y
+
(
G_{4,0}x
-
x
+
\tfrac{2}{5}F_{4,1}^2G_{4,0}u
+
\tfrac{2}{5}G_{4,0}^3u
-
\tfrac{1}{5}F_{4,1}^2u
-
G_{4,0}^2u
+
\tfrac{4}{5}G_{4,0}u
-
\tfrac{1}{5}u
+
F_{4,1}v
)\partial_u,
\endaligned
\]
%%%%%%%%%%%%%%%%%%%%%%%%%%%%%%%%%%%%%%%%%%%%%%%%%%%%%%%%%%%%%%%%%%%%%%
Gr\"obner basis generator of 
moduli space core algebraic variety in 
$\R^3 \ni G_{4,0},F_{4,1},F_{6,0}$:
\[
\B_1 := (2G_{4,0}-1)(G_{4,0}-1)(F_{4,1}^2+G_{4,0}^2-2G_{4,0}+1)(2F_{4,1}-F_{6,0}),
\]
%%%%%%%%%%%%%%%%%%%%%%%%%%%%%%%%%%%%%%%%%%%%%%%%%%%%%%%%%%%%%%%%%%%%%%
\[
\footnotesize
\def\arraystretch{1.25}
\begin{array}{c|cc}
{} & e_1 & e_2 
\\
\hline
e_1 & 0 & (6G_{4,0}-6)(2F_{4,1}-F_{6,0})e_2
\\
e_2 & -(6G_{4,0}-6)(2F_{4,1}-F_{6,0})e_2 & 0
\end{array}
\]

%%%%%%%%%%%%%%%%%%%%%%%%%%%%%%%%%%%%%%%%%%%%%%%%%%%%%%%%%%%%%%%%%%%%%%
%%%%%%%%%%%%%%%%%%%%%%%%%%%%%%%%%%%%%%%%%%%%%%%%%%%%%%%%%%%%%%%%%%%%%%
%%%%%%%%%%%%%%%%%%%%%%%%%%%%%%%%%%%%%%%%%%%%%%%%%%%%%%%%%%%%%%%%%%%%%%

\[
\text{\bf Model 2e3c4d5b}
\ \ \ \ \
\left\{
\aligned
u
&
\,=\,
xy
+
x^3y
+
F_{4,1}x^4y
+
(
\tfrac{6}{5}F_{4,1}^2
+
\tfrac{2}{5}G_{4,0}^2
+
\tfrac{2}{5}G_{4,0}
+
\tfrac{6}{5})x^5y
+
\\
&
\ \ \ \ \
+
(
\tfrac{4}{5}G_{4,0}^2F_{4,1}
+
\tfrac{6}{5}F_{4,1}G_{4,0}
+
\tfrac{8}{5}F_{4,1}^3
+
\tfrac{13}{5}F_{4,1}
)x^6y
+
\cdots,
\\
v
&
\,=\,
x^2
+
G_{4,0}x^4
+
\tfrac{4}{5}F_{4,1}G_{4,0}x^5
+
(
2G_{4,0}^2
+
\tfrac{4}{5}F_{4,1}^2G_{4,0}
)x^6
+
\\
&
\ \ \ \ \
+
(
\tfrac{136}{35}G_{4,0}^2F_{4,1}
+
\tfrac{32}{35}F_{4,1}^3G_{4,0}
)x^7
+
\cdots,
\endaligned\right.
\]
%%%%%%%%%%%%%%%%%%%%%%%%%%%%%%%%%%%%%%%%%%%%%%%%%%%%%%%%%%%%%%%%%%%%%%
\[
\aligned
e_1
&
\,:=\,
-
(
2F_{4,1}x
+
2G_{4,0}v
-
1
)\partial_x
+
\tfrac{
(
2G_{4,0}^2u
-
5G_{4,0}u
+
3u
)
}{G_{4,0}-1}\partial_y
+
\tfrac{
(
-
2F_{4,1}G_{4,0}u
+
2F_{4,1}u
+
G_{4,0}y
-
y
)
}
{G_{4,0}-1}\partial_u
+
\\
&
\ \ \ \ \
-
(
4F_{4,1}v
-
2x
)\partial_v
,
\\
e_2
&
\,:=\,
\tfrac{
(
F_{4,1}x
-
G_{4,0}v
+
G_{4,0}
+
v
-
1
)
}
{G_{4,0}-1}\partial_y
+
\tfrac{(F_{4,1}v+G_{4,0}x-x)}{G_{4,0}-1}\partial_u
, 
\\
e_3
&
\,:=\,
\tfrac{(G_{4,0}y-y)}{G_{4,0}-1}\partial_y
+
\tfrac{(G_{4,0}u-u)}{G_{4,0}-1}\partial_u,
\endaligned
\]
%%%%%%%%%%%%%%%%%%%%%%%%%%%%%%%%%%%%%%%%%%%%%%%%%%%%%%%%%%%%%%%%%%%%%%
Gr\"obner basis generator of 
moduli space core algebraic variety in 
$\R^2 \ni G_{4,0},F_{4,1}$:
\[
\B_1 := 2F_{4,1}^2G_{4,0}+2G_{4,0}^3-F_{4,1}^2-5G_{4,0}^2+4G_{4,0}-1,
\]
%%%%%%%%%%%%%%%%%%%%%%%%%%%%%%%%%%%%%%%%%%%%%%%%%%%%%%%%%%%%%%%%%%%%%%
\[
\footnotesize
\def\arraystretch{1.25}
\begin{array}{c|ccc}
{} & e_1 & e_2 & e_3
\\
\hline
e_1 & 0 & \tfrac{F_{4,1}e_2}{G_{4,0}-1} & 0
\\
e_2 & -\tfrac{F_{4,1}e_2}{G_{4,0}-1} & 0 & e_2
\\
e_3 & -\tfrac{F_{4,1}e_2}{G_{4,0}-1} & 0 & e_2
\end{array}
\]

%%%%%%%%%%%%%%%%%%%%%%%%%%%%%%%%%%%%%%%%%%%%%%%%%%%%%%%%%%%%%%%%%%%%%%
%%%%%%%%%%%%%%%%%%%%%%%%%%%%%%%%%%%%%%%%%%%%%%%%%%%%%%%%%%%%%%%%%%%%%%
%%%%%%%%%%%%%%%%%%%%%%%%%%%%%%%%%%%%%%%%%%%%%%%%%%%%%%%%%%%%%%%%%%%%%%

\[
\text{\bf Model 2e3c4e5a}
\ \ \ \ \
\left\{
\aligned
u
&
\,=\,
xy
-
x^3y
+
F_{4,1}x^4y
+
x^5
+
(
-
\tfrac{6}{5}F_{4,1}^2
+
\tfrac{2}{5}G_{4,0}^2
-
\tfrac{2}{5}G_{4,0}+\tfrac{6}{5}
)x^5y
+
\\
&
\ \ \ \ \
+
F_{6,0}x^6
+
(
-
\tfrac{13}{5}F_{4,1}
-
\tfrac{4}{5}G_{4,0}^2F_{4,1}
+
\tfrac{6}{5}F_{4,1}G_{4,0}
+
\tfrac{8}{5}F_{4,1}^3
)x^6y
+
F_{7,0}x^7
+
\\
&
\ \ \ \ \
+
\cdots
,
\\
v
&
\,=\,
x^2
+
G_{4,0}x^4
-
\tfrac{4}{5}F_{4,1}G_{4,0}x^5
+
(
2G_{4,0}^2
+
\tfrac{4}{5}F_{4,1}^2G_{4,0}
)x^6
+
\\
&
\ \ \ \ \
+
(
-
\tfrac{136}{35}G_{4,0}^2F_{4,1}
-
\tfrac{32}{35}F_{4,1}^3G_{4,0}
)x^7
+
\cdots,
\endaligned\right.
\]
%%%%%%%%%%%%%%%%%%%%%%%%%%%%%%%%%%%%%%%%%%%%%%%%%%%%%%%%%%%%%%%%%%%%%%
\[
\aligned
e_1
&
\,:=\,
(
2F_{4,1}x
-
2G_{4,0}v
+
1
)
\partial_x
+
\\
&
\ \ \ \ \
-
\tfrac{(-60F_{4,1}G_{4,0}y-30F_{6,0}G_{4,0}y-10G_{4,0}^2u-65F_{4,1}y-30F_{6,0}y-25G_{4,0}u-15u-25x)}{5G_{4,0}+5}\partial_y
+
\\
&
\ \ \ \ \
-
\tfrac{(-70F_{4,1}G_{4,0}u-30F_{6,0}G_{4,0}u-75F_{4,1}u-30F_{6,0}u-5G_{4,0}y-25v-5y)}{5G_{4,0}+5}\partial_u+(4F_{4,1}v+2x)\partial_v
,
\\
e_2
&
\,:=\,
-
\tfrac{(2F_{4,1}^2G_{4,0}y-2G_{4,0}^3y+F_{4,1}^2y-5G_{4,0}^2y-5F_{4,1}x-5G_{4,0}v-4G_{4,0}y-5G_{4,0}-5v-y-5)}{5G_{4,0}+5}\partial_y
+
\\
&
\ \ \ \ \
-
\tfrac{(2F_{4,1}^2G_{4,0}u-2G_{4,0}^3u+F_{4,1}^2u-5G_{4,0}^2u-5F_{4,1}v-4G_{4,0}u-5G_{4,0}x-u-5x)}{5G_{4,0}+5}\partial_u,
\endaligned
\]
%%%%%%%%%%%%%%%%%%%%%%%%%%%%%%%%%%%%%%%%%%%%%%%%%%%%%%%%%%%%%%%%%%%%%%
Gr\"obner basis generator of 
moduli space core algebraic variety in 
$\R^3 \ni G_{4,0},F_{4,1},F_{6,0}$:
\[
\B_1 := -\tfrac{6}{5}(2G_{4,0}+1)(F_{4,1}+1+G_{4,0})(F_{4,1}-1-G_{4,0})(2F_{4,1}+F_{6,0}),
\]
%%%%%%%%%%%%%%%%%%%%%%%%%%%%%%%%%%%%%%%%%%%%%%%%%%%%%%%%%%%%%%%%%%%%%%
\[
\footnotesize
\def\arraystretch{1.25}
\begin{array}{c|cc}
{} & e_1 & e_2 
\\
\hline
e_1 & 0 & -\tfrac{(12F_{4,1}G_{4,0}+6F_{6,0}G_{4,0}+12F_{4,1}+6F_{6,0})}{G_{4,0}+1}e_2
\\
e_2 & \tfrac{(12F_{4,1}G_{4,0}+6F_{6,0}G_{4,0}+12F_{4,1}+6F_{6,0})}{G_{4,0}+1}e_2 & 0
\end{array}
\]

%%%%%%%%%%%%%%%%%%%%%%%%%%%%%%%%%%%%%%%%%%%%%%%%%%%%%%%%%%%%%%%%%%%%%%
%%%%%%%%%%%%%%%%%%%%%%%%%%%%%%%%%%%%%%%%%%%%%%%%%%%%%%%%%%%%%%%%%%%%%%
%%%%%%%%%%%%%%%%%%%%%%%%%%%%%%%%%%%%%%%%%%%%%%%%%%%%%%%%%%%%%%%%%%%%%%

\[
\text{\bf Model 2e3c4e5b}
\ \ \ \ \
\left\{
\aligned
u
&
\,=\,
xy
-
x^3y
+
F_{4,1}x^4y
+
(
-
\tfrac{6}{5}F_{4,1}^2
+
\tfrac{2}{5}G_{4,0}^2
-
\tfrac{2}{5}G_{4,0}
+
\tfrac{6}{5})x^5y
+
\\
&
\ \ \ \ \
+
(
-
\tfrac{4}{5}G_{4,0}^2F_{4,1}
+
\tfrac{6}{5}F_{4,1}G_{4,0}
+
\tfrac{8}{5}F_{4,1}^3
-
\tfrac{13}{5}F_{4,1}
)x^6y
+
\cdots
,
\\
v
&
\,=\,
x^2
+
G_{4,0}x^4
-
\tfrac{4}{5}F_{4,1}G_{4,0}x^5
+
(
2G_{4,0}^2
+
\tfrac{4}{5}F_{4,1}^2G_{4,0}
)x^6
+
\\
&
\ \ \ \ \
+
(
-
\tfrac{136}{35}G_{4,0}^2F_{4,1}
-
\tfrac{32}{35}F_{4,1}^3G_{4,0}
)x^7
+
\cdots,
\endaligned\right.
\]
%%%%%%%%%%%%%%%%%%%%%%%%%%%%%%%%%%%%%%%%%%%%%%%%%%%%%%%%%%%%%%%%%%%%%%
\[
\aligned
e_1
&
\,:=\,
2F_{4,1}u\partial_u
+
4F_{4,1}v\partial_v
+
2F_{4,1}x\partial_x
-
2G_{4,0}v\partial_x
+
2G_{4,0}u\partial_y
+
y\partial_u
+
2u\partial_v
+
3u\partial_y
+
\partial_x
, 
\\
e_2
&
\,:=\,
\tfrac{(F_{4,1}x+G_{4,0}v+G_{4,0}+v+1)}{G_{4,0}+1}\partial_y+
\tfrac{(F_{4,1}v+G_{4,0}x+x)}{G_{4,0}+1}\partial_u
, 
\\
e_3
&
\,:=\,
\tfrac{(G_{4,0}y+y)}{G_{4,0}+1}\partial_y+\tfrac{(G_{4,0}u+u)}{G_{4,0}+1}\partial_u,
\endaligned
\]
%%%%%%%%%%%%%%%%%%%%%%%%%%%%%%%%%%%%%%%%%%%%%%%%%%%%%%%%%%%%%%%%%%%%%%
Gr\"obner basis generator of 
moduli space core algebraic variety in 
$\R^2 \ni G_{4,0},F_{4,1}$:
\[
\B_1 := 2F_{4,1}^2G_{4,0}-2G_{4,0}^3+F_{4,1}^2-5G_{4,0}^2-4G_{4,0}-1,
\]
%%%%%%%%%%%%%%%%%%%%%%%%%%%%%%%%%%%%%%%%%%%%%%%%%%%%%%%%%%%%%%%%%%%%%%
\[
\footnotesize
\def\arraystretch{1.25}
\begin{array}{c|ccc}
{} & e_1 & e_2 & e_3
\\
\hline
e_1 & 0 & \tfrac{F_{4,1}}{G_{4,0}+1}e_2
\\
e_2 & -\tfrac{F_{4,1}}{G_{4,0}+1}e_2 & 0 & e_2
\\
e_3 & -\tfrac{F_{4,1}}{G_{4,0}+1}e_2 & 0 & e_2
\end{array}
\]

%%%%%%%%%%%%%%%%%%%%%%%%%%%%%%%%%%%%%%%%%%%%%%%%%%%%%%%%%%%%%%%%%%%%%%
%%%%%%%%%%%%%%%%%%%%%%%%%%%%%%%%%%%%%%%%%%%%%%%%%%%%%%%%%%%%%%%%%%%%%%
%%%%%%%%%%%%%%%%%%%%%%%%%%%%%%%%%%%%%%%%%%%%%%%%%%%%%%%%%%%%%%%%%%%%%%

\[
\text{\bf Model 2e3c4f5a}
\ \ \ \ \
\left\{
\aligned
u
&
\,=\,
xy
+
x^3y
+
x^4
+
\tfrac{5}{4}G_{5,0}x^4y
+
(
2
+
\tfrac{1}{8}G_{5,0}^2
)x^5y
+
(
-
\tfrac{5}{4}G_{5,0}^2
+
\tfrac{14}{5}
)x^6
\\
&
\ \ \ \ \
+
(
\tfrac{23}{4}G_{5,0}
+
\tfrac{25}{8}G_{5,0}^3
)x^6y
+
(
\tfrac{20}{7}(
-
\tfrac{5}{4}G_{5,0}^2
+
\tfrac{14}{5}
)G_{5,0}
-
\tfrac{67}{14}G_{5,0}
)x^7
+
\cdots
,
\\
v
&
\,=\,
x^2+x^4+G_{5,0}x^5+(2+(5}{4)G_{5,0}^2)x^6+((34}{7)G_{5,0}+(25}{14)G_{5,0}^3)x^7 
+
\cdots
\endaligned\right.
\]
%%%%%%%%%%%%%%%%%%%%%%%%%%%%%%%%%%%%%%%%%%%%%%%%%%%%%%%%%%%%%%%%%%%%%%
\[
\aligned
e_1
&
\,:=\,
-(
\tfrac{5}{2}G_{5,0}x
+
2v
-
1
)\partial_x
+
\tfrac{1}{20}\tfrac{(600G_{5,0}^2x-150G_{5,0}^2y-20G_{5,0}u-80G_{5,0}v+16x)}{G_{5,0}}\partial_y
\\
&
-
\tfrac{1}{20}\tfrac{(200G_{5,0}^2u-600G_{5,0}^2v-20G_{5,0}y-16v)}{G_{5,0}}\partial_u
-
(5G_{5,0}v-2x)\partial_v
,
\\
e_2
&
\,:=\,
-
\tfrac{1}{4}(30G_{5,0}x-5G_{5,0}y+4v-4)\partial_y
+
\tfrac{1}{4}(5G_{5,0}u-30G_{5,0}v+4x)\partial_u 
.
\endaligned
\]
%%%%%%%%%%%%%%%%%%%%%%%%%%%%%%%%%%%%%%%%%%%%%%%%%%%%%%%%%%%%%%%%%%%%%%
\[
\footnotesize
\def\arraystretch{1.25}
\begin{array}{c|cc}
{} & e_1 & e_2 
\\
\hline
e_1 & 0 & 0
\\
e_2 & 0 & 0
\end{array}
\]

%%%%%%%%%%%%%%%%%%%%%%%%%%%%%%%%%%%%%%%%%%%%%%%%%%%%%%%%%%%%%%%%%%%%%%
%%%%%%%%%%%%%%%%%%%%%%%%%%%%%%%%%%%%%%%%%%%%%%%%%%%%%%%%%%%%%%%%%%%%%%
%%%%%%%%%%%%%%%%%%%%%%%%%%%%%%%%%%%%%%%%%%%%%%%%%%%%%%%%%%%%%%%%%%%%%%

\[
\text{\bf Model 2e3c4f5b}
\ \ \ \ \
\left\{
\aligned
u
&
\,=\,
xy+x^3y+x^4
+
F_{5,0}x^5
+
2x^5y
+
(
\tfrac{5}{6}F_{5,0}^2+\tfrac{17}{6}
)x^6
+
\\
&
\ \ \ \ \
+
(
\tfrac{151}{42}F_{5,0}
+
\tfrac{25}{42}F_{5,0}^3
)x^7
+
\cdots
,
\\
v
&
\,=\,
x^2+x^4+2x^6
+
\cdots
,
\endaligned\right.
\]
for any value of $F_{5,0}$,
%%%%%%%%%%%%%%%%%%%%%%%%%%%%%%%%%%%%%%%%%%%%%%%%%%%%%%%%%%%%%%%%%%%%%%
\[
\aligned
e_1
&
\,:=\,
-(2v-1)\partial_x+(5F_{5,0}y-u-4v)\partial_y+(5F_{5,0}u+y)\partial_u+2x\partial_v
,
\\
e_2
&
\,:=\,
-(v-1)\partial_y+x\partial_u
,
\\
e_3
&
\,:=\,
v\partial_u+x\partial_y
.
\endaligned
\]

%%%%%%%%%%%%%%%%%%%%%%%%%%%%%%%%%%%%%%%%%%%%%%%%%%%%%%%%%%%%%%%%%%%%%%
\[
\footnotesize
\def\arraystretch{1.25}
\begin{array}{c|ccc}
{} & e_1 & e_2 & e_3
\\
\hline
e_1 & 0 & -5F_{5,0}e_2-e_3 & -5F_{5,0}e_3+e_2
\\
e_2 & 5F_{5,0}e_2+e_3 & 0 & 0
\\
e_3 & 5F_{5,0}e_2+e_3 & 0 & 0
\end{array}
\]

%%%%%%%%%%%%%%%%%%%%%%%%%%%%%%%%%%%%%%%%%%%%%%%%%%%%%%%%%%%%%%%%%%%%%%
%%%%%%%%%%%%%%%%%%%%%%%%%%%%%%%%%%%%%%%%%%%%%%%%%%%%%%%%%%%%%%%%%%%%%%
%%%%%%%%%%%%%%%%%%%%%%%%%%%%%%%%%%%%%%%%%%%%%%%%%%%%%%%%%%%%%%%%%%%%%%

\[
\text{\bf Model 2e3c4g}
\ \ \ \ \
\left\{
\aligned
u
&
\,=\,
xy+x^3y+2x^5y
+
\cdots
,
\\
v
&
\,=\,
2x^6+x^4+x^2 
+
\cdots,
\endaligned\right.
\]
%%%%%%%%%%%%%%%%%%%%%%%%%%%%%%%%%%%%%%%%%%%%%%%%%%%%%%%%%%%%%%%%%%%%%%
\[
\aligned
e_1
&
\,:=\,
-(2v-1)\partial_x-u\partial_y+y\partial_u+2x\partial_v 
,
\ \ \ \ \
e_2
\,:=\,
-(v-1)\partial_y+x\partial_u
, 
\\
e_3
&
\,:=\,
v\partial_u+x\partial_y
,
\ \ \ \ \ \ \ \ \ \ \ \ \ \ \ \ \ \ \ \ \ \ \ \ \ \ \ \ \ \ \ \ \ \ \ \ \ \ \ \ \ \ \ \ 
e_4
\,:=\,
u\partial_u+y\partial_y,
\endaligned
\]

%%%%%%%%%%%%%%%%%%%%%%%%%%%%%%%%%%%%%%%%%%%%%%%%%%%%%%%%%%%%%%%%%%%%%%
\[
\footnotesize
\def\arraystretch{1.25}
\begin{array}{c|cccc}
{} & e_1 & e_2 & e_3 & e_4
\\
\hline
e_1 & 0 & -e_3 & e_2 & 0
\\
e_2 & e_3 & 0 & 0 & e_2
\\
e_3 & -e_2 & 0 & 0 & e_3
\\
e_4 & 0 & -e_2 & -e_3 & 0
\end{array}
\]

%%%%%%%%%%%%%%%%%%%%%%%%%%%%%%%%%%%%%%%%%%%%%%%%%%%%%%%%%%%%%%%%%%%%%%
%%%%%%%%%%%%%%%%%%%%%%%%%%%%%%%%%%%%%%%%%%%%%%%%%%%%%%%%%%%%%%%%%%%%%%
%%%%%%%%%%%%%%%%%%%%%%%%%%%%%%%%%%%%%%%%%%%%%%%%%%%%%%%%%%%%%%%%%%%%%%

\[
\text{\bf Model 2e3c4h}
\ \ \ \ \
\left\{
\aligned
u
&
\,=\,
xy-x^3y+x^4+F_{5,0}x^5+2x^5y
+
(
\tfrac{5}{6}F_{5,0}^2
-
\tfrac{17}{6}
)x^6
+
\\
&
\ \ \ \ \
+
(
-
\tfrac{151}{42}F_{5,0}
+
\tfrac{25}{42}F_{5,0}^3
)x^7
+
\cdots,
\\
v
&
\,=\,
x^2-x^4+2x^6\cdots
,
\endaligned\right.
\]

%%%%%%%%%%%%%%%%%%%%%%%%%%%%%%%%%%%%%%%%%%%%%%%%%%%%%%%%%%%%%%%%%%%%%%
\[
\aligned
e_1
&
\,:=\,
(2v+1)\partial_x+(5F_{5,0}y+u-4v)\partial_y+(5F_{5,0}u+y)\partial_u+2x\partial_v
,
\\
e_2
&
\,:=\,
(v+1)\partial_y+x\partial_u
,
& 
\\
e_3
&
\,:=\,
v\partial_u+x\partial_y
,
\endaligned
\]
for any value of $F_{5,0}$,
%%%%%%%%%%%%%%%%%%%%%%%%%%%%%%%%%%%%%%%%%%%%%%%%%%%%%%%%%%%%%%%%%%%%%%
\[
\footnotesize
\def\arraystretch{1.25}
\begin{array}{c|ccc}
{} & e_1 & e_2 & e_3
\\
\hline
e_1 & 0 & -5F_{5,0}e_2+e_3 & -5F_{5,0}e_3+e_2
\\
e_2 & 5F_{5,0}e_2-e_3 & 0 & 0
\\
e_3 & 5F_{5,0}e_3-e_2 & 0 & 0
\end{array}
\]
%%%%%%%%%%%%%%%%%%%%%%%%%%%%%%%%%%%%%%%%%%%%%%%%%%%%%%%%%%%%%%%%%%%%%%
%%%%%%%%%%%%%%%%%%%%%%%%%%%%%%%%%%%%%%%%%%%%%%%%%%%%%%%%%%%%%%%%%%%%%%
%%%%%%%%%%%%%%%%%%%%%%%%%%%%%%%%%%%%%%%%%%%%%%%%%%%%%%%%%%%%%%%%%%%%%%

\[
\text{\bf Model 2e3c4j5a}
\ \ \ \ \
\left\{
\aligned
u
&
\,=\,
xy+x^5
+
\tfrac{2}{5}x^5y
+
F_{6,0}x^6
+
\tfrac{2}{5}F_{6,0}x^6y
+
(\tfrac{16}{7}+\tfrac{6}{7}F_{6,0}^2)x^7
+
\cdots
,
\\
v
&
\,=\,
x^2+x^4+\tfrac{2}{5}F_{6,0}x^5+(2+\tfrac{1}{5}F_{6,0}^2)x^6+(\tfrac{68}{35}F_{6,0}+\tfrac{4}{35}F_{6,0}^3)x^7
+\cdots,
\endaligned\right.
\]
%%%%%%%%%%%%%%%%%%%%%%%%%%%%%%%%%%%%%%%%%%%%%%%%%%%%%%%%%%%%%%%%%%%%%%
\[
\aligned
e_1
&
\,:=\,
-(F_{6,0}x+2v-1)\partial_x+(5x+2u)\partial_y-(F_{6,0}u-5v-y)\partial_u-(2F_{6,0}v-2x)\partial_v
, 
\\
e_2
&
\,:=\,
(1+\tfrac{2}{5}y)\partial_y+(x+\tfrac{2}{5}u)\partial_u,
\endaligned
\]
for any value of $F_{6,0}$,
%%%%%%%%%%%%%%%%%%%%%%%%%%%%%%%%%%%%%%%%%%%%%%%%%%%%%%%%%%%%%%%%%%%%%%
\[
\footnotesize
\def\arraystretch{1.25}
\begin{array}{c|cc}
{} & e_1 & e_2 
\\
\hline
e_1 & 0 & 0
\\
e_2 & 0 & 0
\end{array}
\]

%%%%%%%%%%%%%%%%%%%%%%%%%%%%%%%%%%%%%%%%%%%%%%%%%%%%%%%%%%%%%%%%%%%%%%
%%%%%%%%%%%%%%%%%%%%%%%%%%%%%%%%%%%%%%%%%%%%%%%%%%%%%%%%%%%%%%%%%%%%%%
%%%%%%%%%%%%%%%%%%%%%%%%%%%%%%%%%%%%%%%%%%%%%%%%%%%%%%%%%%%%%%%%%%%%%%

\[
\text{\bf Model 2e3c4k5a}
\ \ \ \ \
\left\{
\aligned
u
&
\,=\,
xy+x^5+\tfrac{2}{5}x^5y+F_{6,0}x^6+\tfrac{2}{5}F_{6,0}x^6y+(-\tfrac{16}{7}+\tfrac{6}{7}F_{6,0}^2)x^7
+
\cdots
,
\\
v
&
\,=\,
x^2-x^4
-
\tfrac{2}{5}F_{6,0}x^5
+
(
2
-
\tfrac{1}{5}F_{6,0}^2
)x^6
+
(
\tfrac{68}{35}F_{6,0}
-
\tfrac{4}{35}F_{6,0}^3)x^7
+
\cdots
,
\endaligned\right.
\]
%%%%%%%%%%%%%%%%%%%%%%%%%%%%%%%%%%%%%%%%%%%%%%%%%%%%%%%%%%%%%%%%%%%%%%
\[
\aligned
e_1
&
\,:=\,
-(F_{6,0}x-2v-1)\partial_x-(5x+2u)\partial_y-(F_{6,0}u+5v-y)\partial_u-(2F_{6,0}v-2x)\partial_v 
, 
\\
e_2
&
\,:=\,
(1+\tfrac{2}{5}y)\partial_y+(x+\tfrac{2}{5}u)\partial_u
,
\endaligned
\]
for any value of $F_{6,0}$.
%%%%%%%%%%%%%%%%%%%%%%%%%%%%%%%%%%%%%%%%%%%%%%%%%%%%%%%%%%%%%%%%%%%%%%
\[
\footnotesize
\def\arraystretch{1.25}
\begin{array}{c|cc}
{} & e_1 & e_2 
\\
\hline
e_1 & 0 & 0 
\\
e_2 & 0 & 0
\end{array}
\]

%%%%%%%%%%%%%%%%%%%%%%%%%%%%%%%%%%%%%%%%%%%%%%%%%%%%%%%%%%%%%%%%%%%%%%
%%%%%%%%%%%%%%%%%%%%%%%%%%%%%%%%%%%%%%%%%%%%%%%%%%%%%%%%%%%%%%%%%%%%%%
%%%%%%%%%%%%%%%%%%%%%%%%%%%%%%%%%%%%%%%%%%%%%%%%%%%%%%%%%%%%%%%%%%%%%%

\[
\text{\bf Model 2e3c4l5a}
\ \ \ \ \
\left\{
\aligned
u
&
\,=\,
xy+x^4+x^5+F_{6,0}x^6+(\tfrac{12}{7}F_{6,0}^2-\tfrac{5}{7}F_{6,0})x^7
+
\cdots
,
\\
v
&
\,=\,
x^2,
\endaligned\right.
\]
%%%%%%%%%%%%%%%%%%%%%%%%%%%%%%%%%%%%%%%%%%%%%%%%%%%%%%%%%%%%%%%%%%%%%%
\[
\aligned
e_1
&
\,:=\,
-(6F_{6,0}x-5x-1)\partial_x-(18F_{6,0}y+4v-20y)\partial_y-(24F_{6,0}u-25u-y)\partial_u
+
\\
&
\ \ \ \ \
-(12F_{6,0}v-10v-2x)\partial_v
,
\\
e_2
&
\,:=\,
x\partial_u+\partial_y, 
\\
e_3
&
\,:=\,
v\partial_u+x\partial_y,
\endaligned
\]

%%%%%%%%%%%%%%%%%%%%%%%%%%%%%%%%%%%%%%%%%%%%%%%%%%%%%%%%%%%%%%%%%%%%%%
\[
\footnotesize
\def\arraystretch{1.25}
\begin{array}{c|ccc}
{} & e_1 & e_2 & e_3
\\
\hline
e_1 & 0 & (18F_{6,0}-20)e_2 & e_2+(12F_{6,0}-15)e_3
\\
e_2 & -(18F_{6,0}-20)e_2 & 0 & 0 
\\
e_3 & -e_2+(12F_{6,0}-15)e_3 & 0 & 0
\end{array}
\]

%%%%%%%%%%%%%%%%%%%%%%%%%%%%%%%%%%%%%%%%%%%%%%%%%%%%%%%%%%%%%%%%%%%%%%
%%%%%%%%%%%%%%%%%%%%%%%%%%%%%%%%%%%%%%%%%%%%%%%%%%%%%%%%%%%%%%%%%%%%%%
%%%%%%%%%%%%%%%%%%%%%%%%%%%%%%%%%%%%%%%%%%%%%%%%%%%%%%%%%%%%%%%%%%%%%%

\[
\text{\bf Model 2e3c4m}
\ \ \ \ \
\left\{
\aligned
u
&
\,=\,
xy
,
\\
v
&
\,=\,
x^2,
\endaligned\right.
\]
%%%%%%%%%%%%%%%%%%%%%%%%%%%%%%%%%%%%%%%%%%%%%%%%%%%%%%%%%%%%%%%%%%%%%%
\[
\def\arraystretch{1.25}
\begin{array}{llll}
e_1
\,:=\,
y\partial_u+2x\partial_v+\partial_x
, &
e_2
\,:=\,
x\partial_u+\partial_y
, &
& 
\\
e_3
\,:=\,
u\partial_u+2v\partial_v+x\partial_x
, &
e_4
\,:=\,
v\partial_u+x\partial_y
, &
e_5
\,:=\,
u\partial_u+y\partial_y,
\end{array}
\]

%%%%%%%%%%%%%%%%%%%%%%%%%%%%%%%%%%%%%%%%%%%%%%%%%%%%%%%%%%%%%%%%%%%%%%
\[
\footnotesize
\def\arraystretch{1.25}
\begin{array}{c|ccccc}
{} & e_1 & e_2 & e_3 & e_4 & e_5 
\\
\hline
e_1 & 0 & 0 & e_1 & e_2 & 0
\\
e_2 & 0 & 0 & 0 & 0 & e_2
\\
e_3 & -e_1 & 0 & 0 & e_4 & 0
\\
e_4 & -e_2 & 0 & -e_4 & 0 & e_4
\\
e_5 & 0 & -e_2 & 0 & -e_4 & 0
\end{array}
\]
%%%%%%%%%%%%%%%%%%%%%%%%%%%%%%%%%%%%%%%%%%%%%%%%%%%%%%%%%%%%%%%%%%%%%%
%%%%%%%%%%%%%%%%%%%%%%%%%%%%%%%%%%%%%%%%%%%%%%%%%%%%%%%%%%%%%%%%%%%%%%
%%%%%%%%%%%%%%%%%%%%%%%%%%%%%%%%%%%%%%%%%%%%%%%%%%%%%%%%%%%%%%%%%%%%%%

%%%%%%%%%%%%%%%%%%%%%%%%%%%%%%%%%%%%%%%%%%%%%%%%%%%%%%%%%%%%%%%%%%%%%%
\SectionHead{2f Models}
{2f-models}
%%%%%%%%%%%%%%%%%%%%%%%%%%%%%%%%%%%%%%%%%%%%%%%%%%%%%%%%%%%%%%%%%%%%%%

%%%%%%%%%%%%%%%%%%%%%%%%%%%%%%%%%%%%%%%%%%%%%%%%%%%%%%%%%%%%%%%%%%%%%%
%%%%%%%%%%%%%%%%%%%%%%%%%%%%%%%%%%%%%%%%%%%%%%%%%%%%%%%%%%%%%%%%%%%%%%
%%%%%%%%%%%%%%%%%%%%%%%%%%%%%%%%%%%%%%%%%%%%%%%%%%%%%%%%%%%%%%%%%%%%%%

\[
\text{\bf Model 2f3a}
\ \ \ \ \
\left\{
\aligned
u
&
\,=\,
xy+F_{0,3}y^3+(4-F_{0,3})x^3+F_{0,4}y^4+F_{1,3}xy^3+F_{2,2}x^2y^2+F_{3,1}x^3y
+
\\
&
\ \ \ \ \
+
F_{4,0}x^4
+
\cdots,
\\
v
&
\,=\,
y^2+x^2+(-8-G_{2,1})xy^2+G_{2,1}x^2y+G_{0,4}y^4+G_{1,3}xy^3
+
\cdots.
\endaligned\right.
\]

%%%%%%%%%%%%%%%%%%%%%%%%%%%%%%%%%%%%%%%%%%%%%%%%%%%%%%%%%%%%%%%%%%%%%%
\[
\aligned
e_1
&
\,:=\,
\Big(
1-\tfrac{3}{8}uG_{2,1}F_{4,0}+\tfrac{3}{64}uG_{2,1}G_{1,3}+\tfrac{1}{16}vG_{2,1}F_{2,2}+\tfrac{15}{16}uF_{0,3}F_{3,1}-\tfrac{3}{2}vF_{3,1}-\tfrac{9}{32}xG_{3,1}
+
\\
&
\ \ \ \ \
+\tfrac{3}{64}xG_{2,1}^2+\tfrac{1}{16}vG_{2,1}G_{4,0}-3vG_{4,0}+\tfrac{27}{16}xF_{0,3}^2-\tfrac{9}{8}uG_{3,1}-\tfrac{3}{16}yF_{3,1}+\tfrac{1}{16}yG_{2,2}-\tfrac{3}{4}uG_{2,2}
+
\\
&
\ \ \ \ \
-\tfrac{3}{32}uG_{2,1}G_{2,2}-\tfrac{1}{8}vF_{0,3}F_{1,3}+\tfrac{57}{32}uF_{0,3}^2G_{2,1}+\tfrac{3}{8}vF_{0,3}F_{3,1}+\tfrac{5}{8}uF_{0,3}F_{2,2}+\tfrac{9}{32}uG_{2,1}F_{3,1}
+
\\
&
\ \ \ \ \
+17vG_{2,1}+21uF_{0,3}+42vF_{0,3}-\tfrac{21}{2}uG_{2,1}+\tfrac{3}{16}uG_{2,1}G_{4,0}+\tfrac{9}{16}xF_{3,1}-\tfrac{3}{16}xG_{2,2}
+
\\
&
\ \ \ \ \
-\tfrac{55}{8}vF_{0,3}G_{2,1}+\tfrac{5}{32}uF_{0,3}G_{1,3}+\tfrac{7}{8}uG_{1,3}+\tfrac{1}{64}vG_{2,1}G_{1,3}-\tfrac{1}{8}yF_{2,2}-\tfrac{7}{4}xF_{4,0}-\tfrac{1}{2}uF_{2,2}
+
\\
&
\ \ \ \ \
-\tfrac{51}{4}xF_{0,3}+\tfrac{9}{4}yF_{0,3}+\tfrac{7}{32}xG_{1,3}+\tfrac{7}{8}yG_{2,1}+\tfrac{1}{128}vG_{2,1}^3+\tfrac{1}{4}yF_{4,0}+\tfrac{1}{16}yF_{1,3}+\tfrac{1}{2}vG_{2,2}
+
\\
&
\ \ \ \ \
-\tfrac{1}{4}vG_{1,3}+\tfrac{3}{128}uG_{2,1}^3+\tfrac{45}{16}uF_{0,3}^3-12vF_{0,3}^2-\tfrac{1}{32}yG_{1,3}-\tfrac{1}{8}yG_{4,0}-\tfrac{33}{2}uF_{0,3}^2-\tfrac{1}{32}vG_{2,1}G_{2,2}
+
\\
&
\ \ \ \ \
-\tfrac{3}{16}vF_{0,3}G_{3,1}-\tfrac{7}{16}xF_{1,3}+\tfrac{9}{8}vF_{0,3}^3+\tfrac{3}{32}yG_{3,1}-\tfrac{1}{64}yG_{2,1}^2-\tfrac{9}{16}yF_{0,3}^2+\tfrac{5}{8}uF_{0,3}G_{4,0}+\tfrac{3}{4}vG_{3,1}
+
\\
&
\ \ \ \ \
-\tfrac{9}{16}vG_{2,1}^2+\tfrac{3}{32}vG_{2,1}F_{3,1}-\tfrac{1}{8}vF_{0,3}G_{2,2}+\tfrac{21}{32}vF_{0,3}^2G_{2,1}-\tfrac{5}{16}uF_{0,3}F_{1,3}+\tfrac{3}{16}uG_{2,1}F_{2,2}
+
\\
&
\ \ \ \ \
-\tfrac{37}{8}xG_{2,1}
-\tfrac{9}{64}uG_{2,1}G_{3,1}+\tfrac{1}{4}vF_{0,3}G_{4,0}-\tfrac{1}{32}vG_{2,1}F_{1,3}-3uF_{4,0}+\tfrac{1}{2}vF_{1,3}+\tfrac{1}{4}vF_{0,3}F_{2,2}
+
\\
&
\ \ \ \ \
-\tfrac{9}{8}uG_{2,1}^2-\tfrac{1}{8}vG_{2,1}F_{4,0}-\tfrac{3}{4}uF_{1,3}+\tfrac{7}{8}xG_{4,0}-vF_{2,2}+\tfrac{9}{4}uF_{3,1}-\tfrac{13}{2}uF_{0,3}G_{2,1}
+\tfrac{9}{16}xF_{0,3}G_{2,1}
+
\\
&
\ \ \ \ \
+\tfrac{1}{16}vF_{0,3}G_{1,3}+2vF_{4,0}-\tfrac{3}{16}yF_{0,3}G_{2,1}-\tfrac{3}{32}uG_{2,1}F_{1,3}+\tfrac{3}{8}xF_{2,2}-\tfrac{5}{4}uF_{0,3}F_{4,0}+\tfrac{23}{64}uF_{0,3}G_{2,1}^2
+
\\
&
\ \ \ \ \
+\tfrac{3}{2}uG_{4,0}+16x-\tfrac{3}{64}vG_{2,1}G_{3,1}-\tfrac{5}{16}uF_{0,3}G_{2,2}-\tfrac{1}{2}vF_{0,3}F_{4,0}+\tfrac{1}{8}vF_{0,3}G_{2,1}^2
+
\\
&
\ \ \ \ \
-\tfrac{15}{32}uF_{0,3}G_{3,1}-48v
\Big)\,\partial_x
-
\endaligned
\]
\[
\aligned
&
-
\Big(
-\tfrac{3}{8}uG_{2,1}F_{4,0}+\tfrac{3}{64}uG_{2,1}G_{1,3}+\tfrac{1}{16}vG_{2,1}F_{2,2}+\tfrac{15}{16}uF_{0,3}F_{3,1}+\tfrac{3}{4}vF_{3,1}-\tfrac{3}{32}xG_{3,1}+
\\
&
\ \ \ \ \
+\tfrac{1}{64}xG_{2,1}^2+\tfrac{1}{16}vG_{2,1}G_{4,0}+\tfrac{1}{2}vG_{4,0}+\tfrac{9}{16}xF_{0,3}^2+\tfrac{15}{8}uG_{3,1}-\tfrac{9}{16}yF_{3,1}+\tfrac{3}{16}yG_{2,2}+\tfrac{5}{4}uG_{2,2}
+
\\
&
\ \ \ \ \
-\tfrac{3}{32}uG_{2,1}G_{2,2}-\tfrac{1}{8}vF_{0,3}F_{1,3}+\tfrac{57}{32}uF_{0,3}^2G_{2,1}+\tfrac{3}{8}vF_{0,3}F_{3,1}+\tfrac{5}{8}uF_{0,3}F_{2,2}+\tfrac{9}{32}uG_{2,1}F_{3,1}
+
\\
&
\ \ \ \ \
-\tfrac{7}{2}vG_{2,1}+69uF_{0,3}+3vF_{0,3}+\tfrac{55}{2}uG_{2,1}+\tfrac{3}{16}uG_{2,1}G_{4,0}+\tfrac{3}{16}xF_{3,1}-\tfrac{1}{16}xG_{2,2}-\tfrac{25}{8}vF_{0,3}G_{2,1}
+
\\
&
\ \ \ \ \
+\tfrac{5}{32}uF_{0,3}G_{1,3}-\tfrac{5}{8}uG_{1,3}+\tfrac{1}{64}vG_{2,1}G_{1,3}-\tfrac{3}{8}yF_{2,2}-\tfrac{1}{4}xF_{4,0}-\tfrac{5}{2}uF_{2,2}-\tfrac{21}{4}xF_{0,3}
+
\\
&
\ \ \ \ \
+\tfrac{51}{4}yF_{0,3}+\tfrac{1}{32}xG_{1,3}+\tfrac{33}{8}yG_{2,1}+\tfrac{1}{128}vG_{2,1}^3+\tfrac{7}{4}yF_{4,0}+\tfrac{7}{16}yF_{1,3}-\tfrac{1}{4}vG_{2,2}+\tfrac{5}{8}vG_{1,3}
+
\\
&
\ \ \ \ \
+\tfrac{3}{128}uG_{2,1}^3+\tfrac{45}{16}uF_{0,3}^3-\tfrac{21}{4}vF_{0,3}^2-\tfrac{7}{32}yG_{1,3}-\tfrac{7}{8}yG_{4,0}-\tfrac{51}{2}uF_{0,3}^2-\tfrac{1}{32}vG_{2,1}G_{2,2}
+
\\
&
\ \ \ \ \
-\tfrac{3}{16}vF_{0,3}G_{3,1}-\tfrac{1}{16}xF_{1,3}+\tfrac{9}{8}vF_{0,3}^3+\tfrac{9}{32}yG_{3,1}-\tfrac{3}{64}yG_{2,1}^2-\tfrac{27}{16}yF_{0,3}^2+\tfrac{5}{8}uF_{0,3}G_{4,0}
+
\\
&
\ \ \ \ \
-\tfrac{3}{8}vG_{3,1}-\tfrac{3}{8}vG_{2,1}^2+\tfrac{3}{32}vG_{2,1}F_{3,1}-\tfrac{1}{8}vF_{0,3}G_{2,2}+\tfrac{21}{32}vF_{0,3}^2G_{2,1}-\tfrac{5}{16}uF_{0,3}F_{1,3}
+
\\
&
\ \ \ \ \
+\tfrac{3}{16}uG_{2,1}F_{2,2}-\tfrac{7}{8}xG_{2,1}-\tfrac{9}{64}uG_{2,1}G_{3,1}+\tfrac{1}{4}vF_{0,3}G_{4,0}-\tfrac{1}{32}vG_{2,1}F_{1,3}+5uF_{4,0}
-\tfrac{1}{4}vF_{1,3}
+
\\
&
\ \ \ \ \
+\tfrac{1}{4}vF_{0,3}F_{2,2}-\tfrac{13}{8}uG_{2,1}^2-\tfrac{1}{8}vG_{2,1}F_{4,0}+\tfrac{5}{4}uF_{1,3}+\tfrac{1}{8}xG_{4,0}+\tfrac{1}{2}vF_{2,2}-\tfrac{3}{4}uF_{3,1}
+
\\
&
\ \ \ \ \
-14uF_{0,3}G_{2,1}+\tfrac{3}{16}xF_{0,3}G_{2,1}+\tfrac{1}{16}vF_{0,3}G_{1,3}-vF_{4,0}-\tfrac{9}{16}yF_{0,3}G_{2,1}-\tfrac{3}{32}uG_{2,1}F_{1,3}
+
\\
&
\ \ \ \ \
+\tfrac{1}{8}xF_{2,2}-\tfrac{5}{4}uF_{0,3}F_{4,0}+\tfrac{23}{64}uF_{0,3}G_{2,1}^2-48u-\tfrac{9}{2}uG_{4,0}-20y+12x-\tfrac{3}{64}vG_{2,1}G_{3,1}
+
\\
&
\ \ \ \ \
-\tfrac{5}{16}uF_{0,3}G_{2,2}-\tfrac{1}{2}vF_{0,3}F_{4,0}+\tfrac{1}{8}vF_{0,3}G_{2,1}^2-\tfrac{15}{32}uF_{0,3}G_{3,1}
\Big)\partial_y
+
\endaligned
\]
\[
\aligned
&
+
\Big(
-\tfrac{7}{2}uF_{4,0}-\tfrac{7}{8}uF_{1,3}-\tfrac{35}{4}uG_{2,1}-\tfrac{9}{16}uG_{3,1}+\tfrac{3}{32}uG_{2,1}^2-\tfrac{1}{8}vF_{2,2}-\tfrac{51}{2}uF_{0,3}
+
\\
&
\ \ \ \ \
+\tfrac{9}{4}vF_{0,3}+\tfrac{7}{8}vG_{2,1}+\tfrac{3}{32}vG_{3,1}-\tfrac{1}{64}vG_{2,1}^2-\tfrac{9}{16}vF_{0,3}^2-\tfrac{3}{16}vF_{3,1}+\tfrac{1}{16}vG_{2,2}-\tfrac{1}{32}vG_{1,3}
+
\\
&
\ \ \ \ \
-\tfrac{1}{8}vG_{4,0}-\tfrac{3}{16}vF_{0,3}G_{2,1}+\tfrac{9}{8}uF_{0,3}G_{2,1}+\tfrac{7}{16}uG_{1,3}+\tfrac{7}{4}uG_{4,0}+\tfrac{27}{8}uF_{0,3}^2+\tfrac{9}{8}uF_{3,1}
+
\\
&
\ \ \ \ \
-\tfrac{3}{8}uG_{2,2}+36u+y+\tfrac{1}{4}vF_{4,0}+\tfrac{1}{16}vF_{1,3}+\tfrac{3}{4}uF_{2,2}
\Big)\partial_u
-
\endaligned
\]
\[
\aligned
&
-
\Big(
-uF_{4,0}-\tfrac{1}{4}uF_{1,3}-\tfrac{11}{2}uG_{2,1}-\tfrac{3}{8}uG_{3,1}+\tfrac{1}{16}uG_{2,1}^2-\tfrac{3}{4}vF_{2,2}-15uF_{0,3}+\tfrac{51}{2}vF_{0,3}
+
\\
&
\ \ \ \ \
+\tfrac{37}{4}vG_{2,1}+\tfrac{9}{16}vG_{3,1}-\tfrac{3}{32}vG_{2,1}^2-\tfrac{27}{8}vF_{0,3}^2
-\tfrac{9}{8}vF_{3,1}+\tfrac{3}{8}vG_{2,2}-\tfrac{7}{16}vG_{1,3}-\tfrac{7}{4}vG_{4,0}
+
\\
&
\ \ \ \ \
-\tfrac{9}{8}vF_{0,3}G_{2,1}+\tfrac{3}{4}uF_{0,3}G_{2,1}+\tfrac{1}{8}uG_{1,3}+\tfrac{1}{2}uG_{4,0}+\tfrac{9}{4}uF_{0,3}^2+\tfrac{3}{4}uF_{3,1}-\tfrac{1}{4}uG_{2,2}+24u
+
\\
&
\ \ \ \ \
-2x-32v+\tfrac{7}{2}vF_{4,0}+\tfrac{7}{8}vF_{1,3}+\tfrac{1}{2}uF_{2,2}
\Big)\partial_v
+
\endaligned
\]
\[
\aligned
\\
e_2
&
\,:=\,
\Big(
\tfrac{7}{2}uG_{0,4}-3uF_{0,4}-\tfrac{1}{8}yG_{0,4}+\tfrac{1}{4}yF_{0,4}+\tfrac{45}{16}uF_{0,3}^3-\tfrac{3}{16}yF_{1,3}+\tfrac{3}{128}uG_{2,1}^3+\tfrac{3}{32}yG_{1,3}
+
\\
&
\ \ \ \ \
+\tfrac{1}{16}yF_{3,1}+\tfrac{1}{16}yG_{2,2}-\tfrac{1}{64}yG_{2,1}^2-\tfrac{9}{16}yF_{0,3}^2-\tfrac{3}{4}uF_{1,3}-\tfrac{1}{32}yG_{3,1}-\tfrac{3}{32}uG_{2,1}F_{3,1}
+
\\
&
\ \ \ \ \
-\tfrac{1}{8}vF_{0,3}F_{3,1}+\tfrac{57}{32}uF_{0,3}^2G_{2,1}+\tfrac{5}{8}uF_{0,3}F_{2,2}+\tfrac{1}{8}vF_{0,3}G_{2,1}^2+\tfrac{5}{32}uF_{0,3}G_{3,1}-\tfrac{5}{16}uF_{0,3}F_{3,1}
+
\\
&
\ \ \ \ \
+\tfrac{1}{16}vG_{2,1}F_{2,2}-\tfrac{9}{64}uG_{2,1}G_{1,3}+2vF_{0,4}+\tfrac{7}{8}xG_{0,4}-\tfrac{7}{4}xF_{0,4}-vG_{0,4}-\tfrac{15}{32}uF_{0,3}G_{1,3}
+
\\
&
\ \ \ \ \
-\tfrac{3}{16}yF_{0,3}G_{2,1}+\tfrac{9}{32}uG_{2,1}F_{1,3}+\tfrac{9}{16}xF_{0,3}G_{2,1}-\tfrac{1}{32}vG_{2,1}F_{3,1}+\tfrac{3}{32}vG_{2,1}F_{1,3}+\tfrac{23}{64}uF_{0,3}G_{2,1}^2
+
\\
&
\ \ \ \ \
+\tfrac{15}{2}uG_{2,1}+\tfrac{3}{8}uG_{3,1}+\tfrac{3}{4}uG_{2,1}^2-vF_{2,2}-\tfrac{9}{32}xG_{1,3}-\tfrac{3}{8}yG_{2,1}+\tfrac{1}{128}vG_{2,1}^3+\tfrac{21}{32}vF_{0,3}^2G_{2,1}
+
\\
&
\ \ \ \ \
-\tfrac{1}{8}vF_{0,3}G_{2,2}+\tfrac{7}{32}xG_{3,1}+\tfrac{3}{64}xG_{2,1}^2+\tfrac{3}{16}uG_{2,1}F_{2,2}+\tfrac{3}{64}uG_{2,1}G_{3,1}+21uF_{0,3}+6vF_{0,3}
+
\\
&
\ \ \ \ \
-6vG_{2,1}-\tfrac{3}{4}vG_{3,1}+\tfrac{1}{16}vG_{2,1}^2-3vF_{0,3}^2+\tfrac{1}{2}vF_{3,1}+\tfrac{1}{2}vG_{2,2}+\tfrac{3}{4}vG_{1,3}-\tfrac{1}{8}vF_{0,3}G_{2,1}
+
\\
&
\ \ \ \ \
+\tfrac{11}{2}uF_{0,3}G_{2,1}+\tfrac{15}{16}uF_{0,3}F_{1,3}+\tfrac{15}{4}xF_{0,3}-\tfrac{9}{4}yF_{0,3}+\tfrac{3}{8}vF_{0,3}F_{1,3}-\tfrac{3}{32}uG_{2,1}G_{2,2}
+
\\
&
\ \ \ \ \
-\tfrac{9}{8}uG_{1,3}+\tfrac{1}{4}vF_{0,3}F_{2,2}+\tfrac{21}{8}xG_{2,1}-\tfrac{3}{16}xG_{2,2}+\tfrac{27}{16}xF_{0,3}^2+\tfrac{1}{4}vF_{0,3}G_{0,4}-\tfrac{1}{2}vF_{0,3}F_{0,4}
+
\\
&
\ \ \ \ \
+\tfrac{1}{16}vF_{0,3}G_{3,1}+\tfrac{1}{16}vG_{2,1}G_{0,4}
-\tfrac{5}{4}uF_{0,3}F_{0,4}-\tfrac{3}{8}uG_{2,1}F_{0,4}+\tfrac{5}{8}uF_{0,3}G_{0,4}-\tfrac{1}{8}vG_{2,1}F_{0,4}
+
\\
&
\ \ \ \ \
+\tfrac{1}{64}vG_{2,1}G_{3,1}-\tfrac{5}{16}uF_{0,3}G_{2,2}+\tfrac{3}{16}uG_{2,1}G_{0,4}-\tfrac{3}{16}vF_{0,3}G_{1,3}-\tfrac{7}{16}xF_{3,1}-\tfrac{3}{64}vG_{2,1}G_{1,3}
+
\\
&
\ \ \ \ \
+\tfrac{9}{8}vF_{0,3}^3+6uF_{0,3}^2-\tfrac{3}{4}uF_{3,1}-\tfrac{3}{4}uG_{2,2}-\tfrac{1}{32}vG_{2,1}G_{2,2}-\tfrac{1}{8}yF_{2,2}+\tfrac{9}{16}xF_{1,3}+24u
+
\\
&
\ \ \ \ \
-2y+14x-48v+\tfrac{3}{8}xF_{2,2}-\tfrac{3}{2}vF_{1,3}+\tfrac{3}{2}uF_{2,2}
\Big)\partial_x
-
\endaligned
\]
\[
\aligned
&
-
\Big(
-1-\tfrac{5}{2}uG_{0,4}+5uF_{0,4}-\tfrac{7}{8}yG_{0,4}+\tfrac{7}{4}yF_{0,4}+\tfrac{45}{16}uF_{0,3}^3-\tfrac{9}{16}yF_{1,3}+\tfrac{3}{128}uG_{2,1}^3
+
\\
&
\ \ \ \ \
+\tfrac{9}{32}yG_{1,3}+\tfrac{7}{16}yF_{3,1}+\tfrac{3}{16}yG_{2,2}-\tfrac{3}{64}yG_{2,1}^2-\tfrac{27}{16}yF_{0,3}^2-\tfrac{15}{4}uF_{1,3}-\tfrac{7}{32}yG_{3,1}
+
\\
&
\ \ \ \ \
-\tfrac{3}{32}uG_{2,1}F_{3,1}-\tfrac{1}{8}vF_{0,3}F_{3,1}+\tfrac{57}{32}uF_{0,3}^2G_{2,1}+\tfrac{5}{8}uF_{0,3}F_{2,2}+\tfrac{1}{8}vF_{0,3}G_{2,1}^2+\tfrac{5}{32}uF_{0,3}G_{3,1}
+
\\
&
\ \ \ \ \
-\tfrac{5}{16}uF_{0,3}F_{3,1}+\tfrac{1}{16}vG_{2,1}F_{2,2}-\tfrac{9}{64}uG_{2,1}G_{1,3}-vF_{0,4}+\tfrac{1}{8}xG_{0,4}-\tfrac{1}{4}xF_{0,4}+\tfrac{5}{2}vG_{0,4}
+
\\
&
\ \ \ \ \
-\tfrac{15}{32}uF_{0,3}G_{1,3}-\tfrac{9}{16}yF_{0,3}G_{2,1}+\tfrac{9}{32}uG_{2,1}F_{1,3}+\tfrac{3}{16}xF_{0,3}G_{2,1}-\tfrac{1}{32}vG_{2,1}F_{3,1}+\tfrac{3}{32}vG_{2,1}F_{1,3}
+
\\
&
\ \ \ \ \
+\tfrac{23}{64}uF_{0,3}G_{2,1}^2-\tfrac{17}{2}uG_{2,1}-\tfrac{9}{8}uG_{3,1}+\tfrac{1}{4}uG_{2,1}^2+\tfrac{1}{2}vF_{2,2}-\tfrac{3}{32}xG_{1,3}-\tfrac{25}{8}yG_{2,1}+\tfrac{1}{128}vG_{2,1}^3
+
\\
&
\ \ \ \ \
+\tfrac{21}{32}vF_{0,3}^2G_{2,1}-\tfrac{1}{8}vF_{0,3}G_{2,2}+\tfrac{1}{32}xG_{3,1}+\tfrac{1}{64}xG_{2,1}^2+\tfrac{3}{16}uG_{2,1}F_{2,2}+\tfrac{3}{64}uG_{2,1}G_{3,1}
+
\\
&
\ \ \ \ \
-15uF_{0,3}+21vF_{0,3}+\tfrac{5}{2}vG_{2,1}+\tfrac{1}{8}vG_{3,1}+\tfrac{1}{4}vG_{2,1}^2+\tfrac{15}{4}vF_{0,3}^2-\tfrac{1}{4}vF_{3,1}-\tfrac{1}{4}vG_{2,2}
+
\\
&
\ \ \ \ \
-\tfrac{3}{8}vG_{1,3}+\tfrac{29}{8}vF_{0,3}G_{2,1}-2uF_{0,3}G_{2,1}+\tfrac{15}{16}uF_{0,3}F_{1,3}-\tfrac{3}{4}xF_{0,3}-\tfrac{15}{4}yF_{0,3}+\tfrac{3}{8}vF_{0,3}F_{1,3}
+
\\
&
\ \ \ \ \
-\tfrac{3}{32}uG_{2,1}G_{2,2}+\tfrac{15}{8}uG_{1,3}+\tfrac{1}{4}vF_{0,3}F_{2,2}+\tfrac{3}{8}xG_{2,1}-\tfrac{1}{16}xG_{2,2}+\tfrac{9}{16}xF_{0,3}^2+\tfrac{1}{4}vF_{0,3}G_{0,4}
+
\\
&
\ \ \ \ \
-\tfrac{1}{2}vF_{0,3}F_{0,4}+\tfrac{1}{16}vF_{0,3}G_{3,1}+\tfrac{1}{16}vG_{2,1}G_{0,4}-\tfrac{5}{4}uF_{0,3}F_{0,4}-\tfrac{3}{8}uG_{2,1}F_{0,4}+\tfrac{5}{8}uF_{0,3}G_{0,4}
+
\\
&
\ \ \ \ \
-\tfrac{1}{8}vG_{2,1}F_{0,4}+\tfrac{1}{64}vG_{2,1}G_{3,1}-\tfrac{5}{16}uF_{0,3}G_{2,2}+\tfrac{3}{16}uG_{2,1}G_{0,4}-\tfrac{3}{16}vF_{0,3}G_{1,3}-\tfrac{1}{16}xF_{3,1}
+
\\
&
\ \ \ \ \
-\tfrac{3}{64}vG_{2,1}G_{1,3}+\tfrac{9}{8}vF_{0,3}^3-3uF_{0,3}^2+\tfrac{5}{4}uF_{3,1}+\tfrac{5}{4}uG_{2,2}-\tfrac{1}{32}vG_{2,1}G_{2,2}-\tfrac{3}{8}yF_{2,2}+\tfrac{3}{16}xF_{1,3}
+
\\
&
\ \ \ \ \
-72u-14y+2x+8v+\tfrac{1}{8}xF_{2,2}+\tfrac{3}{4}vF_{1,3}-\tfrac{1}{2}uF_{2,2}
\Big)\partial_y
+
\endaligned
\]
\[
\aligned
&
+
\Big(
\tfrac{7}{4}uG_{0,4}-\tfrac{7}{2}uF_{0,4}+\tfrac{9}{8}uF_{1,3}+\tfrac{1}{4}vF_{0,4}-\tfrac{1}{8}vG_{0,4}+\tfrac{23}{4}uG_{2,1}+\tfrac{7}{16}uG_{3,1}+\tfrac{3}{32}uG_{2,1}^2
+
\\
&
\ \ \ \ \
-\tfrac{1}{8}vF_{2,2}+\tfrac{15}{2}uF_{0,3}+\tfrac{3}{4}vF_{0,3}-\tfrac{3}{8}vG_{2,1}-\tfrac{1}{32}vG_{3,1}-\tfrac{1}{64}vG_{2,1}^2-\tfrac{9}{16}vF_{0,3}^2+\tfrac{1}{16}vF_{3,1}
+
\\
&
\ \ \ \ \
+\tfrac{1}{16}vG_{2,2}+\tfrac{3}{32}vG_{1,3}-\tfrac{3}{16}vF_{0,3}G_{2,1}+\tfrac{9}{8}uF_{0,3}G_{2,1}-\tfrac{9}{16}uG_{1,3}+\tfrac{27}{8}uF_{0,3}^2-\tfrac{7}{8}uF_{3,1}
+
\\
&
\ \ \ \ \
-\tfrac{3}{8}uG_{2,2}+28u+x-2v-\tfrac{3}{16}vF_{1,3}+\tfrac{3}{4}uF_{2,2}
\Big)\partial_u
+
\endaligned
\]
\[
\aligned
&
\ \ \ \ \
-
\Big(
\tfrac{1}{2}uG_{0,4}-uF_{0,4}+\tfrac{3}{4}uF_{1,3}+\tfrac{7}{2}vF_{0,4}-\tfrac{7}{4}vG_{0,4}+\tfrac{7}{2}uG_{2,1}+\tfrac{1}{8}uG_{3,1}+\tfrac{1}{16}uG_{2,1}^2
+
\\
&
\ \ \ \ \
-\tfrac{3}{4}vF_{2,2}+3uF_{0,3}-\tfrac{15}{2}vF_{0,3}-\tfrac{25}{4}vG_{2,1}-\tfrac{7}{16}vG_{3,1}-\tfrac{3}{32}vG_{2,1}^2-\tfrac{27}{8}vF_{0,3}^2+\tfrac{7}{8}vF_{3,1}
+
\\
&
\ \ \ \ \
+\tfrac{3}{8}vG_{2,2}+\tfrac{9}{16}vG_{1,3}-\tfrac{9}{8}vF_{0,3}G_{2,1}+\tfrac{3}{4}uF_{0,3}G_{2,1}-\tfrac{3}{8}uG_{1,3}+\tfrac{9}{4}uF_{0,3}^2
+
\\
&
\ \ \ \ \
-\tfrac{1}{4}uF_{3,1}-\tfrac{1}{4}uG_{2,2}+24u-2y-28v-\tfrac{9}{8}vF_{1,3}+\tfrac{1}{2}uF_{2,2}
\Big)\partial_v.
\endaligned
\]

%%%%%%%%%%%%%%%%%%%%%%%%%%%%%%%%%%%%%%%%%%%%%%%%%%%%%%%%%%%%%%%%%%%%%%
\noindent
A Gr\"obner basis has a lot of generators to describe the
moduli space core algebraic variety in 
$\R^{12} \ni F_{0,3},G_{2,1},F_{0,4},G_{0,4},F_{4,0},G_{4,0},F_{1,3},G_{1,3},F_{3,1},G_{3,1},F_{2,2},G_{2,2}$.

%%%%%%%%%%%%%%%%%%%%%%%%%%%%%%%%%%%%%%%%%%%%%%%%%%%%%%%%%%%%%%%%%%%%%%
\[
\footnotesize
\def\arraystretch{1.25}
\begin{array}{c|cc}
{} & e_1 & e_2 
\\
\hline
e_1 &
0
&
\rotatebox[origin=c]{0}{
\begin{tabular}{p{7cm}}
$
(
14
+
\tfrac{3}{4}F_{0,3}G_{2,1}
+
\tfrac{1}{16}G_{2,1}^2
+
\tfrac{9}{4}F_{0,3}^2
+
\tfrac{1}{2}F_{1,3}
+
\tfrac{7}{4}G_{2,1}
-
\tfrac{1}{4}G_{2,2}
+
\tfrac{1}{2}F_{2,2}
+
\tfrac{3}{2}F_{0,3}
+
\tfrac{1}{8}G_{3,1}
-
\tfrac{1}{4}F_{3,1}
-
\tfrac{1}{4}G_{1,3}
+
\tfrac{7}{8}G_{0,4}
+
\tfrac{1}{8}G_{4,0}
-
\tfrac{1}{4}F_{4,0}
-
\tfrac{7}{4}F_{0,4}
)e_1
+
(
-
22
-
\tfrac{3}{4}F_{0,3}G_{2,1}
-
\tfrac{1}{16}G_{2,1}^2
-
\tfrac{9}{4}F_{0,3}^2
+
\tfrac{1}{4}F_{1,3}
+
\tfrac{15}{4}G_{2,1}
+
\tfrac{1}{4}G_{2,2}
-
\tfrac{1}{2}F_{2,2}
+
\tfrac{27}{2}F_{0,3}
+
\tfrac{1}{4}G_{3,1}
-
\tfrac{1}{2}F_{3,1}
-
\tfrac{1}{8}G_{1,3}
-
\tfrac{1}{8}G_{0,4}
-
\tfrac{7}{8}G_{4,0}
+
\tfrac{7}{4}F_{4,0}
+
\tfrac{1}{4}F_{0,4}
)e_2$
\end{tabular}}
\\
e_2 &
\rotatebox[origin=c]{0}{
	\begin{tabular}{p{7cm}}
	$
	-
	(
	14
	+
	\tfrac{3}{4}F_{0,3}G_{2,1}
	+
	\tfrac{1}{16}G_{2,1}^2
	+
	\tfrac{9}{4}F_{0,3}^2
	+
	\tfrac{1}{2}F_{1,3}
	+
	\tfrac{7}{4}G_{2,1}
	-
	\tfrac{1}{4}G_{2,2}
	+
	\tfrac{1}{2}F_{2,2}
	+
	\tfrac{3}{2}F_{0,3}
	+
	\tfrac{1}{8}G_{3,1}
	-
	\tfrac{1}{4}F_{3,1}
	-
	\tfrac{1}{4}G_{1,3}
	+
	\tfrac{7}{8}G_{0,4}
	+
	\tfrac{1}{8}G_{4,0}
	-
	\tfrac{1}{4}F_{4,0}
	-
	\tfrac{7}{4}F_{0,4}
	)e_1
	+
	(
	-
	22
	-
	\tfrac{3}{4}F_{0,3}G_{2,1}
	-
	\tfrac{1}{16}G_{2,1}^2
	-
	\tfrac{9}{4}F_{0,3}^2
	+
	\tfrac{1}{4}F_{1,3}
	+
	\tfrac{15}{4}G_{2,1}
	+
	\tfrac{1}{4}G_{2,2}
	-
	\tfrac{1}{2}F_{2,2}
	+
	\tfrac{27}{2}F_{0,3}
	+
	\tfrac{1}{4}G_{3,1}
	-
	\tfrac{1}{2}F_{3,1}
	-
	\tfrac{1}{8}G_{1,3}
	-
	\tfrac{1}{8}G_{0,4}
	-
	\tfrac{7}{8}G_{4,0}
	+
	\tfrac{7}{4}F_{4,0}
	+
	\tfrac{1}{4}F_{0,4}
	)e_2$
	\end{tabular}}
& 0
\end{array}
\]

%%%%%%%%%%%%%%%%%%%%%%%%%%%%%%%%%%%%%%%%%%%%%%%%%%%%%%%%%%%%%%%%%%%%%%
%%%%%%%%%%%%%%%%%%%%%%%%%%%%%%%%%%%%%%%%%%%%%%%%%%%%%%%%%%%%%%%%%%%%%%
%%%%%%%%%%%%%%%%%%%%%%%%%%%%%%%%%%%%%%%%%%%%%%%%%%%%%%%%%%%%%%%%%%%%%%

\[
\text{\bf Model 2f3b}
\ \ \ \ \
\left\{
\aligned
u
&
\,=\,
\text{A lot of terms at order 5}
,
\\
v
&
\,=\,
\text{A lot of terms at order 5},
\endaligned\right.
\]
%%%%%%%%%%%%%%%%%%%%%%%%%%%%%%%%%%%%%%%%%%%%%%%%%%%%%%%%%%%%%%%%%%%%%%
\[
\def\arraystretch{1.25}
\begin{array}{ll}
e_1
\,:=\,
\text{A lot of terms}
, &
e_2
\,:=\,
\text{A lot of terms}.
\end{array}
\]

%%%%%%%%%%%%%%%%%%%%%%%%%%%%%%%%%%%%%%%%%%%%%%%%%%%%%%%%%%%%%%%%%%%%%%
\noindent
A Gr\"obner basis has many generators to describe the
moduli space core algebraic variety in 
$\R^{12} \ni F_{0,3},G_{2,1},G_{0,4},F_{0,4},
F_{4,0},G_{4,0},F_{1,3},G_{1,3},F_{3,1},G_{3,1},F_{2,2},G_{2,2}$.
%%%%%%%%%%%%%%%%%%%%%%%%%%%%%%%%%%%%%%%%%%%%%%%%%%%%%%%%%%%%%%%%%%%%%%
\[
\footnotesize
\def\arraystretch{1.25}
\begin{array}{c|cc}
{} & e_1 & e_2 
\\
\hline
e_1 &
0
&
\rotatebox[origin=c]{0}{
\begin{tabular}{p{7cm}}
$
(\tfrac{9}{8}F_{0,3}^2
+
\tfrac{1}{32}G_{2,1}^2
+
\tfrac{3}{8}F_{0,3}G_{2,1}
+
G_{0,4}
+
\tfrac{1}{4}G_{4,0}
+
\tfrac{1}{4}F_{1,3}
-
\tfrac{1}{8}G_{2,2}
+
\tfrac{1}{4}F_{2,2}
-
9F_{0,3}
+
\tfrac{1}{4}G_{3,1}
-
\tfrac{1}{2}F_{3,1}
-
\tfrac{1}{8}G_{1,3}
-
\tfrac{1}{2}F_{4,0}
-
2F_{0,4}
)e_1
+
(
-
8
-
\tfrac{9}{8}F_{0,3}^2
-
\tfrac{1}{32}G_{2,1}^2
-
\tfrac{3}{8}F_{0,3}G_{2,1}
-
\tfrac{1}{4}G_{0,4}
-
G_{4,0}
+
\tfrac{1}{2}F_{1,3}
+
G_{2,1}
+
\tfrac{1}{8}G_{2,2}
-
\tfrac{1}{4}F_{2,2}
-
3F_{0,3}
+
\tfrac{1}{8}G_{3,1}
-
\tfrac{1}{4}F_{3,1}
-
\tfrac{1}{4}G_{1,3}
+
2F_{4,0}
+
\tfrac{1}{2}F_{0,4})e_2
$
\end{tabular}}
\\
e_2 &
\rotatebox[origin=c]{0}{
\begin{tabular}{p{7cm}}
$
-
(
\tfrac{9}{8}F_{0,3}^2
+
\tfrac{1}{32}G_{2,1}^2
+
\tfrac{3}{8}F_{0,3}G_{2,1}
+
G_{0,4}
+
\tfrac{1}{4}G_{4,0}
+
\tfrac{1}{4}F_{1,3}
-
\tfrac{1}{8}G_{2,2}
+
\tfrac{1}{4}F_{2,2}
-
9F_{0,3}
+
\tfrac{1}{4}G_{3,1}
-
\tfrac{1}{2}F_{3,1}
-
\tfrac{1}{8}G_{1,3}
-
\tfrac{1}{2}F_{4,0}
-
2F_{0,4}
)e_1
+
(
-
8
-
\tfrac{9}{8}F_{0,3}^2
-
\tfrac{1}{32}G_{2,1}^2
-
\tfrac{3}{8}F_{0,3}G_{2,1}
-
\tfrac{1}{4}G_{0,4}
-
G_{4,0}
+
\tfrac{1}{2}F_{1,3}
+
G_{2,1}
+
\tfrac{1}{8}G_{2,2}
-
\tfrac{1}{4}F_{2,2}
-
3F_{0,3}
+
\tfrac{1}{8}G_{3,1}
-
\tfrac{1}{4}F_{3,1}
-
\tfrac{1}{4}G_{1,3}
+
2F_{4,0}
+
\tfrac{1}{2}F_{0,4})e_2
$
\end{tabular}}
& 0
\end{array}
\]

%%%%%%%%%%%%%%%%%%%%%%%%%%%%%%%%%%%%%%%%%%%%%%%%%%%%%%%%%%%%%%%%%%%%%%
%%%%%%%%%%%%%%%%%%%%%%%%%%%%%%%%%%%%%%%%%%%%%%%%%%%%%%%%%%%%%%%%%%%%%%
%%%%%%%%%%%%%%%%%%%%%%%%%%%%%%%%%%%%%%%%%%%%%%%%%%%%%%%%%%%%%%%%%%%%%%

\[
\text{\bf Model 2f3c}
\ \ \ \ \
\left\{
\aligned
u
&
\,=\,
\text{A lot of terms at order 5}
,
\\
v
&
\,=\,
\text{A lot of terms at order 5}
.
\endaligned\right.
\]
%%%%%%%%%%%%%%%%%%%%%%%%%%%%%%%%%%%%%%%%%%%%%%%%%%%%%%%%%%%%%%%%%%%%%%
\[
\def\arraystretch{1.25}
\begin{array}{ll}
e_1
\,:=\,
\text{A lot of terms}
, &
e_2
\,:=\,
\text{A lot of terms}.
\end{array}
\]

%%%%%%%%%%%%%%%%%%%%%%%%%%%%%%%%%%%%%%%%%%%%%%%%%%%%%%%%%%%%%%%%%%%%%%
\noindent
A Gr\"obner basis has 92 generators to describe the
moduli space core algebraic variety in 
$\R^9 \ni _{0,3},F_{0,4},G_{0,4},
G_{4,0},F_{1,3},G_{1,3},G_{3,1},F_{2,2},G_{2,2}$.
 
%%%%%%%%%%%%%%%%%%%%%%%%%%%%%%%%%%%%%%%%%%%%%%%%%%%%%%%%%%%%%%%%%%%%%%
\[
\footnotesize
\def\arraystretch{1.25}
\begin{array}{c|cc}
{} & e_1 & e_2 
\\
\hline
e_1 &
0
&
\rotatebox[origin=c]{0}{
\begin{tabular}{p{7cm}}
$
(
-
\tfrac{577}{36}
+
\tfrac{11}{432}F_{0,3}F_{2,2}
-
\tfrac{5}{24}F_{0,3}F_{0,4}
+
\tfrac{5}{192}F_{0,3}G_{3,1}
-
\tfrac{11}{576}F_{0,3}G_{1,3}
-
\tfrac{1}{72}F_{0,3}F_{1,3}
+
\tfrac{1}{216}F_{0,3}G_{2,2}
-
\tfrac{23}{48}F_{0,3}^3
+
\tfrac{175}{144}F_{0,3}^2
+
\tfrac{233}{54}F_{0,3}
+
\tfrac{1}{48}F_{1,3}
-
\tfrac{65}{144}F_{2,2}
+
\tfrac{73}{192}G_{1,3}
+
\tfrac{7}{144}G_{2,2}
-
\tfrac{71}{192}G_{3,1}
+
\tfrac{21}{8}F_{0,4}
-
\tfrac{1}{12}G_{0,4}
-
\tfrac{1}{4}G_{4,0}
)e_1
+
(
\tfrac{1097}{36}
-
\tfrac{187}{432}F_{0,3}F_{2,2}
+
\tfrac{85}{24}F_{0,3}F_{0,4}
-
\tfrac{85}{192}F_{0,3}G_{3,1}
+
\tfrac{187}{576}F_{0,3}G_{1,3}
+
\tfrac{17}{72}F_{0,3}F_{1,3}
-
\tfrac{17}{216}F_{0,3}G_{2,2}
+
\tfrac{391}{48}F_{0,3}^3
+
\tfrac{2209}{144}F_{0,3}^2
-
\tfrac{1369}{54}F_{0,3}
-
\tfrac{17}{48}F_{1,3}
+
\tfrac{97}{144}F_{2,2}
-
\tfrac{89}{192}G_{1,3}
-
\tfrac{47}{144}G_{2,2}
+
\tfrac{55}{192}G_{3,1}
-
\tfrac{21}{8}F_{0,4}
+
\tfrac{17}{12}G_{0,4}
+
\tfrac{5}{4}G_{4,0}
)e_2
$
\end{tabular}}
\\
e_2 &
\rotatebox[origin=c]{0}{
\begin{tabular}{p{7cm}}
$
-
(
-
\tfrac{577}{36}
+
\tfrac{11}{432}F_{0,3}F_{2,2}
-
\tfrac{5}{24}F_{0,3}F_{0,4}
+
\tfrac{5}{192}F_{0,3}G_{3,1}
-
\tfrac{11}{576}F_{0,3}G_{1,3}
-
\tfrac{1}{72}F_{0,3}F_{1,3}
+
\tfrac{1}{216}F_{0,3}G_{2,2}
-
\tfrac{23}{48}F_{0,3}^3
+
\tfrac{175}{144}F_{0,3}^2
+
\tfrac{233}{54}F_{0,3}
+
\tfrac{1}{48}F_{1,3}
-
\tfrac{65}{144}F_{2,2}
+
\tfrac{73}{192}G_{1,3}
+
\tfrac{7}{144}G_{2,2}
-
\tfrac{71}{192}G_{3,1}
+
\tfrac{21}{8}F_{0,4}
-
\tfrac{1}{12}G_{0,4}
-
\tfrac{1}{4}G_{4,0}
)e_1
+
(
\tfrac{1097}{36}
-
\tfrac{187}{432}F_{0,3}F_{2,2}
+
\tfrac{85}{24}F_{0,3}F_{0,4}
-
\tfrac{85}{192}F_{0,3}G_{3,1}
+
\tfrac{187}{576}F_{0,3}G_{1,3}
+
\tfrac{17}{72}F_{0,3}F_{1,3}
-
\tfrac{17}{216}F_{0,3}G_{2,2}
+
\tfrac{391}{48}F_{0,3}^3
+
\tfrac{2209}{144}F_{0,3}^2
-
\tfrac{1369}{54}F_{0,3}
-
\tfrac{17}{48}F_{1,3}
+
\tfrac{97}{144}F_{2,2}
-
\tfrac{89}{192}G_{1,3}
-
\tfrac{47}{144}G_{2,2}
+
\tfrac{55}{192}G_{3,1}
-
\tfrac{21}{8}F_{0,4}
+
\tfrac{17}{12}G_{0,4}
+
\tfrac{5}{4}G_{4,0}
)e_2
$
\end{tabular}}
& 0
\end{array}
\]

%%%%%%%%%%%%%%%%%%%%%%%%%%%%%%%%%%%%%%%%%%%%%%%%%%%%%%%%%%%%%%%%%%%%%%
%%%%%%%%%%%%%%%%%%%%%%%%%%%%%%%%%%%%%%%%%%%%%%%%%%%%%%%%%%%%%%%%%%%%%%
%%%%%%%%%%%%%%%%%%%%%%%%%%%%%%%%%%%%%%%%%%%%%%%%%%%%%%%%%%%%%%%%%%%%%%

\[
\text{\bf Model 2f3d}
\ \ \ \ \
\left\{
\aligned
u
&
\,=\,
xy+2x^3+
(
-2+\tfrac{1}{2}F_{2,2}-\tfrac{1}{4}G_{3,1}-\tfrac{3}{4}G_{4,0}+\tfrac{1}{2}F_{3,1}-\tfrac{1}{4}G_{0,4}
)y^4
+
(
16
+
\\
&
\ \ \ \ \
-F_{2,2}+2F_{4,0}+\tfrac{1}{2}G_{3,1}+\tfrac{3}{2}G_{4,0}+\tfrac{1}{2}G_{0,4}
)xy^3
+
F_{2,2}x^2y^2+F_{3,1}x^3y
+
\\
&
\ \ \ \ \
+F_{4,0}x^4
+
\cdots
,
\\
v
&
\,=\,
y^2+x^2+8xy^2-4x^2y+G_{0,4}y^4
+
(
G_{3,1}+2G_{4,0}-2G_{0,4}
)xy^3
+
\\
&
\ \ \ \ \
+
(
64+F_{2,2}+6F_{4,0}-\tfrac{3}{2}G_{3,1}-\tfrac{3}{2}G_{4,0}+\tfrac{3}{2}G_{0,4}+3F_{3,1}
)x^2y^2
+
\\
&
\ \ \ \ \
+
G_{3,1}x^3y+G_{4,0}x^4
+
\cdots
.
\endaligned\right.
\]
%%%%%%%%%%%%%%%%%%%%%%%%%%%%%%%%%%%%%%%%%%%%%%%%%%%%%%%%%%%%%%%%%%%%%%
\[
\def\arraystretch{1.25}
\begin{array}{ll}
e_1
&
\,:=\,
(
1+32u-16v+\tfrac{1}{3}xF_{2,2}-\tfrac{8}{3}xF_{4,0}-\tfrac{2}{3}xG_{0,4}+\tfrac{1}{12}xG_{3,1}+\tfrac{2}{3}xG_{4,0}+\tfrac{1}{6}yF_{2,2}-\tfrac{1}{3}yF_{4,0}
+
\\
&
\ \ \ \ \
-\tfrac{1}{12}yG_{0,4}-\tfrac{1}{12}yG_{3,1}+\tfrac{1}{12}yG_{4,0}+vF_{2,2}-2vF_{4,0}-\tfrac{1}{2}vG_{0,4}-\tfrac{1}{2}vG_{3,1}-\tfrac{3}{2}vG_{4,0}
+
\\
&
\ \ \ \ \
-4uF_{4,0}-2uG_{0,4}-\tfrac{1}{2}uG_{3,1}+2uG_{4,0}+\tfrac{8}{3}y+\tfrac{4}{3}x
)\partial_x
-
(
-\tfrac{4}{3}uG_{3,1}-\tfrac{2}{3}uG_{4,0}-\tfrac{1}{6}xF_{2,2}
+
\\
&
\ \ \ \ \
+\tfrac{1}{3}xF_{4,0}+\tfrac{1}{12}xG_{0,4}+\tfrac{1}{12}xG_{3,1}-\tfrac{1}{12}xG_{4,0}-\tfrac{1}{3}yF_{2,2}+\tfrac{8}{3}yF_{4,0}+\tfrac{2}{3}yG_{0,4}-\tfrac{1}{12}yG_{3,1}
+
\\
&
\ \ \ \ \
-\tfrac{2}{3}yG_{4,0}+\tfrac{2}{3}vF_{2,2}-\tfrac{4}{3}vF_{4,0}-\tfrac{4}{3}vG_{0,4}+\tfrac{1}{6}vG_{3,1}+\tfrac{4}{3}vG_{4,0}-\tfrac{16}{3}uF_{4,0}-\tfrac{4}{3}uG_{0,4}
+
\\
&
\ \ \ \ \
+\tfrac{8}{3}uF_{2,2}+3uF_{3,1}+\tfrac{32}{3}u+\tfrac{8}{3}y+\tfrac{10}{3}x+\tfrac{32}{3}v
)\partial_y
+
(
-\tfrac{4}{3}u+\tfrac{8}{3}v+\tfrac{1}{6}vF_{2,2}-\tfrac{1}{3}vF_{4,0}
+
\\
&
\ \ \ \ \
-\tfrac{1}{12}vG_{0,4}-\tfrac{1}{12}vG_{3,1}+\tfrac{1}{12}vG_{4,0}+\tfrac{2}{3}uF_{2,2}-\tfrac{16}{3}uF_{4,0}-\tfrac{4}{3}uG_{0,4}+\tfrac{1}{6}uG_{3,1}+\tfrac{4}{3}uG_{4,0}
+
\\
&
\ \ \ \ \
+y
)\partial_u
+
(
-\tfrac{28}{3}u+\tfrac{8}{3}v+\tfrac{2}{3}vF_{2,2}-\tfrac{16}{3}vF_{4,0}-\tfrac{4}{3}vG_{0,4}+\tfrac{1}{6}vG_{3,1}+\tfrac{4}{3}vG_{4,0}+\tfrac{2}{3}uF_{2,2}
+
\\
&
\ \ \ \ \
-\tfrac{4}{3}uF_{4,0}-\tfrac{1}{3}uG_{0,4}-\tfrac{1}{3}uG_{3,1}+\tfrac{1}{3}uG_{4,0}+2x
)\partial_v
,
\\
e_2
&
\,:=\,
-(
-\tfrac{3}{2}uG_{3,1}-\tfrac{3}{2}uG_{4,0}+\tfrac{2}{3}xF_{2,2}-\tfrac{4}{3}xG_{0,4}-\tfrac{1}{2}xG_{3,1}-xG_{4,0}+\tfrac{1}{3}yF_{2,2}-\tfrac{1}{6}yG_{0,4}
+
\\
&
\ \ \ \ \
-\tfrac{1}{4}yG_{3,1}-\tfrac{1}{2}yG_{4,0}+2vF_{2,2}-vG_{0,4}-vG_{3,1}-3vG_{4,0}+6uF_{4,0}-\tfrac{5}{2}uG_{0,4}+\tfrac{7}{6}xF_{3,1}
+
\\
&
\ \ \ \ \
+\tfrac{1}{3}yF_{3,1}+2vF_{3,1}+uF_{2,2}+4uF_{3,1}+64u+\tfrac{4}{3}y+\tfrac{14}{3}x-8v
)\partial_x
+
(
1-\tfrac{7}{2}uG_{3,1}
+
\\
&
\ \ \ \ \
-8uG_{4,0}-\tfrac{1}{3}xF_{2,2}+\tfrac{1}{6}xG_{0,4}+\tfrac{1}{4}xG_{3,1}+\tfrac{1}{2}xG_{4,0}-\tfrac{2}{3}yF_{2,2}+\tfrac{4}{3}yG_{0,4}+\tfrac{1}{2}yG_{3,1}
+
\\
&
\ \ \ \ \
+yG_{4,0}+\tfrac{4}{3}vF_{2,2}-\tfrac{8}{3}vG_{0,4}-vG_{3,1}-2vG_{4,0}-\tfrac{8}{3}uG_{0,4}-\tfrac{1}{3}xF_{3,1}-\tfrac{7}{6}yF_{3,1}+\tfrac{4}{3}vF_{3,1}
+
\\
&
\ \ \ \ \
+\tfrac{10}{3}uF_{2,2}+\tfrac{16}{3}uF_{3,1}+\tfrac{16}{3}u-\tfrac{20}{3}y-\tfrac{4}{3}x+\tfrac{16}{3}v
)\partial_y
+
\\
&
\ \ \ \ \
-
(
1-\tfrac{7}{2}uG_{3,1}-8uG_{4,0}-\tfrac{1}{3}xF_{2,2}+\tfrac{1}{6}xG_{0,4}+\tfrac{1}{4}xG_{3,1}+\tfrac{1}{2}xG_{4,0}-\tfrac{2}{3}yF_{2,2}+\tfrac{4}{3}yG_{0,4}
+
\\
&
\ \ \ \ \
+\tfrac{1}{2}yG_{3,1}+yG_{4,0}+\tfrac{4}{3}vF_{2,2}-\tfrac{8}{3}vG_{0,4}-vG_{3,1}-2vG_{4,0}-\tfrac{8}{3}uG_{0,4}-\tfrac{1}{3}xF_{3,1}-\tfrac{7}{6}yF_{3,1}
+
\\
&
\ \ \ \ \
+\tfrac{4}{3}vF_{3,1}+\tfrac{10}{3}uF_{2,2}+\tfrac{16}{3}uF_{3,1}+\tfrac{16}{3}u-\tfrac{20}{3}y-\tfrac{4}{3}x+\tfrac{16}{3}v
)\partial_u
+
\\
&
\ \ \ \ \
-
(
-\tfrac{32}{3}u+\tfrac{40}{3}v+\tfrac{7}{3}vF_{3,1}+\tfrac{4}{3}uF_{2,2}-uG_{3,1}-2uG_{4,0}-\tfrac{2}{3}uG_{0,4}+\tfrac{4}{3}uF_{3,1}+\tfrac{4}{3}vF_{2,2}
+
\\
&
\ \ \ \ \
-vG_{3,1}-2vG_{4,0}-\tfrac{8}{3}vG_{0,4}-2y
)\partial_v.
\end{array}
\]

%%%%%%%%%%%%%%%%%%%%%%%%%%%%%%%%%%%%%%%%%%%%%%%%%%%%%%%%%%%%%%%%%%%%%%
\noindent
Gr\"obner basis generators of 
moduli space core algebraic variety in 
$\R^6 \ni F_{2,2},F_{3,1},G_{3,1},G_{0,4},F_{4,0},G_{4,0}$:
\[
\aligned
\B_1 
&
:= 2F_{2,2}+8+2F_{3,1}-4F_{4,0}-3G_{0,4}-G_{3,1}-G_{4,0},
\\
\B_2 
&
:= -2F_{3,1}G_{3,1}+8F_{3,1}G_{4,0}-8F_{4,0}G_{3,1}+32F_{4,0}G_{4,0}+2G_{0,4}^2+G_{0,4}G_{3,1}-8G_{0,4}G_{4,0}
+
\\
&
\ \ \ \ \
+G_{3,1}^2-G_{3,1}G_{4,0}-10G_{4,0}^2+152F_{3,1}+608F_{4,0}-92G_{0,4}-84G_{3,1}-52G_{4,0}+2080,
\\
\B_3 
&
:= 2F_{3,1}G_{0,4}-2F_{3,1}G_{3,1}+6F_{3,1}G_{4,0}+8F_{4,0}G_{0,4}-8F_{4,0}G_{3,1}+24F_{4,0}G_{4,0}-8G_{0,4}G_{4,0}
+
\\
&
\ \ \ \ \
+G_{3,1}^2-8G_{4,0}^2+120F_{3,1}+480F_{4,0}-44G_{0,4}-76G_{3,1}-68G_{4,0}+1568,
\\
\B_4 
&
:= (8F_{4,0}+G_{3,1}-8G_{4,0}-104)(8F_{4,0}+2F_{3,1}-G_{0,4}-3G_{4,0}+24-G_{3,1}),
\\
\B_5 
&
:= 4F_{3,1}^2-8F_{3,1}G_{3,1}+16F_{3,1}G_{4,0}-64F_{4,0}^2+16F_{4,0}G_{0,4}-16F_{4,0}G_{3,1}+112F_{4,0}G_{4,0}
+
\\
&
\ \ \ \ \
+2G_{0,4}G_{3,1}-16G_{0,4}G_{4,0}+3G_{3,1}^2-2G_{3,1}G_{4,0}-32G_{4,0}^2+480F_{3,1}+1536F_{4,0}
+
\\
&
\ \ \ \ \
-200G_{0,4}-280G_{3,1}-344G_{4,0}+5312,
\\
\B_6 &
:= (G_{3,1}^2-8G_{3,1}G_{4,0}+16G_{4,0}^2-124G_{3,1}+560G_{4,0}+4512)
\\
&
\ \ \ \ \ \
(8F_{4,0}+2F_{3,1}-G_{0,4}-3G_{4,0}+24-G_{3,1}).
\endaligned
\]

%%%%%%%%%%%%%%%%%%%%%%%%%%%%%%%%%%%%%%%%%%%%%%%%%%%%%%%%%%%%%%%%%%%%%%
\[
\footnotesize
\def\arraystretch{1.25}
\begin{array}{c|cc}
{} & e_1 & e_2 
\\
\hline
e_1 & 0 
&
\rotatebox[origin=c]{0}{
\begin{tabular}{p{7cm}}
$
(
-\tfrac{22}{3}-\tfrac{5}{6}F_{2,2}+\tfrac{7}{12}G_{3,1}+\tfrac{11}{12}G_{4,0}+\tfrac{17}{12}G_{0,4}-\tfrac{7}{6}F_{3,1}+\tfrac{1}{3}F_{4,0}
)e_1
+
(
\tfrac{4}{3}-\tfrac{2}{3}F_{2,2}+\tfrac{1}{6}G_{3,1}-\tfrac{1}{6}G_{4,0}+\tfrac{5}{6}G_{0,4}-\tfrac{1}{3}F_{3,1}+\tfrac{8}{3}F_{4,0}
)e_2
$
\end{tabular}}
\\
e_2 &
\rotatebox[origin=c]{0}{
\begin{tabular}{p{7cm}}
$
-
(
-\tfrac{22}{3}-\tfrac{5}{6}F_{2,2}+\tfrac{7}{12}G_{3,1}+\tfrac{11}{12}G_{4,0}+\tfrac{17}{12}G_{0,4}-\tfrac{7}{6}F_{3,1}+\tfrac{1}{3}F_{4,0}
)e_1
-
(
\tfrac{4}{3}-\tfrac{2}{3}F_{2,2}+\tfrac{1}{6}G_{3,1}-\tfrac{1}{6}G_{4,0}+\tfrac{5}{6}G_{0,4}-\tfrac{1}{3}F_{3,1}+\tfrac{8}{3}F_{4,0}
)e_2
$
\end{tabular}}
& 0
\end{array}
\]

%%%%%%%%%%%%%%%%%%%%%%%%%%%%%%%%%%%%%%%%%%%%%%%%%%%%%%%%%%%%%%%%%%%%%%
%%%%%%%%%%%%%%%%%%%%%%%%%%%%%%%%%%%%%%%%%%%%%%%%%%%%%%%%%%%%%%%%%%%%%%
%%%%%%%%%%%%%%%%%%%%%%%%%%%%%%%%%%%%%%%%%%%%%%%%%%%%%%%%%%%%%%%%%%%%%%

\[
\text{\bf Model 2f3e}
\ \ \ \ \
\left\{
\aligned
u
&
\,=\,
xy+2y^3
+
(
-8+\tfrac{1}{2}F_{3,1}-\tfrac{1}{4}G_{0,4}+\tfrac{1}{2}F_{2,2}-\tfrac{1}{4}G_{3,1}-\tfrac{3}{4}G_{4,0}
)y^4
+
\\
&
\ \ \ \ \
+
(
4+2F_{4,0}+\tfrac{1}{2}G_{0,4}+\tfrac{1}{2}G_{3,1}+\tfrac{3}{2}G_{4,0}-F_{2,2}
)xy^3
+
F_{2,2}x^2y^2
+
\\
&
\ \ \ \ \
+F_{3,1}x^3y+F_{4,0}x^4
+
\cdots
,
\\
v
&
\,=\,
y^2+x^2-4xy^2+8x^2y+G_{0,4}y^4
+
(
-2G_{0,4}+G_{3,1}+2G_{4,0}
)xy^3
+
\\
&
\ \ \ \ \
+
(
28+6F_{4,0}+\tfrac{3}{2}G_{0,4}-\tfrac{3}{2}G_{3,1}-\tfrac{3}{2}G_{4,0}+F_{2,2}+3F_{3,1}
)x^2y^2
+
\\
&
\ \ \ \ \
+
G_{3,1}x^3y+G_{4,0}x^4
+
\cdots
.
\endaligned\right.
\]
%%%%%%%%%%%%%%%%%%%%%%%%%%%%%%%%%%%%%%%%%%%%%%%%%%%%%%%%%%%%%%%%%%%%%%
\[
\aligned
e_1
&
\,:=\,
-(-1-\tfrac{80}{3}u-\tfrac{32}{3}v-\tfrac{1}{2}xF_{2,2}+\tfrac{1}{12}xG_{3,1}+\tfrac{7}{3}xF_{4,0}+\tfrac{7}{12}xG_{0,4}-\tfrac{7}{12}xG_{4,0}+\tfrac{2}{3}yF_{4,0}
+
\\
&
\ \ \ \ \
+\tfrac{1}{6}yG_{0,4}-\tfrac{1}{12}yG_{3,1}-\tfrac{1}{6}yG_{4,0}+2uF_{2,2}+\tfrac{5}{6}uG_{3,1}-\tfrac{32}{3}uF_{4,0}-\tfrac{5}{3}uG_{0,4}+\tfrac{5}{3}uG_{4,0}
+
\\
&
\ \ \ \ \
-\tfrac{8}{3}vF_{4,0}-\tfrac{2}{3}vG_{0,4}+\tfrac{1}{3}vG_{3,1}+\tfrac{8}{3}vG_{4,0}+\tfrac{8}{3}y+\tfrac{34}{3}x)\partial_x+(-32u+\tfrac{1}{2}yF_{2,2}-3uF_{3,1}
+
\\
&
\ \ \ \ \
+\tfrac{1}{12}xG_{3,1}-\tfrac{2}{3}xF_{4,0}-\tfrac{1}{6}xG_{0,4}+\tfrac{1}{6}xG_{4,0}-\tfrac{7}{3}yF_{4,0}-\tfrac{7}{12}yG_{0,4}-\tfrac{1}{12}yG_{3,1}+\tfrac{7}{12}yG_{4,0}
+
\\
&
\ \ \ \ \
+uG_{3,1}-8uF_{4,0}-2uG_{0,4}+4uG_{4,0}-4vF_{4,0}-\tfrac{28}{3}y-\tfrac{8}{3}x)\partial_y+(uF_{2,2}+y-\tfrac{1}{6}uG_{3,1}
+
\\
&
\ \ \ \ \
-\tfrac{14}{3}uF_{4,0}-\tfrac{62}{3}u-\tfrac{7}{6}uG_{0,4}+\tfrac{7}{6}uG_{4,0}-\tfrac{2}{3}vF_{4,0}-\tfrac{8}{3}v-\tfrac{1}{6}vG_{0,4}+\tfrac{1}{12}vG_{3,1}+\tfrac{1}{6}vG_{4,0})\partial_u
+
\\
&
\ \ \ \ \
+(2x+\tfrac{16}{3}u-\tfrac{8}{3}uF_{4,0}-\tfrac{2}{3}uG_{0,4}+\tfrac{1}{3}uG_{3,1}+\tfrac{2}{3}uG_{4,0}-\tfrac{68}{3}v+vF_{2,2}-\tfrac{1}{6}vG_{3,1}-\tfrac{14}{3}vF_{4,0}
+
\\
&
\ \ \ \ \
-\tfrac{7}{6}vG_{0,4}+\tfrac{7}{6}vG_{4,0})\partial_v
,
\\
e_2
&
\,:=\,
(-\tfrac{196}{3}u-\tfrac{64}{3}v+\tfrac{7}{6}xG_{0,4}+\tfrac{3}{4}xG_{3,1}+\tfrac{3}{2}xG_{4,0}-xF_{2,2}-\tfrac{4}{3}xF_{3,1}-\tfrac{1}{6}yF_{3,1}+\tfrac{1}{3}yG_{0,4}
+
\\
&
\ \ \ \ \
-6uF_{4,0}-\tfrac{1}{2}vG_{3,1}-\tfrac{29}{6}uG_{0,4}-\tfrac{3}{2}uG_{3,1}-\tfrac{9}{2}uG_{4,0}+3uF_{2,2}+\tfrac{8}{3}uF_{3,1}+\tfrac{2}{3}vF_{3,1}-\tfrac{4}{3}vG_{0,4}
+
\\
&
\ \ \ \ \
+\tfrac{56}{3}x-\tfrac{2}{3}y)\partial_x-(-64u-1-\tfrac{3}{4}yG_{3,1}-\tfrac{3}{2}yG_{4,0}+yF_{2,2}-\tfrac{1}{3}xG_{0,4}+\tfrac{1}{6}xF_{3,1}+\tfrac{4}{3}yF_{3,1}
+
\\
&
\ \ \ \ \
-\tfrac{7}{6}yG_{0,4}-4uG_{0,4}-\tfrac{1}{2}uG_{3,1}+2uF_{2,2}+2uF_{3,1}+vF_{3,1}-\tfrac{16}{3}x-\tfrac{68}{3}y)\partial_y-(-x-\tfrac{124}{3}u
+
\\
&
\ \ \ \ \
-\tfrac{7}{3}uG_{0,4}-\tfrac{3}{2}uG_{3,1}-3uG_{4,0}+2uF_{2,2}+\tfrac{8}{3}uF_{3,1}+\tfrac{1}{6}vF_{3,1}-\tfrac{1}{3}vG_{0,4}-\tfrac{16}{3}v)\partial_u-(-2y
+
\\
&
\ \ \ \ \
-\tfrac{4}{3}u+\tfrac{2}{3}uF_{3,1}-\tfrac{4}{3}uG_{0,4}-\tfrac{136}{3}v-\tfrac{7}{3}vG_{0,4}-\tfrac{3}{2}vG_{3,1}-3vG_{4,0}+2vF_{2,2}+\tfrac{8}{3}vF_{3,1})\partial_v
.
\endaligned
\]

%%%%%%%%%%%%%%%%%%%%%%%%%%%%%%%%%%%%%%%%%%%%%%%%%%%%%%%%%%%%%%%%%%%%%%
\noindent
Gr\"obner basis generators of 
moduli space core algebraic variety in 
$\R^6 \ni F_{2,2},F_{3,1},F_{4,0},G_{0,4},G_{3,1},G_{4,0}$:
\[
\aligned
\B_1 & := 2F_{2,2}-48-3G_{0,4}-G_{4,0}+2F_{3,1}-4F_{4,0}-G_{3,1},
\\
\B_2 & := 2F_{3,1}G_{3,1}-8F_{3,1}G_{4,0}+8F_{4,0}G_{3,1}-32F_{4,0}G_{4,0}+4G_{0,4}^2+G_{0,4}G_{3,1}-12G_{0,4}G_{4,0}
+
\\
&
\ \ \ \ \
-G_{3,1}^2-G_{3,1}G_{4,0}+24G_{4,0}^2+40F_{3,1}+160F_{4,0}+116G_{0,4}+28G_{3,1}-132G_{4,0}+832,
\\
\B_3 & := 4F_{3,1}G_{0,4}-2F_{3,1}G_{3,1}+4F_{3,1}G_{4,0}+16F_{4,0}G_{0,4}-8F_{4,0}G_{3,1}+16F_{4,0}G_{4,0}-3G_{0,4}G_{3,1}
+
\\
&
\ \ \ \ \
-4G_{0,4}G_{4,0}+G_{3,1}^2+3G_{3,1}G_{4,0}-12G_{4,0}^2+24F_{3,1}+96F_{4,0}-4G_{0,4}-44G_{3,1}
+
\\
&
\ \ \ \ \
-44G_{4,0}-64,
\\
\B_4 & := (16F_{4,0}-G_{3,1}-4G_{4,0}-12)(8F_{4,0}+2F_{3,1}+G_{0,4}-G_{3,1}-5G_{4,0}+16),
\\
\B_5 & := 4F_{3,1}^2-2F_{3,1}G_{3,1}-8F_{3,1}G_{4,0}-64F_{4,0}^2-16F_{4,0}G_{0,4}+8F_{4,0}G_{3,1}+48F_{4,0}G_{4,0}
+
\\
&
\ \ \ \ \
+G_{0,4}G_{3,1}+4G_{0,4}G_{4,0}-G_{3,1}G_{4,0}-4G_{4,0}^2+56F_{3,1}-32F_{4,0}+20G_{0,4}-4G_{3,1}
+
\\
&
\ \ \ \ \
-36G_{4,0}+320,
\\
\B_6 & := (G_{3,1}^2-8G_{3,1}G_{4,0}+16G_{4,0}^2+48G_{3,1}+64G_{4,0}-80)(8F_{4,0}+2F_{3,1}+G_{0,4}-G_{3,1}-5G_{4,0}+16).
\endaligned
\]

%%%%%%%%%%%%%%%%%%%%%%%%%%%%%%%%%%%%%%%%%%%%%%%%%%%%%%%%%%%%%%%%%%%%%%
\[
\footnotesize
\def\arraystretch{1.25}
\begin{array}{c|cc}
{} & e_1 & e_2 
\\
\hline
e_1 & 0
&
\rotatebox[origin=c]{0}{
\begin{tabular}{p{7cm}}
$
(
\tfrac{64}{3}+\tfrac{4}{3}G_{0,4}+\tfrac{2}{3}G_{3,1}+\tfrac{4}{3}G_{4,0}-F_{2,2}-\tfrac{4}{3}F_{3,1}+\tfrac{2}{3}F_{4,0}
)e_1
+
(
\tfrac{44}{3}+\tfrac{11}{12}G_{0,4}+\tfrac{1}{12}G_{3,1}-\tfrac{7}{12}G_{4,0}-\tfrac{1}{2}F_{2,2}-\tfrac{1}{6}F_{3,1}+\tfrac{7}{3}F_{4,0}
)e_2
$
\end{tabular}}
\\
e_2 & 
\rotatebox[origin=c]{0}{
\begin{tabular}{p{7cm}}
$
-(
\tfrac{64}{3}+\tfrac{4}{3}G_{0,4}+\tfrac{2}{3}G_{3,1}+\tfrac{4}{3}G_{4,0}-F_{2,2}-\tfrac{4}{3}F_{3,1}+\tfrac{2}{3}F_{4,0}
)e_1
+
-(
\tfrac{44}{3}+\tfrac{11}{12}G_{0,4}+\tfrac{1}{12}G_{3,1}-\tfrac{7}{12}G_{4,0}-\tfrac{1}{2}F_{2,2}-\tfrac{1}{6}F_{3,1}+\tfrac{7}{3}F_{4,0}
)e_2
$
\end{tabular}}
& 0
\end{array}
\]

%%%%%%%%%%%%%%%%%%%%%%%%%%%%%%%%%%%%%%%%%%%%%%%%%%%%%%%%%%%%%%%%%%%%%%
%%%%%%%%%%%%%%%%%%%%%%%%%%%%%%%%%%%%%%%%%%%%%%%%%%%%%%%%%%%%%%%%%%%%%%
%%%%%%%%%%%%%%%%%%%%%%%%%%%%%%%%%%%%%%%%%%%%%%%%%%%%%%%%%%%%%%%%%%%%%%

\[
\!\!\!\!\!
\text{\bf Model 2f3f4a}
\ \ \ \ \
\left\{
\aligned
u
&
\,=\,
xy+x^3+y^3
+
(\tfrac{-3}{2}+\tfrac{1}{2}F_{2,2}-G_{4,0})x^4
+
(-4+F_{2,2}-G_{4,0})x^3y+F_{2,2}x^2y^2
+
\\
&
\,\,\,\,\,\,\,\,
+
(-4+F_{2,2}-G_{4,0})xy^3
+
(\tfrac{-3}{2}+\tfrac{1}{2}F_{2,2}-G_{4,0})y^4
+
(\tfrac{-29}{10}-\tfrac{3}{10}F_{2,2}+\tfrac{16}{5}G_{4,0}+
\\
&
\,\,\,\,\,\,\,\,
-F_{2,2}G_{4,0}
+\tfrac{1}{5}F_{2,2}^2+\tfrac{6}{5}G_{4,0}^2)x^5
+
(\tfrac{1}{2}F_{2,2}^2-2F_{2,2}G_{4,0}+\tfrac{3}{2}G_{4,0}^2-F_{2,2}+\tfrac{17}{2}G_{4,0}
+
\\
&
\,\,\,\,\,\,\,\,-\tfrac{1}{2})x^4y
+
(\tfrac{2}{3}F_{2,2}^2-2F_{2,2}G_{4,0}+\tfrac{7}{3}F_{2,2}-\tfrac{4}{3}G_{4,0}-17)x^3y^2
+
(F_{2,2}^2-4F_{2,2}G_{4,0}
+
\\
&
\,\,\,\,\,\,\,\,
+3G_{4,0}^2+2F_{2,2}+G_{4,0}-19)x^2y^3
+
(\tfrac{-7}{2}-\tfrac{3}{2}F_{2,2}+12G_{4,0}-5F_{2,2}G_{4,0}+F_{2,2}^2
+
\\
&
\,\,\,\,\,\,\,\,
+6G_{4,0}^2)xy^4
+
(\tfrac{-29}{10}-\tfrac{3}{10}F_{2,2}+\tfrac{16}{5}G_{4,0}-F_{2,2}G_{4,0}+\tfrac{1}{5}F_{2,2}^2+\tfrac{6}{5}G_{4,0}^2)y^5
+
\cdots
,\\
v
& 
\,=\,
x^2+y^2+2x^2y+2xy^2+G_{4,0}x^4
+
(2F_{2,2}-4-2G_{4,0})x^3y
+
(-14+4F_{2,2}
+
\\
&
\,\,\,\,\,\,\,\,
-6G_{4,0})x^2y^2
+
(2F_{2,2}-4-2G_{4,0})xy^3
+
G_{4,0}y^4
+
(\tfrac{2}{5}F_{2,2}G_{4,0}-\tfrac{6}{5}G_{4,0}^2
+
\\
&
\,\,\,\,\,\,\,\,
+\tfrac{9}{5}F_{2,2}-\tfrac{12}{5}G_{4,0}-\tfrac{27}{5})x^5
+
(F_{2,2}^2-4F_{2,2}G_{4,0}+3G_{4,0}^2+4F_{2,2}-G_{4,0}
+
\\
&
\,\,\,\,\,\,\,\,
-27)x^4y
+
(\tfrac{8}{3}F_{2,2}^2-12F_{2,2}G_{4,0}+12G_{4,0}^2-\tfrac{14}{3}F_{2,2}+\tfrac{92}{3}G_{4,0}-14)x^3y^2
+
\\
&
\,\,\,\,\,\,\,\,
+
(2F_{2,2}^2-8F_{2,2}G_{4,0}+6G_{4,0}^2-4F_{2,2}+26G_{4,0}-10)x^2y^3
+
(2F_{2,2}G_{4,0}
+
\\
&
\,\,\,\,\,\,\,\,
-6G_{4,0}^2+5F_{2,2}-8G_{4,0}-21)xy^4
+
(\tfrac{2}{5}F_{2,2}G_{4,0}-\tfrac{6}{5}G_{4,0}^2+\tfrac{9}{5}F_{2,2}-\tfrac{12}{5}G_{4,0}
+
\\
&
\,\,\,\,\,\,\,\,
-\tfrac{27}{5})y^5
+
\cdots
.
\endaligned\right.
\]
%%%%%%%%%%%%%%%%%%%%%%%%%%%%%%%%%%%%%%%%%%%%%%%%%%%%%%%%%%%%%%%%%%%%%%
\[
\def\arraystretch{1.25}
\aligned
e_1 & := 
-
(F_{2,2}u+F_{2,2}x+G_{4,0}u+2G_{4,0}v-3G_{4,0}x-2u-5v-x-3y-1)\partial_x
-
(3F_{2,2}u
+
\\
&
\,\,\,\,\,\,\,\,
+F_{2,2}v+F_{2,2}y-5G_{4,0}u-G_{4,0}v-3G_{4,0}y-4u-4v)\partial_y
-
(2F_{2,2}u-6G_{4,0}u-u
+
\\
&
\,\,\,\,\,\,\,\,
-3v-y)\partial_u
-
(2F_{2,2}v-6G_{4,0}v-10u-2v-2x)\partial_v,\\
e_2 & := 
-
(3F_{2,2}u+F_{2,2}v+F_{2,2}x-5G_{4,0}u-G_{4,0}v-3G_{4,0}x-10u-v-3x+3y)\partial_x
+
\\
&
\,\,\,\,\,\,\,\,
-
(F_{2,2}u+F_{2,2}y+G_{4,0}u+2G_{4,0}v-3G_{4,0}y-8u-2v-4y-1)\partial_y
-
(2F_{2,2}u
+
\\
&
\,\,\,\,\,\,\,\,
-6G_{4,0}u-7u-x)\partial_u
-
(2F_{2,2}v-6G_{4,0}v+2u-8v-2y)\partial_v,\\
e_3 & := 
-
(2u-v+x-y)\partial_x
-
(2u-v-x+y)\partial_y
-
(2u-v)\partial_u
+
(4u-2v)\partial_v,
\endaligned
\]
with $F_{2,2},G_{4,0}$ satisfying:
\[
0 = F_{2,2}^2-6F_{2,2}G_{4,0}+9G_{4,0}^2-F_{2,2}+7G_{4,0}-6,
\]

%%%%%%%%%%%%%%%%%%%%%%%%%%%%%%%%%%%%%%%%%%%%%%%%%%%%%%%%%%%%%%%%%%%%%%
\[
\footnotesize
\def\arraystretch{1.25}
\begin{array}{c|ccc}
{} & e_1 & e_2 & e_3
\\
\hline
e_1 & 
0 &
\rotatebox[origin=c]{0}{
\begin{tabular}{p{6cm}}
$(-F_{2,2}+3G_{4,0})e_1+(F_{2,2}-3G_{4,0})e_2+(3F_{2,2}-9G_{4,0})e_3$
\end{tabular}}
&
-e_1+e_2+3e_3
\\
e_2 &
\rotatebox[origin=c]{0}{
\begin{tabular}{p{6cm}}
$-(-F_{2,2}+3G_{4,0})e_1-(F_{2,2}-3G_{4,0})e_2-(3F_{2,2}-9G_{4,0})e_3$
\end{tabular}}
&
e_1-e_2-3e_3
\\
e_3 &
e_1-e_2-3e_3& -e_1+e_2+3e_3& 0
\end{array}
\]

%%%%%%%%%%%%%%%%%%%%%%%%%%%%%%%%%%%%%%%%%%%%%%%%%%%%%%%%%%%%%%%%%%%%%%
%%%%%%%%%%%%%%%%%%%%%%%%%%%%%%%%%%%%%%%%%%%%%%%%%%%%%%%%%%%%%%%%%%%%%%
%%%%%%%%%%%%%%%%%%%%%%%%%%%%%%%%%%%%%%%%%%%%%%%%%%%%%%%%%%%%%%%%%%%%%%

\[
\text{\bf Model 2f3f4b}
\ \ \ \ \
\left\{
\aligned
u
&
\,=\,
xy+x^3+y^3
+
(-\tfrac{1}{40}G_{3,2}-\tfrac{1}{10}G_{4,0}+\tfrac{1}{4}G_{5,0})x^4
+
(\tfrac{-2}{3}-\tfrac{1}{20}G_{3,2}+\tfrac{4}{5}G_{4,0}
+
\\
&
\,\,\,\,\,\,\,\,
+\tfrac{1}{2}G_{5,0})x^3y
+
(\tfrac{7}{3}-\tfrac{1}{20}G_{3,2}+\tfrac{9}{5}G_{4,0}+\tfrac{1}{2}G_{5,0})x^2y^2
+
(\tfrac{-2}{3}-\tfrac{1}{20}G_{3,2}+\tfrac{4}{5}G_{4,0}
+
\\
&
\,\,\,\,\,\,\,\,
+\tfrac{1}{2}G_{5,0})xy^3
+
(-\tfrac{1}{40}G_{3,2}-\tfrac{1}{10}G_{4,0}+\tfrac{1}{4}G_{5,0})y^4
+
(\tfrac{-259}{135}-\tfrac{59}{450}G_{3,2}+\tfrac{63}{25}G_{4,0}
+
\\
&
\,\,\,\,\,\,\,\,
+\tfrac{1}{1200}G_{3,2}^2+\tfrac{1}{12}G_{5,0}^2+\tfrac{1}{50}G_{3,2}G_{4,0}-\tfrac{1}{60}G_{3,2}G_{5,0}-\tfrac{1}{5}G_{4,0}G_{5,0}+\tfrac{73}{90}G_{5,0})x^5
+
\\
&
\,\,\,\,\,\,\,\,
+
(\tfrac{1}{3}+\tfrac{1}{4}G_{3,2}+2G_{4,0})x^4y
+
(\tfrac{-122}{27}-\tfrac{68}{45}G_{3,2}+\tfrac{42}{5}G_{4,0}+\tfrac{1}{120}G_{3,2}^2+\tfrac{5}{6}G_{5,0}^2
+
\\
&
\,\,\,\,\,\,\,\,
+\tfrac{1}{5}G_{3,2}G_{4,0}-\tfrac{1}{6}G_{3,2}G_{5,0}-2G_{4,0}G_{5,0}+\tfrac{91}{9}G_{5,0})x^3y^2
+
(\tfrac{4}{3}+\tfrac{3}{10}G_{3,2}
+
\\
&
\,\,\,\,\,\,\,\,
-\tfrac{24}{5}G_{4,0}+2G_{5,0})x^2y^3
+
(-\tfrac{59}{90}G_{3,2}+\tfrac{43}{5}G_{4,0}+\tfrac{1}{240}G_{3,2}^2+\tfrac{5}{12}G_{5,0}^2
+
\\
&
\,\,\,\,\,\,\,\,
+\tfrac{1}{10}G_{3,2}G_{4,0}-\tfrac{1}{12}G_{3,2}G_{5,0}-G_{4,0}G_{5,0}+\tfrac{73}{18}G_{5,0}-\tfrac{70}{27})xy^4
+
(\tfrac{-4}{3}+\tfrac{1}{20}G_{3,2}
+
\\
&
\,\,\,\,\,\,\,\,
+\tfrac{6}{5}G_{4,0})y^5
+
\cdots
,
\\
v
&
\,=\,
x^2+y^2+2x^2y+2xy^2+G_{4,0}x^4
+
(\tfrac{4}{3}-\tfrac{1}{10}G_{3,2}+\tfrac{8}{5}G_{4,0}+G_{5,0})x^3y
+
\\
&
\,\,\,\,\,\,\,\,
+
(-\tfrac{2}{3}-\tfrac{1}{5}G_{3,2}+\tfrac{6}{5}G_{4,0}+2G_{5,0})x^2y^2
+
(\tfrac{4}{3}-\tfrac{1}{10}G_{3,2}+\tfrac{8}{5}G_{4,0}+G_{5,0})xy^3
+
\\
&
\,\,\,\,\,\,\,\,
+
G_{4,0}y^4+G_{5,0}x^5
+
(\tfrac{1}{120}G_{3,2}^2+\tfrac{5}{6}G_{5,0}^2+\tfrac{1}{5}G_{3,2}G_{4,0}-\tfrac{1}{6}G_{3,2}G_{5,0}-2G_{4,0}G_{5,0}
+
\\
&
\,\,\,\,\,\,\,\,
-\tfrac{29}{18}G_{3,2}+\tfrac{100}{9}G_{5,0}-140}{27+10G_{4,0})x^4y
+
G_{3,2}x^3y^2
+
(\tfrac{-316}{27}-\tfrac{118}{45}G_{3,2}
+
\\
&
\,\,\,\,\,\,\,\,
+\tfrac{132}{5}G_{4,0}+\tfrac{1}{60}G_{3,2}^2+\tfrac{5}{3}G_{5,0}^2+\tfrac{2}{5}G_{3,2}G_{4,0}-\tfrac{1}{3}G_{3,2}G_{5,0}-4G_{4,0}G_{5,0}
+
\\
&
\,\,\,\,\,\,\,\,
+\tfrac{146}{9}G_{5,0})x^2y^3
+
(\tfrac{2}{3}+\tfrac{1}{5}G_{3,2}-\tfrac{16}{5}G_{4,0}+3G_{5,0})xy^4
+
(\tfrac{-158}{135}-\tfrac{163}{450}G_{3,2}
+
\\
&
\,\,\,\,\,\,\,\,
+\tfrac{66}{25}G_{4,0}+\tfrac{1}{600}G_{3,2}^2+\tfrac{1}{6}G_{5,0}^2+\tfrac{1}{25}G_{3,2}G_{4,0}-\tfrac{1}{30}G_{3,2}G_{5,0}-\tfrac{2}{5}G_{4,0}G_{5,0}
+
\\
&
\,\,\,\,\,\,\,\,+\tfrac{118}{45}G_{5,0})y^5
+
\cdots
.
\endaligned\right.
\]
%%%%%%%%%%%%%%%%%%%%%%%%%%%%%%%%%%%%%%%%%%%%%%%%%%%%%%%%%%%%%%%%%%%%%%
\[
\aligned
e_1 & := 
(1-\tfrac{1}{5}uG_{3,2}G_{5,0}-\tfrac{1}{10}xG_{3,2}G_{5,0}+\tfrac{1}{10}vG_{3,2}G_{5,0}+\tfrac{1}{10}yG_{3,2}G_{5,0}-\tfrac{9}{100}yG_{3,2}G_{4,0}+\tfrac{9}{50}uG_{3,2}G_{4,0}
+
\\
&
\,\,\,\,\,\,\,\,
+\tfrac{9}{10}vG_{4,0}G_{5,0}+\tfrac{9}{10}yG_{4,0}G_{5,0}-\tfrac{9}{5}uG_{4,0}G_{5,0}-\tfrac{9}{100}vG_{3,2}G_{4,0}+\tfrac{9}{100}xG_{3,2}G_{4,0}-\tfrac{9}{10}xG_{4,0}G_{5,0}
+
\\
&
\,\,\,\,\,\,\,\,
-\tfrac{32}{9}u+\tfrac{43}{9}y-\tfrac{37}{9}x+\tfrac{61}{9}v+uG_{5,0}^2+\tfrac{149}{120}yG_{3,2}-\tfrac{59}{12}yG_{5,0}-\tfrac{107}{10}yG_{4,0}+\tfrac{18}{25}yG_{4,0}^2+\tfrac{1}{100}uG_{3,2}^2
+
\\
&
\,\,\,\,\,\,\,\,
-\tfrac{73}{30}uG_{3,2}+\tfrac{28}{3}uG_{5,0}+\tfrac{93}{5}uG_{4,0}-\tfrac{36}{25}uG_{4,0}^2-\tfrac{1}{200}vG_{3,2}^2-\tfrac{1}{2}vG_{5,0}^2+\tfrac{149}{120}vG_{3,2}-\tfrac{59}{12}vG_{5,0}
+
\\
&
\,\,\,\,\,\,\,\,
-\tfrac{127}{10}vG_{4,0}+\tfrac{18}{25}vG_{4,0}^2-\tfrac{18}{25}xG_{4,0}^2+\tfrac{53}{12}xG_{5,0}+\tfrac{1}{200}xG_{3,2}^2-\tfrac{1}{200}yG_{3,2}^2-\tfrac{143}{120}xG_{3,2}-\tfrac{1}{2}yG_{5,0}^2
+
\\
&
\,\,\,\,\,\,\,\,
+\tfrac{119}{10}xG_{4,0}+\tfrac{1}{2}xG_{5,0}^2)\partial_x
+
(-\tfrac{1}{5}uG_{3,2}G_{5,0}+\tfrac{1}{10}xG_{3,2}G_{5,0}+\tfrac{1}{10}vG_{3,2}G_{5,0}-\tfrac{1}{10}yG_{3,2}G_{5,0}
+
\\
&
\,\,\,\,\,\,\,\,
+\tfrac{9}{100}yG_{3,2}G_{4,0}+\tfrac{9}{50}uG_{3,2}G_{4,0}+\tfrac{9}{10}vG_{4,0}G_{5,0}-\tfrac{9}{10}yG_{4,0}G_{5,0}-\tfrac{9}{5}uG_{4,0}G_{5,0}-\tfrac{9}{100}vG_{3,2}G_{4,0}
+
\\
&
\,\,\,\,\,\,\,\,
-\tfrac{9}{100}xG_{3,2}G_{4,0}+\tfrac{9}{10}xG_{4,0}G_{5,0}-\tfrac{86}{9}u-\tfrac{46}{9}y+\tfrac{16}{9}x+\tfrac{28}{9}v+uG_{5,0}^2-\tfrac{143}{120}yG_{3,2}+\tfrac{53}{12}yG_{5,0}
+
\\
&
\,\,\,\,\,\,\,\,
+\tfrac{119}{10}yG_{4,0}-\tfrac{18}{25}yG_{4,0}^2+\tfrac{1}{100}uG_{3,2}^2-\tfrac{7}{3}uG_{3,2}+\tfrac{25}{3}uG_{5,0}+21uG_{4,0}-\tfrac{36}{25}uG_{4,0}^2-\tfrac{1}{200}vG_{3,2}^2
+
\\
&
\,\,\,\,\,\,\,\,
-\tfrac{1}{2}vG_{5,0}^2+\tfrac{31}{24}vG_{3,2}-\tfrac{65}{12}vG_{5,0}-\tfrac{23}{2}vG_{4,0}+\tfrac{18}{25}vG_{4,0}^2+\tfrac{18}{25}xG_{4,0}^2-\tfrac{59}{12}xG_{5,0}
+
\\
&
\,\,\,\,\,\,\,\,
-\tfrac{1}{200}xG_{3,2}^2+\tfrac{1}{200}yG_{3,2}^2+\tfrac{149}{120}xG_{3,2}+\tfrac{1}{2}yG_{5,0}^2-\tfrac{107}{10}xG_{4,0}-\tfrac{1}{2}xG_{5,0}^2)\partial_y
+
(-\tfrac{83}{9}u+\tfrac{43}{9}v
+
\\
&
\,\,\,\,\,\,\,\,
+\tfrac{1}{100}uG_{3,2}^2+uG_{5,0}^2-\tfrac{143}{60}uG_{3,2}+\tfrac{53}{6}uG_{5,0}+\tfrac{119}{5}uG_{4,0}-\tfrac{36}{25}uG_{4,0}^2-\tfrac{1}{200}vG_{3,2}^2-\tfrac{1}{2}vG_{5,0}^2
+
\\
&
\,\,\,\,\,\,\,\,
+\tfrac{149}{120}vG_{3,2}-\tfrac{59}{12}vG_{5,0}-\tfrac{107}{10}vG_{4,0}+\tfrac{18}{25}vG_{4,0}^2+\tfrac{9}{50}uG_{3,2}G_{4,0}-\tfrac{1}{5}uG_{3,2}G_{5,0}-\tfrac{9}{5}uG_{4,0}G_{5,0}
+
\\
&
\,\,\,\,\,\,\,\,
-\tfrac{9}{100}vG_{3,2}G_{4,0}+\tfrac{1}{10}vG_{3,2}G_{5,0}+\tfrac{9}{10}vG_{4,0}G_{5,0}+y)\partial_u
-
(-\tfrac{154}{9}u+\tfrac{74}{9}v+\tfrac{1}{50}uG_{3,2}^2+2uG_{5,0}^2
+
\\
&
\,\,\,\,\,\,\,\,
-\tfrac{149}{30}uG_{3,2}+\tfrac{59}{3}uG_{5,0}+\tfrac{214}{5}uG_{4,0}-\tfrac{72}{25}uG_{4,0}^2-\tfrac{1}{100}vG_{3,2}^2-vG_{5,0}^2+\tfrac{143}{60}vG_{3,2}-\tfrac{53}{6}vG_{5,0}
+
\\
&
\,\,\,\,\,\,\,\,
-\tfrac{119}{5}vG_{4,0}+\tfrac{36}{25}vG_{4,0}^2+\tfrac{9}{25}uG_{3,2}G_{4,0}-\tfrac{2}{5}uG_{3,2}G_{5,0}-\tfrac{18}{5}uG_{4,0}G_{5,0}-\tfrac{9}{50}vG_{3,2}G_{4,0}
+
\\
&
\,\,\,\,\,\,\,\,
+\tfrac{1}{5}vG_{3,2}G_{5,0}+\tfrac{9}{5}vG_{4,0}G_{5,0}-2x)\partial_v,\\
e_2 & := 
-
(-\tfrac{3}{10}uG_{3,2}G_{5,0}-\tfrac{3}{20}xG_{3,2}G_{5,0}+\tfrac{3}{20}vG_{3,2}G_{5,0}+\tfrac{3}{20}yG_{3,2}G_{5,0}-\tfrac{21}{100}yG_{3,2}G_{4,0}+\tfrac{21}{50}uG_{3,2}G_{4,0}
+
\\
&
\,\,\,\,\,\,\,\,
+\tfrac{21}{10}vG_{4,0}G_{5,0}+\tfrac{21}{10}yG_{4,0}G_{5,0}-\tfrac{21}{5}uG_{4,0}G_{5,0}-\tfrac{21}{100}vG_{3,2}G_{4,0}+\tfrac{21}{100}xG_{3,2}G_{4,0}-\tfrac{21}{10}xG_{4,0}G_{5,0}
+
\\
&
\,\,\,\,\,\,\,\,
-8u+7y-\tfrac{11}{3}x+\tfrac{17}{3}v+\tfrac{3}{2}uG_{5,0}^2+\tfrac{59}{40}yG_{3,2}-\tfrac{29}{4}yG_{5,0}-\tfrac{91}{10}yG_{4,0}-\tfrac{18}{25}yG_{4,0}^2+\tfrac{3}{200}uG_{3,2}^2
+
\\
&
\,\,\,\,\,\,\,\,
-\tfrac{31}{10}uG_{3,2}+16uG_{5,0}+\tfrac{93}{5}uG_{4,0}+\tfrac{36}{25}uG_{4,0}^2-\tfrac{3}{400}vG_{3,2}^2-\tfrac{3}{4}vG_{5,0}^2+\tfrac{57}{40}vG_{3,2}-\tfrac{27}{4}vG_{5,0}
+
\\
&
\,\,\,\,\,\,\,\,
-\tfrac{83}{10}vG_{4,0}-\tfrac{18}{25}vG_{4,0}^2+\tfrac{18}{25}xG_{4,0}^2+\tfrac{31}{4}xG_{5,0}+\tfrac{3}{400}xG_{3,2}^2-\tfrac{3}{400}yG_{3,2}^2-\tfrac{61}{40}xG_{3,2}-\tfrac{3}{4}yG_{5,0}^2
+
\\
&
\,\,\,\,\,\,\,\,
+\tfrac{79}{10}xG_{4,0}+\tfrac{3}{4}xG_{5,0}^2)\partial_x
-
(-1-\tfrac{3}{10}uG_{3,2}G_{5,0}+\tfrac{3}{20}xG_{3,2}G_{5,0}+\tfrac{3}{20}vG_{3,2}G_{5,0}-\tfrac{3}{20}yG_{3,2}G_{5,0}
+
\\
&
\,\,\,\,\,\,\,\,
+\tfrac{21}{100}yG_{3,2}G_{4,0}+\tfrac{21}{50}uG_{3,2}G_{4,0}+\tfrac{21}{10}vG_{4,0}G_{5,0}-\tfrac{21}{10}yG_{4,0}G_{5,0}-\tfrac{21}{5}uG_{4,0}G_{5,0}-\tfrac{21}{100}vG_{3,2}G_{4,0}
+
\\
&
\,\,\,\,\,\,\,\,
-\tfrac{21}{100}xG_{3,2}G_{4,0}+\tfrac{21}{10}xG_{4,0}G_{5,0}-14u-\tfrac{14}{3}y+4x+2v+\tfrac{3}{2}uG_{5,0}^2-\tfrac{61}{40}yG_{3,2}+\tfrac{31}{4}yG_{5,0}
+
\\
&
\,\,\,\,\,\,\,\,
+\tfrac{79}{10}yG_{4,0}+\tfrac{18}{25}yG_{4,0}^2+\tfrac{3}{200}uG_{3,2}^2-3uG_{3,2}+15uG_{5,0}+21uG_{4,0}+\tfrac{36}{25}uG_{4,0}^2-\tfrac{3}{400}vG_{3,2}^2
+
\\
&
\,\,\,\,\,\,\,\,
-\tfrac{3}{4}vG_{5,0}^2+\tfrac{59}{40}vG_{3,2}-\tfrac{29}{4}vG_{5,0}-\tfrac{71}{10}vG_{4,0}-\tfrac{18}{25}vG_{4,0}^2-\tfrac{18}{25}xG_{4,0}^2-\tfrac{29}{4}xG_{5,0}-\tfrac{3}{400}xG_{3,2}^2
+
\\
&
\,\,\,\,\,\,\,\,
+\tfrac{3}{400}yG_{3,2}^2+\tfrac{59}{40}xG_{3,2}+\tfrac{3}{4}yG_{5,0}^2-\tfrac{91}{10}xG_{4,0}-\tfrac{3}{4}xG_{5,0}^2)\partial_y
-
(-\tfrac{25}{3}u+4v+\tfrac{3}{200}uG_{3,2}^2+\tfrac{3}{2}uG_{5,0}^2
+
\\
&
\,\,\,\,\,\,\,\,
+\tfrac{31}{2}uG_{5,0}-\tfrac{61}{20}uG_{3,2}+\tfrac{79}{5}uG_{4,0}+\tfrac{36}{25}uG_{4,0}^2-\tfrac{3}{400}vG_{3,2}^2-\tfrac{3}{4}vG_{5,0}^2-\tfrac{29}{4}vG_{5,0}+\tfrac{59}{40}vG_{3,2}
+
\\
&
\,\,\,\,\,\,\,\,
-\tfrac{91}{10}vG_{4,0}-\tfrac{18}{25}vG_{4,0}^2+\tfrac{21}{50}uG_{3,2}G_{4,0}-\tfrac{3}{10}uG_{3,2}G_{5,0}-\tfrac{21}{5}uG_{4,0}G_{5,0}-\tfrac{21}{100}vG_{3,2}G_{4,0}
+
\\
&
\,\,\,\,\,\,\,\,
+\tfrac{3}{20}vG_{3,2}G_{5,0}+\tfrac{21}{10}vG_{4,0}G_{5,0}-x)\partial_u
+
(-18u+\tfrac{28}{3}v+\tfrac{3}{100}uG_{3,2}^2+3uG_{5,0}^2+29uG_{5,0}
+
\\
&
\,\,\,\,\,\,\,\,
-\tfrac{59}{10}uG_{3,2}+\tfrac{182}{5}uG_{4,0}+\tfrac{72}{25}uG_{4,0}^2-\tfrac{3}{200}vG_{3,2}^2-\tfrac{3}{2}vG_{5,0}^2-\tfrac{31}{2}vG_{5,0}+\tfrac{61}{20}vG_{3,2}-\tfrac{79}{5}vG_{4,0}
+
\\
&
\,\,\,\,\,\,\,\,
-\tfrac{36}{25}vG_{4,0}^2+\tfrac{21}{25}uG_{3,2}G_{4,0}-\tfrac{3}{5}uG_{3,2}G_{5,0}-\tfrac{42}{5}uG_{4,0}G_{5,0}-\tfrac{21}{50}vG_{3,2}G_{4,0}+\tfrac{3}{10}vG_{3,2}G_{5,0}
+
\\
&
\,\,\,\,\,\,\,\,
+\tfrac{21}{5}vG_{4,0}G_{5,0}+2y)\partial_v
\endaligned
\]
with $G_{3,2},G_{4,0},G_{5,0}$ satisfying:
\[
\aligned
0 & = 9G_{3,2}^2+432G_{3,2}G_{4,0}-180G_{3,2}G_{5,0}+5184G_{4,0}^2-4320G_{4,0}G_{5,0}+900G_{5,0}^2-1020G_{3,2}
+
\\
&
\,\,\,\,\,\,\,\,
-10080G_{4,0}+10200G_{5,0}+4000,
\endaligned
\]
%%%%%%%%%%%%%%%%%%%%%%%%%%%%%%%%%%%%%%%%%%%%%%%%%%%%%%%%%%%%%%%%%%%%%%
\[
\footnotesize
\def\arraystretch{1.25}
\begin{array}{c|cc}
{} & e_1 & e_2 
\\
\hline
e_1 & 
0
&
\rotatebox[origin=c]{0}{
\begin{tabular}{p{7cm}}
$(\tfrac{-10}{9}-\tfrac{3}{25}G_{3,2}G_{4,0}+\tfrac{1}{20}G_{3,2}G_{5,0}+\tfrac{6}{5}G_{4,0}G_{5,0}-\tfrac{1}{400}G_{3,2}^2-\tfrac{1}{4}G_{5,0}^2-\tfrac{36}{25}G_{4,0}^2+\tfrac{17}{60}G_{3,2}+\tfrac{14}{5}G_{4,0}-\tfrac{17}{6}G_{5,0})e_1+(\tfrac{10}{9}+\tfrac{3}{25}G_{3,2}G_{4,0}-\tfrac{1}{20}G_{3,2}G_{5,0}-\tfrac{6}{5}G_{4,0}G_{5,0}+\tfrac{1}{400}G_{3,2}^2+\tfrac{1}{4}G_{5,0}^2+\tfrac{36}{25}G_{4,0}^2-\tfrac{17}{60}G_{3,2}-\tfrac{14}{5}G_{4,0}+\tfrac{17}{6}G_{5,0})e_2$
\end{tabular}}
\\
e_2 &
\rotatebox[origin=c]{0}{
\begin{tabular}{p{7cm}}
$-(\tfrac{-10}{9}-\tfrac{3}{25}G_{3,2}G_{4,0}+\tfrac{1}{20}G_{3,2}G_{5,0}+\tfrac{6}{5}G_{4,0}G_{5,0}-\tfrac{1}{400}G_{3,2}^2-\tfrac{1}{4}G_{5,0}^2-\tfrac{36}{25}G_{4,0}^2+\tfrac{17}{60}G_{3,2}+\tfrac{14}{5}G_{4,0}-\tfrac{17}{6}G_{5,0})e_1-(\tfrac{10}{9}+\tfrac{3}{25}G_{3,2}G_{4,0}-\tfrac{1}{20}G_{3,2}G_{5,0}-\tfrac{6}{5}G_{4,0}G_{5,0}+\tfrac{1}{400}G_{3,2}^2+\tfrac{1}{4}G_{5,0}^2+\tfrac{36}{25}G_{4,0}^2-\tfrac{17}{60}G_{3,2}-\tfrac{14}{5}G_{4,0}+\tfrac{17}{6}G_{5,0})e_2$
\end{tabular}}
&
0
\end{array}
\]

%%%%%%%%%%%%%%%%%%%%%%%%%%%%%%%%%%%%%%%%%%%%%%%%%%%%%%%%%%%%%%%%%%%%%%
%%%%%%%%%%%%%%%%%%%%%%%%%%%%%%%%%%%%%%%%%%%%%%%%%%%%%%%%%%%%%%%%%%%%%%
%%%%%%%%%%%%%%%%%%%%%%%%%%%%%%%%%%%%%%%%%%%%%%%%%%%%%%%%%%%%%%%%%%%%%%

\[
\!\!\!\!\!
\text{\bf Model 2f3f4c}
\ \ \ \ \
\left\{
\aligned
u
&
\,=\,
xy+x^3+y^3
+
(\tfrac{2}{3}+\tfrac{1}{2}F_{3,1}-\tfrac{1}{2}G_{4,0})x^4
+
F_{3,1}x^3y
+
(5+F_{3,1}+G_{4,0})x^2y^2
+
\\
&
\,\,\,\,\,\,\,\,
+
F_{3,1}xy^3
+
(\tfrac{2}{3}+\tfrac{1}{2}F_{3,1}-\tfrac{1}{2}G_{4,0})y^4
+
(\tfrac{11}{15}+\tfrac{34}{15}F_{3,1}+\tfrac{13}{15}G_{4,0}+\tfrac{1}{5}F_{3,1}^2
+
\\
&
\,\,\,\,\,\,\,\,
-\tfrac{2}{5}F_{3,1}G_{4,0}-\tfrac{1}{2}G_{5,0})x^5
+
(-5-5F_{3,1}+6G_{4,0}+\tfrac{5}{2}G_{5,0})x^4y
+
(30+\tfrac{80}{3}F_{3,1}
+
\\
&
\,\,\,\,\,\,\,\,
-\tfrac{34}{3}G_{4,0}+2F_{3,1}^2-4F_{3,1}G_{4,0}-5G_{5,0})x^3y^2
+
(-\tfrac{52}{3}-6F_{3,1}+5G_{5,0})x^2y^3
+
\\
&
\,\,\,\,\,\,\,\,
+
(\tfrac{56}{3}+\tfrac{34}{3}F_{3,1}+\tfrac{1}{3}G_{4,0}+F_{3,1}^2-2F_{3,1}G_{4,0}-\tfrac{5}{2}G_{5,0})xy^4
+
(-4-F_{3,1}+2G_{4,0}
+
\\
&
\,\,\,\,\,\,\,\,
+\tfrac{1}{2}G_{5,0})y^5
+
\cdots
,
\\
v
&
\,=\,
x^2+y^2+2x^2y+2xy^2+G_{4,0}x^4
+
(2F_{3,1}+\tfrac{16}{3})x^3y
+
(4F_{3,1}+2-2G_{4,0})x^2y^2
+
\\
&
\,\,\,\,\,\,\,\,
+
(2F_{3,1}+\tfrac{16}{3})xy^3
+
G_{4,0}y^4+G_{5,0}x^5
+
(\tfrac{88}{3}+\tfrac{86}{3}F_{3,1}-\tfrac{34}{3}G_{4,0}+2F_{3,1}^2
+
\\
&
\,\,\,\,\,\,\,\,
-4F_{3,1}G_{4,0}-5G_{5,0})x^4y
+
(-20F_{3,1}-\tfrac{104}{3}+16G_{4,0}+10G_{5,0})x^3y^2
+
(60
+
\\
&
\,\,\,\,\,\,\,\,
+\tfrac{136}{3}F_{3,1}-\tfrac{20}{3}G_{4,0}+4F_{3,1}^2-8F_{3,1}G_{4,0}-10G_{5,0})x^2y^3
+
(-18-4F_{3,1}
+
\\
&
\,\,\,\,\,\,\,\,
+5G_{5,0})xy^4
+
(\tfrac{142}{15}+\tfrac{98}{15}F_{3,1}-\tfrac{34}{15}G_{4,0}+\tfrac{2}{5}F_{3,1}^2-\tfrac{4}{5}F_{3,1}G_{4,0}-G_{5,0})y^5
+
\cdots
,
\endaligned\right.
\]
for any value for $G_{5,0}$.
%%%%%%%%%%%%%%%%%%%%%%%%%%%%%%%%%%%%%%%%%%%%%%%%%%%%%%%%%%%%%%%%%%%%%%
\[
\aligned
e_1 & := 
-
(-1-3vF_{3,1}G_{4,0}+3xF_{3,1}G_{4,0}+6F_{3,1}G_{4,0}u-3yF_{3,1}G_{4,0}+\tfrac{142}{3}u-25y+25x-27v
+
\\
&
\,\,\,\,\,\,\,\,
-\tfrac{15}{2}yG_{4,0}-6xG_{4,0}^2-\tfrac{11}{2}vG_{4,0}+\tfrac{29}{2}xF_{3,1}-\tfrac{15}{2}xG_{5,0}+\tfrac{11}{2}xG_{4,0}+\tfrac{15}{2}vG_{5,0}+6vG_{4,0}^2-\tfrac{27}{2}yF_{3,1}
+
\\
&
\,\,\,\,\,\,\,\,
+\tfrac{15}{2}yG_{5,0}-12G_{4,0}^2u+28F_{3,1}u+17G_{4,0}u-15G_{5,0}u+6yG_{4,0}^2-\tfrac{27}{2}vF_{3,1})\partial_x
-
(\tfrac{11}{2}yG_{4,0}
+
\\
&
\,\,\,\,\,\,\,\,
+6xG_{4,0}^2-\tfrac{15}{2}vG_{4,0}-\tfrac{27}{2}xF_{3,1}+\tfrac{15}{2}xG_{5,0}-\tfrac{15}{2}xG_{4,0}+\tfrac{15}{2}vG_{5,0}+6vG_{4,0}^2+\tfrac{29}{2}yF_{3,1}-\tfrac{15}{2}yG_{5,0}
+
\\
&
\,\,\,\,\,\,\,\,
-12G_{4,0}^2u+30F_{3,1}u+13G_{4,0}u-15G_{5,0}u-6yG_{4,0}^2-\tfrac{25}{2}vF_{3,1}+52u+26y-22x-\tfrac{64}{3}v
+
\\
&
\,\,\,\,\,\,\,\,
-3vF_{3,1}G_{4,0}+6F_{3,1}G_{4,0}u-3xF_{3,1}G_{4,0}+3yF_{3,1}G_{4,0})\partial_y
-
(-25v+51u-\tfrac{27}{2}vF_{3,1}-\tfrac{15}{2}vG_{4,0}
+
\\
&
\,\,\,\,\,\,\,\,
+\tfrac{15}{2}vG_{5,0}+6vG_{4,0}^2-12G_{4,0}^2u+29F_{3,1}u+11G_{4,0}u-15G_{5,0}u-3vF_{3,1}G_{4,0}+6F_{3,1}G_{4,0}u
+
\\
&
\,\,\,\,\,\,\,\,
-y)\partial_u
+
(12F_{3,1}G_{4,0}u-6F_{3,1}G_{4,0}v-24G_{4,0}^2u+12G_{4,0}^2v+54F_{3,1}u-29F_{3,1}v+30G_{4,0}u
+
\\
&
\,\,\,\,\,\,\,\,
-11G_{4,0}v-30G_{5,0}u+15G_{5,0}v+98u-50v+2x)\partial_v,\\
e_2 & := 
(-3yF_{3,1}^2+3xF_{3,1}^2+6F_{3,1}^2u-3vF_{3,1}^2+9vF_{3,1}G_{4,0}-9xF_{3,1}G_{4,0}-18F_{3,1}G_{4,0}u+9yF_{3,1}G_{4,0}
+
\\
&
\,\,\,\,\,\,\,\,
+90u-49y+45x-\tfrac{149}{3}v+\tfrac{49}{2}yG_{4,0}+6xG_{4,0}^2+\tfrac{49}{2}vG_{4,0}+\tfrac{69}{2}xF_{3,1}-\tfrac{15}{2}xG_{5,0}-\tfrac{45}{2}xG_{4,0}
+
\\
&
\,\,\,\,\,\,\,\,
+\tfrac{15}{2}vG_{5,0}-6vG_{4,0}^2-\tfrac{71}{2}yF_{3,1}+\tfrac{15}{2}yG_{5,0}+12G_{4,0}^2u+68F_{3,1}u-47G_{4,0}u-15G_{5,0}u-6yG_{4,0}^2
+
\\
&
\,\,\,\,\,\,\,\,
-\tfrac{73}{2}vF_{3,1})\partial_x
+
(1+3yF_{3,1}^2-3xF_{3,1}^2+6F_{3,1}^2u-3vF_{3,1}^2+9vF_{3,1}G_{4,0}+9xF_{3,1}G_{4,0}
+
\\
&
\,\,\,\,\,\,\,\,
-18F_{3,1}G_{4,0}u-9yF_{3,1}G_{4,0}+\tfrac{284}{3}u+46y-46x-44v-\tfrac{45}{2}yG_{4,0}-6xG_{4,0}^2+\tfrac{45}{2}vG_{4,0}
+
\\
&
\,\,\,\,\,\,\,\,
-\tfrac{71}{2}xF_{3,1}+\tfrac{15}{2}xG_{5,0}+\tfrac{49}{2}xG_{4,0}+\tfrac{15}{2}vG_{5,0}-6vG_{4,0}^2+\tfrac{69}{2}yF_{3,1}-\tfrac{15}{2}yG_{5,0}+12G_{4,0}^2u+70F_{3,1}u
+
\\
&
\,\,\,\,\,\,\,\,
-51G_{4,0}u-15G_{5,0}u+6yG_{4,0}^2-\tfrac{71}{2}vF_{3,1})\partial_y
+
(-46v+91u-\tfrac{71}{2}vF_{3,1}+\tfrac{49}{2}vG_{4,0}+\tfrac{15}{2}vG_{5,0}
+
\\
&
\,\,\,\,\,\,\,\,
-3vF_{3,1}^2-6vG_{4,0}^2+6F_{3,1}^2u+12G_{4,0}^2u+69F_{3,1}u-45G_{4,0}u-15G_{5,0}u+9vF_{3,1}G_{4,0}
+
\\
&
\,\,\,\,\,\,\,\,
-18F_{3,1}G_{4,0}u+x)\partial_u
-
(12F_{3,1}^2u-6F_{3,1}^2v-36F_{3,1}G_{4,0}u+18F_{3,1}G_{4,0}v+24G_{4,0}^2u
+
\\
&
\,\,\,\,\,\,\,\,
-12G_{4,0}^2v+142F_{3,1}u-69F_{3,1}v-98G_{4,0}u+45G_{4,0}v-30G_{5,0}u+15G_{5,0}v+186u
+
\\
&
\,\,\,\,\,\,\,\,
-92v-2y)\partial_v,
\endaligned
\]
with $F_{3,1},G_{4,0}$ satisfying:
\[
\aligned
0 = 3F_{3,1}^2-12F_{3,1}G_{4,0}+12G_{4,0}^2+21F_{3,1}-30G_{4,0}+20,
\endaligned
\]
%%%%%%%%%%%%%%%%%%%%%%%%%%%%%%%%%%%%%%%%%%%%%%%%%%%%%%%%%%%%%%%%%%%%%%
\[
\footnotesize
\def\arraystretch{1.25}
\begin{array}{c|cc}
{} & e_1 & e_2 
\\
\hline
e_1 &
0
&
\rotatebox[origin=c]{0}{
\begin{tabular}{p{7cm}}
$(3F_{3,1}^2-12F_{3,1}G_{4,0}+12G_{4,0}^2+21F_{3,1}-30G_{4,0}+20)e_1+(-3F_{3,1}^2+12F_{3,1}G_{4,0}-12G_{4,0}^2-21F_{3,1}+30G_{4,0}-20)e_2$
\end{tabular}}
\\
e_2 &
\rotatebox[origin=c]{0}{
\begin{tabular}{p{7cm}}
$-(3F_{3,1}^2-12F_{3,1}G_{4,0}+12G_{4,0}^2+21F_{3,1}-30G_{4,0}+20)e_1-(-3F_{3,1}^2+12F_{3,1}G_{4,0}-12G_{4,0}^2-21F_{3,1}+30G_{4,0}-20)e_2$
\end{tabular}}
&
0
\end{array}
\]
%%%%%%%%%%%%%%%%%%%%%%%%%%%%%%%%%%%%%%%%%%%%%%%%%%%%%%%%%%%%%%%%%%%%%%
%%%%%%%%%%%%%%%%%%%%%%%%%%%%%%%%%%%%%%%%%%%%%%%%%%%%%%%%%%%%%%%%%%%%%%
%%%%%%%%%%%%%%%%%%%%%%%%%%%%%%%%%%%%%%%%%%%%%%%%%%%%%%%%%%%%%%%%%%%%%%

\[
\text{\bf Model 2f3f4da}
\ \ \ \ \
\left\{
\aligned
u
&
\,=\,
yx+x^3+y^3-\tfrac{1}{16}x^4-\tfrac{1}{4}yx^3
+
\tfrac{41}{8}y^2x^2+\tfrac{3}{4}y^3x-\tfrac{9}{16}y^4-\tfrac{1}{8}x^5+8yx^4
+
\\
&
\,\,\,\,\,\,\,\,
+\tfrac{11}{4}y^2x^3-5y^3x^2+\tfrac{43}{8}y^4x+3y^5+\tfrac{249}{64}x^6+\tfrac{89}{32}yx^5-\tfrac{921}{64}y^2x^4+\tfrac{495}{16}y^3x^3
+
\\
&
\,\,\,\,\,\,\,\,
+\tfrac{1879}{64}y^4x^2-\tfrac{167}{32}y^5x-\tfrac{343}{64}y^6
+
\cdots
,
\\
v
&
\,=\,x^2+y^2+2yx^2+2y^2x+\tfrac{9}{8}x^4+\tfrac{9}{2}yx^3-\tfrac{1}{4}y^2x^2+\tfrac{5}{2}y^3x+\tfrac{17}{8}y^4+\tfrac{5}{2}x^5
+
\\
&
\,\,\,\,\,\,\,\,
-\tfrac{3}{4}yx^4+7y^2x^3+\tfrac{37}{2}y^3x^2+\tfrac{5}{2}y^4x-\tfrac{7}{4}y^5-\tfrac{17}{32}x^6+\tfrac{111}{16}yx^5+\tfrac{1505}{32}y^2x^4
+
\\
&
\,\,\,\,\,\,\,\,
+\tfrac{89}{8}y^3x^3-\tfrac{287}{32}y^4x^2+\tfrac{271}{16}y^5x+\tfrac{367}{32}y^6
+
\cdots,
\endaligned\right.
\]
%%%%%%%%%%%%%%%%%%%%%%%%%%%%%%%%%%%%%%%%%%%%%%%%%%%%%%%%%%%%%%%%%%%%%%
\[
\aligned
e_1 & := 
(\tfrac{1}{2}y+\tfrac{1}{4}v+1+x)\partial_x
-
(7v\tfrac{1}{4}+5x\tfrac{1}{2})\partial_y
+
(y+u+\tfrac{1}{2}v)\partial_u
+
(2x+2v)\partial_v,\\
e_2 & := 
(u-7y\tfrac{1}{2}-15v\tfrac{1}{4}+4x)\partial_x
+
(1+3u+5y-11v\tfrac{1}{4}-\tfrac{1}{2}x)\partial_y
+
\\
&
\,\,\,\,\,\,\,\,
+
(x+9u-\tfrac{1}{2}v)\partial_u
-
(-2y+4u-10v)\partial_v,
\endaligned
\]
%%%%%%%%%%%%%%%%%%%%%%%%%%%%%%%%%%%%%%%%%%%%%%%%%%%%%%%%%%%%%%%%%%%%%%
\[
\footnotesize
\def\arraystretch{1.25}
\begin{array}{c|cc}
{} & e_1 & e_2 
\\
\hline
e_1 & 
0 & \tfrac{7}{2}e_1-\tfrac{1}{2}e_2
\\
e_2 &
-\tfrac{7}{2}e_1+\tfrac{1}{2}e_2 & 0
\end{array}
\]

%%%%%%%%%%%%%%%%%%%%%%%%%%%%%%%%%%%%%%%%%%%%%%%%%%%%%%%%%%%%%%%%%%%%%%
%%%%%%%%%%%%%%%%%%%%%%%%%%%%%%%%%%%%%%%%%%%%%%%%%%%%%%%%%%%%%%%%%%%%%%
%%%%%%%%%%%%%%%%%%%%%%%%%%%%%%%%%%%%%%%%%%%%%%%%%%%%%%%%%%%%%%%%%%%%%%

\[
\text{\bf Model 2f3f4db}
\ \ \ \ \
\left\{
\aligned
u
&
\,=\,
yx+x^3+y^3+\tfrac{3}{32}x^4+\tfrac{1}{8}yx^3+\tfrac{65}{16}y^2x^2+\tfrac{9}{8}y^3x-\tfrac{13}{32}y^4-\tfrac{1}{16}x^5+\tfrac{37}{8}yx^4
+
\\
&
\,\,\,\,\,\,\,\,
+\tfrac{37}{8}y^2x^3+y^3x^2+\tfrac{15}{16}y^4x+\tfrac{23}{8}y^5+\tfrac{99}{64}x^6+\tfrac{195}{32}yx^5+\tfrac{417}{64}y^2x^4+\tfrac{81}{16}y^3x^3
+
\\
&
\,\,\,\,\,\,\,\,
+\tfrac{933}{64}y^4x^2+\tfrac{507}{32}y^5x-\tfrac{489}{64}y^6
+
\cdots
,
\\
v
&
\,=\,
x^2+y^2+2yx^2+2y^2x+\tfrac{13}{16}x^4+\tfrac{15}{4}yx^3+\tfrac{15}{8}y^2x^2+\tfrac{7}{4}y^3x+\tfrac{29}{16}y^4+\tfrac{7}{4}x^5
+
\\
&
\,\,\,\,\,\,\,\,
+\tfrac{31}{8}yx^4+6y^2x^3+\tfrac{37}{4}y^3x^2+\tfrac{37}{4}y^4x-\tfrac{17}{8}y^5+\tfrac{65}{32}x^6+\tfrac{105}{16}yx^5+\tfrac{515}{32}y^2x^4
+
\\
&
\,\,\,\,\,\,\,\,
+\tfrac{315}{8}y^3x^3+\tfrac{615}{32}y^4x^2-\tfrac{223}{16}y^5x+\tfrac{469}{32}y^6
+
\cdots,
\endaligned\right.
\]
%%%%%%%%%%%%%%%%%%%%%%%%%%%%%%%%%%%%%%%%%%%%%%%%%%%%%%%%%%%%%%%%%%%%%%
\[
\def\arraystretch{1.25}
\aligned
e_1 & := (y+\tfrac{11}{8}v+1-\tfrac{1}{2}x+\tfrac{3}{4}u)\partial_x-(\tfrac{3}{2}y+\tfrac{7}{8}v+2x+\tfrac{11}{4}u)\partial_y-(-y+2u-v)\partial_u+(2x+2u-v)\partial_v,\\
e_2 & := (\tfrac{1}{4}u-4y-\tfrac{31}{8}v+\tfrac{7}{2}x)\partial_x+(1+\tfrac{23}{4}u+\tfrac{9}{2}y-\tfrac{21}{8}v-x)\partial_y+(x+8u-v)\partial_u-(-2y+6u-9v)\partial_v.
\endaligned
\]
%%%%%%%%%%%%%%%%%%%%%%%%%%%%%%%%%%%%%%%%%%%%%%%%%%%%%%%%%%%%%%%%%%%%%%
\[
\footnotesize
\def\arraystretch{1.25}
\begin{array}{c|cc}
{} & e_1 & e_2 
\\
\hline
e_1 &
0 & \tfrac{5}{2}e_1+\tfrac{1}{2}e_2
\\
e_2 &
-\tfrac{5}{2}e_1-\tfrac{1}{2}e_2 & 0
\end{array}
\]

%%%%%%%%%%%%%%%%%%%%%%%%%%%%%%%%%%%%%%%%%%%%%%%%%%%%%%%%%%%%%%%%%%%%%%
%%%%%%%%%%%%%%%%%%%%%%%%%%%%%%%%%%%%%%%%%%%%%%%%%%%%%%%%%%%%%%%%%%%%%%
%%%%%%%%%%%%%%%%%%%%%%%%%%%%%%%%%%%%%%%%%%%%%%%%%%%%%%%%%%%%%%%%%%%%%%

\[
\text{\bf Model 2f3f4ea}
\ \ \ \ \
\left\{
\aligned
u
&
\,=\,
yx+x^3+y^3-\tfrac{9}{16}x^4+\tfrac{3}{4}yx^3+\tfrac{41}{8}y^2x^2-\tfrac{1}{4}y^3x-\tfrac{1}{16}y^4+3x^5+\tfrac{43}{8}yx^4
+
\\
&
\,\,\,\,\,\,\,\,
-5y^2x^3+\tfrac{11}{4}y^3x^2+8y^4x-\tfrac{1}{8}y^5-\tfrac{343}{64}x^6-\tfrac{167}{32}yx^5+\tfrac{1879}{64}y^2x^4+\tfrac{495}{16}y^3x^3
+
\\
&
\,\,\,\,\,\,\,\,
-\tfrac{921}{64}y^4x^2+\tfrac{89}{32}y^5x+\tfrac{249}{64}y^6
+
\cdots
,
\\
v
&
\,=\,
x^2+y^2+2yx^2+2y^2x+\tfrac{17}{8}x^4+\tfrac{5}{2}yx^3-\tfrac{1}{4}y^2x^2+\tfrac{9}{2}y^3x+\tfrac{9}{8}y^4-\tfrac{7}{4}x^5
+
\\
&
\,\,\,\,\,\,\,\,
+\tfrac{5}{2}yx^4+\tfrac{37}{2}y^2x^3+7y^3x^2-\tfrac{3}{4}y^4x+\tfrac{5}{2}y^5+\tfrac{367}{32}x^6+\tfrac{271}{16}yx^5-\tfrac{287}{32}y^2x^4
+
\\
&
\,\,\,\,\,\,\,\,
+\tfrac{89}{8}y^3x^3+\tfrac{1505}{32}y^4x^2+\tfrac{111}{16}y^5x-\tfrac{17}{32}y^6
+
\cdots,
\endaligned\right.
\]
%%%%%%%%%%%%%%%%%%%%%%%%%%%%%%%%%%%%%%%%%%%%%%%%%%%%%%%%%%%%%%%%%%%%%%
\[
\def\arraystretch{1.25}
\aligned
e_1 & := (5x+1-\tfrac{1}{2}y+3u-11v\tfrac{1}{4})\partial_x+(-7x\tfrac{1}{2}+4y+u-15v\tfrac{1}{4})\partial_y+(y+9u-\tfrac{1}{2}v)\partial_u
+
\\
&
\,\,\,\,\,\,\,\,
-(-2x+4u-10v)\partial_v,\\
e_2 & := -(5y\tfrac{1}{2}+7v\tfrac{1}{4})\partial_x+(1+\tfrac{1}{2}x+y+\tfrac{1}{4}v)\partial_y+(x+u+\tfrac{1}{2}v)\partial_u+(2y+2v)\partial_v,
\endaligned
\]
%%%%%%%%%%%%%%%%%%%%%%%%%%%%%%%%%%%%%%%%%%%%%%%%%%%%%%%%%%%%%%%%%%%%%%
\[
\footnotesize
\def\arraystretch{1.25}
\begin{array}{c|cc}
{} & e_1 & e_2 
\\
\hline
e_1 & 
0 & \tfrac{1}{2}e_1-\tfrac{7}{2}e_2
\\
e_2 &
-\tfrac{1}{2}e1+\tfrac{7}{2}e_2 & 0
\end{array}
\]
%%%%%%%%%%%%%%%%%%%%%%%%%%%%%%%%%%%%%%%%%%%%%%%%%%%%%%%%%%%%%%%%%%%%%%
%%%%%%%%%%%%%%%%%%%%%%%%%%%%%%%%%%%%%%%%%%%%%%%%%%%%%%%%%%%%%%%%%%%%%%
%%%%%%%%%%%%%%%%%%%%%%%%%%%%%%%%%%%%%%%%%%%%%%%%%%%%%%%%%%%%%%%%%%%%%%

\[
\text{\bf Model 2f3f4eb}
\ \ \ \ \
\left\{
\aligned
u
&
\,=\,
yx+x^3+y^3-\tfrac{13}{32}x^4+\tfrac{9}{8}yx^3+\tfrac{65}{16}y^2x^2+\tfrac{1}{8}y^3x+\tfrac{3}{32}y^4+\tfrac{23}{8}x^5+\tfrac{15}{16}yx^4
+
\\
&
\,\,\,\,\,\,\,\,
+y^2x^3+\tfrac{37}{8}y^3x^2+\tfrac{37}{8}y^4x-\tfrac{1}{16}y^5-\tfrac{489}{64}x^6+\tfrac{507}{32}yx^5+\tfrac{933}{64}y^2x^4+\tfrac{81}{16}y^3x^3
+
\\
&
\,\,\,\,\,\,\,\,
+\tfrac{417}{64}y^4x^2+\tfrac{195}{32}y^5x+\tfrac{99}{64}y^6
+
\cdots
,
\\
v
&
\,=\,
x^2+y^2+2yx^2+2y^2x+\tfrac{29}{16}x^4+\tfrac{7}{4}yx^3+\tfrac{15}{8}y^2x^2+\tfrac{15}{4}y^3x+\tfrac{13}{16}y^4-\tfrac{17}{8}x^5
+
\\
&
\,\,\,\,\,\,\,\,
+\tfrac{37}{4}yx^4+\tfrac{37}{4}y^2x^3+6y^3x^2+\tfrac{31}{8}y^4x+\tfrac{7}{4}y^5+\tfrac{469}{32}x^6-\tfrac{223}{16}yx^5+\tfrac{615}{32}y^2x^4
+
\\
&
\,\,\,\,\,\,\,\,
+\tfrac{315}{8}y^3x^3+\tfrac{515}{32}y^4x^2+\tfrac{105}{16}y^5x+\tfrac{65}{32}y^6
+
\cdot,
\endaligned\right.
\]
%%%%%%%%%%%%%%%%%%%%%%%%%%%%%%%%%%%%%%%%%%%%%%%%%%%%%%%%%%%%%%%%%%%%%%
\[
\def\arraystretch{1.25}
\aligned
e_1 & := (9x\tfrac{1}{2}+1-y+23u\tfrac{1}{4}-21v\tfrac{1}{8})\partial_x+(-4x+7y\tfrac{1}{2}+\tfrac{1}{4}u-31v\tfrac{1}{8})\partial_y+(y+8u-v)\partial_u
+
\\
&
\,\,\,\,\,\,\,\,-(-2x+6u-9v)\partial_v,\\
e_2 & := -(3x\tfrac{1}{2}+2y+11u\tfrac{1}{4}+7v\tfrac{1}{8})\partial_x+(1+x-\tfrac{1}{2}y+3u\tfrac{1}{4}+11v\tfrac{1}{8})\partial_y-(-x+2u-v)\partial_u
+
\\
&
\,\,\,\,\,\,\,\,+(2y+2u-v)\partial_v,
\endaligned
\]
%%%%%%%%%%%%%%%%%%%%%%%%%%%%%%%%%%%%%%%%%%%%%%%%%%%%%%%%%%%%%%%%%%%%%%
\[
\footnotesize
\def\arraystretch{1.25}
\begin{array}{c|cc}
{} & e_1 & e_2 
\\
\hline
e_1 & 
0 & -\tfrac{1}{2}e_1-\tfrac{5}{2}e_2
\\
e_2 &
\tfrac{1}{2}e_1+\tfrac{5}{2}e_2 & 0
\end{array}
\]

%%%%%%%%%%%%%%%%%%%%%%%%%%%%%%%%%%%%%%%%%%%%%%%%%%%%%%%%%%%%%%%%%%%%%%
%%%%%%%%%%%%%%%%%%%%%%%%%%%%%%%%%%%%%%%%%%%%%%%%%%%%%%%%%%%%%%%%%%%%%%
%%%%%%%%%%%%%%%%%%%%%%%%%%%%%%%%%%%%%%%%%%%%%%%%%%%%%%%%%%%%%%%%%%%%%%

\[
\text{\bf Model 2f3g}
\ \ \ \ \
\left\{
\aligned
u
&
\,=\,
xy+F_{0,3}y^3
+
(2-F_{0,3})x^3
+
F_{0,4}y^4
+
F_{1,3}xy^3
+
F_{2,2}x^2y^2
+
(\tfrac{26}{3}
+
\\
&
\ \ \ \ \
+\tfrac{41}{6}F_{0,3}^2+\tfrac{40}{9}F_{0,3}+\tfrac{1}{2}F_{1,3}+\tfrac{5}{6}F_{2,2}-\tfrac{5}{8}G_{1,3}-\tfrac{1}{6}G_{2,2}+\tfrac{3}{8}G_{3,1}-3F_{0,4}
+
\\
&
\ \ \ \ \
+2G_{0,4}+\tfrac{23}{2}F_{0,3}^3+\tfrac{1}{3}F_{0,3}F_{1,3}-\tfrac{11}{18}F_{0,3}F_{2,2}+\tfrac{11}{24}F_{0,3}G_{1,3}-\tfrac{1}{9}F_{0,3}G_{2,2}
+
\\
&
\ \ \ \ \
-\tfrac{5}{8}F_{0,3}G_{3,1}+5F_{0,3}F_{0,4})x^3y
+
(\tfrac{99}{8}-\tfrac{253}{96}F_{0,3}^4+\tfrac{14753}{864}F_{0,3}^2-\tfrac{187}{72}F_{0,3}^3
+
\\
&
\ \ \ \ \
+\tfrac{55}{192}F_{0,3}G_{1,3}-\tfrac{11}{24}F_{0,3}G_{0,4}-\tfrac{11}{144}F_{0,3}^2F_{1,3}+\tfrac{121}{864}F_{0,3}^2F_{2,2}-\tfrac{121}{1152}F_{0,3}^2G_{1,3}
+
\\
&
\ \ \ \ \
+\tfrac{11}{432}F_{0,3}^2G_{2,2}+\tfrac{55}{384}F_{0,3}^2G_{3,1}-\tfrac{55}{48}F_{0,3}^2F_{0,4}+\tfrac{29}{96}F_{0,3}F_{1,3}-\tfrac{47}{144}F_{0,3}F_{2,2}
+
\\
&
\ \ \ \ \
+\tfrac{9}{4}F_{0,3}F_{0,4}+\tfrac{23}{288}F_{0,3}G_{2,2}-\tfrac{31}{96}F_{0,3}G_{3,1}-\tfrac{5}{8}F_{0,3}G_{4,0}-\tfrac{1361}{72}F_{0,3}-\tfrac{15}{16}F_{0,4}
+
\\
&
\ \ \ \ \
-\tfrac{5}{32}G_{2,2}+\tfrac{7}{128}G_{3,1}-\tfrac{1}{32}F_{1,3}+\tfrac{11}{32}F_{2,2}+\tfrac{5}{8}G_{0,4}-\tfrac{9}{128}G_{1,3}+\tfrac{3}{8}G_{4,0})x^4
+
\cdots
,
\\
v
&
\,=\,
y^2
+
x^2
+
(-8-6F_{0,3})xy^2
+
(-4+6F_{0,3})x^2y
+
G_{0,4}y^4
+
G_{1,3}xy^3
+
\\
&
\ \ \ \ \
+
G_{2,2}x^2y^2
+
G_{3,1}x^3y
+
G_{4,0}x^4
+
\cdots
.
\endaligned\right.
\]
%%%%%%%%%%%%%%%%%%%%%%%%%%%%%%%%%%%%%%%%%%%%%%%%%%%%%%%%%%%%%%%%%%%%%%
\[
\def\arraystretch{1.25}
\begin{array}{ll}
e_1
&
\,:=\,
-(-1+3uG_{4,0}-\tfrac{1}{4}uG_{2,2}+\tfrac{11}{16}uG_{3,1}-\tfrac{233}{18}uF_{0,3}+\tfrac{7}{4}uF_{2,2}-\tfrac{221}{144}yF_{0,3}^2-\tfrac{1}{12}yG_{0,4}
+
\\
&
\ \ \ \ \
-\tfrac{23}{48}yF_{0,3}^3-\tfrac{13}{48}xF_{1,3}-\tfrac{85}{192}xG_{1,3}+\tfrac{1}{8}yF_{0,4}+\tfrac{43}{192}xG_{3,1}+\tfrac{5}{4}xG_{4,0}-\tfrac{17}{48}xG_{2,2}-\tfrac{1297}{54}xF_{0,3}
+
\\
&
\ \ \ \ \
+\tfrac{305}{54}yF_{0,3}+\tfrac{115}{48}vF_{0,3}^4+\tfrac{1}{48}yF_{2,2}+\tfrac{5}{48}yF_{1,3}-\tfrac{17}{8}xF_{0,4}+\tfrac{17}{12}xG_{0,4}+\tfrac{391}{48}xF_{0,3}^3+\tfrac{2461}{144}xF_{0,3}^2
+
\\
&
\ \ \ \ \
+\tfrac{31}{48}xF_{2,2}-\tfrac{1}{4}yG_{4,0}+uG_{0,4}-\tfrac{11}{192}yG_{3,1}-\tfrac{25}{192}vF_{0,3}^2G_{3,1}-\tfrac{11}{576}yF_{0,3}G_{1,3}-\tfrac{25}{72}uF_{0,3}G_{2,2}
+
\\
&
\ \ \ \ \
-\tfrac{43}{72}uF_{0,3}F_{2,2}+\tfrac{85}{24}xF_{0,3}F_{0,4}-\tfrac{85}{192}xF_{0,3}G_{3,1}-\tfrac{47}{192}vF_{0,3}G_{1,3}+\tfrac{3}{4}uF_{0,3}F_{0,4}+\tfrac{5}{192}yF_{0,3}G_{3,1}
+
\\
&
\ \ \ \ \
-\tfrac{17}{216}xF_{0,3}G_{2,2}+\tfrac{1}{216}yF_{0,3}G_{2,2}-\tfrac{29}{48}vF_{0,3}F_{1,3}-\tfrac{55}{432}vF_{0,3}^2F_{2,2}+\tfrac{25}{24}vF_{0,3}^2F_{0,4}+\tfrac{5}{4}vF_{0,3}G_{4,0}
+
\\
&
\ \ \ \ \
+\tfrac{7}{6}uF_{0,3}G_{0,4}+\tfrac{85}{192}vF_{0,3}G_{3,1}+\tfrac{187}{576}xF_{0,3}G_{1,3}+\tfrac{5}{12}vF_{0,3}G_{0,4}-\tfrac{1}{72}yF_{0,3}F_{1,3}-\tfrac{5}{24}yF_{0,3}F_{0,4}
+
\\
&
\ \ \ \ \
+\tfrac{161}{24}uF_{0,3}^4+\tfrac{5}{192}yG_{1,3}+\tfrac{1}{48}yG_{2,2}-\tfrac{4607}{108}uF_{0,3}^2+\tfrac{1961}{72}uF_{0,3}^3+\tfrac{7}{36}uF_{0,3}^2F_{1,3}+\tfrac{35}{12}uF_{0,3}^2F_{0,4}
+
\\
&
\ \ \ \ \
-\tfrac{35}{96}uF_{0,3}^2G_{3,1}+\tfrac{7}{2}uF_{0,3}G_{4,0}-\tfrac{3}{2}uF_{0,4}+\tfrac{5}{32}vG_{1,3}+\tfrac{5}{8}vF_{1,3}+\tfrac{1}{2}vG_{4,0}-\tfrac{11}{32}vG_{3,1}
+
\\
&
\ \ \ \ \
-\tfrac{77}{216}uF_{0,3}^2F_{2,2}-\tfrac{5}{216}vF_{0,3}^2G_{2,2}-\tfrac{7}{108}uF_{0,3}^2G_{2,2}+\tfrac{55}{576}vF_{0,3}^2G_{1,3}-\tfrac{187}{432}xF_{0,3}F_{2,2}
+
\\
&
\ \ \ \ \
+\tfrac{11}{432}yF_{0,3}F_{2,2}+\tfrac{47}{96}uF_{0,3}G_{3,1}-\tfrac{11}{144}vF_{0,3}G_{2,2}-\tfrac{31}{24}uF_{0,3}F_{1,3}-\tfrac{13}{96}uF_{0,3}G_{1,3}
+
\\
&
\ \ \ \ \
+\tfrac{7}{144}vF_{0,3}F_{2,2}+\tfrac{77}{288}uF_{0,3}^2G_{1,3}+\tfrac{5}{72}vF_{0,3}^2F_{1,3}+\tfrac{17}{72}xF_{0,3}F_{1,3}-\tfrac{15}{8}vF_{0,3}F_{0,4}+71u
+
\\
&
\ \ \ \ \
-\tfrac{71}{12}y+\tfrac{391}{12}x-\tfrac{7}{2}v-\tfrac{5}{4}uF_{1,3}-\tfrac{13}{16}uG_{1,3}+\tfrac{3}{4}vF_{0,4}-\tfrac{1}{2}vG_{0,4}+\tfrac{691}{144}vF_{0,3}^3-\tfrac{4201}{216}vF_{0,3}^2
+
\\
&
\ \ \ \ \
+\tfrac{413}{36}vF_{0,3}+\tfrac{1}{8}vF_{2,2}+\tfrac{1}{8}vG_{2,2})\partial_x
+
\\
&
\ \ \ \ \
+
(-2uG_{4,0}+\tfrac{5}{6}uG_{2,2}-\tfrac{49}{24}uG_{3,1}+\tfrac{5825}{54}uF_{0,3}-\tfrac{13}{6}uF_{2,2}-\tfrac{2461}{144}yF_{0,3}^2-\tfrac{17}{12}yG_{0,4}
+
\\
&
\ \ \ \ \
-\tfrac{391}{48}yF_{0,3}^3-\tfrac{5}{48}xF_{1,3}-\tfrac{5}{192}xG_{1,3}+\tfrac{17}{8}yF_{0,4}+\tfrac{11}{192}xG_{3,1}+\tfrac{1}{4}xG_{4,0}-\tfrac{1}{48}xG_{2,2}-\tfrac{143}{54}xF_{0,3}
+
\\
&
\ \ \ \ \
+\tfrac{1459}{54}yF_{0,3}+\tfrac{115}{48}vF_{0,3}^4-\tfrac{31}{48}yF_{2,2}+\tfrac{13}{48}yF_{1,3}-\tfrac{1}{8}xF_{0,4}+\tfrac{1}{12}xG_{0,4}+\tfrac{23}{48}xF_{0,3}^3+\tfrac{221}{144}xF_{0,3}^2
+
\\
&
\ \ \ \ \
-\tfrac{1}{48}xF_{2,2}-\tfrac{5}{4}yG_{4,0}-\tfrac{22}{3}uG_{0,4}-\tfrac{43}{192}yG_{3,1}-\tfrac{25}{192}vF_{0,3}^2G_{3,1}-\tfrac{187}{576}yF_{0,3}G_{1,3}
+
\\
&
\ \ \ \ \
+\tfrac{25}{216}uF_{0,3}G_{2,2}+\tfrac{421}{216}uF_{0,3}F_{2,2}+\tfrac{5}{24}xF_{0,3}F_{0,4}-\tfrac{5}{192}xF_{0,3}G_{3,1}-\tfrac{31}{576}vF_{0,3}G_{1,3}
+
\\
&
\ \ \ \ \
-\tfrac{241}{12}uF_{0,3}F_{0,4}+\tfrac{85}{192}yF_{0,3}G_{3,1}-\tfrac{1}{216}xF_{0,3}G_{2,2}+\tfrac{17}{216}yF_{0,3}G_{2,2}-\tfrac{67}{144}vF_{0,3}F_{1,3}
+
\\
&
\ \ \ \ \
-\tfrac{55}{432}vF_{0,3}^2F_{2,2}+\tfrac{25}{24}vF_{0,3}^2F_{0,4}+\tfrac{5}{4}vF_{0,3}G_{4,0}+\tfrac{7}{6}uF_{0,3}G_{0,4}+\tfrac{35}{192}vF_{0,3}G_{3,1}+\tfrac{11}{576}xF_{0,3}G_{1,3}
+
\\
&
\ \ \ \ \
+\tfrac{5}{12}vF_{0,3}G_{0,4}-\tfrac{17}{72}yF_{0,3}F_{1,3}-\tfrac{85}{24}yF_{0,3}F_{0,4}+\tfrac{161}{24}uF_{0,3}^4+\tfrac{85}{192}yG_{1,3}+\tfrac{17}{48}yG_{2,2}-\tfrac{5731}{54}uF_{0,3}^2
+
\end{array}
\]

\[
\def\arraystretch{1.25}
\begin{array}{ll}
&
\ \ \ \ \
-\tfrac{1489}{72}uF_{0,3}^3+\tfrac{7}{36}uF_{0,3}^2F_{1,3}+\tfrac{35}{12}uF_{0,3}^2F_{0,4}-\tfrac{35}{96}uF_{0,3}^2G_{3,1}+\tfrac{7}{2}uF_{0,3}G_{4,0}+11uF_{0,4}
+
\\
&
\ \ \ \ \
-\tfrac{29}{48}vG_{1,3}-\tfrac{5}{12}vF_{1,3}+vG_{4,0}+\tfrac{11}{48}vG_{3,1}-\tfrac{77}{216}uF_{0,3}^2F_{2,2}-\tfrac{5}{216}vF_{0,3}^2G_{2,2}
+
\\
&
\ \ \ \ \
-\tfrac{7}{108}uF_{0,3}^2G_{2,2}+\tfrac{55}{576}vF_{0,3}^2G_{1,3}-\tfrac{11}{432}xF_{0,3}F_{2,2}+\tfrac{187}{432}yF_{0,3}F_{2,2}+\tfrac{99}{32}uF_{0,3}G_{3,1}
+
\\
&
\ \ \ \ \
-\tfrac{53}{432}vF_{0,3}G_{2,2}-\tfrac{193}{72}uF_{0,3}F_{1,3}-\tfrac{589}{288}uF_{0,3}G_{1,3}-\tfrac{89}{432}vF_{0,3}F_{2,2}+\tfrac{77}{288}uF_{0,3}^2G_{1,3}
+
\\
&
\ \ \ \ \
+\tfrac{5}{72}vF_{0,3}^2F_{1,3}+\tfrac{1}{72}xF_{0,3}F_{1,3}+\tfrac{5}{24}vF_{0,3}F_{0,4}-\tfrac{266}{3}u-\tfrac{343}{12}y-\tfrac{1}{12}x+\tfrac{71}{3}v+\tfrac{1}{6}uF_{1,3}
+
\\
&
\ \ \ \ \
+\tfrac{55}{24}uG_{1,3}-\tfrac{1}{2}vF_{0,4}+\tfrac{1}{3}vG_{0,4}+\tfrac{1381}{144}vF_{0,3}^3-\tfrac{1415}{108}vF_{0,3}^2-\tfrac{325}{108}vF_{0,3}-\tfrac{1}{12}vF_{2,2}
+
\\
&
\ \ \ \ \
-\tfrac{1}{12}vG_{2,2})\partial_y
+
\\
&
\ \ \ \ \
-
(\tfrac{5}{2}uG_{4,0}-\tfrac{17}{24}uG_{2,2}+\tfrac{43}{96}uG_{3,1}-\tfrac{1378}{27}uF_{0,3}+\tfrac{31}{24}uF_{2,2}+\tfrac{17}{6}uG_{0,4}-\tfrac{17}{108}uF_{0,3}G_{2,2}
+
\\
&
\ \ \ \ \
-\tfrac{187}{216}uF_{0,3}F_{2,2}-\tfrac{11}{576}vF_{0,3}G_{1,3}+\tfrac{85}{12}uF_{0,3}F_{0,4}-\tfrac{1}{72}vF_{0,3}F_{1,3}+\tfrac{5}{192}vF_{0,3}G_{3,1}+\tfrac{2461}{72}uF_{0,3}^2
+
\\
&
\ \ \ \ \
+\tfrac{391}{24}uF_{0,3}^3-\tfrac{17}{4}uF_{0,4}+\tfrac{5}{192}vG_{1,3}+\tfrac{5}{48}vF_{1,3}-\tfrac{1}{4}vG_{4,0}-\tfrac{11}{192}vG_{3,1}-\tfrac{85}{96}uF_{0,3}G_{3,1}
+
\\
&
\ \ \ \ \
+\tfrac{1}{216}vF_{0,3}G_{2,2}+\tfrac{17}{36}uF_{0,3}F_{1,3}+\tfrac{187}{288}uF_{0,3}G_{1,3}+\tfrac{11}{432}vF_{0,3}F_{2,2}-\tfrac{5}{24}vF_{0,3}F_{0,4}+\tfrac{367}{6}u
+
\\
&
\ \ \ \ \
-y-\tfrac{71}{12}v-\tfrac{13}{24}uF_{1,3}-\tfrac{85}{96}uG_{1,3}+\tfrac{1}{8}vF_{0,4}-\tfrac{1}{12}vG_{0,4}-\tfrac{23}{48}vF_{0,3}^3-\tfrac{221}{144}vF_{0,3}^2+\tfrac{305}{54}vF_{0,3}
+
\\
&
\ \ \ \ \
+\tfrac{1}{48}vF_{2,2}+\tfrac{1}{48}vG_{2,2})\partial_u
+
(uG_{4,0}-\tfrac{1}{12}uG_{2,2}+\tfrac{11}{48}uG_{3,1}-\tfrac{124}{27}uF_{0,3}-\tfrac{1}{12}uF_{2,2}+\tfrac{1}{3}uG_{0,4}
+
\\
&
\ \ \ \ \
-\tfrac{1}{54}uF_{0,3}G_{2,2}-\tfrac{11}{108}uF_{0,3}F_{2,2}-\tfrac{187}{288}vF_{0,3}G_{1,3}+\tfrac{5}{6}uF_{0,3}F_{0,4}-\tfrac{17}{36}vF_{0,3}F_{1,3}+\tfrac{85}{96}vF_{0,3}G_{3,1}
+
\\
&
\ \ \ \ \
+\tfrac{221}{36}uF_{0,3}^2+\tfrac{23}{12}uF_{0,3}^3-\tfrac{1}{2}uF_{0,4}+\tfrac{85}{96}vG_{1,3}+\tfrac{13}{24}vF_{1,3}-\tfrac{5}{2}vG_{4,0}-\tfrac{43}{96}vG_{3,1}-\tfrac{5}{48}uF_{0,3}G_{3,1}
+
\\
&
\ \ \ \ \
+\tfrac{17}{108}vF_{0,3}G_{2,2}+\tfrac{1}{18}uF_{0,3}F_{1,3}+\tfrac{11}{144}uF_{0,3}G_{1,3}+\tfrac{187}{216}vF_{0,3}F_{2,2}-\tfrac{85}{12}vF_{0,3}F_{0,4}+\tfrac{11}{3}u+2x
+
\\
&
\ \ \ \ \
-\tfrac{391}{6}v-\tfrac{5}{12}uF_{1,3}-\tfrac{5}{48}uG_{1,3}+\tfrac{17}{4}vF_{0,4}-\tfrac{17}{6}vG_{0,4}-\tfrac{391}{24}vF_{0,3}^3-\tfrac{2461}{72}vF_{0,3}^2+\tfrac{1297}{27}vF_{0,3}
+
\\
&
\ \ \ \ \
-\tfrac{31}{24}vF_{2,2}+\tfrac{17}{24}vG_{2,2})\partial_v
,
\\
e_2
&
\,:=\,
(-\tfrac{1}{3}uG_{2,2}-\tfrac{3}{4}uG_{3,1}-\tfrac{238}{9}uF_{0,3}-\tfrac{1}{3}uF_{2,2}-\tfrac{7}{4}yF_{0,3}^2-\tfrac{1}{12}xF_{1,3}+\tfrac{17}{48}xG_{1,3}-\tfrac{1}{2}yF_{0,4}
+
\\
&
\ \ \ \ \
-\tfrac{5}{16}xG_{3,1}+\tfrac{1}{36}xG_{2,2}-\tfrac{4}{3}xF_{0,3}-\tfrac{4}{3}yF_{0,3}+\tfrac{1}{36}yF_{2,2}-\tfrac{1}{12}yF_{1,3}+\tfrac{5}{2}xF_{0,4}+\tfrac{11}{4}xF_{0,3}^2
+
\\
&
\ \ \ \ \
-\tfrac{17}{36}xF_{2,2}+2uG_{0,4}+\tfrac{1}{16}yG_{3,1}-\tfrac{7}{18}uF_{0,3}G_{2,2}-\tfrac{7}{18}uF_{0,3}F_{2,2}+\tfrac{5}{48}vF_{0,3}G_{1,3}
+
\\
&
\ \ \ \ \
+7uF_{0,3}F_{0,4}+\tfrac{5}{12}vF_{0,3}F_{1,3}-\tfrac{5}{16}vF_{0,3}G_{3,1}-\tfrac{1}{48}yG_{1,3}+\tfrac{1}{36}yG_{2,2}+\tfrac{47}{3}uF_{0,3}^2+\tfrac{49}{2}uF_{0,3}^3
+
\\
&
\ \ \ \ \
+6uF_{0,4}-\tfrac{1}{8}vG_{1,3}-\tfrac{1}{2}vF_{1,3}-\tfrac{1}{8}vG_{3,1}-\tfrac{7}{8}uF_{0,3}G_{3,1}-\tfrac{5}{36}vF_{0,3}G_{2,2}+\tfrac{7}{6}uF_{0,3}F_{1,3}
+
\\
&
\ \ \ \ \
+\tfrac{7}{24}uF_{0,3}G_{1,3}-\tfrac{5}{36}vF_{0,3}F_{2,2}+\tfrac{5}{2}vF_{0,3}F_{0,4}-\tfrac{68}{3}u+\tfrac{17}{9}y-\tfrac{91}{9}x-\tfrac{14}{3}v-2uF_{1,3}+\tfrac{1}{4}uG_{1,3}
+
\\
&
\ \ \ \ \
-3vF_{0,4}+\tfrac{35}{4}vF_{0,3}^3-\tfrac{59}{6}vF_{0,3}^2-\tfrac{85}{9}vF_{0,3}+\tfrac{1}{6}vF_{2,2}+\tfrac{1}{6}vG_{2,2})\partial_x
+
\\
&
\ \ \ \ \
-
(-1+\tfrac{4}{9}uG_{2,2}+\tfrac{1}{2}uG_{3,1}-\tfrac{250}{9}uF_{0,3}+\tfrac{22}{9}uF_{2,2}-\tfrac{11}{4}yF_{0,3}^2+\tfrac{1}{12}xF_{1,3}+\tfrac{1}{48}xG_{1,3}
+
\\
&
\ \ \ \ \
-\tfrac{5}{2}yF_{0,4}-\tfrac{1}{16}xG_{3,1}-\tfrac{1}{36}xG_{2,2}-\tfrac{5}{3}xF_{0,3}-\tfrac{5}{3}yF_{0,3}+\tfrac{17}{36}yF_{2,2}+\tfrac{1}{12}yF_{1,3}+\tfrac{1}{2}xF_{0,4}
+
\\
&
\ \ \ \ \
+\tfrac{7}{4}xF_{0,3}^2-\tfrac{1}{36}xF_{2,2}+\tfrac{5}{16}yG_{3,1}-\tfrac{7}{18}uF_{0,3}G_{2,2}-\tfrac{7}{18}uF_{0,3}F_{2,2}+\tfrac{5}{48}vF_{0,3}G_{1,3}+7uF_{0,3}F_{0,4}
+
\\
&
\ \ \ \ \
+\tfrac{5}{12}vF_{0,3}F_{1,3}-\tfrac{5}{16}vF_{0,3}G_{3,1}-\tfrac{17}{48}yG_{1,3}-\tfrac{1}{36}yG_{2,2}-\tfrac{100}{3}uF_{0,3}^2+\tfrac{49}{2}uF_{0,3}^3-8uF_{0,4}
+
\\
&
\ \ \ \ \
+\tfrac{1}{12}vG_{1,3}+\tfrac{1}{3}vF_{1,3}-\tfrac{1}{4}vG_{3,1}-\tfrac{7}{8}uF_{0,3}G_{3,1}-\tfrac{5}{36}vF_{0,3}G_{2,2}+\tfrac{7}{6}uF_{0,3}F_{1,3}+\tfrac{7}{24}uF_{0,3}G_{1,3}
+
\\
&
\ \ \ \ \
-\tfrac{5}{36}vF_{0,3}F_{2,2}+\tfrac{5}{2}vF_{0,3}F_{0,4}+\tfrac{128}{9}u+\tfrac{109}{9}y-\tfrac{17}{9}x-\tfrac{68}{9}v-\tfrac{4}{3}uF_{1,3}-\tfrac{1}{3}uG_{1,3}+2vF_{0,4}
+
\\
&
\ \ \ \ \
+2vG_{0,4}+\tfrac{35}{4}vF_{0,3}^3+\tfrac{50}{3}vF_{0,3}^2+\tfrac{35}{9}vF_{0,3}-\tfrac{1}{9}vF_{2,2}-\tfrac{1}{9}vG_{2,2})\partial_y
+
\\
&
\ \ \ \ \
+
(\tfrac{11}{2}uF_{0,3}^2+\tfrac{1}{3}uF_{0,3}-\tfrac{1}{6}uF_{1,3}-\tfrac{17}{18}uF_{2,2}+\tfrac{17}{24}uG_{1,3}+\tfrac{1}{18}uG_{2,2}-\tfrac{5}{8}uG_{3,1}+5uF_{0,4}
+
\\
&
\ \ \ \ \
-\tfrac{7}{4}vF_{0,3}^2+\tfrac{5}{3}vF_{0,3}-\tfrac{1}{12}vF_{1,3}+\tfrac{1}{36}vF_{2,2}-\tfrac{1}{48}vG_{1,3}+\tfrac{1}{36}vG_{2,2}+\tfrac{1}{16}vG_{3,1}-\tfrac{1}{2}vF_{0,4}
+
\\
&
\ \ \ \ \
-\tfrac{200}{9}u+\tfrac{17}{9}v+x)\partial_u
-
(7uF_{0,3}^2+\tfrac{34}{3}uF_{0,3}+\tfrac{1}{3}uF_{1,3}-\tfrac{1}{9}uF_{2,2}+\tfrac{1}{12}uG_{1,3}-\tfrac{1}{9}uG_{2,2}
+
\\
&
\ \ \ \ \
-\tfrac{1}{4}uG_{3,1}+2uF_{0,4}-\tfrac{11}{2}vF_{0,3}^2-\tfrac{10}{3}vF_{0,3}+\tfrac{1}{6}vF_{1,3}+\tfrac{17}{18}vF_{2,2}-\tfrac{17}{24}vG_{1,3}-\tfrac{1}{18}vG_{2,2}
+
\\
&
\ \ \ \ \
+\tfrac{5}{8}vG_{3,1}-5vF_{0,4}+\tfrac{76}{9}u+\tfrac{218}{9}v-2y)\partial_v.
\end{array}
\]

%%%%%%%%%%%%%%%%%%%%%%%%%%%%%%%%%%%%%%%%%%%%%%%%%%%%%%%%%%%%%%%%%%%%%%
\noindent
A Gr\"obner basis has 92 generators to describe the 
moduli space core algebraic variety in 
$\R^7 \ni F_{0,3},F_{0,4},G_{0,4},G_{1,3},
G_{2,2},G_{3,1},G_{4,0}$.

%%%%%%%%%%%%%%%%%%%%%%%%%%%%%%%%%%%%%%%%%%%%%%%%%%%%%%%%%%%%%%%%%%%%%%
\[
\footnotesize
\def\arraystretch{1.25}
\begin{array}{c|cc}
{} & e_1 & e_2 
\\
\hline
e_1 & 0 & 0
\\
e_2 & 0 & 0
\end{array}
\]

%%%%%%%%%%%%%%%%%%%%%%%%%%%%%%%%%%%%%%%%%%%%%%%%%%%%%%%%%%%%%%%%%%%%%%
%%%%%%%%%%%%%%%%%%%%%%%%%%%%%%%%%%%%%%%%%%%%%%%%%%%%%%%%%%%%%%%%%%%%%%
%%%%%%%%%%%%%%%%%%%%%%%%%%%%%%%%%%%%%%%%%%%%%%%%%%%%%%%%%%%%%%%%%%%%%%

\[
\text{\bf Model 2f3h}
\ \ \ \ \
\left\{
\aligned
u
&
\,=\,
xy+2x^3
+
(\tfrac{1}{3}F_{2,2}-F_{4,0})y^4 
+
(36+F_{3,1}+G_{4,0}-G_{0,4})xy^3
+
F_{2,2}x^2y^2
+
\\
&
\ \ \ \ \
+
F_{3,1}x^3y
+
F_{4,0}x^4
+
\cdots,
\\
v
&
\,=\,
y^2+x^2-12x^2y
+
G_{0,4}y^4+G_{1,3}xy^3
+
(108+3G_{4,0}+3G_{0,4})x^2y^2
+
\\
&
\ \ \ \ \
+
(\tfrac{4}{3}F_{2,2}-8F_{4,0}+G_{1,3})x^3y
+
G_{4,0}x^4
+
\cdots,
\endaligned\right.
\]
%%%%%%%%%%%%%%%%%%%%%%%%%%%%%%%%%%%%%%%%%%%%%%%%%%%%%%%%%%%%%%%%%%%%%%
\[
\def\arraystretch{1.25}
\begin{array}{ll}
e_1
&
\,:=\,
-(2xF_{4,0}-1-\tfrac{1}{4}G_{1,3}x+\tfrac{1}{10}yF_{3,1}-\tfrac{1}{10}yG_{0,4}-\tfrac{1}{10}yG_{4,0}-\tfrac{18}{5}y+2uF_{2,2}-\tfrac{1}{2}uG_{1,3}
+
\\
&
\ \ \ \ \
+36v+vF_{3,1}-vG_{0,4}+vG_{4,0})\partial_x
-(\tfrac{1}{10}xF_{3,1}-\tfrac{1}{10}xG_{0,4}-\tfrac{1}{10}xG_{4,0}+\tfrac{12}{5}x+2yF_{4,0}
+
\\
&
\ \ \ \ \
-\tfrac{1}{4}yG_{1,3}+\tfrac{1}{5}uF_{3,1}+\tfrac{14}{5}uG_{0,4}+\tfrac{4}{5}uG_{4,0}+\tfrac{144}{5}u+\tfrac{1}{2}vG_{1,3})\partial_y
-(4uF_{4,0}-y-\tfrac{1}{2}uG_{1,3}
+
\\
&
\ \ \ \ \
+\tfrac{1}{10}vF_{3,1}-\tfrac{1}{10}vG_{0,4}-\tfrac{1}{10}vG_{4,0}-\tfrac{18}{5}v)\partial_u
-(-2x+\tfrac{108}{5}u+\tfrac{2}{5}uF_{3,1}-\tfrac{2}{5}uG_{0,4}-\tfrac{2}{5}uG_{4,0}
+
\\
&
\ \ \ \ \
+4vF_{4,0}-\tfrac{1}{2}vG_{1,3})\partial_v
, 
\\
e_2
&
\,:=\,
-(\tfrac{1}{2}xF_{3,1}-xG_{0,4}+6x+\tfrac{1}{15}yF_{2,2}-\tfrac{1}{20}yG_{1,3}-5uG_{0,4}+3uF_{3,1}+3uG_{4,0}+108u
+
\\
&
\ \ \ \ \
+\tfrac{4}{3}vF_{2,2}-4vF_{4,0})\partial_x
+(1-\tfrac{1}{15}xF_{2,2}+\tfrac{1}{20}xG_{1,3}-12y-\tfrac{1}{2}yF_{3,1}+yG_{0,4}+\tfrac{8}{15}uF_{2,2}
+
\\
&
\ \ \ \ \
-\tfrac{9}{10}uG_{1,3}-4uF_{4,0}-2vG_{0,4})\partial_y
-(uF_{3,1}-2uG_{0,4}+18u-x+\tfrac{1}{15}vF_{2,2}
+
\\
&
\ \ \ \ \
-\tfrac{1}{20}vG_{1,3})\partial_u
-(-2y+\tfrac{4}{15}uF_{2,2}-\tfrac{1}{5}uG_{1,3}+24v+vF_{3,1}-2vG_{0,4})\partial_v
.
\end{array}
\]

%%%%%%%%%%%%%%%%%%%%%%%%%%%%%%%%%%%%%%%%%%%%%%%%%%%%%%%%%%%%%%%%%%%%%%
\noindent
Gr\"obner basis generators of 
moduli space core algebraic variety in 
$\R^6 \ni F_{3,1},F_{2,2},F_{4,0},G_{0,4},G_{1,3},G_{4,0}$:
\[
\aligned
\B_1 & := 96F_{2,2}G_{4,0}+216F_{4,0}G_{0,4}-792F_{4,0}G_{4,0}-63G_{0,4}G_{1,3}+63G_{1,3}G_{4,0}+1736F_{2,2}
+
\\
&
\ \ \ \ \
-14352F_{4,0}+708G_{1,3},
\\
\B_2 & := 36F_{2,2}G_{0,4}+204F_{2,2}G_{4,0}-1440F_{4,0}G_{4,0}-117G_{0,4}G_{1,3}+117G_{1,3}G_{4,0}+4304F_{2,2}
+
\\
&
\ \ \ \ \
-28320F_{4,0}+1392G_{1,3},
\\
\B_3 & := 56F_{2,2}G_{4,0}+24F_{3,1}F_{4,0}-432F_{4,0}G_{4,0}-27G_{0,4}G_{1,3}+33G_{1,3}G_{4,0}+1176F_{2,2}
+
\\
&
\ \ \ \ \
-7632F_{4,0}+468G_{1,3},
\\
\B_4 & := 5832F_{3,1}G_{4,0}-5832F_{4,0}G_{1,3}+1536G_{0,4}^2-10224G_{0,4}G_{4,0}+165G_{1,3}^2+6048G_{4,0}^2
+
\\
&
\ \ \ \ \
+147872F_{3,1}-227904G_{0,4}+394496G_{4,0}+6363648,
\\
\B_5 & := 10344F_{3,1}G_{4,0}+24576F_{4,0}^2-16488F_{4,0}G_{1,3}-11952G_{0,4}G_{4,0}+585G_{1,3}^2+8736G_{4,0}^2
+
\\
&
\ \ \ \ \
+287264F_{3,1}-282432G_{0,4}+680192G_{4,0}+13165056,
\\
\B_6 & := 1536F_{3,1}G_{0,4}+14616F_{3,1}G_{4,0}-16152F_{4,0}G_{1,3}-24144G_{0,4}G_{4,0}+615G_{1,3}^2
+
\\
&
\ \ \ \ \
+14304G_{4,0}^2+375776F_{3,1}-556224G_{0,4}+958208G_{4,0}+15957504,
\\
\B_7 & := 512F_{3,1}^2+14776F_{3,1}G_{4,0}-16824F_{4,0}G_{1,3}-21264G_{0,4}G_{4,0}+555G_{1,3}^2+14432G_{4,0}^2
+
\\
&
\ \ \ \ \
+409696F_{3,1}-537024G_{0,4}+987904G_{4,0}+16860672.
\endaligned
\]
%%%%%%%%%%%%%%%%%%%%%%%%%%%%%%%%%%%%%%%%%%%%%%%%%%%%%%%%%%%%%%%%%%%%%%
\[
\footnotesize
\def\arraystretch{1.25}
\begin{array}{c|cc}
{} & e_1 & e_2 
\\
\hline
e_1 & 0 &
\rotatebox[origin=c]{0}{
\begin{tabular}{p{7cm}}
$ (\tfrac{-48}{5}-\tfrac{2}{5}F_{3,1}+\tfrac{9}{10}G_{0,4}-\tfrac{1}{10}G_{4,0})e_1+(2F_{4,0}-\tfrac{1}{5}G_{1,3}-\tfrac{1}{15}F_{2,2})e_2$
\end{tabular}}
\\
e_2 &
\rotatebox[origin=c]{0}{
\begin{tabular}{p{7cm}} $
-(\tfrac{-48}{5}-\tfrac{2}{5}F_{3,1}+\tfrac{9}{10}G_{0,4}-\tfrac{1}{10}G_{4,0})e_1-(2F_{4,0}-\tfrac{1}{5}G_{1,3}-\tfrac{1}{15}F_{2,2})e_2$
\end{tabular}}
& 0
\end{array}
\]

%%%%%%%%%%%%%%%%%%%%%%%%%%%%%%%%%%%%%%%%%%%%%%%%%%%%%%%%%%%%%%%%%%%%%%
%%%%%%%%%%%%%%%%%%%%%%%%%%%%%%%%%%%%%%%%%%%%%%%%%%%%%%%%%%%%%%%%%%%%%%
%%%%%%%%%%%%%%%%%%%%%%%%%%%%%%%%%%%%%%%%%%%%%%%%%%%%%%%%%%%%%%%%%%%%%%

\[
\text{\bf Model 2f3i}
\ \ \ \ \
\left\{
\aligned
u
&
\,=\,
xy+2y^3
+
(\tfrac{1}{3}F_{2,2}-F_{4,0})y^4+(G_{4,0}+F_{3,1}-G_{0,4}-36)xy^3
+
F_{2,2}x^2y^2
+
\\
&
\ \ \ \ \
+F_{3,1}x^3y+F_{4,0}x^4
+
\cdots
,
\\
v
&
\,=\,
y^2+x^2-12xy^2+G_{0,4}y^4+G_{1,3}xy^3
+
(3G_{4,0}+3G_{0,4}+108)x^2y^2
+
\\
&
\ \ \ \ \
+
(\tfrac{4}{3}F_{2,2}-8F_{4,0}+G_{1,3})x^3y
+G_{4,0}x^4
+
\cdots
.
\endaligned\right.
\]
%%%%%%%%%%%%%%%%%%%%%%%%%%%%%%%%%%%%%%%%%%%%%%%%%%%%%%%%%%%%%%%%%%%%%%
\[
\def\arraystretch{1.25}
\begin{array}{ll}
e_1
&
\,:=\,
-(-1+\tfrac{1}{2}xF_{3,1}-\tfrac{1}{2}xG_{0,4}-\tfrac{1}{2}xG_{4,0}-6x+\tfrac{2}{5}yF_{4,0}-\tfrac{1}{20}yG_{1,3}-\tfrac{56}{5}uF_{4,0}+\tfrac{9}{10}uG_{1,3}
+
\\
&
\ \ \ \ \
+2uF_{2,2}+2vG_{4,0})\partial_x
-(\tfrac{2}{5}xF_{4,0}-\tfrac{1}{20}G_{1,3}x-12y+\tfrac{1}{2}yF_{3,1}-\tfrac{1}{2}yG_{0,4}-\tfrac{1}{2}yG_{4,0}
+
\\
&
\ \ \ \ \
+3uF_{3,1}-2uG_{4,0}+4vF_{4,0})\partial_y
-(uF_{3,1}-uG_{0,4}-uG_{4,0}-18u-y+\tfrac{2}{5}vF_{4,0}
+
\\
&
\ \ \ \ \
-\tfrac{1}{20}vG_{1,3})\partial_u
-(-2x+\tfrac{8}{5}uF_{4,0}-\tfrac{1}{5}uG_{1,3}-12v+vF_{3,1}-vG_{0,4}-vG_{4,0})\partial_v 
,
\\
e_2
&
\,:=\,
-(\tfrac{1}{3}xF_{2,2}-\tfrac{1}{4}xG_{1,3}-\tfrac{6}{5}y+\tfrac{1}{10}yF_{3,1}-\tfrac{1}{5}yG_{0,4}+\tfrac{1}{5}uF_{3,1}+\tfrac{3}{5}uG_{0,4}+\tfrac{108}{5}u+3uG_{4,0}
+
\\
&
\ \ \ \ \
+\tfrac{2}{3}vF_{2,2}-4vF_{4,0}+\tfrac{1}{2}vG_{1,3})\partial_x
-(-1+\tfrac{1}{10}xF_{3,1}-\tfrac{1}{5}xG_{0,4}-\tfrac{36}{5}x+\tfrac{1}{3}yF_{2,2}-\tfrac{1}{4}yG_{1,3}
+
\\
&
\ \ \ \ \
+\tfrac{4}{3}uF_{2,2}+4uF_{4,0}-\tfrac{1}{2}uG_{1,3}+vF_{3,1})\partial_y
-(-x+\tfrac{2}{3}uF_{2,2}-\tfrac{1}{2}uG_{1,3}+\tfrac{1}{10}vF_{3,1}-\tfrac{1}{5}vG_{0,4}
+
\\
&
\ \ \ \ \
-\tfrac{36}{5}v)\partial_u
-(-2y+\tfrac{36}{5}u+\tfrac{2}{5}uF_{3,1}-\tfrac{4}{5}uG_{0,4}+\tfrac{2}{3}vF_{2,2}-\tfrac{1}{2}vG_{1,3})\partial_v.
\end{array}
\]

%%%%%%%%%%%%%%%%%%%%%%%%%%%%%%%%%%%%%%%%%%%%%%%%%%%%%%%%%%%%%%%%%%%%%%
\noindent
Gr\"obner basis generators of 
moduli space core algebraic variety in 
$\R^6 \ni F_{2,2},F_{3,1},F_{4,0},G_{0,4},G_{1,3},G_{4,0}$:
\[
\aligned
\B_1 & := 36F_{3,1}F_{4,0}-144F_{4,0}G_{4,0}-9G_{0,4}G_{1,3}+18G_{1,3}G_{4,0}+280F_{2,2}-1200F_{4,0}-24G_{1,3},
\\
\B_2 & := -240F_{2,2}G_{4,0}+216F_{4,0}G_{0,4}+1224F_{4,0}G_{4,0}+63G_{0,4}G_{1,3}-63G_{1,3}G_{4,0}-640F_{2,2}
+
\\
&
\ \ \ \ \
+11136F_{4,0}+888G_{1,3},
\\
\B_3 & := -120F_{2,2}G_{4,0}+27F_{3,1}G_{1,3}+288F_{4,0}G_{4,0}-9G_{0,4}G_{1,3}+9G_{1,3}G_{4,0}+280F_{2,2}
+
\\
&
\ \ \ \ \
-1200F_{4,0}-24G_{1,3},
\\
\B_4 & := -3672F_{3,1}G_{4,0}+3672F_{4,0}G_{1,3}+1200G_{0,4}^2+2304G_{0,4}G_{4,0}-129G_{1,3}^2-1440G_{4,0}^2
+
\\
&
\ \ \ \ \
-7360F_{3,1}+49920G_{0,4}+11744G_{4,0}+241920,
\\
\B_5 & := -192F_{3,1}G_{4,0}+192F_{4,0}^2+144F_{4,0}G_{1,3}+144G_{0,4}G_{4,0}+9G_{1,3}^2-240G_{4,0}^2+80F_{3,1}
+
\\
&
\ \ \ \ \
+240G_{0,4}-1456G_{4,0}+8640,
\\
\B_6 & := 300F_{3,1}G_{0,4}-756F_{3,1}G_{4,0}+456F_{4,0}G_{1,3}-408G_{0,4}G_{4,0}+33G_{1,3}^2-120G_{4,0}^2+7520F_{3,1}
+
\\
&
\ \ \ \ \
+960G_{0,4}-16288G_{4,0}+34560,
\\
\B_7 & := 200F_{3,1}^2-2168F_{3,1}G_{4,0}+1368F_{4,0}G_{1,3}+576G_{0,4}G_{4,0}+99G_{1,3}^2-1360G_{4,0}^2+8160F_{3,1}
+
\\
&
\ \ \ \ \
+2880G_{0,4}-20064G_{4,0}+103680.
\endaligned
\]
%%%%%%%%%%%%%%%%%%%%%%%%%%%%%%%%%%%%%%%%%%%%%%%%%%%%%%%%%%%%%%%%%%%%%%
\[
\footnotesize
\def\arraystretch{1.25}
\begin{array}{c|cc}
{} & e_1 & e_2 
\\
\hline
e_1 & 
0
&
\rotatebox[origin=c]{0}{
\begin{tabular}{p{7cm}}
$
(-\tfrac{1}{3}F_{2,2}+\tfrac{1}{5}G_{1,3}+\tfrac{2}{5}F_{4,0})e_1+(\tfrac{-24}{5}+\tfrac{2}{5}F_{3,1}-\tfrac{3}{10}G_{0,4}-\tfrac{1}{2}G_{4,0})e_2
$
\end{tabular}}
\\
e_2 & \rotatebox[origin=c]{0}{
\begin{tabular}{p{7cm}}
$
(-\tfrac{1}{3}F_{2,2}+\tfrac{1}{5}G_{1,3}+\tfrac{2}{5}F_{4,0})e_1+(\tfrac{-24}{5}+\tfrac{2}{5}F_{3,1}-\tfrac{3}{10}G_{0,4}-\tfrac{1}{2}G_{4,0})e_2
$
\end{tabular}}
& 0
\end{array}
\]

%%%%%%%%%%%%%%%%%%%%%%%%%%%%%%%%%%%%%%%%%%%%%%%%%%%%%%%%%%%%%%%%%%%%%%
%%%%%%%%%%%%%%%%%%%%%%%%%%%%%%%%%%%%%%%%%%%%%%%%%%%%%%%%%%%%%%%%%%%%%%
%%%%%%%%%%%%%%%%%%%%%%%%%%%%%%%%%%%%%%%%%%%%%%%%%%%%%%%%%%%%%%%%%%%%%%

\[
\text{\bf Model 2f3j4a}
\ \ \ \ \
\left\{
\aligned
u
&
\,=\,
xy+x^3+y^3+F_{4,0}x^4+F_{3,1}x^3y
+
(6F_{4,0}-\tfrac{21}{2})x^2y^2
+
(-1+F_{3,1})xy^3
+
\\
&
\,\,\,\,\,\,\,\,
+
(-\tfrac{1}{2}+F_{4,0})y^4
+
\cdots
,
\\
v
&
\,=\,
x^2+y^2-6x^2y-6xy^2+G_{4,0}x^4
+
(8F_{4,0}-2F_{3,1}+4G_{4,0})x^3y
+
\\
&
\,\,\,\,\,\,\,\,
+
(57+6G_{4,0})x^2y^2
+
(8F_{4,0}+2-2F_{3,1}+4G_{4,0})xy^3
+
(G_{4,0}+1)y^4
+
\cdots
.
\endaligned\right.
\]
%%%%%%%%%%%%%%%%%%%%%%%%%%%%%%%%%%%%%%%%%%%%%%%%%%%%%%%%%%%%%%%%%%%%%%
\[
\aligned
e_1 & := 
(1-\tfrac{40}{3}xF_{4,0}+\tfrac{19}{6}xF_{3,1}+\tfrac{1}{3}xG_{4,0}+\tfrac{8}{3}yF_{4,0}-\tfrac{5}{6}yF_{3,1}+\tfrac{1}{3}yG_{4,0}+\tfrac{185}{2}u-\tfrac{147}{4}v+\tfrac{127}{4}x
+
\\
&
\,\,\,\,\,\,\,\,
-\tfrac{136}{3}uF_{4,0}+\tfrac{32}{3}uF_{3,1}-\tfrac{8}{3}uG_{4,0}+\tfrac{40}{3}vF_{4,0}-\tfrac{25}{6}vF_{3,1}-\tfrac{1}{3}vG_{4,0}-\tfrac{15}{4}y)\partial_x
+
(\tfrac{8}{3}xF_{4,0}
+
\\
&
\,\,\,\,\,\,\,\,
-\tfrac{5}{6}xF_{3,1}+\tfrac{1}{3}xG_{4,0}-\tfrac{40}{3}yF_{4,0}+\tfrac{19}{6}yF_{3,1}+\tfrac{1}{3}yG_{4,0}+\tfrac{141}{2}u-\tfrac{115}{4}v-\tfrac{27}{4}x-\tfrac{112}{3}uF_{4,0}
+
\\
&
\,\,\,\,\,\,\,\,
+\tfrac{26}{3}uF_{3,1}-\tfrac{8}{3}uG_{4,0}+\tfrac{28}{3}vF_{4,0}-\tfrac{19}{6}vF_{3,1}-\tfrac{1}{3}vG_{4,0}+\tfrac{139}{4}y)\partial_y
+
(y-\tfrac{80}{3}uF_{4,0}+\tfrac{19}{3}uF_{3,1}
+
\\
&
\,\,\,\,\,\,\,\,
+\tfrac{2}{3}uG_{4,0}+\tfrac{133}{2}u+\tfrac{8}{3}vF_{4,0}-\tfrac{5}{6}vF_{3,1}+\tfrac{1}{3}vG_{4,0}-\tfrac{15}{4}v)\partial_u
-
(-2x-\tfrac{32}{3}uF_{4,0}+\tfrac{10}{3}uF_{3,1}
+
\\
&
\,\,\,\,\,\,\,\,
-\tfrac{4}{3}uG_{4,0}+33u+\tfrac{80}{3}vF_{4,0}-\tfrac{19}{3}vF_{3,1}-\tfrac{2}{3}vG_{4,0}-\tfrac{127}{2}v)\partial_v,\\
e_2 & := 
-
(\tfrac{7}{2}xF_{3,1}-\tfrac{1}{3}xG_{4,0}-\tfrac{40}{3}xF_{4,0}-\tfrac{1}{2}yF_{3,1}-\tfrac{1}{3}yG_{4,0}+\tfrac{8}{3}yF_{4,0}+\tfrac{213}{2}u-\tfrac{149}{4}v-\tfrac{25}{4}y
+
\\
&
\,\,\,\,\,\,\,\,
+10uF_{3,1}+\tfrac{8}{3}uG_{4,0}-\tfrac{112}{3}uF_{4,0}-\tfrac{7}{2}vF_{3,1}+\tfrac{1}{3}vG_{4,0}+\tfrac{52}{3}vF_{4,0}+\tfrac{117}{4}x)\partial_x
-
(-\tfrac{1}{2}xF_{3,1}
+
\\
&
\,\,\,\,\,\,\,\,
-1-\tfrac{1}{3}xG_{4,0}+\tfrac{8}{3}xF_{4,0}+\tfrac{7}{2}yF_{3,1}-\tfrac{1}{3}yG_{4,0}-\tfrac{40}{3}yF_{4,0}+\tfrac{181}{2}u-\tfrac{105}{4}v+\tfrac{129}{4}y+8uF_{3,1}
+
\\
&
\,\,\,\,\,\,\,\,
+\tfrac{8}{3}uG_{4,0}-\tfrac{88}{3}uF_{4,0}-\tfrac{5}{2}vF_{3,1}+\tfrac{1}{3}vG_{4,0}+\tfrac{40}{3}vF_{4,0}-\tfrac{37}{4}x)\partial_y
-
(-x+\tfrac{123}{2}u+7uF_{3,1}
+
\\
&
\,\,\,\,\,\,\,\,
-\tfrac{2}{3}uG_{4,0}-\tfrac{80}{3}uF_{4,0}-\tfrac{37}{4}v-\tfrac{1}{2}vF_{3,1}-\tfrac{1}{3}vG_{4,0}+\tfrac{8}{3}vF_{4,0})\partial_u
+
(2y+19u+2uF_{3,1}
+
\\
&
\,\,\,\,\,\,\,\,
+\tfrac{4}{3}uG_{4,0}-\tfrac{32}{3}uF_{4,0}-\tfrac{129}{2}v-7vF_{3,1}+\tfrac{2}{3}vG_{4,0}+\tfrac{80}{3}vF_{4,0})\partial_v,
\endaligned
\]
with invariants $F_{3,1},F_{4,0},G_{4,0}$ satisfying:
\[
\aligned
0
\,=\,
\B_1 & := (2F_{3,1}+19+8F_{4,0})(153+16F_{3,1}-64F_{4,0}),
\\
0
\,=\,
\B_2 & := (2F_{4,0}+G_{4,0}+7)(153+16F_{3,1}-64F_{4,0}).
\endaligned
\]
%%%%%%%%%%%%%%%%%%%%%%%%%%%%%%%%%%%%%%%%%%%%%%%%%%%%%%%%%%%%%%%%%%%%%%
\[
\footnotesize
\def\arraystretch{1.25}
\begin{array}{c|cc}
{} & e_1 & e_2 
\\
\hline
e_1 & 
0
&
\rotatebox[origin=c]{0}{
\begin{tabular}{p{7cm}}
$(-\tfrac{51}{2}-\tfrac{8}{3}F_{3,1}+\tfrac{32}{3}F_{4,0})e_1+(-\tfrac{51}{2}-\tfrac{8}{3}F_{3,1}+\tfrac{32}{3}F_{4,0})e_2$
\end{tabular}}
\\
e_2 &
\rotatebox[origin=c]{0}{
\begin{tabular}{p{7cm}}
$-(-\tfrac{51}{2}-\tfrac{8}{3}F_{3,1}+\tfrac{32}{3}F_{4,0})e_1-(-\tfrac{51}{2}-\tfrac{8}{3}F_{3,1}+\tfrac{32}{3}F_{4,0})e_2$
\end{tabular}}
& 
0
\end{array}
\]

%%%%%%%%%%%%%%%%%%%%%%%%%%%%%%%%%%%%%%%%%%%%%%%%%%%%%%%%%%%%%%%%%%%%%%
%%%%%%%%%%%%%%%%%%%%%%%%%%%%%%%%%%%%%%%%%%%%%%%%%%%%%%%%%%%%%%%%%%%%%%
%%%%%%%%%%%%%%%%%%%%%%%%%%%%%%%%%%%%%%%%%%%%%%%%%%%%%%%%%%%%%%%%%%%%%%

\[
\text{\bf Model 2f3j4b}
\ \ \ \ \
\left\{
\aligned
u
&
\,=\,
xy+x^3+y^3+F_{4,0}x^4+F_{3,1}x^3y
+
(6F_{4,0}-\tfrac{15}{2})x^2y^2
+
(1+F_{3,1})xy^3
+
\\
&
\,\,\,\,\,\,\,\,
+
(\tfrac{1}{2}+F_{4,0})y^4
+
\cdots
,
\\
v
&
\,=\,
x^2+y^2-6x^2y-6xy^2+G_{4,0}x^4
+
(8F_{4,0}-2F_{3,1}+4G_{4,0})x^3y
+
\\
&
\,\,\,\,\,\,\,\,
+
(51+6G_{4,0})x^2y^2
+
(8F_{4,0}-2-2F_{3,1}+4G_{4,0})xy^3
+
(G_{4,0}-1)y^4
+
\cdots,
\endaligned\right.
\]
%%%%%%%%%%%%%%%%%%%%%%%%%%%%%%%%%%%%%%%%%%%%%%%%%%%%%%%%%%%%%%%%%%%%%%
\[
\def\arraystretch{1.25}
\aligned
e_1 & := 
-
(-1+\tfrac{487}{6}u-\tfrac{269}{12}v+\tfrac{353}{12}x+\tfrac{8}{3}yF_{4,0}-\tfrac{1}{2}yF_{3,1}-\tfrac{1}{3}yG_{4,0}-\tfrac{40}{3}xF_{4,0}+\tfrac{7}{2}xF_{3,1}-\tfrac{1}{3}xG_{4,0}
+
\\
&
\,\,\,\,\,\,\,\,
+\tfrac{40}{3}vF_{4,0}-\tfrac{5}{2}vF_{3,1}+\tfrac{1}{3}vG_{4,0}-\tfrac{88}{3}uF_{4,0}+8uF_{3,1}+\tfrac{8}{3}uG_{4,0}-\tfrac{97}{12}y)\partial_x
-
(\tfrac{571}{6}u-\tfrac{389}{12}v-\tfrac{61}{12}x
+
\\
&
\,\,\,\,\,\,\,\,
-\tfrac{40}{3}yF_{4,0}+\tfrac{7}{2}yF_{3,1}-\tfrac{1}{3}yG_{4,0}+\tfrac{8}{3}xF_{4,0}-\tfrac{1}{2}xF_{3,1}-\tfrac{1}{3}xG_{4,0}+\tfrac{52}{3}vF_{4,0}-\tfrac{7}{2}vF_{3,1}+\tfrac{1}{3}vG_{4,0}
+
\\
&
\,\,\,\,\,\,\,\,
-\tfrac{112}{3}uF_{4,0}+10uF_{3,1}+\tfrac{8}{3}uG_{4,0}+\tfrac{317}{12}y)\partial_y
-
(-y-\tfrac{80}{3}uF_{4,0}+7uF_{3,1}-\tfrac{2}{3}uG_{4,0}+\tfrac{335}{6}u
+
\\
&
\,\,\,\,\,\,\,\,
+\tfrac{8}{3}vF_{4,0}-\tfrac{1}{2}vF_{3,1}-\tfrac{1}{3}vG_{4,0}-\tfrac{97}{12}v)\partial_u
+
(2x-\tfrac{32}{3}uF_{4,0}+2uF_{3,1}+\tfrac{4}{3}uG_{4,0}+\tfrac{43}{3}u+\tfrac{80}{3}vF_{4,0}
+
\\
&
\,\,\,\,\,\,\,\,
-7vF_{3,1}+\tfrac{2}{3}vG_{4,0}-\tfrac{353}{6}v)\partial_v,\\
e_2 & := 
(\tfrac{379}{6}u-\tfrac{323}{12}v-\tfrac{79}{12}y+\tfrac{1}{3}yG_{4,0}+\tfrac{8}{3}yF_{4,0}+\tfrac{19}{6}xF_{3,1}+\tfrac{1}{3}xG_{4,0}-\tfrac{40}{3}xF_{4,0}-\tfrac{5}{6}yF_{3,1}-\tfrac{19}{6}vF_{3,1}
+
\\
&
\,\,\,\,\,\,\,\,
-\tfrac{1}{3}vG_{4,0}+\tfrac{28}{3}vF_{4,0}+\tfrac{26}{3}uF_{3,1}-\tfrac{8}{3}uG_{4,0}-\tfrac{112}{3}uF_{4,0}+\tfrac{371}{12}x)\partial_x
+
(\tfrac{499}{6}u+1-\tfrac{407}{12}v+\tfrac{335}{12}y
+
\\
&
\,\,\,\,\,\,\,\,
+\tfrac{1}{3}yG_{4,0}-\tfrac{40}{3}yF_{4,0}-\tfrac{5}{6}xF_{3,1}+\tfrac{1}{3}xG_{4,0}+\tfrac{8}{3}xF_{4,0}+\tfrac{19}{6}yF_{3,1}-\tfrac{25}{6}vF_{3,1}-\tfrac{1}{3}vG_{4,0}+\tfrac{40}{3}vF_{4,0}
+
\\
&
\,\,\,\,\,\,\,\,
+\tfrac{32}{3}uF_{3,1}-\tfrac{8}{3}uG_{4,0}-\tfrac{136}{3}uF_{4,0}-\tfrac{43}{12}x)\partial_y
+
(x+\tfrac{353}{6}u+\tfrac{19}{3}uF_{3,1}+\tfrac{2}{3}uG_{4,0}-\tfrac{80}{3}uF_{4,0}-\tfrac{43}{12}v
+
\\
&
\,\,\,\,\,\,\,\,
-\tfrac{5}{6}vF_{3,1}+\tfrac{1}{3}vG_{4,0}+\tfrac{8}{3}vF_{4,0})\partial_u
-
(-2y+\tfrac{97}{3}u+\tfrac{10}{3}uF_{3,1}-\tfrac{4}{3}uG_{4,0}-\tfrac{32}{3}uF_{4,0}-\tfrac{335}{6}v-\tfrac{19}{3}vF_{3,1}
+
\\
&
\,\,\,\,\,\,\,\,
-\tfrac{2}{3}vG_{4,0}+\tfrac{80}{3}vF_{4,0})\partial_v,
\endaligned
\]
with invariants $F_{4,0},F_{3,1},G_{4,0}$ satisfying:
\[
\aligned
0
\,=\,
\B_1 & := (2F_{3,1}+25+8F_{4,0})(137-64F_{4,0}+16F_{3,1}),
\\
0
\,=\,
\B_2 & := (2F_{4,0}+G_{4,0}+7)(137-64F_{4,0}+16F_{3,1}),
\endaligned
\]
%%%%%%%%%%%%%%%%%%%%%%%%%%%%%%%%%%%%%%%%%%%%%%%%%%%%%%%%%%%%%%%%%%%%%%
\[
\footnotesize
\def\arraystretch{1.25}
\begin{array}{c|cc}
{} & e_1 & e_2 
\\
\hline
e_1 & 
0
&
\rotatebox[origin=c]{0}{
\begin{tabular}{p{7cm}}
$(\tfrac{137}{6}-\tfrac{32}{3}F_{4,0}+\tfrac{8}{3}F_{3,1})e_1+(\tfrac{137}{6}-\tfrac{32}{3}F_{4,0}+\tfrac{8}{3}F_{3,1})e_2$
\end{tabular}}
\\
e_2 &
\rotatebox[origin=c]{0}{
\begin{tabular}{p{7cm}}
$-(\tfrac{137}{6}-\tfrac{32}{3}F_{4,0}+\tfrac{8}{3}F_{3,1})e_1-^(\tfrac{137}{6}-\tfrac{32}{3}F_{4,0}+\tfrac{8}{3}F_{3,1})e_2$
\end{tabular}}
&
0
\end{array}
\]

\rotatebox[origin=c]{90}{
\begin{tabular}{p{3cm}}
\end{tabular}}

%%%%%%%%%%%%%%%%%%%%%%%%%%%%%%%%%%%%%%%%%%%%%%%%%%%%%%%%%%%%%%%%%%%%%%
%%%%%%%%%%%%%%%%%%%%%%%%%%%%%%%%%%%%%%%%%%%%%%%%%%%%%%%%%%%%%%%%%%%%%%
%%%%%%%%%%%%%%%%%%%%%%%%%%%%%%%%%%%%%%%%%%%%%%%%%%%%%%%%%%%%%%%%%%%%%%

\[
\text{\bf Model 2f3j4c}
\ \ \ \ \
\left\{
\aligned
u
&
\,=\,
xy+x^3+y^3+F_{4,0}x^4
+
(4F_{4,0}-9)x^3y
+
(6F_{4,0}-9)x^2y^2
+
(4F_{4,0}-9)xy^3
+
\\
&
\,\,\,\,\,\,\,\,
+
F_{4,0}y^4+F_{5,0}x^5
+
(-24F_{4,0}+5F_{5,0}+27)x^4y
+
(108-72F_{4,0}+10F_{5,0})x^3y^2
+
\\
&
\,\,\,\,\,\,\,\,
+
(108-72F_{4,0}+10F_{5,0})x^2y^3
+
(-24F_{4,0}+5F_{5,0}+27)xy^4+F_{5,0}y^5
+
\cdots
,
\\
v
&
\,=\,
x^2+y^2-6x^2y-6xy^2
+
(-2F_{4,0}-6)x^4
+
(-8F_{4,0}-6)x^3y
+
\\
&
\,\,\,\,\,\,\,\,
+
(18-12F_{4,0})x^2y^2
+
(-8F_{4,0}-6)xy^3
+
(-2F_{4,0}-6)y^4
+
G_{5,0}x^5
+
\\
&
\,\,\,\,\,\,\,\,
+
(90+48F_{4,0}+5G_{5,0})x^4y
+
(108+144F_{4,0}+10G_{5,0})x^3y^2
+
\\
&
\,\,\,\,\,\,\,\,
+
(108+144F_{4,0}+10G_{5,0})x^2y^3
+
(90+48F_{4,0}+5G_{5,0})xy^4
+
G_{5,0}y^5
+
\cdots
.
\endaligned\right.
\]
%%%%%%%%%%%%%%%%%%%%%%%%%%%%%%%%%%%%%%%%%%%%%%%%%%%%%%%%%%%%%%%%%%%%%%
\[
\def\arraystretch{1.25}
\aligned
e_1 & := 
-
(-1+\tfrac{71}{2}xF_{4,0}+\tfrac{25}{4}xF_{5,0}+\tfrac{25}{8}xG_{5,0}-\tfrac{11}{2}yF_{4,0}-\tfrac{5}{4}yF_{5,0}-\tfrac{5}{8}yG_{5,0}+93uF_{4,0}+\tfrac{35}{2}uF_{5,0}
+
\\
&
\,\,\,\,\,\,\,\,
+\tfrac{35}{4}uG_{5,0}-\tfrac{63}{2}vF_{4,0}-\tfrac{25}{4}vF_{5,0}-\tfrac{25}{8}vG_{5,0}-12y+51x+135u-54v)\partial_x
-
(-\tfrac{11}{2}xF_{4,0}
+
\\
&
\,\,\,\,\,\,\,\,
-\tfrac{5}{4}xF_{5,0}-\tfrac{5}{8}xG_{5,0}+\tfrac{71}{2}yF_{4,0}+\tfrac{25}{4}yF_{5,0}+\tfrac{25}{8}yG_{5,0}+93uF_{4,0}+\tfrac{35}{2}uF_{5,0}+\tfrac{35}{4}uG_{5,0}
+
\\
&
\,\,\,\,\,\,\,\,
-\tfrac{63}{2}vF_{4,0}-\tfrac{25}{4}vF_{5,0}-\tfrac{25}{8}vG_{5,0}+48y-9x+135u-54v)\partial_y
-
(71uF_{4,0}+99u-y+\tfrac{25}{2}uF_{5,0}
+
\\
&
\,\,\,\,\,\,\,\,
+\tfrac{25}{4}uG_{5,0}-12v-\tfrac{11}{2}vF_{4,0}-\tfrac{5}{4}vF_{5,0}-\tfrac{5}{8}vG_{5,0})\partial_u
+
(2x+30u+22uF_{4,0}+5uF_{5,0}+\tfrac{5}{2}uG_{5,0}
+
\\
&
\,\,\,\,\,\,\,\,
-102v-71vF_{4,0}-\tfrac{25}{2}vF_{5,0}-\tfrac{25}{4}vG_{5,0})\partial_v,\\
e_2 & := 
-
(\tfrac{71}{2}xF_{4,0}+\tfrac{25}{4}xF_{5,0}+\tfrac{25}{8}xG_{5,0}-\tfrac{11}{2}yF_{4,0}-\tfrac{5}{4}yF_{5,0}-\tfrac{5}{8}yG_{5,0}+93uF_{4,0}+\tfrac{35}{2}uF_{5,0}
+
\\
&
\,\,\,\,\,\,\,\,
+\tfrac{35}{4}uG_{5,0}-\tfrac{63}{2}vF_{4,0}-\tfrac{25}{4}vF_{5,0}-\tfrac{25}{8}vG_{5,0}+48x-9y+135u-54v)\partial_x
-
(-\tfrac{11}{2}xF_{4,0}
+
\\
&
\,\,\,\,\,\,\,\,
-\tfrac{5}{4}xF_{5,0}-1-\tfrac{5}{8}xG_{5,0}+\tfrac{71}{2}yF_{4,0}+\tfrac{25}{4}yF_{5,0}+\tfrac{25}{8}yG_{5,0}+93uF_{4,0}+\tfrac{35}{2}uF_{5,0}+\tfrac{35}{4}uG_{5,0}
+
\\
&
\,\,\,\,\,\,\,\,
-\tfrac{63}{2}vF_{4,0}-\tfrac{25}{4}vF_{5,0}-\tfrac{25}{8}vG_{5,0}-12x+51y+135u-54v)\partial_y
-
(71uF_{4,0}+99u-x
+
\\
&
\,\,\,\,\,\,\,\,
+\tfrac{25}{2}uF_{5,0}+\tfrac{25}{4}uG_{5,0}-12v-\tfrac{11}{2}vF_{4,0}-\tfrac{5}{4}vF_{5,0}-\tfrac{5}{8}vG_{5,0})\partial_u
+
(2y+30u+22uF_{4,0}
+
\\
&
\,\,\,\,\,\,\,\,
+5uF_{5,0}+\tfrac{5}{2}uG_{5,0}-102v-71vF_{4,0}-\tfrac{25}{2}vF_{5,0}-\tfrac{25}{4}vG_{5,0})\partial_v,
\endaligned
\]
with invariants $F_{4,0},F_{5,0},G_{5,0}$ satisfying:
\[
0
\,=\,
\B_1
\,:=\,
152F_{4,0}^2+20F_{4,0}F_{5,0}+10F_{4,0}G_{5,0}
+240F_{4,0}-10F_{5,0}+5G_{5,0}-72,
\]
%%%%%%%%%%%%%%%%%%%%%%%%%%%%%%%%%%%%%%%%%%%%%%%%%%%%%%%%%%%%%%%%%%%%%%
\[
\footnotesize
\def\arraystretch{1.25}
\begin{array}{c|cc}
{} & e_1 & e_2 
\\
\hline
e_1 & 
0
&
\rotatebox[origin=c]{0}{
\begin{tabular}{p{7cm}}
$(-60-41F_{4,0}-\tfrac{15}{2}F_{5,0}-\tfrac{15}{4}G_{5,0})e_1+(60+41F_{4,0}+\tfrac{15}{2}F_{5,0}+\tfrac{15}{4}G_{5,0})e_2$
\end{tabular}}
\\
e_2 &
\rotatebox[origin=c]{0}{
\begin{tabular}{p{7cm}}
$-(-60-41F_{4,0}-\tfrac{15}{2}F_{5,0}-\tfrac{15}{4}G_{5,0})e_1-(60+41F_{4,0}+\tfrac{15}{2}F_{5,0}+\tfrac{15}{4}G_{5,0})e_2$
\end{tabular}}
&
0
\end{array}
\]

\rotatebox[origin=c]{0}{
\begin{tabular}{p{7cm}}
\end{tabular}}

%%%%%%%%%%%%%%%%%%%%%%%%%%%%%%%%%%%%%%%%%%%%%%%%%%%%%%%%%%%%%%%%%%%%%%
%%%%%%%%%%%%%%%%%%%%%%%%%%%%%%%%%%%%%%%%%%%%%%%%%%%%%%%%%%%%%%%%%%%%%%
%%%%%%%%%%%%%%%%%%%%%%%%%%%%%%%%%%%%%%%%%%%%%%%%%%%%%%%%%%%%%%%%%%%%%%

\[
\text{\bf Model 2f3j4d}
\ \ \ \ \
\left\{
\aligned
u
&
\,=\,
xy+x^3+y^3+F_{4,0}x^4
+
(4F_{4,0}-9)x^3y
+
(6F_{4,0}-9)x^2y^2
+
(4F_{4,0}-9)xy^3
+
\\
&
\,\,\,\,\,\,\,\,
+
F_{4,0}y^4+F_{5,0}x^5
+
(-24F_{4,0}+5F_{5,0}+27)x^4y
+
(108-72F_{4,0}
+
\\
&
\,\,\,\,\,\,\,\,
+10F_{5,0})x^3y^2
+
(108-72F_{4,0}+10F_{5,0})x^2y^3
+
(-24F_{4,0}+5F_{5,0}
+
\\
&
\,\,\,\,\,\,\,\,
+27)xy^4
+
F_{5,0}y^5
+
\cdots
,
\\
v
&
\,=\,
x^2+y^2-6x^2y-6xy^2
+
(-2F_{4,0}-8)x^4
+
(-8F_{4,0}-14)x^3y
+
\\
&
\,\,\,\,\,\,\,\,
+
(6-12F_{4,0})x^2y^2
+
(-8F_{4,0}-14)xy^3
+
(-2F_{4,0}-8)y^4
+
G_{5,0}x^5
+
\\
&
\,\,\,\,\,\,\,\,
+
(138+48F_{4,0}+5G_{5,0})x^4y
+
(252+144F_{4,0}+10G_{5,0})x^3y^2
+
\\
&
\,\,\,\,\,\,\,\,
+
(252+144F_{4,0}+10G_{5,0})x^2y^3
+
(138+48F_{4,0}+5G_{5,0})xy^4
+
G_{5,0}y^5
+
\cdots,
\endaligned\right.
\]
%%%%%%%%%%%%%%%%%%%%%%%%%%%%%%%%%%%%%%%%%%%%%%%%%%%%%%%%%%%%%%%%%%%%%%
\[
\def\arraystretch{1.25}
\aligned
e_1 & := 
(-\tfrac{19}{2}yF_{4,0}-\tfrac{5}{4}yF_{5,0}+1-\tfrac{5}{8}yG_{5,0}+\tfrac{79}{2}xF_{4,0}+\tfrac{25}{4}xF_{5,0}+\tfrac{25}{8}xG_{5,0}+117uF_{4,0}+\tfrac{35}{2}uF_{5,0}
+
\\
&
\,\,\,\,\,\,\,\,
+\tfrac{35}{4}uG_{5,0}-\tfrac{87}{2}vF_{4,0}-\tfrac{25}{4}vF_{5,0}-\tfrac{25}{8}vG_{5,0}-\tfrac{33}{2}y+\tfrac{175}{2}x+260u-\tfrac{169}{2}v)\partial_x
+
(\tfrac{79}{2}yF_{4,0}
+
\\
&
\,\,\,\,\,\,\,\,
+\tfrac{25}{4}yF_{5,0}+\tfrac{25}{8}yG_{5,0}-\tfrac{19}{2}xF_{4,0}-\tfrac{5}{4}xF_{5,0}-\tfrac{5}{8}xG_{5,0}+117uF_{4,0}+\tfrac{35}{2}uF_{5,0}+\tfrac{35}{4}uG_{5,0}
+
\\
&
\,\,\,\,\,\,\,\,
-\tfrac{87}{2}vF_{4,0}-\tfrac{25}{4}vF_{5,0}-\tfrac{25}{8}vG_{5,0}+\tfrac{181}{2}y-\tfrac{39}{2}x+260u-\tfrac{169}{2}v)\partial_y
+
(y+178u+79uF_{4,0}
+
\\
&
\,\,\,\,\,\,\,\,
+\tfrac{25}{2}uF_{5,0}+\tfrac{25}{4}uG_{5,0}-\tfrac{33}{2}v-\tfrac{19}{2}vF_{4,0}-\tfrac{5}{4}vF_{5,0}-\tfrac{5}{8}vG_{5,0})\partial_u
-
(38uF_{4,0}+5uF_{5,0}+84u
+
\\
&
\,\,\,\,\,\,\,\,
-2x+\tfrac{5}{2}uG_{5,0}-175v-79vF_{4,0}-\tfrac{25}{2}vF_{5,0}-\tfrac{25}{4}vG_{5,0})\partial_v,\\
e_2 & := 
(-\tfrac{5}{8}yG_{5,0}+\tfrac{79}{2}xF_{4,0}+\tfrac{25}{4}xF_{5,0}+\tfrac{25}{8}xG_{5,0}-\tfrac{19}{2}yF_{4,0}-\tfrac{5}{4}yF_{5,0}+117uF_{4,0}+\tfrac{35}{2}uF_{5,0}
+
\\
&
\,\,\,\,\,\,\,\,
+\tfrac{35}{4}uG_{5,0}-\tfrac{87}{2}vF_{4,0}-\tfrac{25}{4}vF_{5,0}-\tfrac{25}{8}vG_{5,0}+\tfrac{181}{2}x-\tfrac{39}{2}y+260u-\tfrac{169}{2}v)\partial_x
+
(1+\tfrac{25}{8}yG_{5,0}
+
\\
&
\,\,\,\,\,\,\,\,
-\tfrac{19}{2}xF_{4,0}-\tfrac{5}{4}xF_{5,0}-\tfrac{5}{8}xG_{5,0}+\tfrac{79}{2}yF_{4,0}+\tfrac{25}{4}yF_{5,0}+117uF_{4,0}+\tfrac{35}{2}uF_{5,0}+\tfrac{35}{4}uG_{5,0}
+
\\
&
\,\,\,\,\,\,\,\,
-\tfrac{87}{2}vF_{4,0}-\tfrac{25}{4}vF_{5,0}-\tfrac{25}{8}vG_{5,0}-\tfrac{33}{2}x+\tfrac{175}{2}y+260u-\tfrac{169}{2}v)\partial_y
+
(x+178u+79uF_{4,0}
+
\\
&
\,\,\,\,\,\,\,\,
+\tfrac{25}{2}uF_{5,0}
+\tfrac{25}{4}uG_{5,0}-\tfrac{33}{2}v-\tfrac{19}{2}vF_{4,0}-\tfrac{5}{4}vF_{5,0}-\tfrac{5}{8}vG_{5,0})\partial_u
-
(38uF_{4,0}+5uF_{5,0}+84u
+
\\
&
\,\,\,\,\,\,\,\,
-2y+\tfrac{5}{2}uG_{5,0}-175v-79vF_{4,0}-\tfrac{25}{2}vF_{5,0}-\tfrac{25}{4}vG_{5,0})\partial_v,
\endaligned
\]
with invariants $F_{4,0},F_{5,0},G_{5,0}$ satisfying:
\[
0
\,=\,
\B_1
\,:=\,
44F_{4,0}^2+10F_{4,0}F_{5,0}+5F_{4,0}G_{5,0}+114F_{4,0}+10F_{5,0}
+100,
\]
%%%%%%%%%%%%%%%%%%%%%%%%%%%%%%%%%%%%%%%%%%%%%%%%%%%%%%%%%%%%%%%%%%%%%%
\[
\footnotesize
\def\arraystretch{1.25}
\begin{array}{c|cc}
{} & e_1 & e_2 
\\
\hline
e_1 & 
0
&
\rotatebox[origin=c]{0}{
\begin{tabular}{p{7cm}}
$(107+49F_{4,0}+\tfrac{15}{2}F_{5,0}+\tfrac{15}{4}G_{5,0})e_1+(-107-49F_{4,0}-\tfrac{15}{2}F_{5,0}-\tfrac{15}{4}G_{5,0})e_2$
\end{tabular}}
\\
e_2 &
\rotatebox[origin=c]{0}{
\begin{tabular}{p{7cm}}
$-(107+49F_{4,0}+\tfrac{15}{2}F_{5,0}+\tfrac{15}{4}G_{5,0})e_1-(-107-49F_{4,0}-\tfrac{15}{2}F_{5,0}-\tfrac{15}{4}G_{5,0})e_2$
\end{tabular}} & 0
\end{array}
\]

%%%%%%%%%%%%%%%%%%%%%%%%%%%%%%%%%%%%%%%%%%%%%%%%%%%%%%%%%%%%%%%%%%%%%%
%%%%%%%%%%%%%%%%%%%%%%%%%%%%%%%%%%%%%%%%%%%%%%%%%%%%%%%%%%%%%%%%%%%%%%
%%%%%%%%%%%%%%%%%%%%%%%%%%%%%%%%%%%%%%%%%%%%%%%%%%%%%%%%%%%%%%%%%%%%%%

\[
\text{\bf Model 2f3j4e}
\ \ \ \ \
\left\{
\aligned
u
&
\,=\,
xy+x^3+y^3-\tfrac{1}{2}x^4-11x^3y-12x^2y^2-11xy^3-\tfrac{1}{2}y^4+F_{5,0}x^5
+
\\
&
\,\,\,\,\,\,\,\,
+
(39+5F_{5,0})x^4y
+
(144+10F_{5,0})x^3y^2
+
(144+10F_{5,0})x^2y^3
+
\\
&
\,\,\,\,\,\,\,\,
+
(39+5F_{5,0})xy^4
+
F_{5,0}y^5
+
(-\tfrac{20}{3}F_{5,0}^2-\tfrac{373}{3}F_{5,0}-572)x^6
+
\\
&
\,\,\,\,\,\,\,\,
+
(-40F_{5,0}^2-776F_{5,0}-3549)x^5y
+
(-100F_{5,0}^2-1985F_{5,0}-9480)x^4y^2
+
\\
&
\,\,\,\,\,\,\,\,
+
(-\tfrac{400}{3}F_{5,0}^2-\tfrac{8000}{3}F_{5,0}-13060)x^3y^3
+
(-100F_{5,0}^2-1985F_{5,0}-9480)x^2y^4
+
\\
&
\,\,\,\,\,\,\,\,
+
(-40F_{5,0}^2-776F_{5,0}-3549)xy^5
+
(-\tfrac{20}{3}F_{5,0}^2-\tfrac{373}{3}F_{5,0}-572)y^6
+
\cdots
,
\\
v
&
\,=\,
x^2+y^2-6x^2y-6xy^2-6x^4-6x^3y+18x^2y^2-6xy^3-6y^4
+
\\
&
\,\,\,\,\,\,\,\,
+
(-2F_{5,0}-16)x^5
+
(10-10F_{5,0})x^4y
+
(-52-20F_{5,0})x^3y^2
+
\\
&
\,\,\,\,\,\,\,\,
+
(-52-20F_{5,0})x^2y^3
+
(10-10F_{5,0})xy^4
+
(-2F_{5,0}-16)y^5
+
\\
&
\,\,\,\,\,\,\,\,
+
(\tfrac{40}{3}F_{5,0}^2+\tfrac{710}{3}F_{5,0}+1120)x^6+(80F_{5,0}^2+1480F_{5,0}+6930)x^5y
+
\\
&
\,\,\,\,\,\,\,\,
+
(200F_{5,0}^2+3790F_{5,0}+17316)x^4y^2
+
(\tfrac{800}{3}F_{5,0}^2+\tfrac{15280}{3}F_{5,0}+23336)x^3y^3
+
\\
&
\,\,\,\,\,\,\,\,
+
(200F_{5,0}^2+3790F_{5,0}+17316)x^2y^4
+
(80F_{5,0}^2+1480F_{5,0}+6930)xy^5
+
\\
&
\,\,\,\,\,\,\,\,
+
(\tfrac{40}{3}F_{5,0}^2+\tfrac{710}{3}F_{5,0}+1120)y^6
+
\cdots
.
\endaligned\right.
\]
%%%%%%%%%%%%%%%%%%%%%%%%%%%%%%%%%%%%%%%%%%%%%%%%%%%%%%%%%%%%%%%%%%%%%%
\[
\def\arraystretch{1.25}
\aligned
e_1 & := 
(1+\tfrac{467}{2}x+25xF_{5,0}-\tfrac{89}{2}y-5yF_{5,0}+70F_{5,0}u+662u-\tfrac{457}{2}v-25vF_{5,0})\partial_x
+
(-\tfrac{95}{2}x
+
\\
&
\,\,\,\,\,\,\,\,
-5xF_{5,0}+\tfrac{473}{2}y+25yF_{5,0}+70F_{5,0}u+662u-\tfrac{457}{2}v-25vF_{5,0})\partial_y
+
(y+50F_{5,0}u+470u
+
\\
&
\,\,\,\,\,\,\,\,
-\tfrac{89}{2}v-5vF_{5,0})\partial_u
-
(20F_{5,0}u-50F_{5,0}v+196u-467v-2x)\partial_v,\\
e_2 & := 
(\tfrac{473}{2}x+25xF_{5,0}-\tfrac{95}{2}y-5yF_{5,0}+70F_{5,0}u+662u-\tfrac{457}{2}v-25vF_{5,0})\partial_x
+
(-\tfrac{89}{2}x
+
\\
&
\,\,\,\,\,\,\,\,
-5xF_{5,0}+1+\tfrac{467}{2}y+25yF_{5,0}+70F_{5,0}u+662u-\tfrac{457}{2}v-25vF_{5,0})\partial_y
+
(x+50F_{5,0}u
+
\\
&
\,\,\,\,\,\,\,\,
+470u-\tfrac{89}{2}v-5vF_{5,0})\partial_u
-
(20F_{5,0}u-50F_{5,0}v+196u-467v-2y)\partial_v.
\endaligned
\]
%%%%%%%%%%%%%%%%%%%%%%%%%%%%%%%%%%%%%%%%%%%%%%%%%%%%%%%%%%%%%%%%%%%%%%
\[
\footnotesize
\def\arraystretch{1.25}
\begin{array}{c|cc}
{} & e_1 & e_2 
\\
\hline
e_1 & 
0 & (30F_{5,0}+281)e_1-(30F_{5,0}+281)e_2
\\
e_2 &
-(30F_{5,0}+281)e_1+(30F_{5,0}+281)e_2 & 0
\end{array}
\]

\rotatebox[origin=c]{90}{
\begin{tabular}{p{3cm}}
\end{tabular}}

%%%%%%%%%%%%%%%%%%%%%%%%%%%%%%%%%%%%%%%%%%%%%%%%%%%%%%%%%%%%%%%%%%%%%%
%%%%%%%%%%%%%%%%%%%%%%%%%%%%%%%%%%%%%%%%%%%%%%%%%%%%%%%%%%%%%%%%%%%%%%
%%%%%%%%%%%%%%%%%%%%%%%%%%%%%%%%%%%%%%%%%%%%%%%%%%%%%%%%%%%%%%%%%%%%%%

\[
\text{\bf Model 2f3j4f}
\ \ \ \ \
\left\{
\aligned
u
&
\,=\,
xy+x^3+y^3-9x^3y-9x^2y^2-9xy^3+F_{5,0}x^5
+
(27+5F_{5,0})x^4y
+
\\
&
\,\,\,\,\,\,\,\,
+
(108+10F_{5,0})x^3y^2
+
(108+10F_{5,0})x^2y^3
+
(27+5F_{5,0})xy^4
+
\\
&
\,\,\,\,\,\,\,\,
+
F_{5,0}y^5
+
(\tfrac{20}{3}F_{5,0}^2+\tfrac{365}{3}F_{5,0}+547)x^6
+
(40F_{5,0}^2+700F_{5,0}+3201)x^5y
+
\\
&
\,\,\,\,\,\,\,\,
+
(100F_{5,0}^2+1705F_{5,0}+7557)x^4y^2
+
(\tfrac{400}{3}F_{5,0}^2+\tfrac{6760}{3}F_{5,0}+9752)x^3y^3
+
\\
&
\,\,\,\,\,\,\,\,
+
(100F_{5,0}^2+1705F_{5,0}+7557)x^2y^4
+
(40F_{5,0}^2+700F_{5,0}+3201)xy^5
+
\\
&
\,\,\,\,\,\,\,\,
+
(\tfrac{20}{3}F_{5,0}^2+\tfrac{365}{3}F_{5,0}+547)y^6
+
\cdots
,
\\
v
&
\,=\,
x^2+y^2-6x^2y-6xy^2-7x^4-10x^3y+12x^2y^2-10xy^3-7y^4
+
\\
&
\,\,\,\,\,\,\,\,
+
(-2F_{5,0}-22)x^5
+
(4-10F_{5,0})x^4y
+
(-40-20F_{5,0})x^3y^2
+
\\
&
\,\,\,\,\,\,\,\,
+
(-40-20F_{5,0})x^2y^3
+
(4-10F_{5,0})xy^4
+
(-2F_{5,0}-22)y^5
+
\\
&
\,\,\,\,\,\,\,\,
+
(-\tfrac{40}{3}F_{5,0}^2-\tfrac{766}{3}F_{5,0}-1136)x^6
+
(-80F_{5,0}^2-1472F_{5,0}-6498)x^5y
+
\\
&
\,\,\,\,\,\,\,\,
+
(-200F_{5,0}^2-3590F_{5,0}-16308)x^4y^2
+
(-\tfrac{800}{3}F_{5,0}^2-\tfrac{14240}{3}F_{5,0}-21568)x^3y^3
+
\\
&
\,\,\,\,\,\,\,\,
+
(-200F_{5,0}^2-3590F_{5,0}-16308)x^2y^4
+
(-80F_{5,0}^2-1472F_{5,0}-6498)xy^5
+
\\
&
\,\,\,\,\,\,\,\,
+
(-\tfrac{40}{3}F_{5,0}^2-\tfrac{766}{3}F_{5,0}-1136)y^6
+
\cdots,
\endaligned\right.
\]
for any value for $F_{5,0}$,
%%%%%%%%%%%%%%%%%%%%%%%%%%%%%%%%%%%%%%%%%%%%%%%%%%%%%%%%%%%%%%%%%%%%%%
\[
\def\arraystretch{1.25}
\aligned
e_1 & := 
-
(70F_{5,0}u-25F_{5,0}v+25F_{5,0}x-5F_{5,0}y+613u-226v+223x-46y-1)\partial_x
+
\\
&
\,\,\,\,\,\,\,\,
-
(70F_{5,0}u-25F_{5,0}v-5F_{5,0}x+25F_{5,0}y+613u-226v-43x+220y)\partial_y
+
\\
&
\,\,\,\,\,\,\,\,
-
(50F_{5,0}u-5F_{5,0}v+443u-46v-y)\partial_u
+
(20F_{5,0}u-50F_{5,0}v+166u-446v+2x)\partial_v,\\
e_2 & := 
-
(70F_{5,0}u-25F_{5,0}v+25F_{5,0}x-5F_{5,0}y+613u-226v+220x-43y)\partial_x
+
\\
&
\,\,\,\,\,\,\,\,
-
(70F_{5,0}u-25F_{5,0}v-5F_{5,0}x+25F_{5,0}y+613u-226v-46x+223y-1)\partial_y
+
\\
&
\,\,\,\,\,\,\,\,
-
(50F_{5,0}u-5F_{5,0}v+443u-46v-x)\partial_u
+
(20F_{5,0}u-50F_{5,0}v+166u-446v+2y)\partial_v.
\endaligned
\]
%%%%%%%%%%%%%%%%%%%%%%%%%%%%%%%%%%%%%%%%%%%%%%%%%%%%%%%%%%%%%%%%%%%%%%
\[
\footnotesize
\def\arraystretch{1.25}
\begin{array}{c|cc}
{} & e_1 & e_2 
\\
\hline
e_1 & 
0 & -(30F_{5,0}+266)e_1+(30F_{5,0}+266)e_2
\\
e_2 &
(30F_{5,0}+266)e_1-(30F_{5,0}+266)e_2 & 0
\end{array}
\]

\rotatebox[origin=c]{90}{
\begin{tabular}{p{3cm}}
\end{tabular}}

%%%%%%%%%%%%%%%%%%%%%%%%%%%%%%%%%%%%%%%%%%%%%%%%%%%%%%%%%%%%%%%%%%%%%%
%%%%%%%%%%%%%%%%%%%%%%%%%%%%%%%%%%%%%%%%%%%%%%%%%%%%%%%%%%%%%%%%%%%%%%
%%%%%%%%%%%%%%%%%%%%%%%%%%%%%%%%%%%%%%%%%%%%%%%%%%%%%%%%%%%%%%%%%%%%%%

\[
\text{\bf Model 2f3j4g}
\ \ \ \ \
\left\{
\aligned
u
&
\,=\,
xy+x^3+y^3-\tfrac{1}{4}x^4-10x^3y-\tfrac{21}{2}x^2y^2-10xy^3-\tfrac{1}{4}y^4-\tfrac{37}{4}x^5-\tfrac{53}{4}x^4y
+
\\
&
\,\,\,\,\,\,\,\,
+\tfrac{67}{2}x^3y^2+\tfrac{67}{2}x^2y^3-\tfrac{53}{4}xy^4-\tfrac{37}{4}y^5
+
\cdots
,
\\
v
&
\,=\,
x^2+y^2-6x^2y-6xy^2-\tfrac{13}{2}x^4-8x^3y+15x^2y^2-8xy^3-\tfrac{13}{2}y^4-\tfrac{1}{2}x^5
+
\\
&
\,\,\,\,\,\,\,\,
+\tfrac{199}{2}x^4y+139x^3y^2+139x^2y^3+\tfrac{199}{2}xy^4-\tfrac{1}{2}y^5
+
\cdots,
\endaligned\right.
\]
%%%%%%%%%%%%%%%%%%%%%%%%%%%%%%%%%%%%%%%%%%%%%%%%%%%%%%%%%%%%%%%%%%%%%%
\[
\def\arraystretch{1.25}
\aligned
e_1 & := 
(1+\tfrac{9}{5}y+\tfrac{49}{5}u+4v)\partial_x
+
(-\tfrac{6}{5}x+3y+\tfrac{49}{5}u+4v)\partial_y
+
(y+3u+\tfrac{9}{5}v)\partial_u
-
(-2x+\tfrac{54}{5}u)\partial_v,\\
e_2 & := 
(-\tfrac{3}{5}y+\tfrac{7}{5}u+7v)\partial_x
+
(1+\tfrac{12}{5}x-3y+\tfrac{7}{5}u+7v)\partial_y
-
(-x+3u-\tfrac{12}{5}v)\partial_u
-
(-2y+\tfrac{42}{5}u+6v)\partial_v,\\
e_3 & := 
(x-\tfrac{1}{5}y+\tfrac{14}{5}u-v)\partial_x
+
(-\tfrac{1}{5}x+y+\tfrac{14}{5}u-v)\partial_y
+
(2u-\tfrac{1}{5}v)\partial_u
-
(\tfrac{4}{5}u-2v)\partial_v,
\endaligned
\]

%%%%%%%%%%%%%%%%%%%%%%%%%%%%%%%%%%%%%%%%%%%%%%%%%%%%%%%%%%%%%%%%%%%%%%
\[
\footnotesize
\def\arraystretch{1.25}
\begin{array}{c|ccc}
{} & e_1 & e_2 & e_3
\\
\hline
e_1 & 
0 
&
-\tfrac{9}{5}e_1-\tfrac{3}{5}e_2+\tfrac{3}{5}e_3
&
e_1-\tfrac{1}{5}e_2-\tfrac{7}{5}e_3
\\
e_2 &
\tfrac{9}{5}e_1+\tfrac{3}{5}e_2-\tfrac{3}{5}e_3 & 0
& 
-\tfrac{1}{5}e_1+e_2+\tfrac{11}{5}e_3
\\
e_3 &
-e_1+\tfrac{1}{5}e_2+\tfrac{7}{5}e_3
& 
\tfrac{1}{5}e_1-e_2-\tfrac{11}{5}e_3
& 
0
\end{array}
\]

%%%%%%%%%%%%%%%%%%%%%%%%%%%%%%%%%%%%%%%%%%%%%%%%%%%%%%%%%%%%%%%%%%%%%%
%%%%%%%%%%%%%%%%%%%%%%%%%%%%%%%%%%%%%%%%%%%%%%%%%%%%%%%%%%%%%%%%%%%%%%
%%%%%%%%%%%%%%%%%%%%%%%%%%%%%%%%%%%%%%%%%%%%%%%%%%%%%%%%%%%%%%%%%%%%%%

\[
\text{\bf Model 2f3k}
\ \ \ \ \
\left\{
\aligned
u
&
\,=\,
xy-2y^3+2x^3
+
(-\tfrac{1}{4}F_{3,1}+\tfrac{1}{2}G_{0,4}+\tfrac{1}{8}G_{3,1})y^4
+
\\
&
\ \ \ \ \
+
(\tfrac{1}{2}G_{1,3}+F_{3,1}-\tfrac{1}{2}G_{3,1})xy^3
+
F_{2,2}x^2y^2+F_{3,1}x^3y+F_{4,0}x^4+\cdots
,
\\
v
&
\,=\,
y^2+x^2+4xy^2-4x^2y+G_{0,4}y^4+G_{1,3}xy^3
+
\\
&
\ \ \ \ \
+
(2F_{2,2}+3F_{3,1}-\tfrac{3}{2}G_{3,1}+64)x^2y^2
+
G_{3,1}x^3y
+
\\
&
\ \ \ \ \
+
(\tfrac{1}{2}F_{3,1}+2F_{4,0}-\tfrac{1}{4}G_{3,1})x^4
+
\cdots,
\endaligned\right.
\]
%%%%%%%%%%%%%%%%%%%%%%%%%%%%%%%%%%%%%%%%%%%%%%%%%%%%%%%%%%%%%%%%%%%%%%
\[
\def\arraystretch{1.25}
\begin{array}{ll}
e_1
&
\,:=\,
(1+8u-8v+\tfrac{3}{4}xF_{2,2}-\tfrac{7}{2}xF_{4,0}+\tfrac{7}{16}G_{1,3}x-\tfrac{9}{16}xG_{3,1}+\tfrac{1}{4}yF_{2,2}-\tfrac{1}{2}yF_{4,0}+\tfrac{1}{16}yG_{1,3}
+
\\
&
\ \ \ \ \
-\tfrac{3}{16}yG_{3,1}-vF_{3,1}+2uF_{2,2}-8uF_{4,0}+\tfrac{3}{2}uG_{1,3}-3uG_{3,1}+\tfrac{3}{2}vF_{2,2}-7vF_{4,0}+\tfrac{3}{8}vG_{1,3}
+
\\
&
\ \ \ \ \
-\tfrac{5}{8}vG_{3,1}+4y)\partial_x
-
(2uF_{3,1}+32u+4v-\tfrac{1}{4}xF_{2,2}+\tfrac{1}{2}xF_{4,0}-\tfrac{1}{16}G_{1,3}x+\tfrac{3}{16}xG_{3,1}
+
\\
&
\ \ \ \ \
-\tfrac{3}{4}yF_{2,2}+\tfrac{7}{2}yF_{4,0}-\tfrac{7}{16}yG_{1,3}+\tfrac{9}{16}yG_{3,1}+4uF_{2,2}-12uF_{4,0}+uG_{1,3}-\tfrac{5}{2}uG_{3,1}
+
\\
&
\ \ \ \ \
+\tfrac{3}{2}vF_{2,2}-3vF_{4,0}+\tfrac{7}{8}vG_{1,3}-\tfrac{9}{8}vG_{3,1}+2y+2x)\partial_y
+
(y+\tfrac{3}{2}uF_{2,2}-7uF_{4,0}+\tfrac{7}{8}uG_{1,3}
+
\\
&
\ \ \ \ \
-\tfrac{9}{8}uG_{3,1}-2u+4v+\tfrac{1}{4}vF_{2,2}-\tfrac{1}{2}vF_{4,0}+\tfrac{1}{16}vG_{1,3}-\tfrac{3}{16}vG_{3,1})\partial_u
+
\\
&
\ \ \ \ \
+
(uF_{2,2}-2uF_{4,0}-4u+2x+\tfrac{1}{4}uG_{1,3}-\tfrac{3}{4}uG_{3,1}+\tfrac{3}{2}vF_{2,2}-7vF_{4,0}+\tfrac{7}{8}vG_{1,3}-\tfrac{9}{8}vG_{3,1})\partial_v
, 
\\
e_2
&
\,:=\,
-
(32u+4v-\tfrac{9}{16}xG_{1,3}+\tfrac{7}{8}xF_{3,1}-\tfrac{7}{4}xG_{0,4}+\tfrac{3}{4}xF_{2,2}-\tfrac{3}{16}yG_{1,3}+\tfrac{1}{8}yF_{3,1}-\tfrac{1}{4}yG_{0,4}
+
\\
&
\ \ \ \ \
+\tfrac{1}{4}yF_{2,2}-\tfrac{3}{2}uG_{3,1}+\tfrac{1}{2}vG_{3,1}-\tfrac{3}{2}uG_{1,3}+5uF_{3,1}-6uG_{0,4}+4uF_{2,2}-\tfrac{9}{8}vG_{1,3}
+
\\
&
\ \ \ \ \
+\tfrac{3}{4}vF_{3,1}-\tfrac{3}{2}vG_{0,4}+\tfrac{3}{2}vF_{2,2}-2x-2y)\partial_x
+
(8u+1-8v+\tfrac{3}{16}xG_{1,3}-\tfrac{1}{8}xF_{3,1}
+
\\
&
\ \ \ \ \
+\tfrac{1}{4}xG_{0,4}-\tfrac{1}{4}xF_{2,2}+\tfrac{9}{16}yG_{1,3}-\tfrac{7}{8}yF_{3,1}+\tfrac{7}{4}yG_{0,4}-\tfrac{3}{4}yF_{2,2}+\tfrac{1}{2}uG_{3,1}-3uG_{1,3}
+
\\
&
\ \ \ \ \
+2uF_{3,1}-4uG_{0,4}+2uF_{2,2}-\tfrac{9}{8}vG_{1,3}+\tfrac{3}{4}vF_{3,1}-\tfrac{7}{2}vG_{0,4}+\tfrac{3}{2}vF_{2,2}-4x)\partial_y
+
\\
&
\ \ \ \ \
-
(-x-\tfrac{9}{8}uG_{1,3}+\tfrac{7}{4}uF_{3,1}-\tfrac{7}{2}uG_{0,4}+\tfrac{3}{2}uF_{2,2}-2u-\tfrac{3}{16}vG_{1,3}+\tfrac{1}{8}vF_{3,1}-\tfrac{1}{4}vG_{0,4}
+
\\
&
\ \ \ \ \
+\tfrac{1}{4}vF_{2,2}+4v)\partial_u
-
(-2y-4u-\tfrac{3}{4}uG_{1,3}+\tfrac{1}{2}uF_{3,1}-uG_{0,4}+uF_{2,2}-\tfrac{9}{8}vG_{1,3}
+
\\
&
\ \ \ \ \
+\tfrac{7}{4}vF_{3,1}-\tfrac{7}{2}vG_{0,4}+\tfrac{3}{2}vF_{2,2})\partial_v,
\end{array}
\]
%%%%%%%%%%%%%%%%%%%%%%%%%%%%%%%%%%%%%%%%%%%%%%%%%%%%%%%%%%%%%%%%%%%%%%
Gr\"obner basis generators of 
moduli space core algebraic variety in 
$\R^6 \ni F_{2,2},F_{3,1},F_{4,0},G_{0,4},G_{1,3},G_{3,1}$:
\[
\aligned
\B_1 & := -16F_{2,2}G_{1,3}+16F_{2,2}G_{3,1}+4F_{3,1}^2-16F_{3,1}G_{0,4}-10F_{3,1}G_{1,3}+10F_{3,1}G_{3,1}-64F_{4,0}^2
+
\\
&
\ \ \ \ \
+56F_{4,0}G_{1,3}-40F_{4,0}G_{3,1}+16G_{0,4}^2+20G_{0,4}G_{1,3}-20G_{0,4}G_{3,1}+5G_{1,3}G_{3,1}-6G_{3,1}^2
+
\\
&
\ \ \ \ \
-32F_{3,1}-128F_{4,0}+64G_{0,4}+16G_{3,1};
\\
\B_2 & := 16F_{2,2}G_{1,3}-16F_{2,2}G_{3,1}+16F_{3,1}F_{4,0}+8F_{3,1}G_{0,4}+6F_{3,1}G_{1,3}-10F_{3,1}G_{3,1}+64F_{4,0}^2
+
\\
&
\ \ \ \ \
-56F_{4,0}G_{1,3}+32F_{4,0}G_{3,1}-16G_{0,4}^2-20G_{0,4}G_{1,3}+24G_{0,4}G_{3,1}-3G_{1,3}G_{3,1}+5G_{3,1}^2
+
\\
&
\ \ \ \ \
+128F_{3,1}+512F_{4,0}-256G_{0,4}-96G_{1,3}+32G_{3,1},
\\
\B_3 & := 128F_{2,2}F_{3,1}-1280F_{2,2}G_{0,4}+560F_{2,2}G_{1,3}-816F_{2,2}G_{3,1}+272F_{3,1}G_{0,4}+86F_{3,1}G_{1,3}
+
\\
&
\ \ \ \ \
-506F_{3,1}G_{3,1}+3072F_{4,0}^2+5184F_{4,0}G_{0,4}-2984F_{4,0}G_{1,3}+1944F_{4,0}G_{3,1}-544G_{0,4}^2
+
\\
&
\ \ \ \ \
-828G_{0,4}G_{1,3}+1788G_{0,4}G_{3,1}+160G_{1,3}^2-139G_{1,3}G_{3,1}+253G_{3,1}^2-10752F_{2,2}
+
\\
&
\ \ \ \ \
+6880F_{3,1}+71552F_{4,0}-19904G_{0,4}-11968G_{1,3}+8528G_{3,1}-40960,
\\
\B_4 & := 64F_{2,2}^2+1536F_{2,2}G_{0,4}-288F_{2,2}G_{1,3}+576F_{2,2}G_{3,1}-144F_{3,1}G_{0,4}+90F_{3,1}G_{1,3}
+
\\
&
\ \ \ \ \
+330F_{3,1}G_{3,1}-2304F_{4,0}^2-6720F_{4,0}G_{0,4}+2184F_{4,0}G_{1,3}-1464F_{4,0}G_{3,1}+288G_{0,4}^2
+
\\
&
\ \ \ \ \
+540G_{0,4}G_{1,3}-1692G_{0,4}G_{3,1}-192G_{1,3}^2+51G_{1,3}G_{3,1}-165G_{3,1}^2+15616F_{2,2}
+
\\
&
\ \ \ \ \
-5600F_{3,1}-86400F_{4,0}+19392G_{0,4}+13824G_{1,3}-11408G_{3,1}+55296,
\\
\B_5 & := 28128F_{2,2}G_{1,3}^2-19328F_{2,2}G_{1,3}G_{3,1}-8800F_{2,2}G_{3,1}^2-5522F_{3,1}G_{1,3}^2
+
\\
&
\ \ \ \ \
+29508F_{3,1}G_{1,3}G_{3,1}-23986F_{3,1}G_{3,1}^2+77312F_{4,0}^2G_{1,3}-77312F_{4,0}^2G_{3,1}
+
\\
&
\ \ \ \ \
-163592F_{4,0}G_{1,3}^2+234000F_{4,0}G_{1,3}G_{3,1}-70408F_{4,0}G_{3,1}^2+11044G_{0,4}G_{1,3}^2
+
\\
&
\ \ \ \ \
-59016G_{0,4}G_{1,3}G_{3,1}+47972G_{0,4}G_{3,1}^2+15164G_{1,3}^3-50755G_{1,3}^2G_{3,1}+38762G_{1,3}G_{3,1}^2
+
\\
&
\ \ \ \ \
-3171G_{3,1}^3-887296F_{2,2}G_{1,3}+887296F_{2,2}G_{3,1}-399136F_{3,1}G_{1,3}+642784F_{3,1}G_{3,1}
+
\\
&
\ \ \ \ \
+2107264F_{4,0}G_{1,3}-1132672F_{4,0}G_{3,1}+489024G_{0,4}G_{1,3}-976320G_{0,4}G_{3,1}-60960G_{1,3}^2
+
\\
&
\ \ \ \ \
+199568G_{1,3}G_{3,1}-260432G_{3,1}^2-3537408F_{3,1}-14149632F_{4,0}+7074816G_{0,4}
+
\\
&
\ \ \ \ \
+293888G_{1,3}+1474816G_{3,1},
\\
\B_6 & := 19328F_{2,2}G_{0,4}G_{1,3}-19328F_{2,2}G_{0,4}G_{3,1}+5232F_{2,2}G_{1,3}^2-5232F_{2,2}G_{3,1}^2+1410F_{3,1}G_{1,3}^2
+
\\
&
\ \ \ \ \
+4828F_{3,1}G_{1,3}G_{3,1}-6238F_{3,1}G_{3,1}^2-77312F_{4,0}G_{0,4}G_{1,3}+77312F_{4,0}G_{0,4}G_{3,1}+
\\
&
\ \ \ \ \
-15288F_{4,0}G_{1,3}^2+19312F_{4,0}G_{1,3}G_{3,1}-4024F_{4,0}G_{3,1}^2+4428G_{0,4}G_{1,3}^2
+
\\
&
\ \ \ \ \
-28984G_{0,4}G_{1,3}G_{3,1}+24556G_{0,4}G_{3,1}^2+552G_{1,3}^3-7393G_{1,3}^2G_{3,1}+6690G_{1,3}G_{3,1}^2
+
\\
&
\ \ \ \ \
+151G_{3,1}^3+29952F_{2,2}G_{1,3}-29952F_{2,2}G_{3,1}-116320F_{3,1}G_{1,3}+162976F_{3,1}G_{3,1}
+
\\
&
\ \ \ \ \
-585088F_{4,0}G_{1,3}+771712F_{4,0}G_{3,1}+271296G_{0,4}G_{1,3}-364608G_{0,4}G_{3,1}+138016G_{1,3}^2
+
\\
&
\ \ \ \ \
-251088G_{1,3}G_{3,1}+89744G_{3,1}^2-677376F_{3,1}-2709504F_{4,0}+1354752G_{0,4}
+
\\
&
\ \ \ \ \
+737280G_{1,3}-398592G_{3,1},
\\
\B_7 & := -7392F_{2,2}G_{1,3}^2+7248F_{2,2}G_{1,3}G_{3,1}+144F_{2,2}G_{3,1}^2+4832F_{3,1}G_{0,4}^2-2416F_{3,1}G_{0,4}G_{3,1}
+
\\
&
\ \ \ \ \
-3588F_{3,1}G_{1,3}^2+3408F_{3,1}G_{1,3}G_{3,1}+482F_{3,1}G_{3,1}^2+19328F_{4,0}G_{0,4}^2-9664F_{4,0}G_{0,4}G_{3,1}
+
\\
&
\ \ \ \ \
+15216F_{4,0}G_{1,3}^2-15360F_{4,0}G_{1,3}G_{3,1}+1352F_{4,0}G_{3,1}^2-9664G_{0,4}^3-4832G_{0,4}^2G_{1,3}
+
\endaligned
\]

\[
\aligned
&
\ \ \ \ \
+7248G_{0,4}^2G_{3,1}+7176G_{0,4}G_{1,3}^2-4400G_{0,4}G_{1,3}G_{3,1}-2172G_{0,4}G_{3,1}^2+816G_{1,3}^3
+
\\
&
\ \ \ \ \
+2136G_{1,3}^2G_{3,1}-3254G_{1,3}G_{3,1}^2+151G_{3,1}^3+109824F_{2,2}G_{1,3}-109824F_{2,2}G_{3,1}
+
\\
&
\ \ \ \ \
+135296F_{3,1}G_{0,4}+69176F_{3,1}G_{1,3}-93800F_{3,1}G_{3,1}+541184F_{4,0}G_{0,4}-162592F_{4,0}G_{1,3}
+
\\
&
\ \ \ \ \
+64096F_{4,0}G_{3,1}-270592G_{0,4}^2-273648G_{0,4}G_{1,3}+255248G_{0,4}G_{3,1}-42776G_{1,3}^2
+
\\
&
\ \ \ \ \
+33060G_{1,3}G_{3,1}+22028G_{3,1}^2+802048F_{3,1}+3208192F_{4,0}-1604096G_{0,4}
+
\\
&
\ \ \ \ \
-582400G_{1,3}+181376G_{3,1},
\endaligned
\]
%%%%%%%%%%%%%%%%%%%%%%%%%%%%%%%%%%%%%%%%%%%%%%%%%%%%%%%%%%%%%%%%%%%%%%
\[
\footnotesize
\def\arraystretch{1.25}
\begin{array}{c|cc}
{} & e_1 & e_2 
\\
\hline
e_1 & 0 & \rotatebox[origin=c]{0}{
\begin{tabular}{p{7cm}}
$(-2+\tfrac{1}{2}G_{1,3}-\tfrac{7}{8}F_{3,1}-F_{2,2}+\tfrac{1}{2}F_{4,0}+\tfrac{7}{4}G_{0,4}+\tfrac{3}{16}G_{3,1})e_1+(-2-\tfrac{1}{4}G_{1,3}-\tfrac{1}{8}F_{3,1}-F_{2,2}+\tfrac{7}{2}F_{4,0}+\tfrac{1}{4}G_{0,4}+\tfrac{9}{16}G_{3,1})e_2$
\end{tabular}}
\\
e_2 &
\rotatebox[origin=c]{0}{
\begin{tabular}{p{7cm}}
$-(-2+\tfrac{1}{2}G_{1,3}-\tfrac{7}{8}F_{3,1}-F_{2,2}+\tfrac{1}{2}F_{4,0}+\tfrac{7}{4}G_{0,4}+\tfrac{3}{16}G_{3,1})e_1-(-2-\tfrac{1}{4}G_{1,3}-\tfrac{1}{8}F_{3,1}-F_{2,2}+\tfrac{7}{2}F_{4,0}+\tfrac{1}{4}G_{0,4}+\tfrac{9}{16}G_{3,1})e_2$
\end{tabular}}
& 0
\end{array}
\]

%%%%%%%%%%%%%%%%%%%%%%%%%%%%%%%%%%%%%%%%%%%%%%%%%%%%%%%%%%%%%%%%%%%%%%
%%%%%%%%%%%%%%%%%%%%%%%%%%%%%%%%%%%%%%%%%%%%%%%%%%%%%%%%%%%%%%%%%%%%%%
%%%%%%%%%%%%%%%%%%%%%%%%%%%%%%%%%%%%%%%%%%%%%%%%%%%%%%%%%%%%%%%%%%%%%%

\[
\text{\bf Model 2f3n4a}
\ \ \ \ \
\left\{
\aligned
u
&
\,=\,
xy+4xy^3+4x^3y
+
\tfrac{1}{2}G_{5,0}y^5+5F_{5,0}xy^4+5G_{5,0}x^2y^3+10F_{5,0}x^3y^2
+
\\
&
\ \ \ \ \
+\tfrac{5}{2}G_{5,0}x^4y+F_{5,0}x^5+\tfrac{5}{4}F_{5,0}G_{5,0}y^6+(24+\tfrac{15}{8}G_{5,0}^2+\tfrac{15}{2}F_{5,0}^2)xy^5
+
\\
&
\ \ \ \ \
+\tfrac{75}{4}F_{5,0}G_{5,0}x^2y^4+(80+25F_{5,0}^2+\tfrac{25}{4}G_{5,0}^2)x^3y^3+\tfrac{75}{4}F_{5,0}G_{5,0}x^4y^2
+
\\
&
\ \ \ \ \
+(24+\tfrac{15}{8}G_{5,0}^2+\tfrac{15}{2}F_{5,0}^2)x^5y+\tfrac{5}{4}F_{5,0}G_{5,0}x^6+(\tfrac{34}{7}G_{5,0}+\tfrac{25}{112}G_{5,0}^3
+
\\
&
\ \ \ \ \
+\tfrac{75}{28}F_{5,0}^2G_{5,0})y^7+(68F_{5,0}+\tfrac{75}{8}F_{5,0}G_{5,0}^2+\tfrac{25}{2}F_{5,0}^3)xy^6+(102G_{5,0}
+
\\
&
\ \ \ \ \
+\tfrac{75}{16}G_{5,0}^3+\tfrac{225}{4}F_{5,0}^2G_{5,0})x^2y^5+(340F_{5,0}+\tfrac{125}{2}F_{5,0}^3+\tfrac{375}{8}F_{5,0}G_{5,0}^2)x^3y^4
+
\\
&
\ \ \ \ \
+(170G_{5,0}+\tfrac{375}{4}F_{5,0}^2G_{5,0}+\tfrac{125}{16}G_{5,0}^3)x^4y^3+(204F_{5,0}+\tfrac{75}{2}F_{5,0}^3
+
\\
&
\ \ \ \ \
+\tfrac{225}{8}F_{5,0}G_{5,0}^2)x^5y^2+(34G_{5,0}+\tfrac{25}{16}G_{5,0}^3+\tfrac{75}{4}F_{5,0}^2G_{5,0})x^6y+(\tfrac{68}{7}F_{5,0}
+
\\
&
\ \ \ \ \
+\tfrac{75}{56}F_{5,0}G_{5,0}^2+\tfrac{25}{14}F_{5,0}^3)x^7
+
\cdots
,
\\
v
&
\,=\,
y^2+x^2+2y^4+12x^2y^2+2x^4+2F_{5,0}y^5+5G_{5,0}xy^4+20F_{5,0}x^2y^3
+
\\
&
\ \ \ \ \
+10G_{5,0}x^3y^2+10F_{5,0}x^4y+G_{5,0}x^5+(8+\tfrac{5}{8}G_{5,0}^2+\tfrac{5}{2}F_{5,0}^2)y^6
+
\\
&
\ \ \ \ \
+15F_{5,0}G_{5,0}xy^5+(120+\tfrac{75}{8}G_{5,0}^2+\tfrac{75}{2}F_{5,0}^2)x^2y^4+50F_{5,0}G_{5,0}x^3y^3
+
\\
&
\ \ \ \ \
+(120+\tfrac{75}{8}G_{5,0}^2+\tfrac{75}{2}F_{5,0}^2)x^4y^2+15F_{5,0}G_{5,0}x^5y+(8+\tfrac{5}{8}G_{5,0}^2+\tfrac{5}{2}F_{5,0}^2)x^6
+
\\
&
\ \ \ \ \
+(\tfrac{136}{7}F_{5,0}+\tfrac{25}{7}F_{5,0}^3+\tfrac{75}{28}F_{5,0}G_{5,0}^2)y^7+(68G_{5,0}+\tfrac{25}{8}G_{5,0}^3+\tfrac{75}{2}F_{5,0}^2G_{5,0})xy^6
+
\\
&
\ \ \ \ \
+(408F_{5,0}+75F_{5,0}^3+\tfrac{225}{4}F_{5,0}G_{5,0}^2)x^2y^5+(340G_{5,0}+\tfrac{125}{8}G_{5,0}^3
+
\\
&
\ \ \ \ \
+\tfrac{375}{2}F_{5,0}^2G_{5,0})x^3y^4+(125F_{5,0}^3+\tfrac{375}{4}F_{5,0}G_{5,0}^2+680F_{5,0})x^4y^3
+
\\
&
\ \ \ \ \
+(204G_{5,0}+\tfrac{75}{8}G_{5,0}^3+\tfrac{225}{2}F_{5,0}^2G_{5,0})x^5y^2+(136F_{5,0}+25F_{5,0}^3
+
\\
&
\ \ \ \ \
+\tfrac{75}{4}F_{5,0}G_{5,0}^2)x^6y+(\tfrac{68}{7}G_{5,0}+\tfrac{25}{56}G_{5,0}^3+\tfrac{75}{14}F_{5,0}^2G_{5,0})x^7
+
\cdots
.
\endaligned\right.
\]
%%%%%%%%%%%%%%%%%%%%%%%%%%%%%%%%%%%%%%%%%%%%%%%%%%%%%%%%%%%%%%%%%%%%%%
\[
\def\arraystretch{1.25}
\begin{array}{ll}
e_1
&
\,:=\,
-(-1+\tfrac{5}{4}xG_{5,0}+\tfrac{5}{2}yF_{5,0}+4v)\partial_x-(\tfrac{5}{2}xF_{5,0}+\tfrac{5}{4}yG_{5,0}+8u)\partial_y
+
\\
&
\ \ \ \ \
-(-y+\tfrac{5}{2}uG_{5,0}+\tfrac{5}{2}vF_{5,0})\partial_u-(-2x+10F_{5,0}u+\tfrac{5}{2}vG_{5,0})\partial_v
, 
\\
e_2
&
\,:=\,
-(\tfrac{5}{2}xF_{5,0}+\tfrac{5}{4}yG_{5,0}+8u)\partial_x-(-1+\tfrac{5}{4}xG_{5,0}+\tfrac{5}{2}yF_{5,0}+4v)\partial_y
+
\\
&
\ \ \ \ \
-(-x+5F_{5,0}u+\tfrac{5}{4}vG_{5,0})\partial_u-(5F_{5,0}v+5G_{5,0}u-2y)\partial_v
,
\end{array}
\]
for any value of $F_{5,0},G_{5,0}$,
%%%%%%%%%%%%%%%%%%%%%%%%%%%%%%%%%%%%%%%%%%%%%%%%%%%%%%%%%%%%%%%%%%%%%%
\[
\footnotesize
\def\arraystretch{1.25}
\begin{array}{c|cc}
{} & e_1 & e_2 
\\
\hline
e_1 & 0 & 0
\\
e_2 & 0 & 0
\end{array}
\]
%%%%%%%%%%%%%%%%%%%%%%%%%%%%%%%%%%%%%%%%%%%%%%%%%%%%%%%%%%%%%%%%%%%%%%
%%%%%%%%%%%%%%%%%%%%%%%%%%%%%%%%%%%%%%%%%%%%%%%%%%%%%%%%%%%%%%%%%%%%%%
%%%%%%%%%%%%%%%%%%%%%%%%%%%%%%%%%%%%%%%%%%%%%%%%%%%%%%%%%%%%%%%%%%%%%%

\[
\text{\bf Model 2f3n4b}
\ \ \ \ \
\left\{
\aligned
u
&
\,=\,
xy+y^4+6x^2y^2+x^4+\tfrac{1}{2}G_{5,0}y^5+5F_{5,0}xy^4+5G_{5,0}x^2y^3+10F_{5,0}x^3y^2
+
\\
&
\ \ \ \ \
+\tfrac{5}{2}G_{5,0}x^4y+F_{5,0}x^5+(\tfrac{5}{4}F_{5,0}^2+\tfrac{5}{16}G_{5,0}^2)y^6+(24+\tfrac{15}{2}F_{5,0}G_{5,0})xy^5
+
\\
&
\ \ \ \ \
+(\tfrac{75}{4}F_{5,0}^2+\tfrac{75}{16}G_{5,0}^2)x^2y^4+(25F_{5,0}G_{5,0}+80)x^3y^3+(\tfrac{75}{4}F_{5,0}^2
+
\\
&
\ \ \ \ \
+\tfrac{75}{16}G_{5,0}^2)x^4y^2+(24+\tfrac{15}{2}F_{5,0}G_{5,0})x^5y+(\tfrac{5}{4}F_{5,0}^2+\tfrac{5}{16}G_{5,0}^2)x^6
+
\\
&
\ \ \ \ \
+(\tfrac{68}{7}F_{5,0}+\tfrac{75}{28}F_{5,0}^2G_{5,0}+\tfrac{25}{112}G_{5,0}^3)y^7+(\tfrac{25}{2}F_{5,0}^3+\tfrac{75}{8}G_{5,0}^2F_{5,0}
+
\\
&
\ \ \ \ \
+34G_{5,0})xy^6+(204F_{5,0}+\tfrac{225}{4}F_{5,0}^2G_{5,0}+\tfrac{75}{16}G_{5,0}^3)x^2y^5+(170G_{5,0}
+
\\
&
\ \ \ \ \
+\tfrac{125}{2}F_{5,0}^3+\tfrac{375}{8}G_{5,0}^2F_{5,0})x^3y^4+(340F_{5,0}+\tfrac{375}{4}F_{5,0}^2G_{5,0}+\tfrac{125}{16}G_{5,0}^3)x^4y^3
+
\\
&
\ \ \ \ \
+(102G_{5,0}+\tfrac{75}{2}F_{5,0}^3+\tfrac{225}{8}G_{5,0}^2F_{5,0})x^5y^2+(68F_{5,0}+\tfrac{75}{4}F_{5,0}^2G_{5,0}
+
\\
&
\ \ \ \ \
+\tfrac{25}{16}G_{5,0}^3)x^6y+(\tfrac{25}{14}F_{5,0}^3+\tfrac{75}{56}G_{5,0}^2F_{5,0}+\tfrac{34}{7}G_{5,0})x^7
+
\cdots
,
\\
v
&
\,=\,
y^2+x^2+8xy^3+8x^3y+2F_{5,0}y^5+5G_{5,0}xy^4+20F_{5,0}x^2y^3
+
\\
&
\ \ \ \ \
+10G_{5,0}x^3y^2+10F_{5,0}x^4y+G_{5,0}x^5+(\tfrac{5}{2}F_{5,0}G_{5,0}+8)y^6+(15F_{5,0}^2
+
\\
&
\ \ \ \ \
+\tfrac{15}{4}G_{5,0}^2)xy^5+(120+\tfrac{75}{2}F_{5,0}G_{5,0})x^2y^4+(\tfrac{25}{2}G_{5,0}^2+50F_{5,0}^2)x^3y^3
+
\\
&
\ \ \ \ \
+(120+\tfrac{75}{2}F_{5,0}G_{5,0})x^4y^2+15F_{5,0}^2+\tfrac{15}{4}G_{5,0}^2)x^5y+(\tfrac{5}{2}F_{5,0}G_{5,0}
+
\\
&
\ \ \ \ \
+8)x^6+(\tfrac{68}{7}G_{5,0}+\tfrac{25}{7}F_{5,0}^3+\tfrac{75}{28}G_{5,0}^2F_{5,0})y^7+(136F_{5,0}+\tfrac{75}{2}F_{5,0}^2G_{5,0}
+
\\
&
\ \ \ \ \
+\tfrac{25}{8}G_{5,0}^3)xy^6+(204G_{5,0}+75F_{5,0}^3+\tfrac{225}{4}G_{5,0}^2F_{5,0})x^2y^5+(680F_{5,0}
+
\\
&
\ \ \ \ \
+\tfrac{375}{2}F_{5,0}^2G_{5,0}+\tfrac{125}{8}G_{5,0}^3)x^3y^4+(125F_{5,0}^3+\tfrac{375}{4}G_{5,0}^2F_{5,0}+340G_{5,0})x^4y^3
+
\\
&
\ \ \ \ \
+(408F_{5,0}+\tfrac{225}{2}F_{5,0}^2G_{5,0}+\tfrac{75}{8}G_{5,0}^3)x^5y^2+(68G_{5,0}+25F_{5,0}^3
+
\\
&
\ \ \ \ \
+\tfrac{75}{4}G_{5,0}^2F_{5,0})x^6y+(\tfrac{136}{7}F_{5,0}+\tfrac{75}{14}F_{5,0}^2G_{5,0}+\tfrac{25}{56}G_{5,0}^3)x^7
+
\cdots.
\endaligned\right.
\]
%%%%%%%%%%%%%%%%%%%%%%%%%%%%%%%%%%%%%%%%%%%%%%%%%%%%%%%%%%%%%%%%%%%%%%
\[
\def\arraystretch{1.25}
\begin{array}{ll}
e_1
&
\,:=\,
-(-1+\tfrac{5}{2}xF_{5,0}+\tfrac{5}{4}yG_{5,0}+8u)\partial_x-(\tfrac{5}{4}xG_{5,0}+\tfrac{5}{2}yF_{5,0}+4v)\partial_y
+
\\
&
\ \ \ \ \
-(-y+5uF_{5,0}+\tfrac{5}{4}vG_{5,0})\partial_u-(5F_{5,0}v+5G_{5,0}u-2x)\partial_v 
,
\\
e_2
&
\,:=\,
-(\tfrac{5}{4}xG_{5,0}+\tfrac{5}{2}yF_{5,0}+4v)\partial_x-(-1+\tfrac{5}{2}xF_{5,0}+\tfrac{5}{4}yG_{5,0}+8u)\partial_y
+
\\
&
\ \ \ \ \
-(-x+\tfrac{5}{2}G_{5,0}u+\tfrac{5}{2}F_{5,0}v)\partial_u-(-2y+10uF_{5,0}+\tfrac{5}{2}vG_{5,0})\partial_v,
\end{array}
\]
for any value of $F_{5,0},G_{5,0}$,
%%%%%%%%%%%%%%%%%%%%%%%%%%%%%%%%%%%%%%%%%%%%%%%%%%%%%%%%%%%%%%%%%%%%%%
\[
\footnotesize
\def\arraystretch{1.25}
\begin{array}{c|cc}
{} & e_1 & e_2 
\\
\hline
e_1 & 0 & 0
\\
e_2 & 0 & 0
\end{array}
\]

%%%%%%%%%%%%%%%%%%%%%%%%%%%%%%%%%%%%%%%%%%%%%%%%%%%%%%%%%%%%%%%%%%%%%%
%%%%%%%%%%%%%%%%%%%%%%%%%%%%%%%%%%%%%%%%%%%%%%%%%%%%%%%%%%%%%%%%%%%%%%
%%%%%%%%%%%%%%%%%%%%%%%%%%%%%%%%%%%%%%%%%%%%%%%%%%%%%%%%%%%%%%%%%%%%%%

\[
\!\!\!\!\!\!\!\!\!\!\!\!\!\!
\text{\bf Model 2f3n4c}
\ \ \ \ \
\left\{
\aligned
u
&
\,=\,
xy+\tfrac{1}{2}y^4+2xy^3+3x^2y^2+2x^3y+\tfrac{1}{2}x^4+F_{5,0}y^5+5F_{5,0}xy^4+10F_{5,0}x^2y^3
+
\\
&
\ \ \ \ \
+10F_{5,0}x^3y^2+5F_{5,0}x^4y+F_{5,0}x^5+(\tfrac{5}{2}F_{5,0}^2+2)y^6+(15F_{5,0}^2+12)xy^5
+
\\
&
\ \ \ \ \
+(\tfrac{75}{2}F_{5,0}^2+30)x^2y^4+(50F_{5,0}^2+40)x^3y^3+(\tfrac{75}{2}F_{5,0}^2+30)x^4y^2
+
\\
&
\ \ \ \ \
+(15F_{5,0}^2+12)x^5y+(\tfrac{5}{2}F_{5,0}^2+2)x^6+(\tfrac{50}{7}F_{5,0}^3+\tfrac{68}{7}F_{5,0})y^7+(50F_{5,0}^3
+
\\
&
\ \ \ \ \
+68F_{5,0})xy^6+(150F_{5,0}^3+204F_{5,0})x^2y^5+(250F_{5,0}^3+340F_{5,0})x^3y^4
+
\\
&
\ \ \ \ \
+(250F_{5,0}^3+340F_{5,0})x^4y^3+(150F_{5,0}^3+204F_{5,0})x^5y^2+(50F_{5,0}^3+68F_{5,0})x^6y
+
\\
&
\ \ \ \ \
+(\tfrac{50}{7}F_{5,0}^3+\tfrac{68}{7}F_{5,0})x^7
+
\cdots
,
\\
v
&
\,=\,
y^2+x^2+y^4+4xy^3+6x^2y^2+4x^3y+x^4+2F_{5,0}y^5+10F_{5,0}xy^4+20F_{5,0}x^2y^3
+
\\
&
\ \ \ \ \
+20F_{5,0}x^3y^2+10F_{5,0}x^4y+2F_{5,0}x^5+(5F_{5,0}^2+4)y^6+(30F_{5,0}^2+24)xy^5
+
\\
&
\ \ \ \ \
+(75F_{5,0}^2+60)x^2y^4+(100F_{5,0}^2+80)x^3y^3+(75F_{5,0}^2+60)x^4y^2+(30F_{5,0}^2
+
\\
&
\ \ \ \ \
+24)x^5y+(5F_{5,0}^2+4)x^6+(\tfrac{100}{7}F_{5,0}^3+\tfrac{136}{7}F_{5,0})y^7+(100F_{5,0}^3+136F_{5,0})xy^6
+
\\
&
\ \ \ \ \
+(300F_{5,0}^3+408F_{5,0})x^2y^5+(500F_{5,0}^3+680F_{5,0})x^3y^4+(500F_{5,0}^3+680F_{5,0})x^4y^3
+
\\
&
\ \ \ \ \
+(300F_{5,0}^3+408F_{5,0})x^5y^2+(100F_{5,0}^3+136F_{5,0})x^6y+(\tfrac{100}{7}F_{5,0}^3+\tfrac{136}{7}F_{5,0})x^7
+
\cdots
.
\endaligned\right.
\]
%%%%%%%%%%%%%%%%%%%%%%%%%%%%%%%%%%%%%%%%%%%%%%%%%%%%%%%%%%%%%%%%%%%%%%
\[
\def\arraystretch{1.25}
\begin{array}{ll}
e_1
&
\,:=\,
-(5F_{5,0}x+4u+2v-1)\partial_x-(5F_{5,0}y+4u+2v)\partial_y-(10F_{5,0}u-y)\partial_u-(10F_{5,0}v-2x)\partial_v
, 
\\
e_2
&
\,:=\,
-(5F_{5,0}x+4u+2v)\partial_x-(5F_{5,0}y+4u+2v-1)\partial_y-(10F_{5,0}u-x)\partial_u-(10F_{5,0}v-2y)\partial_v
, 
\\
e_3
&
\,:=\,
-(x-y)\partial_x+(x-y)\partial_y-(2u-v)\partial_u+(4u-2v)\partial_v
.
\end{array}
\]

%%%%%%%%%%%%%%%%%%%%%%%%%%%%%%%%%%%%%%%%%%%%%%%%%%%%%%%%%%%%%%%%%%%%%%
\[
\footnotesize
\def\arraystretch{1.25}
\begin{array}{c|ccc}
{} & e_1 & e_2 & e_3
\\
\hline
e_1 & 0 & -5F_{5,0}e_1+5F_{5,0}e_2 &-e_1+e_2
\\
e_2 & 5F_{5,0}e_1-5F_{5,0}e_2 & 0 & e_1-e_2
\\
e_3 & e_1-e_2 & -e_1+e_2 & 0
\end{array}
\]
for any value of $F_{5,0}$.

%%%%%%%%%%%%%%%%%%%%%%%%%%%%%%%%%%%%%%%%%%%%%%%%%%%%%%%%%%%%%%%%%%%%%%
%%%%%%%%%%%%%%%%%%%%%%%%%%%%%%%%%%%%%%%%%%%%%%%%%%%%%%%%%%%%%%%%%%%%%%
%%%%%%%%%%%%%%%%%%%%%%%%%%%%%%%%%%%%%%%%%%%%%%%%%%%%%%%%%%%%%%%%%%%%%%

\[
\text{\bf Model 2f3n4d}
\ \ \ \ \
\left\{
\aligned
u
&
\,=\,
xy-y^4-6x^2y^2-x^4+\tfrac{1}{2}G_{5,0}y^5+5F_{5,0}xy^4+5G_{5,0}x^2y^3+10F_{5,0}x^3y^2
+
\\
&
\ \ \ \ \
+\tfrac{5}{2}G_{5,0}x^4y+F_{5,0}x^5
+
(-\tfrac{5}{4}F_{5,0}^2-\tfrac{5}{16}G_{5,0}^2)y^6+(24-\tfrac{15}{2}F_{5,0}G_{5,0})xy^5
+
\\
&
\ \ \ \ \
+
(-\tfrac{75}{4}F_{5,0}^2-\tfrac{75}{16}G_{5,0}^2)x^2y^4
+
(-25F_{5,0}G_{5,0}+80)x^3y^3
+
(-\tfrac{75}{4}F_{5,0}^2
+
\\
&
\ \ \ \ \
-\tfrac{75}{16}G_{5,0}^2)x^4y^2
+
(24-\tfrac{15}{2}F_{5,0}G_{5,0})x^5y
+
(-\tfrac{5}{4}F_{5,0}^2-\tfrac{5}{16}G_{5,0}^2)x^6
+
\\
&
\ \ \ \ \
+
(-\tfrac{68}{7}F_{5,0}+\tfrac{75}{28}F_{5,0}^2G_{5,0}+\tfrac{25}{112}G_{5,0}^3)y^7
+
(\tfrac{25}{2}F_{5,0}^3+\tfrac{75}{8}G_{5,0}^2F_{5,0}
+
\\
&
\ \ \ \ \
-34G_{5,0})xy^6
+
(-204F_{5,0}+\tfrac{225}{4}F_{5,0}^2G_{5,0}+\tfrac{75}{16}G_{5,0}^3)x^2y^5
+
(-170G_{5,0}
+
\\
&
\ \ \ \ \
+\tfrac{125}{2}F_{5,0}^3+\tfrac{375}{8}G_{5,0}^2F_{5,0})x^3y^4
+
(-340F_{5,0}+\tfrac{375}{4}F_{5,0}^2G_{5,0}+\tfrac{125}{16}G_{5,0}^3)x^4y^3
+
\\
&
\ \ \ \ \
+
(-102G_{5,0}+\tfrac{75}{2}F_{5,0}^3+\tfrac{225}{8}G_{5,0}^2F_{5,0})x^5y^2
+
(-68F_{5,0}+\tfrac{75}{4}F_{5,0}^2G_{5,0}
+
\\
&
\ \ \ \ \
+\tfrac{25}{16}G_{5,0}^3)x^6y
+
(\tfrac{25}{14}F_{5,0}^3+\tfrac{75}{56}G_{5,0}^2F_{5,0}-\tfrac{34}{7}G_{5,0})x^7
+
\cdots
,
\endaligned\right.
\]

\[
\text{\bf Model 2f3n4d}
\ \ \ \ \
\left\{
\aligned
v
&
\,=\,
y^2+x^2-8xy^3-8x^3y+2F_{5,0}y^5+5G_{5,0}xy^4+20F_{5,0}x^2y^3+10G_{5,0}x^3y^2
+
\\
&
\ \ \ \ \
+10F_{5,0}x^4y+G_{5,0}x^5
+
(-\tfrac{5}{2}F_{5,0}G_{5,0}+8)y^6
+
(-15F_{5,0}^2-\tfrac{15}{4}G_{5,0}^2)xy^5
+
\\
&
\ \ \ \ \
+
(120-\tfrac{75}{2}F_{5,0}G_{5,0})x^2y^4
+
(-\tfrac{25}{2}G_{5,0}^2-50F_{5,0}^2)x^3y^3
+
(120
+
\\
&
\ \ \ \ \
-\tfrac{75}{2}F_{5,0}G_{5,0})x^4y^2
+
(-15F_{5,0}^2-\tfrac{15}{4}G_{5,0}^2)x^5y
+
(-\tfrac{5}{2}F_{5,0}G_{5,0}+8)x^6
+
\\
&
\ \ \ \ \
+
(-\tfrac{68}{7}G_{5,0}+\tfrac{25}{7}F_{5,0}^3+\tfrac{75}{28}G_{5,0}^2F_{5,0})y^7
+
(-136F_{5,0}+\tfrac{75}{2}F_{5,0}^2G_{5,0}
+
\\
&
\ \ \ \ \
+\tfrac{25}{8}G_{5,0}^3)xy^6
+
(-204G_{5,0}+75F_{5,0}^3+\tfrac{225}{4}G_{5,0}^2F_{5,0})x^2y^5
+
(-680F_{5,0}
+
\\
&
\ \ \ \ \
+\tfrac{375}{2}F_{5,0}^2G_{5,0}+\tfrac{125}{8}G_{5,0}^3)x^3y^4
+
(125F_{5,0}^3+\tfrac{375}{4}G_{5,0}^2F_{5,0}-340G_{5,0})x^4y^3
+
\\
&
\ \ \ \ \
+
(-408F_{5,0}+\tfrac{225}{2}F_{5,0}^2G_{5,0}+\tfrac{75}{8}G_{5,0}^3)x^5y^2
+
(-68G_{5,0}+25F_{5,0}^3
+
\\
&
\ \ \ \ \
+\tfrac{75}{4}G_{5,0}^2F_{5,0})x^6y
+
(-\tfrac{136}{7}F_{5,0}+\tfrac{75}{14}F_{5,0}^2G_{5,0}+\tfrac{25}{56}G_{5,0}^3)x^7
+
\cdots
.
\endaligned\right.
\]
%%%%%%%%%%%%%%%%%%%%%%%%%%%%%%%%%%%%%%%%%%%%%%%%%%%%%%%%%%%%%%%%%%%%%%
\[
\def\arraystretch{1.25}
\begin{array}{ll}
e_1
&
\,:=\,
(1+\tfrac{5}{2}xF_{5,0}+\tfrac{5}{4}yG_{5,0}+8u)\partial_x+(\tfrac{5}{4}xG_{5,0}+\tfrac{5}{2}yF_{5,0}+4v)\partial_y
+
\\
&
\ \ \ \ \
+(5uF_{5,0}+y+\tfrac{5}{4}vG_{5,0})\partial_u+(5F_{5,0}v+5G_{5,0}u+2x)\partial_v
, 
\\
e_2
&
\,:=\,
(\tfrac{5}{4}xG_{5,0}+\tfrac{5}{2}yF_{5,0}+4v)\partial_x+(1+\tfrac{5}{2}xF_{5,0}+\tfrac{5}{4}yG_{5,0}+8u)\partial_y
+
\\
&
\ \ \ \ \
+(x+\tfrac{5}{2}G_{5,0}u+\tfrac{5}{2}F_{5,0}v)\partial_u+(10uF_{5,0}+2y+\tfrac{5}{2}vG_{5,0})\partial_v
,
\end{array}
\]
for any value for $F_{5,0},G_{5,0}$.
%%%%%%%%%%%%%%%%%%%%%%%%%%%%%%%%%%%%%%%%%%%%%%%%%%%%%%%%%%%%%%%%%%%%%%
\[
\footnotesize
\def\arraystretch{1.25}
\begin{array}{c|cc}
{} & e_1 & e_2 
\\
\hline
e_1 & 0 & 0
\\
e_2 & 0 & 0
\end{array}
\]

%%%%%%%%%%%%%%%%%%%%%%%%%%%%%%%%%%%%%%%%%%%%%%%%%%%%%%%%%%%%%%%%%%%%%%
%%%%%%%%%%%%%%%%%%%%%%%%%%%%%%%%%%%%%%%%%%%%%%%%%%%%%%%%%%%%%%%%%%%%%%
%%%%%%%%%%%%%%%%%%%%%%%%%%%%%%%%%%%%%%%%%%%%%%%%%%%%%%%%%%%%%%%%%%%%%%

\[
\text{\bf Model 2f3n4e}
\ \ \ \ \
\left\{
\aligned
u
&
\,=\,
xy-4xy^3-4x^3y+\tfrac{1}{2}G_{5,0}y^5+5F_{5,0}xy^4+5G_{5,0}x^2y^3+10F_{5,0}x^3y^2
+
\\
&
\ \ \ \ \
+\tfrac{5}{2}G_{5,0}x^4y+F_{5,0}x^5-\tfrac{5}{4}F_{5,0}G_{5,0}y^6+(24-\tfrac{15}{8}G_{5,0}^2-\tfrac{15}{2}F_{5,0}^2)xy^5
+
\\
&
\ \ \ \ \
-\tfrac{75}{4}F_{5,0}G_{5,0}x^2y^4+(80-25F_{5,0}^2-\tfrac{25}{4}G_{5,0}^2)x^3y^3-\tfrac{75}{4}F_{5,0}G_{5,0}x^4y^2
+
\\
&
\ \ \ \ \
+(24-\tfrac{15}{8}G_{5,0}^2-\tfrac{15}{2}F_{5,0}^2)x^5y-\tfrac{5}{4}F_{5,0}G_{5,0}x^6+(-\tfrac{34}{7}G_{5,0}+\tfrac{25}{112}G_{5,0}^3
+
\\
&
\ \ \ \ \
+\tfrac{75}{28}F_{5,0}^2G_{5,0})y^7+(-68F_{5,0}+\tfrac{75}{8}F_{5,0}G_{5,0}^2+\tfrac{25}{2}F_{5,0}^3)xy^6+(-102G_{5,0}
+
\\
&
\ \ \ \ \
+\tfrac{75}{16}G_{5,0}^3+\tfrac{225}{4}F_{5,0}^2G_{5,0})x^2y^5+(-340F_{5,0}+\tfrac{125}{2}F_{5,0}^3+\tfrac{375}{8}F_{5,0}G_{5,0}^2)x^3y^4
+
\\
&
\ \ \ \ \
+(-170G_{5,0}+\tfrac{375}{4}F_{5,0}^2G_{5,0}+\tfrac{125}{16}G_{5,0}^3)x^4y^3+(-204F_{5,0}+\tfrac{75}{2}F_{5,0}^3
+
\\
&
\ \ \ \ \
+\tfrac{225}{8}F_{5,0}G_{5,0}^2)x^5y^2+(-34G_{5,0}+\tfrac{25}{16}G_{5,0}^3
+\tfrac{75}{4}F_{5,0}^2G_{5,0})x^6y+(-\tfrac{68}{7}F_{5,0}
+
\\
&
\ \ \ \ \
+\tfrac{75}{56}F_{5,0}G_{5,0}^2+\tfrac{25}{14}F_{5,0}^3)x^7
+
\cdots
,
\\
v
&
\,=\,
y^2+x^2-2y^4-12x^2y^2-2x^4+2F_{5,0}y^5+5G_{5,0}xy^4+20F_{5,0}x^2y^3
+
\\
&
\ \ \ \ \
+10G_{5,0}x^3y^2+10F_{5,0}x^4y+G_{5,0}x^5+(8-\tfrac{5}{8}G_{5,0}^2-\tfrac{5}{2}F_{5,0}^2)y^6
+
\\
&
\ \ \ \ \
-15F_{5,0}G_{5,0}xy^5+(120-\tfrac{75}{8}G_{5,0}^2-\tfrac{75}{2}F_{5,0}^2)x^2y^4-50F_{5,0}G_{5,0}x^3y^3
+
\\
&
\ \ \ \ \
+(120-\tfrac{75}{8}G_{5,0}^2-\tfrac{75}{2}F_{5,0}^2)x^4y^2-15F_{5,0}G_{5,0}x^5y+(8-\tfrac{5}{8}G_{5,0}^2-\tfrac{5}{2}F_{5,0}^2)x^6
+
\\
&
\ \ \ \ \
+(-\tfrac{136}{7}F_{5,0}+\tfrac{25}{7}F_{5,0}^3+\tfrac{75}{28}F_{5,0}G_{5,0}^2)y^7+(-68G_{5,0}+\tfrac{25}{8}G_{5,0}^3
+
\\
&
\ \ \ \ \
+\tfrac{75}{2}F_{5,0}^2G_{5,0})xy^6+(-408F_{5,0}+75F_{5,0}^3+\tfrac{225}{4}F_{5,0}G_{5,0}^2)x^2y^5
+
\\
&
\ \ \ \ \
+(-340G_{5,0}+\tfrac{125}{8}G_{5,0}^3+\tfrac{375}{2}F_{5,0}^2G_{5,0})x^3y^4+(125F_{5,0}^3+\tfrac{375}{4}F_{5,0}G_{5,0}^2
+
\\
&
\ \ \ \ \
-680F_{5,0})x^4y^3+(-204G_{5,0}+\tfrac{75}{8}G_{5,0}^3+\tfrac{225}{2}F_{5,0}^2G_{5,0})x^5y^2+(-136F_{5,0}
+
\\
&
\ \ \ \ \
+25F_{5,0}^3+\tfrac{75}{4}F_{5,0}G_{5,0}^2)x^6y+(-\tfrac{68}{7}G_{5,0}+\tfrac{25}{56}G_{5,0}^3+\tfrac{75}{14}F_{5,0}^2G_{5,0})x^7
+
\cdots
.
\endaligned\right.
\]
%%%%%%%%%%%%%%%%%%%%%%%%%%%%%%%%%%%%%%%%%%%%%%%%%%%%%%%%%%%%%%%%%%%%%%
\[
\def\arraystretch{1.25}
\begin{array}{ll}
e_1
&
\,:=\,
(1+\tfrac{5}{4}xG_{5,0}+\tfrac{5}{2}yF_{5,0}+4v)\partial_x+(\tfrac{5}{2}xF_{5,0}+\tfrac{5}{4}yG_{5,0}+8u)\partial_y+(y+\tfrac{5}{2}uG_{5,0}+\tfrac{5}{2}vF_{5,0})\partial_u
+
\\
&
\ \ \ \ \
+(2x+10F_{5,0}u+\tfrac{5}{2}vG_{5,0})\partial_v
,
\\
e_2
&
\,:=\,
(\tfrac{5}{2}xF_{5,0}+\tfrac{5}{4}yG_{5,0}+8u)\partial_x+(1+\tfrac{5}{4}xG_{5,0}+\tfrac{5}{2}yF_{5,0}+4v)\partial_y+(x+5F_{5,0}u+\tfrac{5}{4}vG_{5,0})\partial_u
+
\\
&
\ \ \ \ \
+(5F_{5,0}v+5G_{5,0}u+2y)\partial_v
,
\end{array}
\]
for any value of $F_{5,0},G_{5,0}$.
%%%%%%%%%%%%%%%%%%%%%%%%%%%%%%%%%%%%%%%%%%%%%%%%%%%%%%%%%%%%%%%%%%%%%%
\[
\footnotesize
\def\arraystretch{1.25}
\begin{array}{c|cc}
{} & e_1 & e_2 
\\
\hline
e_1 & 0 & 0
\\
e_2 & 0 & 0
\end{array}
\]

%%%%%%%%%%%%%%%%%%%%%%%%%%%%%%%%%%%%%%%%%%%%%%%%%%%%%%%%%%%%%%%%%%%%%%
%%%%%%%%%%%%%%%%%%%%%%%%%%%%%%%%%%%%%%%%%%%%%%%%%%%%%%%%%%%%%%%%%%%%%%
%%%%%%%%%%%%%%%%%%%%%%%%%%%%%%%%%%%%%%%%%%%%%%%%%%%%%%%%%%%%%%%%%%%%%%

\[
\!\!\!\!\!\!\!\!
\text{\bf Model 2f3n4f}
\ \ \ \ \
\left\{
\aligned
u
&
\,=\,
xy-\tfrac{1}{2}y^4-2xy^3-3x^2y^2-2x^3y-\tfrac{1}{2}x^4+F_{5,0}y^5+5F_{5,0}xy^4+10F_{5,0}x^2y^3
+
\\
&
\ \ \ \ \
+10F_{5,0}x^3y^2+5F_{5,0}x^4y+F_{5,0}x^5+(-\tfrac{5}{2}F_{5,0}^2+2)y^6+(-15F_{5,0}^2+12)xy^5
+
\\
&
\ \ \ \ \
+(-\tfrac{75}{2}F_{5,0}^2+30)x^2y^4+(-50F_{5,0}^2+40)x^3y^3+(-\tfrac{75}{2}F_{5,0}^2+30)x^4y^2
+
\\
&
\ \ \ \ \
+(-15F_{5,0}^2+12)x^5y+(-\tfrac{5}{2}F_{5,0}^2+2)x^6+(\tfrac{50}{7}F_{5,0}^3-\tfrac{68}{7}F_{5,0})y^7+(50F_{5,0}^3
+
\\
&
\ \ \ \ \
-68F_{5,0})xy^6+(150F_{5,0}^3-204F_{5,0})x^2y^5+(250F_{5,0}^3-340F_{5,0})x^3y^4+(250F_{5,0}^3
+
\\
&
\ \ \ \ \
-340F_{5,0})x^4y^3+(150F_{5,0}^3-204F_{5,0})x^5y^2+(50F_{5,0}^3-68F_{5,0})x^6y+(\tfrac{50}{7}F_{5,0}^3
+
\\
&
\ \ \ \ \
-\tfrac{68}{7}F_{5,0})x^7
+
\cdots
,
\\
v
&
\,=\,
y^2+x^2-y^4-4xy^3-6x^2y^2-4x^3y-x^4+2F_{5,0}y^5+10F_{5,0}xy^4
+
\\
&
\ \ \ \ \
+20F_{5,0}x^2y^3+20F_{5,0}x^3y^2
+10F_{5,0}x^4y+2F_{5,0}x^5+(-5F_{5,0}^2+4)y^6
+
\\
&
\ \ \ \ \
+(-30F_{5,0}^2+24)xy^5+(-75F_{5,0}^2+60)x^2y^4
+(-100F_{5,0}^2+80)x^3y^3
+
\\
&
\ \ \ \ \
+(-75F_{5,0}^2+60)x^4y^2
+(-30F_{5,0}^2+24)x^5y+(-5F_{5,0}^2+4)x^6
+(\tfrac{100}{7}F_{5,0}^3
+
\\
&
\ \ \ \ \
-\tfrac{136}{7}F_{5,0})y^7+(100F_{5,0}^3-136F_{5,0})xy^6+(300F_{5,0}^3
-408F_{5,0})x^2y^5+(500F_{5,0}^3
+
\\
&
\ \ \ \ \
-680F_{5,0})x^3y^4+(500F_{5,0}^3-680F_{5,0})x^4y^3+(300F_{5,0}^3
+
\\
&
\ \ \ \ \
-408F_{5,0})x^5y^2+(100F_{5,0}^3-136F_{5,0})x^6y+(\tfrac{100}{7}F_{5,0}^3-\tfrac{136}{7}F_{5,0})x^7
+
\cdots
.
\endaligned\right.
\]
%%%%%%%%%%%%%%%%%%%%%%%%%%%%%%%%%%%%%%%%%%%%%%%%%%%%%%%%%%%%%%%%%%%%%%
\[
\def\arraystretch{1.25}
\begin{array}{ll}
e_1
&
\,:=\,
(5F_{5,0}x+4u+2v+1)\partial_x+(5F_{5,0}y+4u+2v)\partial_y+(10F_{5,0}u+y)\partial_u+(10F_{5,0}v+2x)\partial_v
,
\\
e_2
&
\,:=\,
(5F_{5,0}x+4u+2v)\partial_x+(5F_{5,0}y+4u+2v+1)\partial_y+(10F_{5,0}u+x)\partial_u+(10F_{5,0}v+2y)\partial_v
,
\\
e_3
&
\,:=\,
-(x-y)\partial_x+(x-y)\partial_y-(2u-v)\partial_u+(4u-2v)\partial_v
.
\end{array}
\]

%%%%%%%%%%%%%%%%%%%%%%%%%%%%%%%%%%%%%%%%%%%%%%%%%%%%%%%%%%%%%%%%%%%%%%
\[
\footnotesize
\def\arraystretch{1.25}
\begin{array}{c|ccc}
{} & e_1 & e_2 & e_3
\\
\hline
e_1 & 0 & 5F_{5,0}e_1-5F_{5,0}e_2 & -e_1+e_2
\\
e_2 & -5F_{5,0}e_1+5F_{5,0}e_2 & 0 & e_1-e_2
\\
e_3 & e_1-e_2 & -e_1+e_2 & 0
\end{array}
\]
for any value for $F_{5,0}$.

%%%%%%%%%%%%%%%%%%%%%%%%%%%%%%%%%%%%%%%%%%%%%%%%%%%%%%%%%%%%%%%%%%%%%%
%%%%%%%%%%%%%%%%%%%%%%%%%%%%%%%%%%%%%%%%%%%%%%%%%%%%%%%%%%%%%%%%%%%%%%
%%%%%%%%%%%%%%%%%%%%%%%%%%%%%%%%%%%%%%%%%%%%%%%%%%%%%%%%%%%%%%%%%%%%%%

\[
\text{\bf Model 2f3n4g}
\ \ \ \ \
\left\{
\aligned
u
&
\,=\,
xy-\tfrac{1}{2}y^4+2xy^3-3x^2y^2+2x^3y-\tfrac{1}{2}x^4-F_{5,0}y^5+5F_{5,0}xy^4
+
\\
&
\ \ \ \ \
-10F_{5,0}x^2y^3+10F_{5,0}x^3y^2-5F_{5,0}x^4y+F_{5,0}x^5+(-\tfrac{5}{2}F_{5,0}^2-2)y^6
+
\\
&
\ \ \ \ \
+(15F_{5,0}^2+12)xy^5+(-\tfrac{75}{2}F_{5,0}^2-30)x^2y^4+(50F_{5,0}^2+40)x^3y^3
+
\\
&
\ \ \ \ \
+(-\tfrac{75}{2}F_{5,0}^2-30)x^4y^2+(15F_{5,0}^2+12)x^5y+(-\tfrac{5}{2}F_{5,0}^2-2)x^6
+
\\
&
\ \ \ \ \
+(-\tfrac{50}{7}F_{5,0}^3-\tfrac{68}{7}F_{5,0})y^7+(50F_{5,0}^3+68F_{5,0})xy^6+(-150F_{5,0}^3
+
\\
&
\ \ \ \ \
-204F_{5,0})x^2y^5+(250F_{5,0}^3+340F_{5,0})x^3y^4+(-250F_{5,0}^3-340F_{5,0})x^4y^3
+
\\
&
\ \ \ \ \
+(150F_{5,0}^3+204F_{5,0})x^5y^2+(-50F_{5,0}^3-68F_{5,0})x^6y+(\tfrac{50}{7}F_{5,0}^3
+
\\
&
\ \ \ \ \
+\tfrac{68}{7}F_{5,0})x^7
+
\cdots
,
\\
v
&
\,=\,
y^2+x^2+y^4-4xy^3+6x^2y^2-4x^3y+x^4+2F_{5,0}y^5-10F_{5,0}xy^4
+
\\
&
\ \ \ \ \
+20F_{5,0}x^2y^3
-20F_{5,0}x^3y^2+10F_{5,0}x^4y-2F_{5,0}x^5+(5F_{5,0}^2+4)y^6
+
\\
&
\ \ \ \ \
+(-30F_{5,0}^2-24)xy^5+(75F_{5,0}^2+60)x^2y^4+(-100F_{5,0}^2-80)x^3y^3
+
\\
&
\ \ \ \ \
+(75F_{5,0}^2+60)x^4y^2+(-30F_{5,0}^2-24)x^5y+(5F_{5,0}^2+4)x^6
+
\\
&
\ \ \ \ \
+(\tfrac{100}{7}F_{5,0}^3+\tfrac{136}{7}F_{5,0})y^7+(-100F_{5,0}^3
+
\\
&
\ \ \ \ \
-136F_{5,0})xy^6+(300F_{5,0}^3+408F_{5,0})x^2y^5+(-500F_{5,0}^3-680F_{5,0})x^3y^4
+
\\
&
\ \ \ \ \
v+(500F_{5,0}^3+680F_{5,0})x^4y^3+(-300F_{5,0}^3-408F_{5,0})x^5y^2+(100F_{5,0}^3
+
\\
&
\ \ \ \ \
+136F_{5,0})x^6y+(-\tfrac{100}{7}F_{5,0}^3-\tfrac{136}{7}F_{5,0})x^7
+
\cdots
.
\endaligned\right.
\]
%%%%%%%%%%%%%%%%%%%%%%%%%%%%%%%%%%%%%%%%%%%%%%%%%%%%%%%%%%%%%%%%%%%%%%
\[
\def\arraystretch{1.25}
\begin{array}{ll}
e_1
&
\,:=\,
(5F_{5,0}x+4u-2v+1)\partial_x+(5F_{5,0}y-4u+2v)\partial_y+(10F_{5,0}u+y)\partial_u+(10F_{5,0}v+2x)\partial_v
,
\\
e_2
&
\,:=\,
-(5F_{5,0}x+4u-2v)\partial_x-(5F_{5,0}y-4u+2v-1)\partial_y-(10F_{5,0}u-x)\partial_u-(10F_{5,0}v-2y)\partial_v
, 
\\
e_3
&
\,:=\,
(x+y)\partial_x+(x+y)\partial_y+(2u+v)\partial_u+(4u+2v)\partial_v
.
\end{array}
\]

%%%%%%%%%%%%%%%%%%%%%%%%%%%%%%%%%%%%%%%%%%%%%%%%%%%%%%%%%%%%%%%%%%%%%%
\[
\footnotesize
\def\arraystretch{1.25}
\begin{array}{c|ccc}
{} & e_1 & e_2 & e_3
\\
\hline
e_1 & 0 & -5F_{5,0}e_1-5F_{5,0}e_2 & e_1+e_2
\\
e_2 & 5F_{5,0}e_1+5F_{5,0}e_2 & 0 & e_1+e_2
\\
e_3 & -e_1-e_2 & -e_1-e_2 & 0
\end{array}
\]
for any value for $F_{5,0}$.

%%%%%%%%%%%%%%%%%%%%%%%%%%%%%%%%%%%%%%%%%%%%%%%%%%%%%%%%%%%%%%%%%%%%%%
%%%%%%%%%%%%%%%%%%%%%%%%%%%%%%%%%%%%%%%%%%%%%%%%%%%%%%%%%%%%%%%%%%%%%%
%%%%%%%%%%%%%%%%%%%%%%%%%%%%%%%%%%%%%%%%%%%%%%%%%%%%%%%%%%%%%%%%%%%%%%

\[
\text{\bf Model 2f3n4h}
\ \ \ \ \
\left\{
\aligned
u
&
\,=\,
xy+\tfrac{1}{2}y^4-2xy^3+3x^2y^2-2x^3y+\tfrac{1}{2}x^4-F_{5,0}y^5+5F_{5,0}xy^4
+
\\
&
\ \ \ \ \
-10F_{5,0}x^2y^3+10F_{5,0}x^3y^2-5F_{5,0}x^4y+F_{5,0}x^5+(\tfrac{5}{2}F_{5,0}^2-2)y^6+
\\
&
\ \ \ \ \
+(-15F_{5,0}^2+12)xy^5+(\tfrac{75}{2}F_{5,0}^2-30)x^2y^4+(-50F_{5,0}^2+40)x^3y^3
+
\\
&
\ \ \ \ \
+(\tfrac{75}{2}F_{5,0}^2-30)x^4y^2+(-15F_{5,0}^2+12)x^5y+(\tfrac{5}{2}F_{5,0}^2-2)x^6
+
\\
&
\ \ \ \ \
+(-\tfrac{50}{7}F_{5,0}^3+\tfrac{68}{7}F_{5,0})y^7+(50F_{5,0}^3-68F_{5,0})xy^6+(-150F_{5,0}^3
+
\\
&
\ \ \ \ \
+204F_{5,0})x^2y^5+(250F_{5,0}^3-340F_{5,0})x^3y^4+(-250F_{5,0}^3+340F_{5,0})x^4y^3
+
\\
&
\ \ \ \ \
+(150F_{5,0}^3-204F_{5,0})x^5y^2+(-50F_{5,0}^3+68F_{5,0})x^6y+(\tfrac{50}{7}F_{5,0}^3
+
\\
&
\ \ \ \ \
-\tfrac{68}{7}F_{5,0})x^7
+
\cdots
,
\\
v
&
\,=\,
y^2+x^2-y^4+4xy^3-6x^2y^2+4x^3y-x^4+2F_{5,0}y^5-10F_{5,0}xy^4
+
\\
&
\ \ \ \ \
+20F_{5,0}x^2y^3-20F_{5,0}x^3y^2+10F_{5,0}x^4y-2F_{5,0}x^5+(-5F_{5,0}^2+4)y^6
+
\\
&
\ \ \ \ \
+(30F_{5,0}^2-24)xy^5+(-75F_{5,0}^2+60)x^2y^4+(100F_{5,0}^2-80)x^3y^3
+
\\
&
\ \ \ \ \
+(-75F_{5,0}^2+60)x^4y^2+(30F_{5,0}^2-24)x^5y+(-5F_{5,0}^2+4)x^6
+
\\
&
\ \ \ \ \
+(\tfrac{100}{7}F_{5,0}^3-\tfrac{136}{7}F_{5,0})y^7+(-100F_{5,0}^3+136F_{5,0})xy^6+(300F_{5,0}^3
+
\\
&
\ \ \ \ \
-408F_{5,0})x^2y^5+(-500F_{5,0}^3+680F_{5,0})x^3y^4+(500F_{5,0}^3-680F_{5,0})x^4y^3
+
\\
&
\ \ \ \ \
+(-300F_{5,0}^3+408F_{5,0})x^5y^2+(100F_{5,0}^3-136F_{5,0})x^6y+(-\tfrac{100}{7}F_{5,0}^3
+
\\
&
\ \ \ \ \
+\tfrac{136}{7}F_{5,0})x^7
+
\cdots
.
\endaligned\right.
\]
%%%%%%%%%%%%%%%%%%%%%%%%%%%%%%%%%%%%%%%%%%%%%%%%%%%%%%%%%%%%%%%%%%%%%%
\[
\def\arraystretch{1.25}
\begin{array}{ll}
e_1
&
\,:=\,
-(5F_{5,0}x+4u-2v-1)\partial_x-(5F_{5,0}y-4u+2v)\partial_y-(10F_{5,0}u-y)\partial_u-(10F_{5,0}v-2x)\partial_v
,
\\
e_2
&
\,:=\,
(5F_{5,0}x+4u-2v)\partial_x+(5F_{5,0}y-4u+2v+1)\partial_y+(10F_{5,0}u+x)\partial_u+(10F_{5,0}v+2y)\partial_v
, 
\\
e_3
&
\,:=\,
(x+y)\partial_x+(x+y)\partial_y+(2u+v)\partial_u+(4u+2v)\partial_v
.
\end{array}
\]

\[
\footnotesize
\def\arraystretch{1.25}
\begin{array}{c|ccc}
{} & e_1 & e_2 & e_3
\\
\hline
e_1 & 0 & 5F_{5,0}e_1+5F_{5,0}e_2 & e_1+e_2
\\
e_2 & -5F_{5,0}e_1-5F_{5,0}e_2 & 0 & e_1+e_2
\\
e_3 & -e_1-e_2 & -e_1-e_2 & 0
\end{array}
\]
for any value for $F_{5,0}$.

%%%%%%%%%%%%%%%%%%%%%%%%%%%%%%%%%%%%%%%%%%%%%%%%%%%%%%%%%%%%%%%%%%%%%%
%%%%%%%%%%%%%%%%%%%%%%%%%%%%%%%%%%%%%%%%%%%%%%%%%%%%%%%%%%%%%%%%%%%%%%
%%%%%%%%%%%%%%%%%%%%%%%%%%%%%%%%%%%%%%%%%%%%%%%%%%%%%%%%%%%%%%%%%%%%%%

\[
\text{\bf Model 2f3n4i}
\ \ \ \ \
\left\{
\aligned
u
&
\,=\,
xy
,
\\
v
&
\,=\,
x^2+y^2
.
\endaligned\right.
\]
%%%%%%%%%%%%%%%%%%%%%%%%%%%%%%%%%%%%%%%%%%%%%%%%%%%%%%%%%%%%%%%%%%%%%%
\[
\def\arraystretch{1.25}
\begin{array}{llll}
e_1
\,:=\,
y\partial_u+2x\partial_v+\partial_x
, &
e_2
\,:=\,
x\partial_u+2y\partial_v+\partial_y
, &
& 
\\
e_3
\,:=\,
2u\partial_u+2v\partial_v+x\partial_x+y\partial_y
, &
e_4
\,:=\,
v\partial_u+4u\partial_v+y\partial_x+x\partial_y
.
\end{array}
\]

%%%%%%%%%%%%%%%%%%%%%%%%%%%%%%%%%%%%%%%%%%%%%%%%%%%%%%%%%%%%%%%%%%%%%%
\[
\footnotesize
\def\arraystretch{1.25}
\begin{array}{c|cccc}
{} & e_1 & e_2 & e_3 & e_4
\\
\hline
e_1 & 0 & 0 & e_1 & e_2 
\\
e_2 & 0 & 0 & e_2 & e_1
\\
e_3 & -e_1 & -e_2 & 0 & 0
\\
e_4 & -e_2 & -e_1 & 0 & 0
\end{array}
\]

%%%%%%%%%%%%%%%%%%%%%%%%%%%%%%%%%%%%%%%%%%%%%%%%%%%%%%%%%%%%%%%%%%%%%%
%%%%%%%%%%%%%%%%%%%%%%%%%%%%%%%%%%%%%%%%%%%%%%%%%%%%%%%%%%%%%%%%%%%%%%
%%%%%%%%%%%%%%%%%%%%%%%%%%%%%%%%%%%%%%%%%%%%%%%%%%%%%%%%%%%%%%%%%%%%%%

%%%%%%%%%%%%%%%%%%%%%%%%%%%%%%%%%%%%%%%%%%%%%%%%%%%%%%%%%%%%%%%%%%%%%%
\SectionHead{2g Models}
{2g-models}
%%%%%%%%%%%%%%%%%%%%%%%%%%%%%%%%%%%%%%%%%%%%%%%%%%%%%%%%%%%%%%%%%%%%%%

%%%%%%%%%%%%%%%%%%%%%%%%%%%%%%%%%%%%%%%%%%%%%%%%%%%%%%%%%%%%%%%%%%%%%%
%%%%%%%%%%%%%%%%%%%%%%%%%%%%%%%%%%%%%%%%%%%%%%%%%%%%%%%%%%%%%%%%%%%%%%
%%%%%%%%%%%%%%%%%%%%%%%%%%%%%%%%%%%%%%%%%%%%%%%%%%%%%%%%%%%%%%%%%%%%%%

\[
\!\!\!\!\!\!
\text{\bf Model 2g3a}
\ \ \ \ \
\left\{
\aligned
u
&
\,=\,
xy+F_{0,3}y^3+F_{3,0}x^3+F_{0,4}y^4+F_{1,3}xy^3+F_{2,2}x^2y^2+F_{3,1}x^3y+F_{4,0}x^4
+
\cdots,
\\
v
&
\,=\,
-y^2+x^2-6F_{0,3}xy^2+(6F_{3,0}-8)x^2y+G_{0,4}y^4+G_{1,3}xy^3+G_{2,2}x^2y^2
\\
&
\ \ \ \ \
+
G_{3,1}x^3y+G_{4,0}x^4
+
\cdots,
\endaligned\right.
\]
%%%%%%%%%%%%%%%%%%%%%%%%%%%%%%%%%%%%%%%%%%%%%%%%%%%%%%%%%%%%%%%%%%%%%%
\[
\def\arraystretch{1.25}
\begin{array}{ll}
e_1
&
\,:=\,
(1-\tfrac{9}{2}xF_{0,3}F_{3,0}+\tfrac{15}{4}vF_{3,0}F_{3,1}+\tfrac{21}{2}uF_{0,3}F_{3,1}+\tfrac{15}{4}vF_{3,0}F_{1,3}+21uF_{0,3}G_{4,0}+\tfrac{5}{4}vF_{3,0}G_{2,2}
+
\\
&
\ \ \ \ \
+228uF_{0,3}F_{3,0}+\tfrac{15}{2}vF_{3,0}G_{4,0}+\tfrac{21}{2}uF_{0,3}F_{1,3}-189uF_{0,3}F_{3,0}^2+\tfrac{7}{2}uF_{0,3}G_{2,2}-8uF_{0,3}
+
\\
&
\ \ \ \ \
-80vF_{3,0}+9xF_{0,3}+\tfrac{3}{8}xG_{1,3}-\tfrac{1}{2}xF_{2,2}-3xF_{4,0}-\tfrac{135}{2}vF_{3,0}^3+\tfrac{1}{4}yG_{2,2}+\tfrac{3}{4}yF_{1,3}+\tfrac{3}{2}yG_{4,0}
+
\\
&
\ \ \ \ \
+\tfrac{3}{8}xG_{3,1}-\tfrac{27}{2}yF_{3,0}^2+\tfrac{3}{4}yF_{3,1}+15yF_{3,0}-2uF_{2,2}-\tfrac{1}{2}uG_{1,3}+147vF_{3,0}^2-3vF_{1,3}-3vF_{3,1}
+
\\
&
\ \ \ \ \
-8vG_{4,0}-vG_{2,2})\partial_x+(\tfrac{15}{4}vF_{0,3}F_{1,3}-\tfrac{21}{2}uF_{3,0}F_{3,1}-\tfrac{9}{2}yF_{0,3}F_{3,0}-\tfrac{135}{2}vF_{0,3}F_{3,0}^2
+
\\
&
\ \ \ \ \
+\tfrac{5}{4}vF_{0,3}G_{2,2}-\tfrac{21}{2}uF_{3,0}F_{1,3}-\tfrac{7}{2}uF_{3,0}G_{2,2}-21uF_{3,0}G_{4,0}+84vF_{0,3}F_{3,0}+\tfrac{15}{2}vF_{0,3}G_{4,0}
+
\\
&
\ \ \ \ \
+\tfrac{15}{4}vF_{0,3}F_{3,1}-\tfrac{3}{4}xF_{3,1}-\tfrac{1}{4}xG_{2,2}-\tfrac{3}{2}xG_{4,0}-\tfrac{1}{2}yF_{2,2}-3yF_{4,0}+\tfrac{3}{8}yG_{1,3}+\tfrac{3}{8}yG_{3,1}
+
\\
&
\ \ \ \ \
+189uF_{3,0}^3-\tfrac{1}{2}vG_{1,3}-390uF_{3,0}^2+9uF_{1,3}+6uF_{3,1}+3uG_{2,2}+20uG_{4,0}+\tfrac{27}{2}xF_{3,0}^2
+
\\
&
\ \ \ \ \
-\tfrac{3}{4}xF_{1,3}-18xF_{3,0}+6yF_{0,3}+200uF_{3,0}-8vF_{0,3})\partial_y-(uF_{2,2}+6uF_{4,0}-\tfrac{3}{4}uG_{1,3}
+
\\
&
\ \ \ \ \
-\tfrac{3}{4}uG_{3,1}+15vF_{3,0}-\tfrac{27}{2}vF_{3,0}^2+\tfrac{3}{4}vF_{1,3}+\tfrac{3}{4}vF_{3,1}+\tfrac{1}{4}vG_{2,2}+\tfrac{3}{2}vG_{4,0}-15uF_{0,3}
+
\\
&
\ \ \ \ \
+9uF_{0,3}F_{3,0}-y)\partial_u-(-18vF_{0,3}+vF_{2,2}+6vF_{4,0}-\tfrac{3}{4}vG_{1,3}-\tfrac{3}{4}vG_{3,1}+16u+54uF_{3,0}^2
+
\\
&
\ \ \ \ \
-3uF_{1,3}-78uF_{3,0}-3uF_{3,1}-uG_{2,2}-6uG_{4,0}+9vF_{0,3}F_{3,0}-2x)\partial_v 
, 
\\
e_2
&
\,:=\,
(-3yF_{0,3}-2vF_{2,2}-\tfrac{3}{4}xF_{3,1}+\tfrac{1}{4}xG_{2,2}+\tfrac{1}{2}yF_{2,2}+\tfrac{3}{8}yG_{1,3}+\tfrac{3}{8}yG_{3,1}-\tfrac{3}{2}vG_{1,3}-3uF_{1,3}
+
\\
&
\ \ \ \ \
+63uF_{3,0}F_{0,3}^2+\tfrac{15}{8}vF_{3,0}G_{3,1}+42uF_{0,3}F_{0,4}+\tfrac{21}{4}uF_{0,3}G_{1,3}+\tfrac{5}{2}vF_{3,0}F_{2,2}+15vF_{3,0}F_{0,4}
+
\\
&
\ \ \ \ \
+\tfrac{45}{2}vF_{0,3}F_{3,0}^2+\tfrac{15}{8}vF_{3,0}G_{1,3}+7uF_{0,3}F_{2,2}+\tfrac{21}{4}uF_{0,3}G_{3,1}+\tfrac{9}{2}yF_{0,3}F_{3,0}-27vF_{0,3}F_{3,0}
+
\\
&
\ \ \ \ \
-4x-\tfrac{3}{4}xF_{1,3}+3xF_{3,0}-2vG_{3,1}+3yF_{0,4}+\tfrac{27}{2}xF_{0,3}^2+\tfrac{3}{2}xG_{0,4}-2uG_{0,4}-12vF_{0,4}
+
\\
&
\ \ \ \ \
-18uF_{0,3}^2)\partial_x
+
(1+15vF_{0,3}F_{0,4}-42uF_{3,0}F_{0,4}-7uF_{3,0}F_{2,2}+\tfrac{15}{8}vF_{0,3}G_{3,1}+\tfrac{45}{2}vF_{3,0}F_{0,3}^2
+
\\
&
\ \ \ \ \
+\tfrac{5}{2}vF_{0,3}F_{2,2}+\tfrac{15}{8}vF_{0,3}G_{1,3}-63uF_{0,3}F_{3,0}^2-\tfrac{21}{4}uF_{3,0}G_{3,1}-\tfrac{21}{4}uF_{3,0}G_{1,3}+72uF_{0,3}F_{3,0}
+
\\
&
\ \ \ \ \
-\tfrac{9}{2}xF_{0,3}F_{3,0}-2vG_{0,4}+\tfrac{27}{2}yF_{0,3}^2+36uF_{0,4}-3xF_{0,4}-8y-\tfrac{3}{8}xG_{1,3}-\tfrac{1}{2}xF_{2,2}+\tfrac{1}{4}yG_{2,2}
+
\\
&
\ \ \ \ \
-\tfrac{3}{4}yF_{1,3}-\tfrac{3}{8}xG_{3,1}-\tfrac{3}{4}yF_{3,1}+6yF_{3,0}+4uF_{2,2}+\tfrac{9}{2}uG_{1,3}+5uG_{3,1}+\tfrac{3}{2}yG_{0,4}-18vF_{0,3}^2)\partial_y
+
\\
&
\ \ \ \ \
+
(-12u+9uF_{3,0}-\tfrac{3}{2}uF_{3,1}+3uG_{0,4}+\tfrac{1}{2}uG_{2,2}-3vF_{0,4}-\tfrac{1}{2}vF_{2,2}-\tfrac{3}{8}vG_{1,3}-\tfrac{3}{8}vG_{3,1}
+
\\
&
\ \ \ \ \
+27uF_{0,3}^2-\tfrac{3}{2}uF_{1,3}-\tfrac{9}{2}vF_{0,3}F_{3,0}+x)\partial_u
+
(27vF_{0,3}^2-\tfrac{3}{2}vF_{1,3}-\tfrac{3}{2}vF_{3,1}+3vG_{0,4}+\tfrac{1}{2}vG_{2,2}
+
\\
&
\ \ \ \ \
-16v+12uF_{0,4}+2uF_{2,2}+\tfrac{3}{2}uG_{1,3}+\tfrac{3}{2}uG_{3,1}-18uF_{0,3}+12vF_{3,0}+18uF_{0,3}F_{3,0}-2y)\partial_v.
\end{array}
\]

%%%%%%%%%%%%%%%%%%%%%%%%%%%%%%%%%%%%%%%%%%%%%%%%%%%%%%%%%%%%%%%%%%%%%%
\noindent
A Gr\"obner basis has 433 generators to describe the
moduli space core algebraic variety in 
$\R^{12} \ni F_{0,3},F_{0,4},F_{1,3},F_{2,2},F_{3,0},F_{3,1},F_{4,0},G_{0,4},G_{1,3},G_{2,2},G_{3,1},G_{4,0}$.
%%%%%%%%%%%%%%%%%%%%%%%%%%%%%%%%%%%%%%%%%%%%%%%%%%%%%%%%%%%%%%%%%%%%%%
\[
\footnotesize
\def\arraystretch{1.25}
\begin{array}{c|cc}
{} & e_1 & e_2 
\\
\hline
e_1 & 0 
&
\rotatebox[origin=c]{0}{
\begin{tabular}{p{7cm}}
$(-4+\tfrac{27}{2}F_{3,0}^2+\tfrac{27}{2}F_{0,3}^2-12F_{3,0}-\tfrac{3}{2}F_{1,3}-\tfrac{3}{2}F_{3,1}+\tfrac{3}{2}G_{0,4}-\tfrac{3}{2}G_{4,0})e_1+(-6F_{0,3}+3F_{4,0}-\tfrac{3}{4}G_{1,3}-\tfrac{3}{4}G_{3,1}-3F_{0,4})e_2$
\end{tabular}}
\\
e_2 & \rotatebox[origin=c]{0}{
\begin{tabular}{p{7cm}}
$-(-4+\tfrac{27}{2}F_{3,0}^2+\tfrac{27}{2}F_{0,3}^2-12F_{3,0}-\tfrac{3}{2}F_{1,3}-\tfrac{3}{2}F_{3,1}+\tfrac{3}{2}G_{0,4}-\tfrac{3}{2}G_{4,0})e_1-(-6F_{0,3}+3F_{4,0}-\tfrac{3}{4}G_{1,3}-\tfrac{3}{4}G_{3,1}-3F_{0,4})e_2$
\end{tabular}}
& 0
\end{array}
\]

%%%%%%%%%%%%%%%%%%%%%%%%%%%%%%%%%%%%%%%%%%%%%%%%%%%%%%%%%%%%%%%%%%%%%%
%%%%%%%%%%%%%%%%%%%%%%%%%%%%%%%%%%%%%%%%%%%%%%%%%%%%%%%%%%%%%%%%%%%%%%
%%%%%%%%%%%%%%%%%%%%%%%%%%%%%%%%%%%%%%%%%%%%%%%%%%%%%%%%%%%%%%%%%%%%%%

\[
\text{\bf Model 2g3b}
\ \ \ \ \
\left\{
\aligned
u
&
\,=\,
xy+x^3+F_{0,4}y^4+F_{1,3}xy^3+F_{2,2}x^2y^2+(-F_{1,3}-G_{4,0}+9+G_{0,4})x^3y
+
\\
&
\ \ \ \ \
+(-\tfrac{1}{3}F_{2,2}-F_{0,4})x^4+(-\tfrac{1}{50}F_{1,3}G_{3,1}-\tfrac{14}{25}F_{0,4}F_{1,3}+\tfrac{3}{5}F_{0,4}-\tfrac{14}{25}F_{0,4}G_{0,4}
+
\\
&
\ \ \ \ \
-\tfrac{2}{5}F_{0,4}G_{4,0}-\tfrac{1}{50}G_{0,4}G_{3,1})y^5+(\tfrac{18}{5}F_{1,3}-\tfrac{1}{5}F_{1,3}^2-\tfrac{2}{5}F_{1,3}G_{4,0}-\tfrac{4}{3}F_{0,4}F_{2,2}
+
\\
&
\ \ \ \ \
+F_{0,4}G_{3,1}+\tfrac{18}{5}G_{0,4}-\tfrac{1}{5}G_{0,4}F_{1,3}-\tfrac{2}{5}G_{0,4}G_{4,0})xy^4+(\tfrac{6}{5}F_{2,2}-\tfrac{9}{5}G_{3,1}
+
\\
&
\ \ \ \ \
-\tfrac{78}{5}F_{0,4}-\tfrac{2}{15}F_{2,2}G_{4,0}+\tfrac{1}{5}G_{4,0}G_{3,1}+\tfrac{3}{5}F_{1,3}G_{3,1}-\tfrac{11}{15}F_{2,2}F_{1,3}+\tfrac{6}{5}F_{0,4}F_{1,3}
+
\\
&
\ \ \ \ \
+\tfrac{12}{5}F_{0,4}G_{4,0})x^2y^3+(\tfrac{324}{5}-\tfrac{44}{5}F_{1,3}-\tfrac{72}{5}G_{4,0}-\tfrac{4}{9}F_{2,2}^2+\tfrac{1}{3}F_{2,2}G_{3,1}
+
\\
&
\ \ \ \ \
+\tfrac{6}{5}F_{1,3}G_{4,0}+\tfrac{2}{5}F_{1,3}^2+\tfrac{4}{5}G_{4,0}^2)x^3y^2+(-\tfrac{25}{6}F_{2,2}+\tfrac{19}{10}G_{3,1}-\tfrac{68}{5}F_{0,4}
+
\\
&
\ \ \ \ \
+\tfrac{2}{3}F_{2,2}G_{4,0}-\tfrac{7}{20}G_{4,0}G_{3,1}-\tfrac{3}{10}F_{1,3}G_{3,1}+\tfrac{1}{2}F_{2,2}F_{1,3}+\tfrac{1}{5}F_{0,4}F_{1,3}
+
\\
&
\ \ \ \ \
+\tfrac{2}{5}F_{0,4}G_{4,0}-\tfrac{1}{3}F_{2,2}G_{0,4}+\tfrac{1}{4}G_{0,4}G_{3,1})x^4y+(\tfrac{-54}{25}-\tfrac{4}{25}F_{1,3}+\tfrac{42}{25}G_{4,0}
+
\\
&
\ \ \ \ \
+\tfrac{12}{25}G_{0,4}+\tfrac{4}{45}F_{2,2}^2-\tfrac{1}{15}F_{2,2}G_{3,1}+\tfrac{4}{15}F_{0,4}F_{2,2}-\tfrac{1}{5}F_{0,4}G_{3,1}-\tfrac{4}{25}F_{1,3}G_{4,0}
+
\\
&
\ \ \ \ \
+\tfrac{1}{25}G_{0,4}F_{1,3}+\tfrac{2}{25}G_{0,4}G_{4,0}-\tfrac{1}{25}F_{1,3}^2-\tfrac{4}{25}G_{4,0}^2)x^5
+
\cdots
,
\\
v
&
\,=\,
-y^2+x^2+6x^2y+G_{0,4}y^4+(-8F_{0,4}-\tfrac{4}{3}F_{2,2}-G_{3,1})xy^3+(-3G_{4,0}+27
+
\\
&
\ \ \ \ \
-3G_{0,4})x^2y^2+G_{3,1}x^3y+G_{4,0}x^4+(\tfrac{6}{5}G_{0,4}-\tfrac{2}{5}G_{0,4}G_{4,0}+\tfrac{2}{5}F_{0,4}G_{3,1}
+
\\
&
\ \ \ \ \
+\tfrac{2}{75}F_{2,2}G_{3,1}+\tfrac{16}{75}F_{0,4}F_{2,2}-\tfrac{2}{5}G_{0,4}F_{1,3}+\tfrac{1}{50}G_{3,1}^2+\tfrac{48}{25}F_{0,4}^2-\tfrac{2}{5}G_{0,4}^2)y^5
+
\\
&
\ \ \ \ \
+(-\tfrac{24}{5}F_{2,2}-\tfrac{18}{5}G_{3,1}-\tfrac{126}{5}F_{0,4}+\tfrac{8}{15}F_{2,2}G_{4,0}+\tfrac{2}{5}G_{4,0}G_{3,1}+\tfrac{1}{5}F_{1,3}G_{3,1}
+
\\
&
\ \ \ \ \
+\tfrac{4}{15}F_{2,2}F_{1,3}+\tfrac{12}{5}F_{0,4}F_{1,3}+\tfrac{24}{5}F_{0,4}G_{4,0}-\tfrac{4}{3}F_{2,2}G_{0,4}+G_{0,4}G_{3,1})xy^4
+
\\
&
\ \ \ \ \
+(\tfrac{486}{5}+\tfrac{12}{5}F_{1,3}-\tfrac{108}{5}G_{4,0}-12G_{0,4}+\tfrac{8}{9}F_{2,2}^2+\tfrac{16}{3}F_{0,4}F_{2,2}-4F_{0,4}G_{3,1}
+
\\
&
\ \ \ \ \
+\tfrac{7}{5}F_{1,3}G_{4,0}+G_{0,4}F_{1,3}+2G_{0,4}G_{4,0}+\tfrac{2}{5}F_{1,3}^2+\tfrac{6}{5}G_{4,0}^2-\tfrac{1}{2}G_{3,1}^2)x^2y^3
+
\\
&
\ \ \ \ \
+(-\tfrac{2}{3}F_{2,2}+\tfrac{71}{5}G_{3,1}+\tfrac{64}{5}F_{0,4}+\tfrac{4}{3}F_{2,2}G_{4,0}-\tfrac{9}{5}G_{4,0}G_{3,1}-\tfrac{2}{5}F_{1,3}G_{3,1}
+
\\
&
\ \ \ \ \
-\tfrac{8}{5}F_{0,4}F_{1,3}-\tfrac{16}{5}F_{0,4}G_{4,0}+\tfrac{4}{3}F_{2,2}G_{0,4}-G_{0,4}G_{3,1})x^3y^2+(\tfrac{-27}{5}
+
\\
&
\ \ \ \ \
-\tfrac{13}{5}F_{1,3}+\tfrac{66}{5}G_{4,0}+\tfrac{24}{5}G_{0,4}-\tfrac{1}{3}F_{2,2}G_{3,1}-\tfrac{11}{10}F_{1,3}G_{4,0}-\tfrac{1}{10}G_{0,4}F_{1,3}
+
\\
&
\ \ \ \ \
-\tfrac{1}{5}G_{0,4}G_{4,0}-\tfrac{1}{5}F_{1,3}^2-\tfrac{7}{5}G_{4,0}^2+\tfrac{1}{4}G_{3,1}^2)x^4y+(-\tfrac{38}{75}F_{2,2}-\tfrac{23}{25}G_{3,1}-\tfrac{178}{25}F_{0,4}
+
\\
&
\ \ \ \ \
-\tfrac{28}{75}F_{2,2}G_{4,0}+\tfrac{7}{25}G_{4,0}G_{3,1}+\tfrac{1}{25}F_{1,3}G_{3,1}-\tfrac{4}{75}F_{2,2}F_{1,3}-\tfrac{4}{25}F_{0,4}F_{1,3}
+
\\
&
\ \ \ \ \
-\tfrac{8}{25}F_{0,4}G_{4,0})x^5
+
\cdots
.
\endaligned\right.
\]
%%%%%%%%%%%%%%%%%%%%%%%%%%%%%%%%%%%%%%%%%%%%%%%%%%%%%%%%%%%%%%%%%%%%%%
\[
\def\arraystretch{1.25}
\begin{array}{ll}
e_1
&
\,:=\,
(1+\tfrac{2}{3}xF_{2,2}-\tfrac{1}{2}xG_{3,1}-\tfrac{18}{5}y+\tfrac{1}{5}yF_{1,3}+\tfrac{2}{5}yG_{4,0}-\tfrac{4}{3}uF_{2,2}+4uF_{0,4}+\tfrac{1}{2}uG_{3,1}+vF_{1,3})\partial_x
+
\\
&
\ \ \ \ \
+(-\tfrac{1}{5}xF_{1,3}-\tfrac{2}{5}xG_{4,0}+\tfrac{3}{5}x+\tfrac{2}{3}yF_{2,2}-\tfrac{1}{2}yG_{3,1}+\tfrac{1}{5}uF_{1,3}-\tfrac{3}{5}uG_{4,0}+\tfrac{27}{5}u-3uG_{0,4}
+
\\
&
\ \ \ \ \
+4vF_{0,4}+\tfrac{2}{3}vF_{2,2}+\tfrac{1}{2}vG_{3,1})\partial_y-(-y-\tfrac{4}{3}uF_{2,2}+uG_{3,1}-\tfrac{18}{5}v+\tfrac{1}{5}vF_{1,3}+\tfrac{2}{5}vG_{4,0})\partial_u
+
\\
&
\ \ \ \ \
+(2x+\tfrac{18}{5}u+\tfrac{4}{5}uF_{1,3}+\tfrac{8}{5}uG_{4,0}+\tfrac{4}{3}vF_{2,2}-vG_{3,1})\partial_v
,
\\
e_2
&
\,:=\,
(xF_{1,3}+xG_{0,4}+xG_{4,0}-6x+\tfrac{1}{10}yG_{3,1}+\tfrac{4}{5}yF_{0,4}-2uG_{0,4}-3uF_{1,3}+4vF_{0,4})\partial_x
+
\\
&
\ \ \ \ \
-(-1+\tfrac{1}{10}xG_{3,1}+\tfrac{4}{5}xF_{0,4}+3y-yF_{1,3}-yG_{0,4}-yG_{4,0}+\tfrac{9}{10}uG_{3,1}+\tfrac{56}{5}uF_{0,4}+2uF_{2,2}
+
\\
&
\ \ \ \ \
+2vG_{0,4})\partial_y-(-x-2uF_{1,3}-2uG_{0,4}-2uG_{4,0}+9u+\tfrac{1}{10}vG_{3,1}+\tfrac{4}{5}vF_{0,4})\partial_u+(-2y
+
\\
&
\ \ \ \ \
+\tfrac{2}{5}uG_{3,1}+\tfrac{16}{5}uF_{0,4}-6v+2vF_{1,3}+2vG_{0,4}+2vG_{4,0})\partial_v.
\end{array}
\]

%%%%%%%%%%%%%%%%%%%%%%%%%%%%%%%%%%%%%%%%%%%%%%%%%%%%%%%%%%%%%%%%%%%%%%
\noindent
Gr\"obner basis generators of 
moduli space core algebraic variety in 
$\R^6 \ni F_{0,4},F_{1,3},F_{2,2},G_{0,4},G_{3,1},G_{4,0}$:
\[
\aligned
\B_1 & := 288F_{0,4}G_{4,0}+27F_{1,3}G_{3,1}+120F_{2,2}G_{4,0}+45G_{0,4}G_{3,1}-45G_{3,1}G_{4,0}-1212F_{0,4}-590F_{2,2}
+
\\
&
\ \ \ \ \
+120G_{3,1}; 
\\
\B_2 & := 1656F_{0,4}G_{3,1}+216F_{1,3}G_{4,0}+240F_{2,2}G_{3,1}+2160G_{0,4}G_{4,0}+177G_{3,1}^2+672G_{4,0}^2-2144F_{1,3}
+
\\
&
\ \ \ \ \
-13960G_{0,4}-7536G_{4,0}+13392; 
\\
\B_3 & := -1224F_{1,3}G_{4,0}+204F_{2,2}G_{3,1}+828G_{0,4}^2-2304G_{0,4}G_{4,0}-129G_{3,1}^2-588G_{4,0}^2+8776F_{1,3}
+
\\
&
\ \ \ \ \
+10076G_{0,4}+4524G_{4,0}+6912; 
\\
\B_4 & := -576F_{0,4}G_{4,0}+108F_{2,2}G_{0,4}-204F_{2,2}G_{4,0}-117G_{0,4}G_{3,1}+117G_{3,1}G_{4,0}+3684F_{0,4}
+
\\
&
\ \ \ \ \
+1186F_{2,2}-303G_{3,1};
\\
\B_5 & := 1242F_{1,3}G_{0,4}-882F_{1,3}G_{4,0}-60F_{2,2}G_{3,1}-2196G_{0,4}G_{4,0}-165G_{3,1}^2-444G_{4,0}^2+6148F_{1,3}
+
\\
&
\ \ \ \ \
+15404G_{0,4}+4092G_{4,0}-864;
\\
\B_6 & := 648F_{0,4}G_{0,4}+792F_{0,4}G_{4,0}+240F_{2,2}G_{4,0}+171G_{0,4}G_{3,1}-171G_{3,1}G_{4,0}-5808F_{0,4}
+
\\
&
\ \ \ \ \
-1480F_{2,2}+474G_{3,1};
\\
\B_7 & := 2466F_{1,3}G_{4,0}+276F_{2,2}^2-825F_{2,2}G_{3,1}+2718G_{0,4}G_{4,0}+474G_{3,1}^2+2382G_{4,0}^2-20759F_{1,3}
+
\\
&
\ \ \ \ \
-11617G_{0,4}-17841G_{4,0}-32373;
\\
\B_8 & := 504F_{0,4}G_{4,0}+18F_{1,3}F_{2,2}+156F_{2,2}G_{4,0}+45G_{0,4}G_{3,1}-72G_{3,1}G_{4,0}-2436F_{0,4}-674F_{2,2}
+
\\
&
\ \ \ \ \
+273G_{3,1}; 
\\
\B_9 & := 1104F_{0,4}F_{2,2}-3852F_{1,3}G_{4,0}+918F_{2,2}G_{3,1}-5400G_{0,4}G_{4,0}-615G_{3,1}^2-2784G_{4,0}^2
+
\\
&
\ \ \ \ \
+29878F_{1,3}+29150G_{0,4}+21462G_{4,0}+32346;
\\
\B_{10} & := 621F_{1,3}^2-2016F_{1,3}G_{4,0}+750F_{2,2}G_{3,1}-3600G_{0,4}G_{4,0}-525G_{3,1}^2-2316G_{4,0}^2+23269F_{1,3}
+
\\
&
\ \ \ \ \
+19280G_{0,4}+17988G_{4,0}+25704;
\\
\B_{11} & := 324F_{0,4}F_{1,3}-1584F_{0,4}G_{4,0}-480F_{2,2}G_{4,0}-180G_{0,4}G_{3,1}+261G_{3,1}G_{4,0}+8916F_{0,4}
+
\\
&
\ \ \ \ \
+2330F_{2,2}-894G_{3,1};
\\
\B_{12} & := 39744F_{0,4}^2+49104F_{1,3}G_{4,0}-9840F_{2,2}G_{3,1}+77040G_{0,4}G_{4,0}+7509G_{3,1}^2+33168G_{4,0}^2
+
\\
&
\ \ \ \ \
-372556F_{1,3}-438260G_{0,4}-263724G_{4,0}-313092; \\
\B_{13} & := -1944F_{1,3}G_{4,0}^2+324F_{2,2}G_{3,1}G_{4,0}-4536G_{0,4}G_{4,0}^2-351G_{3,1}^2G_{4,0}-1080G_{4,0}^3
+
\\
&
\ \ \ \ \
+27072F_{1,3}G_{4,0}-1986F_{2,2}G_{3,1}+53118G_{0,4}G_{4,0}+2094G_{3,1}^2+16530G_{4,0}^2-94388F_{1,3}
+
\\
&
\ \ \ \ \
-150154G_{0,4}-57450G_{4,0}-34560;
\\
\B_{14} & := 25920F_{0,4}G_{4,0}^2+4320F_{2,2}G_{4,0}^2+4860G_{0,4}G_{3,1}G_{4,0}+135G_{3,1}^3-2700G_{3,1}G_{4,0}^2
+
\\
&
\ \ \ \ \
-382320F_{0,4}G_{4,0}-70584F_{2,2}G_{4,0}-47376G_{0,4}G_{3,1}+49284G_{3,1}G_{4,0}+1298040F_{0,4}
+
\\
&
\ \ \ \ \
+251516F_{2,2}-146766G_{3,1};
\\
\B_{15} & := 6210G_{0,4}G_{3,1}^2-298080G_{0,4}G_{4,0}^2-22770G_{3,1}^2G_{4,0}+33120G_{4,0}^3-440376F_{1,3}G_{4,0}
+
\\
&
\ \ \ \ \
+141108F_{2,2}G_{3,1}+3240888G_{0,4}G_{4,0}+56847G_{3,1}^2-326472G_{4,0}^2+2925104F_{1,3}
+
\\
&
\ \ \ \ \
-8290352G_{0,4}-1614984G_{4,0}+16834608;
\\
\B_{16} & := 110160F_{0,4}G_{4,0}^2+405F_{2,2}G_{3,1}^2+24840F_{2,2}G_{4,0}^2+14175G_{0,4}G_{3,1}G_{4,0}-14175G_{3,1}G_{4,0}^2
+
\\
&
\ \ \ \ \
-1539900F_{0,4}G_{4,0}-364398F_{2,2}G_{4,0}-151317G_{0,4}G_{3,1}+213498G_{3,1}G_{4,0}
+
\\
&
\ \ \ \ \
+4969920F_{0,4}+1216192F_{2,2}-542742G_{3,1};
\\
\endaligned
\]
%%%%%%%%%%%%%%%%%%%%%%%%%%%%%%%%%%%%%%%%%%%%%%%%%%%%%%%%%%%%%%%%%%%%%%
\[
\footnotesize
\def\arraystretch{1.25}
\begin{array}{c|cc}
{} & e_1 & e_2 
\\
\hline
e_1 & 0 & \rotatebox[origin=c]{0}{
\begin{tabular}{p{7cm}}
$(-\tfrac{12}{5}+\tfrac{4}{5}F_{1,3}+G_{0,4}+\tfrac{3}{5}G_{4,0})e_1+(-\tfrac{2}{3}F_{2,2}+\tfrac{2}{5}G_{3,1}-\tfrac{4}{5}F_{0,4})e_2$
\end{tabular}}
\\
e_2 & \rotatebox[origin=c]{0}{
\begin{tabular}{p{7cm}}
$-(-\tfrac{12}{5}+\tfrac{4}{5}F_{1,3}+G_{0,4}+\tfrac{3}{5}G_{4,0})e_1-(-\tfrac{2}{3}F_{2,2}+\tfrac{2}{5}G_{3,1}-\tfrac{4}{5}F_{0,4})e_2$
\end{tabular}} & 0
\end{array}
\]

%%%%%%%%%%%%%%%%%%%%%%%%%%%%%%%%%%%%%%%%%%%%%%%%%%%%%%%%%%%%%%%%%%%%%%
%%%%%%%%%%%%%%%%%%%%%%%%%%%%%%%%%%%%%%%%%%%%%%%%%%%%%%%%%%%%%%%%%%%%%%
%%%%%%%%%%%%%%%%%%%%%%%%%%%%%%%%%%%%%%%%%%%%%%%%%%%%%%%%%%%%%%%%%%%%%%

\[
\text{\bf Model 2g3c4a}
\ \ \ \ \
\left\{
\aligned
u
&
\,=\,
xy-\tfrac{1}{6}y^4+x^2y^2-\tfrac{1}{6}x^4+F_{0,5}y^5-\tfrac{5}{2}G_{0,5}xy^4-10F_{0,5}x^2y^3+5G_{0,5}x^3y^2
+
\\
&
\ \ \ \ \
+5F_{0,5}x^4y-\tfrac{1}{2}G_{0,5}x^5+(\tfrac{15}{8}G_{0,5}^2-\tfrac{15}{2}F_{0,5}^2)y^6+(\tfrac{-2}{3}+45G_{0,5}F_{0,5})xy^5
+
\\
&
\ \ \ \ \
+(-\tfrac{225}{8}G_{0,5}^2+\tfrac{225}{2}F_{0,5}^2)x^2y^4+(\tfrac{20}{9}-150G_{0,5}F_{0,5})x^3y^3+(\tfrac{225}{8}G_{0,5}^2
+
\\
&
\ \ \ \ \
-\tfrac{225}{2}F_{0,5}^2)x^4y^2+(\tfrac{-2}{3}+45G_{0,5}F_{0,5})x^5y+(-\tfrac{15}{8}G_{0,5}^2+\tfrac{15}{2}F_{0,5}^2)x^6
+
\\
&
\ \ \ \ \
+(\tfrac{450}{7}F_{0,5}^3-\tfrac{675}{14}G_{0,5}^2F_{0,5}+\tfrac{17}{21}G_{0,5})y^7+(\tfrac{225}{4}G_{0,5}^3-675F_{0,5}^2G_{0,5}
+
\\
&
\ \ \ \ \
+\tfrac{34}{3}F_{0,5})xy^6+(-1350F_{0,5}^3+\tfrac{2025}{2}G_{0,5}^2F_{0,5}-17G_{0,5})x^2y^5+(-\tfrac{170}{3}F_{0,5}
+
\\
&
\ \ \ \ \
-\tfrac{1125}{4}G_{0,5}^3+3375F_{0,5}^2G_{0,5})x^3y^4+(\tfrac{85}{3}G_{0,5}-\tfrac{3375}{2}G_{0,5}^2F_{0,5}+2250F_{0,5}^3)x^4y^3
+
\\
&
\ \ \ \ \
+(34F_{0,5}+\tfrac{675}{4}G_{0,5}^3-2025F_{0,5}^2G_{0,5})x^5y^2+(-450F_{0,5}^3+\tfrac{675}{2}G_{0,5}^2F_{0,5}
+
\\
&
\ \ \ \ \
-\tfrac{17}{3}G_{0,5})x^6y+(-\tfrac{225}{28}G_{0,5}^3+\tfrac{675}{7}F_{0,5}^2G_{0,5}-\tfrac{34}{21}F_{0,5})x^7
+
\cdots
,
\\
v
&
\,=\,
-y^2+x^2-\tfrac{4}{3}xy^3+\tfrac{4}{3}x^3y+G_{0,5}y^5+10F_{0,5}xy^4-10G_{0,5}x^2y^3
+
\\
&
\ \ \ \ \
-20F_{0,5}x^3y^2+5G_{0,5}x^4y+2F_{0,5}x^5+(\tfrac{2}{9}-15G_{0,5}F_{0,5})y^6+(\tfrac{45}{2}G_{0,5}^2
+
\\
&
\ \ \ \ \
-90F_{0,5}^2)xy^5+(\tfrac{-10}{3}+225G_{0,5}F_{0,5})x^2y^4+(300F_{0,5}^2-75G_{0,5}^2)x^3y^3
+
\\
&
\ \ \ \ \
+(\tfrac{10}{3}-225G_{0,5}F_{0,5})x^4y^2+(\tfrac{45}{2}G_{0,5}^2-90F_{0,5}^2)x^5y+(\tfrac{-2}{9}+15G_{0,5}F_{0,5})x^6
+
\\
&
\ \ \ \ \
+(\tfrac{1350}{7}F_{0,5}^2G_{0,5}-\tfrac{225}{14}G_{0,5}^3-\tfrac{68}{21}F_{0,5})y^7+(\tfrac{34}{3}G_{0,5}-675G_{0,5}^2F_{0,5}
+
\\
&
\ \ \ \ \
+900F_{0,5}^3)xy^6+(-4050F_{0,5}^2G_{0,5}+\tfrac{675}{2}G_{0,5}^3+68F_{0,5})x^2y^5+(-4500F_{0,5}^3
+
\\
&
\ \ \ \ \
+3375G_{0,5}^2F_{0,5}-\tfrac{170}{3}G_{0,5})x^3y^4+(-\tfrac{340}{3}F_{0,5}+6750F_{0,5}^2G_{0,5}
+
\\
&
\ \ \ \ \
-\tfrac{1125}{2}G_{0,5}^3)x^4y^3+(2700F_{0,5}^3-2025G_{0,5}^2F_{0,5}+34G_{0,5})x^5y^2
+
\\
&
\ \ \ \ \
+(-1350F_{0,5}^2G_{0,5}+\tfrac{225}{2}G_{0,5}^3+\tfrac{68}{3}F_{0,5})x^6y+(-\tfrac{34}{21}G_{0,5}+\tfrac{675}{7}G_{0,5}^2F_{0,5}
+
\\
&
\ \ \ \ \
-\tfrac{900}{7}F_{0,5}^3)x^7
+
\cdots
.
\endaligned\right.
\]
%%%%%%%%%%%%%%%%%%%%%%%%%%%%%%%%%%%%%%%%%%%%%%%%%%%%%%%%%%%%%%%%%%%%%%
\[
\def\arraystretch{1.25}
\begin{array}{ll}
e_1
&
\,:=\,
(1-\tfrac{15}{2}xG_{0,5}+15yF_{0,5}-\tfrac{4}{3}u)\partial_x-(15xF_{0,5}+\tfrac{15}{2}yG_{0,5}-\tfrac{2}{3}v)\partial_y+
\\
&
\ \ \ \ \
-(15F_{0,5}v+15G_{0,5}u-y)\partial_u+(60F_{0,5}u-15G_{0,5}v+2x)\partial_v
,\\
e_2
&
\,:=\,
(15xF_{0,5}+\tfrac{15}{2}yG_{0,5}-\tfrac{2}{3}v)\partial_x+(1-\tfrac{15}{2}xG_{0,5}+15yF_{0,5}-\tfrac{4}{3}u)\partial_y
+
\\
&
\ \ \ \ \
+(x+30F_{0,5}u-\tfrac{15}{2}G_{0,5}v)\partial_u+(30F_{0,5}v+30G_{0,5}u-2y)\partial_v
,
\end{array}
\]
for any value of $F_{0,5},G_{0,5}$.
%%%%%%%%%%%%%%%%%%%%%%%%%%%%%%%%%%%%%%%%%%%%%%%%%%%%%%%%%%%%%%%%%%%%%%
\[
\footnotesize
\def\arraystretch{1.25}
\begin{array}{c|cc}
{} & e_1 & e_2 
\\
\hline
e_1 & 0 & 0
\\
e_2 & 0 & 0
\end{array}
\]

%%%%%%%%%%%%%%%%%%%%%%%%%%%%%%%%%%%%%%%%%%%%%%%%%%%%%%%%%%%%%%%%%%%%%%
%%%%%%%%%%%%%%%%%%%%%%%%%%%%%%%%%%%%%%%%%%%%%%%%%%%%%%%%%%%%%%%%%%%%%%
%%%%%%%%%%%%%%%%%%%%%%%%%%%%%%%%%%%%%%%%%%%%%%%%%%%%%%%%%%%%%%%%%%%%%%

\[
\text{\bf Model 2g3c4b}
\ \ \ \ \
\left\{
\aligned
u
&
\,=\,
xy
,
\\
v
&
\,=\,
x^2-y^2
.
\endaligned\right.
\]
%%%%%%%%%%%%%%%%%%%%%%%%%%%%%%%%%%%%%%%%%%%%%%%%%%%%%%%%%%%%%%%%%%%%%%
\[
\def\arraystretch{1.25}
\begin{array}{llll}
e_1
\,:=\,
\partial_x+y\partial_u+2x\partial_v
, &
e_2
\,:=\,
\partial_y+x\partial_u-2y\partial_v
, &
& 
\\
e_3
\,:=\,
2u\partial_u+2v\partial_v+x\partial_x+y\partial_y
, &
e_4
\,:=\,
v\partial_u-4u\partial_v-y\partial_x+x\partial_y
.
\end{array}
\]

%%%%%%%%%%%%%%%%%%%%%%%%%%%%%%%%%%%%%%%%%%%%%%%%%%%%%%%%%%%%%%%%%%%%%%
\[
\footnotesize
\def\arraystretch{1.25}
\begin{array}{c|cccc}
{} & e_1 & e_2 & e_3 & e_4
\\
\hline
e_1 & 0 & 0 & e_1 & e_2
\\
e_2 & 0 & 0 & e_2 & e_1
\\
e_3 & -e_1 & -e_2 & 0 & 0
\\
e_4 & -e_2 & -e_1 & 0 & 0
\end{array}
\]

%%%%%%%%%%%%%%%%%%%%%%%%%%%%%%%%%%%%%%%%%%%%%%%%%%%%%%%%%%%%%%%%%%%%%%
%%%%%%%%%%%%%%%%%%%%%%%%%%%%%%%%%%%%%%%%%%%%%%%%%%%%%%%%%%%%%%%%%%%%%%
%%%%%%%%%%%%%%%%%%%%%%%%%%%%%%%%%%%%%%%%%%%%%%%%%%%%%%%%%%%%%%%%%%%%%%

%%%%%%%%%%%%%%%%%%%%%%%%%%%%%%%%%%%%%%%%%%%%%%%%%%%%%%%%%%%%%%%%%%%%%%
%%%%%%%%%%%%%%%%%%%%%%%%%%%%%%%%%%%%%%%%%%%%%%%%%%%%%%%%%%%%%%%%%%%%%%
%%%%%%%%%%%%%%%%%%%%%%%%%%%%%%%%%%%%%%%%%%%%%%%%%%%%%%%%%%%%%%%%%%%%%%

\vfill\end{document}